\documentclass[11pt,reqno]{amsart}
\usepackage{amsmath,amsthm,amssymb,mathrsfs,stmaryrd,color}
\usepackage[all]{xy}
\usepackage{url}
\usepackage{slashed}
\usepackage{geometry}

\usepackage[utf8]{inputenc}
\usepackage[T1]{fontenc}

\usepackage{relsize}
\usepackage[bbgreekl]{mathbbol}
\usepackage{amsfonts}

\DeclareMathOperator{\hc}{hc} 
\DeclareMathOperator{\HKR}{HKR}
\DeclareMathOperator{\qOmega}{q\Omega}
\DeclareMathOperator{\CP}{\mathbf{CP}}
\DeclareMathOperator{\THH}{THH}
\DeclareMathOperator{\HH}{HH}
\DeclareMathOperator{\HC}{HC}
\DeclareMathOperator{\HP}{HP}
\DeclareMathOperator{\TC}{TC}
\DeclareMathOperator{\TP}{TP}
\DeclareMathOperator{\slashOmega}{\widetilde{p}\Omega}
\DeclareMathOperator{\slashdR}{\widetilde{p}dR}
\DeclareMathOperator{\Perf}{Perf}
\DeclareMathOperator{\Syn}{syn}
\DeclareMathOperator{\RHom}{RHom}
\DeclareMathOperator{\Hodge}{Hod}
\DeclareMathOperator{\mathet}{\acute{e}t}
\DeclareMathOperator{\DHod}{\slashed{D}}
\DeclareMathOperator{\RGamma}{R \Gamma}
\DeclareMathOperator{\Cart}{Cart}
\DeclareMathOperator{\qdR}{qdR}
\DeclareMathOperator{\WCart}{WCart}
\newcommand{\slashp}{\tilde{p}}
\DeclareMathOperator{\dlog}{dlog}
\DeclareMathOperator{\calN}{\mathcal{N}}
\DeclareMathOperator{\cyc}{cyc}
\DeclareMathOperator{\im}{im}

\DeclareMathOperator{\Vect}{Vect}
\DeclareMathOperator{\Tate}{T}
\DeclareMathOperator{\mot}{M}
\DeclareMathOperator{\BGL}{BGL}
\DeclareMathOperator{\BG}{BG}
\DeclareMathOperator{\BGm}{B\mathbf{G}_{m}}
\DeclareMathOperator{\BP}{BP}
\DeclareMathOperator{\BL}{BL}
\DeclareMathOperator{\GL}{GL}
\DeclareMathOperator{\Gal}{Gal}
\DeclareMathOperator{\SSet}{\mathcal{S}}
\DeclareMathOperator{\cofib}{cofib}
\DeclareSymbolFontAlphabet{\mathbb}{AMSb}
\DeclareSymbolFontAlphabet{\mathbbl}{bbold}
\newcommand{\Prism}{\mathbbl{\Delta}}
\DeclareMathOperator{\End}{End}

\DeclareMathOperator{\fppf}{fppf}
\DeclareMathOperator{\Z}{\mathbf{Z}}
\DeclareMathOperator{\F}{\mathbf{F}}
\DeclareMathOperator{\fib}{fib}

\newcommand{\SynSheaf}[2]{\RGamma_{\mathrm{Syn}}( \Spec(#2), \Z_p(#1) )}
\DeclareMathOperator{\DFilt}{\mathcal{DF}}
\DeclareMathOperator{\DFiltI}{ \widehat{\mathcal{DF}}}
\DeclareMathOperator{\DFiltIComp}{ {\widehat{\mathcal{DF}}}^{\mathrm{c}}}
\DeclareMathOperator{\DFiltComp}{ \mathcal{DF}^{\mathrm{c}}}
\DeclareMathOperator{\Beil}{Beil}
\DeclareMathOperator{\conj}{conj}
\DeclareMathOperator{\id}{id}
\DeclareMathOperator{\Fun}{Fun}
\DeclareMathOperator{\Sphere}{\mathbb{S}}
\DeclareMathOperator{\op}{op}

\DeclareMathOperator{\Nyg}{N}
\DeclareMathOperator{\Ext}{Ext}
\DeclareMathOperator{\Sym}{Sym}

\DeclareMathOperator{\calE}{\mathcal{E}}
\DeclareMathOperator{\calC}{\mathcal{C}}
\DeclareMathOperator{\calD}{\mathcal{D}}
\DeclareMathOperator{\coker}{coker}
\DeclareMathOperator{\Poly}{Poly}
\DeclareMathOperator{\CAlg}{CAlg}

\DeclareMathOperator{\anim}{an}
\DeclareMathOperator{\Aut}{Aut}
\DeclareMathOperator{\Fil}{Fil}

\DeclareMathOperator{\Tor}{Tor}
\DeclareMathOperator{\gr}{gr}
\DeclareMathOperator{\glo}{gl}
\DeclareMathOperator{\Tot}{Tot}
\DeclareMathOperator{\calI}{\mathcal{I}}
\DeclareMathOperator{\calO}{\mathcal{O}}
\DeclareMathOperator{\Crys}{Crys}
\DeclareMathOperator{\Pic}{Pic}

\DeclareMathOperator{\qrsp}{qrsp}
\DeclareMathOperator{\Set}{Set}
\DeclareMathOperator{\QCoh}{QCoh}
\DeclareMathOperator{\calU}{\mathcal{U}}

\DeclareMathOperator{\calL}{\mathcal{L}}
\DeclareMathOperator{\Q}{\mathbf{Q}}

\usepackage{enumitem}
\setlist[enumerate]{itemsep=2pt,parsep=2pt,before={\parskip=2pt}}

\usepackage[colorlinks=true,hyperindex, linkcolor=magenta, pagebackref=false, citecolor=cyan,pdfpagelabels]{hyperref}

\newcommand{\cosimp}[3]{\xymatrix@R=50pt@C=50pt@1{#1 \ar@<.4ex>[r] \ar@<-.4ex>[r] & {\ }#2 \ar@<0.8ex>[r] \ar[r] \ar@<-.8ex>[r] & {\ } #3 \ar@<1.2ex>[r] \ar@<.4ex>[r] \ar@<-.4ex>[r] \ar@<-1.2ex>[r] & \cdots }}

\newcommand{\colim}{\mathop{\mathrm{colim}}}

\newcommand{\adjunction}[4]{\xymatrix@R=50pt@C=50pt@1{#1{\ } \ar@<0.3ex>[r]^-{ {\scriptstyle #2}} & {\ } #3 \ar@<0.3ex>[l]^{ {\scriptstyle #4}}}}

\begin{document}

\newtheorem{theorem}{Theorem}[subsection]
\newtheorem*{theorem*}{Theorem}
\newtheorem*{definition*}{Definition}
\newtheorem{proposition}[theorem]{Proposition}
\newtheorem{lemma}[theorem]{Lemma}
\newtheorem{corollary}[theorem]{Corollary}
\newtheorem{conjecture}[theorem]{Conjecture}

\theoremstyle{definition}
\newtheorem{variant}[theorem]{Variant}
\newtheorem{definition}[theorem]{Definition}
\newtheorem{question}[theorem]{Question}
\newtheorem{remark}[theorem]{Remark}
\newtheorem{warning}[theorem]{Warning}
\newtheorem{example}[theorem]{Example}
\newtheorem{notation}[theorem]{Notation}
\newtheorem{convention}[theorem]{Convention}
\newtheorem{construction}[theorem]{Construction}
\newtheorem{claim}[theorem]{Claim}
\newtheorem{assumption}[theorem]{Assumption}

\newcommand{\Qc}{q-\mathrm{crys}}

\newcommand{\Shv}{\mathrm{Shv}}
\newcommand{\et}{{\acute{e}t}}
\newcommand{\crys}{\mathrm{crys}}
\renewcommand{\inf}{\mathrm{inf}}
\newcommand{\Hom}{\mathrm{Hom}}
\newcommand{\Sch}{\mathrm{Sch}}
\newcommand{\Spf}{\mathrm{Spf}}
\newcommand{\Spa}{\mathrm{Spa}}
\newcommand{\Spec}{\mathrm{Spec}}
\newcommand{\perf}{\mathrm{perf}}
\newcommand{\QSyn}{\mathrm{QSyn}}
\newcommand{\perfd}{\mathrm{perfd}}
\newcommand{\arc}{{\rm arc}}

\newcommand{\rad}{\mathrm{rad}}

\newcommand{\psh}{\mathrm{PShv}}
\newcommand{\scr}{\mathrm{sCAlg}}
\newcommand{\HT}{\mathrm{HT}}
\newcommand{\dR}{\mathrm{dR}}
\newcommand{\LdR}{\mathrm{LdR}}

\setcounter{tocdepth}{2}

\newcommand{\comment}[1]{\textcolor{red}{\footnotesize #1}}

\title{Absolute Prismatic Cohomology}
\author{Bhargav Bhatt}
\author{Jacob Lurie}
\maketitle

\begin{abstract}
The goal of this paper is to study the absolute prismatic cohomology of $p$-adic formal schemes. We do so by recasting the notion of a prismatic crystal on $\mathrm{Spf}(\mathbf{Z}_p)$ in terms of quasicoherent sheaves on a geometric object we call the Cartier-Witt stack. 
\end{abstract}

\tableofcontents

\newpage

\newpage \section{Introduction}\label{section:introduction}

\subsection{Goal}
Let $K$ be a field of characteristic zero and let $X_K$ be a smooth and proper variety over $K$. One can then consider two different cohomological
invariants of $X$:
\begin{itemize}
\item The {\it algebraic de Rham cohomology} $\mathrm{H}^{\ast}_{\dR}( X_K / K)$, defined as the (hyper)cohomology of $X_K$ with coefficients in
the algebraic de Rham complex
$$ \Omega^{0}_{ X_K / K } \rightarrow \Omega^{1}_{X_K /K} \rightarrow \Omega^{2}_{X_K /K} \rightarrow \cdots.$$

\item The {\it $p$-adic \'{e}tale cohomology} $\mathrm{H}^{\ast}_{\mathet}( X_{ \overline{K} }, \Z_p)$; here $p$ denotes
a fixed prime number, $\overline{K}$ denotes an algebraic closure of $K$, and $X_{ \overline{K} }$ denotes the geometric fiber
$\Spec( \overline{K} ) \times_{ \Spec(K) } X_K$.
\end{itemize}
When $K = \mathbf{C}$ is the field of complex numbers, these invariants are closely related: there are canonical isomorphisms
$$ \mathrm{H}^{\ast}_{\dR}( X_{\mathbf{C}} / \mathbf{C} ) \simeq \mathbf{C} \otimes \mathrm{H}^{\ast}_{\mathrm{sing} }( X(\mathbf{C}), \Z) \quad \quad \mathrm{H}^{\ast}_{\mathet}( X_{\mathbf{C}}, \Z_p ) \simeq \Z_p \otimes \mathrm{H}^{\ast}_{\mathrm{sing} }( X(\mathbf{C}), \Z),$$
where $\mathrm{H}^{\ast}_{\mathrm{sing} }( X(\mathbf{C}), \Z)$ denotes the singular cohomology of the complex manifold $X(\mathbf{C})$ of $\mathbf{C}$-valued points of $X_{\mathbf{C}}$. Note that the definition of $\mathrm{H}^{\ast}_{\mathrm{sing} }( X(\mathbf{C}), \Z)$ depends in an essential way on the (archimedean) topology on 
the field of complex numbers.

A central goal of $p$-adic Hodge theory is to obtain an analogous picture in the case where the field $K$ is equipped with a {\em nonarchimedean} topology.
Suppose that $K$ is complete with respect to a $p$-adic absolute value, and let us also assume that the variety $X_{K}$ has good reduction: that is, it is given as the generic fiber of a scheme $X$ which is smooth and proper over the valuation ring $\calO_{K} \subset K$. In this case, the algebraic de Rham cohomology groups of $X_K$ admit an integral refinement: we have
an isomorphism
$$\mathrm{H}^{\ast}_{\dR}( X_K / K) \simeq \mathrm{H}^{\ast}_{\dR}( X / \calO_K )[1/p],$$
where $\mathrm{H}^{\ast}_{\dR}( X / \calO_K )$ denotes the cohomology of the cochain complex
$$ \RGamma_{\dR}( X / \calO_K ) =
\RGamma(X, \Omega^{0}_{ X / \calO_K} \rightarrow \Omega^{1}_{ X/ \calO_K} \rightarrow \cdots ).$$
It is then natural to ask the following:

\begin{question}\label{question:p-adic-Hodge}
What is the relationship between the cochain complexes $\RGamma_{\dR}( X / \calO_K )$ and $\RGamma_{ \mathet}( X_{ \overline{K} }, \Z_p)$?
Is there some more fundamental invariant of $X$ recovering both these complexes, and yielding new structures on them, or new relationships between them?
\end{question}

A rough goal of this paper is to propose an answer this question when $K=\mathbf{Q}_p$ by elaborating on the notion of absolute prismatic cohomology that was implicitly introduced (but not seriously pursued) in the paper \cite{prisms}. 

\subsection{Some recent history}
Recall that \cite{BMS1} gave an answer to Question \ref{question:p-adic-Hodge} in the case where
$K = \overline{K}$ is an algebraically closed field. Let $A_{\mathrm{inf}} = W( \calO_{K}^{\flat} )$ denote Fontaine's
``infinitesimal'' period ring, so that the tautological surjection $\calO_{K}^{\flat} \twoheadrightarrow \calO_{K} / p \calO_{K}$
lifts uniquely to a surjective ring homomorphism $\theta: A_{\mathrm{inf}} \twoheadrightarrow \calO_{K}$. In this situation, we have:

\begin{theorem}[{\cite[Theorem~1.10]{BMS1}}]\label{theorem:BMS1-Summarized}
Let $K$ be an algebraically closed field which is complete with respect to a $p$-adic absolute value.
To every smooth and proper $\calO_{K}$-scheme $X$, one can associate perfect complex of $A_{\mathrm{inf}}$-modules
$\RGamma_{ A_{\mathrm{inf}}}(X)$ with the following features:
\begin{itemize}
\item There is a canonical isomorphism
$$ \RGamma_{ \dR}( X / \calO_K) \simeq \calO_{K} \otimes^{L}_{ A_{\mathrm{inf}} } \RGamma_{ A_{\mathrm{inf}}}( X )$$
in the derived category of $\calO_K$-modules. That is, the algebraic de Rham cohomology of $\mathfrak{X}$ can be recovered from $\RGamma_{ A_{\mathrm{inf}}}(X )$ by extending scalars along $\theta$.

\item Let $\xi$ denote a generator of the (principal) ideal $\ker(\theta) \subset A_{\mathrm{inf}}$. Then there is a canonical isomorphism
$$ \RGamma_{\mathet}( X_{K}, \Z_p) \simeq (\RGamma_{ A_{\mathrm{inf}} }(X)[ 1 / \xi ]^{\wedge})^{\varphi=1}$$
in the derived category of $\Z_p$-modules. That is, the $p$-adic \'{e}tale cohomology of $X_K$ can be recovered from the $p$-adic completion of 
$\RGamma_{ A_{\mathrm{inf}}}( X )[ 1 / \xi ]$ by extracting the fixed points of Frobenius.
\end{itemize}
\end{theorem}

The original definition of the complex $\RGamma_{ A_{\mathrm{inf}}}(X)$ uses perfectoid geometry. To carry out
the constructions of \cite{BMS1}, it is not necessary to assume that $K$ is algebraically closed: it is enough to assume that
$K$ is a perfectoid field containing a system of $p^{n}$th roots of unity. However, this a very strong assumption which
is not satisfied in many cases of arithmetic interest (for example, finite extensions of $\Q_p$ are never perfectoid). To address this point, \cite{prisms} introduced the theory of {\em prismatic cohomology}, which generalizes the construction $X \mapsto \RGamma_{ A_{\mathrm{inf}}}(X )$ to geometric objects which are defined over non-perfectoid ground rings. The starting point of this theory is the following:

\begin{definition}[Prisms: Torsion-Free Case]\label{definition:prism-preliminary}
Let $A$ be a torsion-free commutative ring, let $\varphi_A: A \rightarrow A$ be a ring homomorphism
satisfying $\varphi_{A}(x) \equiv x^{p} \pmod{p}$, and let $I \subseteq A$ be an invertible ideal.
We say that the pair $(A,I)$ is a {\it prism} if $A$ is both $p$-adically and $I$-adically complete,
and the collection of elements
$$ \{ \frac{ \varphi_{A}(x) - x^p}{p}  \}_{x \in I}$$
generate the unit ideal of $A$. In this case, we will generally write $\overline{A}$ for the quotient ring $A/I$.
We say that the prism $(A,I)$ is {\it perfect} if $\varphi_{A}$ is an isomorphism,
and {\it bounded} if the ring $\overline{A}$ has bounded $p$-torsion.
\end{definition}

\begin{warning}
Definition \ref{definition:prism-preliminary} is a special case of the general definition of prism (see Definition \ref{definition:prism}), which does not
require the commutative ring $A$ to be torsion-free.
\end{warning}

\begin{example}
Let $K$ be a perfectoid field, let $A_{\mathrm{inf} } = W( \calO_{K}^{\flat} )$ be Fontaine's infinitesimal period ring, and let
$\theta: A_{\mathrm{inf}} \twoheadrightarrow \calO_{K}$ be the tautological map. Then the pair $( A_{\mathrm{inf}}, \ker(\theta) )$ is a perfect prism.
\end{example}

Let $(A,I)$ be a bounded prism and let $X$ be a smooth and proper $\overline{A}$-scheme. In this setup, the paper \cite{prisms} defined the relative prismatic site $(X/A)_\Prism$ and used it construct a complex of $A$-modules $\RGamma_{\Prism}(X / A)$, whose cohomology we will denote by $\mathrm{H}^{\ast}_{\Prism}( X / A)$ and refer to as the {\it prismatic cohomology of $\mathfrak{X}$ relative to $A$}. In the special case where $(A,I) = ( A_{\mathrm{inf}}, \ker(\theta) )$, this recovers the complex $\RGamma_{ A_{\mathrm{inf}} }( X )$ of \cite{BMS1} (by \cite[Theorem~17.2]{prisms}). Other cohomology theories recovered by this construction include:

\begin{example}[Crystalline cohomology]
\label{example:crystalline-from-relative-prismatic}
Let $k$ be a perfect field of characteristic $p$,  let $(A,I) = (W(k),(p))$, and let $X$ be a smooth and proper $k$-scheme. Then the prismatic cohomology $\mathrm{H}^{\ast}_{\Prism}(X/A)$ can be identified with the crystalline cohomology $\mathrm{H}^{\ast}_{\crys}(X/W(k))$  (see \cite[Theorem~5.2]{prisms}, or
Theorem \ref{theorem:crystalline-comparison} for a more general statement). If $Y$ is a proper smooth $W(k)$-scheme lifting $X$ (i.e., $X = \Spec(k) \times_{\Spec(W(k))} Y$), then one can also identify $\mathrm{H}^{\ast}_{\crys}(X/W(k))$ (and thus $\mathrm{H}^*_{\Prism}(X/A)$) with the algebraic de Rham cohomology
$\mathrm{H}^{\ast}_{\dR}( Y/W(k) )$ by a result of Berthelot \cite{BerthelotCrysCoh}, which we revisit from a ``prismatic'' perspective
in \S \ref{subsection:absolute-de-Rham-comparison}.
\end{example}


The examples discussed thus far concern perfect prisms. However, the general theory of prismatic cohomology offers the flexibility of working over other (not necessarily perfect) prisms. An important such example is:

\begin{example}[$q$-de Rham cohomology]\label{example:q-de-Rham-from-relative}
Let $(A, I)$ denote the $q$-de Rham prism $( \Z_p[[q-1]], ( 1+ q + \cdots + q^{p-1} ) )$ (see Example \ref{example:q-prism}),
so that $\overline{A} = A/I$ is the cyclotomic number ring $\Z_p[ \zeta_p ]$. If $X$ is a proper smooth $\mathbf{Z}_p$-scheme, then the cochain complex $\RGamma_{\qdR}(X) :=\RGamma_{ \Prism}( X_{\overline{A}} / A )$ can be regarded as a $q$-deformation of the algebraic de Rham complex $\RGamma_{\dR}( X/ \Z_p )$, which
computes the {\it $q$-de Rham cohomology} of the scheme $X$ (see \cite[\S~16]{prisms}).
\end{example}

Given a smooth and proper $\Z_p$-scheme $X$ and a bounded prism $(A,I) $, we obtain a smooth and proper $\overline{A}$-scheme $X_{\overline{A}} := \Spec(\overline{A}) \times_{ \Spec(\Z_p) } X$ via base change.  The preceding discussion suggests that the construction carrying a bounded prism $(A,I)$ to the cochain complex $\RGamma_{ \Prism}( X_{\overline{A}} / A )$ of $A$-modules is an important invariant of $X$: it determines both the de Rham cohomology of $X$ (using Example \ref{example:crystalline-from-relative-prismatic}) and the \'{e}tale cohomology of the geometric generic fiber of $X$ (using Theorem \ref{theorem:BMS1-Summarized}). Thus, this invariant can be regarded as an answer to  Question \ref{question:p-adic-Hodge} when $K=\mathbf{Q}_p$. The goal of this paper is to study the finer structure of this invariant and use it to construct more concrete invariants of $X$. (In fact, the results in the paper itself usually have much milder hypotheses on $X$; in the introduction, we often stick to smooth and proper $X$'s to obtain cleaner statements.)

\subsection{Absolute prismatic cohomology and the Cartier-Witt stack}
\label{ss:APCCW}

Recall that the absolute prismatic site $\mathrm{Spf}(\mathbf{Z}_p)_\Prism$ of $\mathrm{Spf}(\mathbf{Z}_p)$ is simply another name for the category of all bounded prisms (see Definition~\ref{definition:absolute-prismatic-site} for the absolute prismatic site in general).  For a fixed (smooth and proper) $\Z_p$-scheme $X$, the construction
$$ (A,I) \mapsto \RGamma_{ \Prism}( X_{\overline{A}} / A )$$
discussed above  is an example of a {\it prismatic crystal on  $\mathrm{Spf}(\mathbf{Z}_p)$}: it associates to each bounded prism $(A,I)$ a perfect complex of $A$-modules, and
to each morphism of bounded prisms $f: (A,I) \rightarrow (B,J)$ a quasi-isomorphism
$$ B \otimes^{L}_{A} \RGamma_{ \Prism}( X_{\overline{A}} / A ) \simeq \RGamma_{ \Prism}( X_{\overline{B}} / B ),$$
whose formation is compatible with composition (up to coherent homotopy).  Our first goal in this paper is to study such objects geometrically, i.e., we shall realize prismatic crystals on $\mathrm{Spf}(\mathbf{Z}_p)$ as quasi-coherent sheaves on an algebro-geometric object.

Let $R$ be a commutative ring in which $p$ is nilpotent. 
In \S \ref{section:cartier-witt}, we will introduce a groupoid $\WCart(R)$ of {\it Cartier-Witt divisors for $R$} (Definition \ref{definition:Cartier-Witt-divisor}), whose
objects are closely related to prism structures on the ring of Witt vectors $W(R)$ (Remark \ref{remark:Cartier-Witt-divisor-prism}). The assignment $R \mapsto \WCart(R)$ can be regarded as a functor, which we will denote by $\WCart$ and refer to as the {\it Cartier-Witt stack}. The stack $\WCart$ is not far from being algebraic:
it can be described as the (stack-theoretic) quotient of an (affine) formal scheme by the action of an (affine) group scheme (see Proposition \ref{proposition:presentation-as-quotient}). To every prism $(A,I)$, one can associate a map $\rho_{A}: \Spf(A) \rightarrow \WCart$ (Construction \ref{construction:point-of-prismatic-stack}).
Consequently, every perfect complex $\mathscr{F}$ on $\WCart$ determines a prismatic crystal, given by the construction
$(A,I) \mapsto \rho^{\ast}_{A}(\mathscr{F})$. We will see that this construction induces an equivalence
\begin{equation}
\label{eq:GeometrizePrismCrys}
 \{ \text{Perfect complexes on $\WCart$} \} \simeq \{ \text{Perfect Prismatic Crystals on } \mathrm{Spf}(\mathbf{Z}_p)  \}
 \end{equation}
(see Proposition \ref{proposition:DWCart-prism-description} for a more general statement). Some important examples of prismatic crystals for this paper are the following:

\begin{example}[The prismatic cohomology sheaf]
\label{example:prismatic-cohomology-sheaf}
Let $X$ be a smooth and proper $\Z_p$-scheme. Then the prismatic crystal $(A,I) \mapsto \RGamma_{\Prism}( X_{\overline{A}} / A )$
determines a perfect complex on $\WCart$, which we will denote by $\mathscr{H}_{\Prism}(X)$ and refer to as the
{\it prismatic cohomology sheaf} of $X$.

More generally, suppose that $\mathfrak{X}$ is a $p$-adic formal scheme
which is quasi-compact, quasi-separated, and that the structure sheaf $\calO_{\mathfrak{X}}$ has bounded
$p$-power torsion (if the last condition is satisfied, we will say that $\mathfrak{X}$ is {\it bounded}; see Definition \ref{definition:bounded-formal-scheme}).
Applying a variant of the preceding construction, we obtain a quasi-coherent complex $\mathscr{H}_{\Prism}( \mathfrak{X} )$ on $\WCart$ (which is generally not a perfect complex); see  Construction~\ref{construction:prismatic-cohomology-sheaves} for a more general construction.
\end{example}

\begin{example}[The Breuil-Kisin twist]
\label{ex:BKTwistIntro}
For every bounded prism $(A,I)$, we let $A\{-1\}$ denote the prismatic cohomology group $\mathrm{H}^{2}_{ \Prism}( \mathbf{P}^{1}_{ \overline{A} } / A )$.
This is an invertible $A$-module, whose inverse we denote by $A\{1\}$ and refer to as the {\it Breuil-Kisin twist of $A$}. In 
\S \ref{section:twist-and-log}, we give a purely algebraic construction of the $A$-module $A\{1\}$ (which does not depend on the theory
of prismatic cohomology); see Definition \ref{definition:twist-general}. The construction $(A,I) \mapsto A\{1\}$ is an invertible
prismatic crystal, which we can identify with a line bundle $\calO_{ \WCart}\{1\}$ on the Cartier-Witt stack. For any quasi-coherent complex $\mathscr{F}$ on $\WCart$, we let $\mathscr{F}\{n\}$ denote the tensor product of $\mathscr{F}$ with the $n$th power of $\calO_{\WCart}\{1\}$.
\end{example}

Let $\mathfrak{X}$ be a bounded $p$-adic formal scheme which is quasi-compact and quasi-separated. For every integer $n$, we let $\RGamma_{\Prism}( \mathfrak{X} )\{n\}$ denote the complex of (derived) global sections $\RGamma( \WCart, \mathscr{H}_{\Prism}( \mathfrak{X})\{n\})$. In the special case $n = 0$, we denote $\RGamma_{\Prism}( \mathfrak{X} )\{n\}$
by $\RGamma_{\Prism}( \mathfrak{X} )$ and refer to it as the {\it absolute prismatic complex}
of $\mathfrak{X}$. We denote the cohomology of this complex by $\mathrm{H}^{\ast}_{\Prism}( \mathfrak{X} )$ and refer to it as
the {\it absolute prismatic cohomology} of $\mathfrak{X}$. In \S \ref{section:absolute-prismatic-cohomology}, we will see that this construction has the following features:

\begin{itemize}
\item If $\mathfrak{X}$ is a quasisyntomic $p$-adic formal scheme, the absolute prismatic cohomology $\mathrm{H}^{\ast}_{\Prism}( \mathfrak{X} )$ admits
a site-theoretic definition: it is the cohomology of the absolute prismatic site of $\mathfrak{X}$ (Theorem \ref{theorem:absolute-compute-with-stack}).

\item If $\mathfrak{X}$ is defined over $A/I$ for a bounded prism $(A,I)$, then there is a comparison map
$$ \mathrm{H}^{\ast}_{\Prism}( \mathfrak{X} ) \rightarrow \mathrm{H}^{\ast}_{\Prism}( \mathfrak{X} / A )$$
from absolute to relative prismatic cohomology. If $(A,I)$ is a perfect prism, then this map is an isomorphism (Proposition \ref{proposition:absolute-vs-relative}).
In particular, absolute prismatic cohomology is a generalization of the $A_{\mathrm{inf}}$-cohomology introduced in \cite{BMS1}.

\item If $X$ is an $\F_p$-scheme, then there is a comparison map
$$ \mathrm{H}^{\ast}_{\Prism}(X) \rightarrow \mathrm{H}^{\ast}_{\crys}( X / \Z_p )$$
from the absolute prismatic cohomology of $X$ to the crystalline cohomology of $X$. If $X$ is quasisyntomic, then
this map is an isomorphism (Theorem \ref{theorem:crystalline-comparison}).
\end{itemize}

By virtue of the first point, it is possible to define the absolute prismatic complex $\RGamma_{\Prism}( \mathfrak{X} )$ (at least for quasisyntomic formal schemes)
without mentioning the Cartier-Witt stack. However, if one wishes to {\em compute} with absolute prismatic cohomology in situations of arithmetic
interest, then it is useful to understand the geometry of $\WCart$.

\begin{example}[$\WCart$ via simple quotient stacks]
The $q$-de Rham prism $( \Z_p[[q-1]], ( 1 + q + \cdots + q^{p-1} ) )$ determines a map $\rho_{\qdR}: \Spf( \Z_p[[q-1]]) \rightarrow \WCart$.
This map is equivariant with respect to the action of the profinite group $\Z_p^{\times}$ (which acts on the $q$-de Rham prism via the
construction $q \mapsto q^{\alpha}$), which is trivial modulo $(q-1)$. We therefore obtain a commutative diagram of stacks
$$ \xymatrix@R=50pt@C=50pt{ [\Spf( \Z_p ) / \Z_p^{\times} ] \ar[r] \ar[d] & [ \Spf( \Z_p[[q-1]] ) / \Z_p^{\times} ] \ar[d]^{ \rho_{\qdR} } \\
\Spf(\Z_p) \ar[r]^-{ \rho_{\dR} } & \WCart. }$$
In \S \ref{subsection:computing-with-WCart}, we show that, for $p > 2$, this diagram behaves like a pushout square
for the purpose of cohomological calculations (see Theorem \ref{theorem:qdR-WCart}).  In particular, if $X$ is a proper smooth $\Z_p$-scheme,
we obtain a homotopy pullback diagram of complexes
$$ \xymatrix@R=50pt@C=50pt{ \RGamma_{\Prism}(X) \ar[r] \ar[d] & \RGamma_{\qdR}( X)^{h \Z_p^{\times} } \ar[d] \\
\RGamma_{\dR}( X / \Z_p) \ar[r] & \RGamma_{\dR}( X / \Z_p)^{h \Z_p^{\times} }, }$$
which allows us to compute the absolute prismatic cohomology of $X$ in terms of the de Rham cohomology of $X$
together with its $q$-deformation (see Remark \ref{remark:concrete-qdR}).
\end{example}

\begin{remark}[Prismatic $F$-crystals]
The structure sheaf $\mathcal{O}_\Prism$ on the prismatic site $\mathrm{Spf}(\mathbf{Z}_p)_\Prism$ comes equipped with an endomorphism that lifts the Frobenius on $\mathcal{O}_\Prism/p$; this structure is reflected in a map $F:\WCart \to \WCart$ which is a lift of the Frobenius on $\WCart \times \mathrm{Spec}(\mathbf{F}_p)$ (see Construction~\ref{construction:Frobenius-on-stack}). Any perfect prismatic crystal $\mathcal{E}$ coming from geometry (such as Example~\ref{example:prismatic-cohomology-sheaf} or $\mathcal{O}_{\WCart}\{-1\}$ from Example~\ref{ex:BKTwistIntro}) is naturally (effective) prismatic $F$-crystal, i.e., it comes equipped with a map $F^* \mathcal{E} \to \mathcal{E}$ satisfying suitable conditions. This structure shall play a crucial role in some of the structures studied in this paper (such as the Nygaard filtration and syntomic complexes). There is also a close relationship between prismatic $F$-crystals and the classical theory of crystalline Galois representations; we refer the interested reader to \cite{bhatt2021prismatic} for more on this connection.
\end{remark}

\begin{remark}[Drinfeld's stacks]
The Cartier-Witt stack $\WCart$ was independently discovered by Drinfeld, albeit with somewhat different motivations; see \cite[\S 4]{drinfeld-prismatic}, where Drinfeld writes $\Sigma$ for what we call $\WCart$. In fact, \cite{drinfeld-prismatic} went further to construct enlargements $\Sigma'$ and $\Sigma''$ of $\Sigma$ which were expected to encode prismatic cohomology with its Nygaard filtration in stack-theoretic terms. (This expectation will be verified in a sequel to the present paper.)
\end{remark}

\begin{remark}[$\WCart_{\mathfrak{X}}$ for a $p$-adic formal scheme $\mathfrak{X}$]
The stack $\WCart$ geometrizes the absolute prismatic site of $\mathrm{Spf}(\mathbf{Z}_p)_\Prism$: perfect complexes on $\WCart$ identify with perfect prismatic crystals on $\mathrm{Spf}(\mathbf{Z}_p)$ via the equivalence \eqref{eq:GeometrizePrismCrys}.  Similarly, for any bounded $p$-adic formal scheme $\mathfrak{X}$ satisfying mild conditions, one can construct a stack $\WCart_{\mathfrak{X}}$ with the property that perfect complexes on $\WCart_{\mathfrak{X}}$ are naturally identified with perfect prismatic crystals on $\mathfrak{X}$.  The definition of $\WCart_{\mathfrak{X}}$ relies on formal derived algebraic geometry (and, in fact, the natural domain of definition of the construction is a derived $p$-adic formal scheme $\mathfrak{X}$). To avoid introducing additional technical complications in an already long paper, we defer the discussion of this generalization to the sequel \cite{BhattLurieAPCsequel}. We remark that these generalizations were also suggested by Drinfeld in \cite[\S 1.1]{drinfeld-prismatic}.
\end{remark}

\subsection{Diffracted Hodge cohomology and the Hodge-Tate divisor}
\label{ss:DiffHodHT}
The global structure of the Cartier-Witt stack is somewhat complicated. However, there is a closed substack $\WCart^{\mathrm{HT}} \subset \WCart$,
which we refer to as the {\it Hodge-Tate divisor}, which admits a relatively simple description: it is the formally completed classifying stack of
a commutative affine group scheme $\mathbf{G}_{m}^{\sharp}$, given by the divided power envelope of $\mathbf{G}_{m}$ along its identity section.
Quasi-coherent complexes on $\WCart^{\mathrm{HT}}$ are easy to analyze: they can be identified with pairs $(M, \Theta)$, where $M$ is a $p$-complete
complex of abelian groups and $\Theta$ is an endomorphism of $M$ satisfying a certain integrality condition (see Theorem \ref{theorem:compute-with-HT}).
We will refer to $\Theta$ as the {\it Sen operator}; it is closely related to the work of Sen on semilinear Galois representations, which we review in
\S \ref{subsection:Sen-theory}. 

\begin{remark}[Cohomological dimension of $\mathfrak{X}_\Prism$]
One consequence of our analysis of the Hodge-Tate divisor is that the Cartier-Witt stack $\WCart$ has cohomological dimension $1$.
For example, if $\mathfrak{X}$ is an affine formal scheme which is smooth of relative dimension $n$ over $\Spf(\Z_p)$, then
the absolute prismatic cohomology groups $\mathrm{H}^{\ast}_{\Prism}( \mathfrak{X} )$ vanish for $\ast > n+1$.
This is consistent with the heuristic that the prismatic cohomology of $\mathfrak{X}$ should be viewed as its de Rham cohomology
relative to ``the field with one element'', over which it has relative dimension $n+1$.
\end{remark}

By applying the preceding analysis to prismatic crystals of algebro-geometric origin, we can obtain a new invariant.
Let $\mathfrak{X}$ be a bounded $p$-adic formal scheme. It follows from the discussion above that we can identify the
restriction $\mathscr{H}_{\Prism}( \mathfrak{X} )|_{ \WCart^{\mathrm{HT}}}$ with a pair $(M, \Theta)$, for some $p$-complete cochain complex $M$.
The complex $M$ can be realized as the (derived) global sections $\RGamma( \mathfrak{X}, \Omega^{\DHod}_{\mathfrak{X}} )$, where
$\Omega^{\DHod}_{\mathfrak{X}} \in \calD( \mathfrak{X} )$ is a certain (quasi-coherent) complex on $\mathfrak{X}$ which we will refer
to as the {\it diffracted Hodge complex} of $\mathfrak{X}$ (Notation \ref{notation:diffracted-Hodge-formal-scheme}). This construction has
the following features:

\begin{itemize}
\item The diffracted Hodge complex $\Omega^{\DHod}_{\mathfrak{X}}$ is equipped with an exhausive filtration
$$ \Fil_{0}^{\conj} \Omega^{\DHod}_{\mathfrak{X}} \rightarrow \Fil_{1}^{\conj} \Omega^{\DHod}_{\mathfrak{X}} \rightarrow
\Fil_{2}^{\conj} \Omega^{\DHod}_{\mathfrak{X}} \rightarrow \cdots,$$
which we will refer to as the {\it conjugate filtration}. Moreover, if $\mathfrak{X}$ is smooth over $\Spf(W(k))$ for some perfect field $k$, then there is a
canonical isomorphism of the associated graded complex $\gr_{\ast}^{\conj}  \Omega^{\DHod}_{\mathfrak{X}}$
with the Hodge complex $\bigoplus_{n \geq 0} \widehat{\Omega}^{n}_{ \mathfrak{X} / W(k) }[-n]$ (Example \ref{example:diffracted-Hodge-of-smooth}).

\item  The Sen operator on $\RGamma( \mathfrak{X}, \Omega^{\DHod}_{\mathfrak{X}} )$ is induced by an
endomorphism of the diffracted Hodge complex $\Omega^{\DHod}_{\mathfrak{X}}$ itself, which we will also denote by
$\Theta$ and refer to as the {\it Sen operator}. Moreover, this endomorphism is compatible with the conjugate filtration
and acts diagonally at the associated graded level: that is, the induced endomorphism of $\gr_{n}^{\conj} \Omega^{\DHod}_{\mathfrak{X} }$ 
is given by multiplication by $-n$ (Notation \ref{notation:Sen-operator-on-diffracted}). In particular, the conjugate filtration on $\Omega^{\DHod}_{\mathfrak{X}}$
splits rationally (into eigenspaces of $\Theta$).

\item Suppose that $\mathfrak{X}$ is smooth over $\Spf(W(k))$ for some perfect field $k$, and let 
$$X_0 = \Spec(k) \times_{\Spf(W(k)) } \mathfrak{X}$$ denote its special fiber. Let $i: X_0 \hookrightarrow \mathfrak{X}$ denote the inclusion map. Then the (derived) pullback $i^{\ast} \Omega^{\DHod}_{\mathfrak{X}}$
can be identified with the de Rham complex of $X_0$ (or, more precisely, with its pushforward along the absolute Frobenius map $\varphi: X_0 \rightarrow X_0$).
This identification carries the conjugate filtration on $\Omega^{\DHod}_{\mathfrak{X}}$ to the usual conjugate filtration on the de Rham complex.
\end{itemize}

\begin{remark}[The Deligne-Illusie theorem via the Sen operator]
\label{rmk:DIintro}
Let $k$ be a perfect field of characteristic $p$ and let $X_0$ be a smooth $k$-scheme. Then the conjugate filtration of the de Rham complex 
$\Omega^{\ast}_{X/k}$ determines a spectral sequence
\begin{equation}\label{conjugate-spectral-sequence} \mathrm{H}^{s}( X_0, \Omega^{t}_{X_0/k} ) \Rightarrow \varphi_{\ast} \mathrm{H}^{s+t}_{\dR}( X_0 / k). \end{equation}
In \cite{MR894379}, Deligne and Illusie showed that if the dimension $\dim(X_0)$ is smaller than $p$ and
$X_0$ can be lifted to a smooth $W(k) / (p^2)$-scheme, then the spectral sequence (\ref{conjugate-spectral-sequence})
degenerates (if $X_0$ is also proper over $k$, it then follows by a dimension-counting argument that the Hodge-to-de-Rham spectral sequence
$\mathrm{H}^{s}( X_0, \Omega^{t}_{X_0/k}) \Rightarrow \mathrm{H}^{s+t}_{\dR}(X_0/k)$ also degenerates). 

Using formal properties of the diffracted Hodge complex, we can immediately deduce a weaker form of this result.
Suppose that $X_0$ arises as the special fiber of a formal scheme $\mathfrak{X}$ which is smooth over $\Spf( W(k) )$. Then the Frobenius pushforward $\varphi_{\ast} \Omega^{\ast}_{X_0/k}$ can be identified with the restriction of the diffracted Hodge complex 
$\Omega^{\DHod}_{\mathfrak{X}}|_{ X_0 }$. If the dimension of $X_0$ is smaller than $p$, then the eigenvalues of the
Sen operator on the associated graded complex $\gr^{\conj}_{\ast} \Omega^{\DHod}_{\mathfrak{X}} \simeq \Omega^{\ast}_{\mathfrak{X} / W(k)}$
are distinct modulo $p$. It follows that the conjugate filtration of $\Omega^{\DHod}_{\mathfrak{X}}$ is canonically split
(into generalized eigenspaces of $\Theta$), which immediately implies the degeneration of the spectral sequence (\ref{conjugate-spectral-sequence}) (see Remark \ref{remark:Deligne-Illusie}).
\end{remark}

\begin{remark}[Extension to other DVRs]
The results and constructions discussed in \S \ref{ss:APCCW} and \S \ref{ss:DiffHodHT} concern the structure of prismatic crystals on $\mathrm{Spf}(\mathbf{Z}_p)$ arising from bounded $p$-adic formal schemes $\mathfrak{X}$ over $\mathbf{Z}_p$. These extend essentially without change if $\mathbf{Z}_p$ is replaced by $W(k)$ for any perfect field $k$ of characteristic $p$. In fact, most of our constructions admit analogues when $\mathbf{Z}_p$ is replaced by $\mathcal{O}_K$ for any complete  discretely valued extension $K/\mathbf{Q}_p$ with perfect residue field.  However, the resulting statements depend on the absolute ramification index of $K$, and the results are strongest when $K$ is unramified. For instance, the condition on the eigenvalues of the (analog of the) Sen operator is weaker in the ramified case, and there is consequently no analog of Remark~\ref{rmk:DIintro} for proper smooth schemes $\mathfrak{X}/\mathcal{O}_K$ for ramified $K$. In the interest of simplicity, we have elected to restriction attention to the unramified case in this paper; the sequel \cite{BhattLurieAPCsequel} contains some results for more general $K$. 
\end{remark}

\subsection{The Nygaard filtration}
For any bounded $p$-adic formal scheme $\mathfrak{X}$, the absolute prismatic complex $\RGamma_{\Prism}( \mathfrak{X} )$ is equipped with
an endomorphism $\varphi: \RGamma_{\Prism}( \mathfrak{X} ) \rightarrow \RGamma_{\Prism}( \mathfrak{X} )$, which we will refer to as
the {\it Frobenius morphism} (Notation \ref{notation:Frobenius-on-absolute}). For the twisted prismatic complexes $\RGamma_{\Prism}(\mathfrak{X})\{n\}$, the situation is more subtle: the Frobenius morphism is only partially defined. More precisely, the complex $\RGamma_{\Prism}(\mathfrak{X})\{n\}$ admits a decreasing filtration
$$ \cdots \rightarrow \Fil^{2}_{\Nyg} \RGamma_{\Prism}(\mathfrak{X})\{n\} \rightarrow \Fil^{1}_{\Nyg} \RGamma_{\Prism}(\mathfrak{X})\{n\}
\rightarrow \Fil^{0}_{\Nyg}  \RGamma_{\Prism}(\mathfrak{X})\{n\} = \RGamma_{\Prism}( \mathfrak{X})\{n\},$$
which we refer to as the {\it absolute Nygaard filtration}, and a naturally defined Frobenius operator $\varphi\{n\}:  \Fil^{n}_{\Nyg} \RGamma_{\Prism}(\mathfrak{X})\{n\} \rightarrow \RGamma_{\Prism}(\mathfrak{X})\{n\}$. In \S \ref{section:Nygaard-filtration}, we introduce the Nygaard filtration and show that it has the following features:
\begin{itemize}
\item Let $(A,I)$ be a perfect prism and let $\mathfrak{X}$ be a $p$-adic formal scheme which is smooth over $\Spf(A/I)$. Then the quasi-isomorphism
$$ \RGamma_{\Prism}( \mathfrak{X} )\{n\} \simeq \RGamma_{\Prism}( \mathfrak{X}/ A )\{n\} \simeq \varphi_{A}^{\ast} \RGamma_{\Prism}( \mathfrak{X}/ A )\{n\}$$
can be promoted to a filtered quasi-isomorphism, where the left side is equipped with the absolute Nygaard filtration and the right
side is equipped with the relative Nygaard filtration introduced in \cite{prisms} (Theorem \ref{theorem:compare-Nygaard-filtrations}).

\item Let $\mathfrak{X} = \Spf(R)$ be an affine formal scheme whose coordinate ring $R$ is quasiregular semiperfectoid. Then
the absolute prismatic complex $\RGamma_{\Prism}( \mathfrak{X} )$ can be identified with the underlying commutative ring of a prism $(A,I)$.
Under this identification, each $\Fil^{n}_{\Nyg} \RGamma_{\Prism}( \mathfrak{X} )$ corresponds to the ideal $\{ x \in A: \varphi_{A}(x) \in I^{n} \}$,
regarded as a chain complex concentrated in degree zero (Corollary \ref{corollary:absolute-Nygaard-qrsp}). Moreover, if $R$ is an $\F_p$-algebra,
we give an explicit set of generators for this ideal (Proposition \ref{proposition:concrete-Nygaard-qrsp}).

\item Let $X$ be an $\F_p$-scheme. If $X$ is regular, then the identification of the absolute prismatic complex $\RGamma_{\Prism}(X)$
with the crystalline cochain complex $\RGamma_{\crys}(X)$ carries the Nygaard filtration on absolute prismatic cohomology to the classical
Nygaard filtration on crystalline cohomology (Proposition \ref{proposition:concrete-Nygaard-regular}).

\item Let $\mathfrak{X}$ be a bounded $p$-adic formal scheme. For every pair of integers $m$ and $n$, we have a canonical fiber sequence
$$ \gr^{m}_{\Nyg} \RGamma_{\Prism}( \mathfrak{X} )\{n\} \rightarrow 
\RGamma( \mathfrak{X}, \Fil_{m}^{\conj} \Omega^{\DHod}_{\mathfrak{X} } ) \xrightarrow{ \Theta+m} 
\RGamma( \mathfrak{X}, \Fil_{m-1}^{\conj} \Omega^{\DHod}_{\mathfrak{X} } )$$
(Remark \ref{remark:Nygaard-associated-graded}). In particular, the complex $ \gr^{m}_{\Nyg} \RGamma_{\Prism}( \mathfrak{X} )\{n\} $
is canonically independent of $n$ (Remark \ref{remark:laurent-in-e}). 
\end{itemize}

For every bounded $p$-adic formal scheme $\mathfrak{X}$, let us write $\RGamma_{ \widehat{\Prism} }( \mathfrak{X} )$ for the homotopy limit of the tower
$$ \cdots \rightarrow \RGamma_{ \Prism}( \mathfrak{X} ) / \Fil^{2}_{\Nyg} \RGamma_{ \Prism}( \mathfrak{X} )
\rightarrow \RGamma_{ \Prism}( \mathfrak{X} ) / \Fil^{1}_{\Nyg} \RGamma_{ \Prism}( \mathfrak{X} )
\rightarrow \RGamma_{ \Prism}( \mathfrak{X} ) / \Fil^{0}_{\Nyg} \RGamma_{ \Prism}( \mathfrak{X} ) \simeq 0.$$
We refer to $\RGamma_{ \widehat{\Prism} }( \mathfrak{X} )$ as the {\it Nygaard-complete} prismatic complex of $\mathfrak{X}$.
Beware that the tautological map $\RGamma_{ \Prism }( \mathfrak{X} ) \rightarrow \RGamma_{ \widehat{\Prism} }( \mathfrak{X} )$
is generally not a quasi-isomorphism. However, we will show that it is a quasi-isomorphism whenever $\mathfrak{X}$
is Noetherian and the special fiber of $\mathfrak{X}$ is regular (Proposition \ref{proposition:smooth-absolute-Nygaard-complete}); these
conditions are satisfied, for example, if $\mathfrak{X}$ is smooth over $\Spf(\Z_p)$.

\begin{remark}[Connections to $K$-theoretic invariants]
The Nygaard-completed prismatic complexes $\RGamma_{ \widehat{\Prism} }( \mathfrak{X} )$ (along with their twisted counterparts $\RGamma_{ \widehat{\Prism} }( \mathfrak{X} )\{n\}$) have a close relationship with $K$-theoretic invariants. Assume for simplicity that $\mathfrak{X} = \Spf(R)$ where
$R$ is a quasisyntomic ring. Let $\TP(R)$ denote the topological periodic cyclic homology spectrum of $R$, and let
$\TP(R)^{\wedge}_{p}$ denote its $p$-completion. In \cite{BMS2}, the first author,
Morrow, and Scholze construct a complete and exhaustive {\em motivic} filtration $\Fil^{\bullet}_{\mot} \TP(R)^{\wedge}_{p}$
on the spectrum $\TP(R)^{\wedge}_{p}$, together with canonical identifications
$$ \gr^{n}_{\mot} \TP(R)^{\wedge}_{p} \simeq \RGamma_{\widehat{\Prism}}( \Spf(R) )\{n\}[2n].$$
In \S \ref{section:periodic-cyclic-homology}, we review the definition of the motivic filtration and
extend the construction to arbitrary commutative rings (Theorem \ref{theorem:BMS2-main}). We also
introduce an integral counterpart of the motivic filtration (Construction \ref{construction:integral-motivic-filtration}),
which is defined on spectrum $\TP(R)$ before completing at the prime $p$ (and with no assumptions on $R$). 
Beware that, in general, the filtrations introduced here are not exhaustive (though they are always complete:
see Corollary \ref{corollary:TP-filtration-complete}).
\end{remark}

\subsection{Syntomic cohomology}
Let $\mathfrak{X}$ be a bounded $p$-adic formal scheme. For each $n \geq 0$, we have a tautological map
$\iota: \Fil^{n}_{\Nyg} \RGamma_{\Prism}( \mathfrak{X} )\{n\} \rightarrow \RGamma_{\Prism}( \mathfrak{X} )\{n\}$. We
let $\RGamma_{\Syn}( \mathfrak{X}, \Z_p(n) )$ denote the fiber of the map
$$ ( \varphi\{n\} - \iota ): \Fil^{n}_{\Nyg} \RGamma_{\Prism}( \mathfrak{X} )\{n\} \rightarrow \RGamma_{\Prism}( \mathfrak{X} )\{n\}.$$
Following \cite{BMS2}, we will refer to $\RGamma_{\Syn}( \mathfrak{X}, \Z_p(n) )$ as the {\it $n$th syntomic complex} of the formal scheme $\mathfrak{X}$.

In the case $n=1$, there is a close connection of the complex $\RGamma_{\Syn}( \mathfrak{X}, \Z_p(n) )$ with the cohomology of the multiplicative group $\mathbf{G}_m$. Let $(A,I)$ be a prism, and let $T_p( \overline{A}^{\times} )$ denote the $p$-adic Tate module of the multiplicative group of the quotient ring $\overline{A} = A/I$.
In \S \ref{section:twist-and-log}, we construct a canonical homomorphism $\log_{\Prism}: T_p( \overline{A}^{\times} ) \rightarrow A\{1\}$, which we refer to as the {\it prismatic logarithm} (Construction \ref{construction:Tate-logarithm}). In \S \ref{section:first-chern}, we exploit the prismatic logarithm to construct a comparison map
$$ c_{1}^{\Syn}: \RGamma_{\mathet}( \Spec(R), \mathbf{G}_m)[-1] \rightarrow \RGamma_{\Syn}( \Spf(R), \Z_p(1) )$$
for every $p$-complete commutative ring $R$, which we will refer to as the {\it syntomic first Chern class} (Proposition \ref{proposition:prismatic-chern-class-construction}). Our main result is that $c_{1}^{\Syn}$ becomes a quasi-isomorphism after (derived) $p$-completion (Theorem \ref{theorem:syntomic-chern-class-isomorphism}).
In the case where $R$ is quasisyntomic, this is essentially Proposition~7.17 of \cite{BMS2}. We give a different proof here, which avoids the machinery of
algebraic $K$-theory.

In \S \ref{section:etale-comparison}, we study the relationship between syntomic cohomology of a formal scheme $\Spf(R)$ with the \'{e}tale cohomology of the affine scheme $\Spec(R[1/p])$. Using the theory of arc descent, we show that there are essentially unique comparison maps
$$  \gamma_{\Syn}^{\mathet}\{n\}: \RGamma_{\Syn}( \Spf(R), \Z_p(n) ) \rightarrow \RGamma_{\mathet}( \Spec(R[1/p]), \Z_p(n) )$$
which are multiplicative, functorial, and compatible with Chern classes (Theorem \ref{theorem:etale-comparison}). In \S \ref{subsection:syntomic-complexes-general},
we exploit this comparison to ``decomplete'' the theory of syntomic complexes. Suppose that $S$ is a commutative ring and assume for simplicity that the $p$-power torsion in $S$ is bounded. We let $\RGamma_{\Syn}( \Spec(S), \Z_p(n) )$ denote the homotopy fiber product of the diagram
$$ \RGamma_{\Syn}( \Spf( \widehat{S} ), \Z_p(n) ) \xrightarrow{
\gamma_{\Syn}^{\mathet}\{n\} } \RGamma_{\mathet}( \Spec( \widehat{S}[1/p] ), \Z_p(n) )
\leftarrow \RGamma_{\mathet}( \Spec( S[1/p] ), \Z_p(n) ),$$
where $\widehat{S}$ denotes the $p$-adic completion of $S$ (Construction \ref{construction:syntomic-complexes-general}). The construction $\Spec(S) \mapsto \RGamma_{\Syn}( \Spec(S), \Z_p(n) )$ satisfies flat descent (Proposition \ref{rhoax}), and therefore admits a natural extension 
$X \mapsto \RGamma_{\Syn}(X, \Z_p(n) )$ to the category of all schemes (Variant \ref{variant:syntomic-sheaf-integrally-globalized}). We denote the cohomology groups of $\RGamma_{\Syn}( X, \Z_p(n))$ by $\mathrm{H}^{\ast}_{\Syn}(X, \Z_p(n) )$, which we refer to as
the {\it syntomic cohomology groups of $X$}. Decompleting the constructions \S \ref{section:first-chern}, we obtain a natural map of complexes
$$ \RGamma_{\mathet}( X, \mathbf{G}_m)[-1] \rightarrow \RGamma_{\Syn}( X, \Z_p(1) )$$
which yields, after passing to cohomology in degree $2$, a homomorphism of abelian groups
$c_{1}^{\Syn}: \Pic(X) \rightarrow \mathrm{H}^{2}_{\Syn}(X, \Z_p(1) )$ which we also refer to as the {\it syntomic first Chern class}.

In \S \ref{section:calculate-syntomic}, we use the existence of Chern classes to reproduce several standard calculations in the setting of syntomic cohomology:
\begin{itemize}
\item Let $X$ be a scheme, let $\mathscr{E}$ be a vector bundle of rank $m$ on $X$, let $\mathbf{P}( \mathscr{E} )$ denote the projectivization of
$\mathscr{E}$, and let $\calO(1)$ denote the tautological line bundle on $\mathbf{P}( \mathscr{E} )$. In \S \ref{subsection:projective-bundle-scheme},
we show that the syntomic cohomology of $\mathbf{P}( \mathscr{E} )$ is a free module over the syntomic cohomology of $X$,
generated by the syntomic cohomology classes $\{ c_1^{\Syn}( \calO(1) )^{i} \}_{0 \leq i < m}$ (Theorem \ref{theorem:projective-bundle-formula}).

\item It follows from Theorem \ref{theorem:projective-bundle-formula} that if $\mathscr{E}$ is a vector bundle of rank $n$ on a scheme $X$,
then there are unique elements $\{ c_i^{\Syn}( \mathscr{E} ) \in \mathrm{H}^{2i}_{\Syn}( X, \Z_p(i) ) \}_{1 \leq i \leq m}$ satisfying the identity
$$ c_{1}^{\Syn}( \calO(1) )^{n} + \sum_{i=1}^{m} c_{i}^{\Syn}( \mathscr{E} ) \cdot c_{1}^{\Syn}( \calO(1) )^{m-i} = 0$$
in $\mathrm{H}^{2m}_{\Syn}( \mathbf{P}(\mathscr{E}), \Z_p(m) )$. We will refer to $c_i^{\Syn}( \mathscr{E} )$ as the
{\it $i$th syntomic Chern class of $\mathscr{E}$} (Construction \ref{construction:higher-chern-classes}). The Chern classes given by
this construction can be regarded as simultaneous refinements of their counterparts in crystalline, de Rham, and \'{e}tale cohomology. 
Moreover, they enjoy all of the expected features: for example, we show that the total Chern class
$$ 1 + c_1^{\Syn}( \mathscr{E} ) + c_{2}^{\Syn}( \mathscr{E} ) + \cdots$$ 
is multiplicative in short exact sequences (Theorem \ref{theorem:additivity-chern}). 

\item Let $\BGL_m$ denote the classifying stack of the affine group scheme $\GL_m$, and let
$\mathscr{E}$ denote the universal vector bundle of rank $m$ over $\BGL_m$. For any scheme $X$ satisfying some mild hypotheses,
the syntomic cohomology ring of the product $\BGL_m \times X$ is a polynomial algebra over the syntomic cohomology ring of $X$,
with generators given by the Chern classes $\{ c_{i}^{\Syn}( \mathscr{E} ) \}_{1 \leq i \leq m}$ (Theorem \ref{theorem:syntomic-cohomology-of-BGL}).
Stated more informally, any cohomological invariant of a vector bundle $\mathscr{E}$ can be expressed in terms of its Chern classes.
\end{itemize}

Moreover, we establish analogous results for other cohomological invariants studied in this paper, such as prismatic cohomology and diffracted Hodge cohomology.

\begin{remark}
One upshot of our construction of syntomic Chern classes is that we also obtain Chern classes for vector bundles in prismatic cohomology, and thus in the $A_{\inf}$-cohomology of \cite{BMS1} for smooth $p$-adic formal schemes over $\mathcal{O}_C$ with $C/\mathbf{Q}_p$ complete and algebraically closed. Such a theory was constructed independently by \cite[\S 6.2]{KubrakPrikhodko} using the interplay between the theory of crystalline Galois representations and Breuil-Kisin modules (see \cite[\S 5.2]{KubrakPrikhodko}).
\end{remark}

\subsection{Cochains and Animation}\label{subsection:cochain-animation}

Throughout most of this paper, we will freely use the language of $\infty$-categories (for various accounts, see \cite{kerodon}, \cite{HTT}, \cite{groth2015short},
and \cite{hinich2018lectures}). For the reader's convenience, we briefly describe two examples of $\infty$-categories and the role they will play in our story.

\begin{example}[Derived $\infty$-Categories]\label{example:derived-infty-category}
Let $A$ be a commutative ring, and let $\mathrm{Ch}(A)$ denote the category of cochain complexes of $A$-modules. 
A morphism $f: M^{\bullet} \rightarrow N^{\bullet}$ in $\mathrm{Ch}(A)$ is said to be a {\it quasi-isomorphism} if it induces
an isomorphism of cohomology groups $\mathrm{H}^{\ast}(M) \rightarrow \mathrm{H}^{\ast}(N)$. Recall that the {\it derived category}
$D(A)$ is obtained from $\mathrm{Ch}(A)$ by formally adjoining inverses to quasi-isomorphisms. More precisely,
there is a functor $T: \mathrm{Ch}(A) \rightarrow D(A)$ with the following universal property: for every category $\calC$,
precomposition with $T$ induces a fully faithful functor
$$ \{ \text{Functors $D(A) \rightarrow \calC$} \} \rightarrow \{ \text{Functors $\mathrm{Ch}(A) \rightarrow \calC$} \},$$
whose essential image consists of functors $F: \mathrm{Ch}(A) \rightarrow \calC$ which carry quasi-isomorphisms of cochain complexes
to isomorphisms in $\calC$.

The derived category $D(A)$ can be realized as the homotopy category of an $\infty$-category $\calD(A)$, which we will refer
to as the {\it derived $\infty$-category} of $A$. This $\infty$-category can be characterized by a similar universal property:
there is a functor $\widetilde{T}: \mathrm{Ch}(A) \rightarrow \calD(A)$ such that, for every $\infty$-category $\calC$, precomposition with
$\widetilde{T}$ induces a fully faithful functor
$$ \{ \text{Functors $\calD(A) \rightarrow \calC$} \} \rightarrow \{ \text{Functors $\mathrm{Ch}(A) \rightarrow \calC$} \},$$
whose essential image consists of those functors $F: \mathrm{Ch}(A) \rightarrow \calC$ which carry quasi-isomorphisms of cochain complexes
to isomorphisms in $\calC$.
\end{example}

\begin{remark}
Let $A$ be a commutative ring. There is a fully faithful functor from the abelian category of $A$-modules to the derived
$\infty$-category $\calD(A)$, whose essential image consists of those complexes $M \in \calD(A)$
for which the cohomology groups $\mathrm{H}^{n}(M)$ vanish for $n \neq 0$. Throughout this paper, we will abuse terminology by identifying the ordinary
category of $A$-modules with its image under functor. In particular, if $M$ is an object of the $\infty$-category $\calD(A)$ for which the 
the cohomology groups $\mathrm{H}^{n}(M)$ vanish for $n \neq 0$, then we will implicitly identify $M$ with the $A$-module $\mathrm{H}^{0}(M)$.
\end{remark}

In this paper, we will study several cohomological invariants of (bounded) $p$-adic formal schemes $\mathfrak{X}$, such as the absolute prismatic
cohomology $\mathrm{H}^{\ast}_{\Prism}( \mathfrak{X} )$. For many purposes, it is essential to work with cochain-level
incarnation of these invariants: that is, to emphasize cochain complexes $\RGamma_{\Prism}( \mathfrak{X} )$
which compute absolute prismatic cohomology, rather than the cohomology groups themselves. These cochain complexes are well-defined
and functorial up to quasi-isomorphism. More precisely, the construction $\mathfrak{X} \mapsto \RGamma_{ \Prism}( \mathfrak{X} )$ determines a contravariant functor
\begin{equation}\label{equation:functor-up-to-homotopy} \{ \text{Bounded $p$-adic formal schemes} \}  \rightarrow  \calD(\Z_p),
\end{equation}
where $\calD(\Z_p)$ is the derived $\infty$-category of Example \ref{example:derived-infty-category}.

\begin{remark}\label{remark:suffra}
By passing to homotopy categories, (\ref{equation:functor-up-to-homotopy}) determines a contravariant functor of ordinary categories
\begin{equation}\label{equation:functor-up-to-homotopy2} \{ \text{Bounded $p$-adic formal schemes} \}  \rightarrow  D(\Z_p).
\end{equation}
However, this passage loses a great deal of essential information. For example, an essential feature of the prismatic complexes
$\RGamma_{\Prism}( \mathfrak{X} )$ (and other invariants we study) is that they are {\em local} in nature. For example,
if $\mathfrak{X}$ is a bounded $p$-adic formal scheme, then the absolute prismatic complex $\RGamma_{\Prism}( \mathfrak{X} )$
can be computed as the limit
$$ \varprojlim_{ \mathfrak{U} \subseteq \mathfrak{X} } \RGamma_{\Prism}( \mathfrak{U} ),$$
in the $\infty$-category $\calD( \Z_p )$, where $\mathfrak{U}$ ranges over the collection of affine open subsets of $\mathfrak{X}$.
Here it is essential to work with the $\infty$-category $\calD(\Z_p)$, rather than its homotopy category $D(\Z_p)$
(where the relevant limit is usually not defined).
\end{remark}

It follows from Remark \ref{remark:suffra} that the functor (\ref{equation:functor-up-to-homotopy}) is completely determined by
its restriction to the category of {\em affine} formal schemes, which we can identify with a covariant functor
$$ \{ \text{$p$-complete rings with bounded $p$-power torsion} \} \rightarrow \calD(\Z_p) \quad \quad R \mapsto \Prism_{R} = \RGamma_{\Prism}( \Spf(R) ).$$
To simplify the statements of theorems (and to slightly strengthen their conclusions), we will enlarge the domain of this functor to include
{\em all} commutative rings. In fact, it will be useful to make a more dramatic enlargement.

\begin{example}[Animated Commutative Rings]
Let $\Poly_{\Z}$ denote the category of finitely generated polynomial rings over $\Z$. An {\it animated commutative ring} is a functor of
$\infty$-categories
$$ \Poly_{\Z}^{\op} \rightarrow \{ \text{Spaces} \}$$
which preserves finite products. The collection of animated commutative rings can be organized into an $\infty$-category
which contains the ordinary category of commutative rings as a full subcategory, where we identify each commutative ring
$R$ with the functor
$$ \Poly_{\Z}^{\op} \rightarrow \{ \text{Sets} \} \quad \quad P \mapsto \{ \text{Ring homomorphisms $P \rightarrow R$} \}.$$
For a more detailed review of the theory of animated commutative rings, we refer the reader to \S \ref{ss:CAlgAnim}.
\end{example}

To every animated commutative ring $R$, we will associate a cochain complex $\Prism_{R}$, which we refer to as the
{\it absolute prismatic complex} of $R$ (Construction \ref{construction:absolute-prismatic-cohomology-general}). The
construction $R \mapsto \Prism_{R}$ determines a functor of $\infty$-categories
\begin{equation}\label{equation:functor-up-to-homotopy3}  \{ \text{Animated commutative rings} \} \rightarrow \calD(\Z_p),
\end{equation}
where neither the domain nor the codomain is an ordinary category. Let us describe one motivation for defining prismatic complexes in this generality:

\begin{remark}[Derived Descent]\label{remark:sugura}
For every animated commutative ring $R$, let us write $R \otimes^{L} \F_p$ for the coproduct of $R$ with the finite field $\F_p$ in the
$\infty$-category of animated commutative rings. We can then form a cosimplicial diagram of animated commutative rings
$$\xymatrix@R=50pt@C=50pt{ R \otimes^{L} \F_p \ar@<.4ex>[r] \ar@<-.4ex>[r] & R \otimes^{L} \F_p \otimes^{L} \F_p
\ar@<.8ex>[r] \ar[r] \ar@<-.8ex>[r] & \cdots }$$
In \S \ref{subsection:absolute-prismatic}, we show that $\Prism_{R}$ can be computed as the limit of the associated cosimplicial diagram
$$\xymatrix@R=50pt@C=50pt{ \Prism_{R \otimes^{L} \F_p} \ar@<.4ex>[r] \ar@<-.4ex>[r] & \Prism_{R \otimes^{L} \F_p \otimes^{L} \F_p}
\ar@<.8ex>[r] \ar[r] \ar@<-.8ex>[r] & \cdots, }$$
formed in the $\infty$-category $\calD(\Z_p)$ (Proposition \ref{proposition:derived-descent-absolute}). 
\end{remark}

\begin{warning}
If $R$ is an ordinary commutative ring which is $p$-torsion-free, then the tensor product $R \otimes^{L} \F_p$ can be identified with the usual quotient ring $R/pR$.
However, the iterated tensor products $R \otimes^{L} \cdots \otimes^{L} \F_p$ will never be ordinary commutative rings (except in the trivial case where $p$ is invertible in $R$). Consequently, to fully capture the phenomenon described in Remark \ref{remark:sugura}, the language of animated commutative rings is essential.
However, Remark \ref{remark:sugura} has a consequence which can be stated more concretely: if $R$ is a $p$-torsion-free commutative ring, then
$\Prism_{R}$ can be identified with the limit of the tower
$$ \cdots \rightarrow \Prism_{ R/p^{3}R} \rightarrow \Prism_{ R / p^2 R} \rightarrow \Prism_{R / pR};$$
see Corollary \ref{corollary:reduce-mod-n}.
\end{warning}

\begin{remark}
Let $R$ be an animated commutative ring. Then the absolute prismatic complex $\Prism_{ R \otimes^{L} \F_p}$ can be identified with
the ($p$-complete) derived de Rham complex $\widehat{\dR}_{R}$ (see Construction \ref{construction:derived-de-Rham} and Theorem \ref{theorem:deduce-comparison}). Consequently, Remark \ref{remark:sugura} provides a mechanism for reducing statements about absolute prismatic cohomology to
questions about (derived) de Rham cohomology. This mechanism will actually be useful in practice: we will apply it in \S \ref{section:first-chern}
to reduce the proof of one of our main results (Theorem \ref{theorem:syntomic-chern-class-isomorphism}) to a concrete calculation with
divided power envelopes (Theorem \ref{theorem:logarithm-sequence}).
\end{remark}
  
\subsection{The Bestiary}

In this paper, we will study several cohomological invariants of bounded $p$-adic formal schemes $\mathfrak{X}$. These invariants are
local in nature, and are therefore determined by their restriction to the case where $\mathfrak{X} = \Spf(R)$ is affine. To simplify the discussion,
we will focus primarily on the affine case and view our invariants as functors of $R$. As in \S \ref{subsection:cochain-animation}, it will be useful
to expand the scope of our definition to allow $R$ to be an arbitrary (animated) commutative ring. For the reader's convenience, we briefly
summarize some of the invariants that we consider and how they are related to one another.

\begin{itemize}
\item To every (animated) commutative ring $R$, we associate an {\it absolute prismatic complex} $\Prism_{R}$ and its twisted
counterparts $\Prism_{R}\{n\}$ (Construction \ref{construction:absolute-prismatic-cohomology-general}), which are objects of the
derived $\infty$-category $\calD( \Z_p)$. Each of these
complexes is equipped with a decreasing {\it Nygaard filtration} $\Fil^{\bullet}_{\Nyg} \Prism_{R}\{n\}$ (Construction \ref{construction:absolute-Nygaard-untwisted})
and a Frobenius morphism $\varphi\{n\}: \Fil^{n}_{\Nyg} \Prism_{R}\{n\} \rightarrow \Prism_{R}\{n\}$.

\item To every (animated) commutative ring $R$, we associate a {\it diffracted Hodge complex}
$\Omega_{R}^{\DHod}$ (Construction \ref{construction:diffracted-Hodge-integral}), which is an object of the derived $\infty$-category $\calD(R)$.
It is equipped with an increasing {\it conjugate filtration} $\Fil_{\bullet}^{\conj} \Omega_{R}^{\DHod}$ and an endomorphism $\Theta$ which we
refer to as the {\it Sen operator} (Remark \ref{remark:Sen-operator-integral}). The diffracted Hodge complex $\Omega_{R}^{\DHod}$ does not depend on the choice of
prime number $p$. However, we will primarily be interested in its $p$-completion, which we denote by 
$\widehat{\Omega}_{R}^{\DHod}$ and refer to as the {\it $p$-complete diffracted Hodge complex of $R$} (Construction \ref{construction:complete-diffracted-Hodge}).

\item To every (animated) commutative ring $R$, we associate an {\it absolute Hodge-Tate complex} $\overline{\Prism}_{R}$ and its
twisted counterparts $\overline{\Prism}_{R}\{n\}$ (Construction \ref{construction:absolute-HT}), which are objects of the derived $\infty$-category
$\calD(R)$. There is a comparison map $\Prism_{R}\{n\} \rightarrow \overline{\Prism}_{R}\{n\}$, given geometrically by
restriction to the Hodge-Tate divisor $\WCart^{\mathrm{HT}} \subset \WCart$. Moreover, there is also a fiber sequence
$\overline{\Prism}_{R}\{n\} \rightarrow \widehat{\Omega}_{R}^{\DHod} \xrightarrow{ \Theta+n} \widehat{\Omega}_{R}^{\DHod}$ (Remark \ref{remark:diffracted-vs-prismatic2}): that is, we can identify $\overline{\Prism}_{R}\{n\}$ with the ``$-n$-eigenspace'' for
the Sen operator on $\widehat{\Omega}_{R}^{\DHod}$.

\item To every (animated) commutative ring $R$ and every integer $n$, we associate {\it syntomic complex} $\RGamma_{\Syn}( \Spec(R), \Z_p(n) )$,
which is an object of the derived $\infty$-category $\calD( \Z_p)$. When $R$ is $p$-complete, $\RGamma_{\Syn}( \Spec(R), \Z_p(n) )$ defined as the fiber of the map
$( \varphi\{n\} - \iota): \Fil^{n}_{\Nyg} \Prism_{R}\{n\} \rightarrow \Prism_{R}\{n\}$ (see Construction \ref{construction:syntomic-complex},
and Construction \ref{construction:syntomic-complexes-general} for the case where $R$ is not assumed to be $p$-complete).

\item Let $S$ be a commutative ring. To every (animated) commutative $S$-algebra $R$, we associate the
($p$-complete) {\it derived de Rham complex} $\widehat{\dR}_{R/S}$ (Construction \ref{construction:derived-de-Rham}), which is
an object of the derived $\infty$-category $\calD(S)$. It is equipped with a decreasing {\it Hodge filtration}
$\Fil^{\bullet}_{\Hodge} \widehat{\dR}_{R/S}$. In the special case $S = \Z$, we denote the complex $\widehat{\dR}_{R/S}$ by
$\widehat{\dR}_{R}$. For every integer $n$, the absolute prismatic complex $\Prism_{R}\{n\}$ is equipped
with a {\it de Rham comparison map} $\gamma_{\Prism}^{\dR}\{n\}: \Prism_{R}\{n\} \rightarrow \widehat{\dR}_{R}$, intertwining
the Nygaard filtration on $\Prism_{R}\{n\}$ with the Hodge filtration on $\widehat{\dR}_{R}$ (see Construction \ref{construction:absolute-Nygaard-untwisted}).
\end{itemize}

The situation is partially summarized by the following diagram:
$$ \xymatrix@R=25pt@C=25pt{ & \RGamma_{\Syn}( \Spec(R), \Z_p(n) ) \ar[d] \ar[dr] \ar[r] & \RGamma_{\mathet}( \Spec(R[1/p]), \Z_p(n))  \\
\Fil^{n}_{\Hodge} \widehat{\dR}_{R} \ar[d] & \Fil^{n}_{\Nyg} \Prism_{R}\{n\} \ar[r]_-{ \varphi\{n\} } \ar[d] \ar[l] & 
\Prism_{R}\{n\} \ar[d] \\
\gr^{n}_{ \Hodge} \widehat{\dR}_{R} \ar[d]^{\sim} & \gr^{n}_{\Nyg} \Prism_{R}\{n\} \ar[r] \ar[l] \ar[d] & 
\overline{\Prism}_{R}\{n\} \ar[d] \\
L\widehat{\Omega}^{n}_{R}[-n] &  \Fil_{n}^{\conj}
\widehat{\Omega}^{\DHod}_{R} \ar[l] \ar[r] & \widehat{\Omega}^{\DHod}_{R}. }$$
Here $L \widehat{\Omega}^{n}_{R}$ denotes the (completed) $n$th exterior power of the absolute cotangent complex of $R$,
which can be identified (after shifting) both with $\gr^{n}_{\Hodge} \widehat{\dR}_{R}$ and $\gr_{n}^{\conj} \widehat{\Omega}^{\DHod}_{R}$.

If $(A,I)$ is a bounded prism with quotient ring $\overline{A} = A/I$, then the theory of relative prismatic cohomology
provides several related invariants. To every (animated) $\overline{A}$-algebra $R$, Construction~7.6 of \cite{prisms}
supplies a {\it relative prismatic complex} $\Prism_{R/A}$, which is an object of the derived $\infty$-category
$\calD(A)$. The Frobenius pullback $\varphi^{\ast} \Prism_{R/A}$ is equipped with a decreasing {\em relative Nygaard filtration} $\Fil^{\bullet}_{\Nyg} \varphi^{\ast} \Prism_{R/A}$
(see \S \ref{subsection:relative-Nygaard}, following \S 12 of \cite{prisms}) and a relative de Rham comparison map
$\varphi^{\ast} \Prism_{R/A} \rightarrow \widehat{\dR}_{R / \overline{A}}$, which intertwines the relative 
Nygaard filtration with the Hodge filtration and induces an isomorphism
$$ \overline{A} \otimes^{L}_{A} \varphi^{\ast} \Prism_{R/A} \simeq \widehat{\dR}_{R / \overline{A} }$$
(Corollary~15.4 of \cite{prisms}, which we recall here as Proposition \ref{proposition:de-Rham-comparison-relative}).
We denote the tensor product $\overline{A} \otimes_{A}^{L} \Prism_{R/A}$ by
$\overline{\Prism}_{R/A}$ and refer to it as the {\it relative Hodge-Tate complex} of $R$; this complex is equipped
with an increasing {\it conjugate filtration} $\Fil_{\bullet}^{\conj} \overline{\Prism}_{R/A}$. The relative and absolute theories
are related by a commutative diagram
$$ \xymatrix@R=25pt@C=20pt{ \Fil^{n}_{\Hodge} \widehat{\dR}_{R} \ar[dd] \ar[dr] & & \Fil^{n}_{\Nyg} \Prism_{R}\{n\} \ar[ll] \ar[rr]^{\varphi\{n\}} \ar[dr]^{\ast} \ar[dd] & & \Prism_{R}\{n\} \ar[dr]^{\ast} \ar[dd] & \\
& \Fil^{n}_{\Hodge} \widehat{\dR}_{R / \overline{A} } \ar[dd] & & \Fil^{n}_{\Nyg} \varphi^{\ast} \Prism_{R/A}\{n\} \ar[rr] \ar[ll] \ar[dd] & & 
\Prism_{R/A}\{n\} \ar[dd] \\
L\widehat{\Omega}^{n}_{R} \ar[dr] & & \gr^{n}_{\Nyg} \Prism_{R}\{n\} \ar[ll] \ar[dr] \ar[rr] & & \overline{\Prism}_{R}\{n\} \ar[dr]^{\ast} &  \\
& L \widehat{\Omega}^{n}_{R/\overline{A}} & & \Fil_{n}^{\conj} \overline{\Prism}_{R/A}\{n\} \ar[ll] \ar[rr] & & \overline{\Prism}_{R/A}\{n\}; }$$
here the diagonal maps with an asterisk are isomorphisms when the prism $(A,I)$ is perfect. 

\subsection{Completions}

Throughout this paper, we will need to consider completions of modules which are not necessarily finitely generated over commutative rings which are not necessarily Noetherian. For this purpose, it will be convenient to adopt the following convention:

\begin{definition}\label{definition:I-complete}
Let $R$ be a commutative ring, let $I \subseteq R$ be a finitely generated ideal, and let $M$ be an $R$-module. Then:
\begin{itemize}
\item[$(a)$] We will say that $M$ is {\it $I$-adically separated} if the intersection $\bigcap_{n \geq 0} I^{n} M$ is equal to $\{0\}$.
\item[$(b)$] We will say that $M$ is {\it $I$-complete} if, for every element $x \in I$, every short exact sequence of $R$-modules
$$ 0 \rightarrow M \rightarrow M' \rightarrow R[1/x] \rightarrow 0$$
admits a unique splitting.
\end{itemize}
In the special case where $I = (p)$ is the ideal generated by $p$, we say that an $R$-module $M$ is {\it $p$-complete} 
if it is $I$-complete, and {\it $p$-adically separated} if it is $I$-adically separated.
\end{definition}

\begin{remark}
Let $R$ be a commutative ring, let $I \subseteq R$ be a finitely generated ideal, and let $M$ be an $R$-module. The following conditions are equivalent:
\begin{itemize}
\item The $R$-module $M$ is both $I$-adically separated and $I$-complete.
\item The canonical map $M \rightarrow \varprojlim M / I^{n} M$ is an isomorphism.
\end{itemize}
See Proposition~3.4.2 of \cite{MR3379634}.
\end{remark}

\begin{warning}
Recall that it is standard to refer an $R$-module $M$ as {\it $I$-adically complete} if the canonical map $M \rightarrow \varprojlim M/I^{n}M$ is an isomorphism: that is, if it is both $I$-adically separated and $I$-complete in the sense of Definition \ref{definition:I-complete}. Beware that this is generally a stronger condition
than the $I$-completeness of $M$.
\end{warning}

\begin{remark}
Let $R$ be a commutative ring, let $I \subseteq R$ be a finitely generated ideal, and let $M$ be an $R$-module.
Condition $(b)$ of Definition \ref{definition:I-complete} is equivalent to the requirement that, for every element $x \in R$,
the abelian groups $\Ext^{n}_{R}( R[1/x], M)$ vanish for every integer $n$ (note that this is automatic for $n \geq 2$, since
the $R$-module $R[1/x]$ has projective dimension $\leq 1$). Moreover, it suffices to verify this condition
for a collection of generators of the ideal $I$.
\end{remark}

\begin{remark}
Let $R$ be a Noetherian commutative ring, let $I \subseteq R$ be an ideal, and let $M$ be a finitely generated $R$-module. If $M$ is
$I$-complete, then it is automatically $I$-adically separated.
\end{remark}

\begin{remark}\label{remark:subcategory-closure}
Let $R$ be a commutative ring and let $I \subseteq R$ be a finitely generated ideal. Then the collection of $I$-complete
$R$-modules is closed under the formation of kernels, cokernels, and extensions. In particular, the category of $I$-complete
$R$-modules is abelian and its inclusion into the category of all $R$-modules is exact. Beware that the analogous assertion is
not true for the subcategory of $I$-complete and $I$-adically separated $R$-modules (if $f: M \rightarrow N$ is a homomorphism
of $R$-modules which are $I$-complete and $I$-adically separated, then the cokernel $\coker(f)$ need not be $I$-adically separated).
\end{remark}

\begin{remark}
Let $f: R \rightarrow S$ be a homomorphism of commutative rings, let $I \subseteq R$ be a finitely generated ideal, and let $J = f(I)S$ be the ideal
of $S$ generated by the image of $I$. Then an $S$-module $M$ is $J$-complete ($J$-adically separated) if and only if it
is $I$-complete ($I$-adically separated) when viewed as an $R$-module via the homomorphism $f$. In particular,
the condition that an $S$-module $M$ is $p$-complete ($p$-adically separated) depends only on the underlying abelian group of $M$, and not on its $S$-module structure.
\end{remark}

\begin{notation}[The Complete Derived $\infty$-Category]\label{notation:complete-derived}
Let $R$ be a commutative ring and let $I \subseteq R$ be a finitely generated ideal. We will say that an object 
$M \in \calD(R)$ is {\it $I$-complete} if the cohomology group $\mathrm{H}^{n}(M)$ is $I$-complete for every integer $n$.
We will generally write $\widehat{\calD}(R)$ for the full subcategory of $\calD(R)$ spanned by those objects which are $I$-complete. In practice, we will almost always
use this notation in the case where $I = (p)$ is the ideal generated by $p$. There is one exception to this rule: if $(A,J)$ is a prism, we write $\widehat{\calD}(A)$
for the full subcategory of $\calD(A)$ spanned by those complexes which are $I$-complete where $I = J + (p)$ is the ideal generated by $p$ and $J$.
\end{notation}

\begin{warning}\label{warning:compare-completions}
Let $R$ be a commutative ring and let $I \subseteq R$ be a finitely generated ideal. It will be important to distinguish between three potentially different notions of $I$-completion.
\begin{itemize}
\item If $M$ is an $R$-module, then we will refer to the inverse limit $\varprojlim_{n \geq 0} M/ I^{n} M$ as the {\it separated $I$-completion of $M$}. It is
$I$-complete, $I$-adically separated, and is characterized by the following universal property: if $N$ is any other $R$-module which is $I$-complete and $I$-adically separated, then the canonical map $\Hom_R( \varprojlim_{n \geq 0} M / I^{n} M, N ) \rightarrow \Hom_{R}(M,N)$ is a bijection.

\item For every object $M \in \calD(R)$, there exists a morphism $f: M \rightarrow \widehat{M}$ where $\widehat{M}$ is $I$-complete
and, for every other object $I$-complete object $N \in \calD(R)$, composition with $f$ induces a bijection
$\Hom_{D(R)}( \widehat{M}, N ) \rightarrow \Hom_{D(R)}(M, N)$. The object $\widehat{M} \in \widehat{\calD}(R)$ is determined (up to canonical isomorphism) by this
universal property and depends functorially on $M$; we refer to $\widehat{M}$ as the {\it $I$-completion of $M$}. 

\item Let $M$ be an $R$-module, which we can regard as an object of $\calD(R)$ whose cohomology groups are concentrated in degree zero.
Beware that the $I$-completion $\widehat{M}$ generally has nonzero cohomology groups in negative degrees. However, the $0$th cohomology group
$\mathrm{H}^{0}( \widehat{M} )$ is an $I$-complete $R$-module, and the tautological map $f: M \rightarrow \widehat{M}$ induces an
$R$-module homomorphism $\mathrm{H}^{0}(f): M \rightarrow \mathrm{H}^{0}( \widehat{M} )$ with the following universal property: for
every $R$-module $N$ which is $I$-complete, composition with $\mathrm{H}^{0}(f)$ induces a bijection
$\Hom_{R}( \mathrm{H}^{0}( \widehat{M} ), N ) \rightarrow \Hom_{R}(M, N)$.
\end{itemize}

If $M$ is an $R$-module, we have tautological comparison maps 
$$ \widehat{M} \xrightarrow{u} \mathrm{H}^{0}( \widehat{M} ) \xrightarrow{v} \varprojlim_{n \geq 0} M/ I^{n} M,$$
where $u$ is a morphism in the derived $\infty$-category $\calD(R)$ and $v$ is an $R$-module homomorphism.
Beware that, in general, neither $u$ nor $v$ is an isomorphism. However, the homomorphism $v$ is always surjective,
and its kernel is equal to the intersection $\bigcap_{n \geq 0} I^{n} \mathrm{H}^{0}( \widehat{M} )$ (measuring
the failure of $\mathrm{H}^{0}( \widehat{M} )$ to be $I$-adically separated).
\end{warning}

\begin{remark}\label{remark:bounded-p-torsion}
In the situation of Warning \ref{warning:compare-completions}, suppose that $I = (p)$ is the ideal generated by $p$ and
that the $R$-module $M$ has bounded $p$-power torsion: that is, that the chain of $p$-power torsion submodules
$$M[p] \subseteq M[p^2] \subseteq M[p^3] \subseteq M[p^4] \subseteq \cdots$$
is eventually stable. In this case, the comparison maps
$$ \widehat{M} \rightarrow \mathrm{H}^{0}( \widehat{M} ) \twoheadrightarrow \varprojlim M/p^{n}M$$
are isomorphisms; that is, the $p$-completion of $M$ as an object of the derived $\infty$-category $\calD(R)$ can be
identified with the usual separated $p$-completion of $M$.
\end{remark}

Let $R$ be a commutative ring which is $p$-complete and $p$-adically separated. We let $\Spf(R) = ( \mathfrak{X}, \calO_{\mathfrak{X}} )$ denote the formal
spectrum of $R$ with respect to its $p$-adic topology. Here $\mathfrak{X}$ is the set of prime ideals $\mathfrak{p} \subset R$ which contain the prime number $p$.
It has a basis consisting of open sets $U_{f} = \{ \mathfrak{p} \in \mathfrak{X}: f \notin \mathfrak{p} \}$ where $f$ is an element of $R$, and its
structure sheaf is given by the formula $\calO_{ \mathfrak{X} }(U_f) = \varprojlim_{n} R[ 1/f ] / p^{n} R[1/f ]$. To guarantee that this construction
is well-behaved, we will generally need to assume that the commutative ring $R$ has bounded $p$-power torsion, so that $\calO_{ \mathfrak{X} }( U_f )$
is also the $p$-completion of $R[1/f]$.

\begin{definition}\label{definition:bounded-formal-scheme}
Let $\mathfrak{X}$ be a $p$-adic formal scheme. We say that $\mathfrak{X}$ is {\it bounded} if the structure sheaf $\calO_{\mathfrak{X}}$
has bounded $p$-power torsion.
\end{definition}

\begin{remark}
Let $\mathfrak{X}$ be a $p$-adic formal scheme. Then:
\begin{itemize}
\item If $\mathfrak{X}$ is locally Noetherian, then it is automatically bounded.
\item If $\mathfrak{X} = \Spf(R)$ is an affine $p$-adic formal scheme, then it is bounded if and only if the coordinate ring $R$ has bounded $p$-power torsion.
\item The condition that $\mathfrak{X}$ is bounded can be tested locally with respect to the Zariski topology.
\item Suppose that $\mathfrak{X} = \Spf(\Z_p) \times X$ is the formal completion of a scheme $X$ along the vanishing locus of $p$. If the structure sheaf $\calO_{X}$ has bounded $p$-power torsion, then $\mathfrak{X}$ is bounded.
\end{itemize}
\end{remark}

\subsection{Acknowledgments} 

Several of the core ideas of this paper (especially the discovery of $\WCart$ and the prismatic logarithm) were discovered and informally announced several years ago, and we apologize for the long delay in their eventual appearance. We are grateful to Johan de Jong, Vladimir Drinfeld,  H\'el\`ene Esnault, Dmitry Kubrak, Lars Hesselholt, Luc Illusie, Shizhang Li, Akhil Mathew, Shubhodip Mondal, Wiesia Niziol, Arthur Ogus, Sasha Petrov, Peter Scholze and Longke Tan for useful conversations and communications. Finally, many thanks to Akhil Mathew for comments on a preliminary version of this paper.

The first author was partially supported by the NSF (\#1801689, \#1952399, \#1840234), the Packard Foundation, and the Simons Foundation (\#622511).

\newpage \section{Breuil-Kisin Twists and the Prismatic Logarithm}\label{section:twist-and-log}

Let $(A,I)$ be a prism and let $\overline{A}$ denote the quotient ring $A/I$. We write $T_p( \overline{A}^{\times})$ for the $p$-adic Tate module of
the multiplicative group $\overline{A}^{\times}$ (Notation \ref{notation:R-flat}). Our goal in this section is to construct an invertible $A$-module $A\{1\}$, which we will
refer to as the {\it Breuil-Kisin twist of $A$}, together with a group homomorphism
$$ \log_{\Prism}: T_p( \overline{A}^{\times} ) \rightarrow A\{1\},$$
which we will refer to as the {\it prismatic logarithm}.

It will be convenient to view our prismatic logarithm as a composition of two more primitive operations. Let $u$ be an element of the Tate module
$T_p(\overline{A}^{\times} )$, which we can identify with a sequence of elements $u_0, u_1, u_2 \cdots \in \overline{A}$ satisfying
$u_0 = 1$ and $u_{n} = u_{n+1}^{p}$ for $n \geq 0$. It follows from an elementary $p$-adic convergence argument that $u$ can be
lifted uniquely to a sequence of elements $\widetilde{u}_0, \widetilde{u}_1, \widetilde{u}_2, \cdots \in A$ satisfying
$\widetilde{u}_{n} = \widetilde{u}_{n+1}^{p}$ for $n \geq 0$ (Proposition \ref{proposition:lift-flat}). Here the element $\widetilde{u}_0 \in A$ need not be equal to $1$.
However it is always a {\em rank one} element of the $\delta$-ring $A$: that is, it satisfies the identity $\delta_{A}( \widetilde{u}_0 ) = 0$
(this immediately implies the slightly weaker identity $\varphi_{A}( \widetilde{u}_0 ) = \widetilde{u}_0^{p}$, and the converse holds when $A$ is $p$-torsion-free).
Let $(1+I)_{\mathrm{rk}=1}$ denote the collection of rank $1$ elements $v \in A$ satisfying $v \equiv 1 \pmod{I}$, which we regard as a subgroup of the unit group
$A^{\times}$. The construction $u \mapsto \widetilde{u}_0$ then determines a group homomorphism $\rho: T_p( \overline{A}^{\times} ) \rightarrow (1 + I )_{\mathrm{rk}=1}$,
which we explain in greater detail in \S\ref{subsection:Tate-modules}. The prismatic logarithm $\log_{\Prism}: T_p( \overline{A}^{\times} ) \rightarrow A\{1\}$
will be defined as the composition of $\rho$ with another group homomorphism $(1+I)_{\mathrm{rk}=1} \rightarrow A\{1\}$, which (by abuse of notation)
we will also denote by $\log_{\Prism}$ and refer to as the {\it prismatic logarithm}.

Heuristically, one can think of the $A$-module $A\{1\}$ as an infinite tensor product
$$ I \otimes_{A} \varphi_{A}^{\ast}(I) \otimes_{A} \varphi_{A}^{2 \ast}(I) \otimes \cdots$$
of iterated pullbacks of $I$ along the Frobenius morphism $\varphi_{A}: A \rightarrow A$. If $v \in A$ is a rank one element
satisfying $v \equiv 1 \mod{I}$, then the prismatic logarithm $\log_{\Prism}(v)$ can be viewed as a suitably normalized version of the
the limit $$\lim_{\alpha \rightarrow 0} \frac{ u^{\alpha} - 1}{\alpha} = \lim_{n \to \infty} \frac{ u^{p^{n} } - 1}{p^{n} }.$$
To make sense of these heuristics, it will be convenient to first assume that the prism $(A,I)$ is {\it transversal}, meaning that
the quotient ring $\overline{A} = A/I$ is $p$-torsion-free (Definition \ref{definition:prism}). We give a brief overview of the theory of transversal prisms in \S\ref{subsection:transversal-prisms}.
Under the assumption that $(A,I)$ is transversal, we define the Breuil-Kisin twist $A\{1\}$ in \S\ref{subsection:BK-twist-transversal}
(Construction \ref{construction:twist-transversal-case}), and the prismatic logarithm $\log_{\Prism}$ in \S\ref{subsection:prismatic-log-transversal} (Construction \ref{construction:prismatic-logarithm}).

To handle the general case, we exploit the fact that every prism can be approximated by transversal prisms. More precisely, we show in \S\ref{subsection:approximation}
that for every prism $(A,I)$ there exists a homomorphism $(B,J) \rightarrow (A,I)$, where the prism $(B,J)$ is transversal (Proposition \ref{proposition:transversal-approximation}). Moreover,
the category of such transversal approximations admits pairwise coproducts, and is therefore sifted (Corollary \ref{siftedness}). In \S\ref{subsection:BK-log-general}, we apply these approximation results
to extend the definition of the Breuil-Kisin twist $A\{1\}$ and the prismatic logarithm $\log_{\Prism}$ to the case of an arbitrary prism $(A,I)$ (Definitions \ref{definition:twist-general} and \ref{definition:prismatic-logarithm-general}). 

In \S\ref{subsection:q-dR}, we specialize our constructions to the case where the prism $(A,I)$ is defined over the $q$-de Rham prism $( \Z_p[[q-1]], ( [p]_{q} ) )$. In this case,
the Breuil-Kisin twist $A\{1\}$ has a canonical generator $e_{A}$ on which the Frobenius acts by the formula $\varphi( e_A ) = e_A/ [p]_{q}$ (Notation \ref{notation:preferred-generator}).
Moreover, if $u \in A$ is a rank one element satisfying $u^{p} \equiv 1 \pmod{I}$, then the prismatic logarithm $\log_{\Prism}$ satisfies the identity
$\log_{\Prism}(u^{p} ) = \log_{q}(u) e_{A}$ where $\log_{q}(u)$ is the $q$-logarithm studied in \cite{cyctrace} (Corollary \ref{corollary:compare-log-bounded}). We will be particularly
interested in the case where the prism $(A,I)$ is crystalline, in which case we have $\log_{\Prism}( u^{p} ) = \log(u) e_A$, where $\log(u) = \sum_{n > 0} (-1)^{n-1} \frac{ (u-1)^{n} }{n}$ is the usual logarithmic series (Corollary \ref{corollary:crystalline-formula-for-log}).

\begin{remark}
The theory of Breuil-Kisin twists is discussed from stacky perspective in \cite[\S 4.9]{drinfeld-prismatic}, while the prismatic logarithm has been discussed from the same perspective in  \cite[Corollary 2.7.11]{DrinfeldFormalGroup}.
\end{remark}

\subsection{Transversal Prisms}\label{subsection:transversal-prisms}

In this section, we review the definition of a {\it transversal prism}, introduced by Ansch\"{u}tz and Le Bras in \cite{cyctrace}. For purposes of giving examples, it will be convenient to introduce a slightly more general notion.

\begin{definition}\label{definition:preprism}
A {\it preprism} is a pair $(A,I)$, where $A$ is a $\delta$-ring and $I \subseteq A$ is an ideal which is invertible as an $A$-module.
If $(A,I)$ and $(B,J)$ are preprisms, then a {\it morphism of preprisms} from $(A,I)$ to $(B,J)$ is a $\delta$-ring homomorphism
$f: A \rightarrow B$ satisfying $J = f(I) B$. We say that a preprism $(A,I)$ is {\it transversal} if the quotient ring $A/I$ is $p$-torsion-free.
\end{definition}

\begin{notation}\label{notation:construct-rho}
Let $(A,I)$ be a preprism. We write $\varphi: A \rightarrow A$ for the Frobenius lift supplied by the $\delta$-structure on $A$, and $\varphi^{\ast}(I)$ for the
invertible $A$-module obtained from $I$ by extending scalars along $\varphi$. Note that the composite map
$I \xrightarrow{\delta} A \twoheadrightarrow A/I$ is additive and $\varphi$-semilinear: the additivity follows from the identity
$$ \delta(x+y) = \delta(x) + \delta(y) - (p-1)! \sum_{i=1}^{p-1} \frac{ x^{i} }{i!} \frac{ y^{p-i} }{(p-i)!}$$
(applied to elements $x,y \in I$), and the $\varphi$-semilinearity from the identity
$$ \delta(xy) = \varphi(x) \delta(y) + \delta(x) y^{p}$$
(applied to elements $x \in A$ and $y \in I$). It follows that there is a unique $A$-linear map $\rho: \varphi^{\ast}(I) \rightarrow A/I$
for which the diagram of sets
$$ \xymatrix@R=50pt@C=50pt{ I \ar[d] \ar[r]^-{\delta} & A \ar@{->>}[d] \\
\varphi^{\ast}(I) \ar[r]^-{\rho} & A/I }$$
is commutative. 
\end{notation}

\begin{definition}\label{definition:prism}
A {\it prism} is a preprism $(A,I)$ which satisfies the following additional conditions:
\begin{itemize}
\item[$(1)$] The $A$-module homomorphism $\rho: \varphi^{\ast}(I) \rightarrow A/I$ of Notation \ref{notation:construct-rho} is surjective.
\item[$(2)$] The commutative ring $A$ is $p$-complete and $I$-complete.
\end{itemize}
We say that a prism $(A,I)$ is {\it transversal} if it is transversal as a preprism: that is, if the commutative ring $A/I$ is $p$-torsion-free.
\end{definition}

\begin{remark}
Recall that a prism $(A,I)$ is {\it bounded} if the quotient ring $A/I$ has bounded $p$-power torsion. Every transversal prism is bounded.
\end{remark}

\begin{remark}
Let $(A,I)$ be a preprism which satisfies condition $(1)$ of Definition \ref{definition:prism}. Then the morphism $\rho$ of
Notation \ref{notation:construct-rho} induces an isomorphism of invertible $(A/I)$-modules
$(A/I) \otimes_{A} \varphi^{\ast}(I) \simeq A/I$.
\end{remark}

\begin{remark}[Morphisms of Prisms]
Let $(A,I)$ and $(B,J)$ be prisms, and let $f: A \rightarrow B$ be a $\delta$-ring homomorphism. The following conditions are equivalent:
\begin{itemize}
\item[$(1)$] The homomorphism $f$ is a morphism of preprisms from $(A,I)$ to $(B,J)$. That is, the ideal $J \subseteq B$ is generated by the image $f(I)$.
\item[$(2)$] The homomorphism $f$ satisfies $f(I) \subseteq J$.
\end{itemize}
The implication $(1) \Rightarrow (2)$ is immediate, and the reverse implication is Proposition~1.5 of \cite{prisms}. If these conditions are satisfied,
then we say that $f$ is a {\it morphism of prisms} from $(A,I)$ to $(B,J)$.
\end{remark}

\begin{remark}\label{albose}
Let $A$ be a commutative ring and let $I \subseteq A$ be an invertible ideal. Suppose that the quotient ring $\overline{A} = A/I$ is $p$-torsion-free.
Then, for every integer $m$, the quotient $I^{m} / I^{m+1} \simeq \overline{A} \otimes_{A} I^{m}$ is an invertible $\overline{A}$-module, and is therefore also
$p$-torsion free. It follows by induction that for each $n \geq 0$, the quotient $A / I^{n}$ is $p$-torsion-free. If $A$ is $I$-complete, then $A$ is also $p$-torsion-free. In particular, if $(A,I)$ is a transversal prism, then the commutative ring $A$ is $p$-torsion-free.
\end{remark}

\begin{remark}\label{remark:transversal-description}
A prism $(A,I)$ is transversal if and only if it satisfies the following conditions:
\begin{itemize}
\item The commutative ring $A$ is $p$-torsion-free.
\item The natural map $u: I/pI \rightarrow A/pA$ is injective.
\end{itemize}
Moreover, if these conditions are satisfied, then the image of $u$ is an invertible ideal in the quotient ring $A/pA$.
\end{remark}

\begin{example}[The $q$-de Rham Prism]\label{example:q-prism}
Let $A = \Z_p[[ q-1]]$ be the completion of the polynomial ring $\mathbf{Z}[q]$ with respect to the ideal $(p, q-1)$.
We regard $A$ as equipped with the $\delta$-structure given by the Frobenius lift
$$ \varphi_{A}: A \rightarrow A \quad \quad q \mapsto q^{p}.$$
Then $(A, ( [p]_q) )$ is a transversal prism; note that the quotient $$A / ( [p]_q) \simeq \Z_p[ q ] / (1 + q + \cdots + q^{p-1} )$$ can be identified with the ring of integers of the local field $\Q_p( \zeta_p )$,
and is therefore $p$-torsion-free.
\end{example}

\begin{proposition}\label{proposition:preprism-to-prism}
Let $(A,I)$ be a transversal preprism. Then there exists a preprism homomorphism $f: (A,I) \rightarrow (B,J)$ with the following universal property:
\begin{itemize}
\item[$(1)$] The pair $(B,J)$ is a prism.
\item[$(2)$] For every prism $(C,K)$, precomposition with $f$ induces a bijection
$$ \{ \textnormal{Prism homomorphisms $(B,J) \rightarrow (C,K)$} \} \rightarrow \{ \textnormal{Preprism homomorphisms $(A,I) \rightarrow (C,K)$} \}$$
\end{itemize}
Moreover, the prism $(B,J)$ is transversal and $f$ induces an open immersion of affine schemes $\Spec( B / J + (p) ) \hookrightarrow \Spec( A / (I+(p)) )$.
\end{proposition}

\begin{proof}
Let $\rho: \varphi^{\ast}(I) \rightarrow A/I$ be as in Notation \ref{notation:construct-rho}. Since $\varphi^{\ast}(I)$ is a projective $A$-module, we can lift $\rho$
to an $A$-module homomorphism $\widetilde{\rho}: \varphi^{\ast}(I) \rightarrow A$. Let $U \subseteq \Spec(A)$ denote the Zariski open subset given by the nonvanishing locus of
$\widetilde{\rho}$. Then $U$ is affine (it can be described locally as the complement of the vanishing locus of single element of $A$). Write $U = \Spec(A')$ and let
$I' = IA'$ denote the ideal of $A'$ generated by $I$. Let $B = \varprojlim_{m,n} A' / (I'^{m} + (p^n) )$ be the separated $(I'+(p))$-completion of $A'$
(which coincides with $(I'+(p))$-completion, since $I'$ is an invertible ideal and $A' / I'$ is $p$-torsion-free). Our assumption that $(A,I)$ is transversal supplies exact sequences 
$$ 0 \rightarrow I \otimes_{A} (A'/ (I'^{m} + (p^n) )) \rightarrow A' /(I'^{m+1} + (p^n)) \rightarrow A'/(I' + (p^n)) \rightarrow 0.$$
Passing to the inverse limit over $m$ and $n$ (and noting that the transition maps are surjective), we deduce that the canonical map
$I \otimes_{A} B \rightarrow B$ is a monomorphism, whose image is the invertible ideal $J = IB$. By construction, the commutative ring $B$ is complete with respect to the ideal $J + (p)$,
and the quotient ring $B/J$ is the $p$-completion of $A'/IA'$, which is a Zariski localization of $A/IA$. In particular, $B/J$ is $p$-torsion-free, so that
$B$ is also $p$-torsion-free by virtue of Remark \ref{albose}.

Let $\varphi_{A}: A \rightarrow A$ denote the Frobenius lift supplied by the $\delta$-structure on $A$, so $\varphi_{A}(I) \subseteq I^{p} + (p)$.
It follows that, for every integer $k \geq 0$, we have $\varphi( I^{2k} + (p^k) ) \subseteq I^{pk} + (p^{k} )$. Consequently, $\varphi_{A}$ induces a morphism of schemes
$\Spec( A / I^{pk} + (p^k) ) \rightarrow \Spec( A / I^{2k} + (p^k) )$. Since $\varphi$ is a lift of Frobenius, this map is the identity at the level of topological spaces,
and therefore preserves the inverse image of the open subscheme $U \subseteq \Spec(A)$. We therefore obtain $\varphi_{A}$-semilinear maps
$$ \varphi_{B,k}: (B / J^{pk} + (p^k)) \rightarrow B / J^{2k} + (p^k).$$
Passing to the inverse limit over $k$, we obtain a map of commutative rings $\varphi_{B}: B \rightarrow B$. This map is a lift of the Frobenius on $B/pB$. Since $B$ is $p$-torsion-free, the morphism $\varphi_{B}$ determines a $\delta$-structure on $B$, which is uniquely characterized by the requirement that the tautological map $f: A \rightarrow B$ is a $\delta$-homomorphism. 

We now claim that the pair $(B,J)$ is a prism. Let $\rho_{B}: \varphi_{B}^{\ast}(J) \rightarrow B/J$ be as in Notation \ref{notation:construct-rho}; we wish to show that $\rho_{B}$ is surjective.
Since $B/J$ is $p$-complete, it will suffice to show that the the composite map $\varphi_{B}^{\ast}(J) \xrightarrow{ \rho_{B} } B/J \rightarrow B / (J+(p))$ is surjective. Equivalently, we must show
that the canonical map $\Spec(B/ J+(p)) \rightarrow \Spec(A)$ factors through the open subscheme $U \subseteq \Spec(A)$, which follows immediately from the construction.

We now prove $(2)$. Suppose we are given a prism $(C,K)$ and a morphism of preprisms $g: (A,I) \rightarrow (C,K)$; we wish to show that $g$ factors uniquely through $f$.
Note that our assumption that $(C,K)$ is a prism guarantees that the map $\Spec(C/K) \rightarrow \Spec(A/I) \hookrightarrow \Spec(A)$ factors through the open subscheme $U \subseteq \Spec(A)$.
Since $C$ is complete with respect to $K$, the ideal $K$ is contained in the Jacobson radical of $C$; it follows that the map $\Spec(C) \rightarrow \Spec(A)$ also factors through $U$.
Consequently, there exists a unique $A$-algebra homomorphism $A' \rightarrow C$. Since $C$ is complete with respect to the ideal
$K + (p) = I'C + (p)$, the ring homomorphism $g: A \rightarrow C$ factors uniquely as a composition $A \xrightarrow{f} B \xrightarrow{ \widehat{g} } C$.
We will complete the proof by showing that $\widehat{g}$ is a homomorphism of $\delta$-rings. Let us identify the $\delta$-structure on $A$ with a ring homomorphism
$s_{A}: A \rightarrow W_2(A)$ (which is a section of the restriction map), and define $s_B: B \rightarrow W_2(B)$ and $s_C: C \rightarrow W_2(C)$ similarly.
To show that $\widehat{g}$ is a homomorphism of $\delta$-rings, we must show that $s_{C} \circ \widehat{g}$ and $W_2( \widehat{g} ) \circ s_B$ coincide
(as ring homomorphisms from $B$ to $W_2(C)$). Let us regard $W_2(C)$ as an $A$-algebra via the ring homomorphism
$s_C \circ g= W_2(g) \circ s_B$ (where the equality follows from our assumption that $g$ is a $\delta$-ring homomorphism), so that we can regard
$s_C \circ \widehat{g}$ and $W_2( \widehat{g}) \circ s_B$ as $A$-algebra homomorphisms. We conclude by observing that an $A$-algebra homomorphism
$B \rightarrow W_2(C)$ is automatically unique, since $W_2(C)$ is complete with respect to the ideal $I + (p)$ and the map
$\Spec( W_2(C) ) \rightarrow \Spec(A)$ factors through the open subscheme $U \subseteq \Spec(A)$.
\end{proof}

\subsection{Breuil-Kisin Twists: Transversal Case}\label{subsection:BK-twist-transversal}

Our goal in this section is to construct the Breuil-Kisin twist $A\{1\}$ in the case of a transversal prism $(A,I)$. We begin by reviewing some elementary properties of transversal prisms (see \cite{cyctrace} for a closely related discussion).

\begin{lemma}[\cite{cyctrace}, Lemma~3.3]\label{lemma:easy-monomorphism}
Let $(A,I)$ be a transversal prism. Then, for every integer $r \geq 0$, the inclusion $I \hookrightarrow A$
induces a monomorphism $\theta: (\varphi_{A}^{r})^{\ast}(I) \rightarrow (\varphi_{A}^{r})^{\ast}(A) \simeq A$, whose cokernel is $p$-torsion-free.
\end{lemma}

\begin{proof}
Since $(\varphi_{A}^{r})^{\ast}(I)$ is an invertible $A$-module, it is $p$-complete and $p$-torsion-free, hence $p$-adically separated.
It will therefore suffice to show that $\theta$ induces a monomorphism $\overline{\theta}: (\varphi_{A}^{r})^{\ast}(I/pI) \rightarrow A/pA$.
This is clear, since $I/pI$ can be identified with an invertible ideal $\overline{I}$ of the quotient ring $A/pA$ (Remark \ref{remark:transversal-description}),
so its Frobenius pullback $(\varphi_{A}^{r})^{\ast}(I/pI)$ maps isomorphically to $\overline{I}^{p^{r}}$.
\end{proof}

\begin{notation}\label{notation:I-r}
Let $(A,I)$ be a transversal prism. For each integer $r \geq 0$, we will abuse notation by identifying
the pullback $(\varphi_{A}^{r})^{\ast}(I)$ with its image under the monomorphism of Lemma \ref{lemma:easy-monomorphism},
which is an invertible ideal of the commutative ring $A$. We let $I_{r}$ denote the product ideal $I \varphi^{\ast}_{A}(I) \cdots (\varphi^{r-1}_{A})^{\ast}(I) \subseteq A$
(so that $I_{r}$ is also an invertible ideal of $A$, abstractly isomorphic to the tensor product $I \otimes_{A} \varphi^{\ast}_{A}(I) \otimes_{A} \cdots \otimes_{A} (\varphi^{r-1}_{A})^{\ast}(I)$).
\end{notation}

\begin{remark}\label{remark:ideal-power-stuff}
Let $(A,I)$ be a transversal prism. Then the ideals $I^{ \frac{p^{r}-1}{p-1} }$ and $I_{r}$ have the same image in the quotient ring $A/pA$.
\end{remark}

\begin{remark}\label{remark:functoriality-image}
Let $(A,I) \rightarrow (B,J)$ be a morphism of prisms. For each $r \geq 0$, $J_r \subseteq B$ is the ideal generated by the image of $I_{r} \subseteq A$.
\end{remark}

\begin{lemma}\label{Iridealtransversal1}
Let $(A,I)$ be a transversal prism, and let $r \geq 0$ be an integer. Then the quotient ring $A/I_{r}$ is $p$-torsion-free.
\end{lemma}

\begin{proof}
We proceed by induction on $r$, the case $r = 0$ being vacuous. To carry out the inductive step, assume that $A/I_{r}$ is $p$-torsion-free.
We have an exact sequence of $A$-modules
$$ 0 \rightarrow I_{r} \otimes_{A} (A/(\varphi_{A}^{\ast})^{r}(I)) \rightarrow A/I_{r+1} \rightarrow A/I_{r} \rightarrow 0,$$
where the first term is $p$-torsion-free by virtue of Lemma \ref{lemma:easy-monomorphism} (and the fact that $I_{r}$ is invertible as an $A$-module).
It follows that $A/I_{r+1}$ is also $p$-torsion-free.
\end{proof}

\begin{remark}[\cite{cyctrace}, Lemma~3.4]\label{remark:inverse-system-compute-A}
Let $(A,I)$ be a transversal prism. Then $A$ is an inverse limit of the diagram of quotient rings
$$ \cdots \twoheadrightarrow A / I_{4} \twoheadrightarrow A / I_{3} \twoheadrightarrow A / I_{2} \twoheadrightarrow A / I_{1} = A / I.$$
Since $A$ and its quotients $A/I_{r}$ are $p$-complete and $p$-torsion free, this is equivalent to the assertion that $A/pA$ is an inverse limit of the diagram
$$\cdots \twoheadrightarrow A / (I_{4} + pA) \twoheadrightarrow A / (I_{3}+pA) \twoheadrightarrow A / (I_{2}+pA) \twoheadrightarrow A / (I+pA).$$
Setting $\overline{A} = A/pA$, taking $\overline{I}$ to be the image of $I$ in $\overline{A}$, and invoking Remark \ref{remark:ideal-power-stuff}, we are reduced to checking that
$\overline{A}$ is the limit of the diagram
$$\cdots \twoheadrightarrow \overline{A}/ \overline{I}^{ \frac{p^4-1}{p-1}}  \twoheadrightarrow \overline{A} / \overline{I}^{ \frac{ p^3-1}{p-1} } \twoheadrightarrow \overline{A} / \overline{I}^{ \frac{ p^2-1}{p-1} } \twoheadrightarrow \overline{A} / \overline{I}.$$
This is clear, since $\overline{A}$ is separated and complete with respect to $\overline{I}$.
\end{remark}

\begin{remark}\label{remark:Frobenius-pullback-system}
Let $(A,I)$ be a transversal prism. Then $A$ can be also realized as the inverse limit of the diagram
$$ \cdots \twoheadrightarrow A / \varphi_{A}^{\ast}(I_4) \twoheadrightarrow A / \varphi_{A}^{\ast}(I_3) \twoheadrightarrow A / \varphi_{A}^{\ast}(I_{2}) \twoheadrightarrow A / \varphi_{A}^{\ast}(I).$$
This follows from a slight variant of the argument supplied in Remark \ref{remark:inverse-system-compute-A}.
\end{remark}

\begin{lemma}\label{Iridealtransversal2}
Let $(A,I)$ be a transversal prism. For each $r > 0$, the natural map
$$ f: (\varphi_{A}^{r})^{\ast}(I) / I_{r+1} \rightarrow A / I_{r}$$
is a monomorphism, whose image is the principal ideal $(p) \subseteq A/I_{r}$.
\end{lemma}

\begin{proof}
We first show that the image of $f$ is divisible by $p$. Equivalently, we show that the ideal $(\varphi_{A}^{r})^{\ast}(I)$ is contained in $(p) + I_r$.
This is clear, since the image of $(\varphi_{A}^{r})^{\ast}(I)$ in the quotient ring $A/pA$ coincides with the image of $I^{p^r}$ (and therefore contains
the image of $I_r$ by virtue of Remark \ref{remark:ideal-power-stuff}). 

Since the quotient ring $A/I_{r}$ is $p$-torsion-free (Lemma \ref{Iridealtransversal1}), we can write $f = pf_0$ for some unique $A$-module homomorphism
$f_0: ( \varphi_{A}^{r} )^{\ast}(I) / I_{r+1} \rightarrow A/I_{r}$. We wish to show that $f_0$ is an isomorphism of $A$-modules. Since the domain and codomain of
$f_0$ are invertible modules over the quotient ring $A/I_{r}$, it will suffice to show that $f_0$ is surjective. We now observe that for each $x \in I$,
we have
\begin{eqnarray*} \varphi_{A}^{r}(x) & = & \varphi_{A}^{r-1}( x^{p} + p \delta(x) ) \\
& = & \varphi_{A}^{r-1}(x) \varphi_{A}^{r-1}(x)^{p-1} + p \varphi_{A}^{r-1}(\delta(x)) \\
& \in & \varphi_A^{r-1}(x) (I_{r-1} + (p) )^{p-1} + p \varphi_{A}^{r-1}(\delta(x)) \\
& \in & I_{r} + p (\varphi_{A}^{r-1})^{\ast}(I) + p\varphi_{A}^{r-1}(\delta(x)),
\end{eqnarray*}
so that $f_0(x) \equiv \varphi_{A}^{r-1}(\delta(x)) \pmod{ (\varphi_{A}^{r-1})^{ \ast}I }$. Since $(A,I)$ is a prism,
the elements $\{ \delta(x) \}_{x \in I}$ generate the unit ideal of $A$. It follows that the elements $\{ f_0(x) \}_{x \in I}$
generate the unit ideal in the quotient ring $A/(\varphi_{A}^{r-1})^{\ast}(I)$, hence also in $A/I_{r}$
(since $A$ is complete with respect to $(\varphi_{A}^{r-1})^{\ast}(I) \subseteq I + (p)$).
\end{proof}

\begin{corollary}[\cite{cyctrace}, Lemma~3.6]
Let $(A,I)$ be a transversal prism. Then, for each $r \geq 0$, the ideal $I_{r} \subseteq A$ is equal to the intersection of the ideals
$\{ (\varphi_{A}^{s})^{\ast}(I) \}_{0 \leq s < r}$. In other words, the canonical map
$$ A / I_{r} \rightarrow \prod_{0 \leq s < r} A / (\varphi_{A}^{s})^{\ast}(I)$$ is injective.
\end{corollary}

\begin{proof}
We proceed by induction on $r$ (the case $r=0$ being trivial). Assume that the conclusion is valid for some integer $r \geq 0$, and let
$x \in A$ belong to each of the ideals $\{ (\varphi_{A}^{s})^{\ast}(I) \}_{0 \leq s \leq r}$. Invoking our inductive hypothesis, we deduce that $x$ belongs to $I_{r}$,
and therefore also to the kernel of the map
$$(\varphi_{A}^{r})^{\ast}(I) \twoheadrightarrow  (\varphi_{A}^{r})^{\ast}(I) / I_{r+1} \xrightarrow{f} A / I_{r},$$
where $f$ is the monomorphism of Lemma \ref{Iridealtransversal2}. It follows that $x$ belongs to $I_{r}$.
\end{proof}

\begin{corollary}\label{corollary:transition-maps}
Let $(A,I)$ be a transversal prism, and let $r \geq 0$. Then:
\begin{itemize}
\item[$(1)$] The $A$-module $I_{r} / I_r^2$ is $p$-torsion-free.
\item[$(2)$] The canonical map $I_{r+1} / I_{r+1}^{2} \rightarrow I_{r} / I_{r}^{2}$ is divisible by $p$.
\item[$(3)$] The induced map $I_{r+1} / I_{r+1}^2 \xrightarrow{ 1/p} I_{r} / I_{r}^2$ is surjective, and therefore induces an isomorphism of invertible
$(A/I_{r})$-modules $I_{r+1} / I_{r} I_{r+1} \xrightarrow{1/p} I_{r} / I_{r}^2)$.
\end{itemize}
\end{corollary}

\begin{proof}
Since $I_{r}$ is an invertible ideal, assertion $(1)$ follows from the fact that $A / I_{r}$ is $p$-torsion-free (Lemma \ref{Iridealtransversal1}), and assertions $(2)$ and $(3)$ follow from Lemma \ref{Iridealtransversal2}.
\end{proof}

\begin{construction}[Breuil-Kisin Twists: Transversal Case]\label{construction:twist-transversal-case}
Let $(A,I)$ be a transversal prism. We let $A\{1\}$ denote the inverse limit of the diagram of surjections
$$  \cdots \twoheadrightarrow I_4/I_4^2 \twoheadrightarrow
I_3/I_3^2  \twoheadrightarrow I_{2} / I_{2}^2 \twoheadrightarrow I/I^2,$$
where the transition maps are given by division by $p$ (Corollary \ref{corollary:transition-maps}). We will refer to $A\{1\}$ as the {\it Breuil-Kisin twist} of $A$.
\end{construction}

\begin{proposition}\label{proposition:invertible-transversal-case}
Let $(A,I)$ be a transversal prism. Then the Breuil-Kisin twist $A\{1\}$ is an invertible $A$-module.
Moreover, for each integer $r \geq 0$, the tautological map $A\{1\} \rightarrow I_{r} / I_{r}^2$ induces an isomorphism
$$ A\{1\} / I_{r} A\{1\} \rightarrow I_{r} / I_{r}^{2}.$$
\end{proposition}

\begin{proof}
Since each $I_{r} / I_{r}^2$ is $p$-complete and $p$-torsion-free, the $A$-module $A\{1\}$ is also $p$-complete and $p$-torsion-free.
It will therefore suffice to show that each quotient $A\{1\} / p^{n} A\{1\}$ is an invertible module over $A/p^{n}A$ (see Lemma~4.11 of \cite{MR3572635}).
Proceeding by induction we can reduce to the case $n=1$. Letting $\overline{I}$ denote the image of $I$ in the quotient ring $A/pA$, we observe that $A\{1\} / p A\{1\}$ is the inverse limit of the tower
$$\{ I_{r} / (p I_{r} + I_{r}^2) \}_{r \geq 0} \simeq \{ \overline{I}^{ \frac{ p^{r} - 1}{p-1} } / \overline{I}^{ 2 \frac{ p^{r} - 1}{p-1} } \}_{r \geq 0}$$
of invertible modules over $(A/pA) / \overline{I}^{ \frac{ p^{r} - 1}{p-1}}$. Since $A/pA$ is $\overline{I}$-adically separated and $I$-complete, it will suffice to show that
each quotient $\overline{I}^{ \frac{ p^{r} - 1}{p-1} } / \overline{I}^{ 2 \frac{ p^{r} - 1}{p-1} }$ is an invertible module over 
$(A/pA) / \overline{I}^{ \frac{p^r - 1}{p-1}}$. This follows immediately from the invertibility of the ideal $\overline{I} \subseteq A/pA$ (Remark \ref{remark:transversal-description}).
\end{proof}

\begin{remark}[Functoriality]
Let $f: (A,I) \rightarrow (B,J)$ be a morphism of transversal prisms. Then there is a unique $A$-module homomorphism
$f\{1\}: A\{1\} \rightarrow B\{1\}$ having the property that, for each $r \geq 0$, the diagram
$$ \xymatrix@R=50pt@C=50pt{ A\{1\} \ar[r]^-{ f\{1\} } \ar@{->>}[d] & B\{1\} \ar@{->>}[d] \\
I_{r} / I_{r}^2 \ar[r]^-{f} & J_{r} / J_{r}^2 }$$
is commutative. Moreover, $f\{1\}$ induces an isomorphism of invertible $B$-modules $B \otimes_{A} A\{1\} \xrightarrow{\sim} B\{1\}$.
\end{remark}

\begin{construction}[The Frobenius on $A\{1\}$]\label{construction:Frobenius-transversal}
Let $(A,I)$ be a transversal prism. For each $A$-module $M$, let us denote the tensor product
$I^{-1} \otimes_{A} M$ by $I^{-1} M$. For each $r \geq 0$, we have a canonical isomorphism of $A$-modules $$\varphi_{A}^{\ast}( I_{r} ) \xrightarrow{\sim} I^{-1} I_{r+1}$$
Reducing modulo the ideal $\varphi_{A}^{\ast}(I_r)$ (which contains $I_{r+1}$) and invoking Proposition \ref{proposition:invertible-transversal-case}, we obtain isomorphisms
$$ \varphi_{A}^{\ast}( A\{1\} ) \otimes_{A} A/\varphi^{\ast}_{A}(I_r) \simeq I^{-1} A\{1\} \otimes_{A} 
A / \varphi^{\ast}_{A}(I_r)$$
Passing to the inverse limit over $r$ (and invoking Remark \ref{remark:Frobenius-pullback-system}), we obtain an isomorphism of invertible $A$-modules $\varphi_{A}^{\ast}( A\{1\}) \simeq I^{-1} A\{1\}$,
which we will identify with a $\varphi_{A}$-semilinear map $\varphi_{ A\{1\} }: A\{1\} \rightarrow I^{-1} A\{1\}$
and refer to as the {\it Frobenius morphism}.
\end{construction}

\subsection{The Prismatic Logarithm: Transversal Case}\label{subsection:prismatic-log-transversal}

Let $(A,I)$ be a transversal prism. Our goal in this section is to construct a group homomorphism
$\log_{\Prism}: (1+I)_{\mathrm{rk}=1} \rightarrow A\{1\}$
which we refer to as the {\it prismatic logarithm}. 
Our construction is motivated by the observation that if $x$ is a positive real number, the classical logarithm $\log(x)$ can be computed as the limit
$$ \lim_{\alpha \rightarrow 0} \frac{ x^{\alpha} - 1}{\alpha}.$$

\begin{proposition}\label{proposition:limit-defining-log}
Let $(A,I)$ be a transversal prism, and let $u \in 1 + I$ be a unit satisfying the identity $\delta(u) = 0$ (or, equivalently, $\varphi_{A}(u) = u^{p}$).
Then:
\begin{itemize}
\item[$(1)$] For each $r \geq 0$, the element $u^{p^{r}}$ is congruent to $1$ modulo the ideal $I_{r+1}$ of Notation \ref{notation:I-r}.
\item[$(2)$] For each $r > 0$, the division by $p$ map
$$\frac{1}{p}: I_{r+1} / I_{r+1}^{2} \twoheadrightarrow I_{r} / I_{r}^2$$
of Corollary \ref{corollary:transition-maps} carries the residue class of $u^{p^{r}}-1$ to the residue class of $u^{p^{r-1}} -1$.
\end{itemize}
\end{proposition}

\begin{proof}
We first prove $(1)$. By virtue of Lemma \ref{Iridealtransversal2}, it will suffice to show that $u^{p^{r}}$ is congruent to $1$ modulo the ideal $(\varphi_{A}^{s})^{\ast}(I)$ for $0 \leq s \leq r$.
This is clear: our assumption that $u$ has rank $1$ guarantees that $u^{p^{r}} = \varphi_{A}^{s}( u^{p^{r-s} } )$, and $u^{p^{r-s}}$ is congruent to $1$ modulo $I$ (since $u$ is congruent to $1$ modulo $I$).

We now prove $(2)$. Set $x = u^{p^{r-1}}$. It follows from $(1)$ that $x$ is congruent to $1$ modulo $I_{r}$, so that $1 + x + \cdots + x^{p-1}$ is congruent to $p$ modulo $I_{r}$.
Multiplying by the element $x-1 \in I_{r}$, we conclude that $x^{p} - 1$ is congruent to $p(x-1)$ modulo $I_{r}^2$.
\end{proof}

\begin{construction}[The Prismatic Logarithm, Transversal Case]\label{construction:prismatic-logarithm}
Let $(A,I)$ be a transversal prism, and let $u \in 1 + I$ satisfy $\delta(u) = 0$. It follows from Proposition \ref{proposition:limit-defining-log} there is a unique element of the Breuil-Kisin twist
$A\{1\}$ having image $u^{p^{r}}-1$ under each of the quotient maps $A\{1\} \twoheadrightarrow I_{r+1} / I_{r+1}^{2}$. We will denote this element by
$\log_{\Prism}(u)$ and refer to it as the {\it prismatic logarithm of $u$}. We regard the construction $u \mapsto \log_{\Prism}(u)$ as a function
$\log_{\Prism}: (1+I)_{\mathrm{rk}=1} \rightarrow A\{1\}$.
\end{construction}

\begin{proposition}\label{proposition:transversal-log-group-homomorphism}
Let $(A,I)$ be a transversal prism. Then the prismatic logarithm map
$$ \log_{\Prism}(u): (1+I)_{\mathrm{rk}=1} \rightarrow A\{1\}$$
is a group homomorphism (where the group structure on the domain is given by multiplication).
\end{proposition}

\begin{proof}
Let $u$ and $v$ be elements of $(1+I)_{\mathrm{rk}=1}$; we wish to verify the identity $\log_{\Prism}(uv) = \log_{\Prism}(u) + \log_{\Prism}(v)$.
To prove this, we must show that for each $r \geq 0$, we have
$$ (uv)^{p^{r}} - 1 \equiv (u^{p^{r}} - 1) + (v^{p^{r}}-1) \pmod{ I_{r+1}^{2} }.$$
We now note that difference between the right and left hand sides is equal to the product $(u^{p^{r}} -1 )(v^{p^r}-1)$, where each factor belongs to $I_{r+1}$ by virtue of
Proposition \ref{proposition:limit-defining-log}.
\end{proof}

\begin{remark}[Functoriality]\label{remark:functoriality-of-log}
Let $f: (B,J) \rightarrow (A,I)$ be a morphism of transversal prisms. Then the diagram of abelian groups
$$ \xymatrix@R=50pt@C=50pt{ (1+J)_{\mathrm{rk}=1} \ar[r]^-{f} \ar[d]^{\log_{\Prism}} & (1+I)_{\mathrm{rk}=1} \ar[d]^{ \log_{\Prism} } \\
B\{1\} \ar[r]^-{ f\{1\} } & A\{1\} }$$
is commutative.
\end{remark}

\subsection{Approximations by Transversal Prisms}\label{subsection:approximation}

We now show that there exists a good supply of transversal prisms.

\begin{proposition}\label{proposition:transversal-approximation}
Let $(A,I)$ be a prism. Then there exists a prism homomorphism $(B,J) \rightarrow (A,I)$, where $(B,J)$ is transversal.
\end{proposition}

The proof of Proposition \ref{proposition:transversal-approximation} will require some preliminaries.

\begin{lemma}\label{lemma:approximate-projective-module}
Let $R$ be a commutative ring and let $M$ be a projective $R$-module of finite rank. Then there exists a ring homomorphism $R_0 \rightarrow R$ and
an isomorphism of $R$-modules $M \simeq M_0 \otimes_{R_0} R$, where $M_0$ is a projective $R_0$-module of finite rank and $R_0$ is a smooth $\Z$-algebra.
\end{lemma}

\begin{proof}
Since $M$ is projective of finite rank, there exists an integer $n \geq 0$ and an idempotent $R$-module homomorphism $e: R^{n} \rightarrow R^{n}$ whose image is isomorphic to $M$.
Let $X$ be the scheme representing the functor
$$ S \mapsto \{ f \in \End_{S}(S^n): f^2 = f \} ).$$
Then $X$ is a closed subscheme of the affine space $\mathbf{A}^{2n}$ (over $\Z$), and is therefore affine. A standard calculation shows that $X$ is smooth over $\Z$,
and can therefore be identified with the spectrum $\Spec(R_0)$, where $R_0$ is a flat $\Z$-algebra. The identity map $\id: R_0 \rightarrow R_0$ corresponds to
an idempotent endomorphism of the module $R_0^{n}$, whose image $M_0$ is a projective $R_0$-module of finite rank. 
The endomorphism $e \in \End_{R}(R^n)$ corresponds to an $R$-valued point of $X$, which we can identify with a ring homomorphism $R_0 \rightarrow R$. By construction, we have an isomorphism
$M \simeq M_0 \otimes_{R_0} R$.
\end{proof}

\begin{lemma}\label{lemma:flatness-lemma}
Let $R$ be a Noetherian ring and let $A$ be the free $p$-complete $\delta$-ring generated by $R$. 
Suppose that $R$ is $p$-torsion-free and that the quotient $R/pR$ is a smooth $\F_p$-algebra. Then $A$ is flat over $R$.
\end{lemma}

\begin{proof}
Choose a surjective $R$-algebra homomorphism $f: B \twoheadrightarrow A$, where $B$ is the $p$-completion of a polynomial ring (possibly on infinitely many generators) over $R$.
Then $B$ is flat over $R$ (since it is $p$-completely flat and $R$ is Noetherian). We will complete the proof by showing that $f$ has a right inverse $g: A \rightarrow B$ in the category of $R$-algebras
(so that $A$ is a retract of $B$). Note that, for any commutative ring $C$, we can identify ring homomorphisms $A \rightarrow C$ with $\delta$-ring homomorphisms $A \rightarrow W(C)$, or equivalently with
ring homomorphisms $R \rightarrow W(C)$. Under this identification, the identity map $\id_{A}$ corresponds to a ring homomorphism $\gamma: R \rightarrow W(A)$, given by a compatible
sequence of ring homomorphisms $\{ \gamma_n: R \rightarrow W_n(A) \}_{n > 0}$. To construct a ring homomorphism $g: A \rightarrow B$ which is right inverse to $f$, we must factor
$\gamma$ a composition $R \xrightarrow{ \widehat{\gamma} } W(B) \xrightarrow{ W(f) } W(A)$: that is, we must provide a compatible sequence of ring homomorphisms
$\{ \overline{\gamma}_n: R \rightarrow W_n(B) \}_{n > 0}$ satisfying $W_n(f) \circ \overline{\gamma}_n = \gamma_n$. We will show that there exists such a sequence for which
$\overline{\gamma}_{1}: R \rightarrow W_1(B) = B$ is the tautological map (given by the $R$-algebra structure on $B$), so that $g$ is a morphism of $R$-algebras.
To prove this, it suffices to solve a sequence of lifting problems
$$ \xymatrix@R=50pt@C=50pt@C=100pt{ & W_{n+1}(B) \ar[d] \\
R \ar[r]_-{ (\gamma_{n+1}, \overline{\gamma}_{n})} \ar@{-->}[ur]^{ \overline{\gamma}_{n+1} } & W_{n+1}(A) \times_{ W_{n}(A) } W_{n}(B). }$$

Let $I \subseteq B$ be the kernel of the ring homomorphism $f$. Then $I$ is $p$-complete and $p$-torsion-free (since $B$ is $p$-torsion-free), and can therefore
be identified with the limit $\varprojlim_{n} I/p^{n}I$. Moreover, the fiber product $W_{n+1}(A) \times_{ W_n(A)} W_n(B)$ can be identified with the quotient of $W_{n+1}(B)$
by the ideal $V^{n}(I)$, where $V$ denotes the Verschiebung. It follows that $W_{n+1}(B)$ can be realized as the limit of a tower of surjective ring homomorphisms
$$ \cdots  \twoheadrightarrow W_{n+1}(B) / V^{n}( p^2 I) \twoheadrightarrow W_{n+1}(B) / V^{n}( pI) \twoheadrightarrow W_{n+1}(B) / V^{n}(I) \simeq W_{n+1}(A) \times_{ W_{n}(A) } W_{n}(B).$$
We are therefore reduced to solving a sequence of lifting problems
$$ \xymatrix@R=50pt@C=50pt{ & W_{n+1}(B) / V^n( p^{m+1} I) \ar[d] \\
R \ar[r] \ar@{-->}[ur] & W_{n+1}(B) / V^{n}( p^m I). }$$
This is possible by virtue of our assumption that $R/pR$ is a smooth $\F_p$-algebra, since the right vertical map exhibits $W_{n+1}(B) / V^n( p^{m+1} I)$ as a square-zero extension of $W_{n+1}(B) / V^n( p^m I)$ by a $p$-torsion module.
\end{proof}

\begin{proof}[Proof of Proposition \ref{proposition:transversal-approximation}]
Let $(A,I)$ be a prism. Then $I$ is a projective $A$-module of rank $1$. By virtue of Lemma \ref{lemma:approximate-projective-module},
there exists a ring homomorphism $u_0: B_0 \rightarrow A$ and an $A$-module isomorphism
$v: A \otimes_{B_0} J_0 \simeq I$, where $B_0$ is a smooth $\Z$-algebra and $J_0$ is an invertible
$B_0$-module. Set $B_1 = \Sym^{\ast}_{B_0}(J_0)$ and let $J_1 \subseteq B_1$ be the invertible ideal given by
$\Sym^{> 0}_{B_0}(J_0)$. Then $v$ is classified by a homomorphism of $B_0$-algebras $u_1: B_1 \rightarrow A$
satisfying $I = J_1 A$. Let $B_2$ denote the free $p$-complete $\delta$-ring generated by $B_1$, and set
$J_2 = J_1 B_2$. Since $B_2$ is a flat $B_1$-algebra (Lemma \ref{lemma:flatness-lemma}), 
the ideal $J_2 \subseteq B_2$ is invertible. Moreover, the quotient $B_2 / J_2$ is flat over over the smooth $\Z$-algebra $B_1 / J_1 \simeq B_0$, and is therefore
$p$-torsion-free. We can therefore regard $(B_2, J_2)$ as a transversal preprism. Invoking the universal property of $B_2$, we deduce that
the ring homomorphism $u_1$ factors uniquely as a composition $B_1 \rightarrow B_2 \xrightarrow{u_2} A$, where $u_2$ is a preprism homomorphism from $(B_2, J_2)$ to $(A,I)$. Applying Proposition \ref{proposition:preprism-to-prism}, we conclude that $u_2$ factors as a composition $(B_2, J_2) \rightarrow (B,J) \rightarrow (A,I)$, where $(B,J)$ is a transversal prism.
\end{proof}

Recall that a morphism of prisms $\theta: (A,I) \rightarrow (B,J)$ is said to be {\it flat} if the underlying ring homomorphism $\overline{\theta}: A/I \rightarrow B/J$ is $p$-completely flat,
and {\it faithfully flat} if the morphism $\overline{\theta}$ is $p$-completely faithfully flat.

\begin{remark}\label{remark:flatness-and-transversality}
Let $\theta: (A,I) \rightarrow (B,J)$ be a flat morphism of prisms. If the prism $(A,I)$ is transversal, then so is $(B,J)$. The converse holds if
$\theta$ is faithfully flat.
\end{remark}

\begin{proposition}\label{proposition:transversal-coproduct}
Let $(A,I)$ be a transversal prism. Then, for every bounded prism $(B,J)$, there exists a coproduct $(A,I) \coprod (B,J)$ in the category of prisms.
Moreover, the canonical map $(B,J) \rightarrow (A,I) \coprod (B,J)$ is flat.
\end{proposition}

\begin{proof}
Since $(A,I)$ is bounded, the category of bounded prisms over $(A,I)$ forms a stack with respect to the
\'{e}tale topology on the affine scheme $\Spec( A / (I + pA) )$. We may therefore assume without loss of generality that $I = (d)$ is a principal ideal generated by a distinguished element $d \in A$. Let $C_0$ denote the completion of the tensor product $A \otimes_{\Z} B$ with respect to the ideal generated by $J$ and $p$ (which coincides with the separated completion,
since the quotient $A \otimes_{\Z} (B/J)$ has bounded $p$-torsion). Since $A/dA$ is flat over $\Z$, the image of $d$ in $C_0$ is $(J,p)$-completely regular relative to $B$.
Invoking Proposition~3.13 of \cite{prisms}, we conclude that there exists a universal $\delta$-homomorphism $f: C_0 \rightarrow C$, where $(C,JC)$ is a prism
for which $JC$ contains the image of $d$. Moreover, the map of prisms $(B,J) \rightarrow (C,JC)$ is flat. It follows immediately
from the construction that $(C,JC)$ is a coproduct of $(A,I)$ and $(B,J)$ in the category of prisms (see Lemma~3.5 of \cite{prisms}). 
\end{proof}

\begin{remark}
In the first conclusion of Proposition \ref{proposition:transversal-coproduct} (i.e., for the existence of coproducts), the hypotheses on the prisms $(A,I)$ and $(B,J)$ can be removed
by working in the more general setting of {\em animated} prisms.
\end{remark}

\begin{corollary}\label{corollary:transversal-coproduct}
Let $(A,I)$ and $(B,J)$ be transversal prisms. Then there exists a coproduct $(A,I) \coprod (B,J)$ in the category of prisms, which is also transversal.
\end{corollary}

\begin{proof}
Combine Proposition \ref{proposition:transversal-coproduct} with Remark \ref{remark:flatness-and-transversality}.
\end{proof}

\begin{corollary}\label{siftedness}
Let $(A,I)$ be a prism and let $\calC$ be the category of transversal prisms $(A_0, I_0)$ equipped with a morphism of prisms $f: (A_0, I_0) \rightarrow (A,I)$. Then $\calC$ is sifted (in the $\infty$-categorical sense). 
In particular, the nerve of $\calC$ is weakly contractible.
\end{corollary}

\begin{proof}
The category $\calC$ is nonempty (Proposition \ref{proposition:transversal-approximation}) and admits pairwise coproducts (Corollary \ref{corollary:transversal-coproduct}).
\end{proof}

\begin{proposition}\label{proposition:characterize-transversal}
Let $(A,I)$ be a transversal prism with $A \neq 0$. Then, for every bounded prism $(B,J)$, the map of prisms
$(B,J) \rightarrow (A,I) \coprod (B,J)$ is faithfully flat.
\end{proposition}

\begin{proof}
Let $(C,K)$ denote the coproduct $(A,I) \coprod (B,J)$. Proposition \ref{proposition:transversal-coproduct}
guarantees that the map $B \rightarrow C$ is $(p,J)$-completely flat; we wish to show that it is $(p,J)$-completely faithfully flat, i.e., that the map $\mathrm{Spec}(C) \to \mathrm{Spec}(B)$ hits every point $x$ of the closed set $\mathrm{Spec}(B/(p,J)) \subset \mathrm{Spec}(B)$. Fix one such point $x$, and let $k$ be a perfect field containing $\kappa(x)$. The map $\overline{b}:B \to \kappa(x) \to k$ refines uniquely to a $\delta$-map $b:B \to W(k)$. Properties of distinguished elements in the $\delta$-ring $W(k)$ (see \cite[Lemma 2.28, Lemma 3.8]{prisms}) show that this map yields a map $\tilde{b}:(B,J) \to (W(k),(p))$ of prisms. On the other hand, as $\mathrm{Spec}(A/(p,I)) \neq \emptyset$ by assumption,  we can, after possibly enlarging $k$, find a map $\overline{a}:A \to A/(p,I) \to k$ of commutative rings. Again, this map refines uniquely to a $\delta$-map $a:A \to W(k)$, which then lifts to a map $\tilde{a}:(A,I) \to (W(k),(p))$ of prisms. The coproduct of $\tilde{a}$ and $\tilde{b}$ then induces a map $(C,K) \to (W(k),(p))$ of prisms. The induced map $\mathrm{Spec}(k) = \mathrm{Spec}(W(k)/(p)) \to \mathrm{Spec}(C/(p,K)) \subset \mathrm{Spec}(C)$ yields a point of $\mathrm{Spec}(C)$ that lifts the given point $x \in \mathrm{Spec}(B)$ by construction, proving the proposition.
\end{proof}

\begin{corollary}\label{corollary:characterize-transversal}
Let $(A,I)$ be a bounded prism. The following conditions are equivalent:
\begin{itemize}
\item[$(a)$] The prism $(A,I)$ is transveral.
\item[$(b)$] For every bounded prism $(B,J)$, there exists a coproduct $(A,I) \coprod (B,J)$ in the category of prisms, which is flat over $(B,J)$.
\item[$(c)$] There exists a nonzero transversal prism $(B,J)$ for which the coproduct $(A,I) \coprod (B,J)$ is flat over $(B,J)$.
\end{itemize}
\end{corollary}

\begin{proof}
The implication $(a) \Rightarrow (b)$ follows from Proposition \ref{proposition:transversal-coproduct} and the implication
$(b) \Rightarrow (c)$ from the existence of a nonzero transversal prism (for example, the $q$-de Rham prism of Example \ref{example:q-prism}).
To prove that $(c)$ implies $(a)$, suppose that $(B,J)$ is a transversal prism. If the coproduct $(A,I) \coprod (B,J)$ is flat over $(B,J)$,
then it is also transversal (Remark \ref{remark:flatness-and-transversality}). If $B \neq 0$, then $(A, I) \coprod (B,J)$ is faithfully flat over $(A,I)$,
so that $(A,I)$ is also transversal (again by Remark \ref{remark:flatness-and-transversality}).
\end{proof}

\begin{remark}
Heuristically, Corollary \ref{corollary:characterize-transversal} asserts that a prism $(A,I)$ is transversal if and only if it is flat over
the (non-existent) initial object of the category of prisms. See Corollary \ref{corollary:flatness-v-transversality} for a more precise formulation of this heuristic.
\end{remark}

\subsection{Extension to General Prisms}\label{subsection:BK-log-general}

We now apply the results of \S\ref{subsection:approximation} to extend the definition of the Breuil-Kisin twist $A\{1\}$ and the prismatic logarithm
$\log_{\Prism}: (1+I)_{\mathrm{rk}=1} \rightarrow A\{1\}$ to the case of an arbitrary prism $(A,I)$.

\begin{proposition}\label{proposition:existence-of-twist}
Let $(A,I)$ be a prism, and let $\calC$ denote the category whose objects are transversal prisms $(A_0, I_0)$ equipped with a morphism of prisms $f: (A_0, I_0) \rightarrow (A,I)$.
Then the diagram
$$ ( (A_0, I_0) \in \calC ) \mapsto A \otimes_{A_0} A_0\{1\}$$
has a direct limit in the category of $A$-modules, which we will denote by $A\{1\}$. Moreover, for every object $(A_0, I_0) \in \calC$, the tautological map $A \otimes_{A_0} A_0\{1\} \rightarrow A\{1\}$ is an isomorphism.
\end{proposition}

\begin{proof}
Since the diagram $\{ A \otimes_{A_0} A_0\{1\} \}_{(A_0, I_0) \in \calC}$ has invertible transition maps, it suffices to observe that the category $\calC$ has weakly contractible nerve, which follows
from Corollary \ref{siftedness}.
\end{proof}

\begin{definition}\label{definition:twist-general}
Let $(A,I)$ be a prism. We will refer to the $A$-module $A\{1\}$ of Proposition \ref{proposition:existence-of-twist} as the {\it Breuil-Kisin twist} of $A$. 
\end{definition}

\begin{remark}
Let $(A,I)$ be a prism. Concretely, the Breuil-Kisin twist $A\{1\}$ can be realized as the tensor product $A \otimes_{A_0} A_0\{1\}$, where $(A_0, I_0)$ is any transversal prism equipped with a map
$f: (A_0, I_0) \rightarrow (A, I)$, and $A_0\{1\}$ denotes the $A_0$-module given by Construction \ref{construction:twist-transversal-case}. In particular, if the prism $(A,I)$ is already transversal,
then we can take $f$ to be the identity morphism $\id_{A}$; in this case, Definition \ref{definition:twist-general} reduces to Construction \ref{construction:twist-transversal-case}
\end{remark}

\begin{notation}\label{notation:general-twists}
Let $(A,I)$ be an arbitrary prism. Then the Breuil-Kisin twist $A\{1\}$ is an invertible $A$-module (this follows immediately from Proposition \ref{proposition:invertible-transversal-case}). For every integer $n$, we let
$A\{n\}$ denote the $n$th tensor power of $A\{1\}$ relative to $A$. More generally, for every $A$-module $M$,
we let $M\{n\}$ denote the tensor product $A\{n\} \otimes_{A} M$.
\end{notation}

\begin{remark}[Functoriality]\label{remark:BK-functoriality}
The Breuil-Kisin twist $A\{1\}$ of Definition \ref{definition:twist-general} depends functorially on the prism $(A,I)$. More precisely, to every morphism of prisms $(A,I) \rightarrow (B,J)$, one
can associate a canonical isomorphism $B \otimes_{A} A\{1\} \simeq B\{1\}$, compatible with composition.
\end{remark}

\begin{remark}
Let $(A,I)$ be a prism. Then the inverse Breuil-Kisin twist $A\{-1\}$ has a geometric interpretation: it can be identified with the relative prismatic cohomology group $\mathrm{H}^2_{\Prism}( \mathbf{P}^1_{A} / A )$. See Variant \ref{variant:projectivelineBKtwist}.
\end{remark}

\begin{remark}\label{remark:twist-modulo-I}
For every prism $(A,I)$, there is a canonical isomorphism of invertible $(A/I)$-modules $\beta: A\{1\} / I A\{1\} \simeq I / I^2$. These isomorphisms are uniquely determined by the requirement that they depend functorially on the prism $(A,I)$ and recover the isomorphism of Proposition \ref{proposition:invertible-transversal-case} in the case where $(A,I)$ is transversal.
\end{remark}

\begin{remark}
Let $(A,I)$ be a prism. It follows from Remark \ref{remark:twist-modulo-I} that the Breuil-Kisin twist $A\{1\}$ is trivial (that is, it is isomorphic to $A$ as an $A$-module) if and only if the ideal $I$ is principal (since
$A$ is $I$-complete, both are equivalent to the triviality of $I/I^2$ as an invertible $(A/I)$-module).
\end{remark}

\begin{remark}[The Frobenius on $A\{1\}$]\label{remark:Frobenius-on-twist}
For every prism $(A,I)$, there is a canonical $\varphi_{A}$-semilinear map $\varphi_{A\{1\} }: A\{1\} \rightarrow I^{-1} A\{1\}$, which induces an $A$-linear isomorphism $\varphi_{A}^{\ast}( A\{1\} ) \simeq I^{-1} A\{1\}$. The morphisms
$\varphi_{A\{1\}}$ are uniquely determined by the requirement that they depend functorially on the prism $(A,I)$
and are given by Construction \ref{construction:Frobenius-transversal} in the case where $(A,I)$ is transversal. Passing to tensor powers, we obtain a $\varphi_{A}$-semilinear map
$\varphi_{A\{n\} }: A\{n\} \rightarrow I^{-n} A\{n\}$ for each integer $n \in \Z$.
\end{remark}


\begin{proposition}\label{proposition:prismatic-logarithm-extension}
Let $(A,I)$ be a prism. Then there is a unique function
$$ \log_{\Prism}: (1+I)_{\mathrm{rk}=1} \rightarrow A\{1\}$$
with the following property: for every morphism of prisms $f: (B,J) \rightarrow (A,I)$ where $(B,J)$ is transversal, the diagram
$$ \xymatrix@R=50pt@C=50pt{ (1+J)_{\mathrm{rk}=1} \ar[r]^-{f} \ar[d]^{\log_{\Prism}} & (1+I)_{\mathrm{rk}=1} \ar[d]^{ \log_{\Prism} } \\
B\{1\} \ar[r]^-{ f\{1\} } & A\{1\} }$$
is commutative; here the left vertical map is given by Construction \ref{construction:prismatic-logarithm}.
\end{proposition}

\begin{definition}\label{definition:prismatic-logarithm-general}
Let $(A,I)$ be a prism. We will refer to the map $\log_{\Prism}: (1+I)_{\mathrm{rk}=1} \rightarrow A\{1\}$
of Proposition \ref{proposition:prismatic-logarithm-extension} as the {\it prismatic logarithm}.
\end{definition}

\begin{example}
Let $(A,I)$ be a transversal prism. Then the prismatic logarithm $\log_{\Prism}: (1+I)_{\mathrm{rk}=1} \rightarrow A\{1\}$
of Construction \ref{construction:prismatic-logarithm} satisfies the conclusion of Proposition \ref{proposition:prismatic-logarithm-extension}
(Remark \ref{remark:functoriality-of-log}), and therefore coincides with the prismatic logarithm $\log_{\Prism}$ of Definition \ref{definition:prismatic-logarithm-general}.
\end{example}

The proof of Proposition \ref{proposition:prismatic-logarithm-extension} will make use of the following:

\begin{lemma}\label{lemma:add-rank-unit}
Let $(B,J)$ be a bounded prism. Then there exists a faithfully flat morphism of prisms $f: (B,J) \rightarrow (C,K)$
and a rank $1$ unit $w \in (1+K)_{\mathrm{rk}=1}$ with the following universal property: for every morphism of prisms
$(B,J) \rightarrow (A,I)$, evaluation on $w$ induces a bijection
$$ \{ \textnormal{$B$-linear $\delta$-homomorphisms $C \rightarrow A$} \} \rightarrow (1+I)_{\mathrm{rk}=1}.$$
\end{lemma}

\begin{proof}
Let $C_0$ denote Laurent polynomial ring $B[w^{\pm 1}]$, equipped with the $\delta$-structure determined by the
$\delta$-structure on $B$ and the relation $\delta(w) = 0$ (so that $\varphi_{B_0}(w) = w^{p}$). Then $C_0$ has the following universal property:
for every $\delta$-algebra $A$ over $B$, evaluation on $u$ induces a bijection
$$ \{ \textnormal{$B$-linear $\delta$-homomorphisms $C_0 \rightarrow A$} \} \rightarrow \{ \textnormal{Rank $1$ units of $A$} \}.$$
Note that $C_0$ is flat as an $B$-module, and the element $w-1$ is $(p,J)$-completely regular relative to $B$.
Let $K_0 \subseteq C_0$ be the ideal generated by $J$ and the element $w-1$, and let
$C = C_0\{ \frac{ K_0}{J} \}^{\wedge}$ be as in Proposition~3.13 of \cite{prisms}. Then the prism $(C,JC)$ has the desired properties.
\end{proof}

\begin{remark}
In the situation of Lemma \ref{lemma:add-rank-unit}, if the prism $(B,J)$ is transversal, then the prism $(C,K)$ is also transversal
(Remark \ref{remark:flatness-and-transversality}).
\end{remark}

\begin{lemma}\label{lemma:approximate-rank-one}
Let $(A,I)$ be a prism and let $u \in (1+I)_{\mathrm{rk}=1}$ be a rank $1$ unit. Then there exists
a morphism of prisms $f: (B,J) \rightarrow (A,I)$ and a rank $1$ unit $v \in (1+J)_{\mathrm{rk}=1}$
where $(B,J)$ is transversal and $f(v) = u$.
\end{lemma}

\begin{proof}
Combine Proposition \ref{proposition:transversal-approximation} with Lemma \ref{lemma:add-rank-unit}.
\end{proof}

\begin{proof}[Proof of Proposition \ref{proposition:prismatic-logarithm-extension}]
Let $(A,I)$ be a prism; we wish to show that there is a unique map
$\log_{\Prism}: (1+I)_{\mathrm{rk}=1} \rightarrow A\{1\}$ satisfying the requirements of Proposition \ref{proposition:prismatic-logarithm-extension}.
We first prove uniqueness. For every rank one unit $u \in (1+I)_{\mathrm{rk=1}}$, Lemma \ref{lemma:approximate-rank-one} guarantees
that we can choose a morphism of prisms $f: (B, J) \rightarrow (A,I)$, where $(B, J)$ is transversal, and an element $v \in (1+J)_{\mathrm{rk}=1}$
satisfying $f(v) = u$. It follows that $\log_{\Prism}(u)$ must be equal to the image of $v$ under the composite map
$$ (1+J)_{\mathrm{rk}=1} \xrightarrow{ \log_{\Prism} } B\{1\} \xrightarrow{ f\{1\} } A\{1\}.$$
To prove existence, we must show that this construction does not depend on the choice of approximation $f$. That is,
we must show that if $f': (B', J') \rightarrow (A,I)$ is another morphism of prisms, where $(B', J')$ is transversal, and
$v' \in (1+J)_{\mathrm{rk}=1}$ is an element satisfying $f'(v') = u$, then we have an equality $f\{1\}( \log_{\Prism}(v)) = f'\{1\}( \log_{\Prism}(v') )$
in $A\{1\}$. Using Proposition \ref{proposition:transversal-coproduct}, we can reduce to the case where $(B,J) = (B',J')$ and $f = f'$.
Choose a flat map of prisms $g: (B,J) \rightarrow (C,K)$ and an element $w \in (1+K)_{\mathrm{rk}=1}$ satisfying the conclusions
of Lemma \ref{lemma:add-rank-unit}. The universal property of $(C,K)$ then guarantees that there are unique
$B$-linear $\delta$-algebra homomorphisms $r,r': C \rightarrow B$ satisfying $r(w) = v$ and $r(w') = v'$.
The equality $f(v) = u = f(v')$ then guarantees that $f \circ r = f \circ r'$. It follows that
$$ f\{1\}( \log_{\Prism}(v) ) = (f \circ r)\{1\}( \log_{\Prism}(w) ) = (f \circ r')\{1\}( \log_{\Prism}(w) ) = f\{1\}( \log_{\Prism}(v') ),$$
as desired.
\end{proof}

\begin{proposition}
\label{prop:PrismLogHom}
Let $(A,I)$ be a prism. Then the prismatic logarithm
$$\log_{\Prism}: (1+I)_{\mathrm{rk}=1} \rightarrow A\{1\}$$ is a group homomorphism.
\end{proposition}

\begin{proof}
Fix a pair of elements $u,v \in (1+I)_{\mathrm{rk}=1}$; we wish to show that $\log_{\Prism}(uv) = \log_{\Prism}(u) + \log_{\Prism}(v)$.
By virtue of Proposition \ref{proposition:transversal-approximation}, there exists a morphism of prisms $f: (B,J) \rightarrow (A,I)$
where $(B,J)$ is transversal. By virtue of Lemma \ref{lemma:add-rank-unit}, we can further arrange that there exist elements $u', v' \in (1+J)_{\mathrm{rk}=1}$
satisfying $f(u') = u$ and $f(v') = v$. In this case, the desired result follows from the identity $\log_{\Prism}(u'v') = \log_{\Prism}(u') + \log_{\Prism}(v')$ in $B\{1\}$ (Proposition \ref{proposition:transversal-log-group-homomorphism}).
\end{proof}

\begin{remark}[Functoriality]\label{remark:functoriality-of-log-2}
Let $f: (B,J) \rightarrow (A,I)$ be a morphism of prisms. Then the diagram of abelian groups
$$ \xymatrix@R=50pt@C=50pt{ (1+J)_{\mathrm{rk}=1} \ar[r]^-{f} \ar[d]^{\log_{\Prism}} & (1+I)_{\mathrm{rk}=1} \ar[d]^{ \log_{\Prism} } \\
B\{1\} \ar[r]^-{ f\{1\} } & A\{1\} }$$
is commutative. To prove this, we can use Lemma \ref{lemma:approximate-rank-one} to reduce to the case where $(B,J)$ is transversal, in which
case it follows immediately from the definition.
\end{remark}

\begin{proposition}\label{proposition:log-Nygaard}
Let $(A,I)$ be a prism and let $u \in (1+I)_{\mathrm{rk}=1}$. Then the prismatic logarithm $\log_{\Prism}(u)$
satisfies the identity $\varphi_{A\{1\}}( \log_{\Prism}(u) ) = \log_{\Prism}(u)$. 
\end{proposition}

\begin{proof}
Using Lemma \ref{lemma:approximate-rank-one}, we can reduce to the case where $(A,I)$ is transversal, in which case it
follows from the identity $\varphi_{A}( u^{p^{r}} - 1) = u^{p^{r+1}} - 1$
for $r \geq 0$.
\end{proof}

\begin{proposition}\label{proposition:normalization-of-logarithm}
Let $(A,I)$ be a prism. Then the diagram of abelian group homomorphisms
$$ \xymatrix@R=50pt@C=50pt{ (1+I)_{\mathrm{rk}=1} \ar[r]^-{ \log_{\Prism} } \ar[d]^{u \mapsto u-1} & A\{1\} \ar[d] \\
I/I^2 \ar[r]^-{ \beta^{-1} }_{\sim} & A\{1\} / I A\{1\} }$$
is commutative, where $\beta$ is the isomorphism of Remark \ref{remark:twist-modulo-I}.
\end{proposition}

\begin{proof}
Using Lemma \ref{lemma:approximate-rank-one}, we can reduce to the case where $(A,I)$ is transversal.
In this case, the desired result follows immediately from the definition of $\log_{\Prism}$ (see Construction \ref{construction:prismatic-logarithm}).
\end{proof}

\subsection{Example: The \texorpdfstring{$q$}{q}-de Rham Prism}\label{subsection:q-dR}

Let $(A,I) = ( \Z_p[[q-1]], ( [p]_q)  )$ denote the $q$-de Rham prism of Example \ref{example:q-prism}. Note that the identity $q^{p} - 1 = (q-1)(1 + \cdots + q^{p-1})$
guarantees that $q^{p} \equiv 1 \pmod{I}$, so that $q^{p}$ is a rank $1$ element of the group $(1+I)^{\times}$. Applying the prismatic logarithm
$\log_{\Prism}$ of Construction \ref{construction:prismatic-logarithm}, we obtain an element $\log_{\Prism}(q^{p}) \in A\{1\}$.

\begin{proposition}\label{proposition:twist-in-q-de-Rham-case}
Let $(A,I) = ( \Z_p[[q-1]], ([p]_{q}) )$ be the $q$-de Rham prism. Then there is a unique element $e_{A} \in A\{1\}$ satisfying
$(q-1) e_{A} = \log_{\Prism}(q^{p})$. Moreover, $A\{1\}$ is the free $A$-module generated by $e_{A}$.
\end{proposition}

\begin{proof}
For each $r \geq 0$, the ideal $I_{r} \subseteq A$ of Notation \ref{notation:I-r} is generated by
the element $[ p^{r} ]_{q} = \frac{ q^{p^{r}} - 1}{q-1} = 1 + q + \cdots + q^{ p^{r} - 1}$.
In particular, we have $q^{p^{r}} \equiv 1 \pmod{I_r}$, which yields a congruence
$$ [ p^{r+1} ]_{q} = [p]_{q^{p^{r}}} [p^{r}]_{q}  \equiv p [p^{r}]_{q} \pmod{ I_{r}^{2} }.$$
It follows that there is a unique element $e_{A} \in A\{1\}$ whose image in each quotient $I_{r} / I_{r}^{2}$
is given by the residue class of $[ p^{r} ]_{q}$, and that $e_{A}$ is a generator of $A\{1\}$.
The identities $(q-1) [p^{r}]_{q} = q^{p^{r}}-1$ imply that $(q-1) e_{A} = \log_{\Prism}(q^{p} )$.
Since $q-1$ is not a zero divisor in $A$, the element $e_{A}$ is uniquely determined by this property.
\end{proof}

\begin{remark}\label{remark:preferred-generator-universal-case}
Stated more informally, Proposition \ref{proposition:twist-in-q-de-Rham-case} asserts that
the invertible $A$-module $A\{1\}$ is freely generated by the element $e_{A} = \frac{ \log_{\Prism}( q^{p} ) }{q-1}$.
Note that the Frobenius morphism $\varphi_{A\{1\}}: A\{1\} \rightarrow I^{-1}\{1\}$ of 
Remark \ref{remark:Frobenius-on-twist} carries the element $\log_{\Prism}(q^{p})$ to itself (Proposition \ref{proposition:log-Nygaard}),
and therefore satisfies the identity
$$ \varphi_{A\{1\} }( e_A ) = \varphi_{ A\{1\} }( \frac{ \log_{\Prism}(q^{p})}{q-1} ) =
\frac{ \log_{\Prism}(q^{p})}{q^{p}-1} = [p]_{q}^{-1} e_{A}.$$
\end{remark}

\begin{notation}\label{notation:preferred-generator}
Let $f: ( \Z_p[[q-1]], ( [p]_{q} ) ) \rightarrow ( A, I)$ be a morphism of prisms. We will write
$e_{A}$ for the image of the generator $\frac{ \log_{\Prism}(q^{p})}{q-1}$
under the homomorphism $f\{1\}: \Z_p[[q-1]]\{1\} \rightarrow A\{1\}$. It follows from
Remark \ref{remark:preferred-generator-universal-case} that the $A$-module $A\{1\}$
is freely generated by the element $e_{A}$, and that the Frobenius map $\varphi_{A\{1\}}: A\{1\} \rightarrow I^{-1}\{1\}$ is given by concretely
by the formula
$$ \varphi_{A\{1\}}( \lambda e_{A} ) = \frac{ \varphi_{A}(\lambda)}{ f( [p]_q) } e_{A}.$$
Beware that the element $e_{A} \in A\{1\}$ depends not only on the prism $(A,I)$, but also on the choice of morphism $f$.
\end{notation}

\begin{example}[The Breuil-Kisin twist over crystalline prisms]\label{BKcrystalline}
Let $(A,I)$ be a crystalline prism (so that $I = (p)$ is the principal ideal generated by $p$). Then
there is a canonical map of prisms
$$ f: ( \Z_p[[q-1]], ( [p]_{q} ) ) \rightarrow (A,I) \quad \quad q \mapsto 1$$
satisfying $f( [p]_q) = p$. It follows that Notation \ref{notation:preferred-generator} supplies a canonical generator $e_{A} \in A\{1\}$
which satisfies the identity $\varphi_{A\{1\}}( e_A ) = \frac{ e_A }{p}$.
\end{example}

We now compare the prismatic logarithm of Definition \ref{definition:prismatic-logarithm-general} with the $q$-logarithm of \cite{cyctrace}. 

\begin{proposition}\label{proposition:universal-divided-power-element}
There exists a morphism of prisms $( \Z_p[[q-1]], ([p]_q) ) \rightarrow (A, I)$ and a rank
$1$ unit $u \in A$ with the following universal property:
\begin{itemize}
\item[$(1)$] The element $\varphi_{A}(u)$ is congruent to $1$ modulo $I$.

\item[$(2)$] For every morphism of prisms $(\Z_p[[q-1]], ([p]_q) ) \rightarrow (B,J)$
and every rank $1$ unit $v \in B$ satisfying $\varphi_{B}(v) \equiv 1 \pmod{J}$, there is a unique $\Z_p[[q-1]]$-linear $\delta$-algebra homomorphism
$f: A \rightarrow B$ satisfying $f(u) = v$.
\end{itemize}
Moreover, the commutative ring $A$ is flat over $\Z_p[[q-1]]$.
\end{proposition}

\begin{proof}
Let $A_0$ denote the Laurent polynomial ring $\Z_p[[q-1]][ u^{\pm 1} ]$, which we regard as a $\delta$-algebra over $\Z_p[[q-1]]$ by setting
$\delta(u) = 0$ (so that $u$ is a rank $1$ unit of $A_0$). Note that the quotient $A_0 / (u^{p} - 1)$ is free of rank $p$ as a $\Z_p[[q-1]]$-module, so that the element $u^{p}-1 \in A_0$ is
$( p, q-1 )$-completely regular relative to $\Z_p[[q-1]]$. Let $I_0 \subseteq A_0$ be the ideal generated by $[p]_q$ and $u^{p} -1$ and let
$A$ denote the $\delta$-ring $A_0\{ \frac{I_0}{[p]_q}  \}^{\wedge} = A_0\{ \frac{ \varphi(u-1)}{ [p]_q} \}^{\wedge}$ given by Proposition~3.13 of \cite{prisms}. Setting
$I = [p]_q A$, we have a morphism of prisms $(\Z_p[[q-1]], ( [p]_q) ) \rightarrow ( A, I )$ which satisfies conditions $(1)$ and $(2)$ by construction. Moreover, Proposition~3.13 of \cite{prisms}
guarantees that $A$ is $( p, [p]_q)$-completely flat over $\Z_p[[q-1]]$, and therefore flat over $\Z_p[[q-1]]$ (since $A$ is complete with respect to the ideal
$(p, [p]_{q} )$ and the commutative ring $\Z_p[[q-1]]$ is Noetherian).
\end{proof}

\begin{remark}\label{remark:divided-power-basis}
Let $A$ be as in Proposition \ref{proposition:universal-divided-power-element}. Then the element $q-1 \in A$ is a non-zero-divisor, and the quotient ring
$A / (q-1) A \simeq \Z_p\{ \frac{ \varphi_{A}( u-1)}{p} \}^{\wedge}$ can be identified with the $p$-completion of the divided power algebra
$\bigoplus_{n \geq 0} \frac{ (u -1)^{n} }{n!} \Z_p$. See Corollary~2.38 of \cite{prisms}.
\end{remark}

For each nonnegative integer $n \geq 0$, we let $[n]_{q}!$ denote the product $\prod_{m=1}^{n} [m]_q$.

\begin{proposition}\label{proposition:q-divided-power-structure}
Let $A$ be the $\delta$-ring of Proposition \ref{proposition:universal-divided-power-element}. For each integer $n \geq 0$, there is a unique element
$\gamma_{n,q}(u - 1)  \in A$ satisfying the identity $$[n]_q! \gamma_{n,q}(u -1) = (u-1) (u - q) \cdots (u - q^{n-1} ).$$
Moreover, as a module over $\Z_p[[q-1]]$, $A$ can be identified with the $(p,q-1)$-completion of the free module
$$ \bigoplus_{n \geq 0} \gamma_{n,q}(u-1) \Z_p[[q-1]].$$
\end{proposition}

\begin{proof}
Since $A$ is flat over $\Z_p[[q-1]]$, the elements $\gamma_{n,q}(u-1)$ are automatically unique if they exist. For existence, we refer to
Proposition~4.9 of \cite{cyctrace}. Note that reduction modulo $q-1$ carries each $\gamma_{n,q}( u-1)$ to the divided power
$\frac{ (u-1)^n}{n!} \in A/(q-1)$. Using Remark \ref{remark:divided-power-basis}, we conclude that the natural map
$$ (\bigoplus_{n \geq 0} \gamma_{n,q}(u-1) \Z_p[[q-1]])^{\wedge}_{(p,q-1)} \rightarrow A$$
induces an isomorphism after reduction modulo $(q-1)$, and is therefore an isomorphism (since both sides are $(q-1)$-complete and
$(q-1)$-torsion-free).
\end{proof}

\begin{definition}[\cite{cyctrace}]
Let $A$ be the $\delta$-ring of Proposition \ref{proposition:universal-divided-power-element}. We let $\log_{q}(u)$ denote the element of $A$ given by
the infinite sum
$$ \sum_{n=1}^{\infty} (-1)^{n-1} q^{ -n(n-1)/2} \frac{ (u-1) (u-q) \cdots (u- q^{n-1} ) }{ 1 + q + \cdots + q^{n-1} }.$$
Note that this sum is convergent in $A$, since its $n$th term is divisible by $[n-1]_{q}!$.
\end{definition}

\begin{proposition}[\cite{cyctrace}]\label{proposition:characterize-q-log}
The element $x = \log_{q}(u) \in A$ is characterized by the following properties:
\begin{itemize}
\item[$(1)$] Let $\rho: A \rightarrow \Z_p[[q-1]]$ denote the unique $\Z_p[[q-1]]$-linear $\delta$-homomorphism satisfying $\rho(u) = 1$. Then
$\rho(x) = 0$.
\item[$(2)$] Let $\tau: A \rightarrow A$ denote the unique $\Z_p[[q-1]]$-linear $\delta$-homomorphism satisfying $\tau(u) = qu$. Then
$\tau(x) - x = q-1$.
\end{itemize}
\end{proposition}

\begin{proof}
See Lemma~4.6 of \cite{cyctrace}; note that the identities $(1)$ and $(2)$ can be written more informally as $\log_{q}(1) = 0$ and $\frac{ \log_{q}(qu) - \log_{q}(u)}{ qu-u} = \frac{1}{u}$,
respectively. 
\end{proof}

\begin{proposition}\label{proposition:qlog-agreement}
Let $(A, ( [p]_q ) )$ be the prism of Proposition \ref{proposition:universal-divided-power-element}, and let $e_{A} \in A\{1\}$ be the generator of
Notation \ref{notation:preferred-generator}. Then the prismatic logarithm $\log_{\Prism}$ of Construction \ref{construction:prismatic-logarithm} satisfies the identity
$$ \log_{\Prism}( u^{p} ) = \log_{q}(u) e_{A}.$$
\end{proposition}

\begin{proof}
Since $A\{1\}$ is a free module of rank $1$ generated by $e_{A}$, there is a unique element $f \in A$ satisfying $\log_{\Prism}( u^p ) = f e_{A}$.
To show that $f = \log_{q}(u)$, it will suffice to verify that $f$ satisfies conditions $(1)$ and $(2)$ of Proposition \ref{proposition:characterize-q-log}. Condition $(1)$ follows from the identity $\log_{\Prism}(1) = 0$, and $(2)$ from the identity $\log_{\Prism}( q^{p} u^{p} ) - \log_{\Prism}(u^{p} ) = \log_{\Prism}(q^{p}) = (q-1) e_{A}$.
\end{proof}

\begin{corollary}\label{corollary:compare-log-bounded}
Let $f: (\Z_p[[q-1]], ( [p]_q) ) \rightarrow (A,I)$ be any morphism of prisms, and let $e_{A}$ be the generator of $A\{1\}$ described in
Notation \ref{notation:preferred-generator}. Assume that $(A,I)$ is bounded and that, for every positive integer $n$, the element
$[n]_{q}$ is a non-zero-divisor in $(A,I)$. Let $u \in A$ be a rank $1$ unit satisfying $\varphi_{A}(u-1) = u^{p} -1  \in J$. Then:
\begin{itemize}
\item[$(1)$] For every nonnegative integer $n$, the product $(u-1)(u-q)\cdots( u-q^{n-1})$ is uniquely divisible by $[n]_{q}!$.
\item[$(2)$] The sum 
$$ \log_{q}(u) = \sum_{n > 0} (-1)^{n-1} q^{ -n(n-1)/2} \frac{ (u-1) (u-q) \cdots (u-q^{n-1} ) }{ [n]_q }$$
converges with respect to the $(p, q-1)$-adic topology on $A$.

\item[$(3)$] We have $\log_{\Prism}(u^{p} ) = \log_{q}(u) e_{A}$ in $A\{1\}$.
\end{itemize}
\end{corollary}

\begin{proof}
The assertion reduces immediately to the universal case where $(A,I)$ is the prism of Proposition \ref{proposition:universal-divided-power-element}
(note that the technical assumption that $(A,I)$ is bounded guarantees that $A$ is separated and complete with respect to $(p,q-1)$-adic topology, so that
the meaning of convergence is unambiguous).
\end{proof}

\begin{corollary}\label{corollary:crystalline-formula-for-log}
Let $A$ be a $\delta$-ring which is $p$-complete and $p$-torsion-free, so that we can regard $(A, (p))$ as a crystalline prism,
and let $e_{A} \in A\{1\}$ be the generator of Remark \ref{BKcrystalline}. Let $u \in A$ be an element satisfying $\varphi_{A}(u) = u^{p}$
and $u^{p} \equiv 1 \pmod{p}$. Then:
\begin{itemize}
\item[$(1)$] The element $u-1 \in A$ has divided powers: that is, for every positive integer $n$, the element
$(u-1)^{n}$ is (uniquely) divisible by $n!$.

\item[$(2)$] The sum 
$$ \log(u) = \sum_{n > 0} (-1)^{n-1}  \frac{ (u-1)^{n} }{ n }$$
converges with respect to the $p$-adic topology

\item[$(3)$] We have $\log_{\Prism}(u^{p} ) = \log(u) e_{A}$ in $A\{1\}$.
\end{itemize}
\end{corollary}

\begin{proof}
Apply Corollary \ref{corollary:compare-log-bounded} to the morphism of prisms
$$ ( \Z_p[[q-1]], ([p]_q) ) \rightarrow (A, (p)) \quad \quad q \mapsto 1.$$
\end{proof}

\subsection{Tate Modules}\label{subsection:Tate-modules}

The prismatic logarithm of Definition \ref{definition:prismatic-logarithm-general} associates to every prism $(A,I)$ a group homomorphism
$$ \log_{\Prism}: ( 1 + I )_{\mathrm{rk}=1} \rightarrow A\{1\}.$$
In practice, the abelian group $(1 + I )_{\mathrm{rk} = 1}$ can be cumbersome to work with, since it depends on the $\delta$-structure on $A$.
In this section, we introduce a variant of the prismatic logarithm (Construction \ref{construction:Tate-logarithm}) which does not not share this defect (and contains essentially the same information:
see Remark \ref{remark:no-info-lost}).

\begin{notation}[The Tate Module]\label{notation:R-flat}
Let $R$ be a commutative ring. We let $R^{\flat}$ denote the inverse limit of the diagram
$$ \cdots \xrightarrow{x \mapsto x^{p}} R \xrightarrow{x \mapsto x^{p}} R \xrightarrow{x \mapsto x^{p}} R,$$
which we regard as a commutative monoid under multiplication. By definition, an element of $R^{\flat}$ is given by a sequence
$\{ x_n \}_{n \geq 0}$ of elements of $R$ which satisfy satisfying $x_{n+1}^{p} = x_n$ for each $n \geq 0$. 
We let $R^{\flat \times}$ denote the group of invertible elements of $R^{\flat}$: that is, the inverse limit of the diagram
of unit groups
$$ \cdots \xrightarrow{x \mapsto x^{p}} R^{\times} \xrightarrow{x \mapsto x^{p}} R^{\times} \xrightarrow{x \mapsto x^{p}} R^{\times}.$$
The construction $\{ x_n \}_{n \geq 0} \mapsto x_0$ determines a group homomorphism $R^{\flat \times} \rightarrow R^{\times}$,
which we denote by $x \mapsto x^{\sharp}$. We denote the kernel of this homomorphism by $T_p(R^{\times} )$, and refer to it
as the {\it Tate module of $R$}. By definition, $T_p(R^{\times} )$ is given by the inverse limit of the diagram of abelian groups
$$ \cdots \xrightarrow{x \mapsto x^{p}} \mu_{p^2}(R) \xrightarrow{x \mapsto x^{p}} \mu_{p}(R) \xrightarrow{x \mapsto x^{p}} \{1\},$$
where $\mu_{p^{n}}(R) = \{ u \in R: u^{p^{n}} = 1 \}$ denotes the group of $p^n$th roots of unity in $R$.
By construction, we have an exact sequence of abelian groups
$$ 0 \rightarrow T_p( R^{\times} ) \rightarrow R^{\flat \times} \xrightarrow{ x \mapsto x^{\sharp} } R^{\times},$$
which is exact on the right if and only if every invertible element of $R$ admits a $p$th root.
\end{notation}

\begin{proposition}\label{proposition:pth-roots-rank-1}
Let $A$ be a $p$-complete $\delta$-ring. Then, for every element $x \in A^{\flat}$, the image $x^{\sharp} \in A$ is
rank one: that is, it satisfies $\delta( x^{\sharp} ) = 0$.
\end{proposition} 

\begin{proof}
If $A$ is $p$-adically separated, then the result follows from Lemma~2.32 of \cite{prisms}. To prove it in general, let us identify
the element $x \in A^{\flat}$ with a homomorphism of commutative rings $f: \Z[ X^{1/p^{\infty}} ] \rightarrow A$. Then $f$ factors
canonically as a composition $$\Z[ X^{1/p^{\infty}} ] \rightarrow B \xrightarrow{f'} A,$$
where $B$ is the free $\delta$-ring generated by $\Z[ X^{1/p^{\infty} } ]$. We may therefore
assume without loss of generality that $A$ is the $p$-completion of $B$. 
Writing $\Z [ X^{1/p^{\infty}} ]$ as a filtered colimit of subrings of the form $\Z[ X^{1/p^{n}} ]$, we deduce that
$B$ can be realized as a filtered colimit of free $\delta$-rings, and is therefore torsion-free as an abelian group.
In particular, the $p$-completion $\widehat{B}$ is $p$-adically separated, so the desired result follows from \cite{prisms}.
\end{proof}

\begin{proposition}\label{proposition:lift-flat}
Let $(A,I)$ be a prism. Then the quotient map $A \twoheadrightarrow A/I$ induces an isomorphism of abelian groups
$A^{\flat \times} \xrightarrow{\sim} (A/I)^{\flat \times}$.
\end{proposition}

\begin{proof}
The abelian group $A^{\flat \times}$ can be identified with the inverse limit of the system
$$ \cdots \rightarrow (A/I^{3})^{\flat \times} \rightarrow (A/I^2)^{\flat \times} \rightarrow (A/I)^{\flat \times}.$$
It will therefore suffice to show that each of the transition maps $( A/I^{n+1})^{\flat \times} \rightarrow (A/I^{n} )^{\flat \times}$ is an isomorphism.
Let $\overrightarrow{M}$ denote the inverse system of abelian groups
$$ \cdots \rightarrow I^{n} / I^{n+1} \xrightarrow{ p} I^{n} / I^{n+1} \xrightarrow{ p } I^{n} / I^{n+1}$$
where each term is equal to the quotient $I^{n} / I^{n+1}$ and each transition map is given by multiplication by $p$.
The short exact sequence 
$$0 \rightarrow I^{n} / I^{n+1} \xrightarrow{ x \mapsto 1 + x } (A / I^{n+1} )^{\times} \rightarrow (A/I^n)^{\flat \times} \rightarrow 0$$
yields a long exact sequence
$$ 0 \rightarrow {\lim}^{0} \overrightarrow{M} \rightarrow (A / I^{n+1} )^{\times \flat} \rightarrow (A /I^{n} )^{\flat \times} \rightarrow {\lim}^{1} \overrightarrow{M}$$
We conclude by observing that the groups $\lim^{0} \overrightarrow{M}$ and $\lim^{1} \overrightarrow{M}$ vanish, since the abelian group $I^{n} / I^{n+1}$ is $p$-complete.
\end{proof}

\begin{construction}\label{construction:Tate-logarithm}
Let $(A,I)$ be a prism and let $\overline{A}$ denote the quotient ring $A/I$. Proposition \ref{proposition:lift-flat} guarantees that every element $x \in \overline{A}^{\flat \times}$ can be lifted uniquely
to an element $\widetilde{x} \in A^{\flat \times}$, and Proposition \ref{proposition:pth-roots-rank-1} guarantees that $\widetilde{x}^{\sharp}$ is
a rank one element of $A^{\times}$. The construction $x \mapsto \widetilde{x}^{\sharp}$ determines a group homomorphism
$\rho: (A/I)^{\flat \times} \rightarrow A^{\times}_{\mathrm{rk}=1}$ which fits into a commutative diagram
$$ \xymatrix@R=50pt@C=50pt{ 0 \ar[r] & T_p( \overline{A}^{\times} ) \ar[r] \ar[d]^{\rho_0} & \overline{A}^{\flat \times} \ar[r]^-{\sharp} \ar[d]^{\rho} & \overline{A}^{\times} \ar@{=}[d] \\
0 \ar[r] & (1+I)_{\mathrm{rk}=1} \ar[r] & A^{\times}_{\mathrm{rk}=1} \ar[r] & \overline{A}^{\times}. }$$
Let $\log_{\Prism}: (1+I)_{\mathrm{rk}=1} \rightarrow A\{1\}$ be the prismatic logarithm of Definition \ref{definition:prismatic-logarithm-general}.
We will abuse terminology by referring to the composite map
$$ T_p( \overline{A}^{\times} ) \xrightarrow{\rho_0} (1+I)_{\mathrm{rk}=1} \xrightarrow{ \log_{\Prism} } A\{1\}$$
as the {\it prismatic logarithm}, and denoting it also by $\log_{\Prism}: T_p( \overline{A}^{\times} ) \rightarrow A\{1\}$.
\end{construction}

\begin{example}\label{example:perfected-q-prism}
Let $A$ denote the completion of the ring $\Z[ q^{1/p^{\infty}} ]$ with respect to the ideal $(p,q-1)$, equipped with the $\delta$-structure given by
$\varphi_{A}( q^{ 1/p^n} ) = q^{1 / p^{n-1} }$, and let $I \subseteq A$ be the ideal generated by $[p]_q = \frac{ q^{p}- 1}{q-1}$.
Then $\epsilon = ( 1, q, q^{1/p}, q^{1/p^2}, \cdots )$ is a compatible sequence of $p$th power roots of unity in the quotient ring $\overline{A} = A/I$,
which we can regard as an element of the Tate module $T_p( \overline{A}^{\times} )$. By construction, the
homomorphism $\rho_0: T_p( \overline{A}^{\times} ) \rightarrow (1+I)_{\mathrm{rk}=1}$
of Construction \ref{construction:Tate-logarithm} satisfies $\rho_0(\epsilon) = q^{p}$. We therefore have $\log_{\Prism}(\epsilon) = (q-1) e_{A}$, where $e_{A}$ is the generator of $A\{1\}$ described in Proposition \ref{proposition:twist-in-q-de-Rham-case}.
\end{example}

\begin{remark}\label{remark:no-info-lost}
Replacing the prismatic logarithm $\log_{\Prism}: (1+I)_{\mathrm{rk}=1} \rightarrow A\{1\}$ of Definition \ref{definition:prismatic-logarithm-general} by
the prismatic logarithm $\log_{\Prism}: T_p( \overline{A}^{\times} ) \rightarrow A\{1\}$ of Construction \ref{construction:Tate-logarithm} does not lose any information,
provided that the prism $(A,I)$ is allowed to vary.
If $(A,I)$ is any bounded prism and $u$ is a rank one element of $1+I$, then we can form a new bounded prism $(A',IA')$, where
$A'$ is the $(p,I)$-completion of the tensor product $A \otimes_{\Z[u] } \Z[ u^{1/p^{\infty}}]$ (endowed with the $\delta$-structure where
each $u^{1/p^{n}}$ is annihilated by $\delta$). Setting $\overline{A}' = A' / IA'$ and letting $\overline{u}^{1/p^n}$ denote the image of $u^{1/p^n}$ in
$\overline{A}'$, we note that the sequence $\{ \overline{u}^{1/p^{n}} \}_{n \geq 0}$ determines an element $e \in T_p( \overline{A}'^{\times} )$,
and the tautological map $A\{1\} \rightarrow A'\{1\}$ is an injection which carries $\log_{\Prism}( u )$ (in the sense of Definition \ref{definition:prismatic-logarithm-general})
to $\log_{\Prism}(e)$ (in the sense of Construction \ref{construction:Tate-logarithm}).
\end{remark}

\begin{remark}[Comparison with Prismatic Dieudonn\'{e} Theory]\label{remark:Dieudonne-comparison}
Let $(A,I)$ be a bounded prism and let $\calC$ denote the category of bounded prisms $(B,J)$ equipped with a morphism $(A,I) \rightarrow (B,J)$,
The construction $(B,J) \mapsto T_p( (B/J)^{\times} )$ determines a functor from $\calC$ to the category of abelian groups, which we will denote by $T$.
In \cite{prismaticdieudonne}, Ansch\"{u}tz and Le Bras define the Breuil-Kisin twist $A\{-1\}$ to be the set of natural transformations from $T$ to the forgetful functor
$$ U: \calC \rightarrow \{ \text{Abelian groups} \} \quad \quad (B,J) \mapsto B.$$
To avoid confusion, let us temporarily denote this set of natural transformations by $A\{-1\}'$; it is an invertible $A$-module which depends functorially on the bounded
prism $(A,I)$ (see Definition~4.7.2 of \cite{prismaticdieudonne}). 

Let $A\{-1\}$ denote the inverse of the invertible $A$-module $A\{1\}$ introduced in Definition \ref{definition:twist-general}. Every element $x \in A\{-1\}$
determines a natural transformation of functors $\lambda(x): T \rightarrow U$, given by the collection of group homomorphisms
$$ \lambda(x)_{(B,J)}: T_p( (B/J)^{\times} ) \rightarrow B \quad \quad y \mapsto x \cdot \log_{\Prism}(y).$$
The construction $x \mapsto \lambda(x)$ determines an $A$-module homomorphism $\lambda: A\{-1\} \rightarrow A\{-1\}'$, depending functorially on $A$.
We claim that $\lambda$ is an isomorphism: in other words, Definition \ref{definition:twist-general} agrees with the Breuil-Kisin twist defined in \cite{prismaticdieudonne}.
To prove this, we may assume without loss of generality that $(A,I)$ is the $q$-de Rham prism $( \Z_p[[q-1]], ( [p]_q ) )$. In this case, the invertible $A$-module $A\{1\}$ is generated by the element $e_{A} = \frac{ \log_{\Prism}(q^p) }{ q-1 }$ (Proposition \ref{proposition:twist-in-q-de-Rham-case}). It will therefore suffice to show that
$\lambda$ carries $e_{A}^{-1}$ to a generator of the invertible $A$-module $A\{-1\}'$, which follows by combining 
Proposition \ref{proposition:qlog-agreement} with Proposition~4.7.3 of \cite{prismaticdieudonne}.
\end{remark}

\newpage \section{The Cartier-Witt Stack}\label{section:cartier-witt}
\label{sec:WCart}

Let $R$ be a commutative ring in which $p$ is nilpotent, so that the ring of Witt vectors $W(R)$ has the structure of a $p$-complete $\delta$-ring.
In \S\ref{subsection:prismatic-stack}, we introduce the notion of a {\it Cartier-Witt divisor for $R$} (Definition \ref{definition:Cartier-Witt-divisor}). The collection of Cartier-Witt divisors for $R$ comprise a groupoid $\WCart(R)$, which can be thought of as an enlargement of the set of prism structures on $W(R)$
(see Remark \ref{remark:Cartier-Witt-divisor-prism} for a precise statement). The construction $R \mapsto \WCart(R)$ determines a groupoid-valued functor
on the category of commutative rings, which we will denote by $\WCart$ and refer to as the {\it Cartier-Witt stack}.

Heuristically, one can think of the Cartier-Witt stack $\WCart$ as playing the role of the formal spectrum $\Spf( A_{\mathrm{init}} )$, where $(A_{\mathrm{init}}, I_{\mathrm{init}} )$ is an
initial object of the category of prisms. Beware that this is not literally true: the category of prisms does not have an initial object, and the Cartier-Witt stack $\WCart$
is not a formal scheme. Nevertheless, in \S\ref{subsection:cartier-witt-atlas} we support this heuristic by associating to every prism $(A,I)$ a morphism
$\rho_{A}: \Spf(A) \rightarrow \WCart$, which depends functorially on $(A,I)$. Under mild hypotheses, we show that if $(B,J)$ is another prism,
then the fiber product $\Spf(A) \times_{\WCart} \Spf(B)$ can be identified with the affine formal scheme $\Spf(C)$, where $(C,K)$ is a coproduct of $(A,I)$ and $(B,J)$ in the category of prisms (Proposition \ref{proposition:pullback-diagram-prisms}).

Let $\calD( \WCart )$ denote the derived $\infty$-category of quasi-coherent complexes on $\WCart$ (see Definition \ref{definition:DWCart}).
For every prism $(A,I)$, the morphism $\rho_{A}: \Spf(A) \rightarrow \WCart$ determines a pullback functor $$\rho_{A}^{\ast}: \calD( \WCart ) \rightarrow \calD( \Spf(A) ) \simeq \widehat{\calD}(A),$$
where $\widehat{\calD}(A)$ denotes the full subcategory of $\calD(A)$ spanned by those complexes which are $(p,I)$-complete.
In \S\ref{subsection:complexes-on-WCart}, we show that these maps supply an equivalence of $\infty$-categories
$\calD( \WCart ) \rightarrow \varprojlim_{ (A,I) } \widehat{\calD}(A)$, where limit is indexed by the category of bounded prisms (Proposition \ref{proposition:DWCart-prism-description}). In other words, to supply a quasi-coherent complex on $\WCart$, one must give a rule which associates to each bounded prism $(A,I)$
a $(p,I)$-complete complex of $A$-modules, which is compatible with $(p,I)$-completed extension of scalars. For this reason, the Cartier-Witt stack
plays an essential role in organizing the information supplied by the (relative) prismatic cohomology of \cite{prisms}.

The structure of the Cartier-Witt stack is somewhat subtle. However, $\WCart$ contains a divisor $\WCart^{\mathrm{HT}}$, which we will refer to as the {\it Hodge-Tate divisor}, which is much easier to analyze. In \S\ref{subsection:HT-divisor}, we show that it can be identified with ($p$-complete) classifying stack of the affine group scheme $\mathbf{G}_{m}^{\sharp}$ obtained as the divided power envelope of the multiplicative group $\mathbf{G}_{m}$ along its identity section (Theorem \ref{theorem:describe-HT}). In \S\ref{subsection:HT-locus}, we apply this result to classify quasi-coherent complexes $\mathscr{E}$ on $\WCart^{\mathrm{HT}}$: they
can be identified with pairs $( \mathscr{E}_{\eta}, \Theta_{\mathscr{E} } )$, where $\mathscr{E}_{\eta}$ is a $p$-complete complex of abelian groups and $\Theta_{\mathscr{E}}$ is an endomorphism of $\mathscr{E}_{\eta}$ satisfying a suitable integrality condition (Theorem \ref{theorem:compute-with-HT}). 
In \S\ref{subsection:Sen-theory}, we show that $\Theta_{\mathscr{E}}$ can be regarded as an analogue of the classical {\it Sen operator} appearing in the theory of semilinear Galois representations (Theorem \ref{theorem:Sen-theory-refined}).

In \S\ref{subsection:computing-with-WCart}, we extend the preceding ideas to understand quasi-coherent complexes on the {\em entire} Cartier-Witt stack $\WCart$.
Applying the construction of \S\ref{subsection:cartier-witt-atlas} to the $q$-de Rham prism $( \Z_p[[q-1]], [ p]_{q} )$ of Example \ref{example:q-prism}, we obtain
a comparison map $\rho_{\qdR}: \Spf( \Z_p[[q-1]] ) \rightarrow \WCart$. This map is invariant under the action of the profinite group
$\Z_p^{\times}$ (acting on the the $q$-de Rham prism $\Z_p[[q-1]]$ via the construction $q \mapsto q^{u}$ for $u \in \Z_p^{\times}$), and therefore
descends to a map of stacks $[\Spf( \Z_p[[q-1]] ) / \Z_p^{\times}] \rightarrow \WCart$ (Proposition \ref{proposition:continuous-action}).
This morphism fits into a commutative diagram
$$ \xymatrix@R=50pt@C=50pt{ [\Spf(\Z_p) / \Z_p^{\times}] \ar[r]^-{q \mapsto 1} \ar[d] & [\Spf( \Z_p[[q-1]] ) / \Z_p^{\times}] \ar[d]^{\rho_{\qdR}} \\
\Spf(\Z_p) \ar[r]^-{ \rho_{\dR} } & \WCart. }$$
In \S\ref{subsection:computing-with-WCart}, we show that when $p$ is an odd prime, then this diagram is close to being a pushout square in the following sense: for
every quasi-coherent complex $\mathscr{E}$ on $\WCart$, the induced diagram of complexes
$$ \xymatrix@R=50pt@C=50pt{ \RGamma( \Z_p^{\times}, \rho_{\dR}^{\ast} \mathscr{E} ) & \RGamma( \Z_p^{\times}, \rho_{\qdR}^{\ast} \mathscr{E} ) \ar[l] \\
\rho_{\dR}^{\ast}( \mathscr{E} ) \ar[u] & \RGamma( \WCart, \mathscr{E} ) \ar[l] \ar[u]  }$$
is a pullback square in $\widehat{\calD}(\Z_p)$ (see Theorem \ref{theorem:qdR-WCart}). We will deduce this from an analogous statement for the
Hodge-Tate divisor $\WCart^{\mathrm{HT}}$, which we carry out in \S\ref{subsection:exp-theta} (see Proposition \ref{proposition:clean-statement-of-fiber}).

The Cartier-Witt stack $\WCart$ is equipped with an endomorphism $F: \WCart \rightarrow \WCart$, which we will refer to as the {\it Frobenius morphism} (Construction \ref{construction:Frobenius-on-stack}). A crucial property of the Frobenius morphism is that it {\em contracts} the Hodge-Tate divisor $\WCart^{\mathrm{HT}}$, in the sense that there is a commutative diagram $$ \xymatrix@R=50pt@C=50pt{ \WCart^{\mathrm{HT}} \ar@{^{(}->}[r] \ar[d] & \WCart \ar[d]^{F} \\
\Spf(\Z_p) \ar[r]^-{ \rho_{\dR}} & \WCart. }$$
In \S\ref{subsection:Frobenius-on-stack}, we show that this diagram is also close to being a pushout square: for every quasi-coherent complex
$\mathscr{E}$ on $\WCart$, the induced diagram
$$ \xymatrix@R=50pt@C=50pt{ \RGamma( \WCart^{\mathrm{HT}}, F^{\ast}(\mathscr{E})|_{ \WCart^{\mathrm{HT}}} )& \RGamma( \WCart, F^{\ast}( \mathscr{E} )) \ar[l] \\
\RGamma( \Spf(\Z_p), \rho_{\dR}^{\ast}(\mathscr{E}) ) \ar[u] & \RGamma( \WCart, \mathscr{E} ) \ar[l] \ar[u] }$$
is a pullback square (Theorem \ref{theorem:Frobenius-pushout-square}). Heuristically, one can think of this result as quantifying the extent to which
the (non-existent) initial prism $(A_{\mathrm{init}}, I_{\mathrm{init}} )$ fails to be perfect.

\begin{remark}
The Cartier-Witt stack $\WCart$ was introduced independently by Drinfeld, who denotes it by $\Sigma$. We refer to the reader to \cite{drinfeld-prismatic} for more details.
\end{remark}

\subsection{Cartier-Witt Divisors}\label{subsection:prismatic-stack}

Let $X$ be a scheme. Recall that a {\it Cartier divisor in $X$} is a closed subscheme $D \subseteq X$ with the property that the ideal sheaf
$\calO(-D) \subseteq \calO_{X}$ is an invertible $\calO_{X}$-module. Beware that this condition generally does not behave well with respect to base change:
if $f: Y \rightarrow X$ is a morphism of schemes and $D \subseteq X$ is a Cartier divisor, then the inverse image $f^{-1}(D) \subseteq Y$ need not be a Cartier divisor.
For this reason, it will be convenient to work with the following more general notion:

\begin{definition}\label{definition:generalized-Cartier-divisor}
Let $X$ be a scheme. A {\it generalized Cartier divisor of $X$} is a pair $( \mathcal{I}, \alpha )$, where $\mathcal{I}$ is an invertible $\calO_{X}$-module
and $\alpha: \mathcal{I} \rightarrow \calO_{X}$ is a a morphism of $\calO_{X}$-modules. We let $\Cart(X)$ denote the category whose objects
are generalized Cartier divisors $(\mathcal{I}, \alpha)$, where a morphism from $(\mathcal{I}, \alpha)$ to $(\mathcal{I}', \alpha')$ is an
isomorphism of $\calO_{X}$-modules $\rho: \mathcal{I} \xrightarrow{\sim} \mathcal{I}'$ satisfying $\alpha = \alpha' \circ \rho$.

In the special case where $X = \Spec(R)$ is an affine scheme, we will denote the category $\Cart(X)$ by $\Cart(R)$; we will identify
its objects with pairs $(I, \alpha)$ where $I$ is an invertible $R$-module and $\alpha: I \rightarrow R$ is an $R$-linear map.
\end{definition}

\begin{remark}
Let $X$ be a scheme and let $D \subseteq X$ be a Cartier divisor, let $\calO(-D)$ denote its ideal sheaf, and let $\iota: \calO(-D) \hookrightarrow \calO_{X}$ be the inclusion map.
Then the pair $(\calO(-D), \iota)$ is a generalized Cartier divisor of $X$. This construction determines a fully faithful embedding from
the set of Cartier divisors of $X$ (regarded as a category having only identity morphisms) to the groupoid of generalized Cartier divisors in $X$.
The essential image of this embedding consists of those generalized Cartier divisors $( \calI, \alpha )$ for which the map $\alpha: \calI \rightarrow \calO_{X}$ is a monomorphism.
\end{remark}

\begin{remark}[Functoriality]
Let $f: X \rightarrow Y$ be a morphism of schemes. Then $f$ determines a pullback functor $\Cart(Y) \rightarrow \Cart(X)$, which carries each
generalized Cartier divisor $( \calI, \alpha)$ of $Y$ to the pullback $( f^{\ast}(\calI), f^{\ast}(\alpha) )$, regarded as a generalized Cartier divisor of $X$.
Beware that if $\alpha$ is a monomorphism (so that $(\calI, \alpha)$ corresponds to a Cartier divisor of $Y$), then $f^{\ast}(\alpha)$ need not be a monomorphism
(so $(f^{\ast}(\calI), f^{\ast}(\alpha) )$ need not correspond to a Cartier divisor of $X$).
\end{remark}

The construction $R \mapsto \Cart(R)$ determines a functor from the category of commutative rings to the $2$-category of groupoids, which is a stack for the fpqc topology.
We will denote this stack by $\Cart$ and refer to it as the {\it moduli stack of generalized Cartier divisors}. In fact, it is an Artin stack, which can be identified with the stack-theoretic quotient
$[\mathbf{A}^{1} / \mathbf{G}_m]$.

\begin{definition}\label{definition:Cartier-Witt-divisor}
Let $R$ be a commutative ring in which $p$ is nilpotent, let $W(R)$ denote the ring of Witt vectors of $R$, and let $(I, \alpha) \in \Cart(W(R)))$ be a generalized Cartier divisor
of the affine scheme $\Spec( W(R) )$ (so that $I$ is an invertible $W(R)$-module and $\alpha: I \rightarrow W(R)$ is a morphism of $W(R)$-modules)). We will say that the pair $(I,\alpha)$ is a {\it Cartier-Witt divisor of $R$} if the following conditions are satisfied:
\begin{itemize}
\item[$(a)$] The image of the map $I \xrightarrow{\alpha} W(R) \twoheadrightarrow R$ is a nilpotent ideal of $R$.
\item[$(b)$] The image of the map $I \xrightarrow{\alpha} W(R) \xrightarrow{\delta} W(R)$ generates the unit ideal of $W(R)$.
\end{itemize}
We let $\WCart(R)$ denote the full subcategory of $\Cart(W(R))$ spanned by the Cartier-Witt divisors of $R$. By convention,
we define $\WCart(R) = \emptyset$ when $p$ is not nilpotent in $R$.
\end{definition}

\begin{remark}[From Cartier-Witt Divisors to Prisms]\label{remark:Cartier-Witt-divisor-prism}
Let $R$ be a commutative ring and let $I \subseteq W(R)$ be an invertible ideal. Then the following conditions are equivalent:
\begin{itemize}
\item The pair $(W(R), I)$ is a prism. Moreover, the $(p,I)$-adic topology on $W(R)$ is a refinement of the $V$-adic topology (or equivalently that the image $p$ and $I$ under $W(R) \twoheadrightarrow R$ is nilpotent).

\item The pair $(I, \iota)$ is a Cartier-Witt divisor of $R$, where $\iota: I \hookrightarrow W(R)$ is the inclusion map
(by convention, this includes the condition that $p$ is nilpotent in $R$).
\end{itemize}
Consequently, there is a close relationship between prism structures on $W(R)$ and Cartier-Witt divisors of $R$.
\end{remark}

\begin{remark}\label{remark:arithmetic-lambda-line}
Let $R$ be a commutative ring. Then the restriction map $W(R) \twoheadrightarrow R$ induces a pullback functor
$\mu_{R}: \Cart( W(R) ) \rightarrow \Cart(R)$. In particular, if $(I, \alpha)$ is a Cartier-Witt divisor of $R$,
then its image under $\mu_{R}$ is a generalized Cartier divisor of $R$. This construction depends functorially on $R$,
and therefore determines a morphism of stacks
$$ \mu: \WCart \rightarrow \Cart = [\mathbf{A}^{1} / \mathbf{G}_{m}].$$

Beware that, if $(I,\alpha)$ is a Cartier-Witt divisor of $R$, then its image under $\mu_R$ is {\em never}
a Cartier divisor of $\Spec(R)$ (except in the trivial case $R \simeq 0$), since the induced map
$R \otimes_{W(R)} I \rightarrow R$ takes values in the nilradical of $R$. In other words, the morphism
$\mu: \WCart \rightarrow [\mathbf{A}^{1} / \mathbf{G}_m]$ factors through the substack
$[\widehat{\mathbf{A}}^{1} / \mathbf{G}_m]$, where $\widehat{\mathbf{A}}^{1}$ denotes the formal affine line.
\end{remark}

\begin{remark}[Functoriality]\label{remark:functoriality-WCart}
Let $f: R \rightarrow S$ be a homomorphism of commutative rings. Then $f$ induces a ring homomorphism
$W(f): W(R) \rightarrow W(S)$, which in turn induces a pullback functor
$$ W(f)^{\ast}: \Cart( W(S) ) \rightarrow \Cart( W(R) ).$$
It is not difficult to see that this pullback functor carries Cartier-Witt divisors of $\Spec(S)$ to Cartier-Witt divisors of $R$,
and therefore restricts to a functor $\WCart(S) \rightarrow \WCart(R)$. We can therefore view the construction
$R \mapsto \WCart(R)$ as a functor from the category of commutative rings to the $2$-category of groupoids.
We will denote this functor by $\WCart$ and refer to it as the {\it Cartier-Witt stack}. It is not difficult to see that
$\WCart$ satisfies descent for the fpqc topology (and, in particular, for the \'{e}tale topology).
\end{remark}

\begin{warning}
Let $R$ be a commutative ring and let $I \subseteq W(R)$ define a prism structure on the ring of Witt vectors $W(R)$.
Assume that the $(p,I)$-adic topology on $W(R)$ is a refinement of the $V$-adic topology, so that the inclusion map
$\iota: I \hookrightarrow W(R)$ determines a Cartier-Witt divisor of $R$. If $f: R \rightarrow S$ is a ring homomorphism,
then the pullback $W(f)^{\ast}(I, \iota)$ is a Cartier-Witt divisor of $\Spec(S)$. Beware that this Cartier-Witt divisor need not
correspond to a prism structure on $W(S)$, because the induced map $I \otimes_{ W(R) } W(S) \rightarrow W(S)$ need not be injective.
\end{warning}

\subsection{Relationship with Prisms}\label{subsection:cartier-witt-atlas}

We now describe an explicit presentation of the Cartier-Witt stack $\WCart$ as the quotient of an affine formal scheme by the action of affine group scheme.

\begin{construction}\label{construction:Cartier-Witt-atlas}
Let $R$ be a commutative ring in which $p$ is nilpotent and let $W(R)$ denote the ring of Witt vectors of $R$. Every element of $W(R)$ has a unique Teichm\"{u}ller expansion
$\sum_{n \geq 0} V^n [a_n]$, where each $a_n$ is an element of $R$. We let $\WCart_{0}(R)$ denote the subset of $W(R)$ consisting of those Witt vectors
$\sum_{n \geq 0} V^n [a_n]$ where $a_0 \in R$ is nilpotent and $a_1 \in R$ is a unit. By convention, we define $\WCart_{0}(R) = \emptyset$ when $p$ is not nilpotent in $R$.
Note that the functor 
$$\WCart_0: \{ \textnormal{Commutative rings} \} \rightarrow \{ \textnormal{Sets} \}$$ is (representable by) an affine formal scheme $\Spf(A^{0})$, where $A^{0}$ denotes the completion of the localized polynomial ring $\Z[ a_0, a_1^{\pm 1}, a_2, a_3, \cdots ]$ with respect to the ideal $(p, a_0)$.

For every commutative ring $R$, let $W(R)^{\times}$ denote the group of units of the commutative ring $W(R)$. The functor $R \mapsto W(R)^{\times}$ is (representable by) an affine group scheme over $\Z$, which we will denote by $W^{\times}$. For every commutative ring $R$, the action of $W(R)^{\times}$ on $W(R)$ (by multiplication) preserves the subset $\WCart_0(R) \subseteq W(R)$ appearing in Construction \ref{construction:Cartier-Witt-atlas}; this determines an action of the affine group scheme $W^{\times}$ on the affine formal scheme $\WCart_{0} \simeq \Spf(A^{0})$.
\end{construction}

\begin{remark}\label{remark:covering-of-stack}
Let $R$ be a commutative ring and let $\alpha: I \rightarrow W(R)$ be a generalized Cartier divisor of $\Spec( W(R) )$. We will say that
the pair $(I, \alpha)$ is {\it principalized} if $I$ is equal to $W(R)$, so that $\alpha$ is given by multiplication by some element $f = \alpha(1) \in W(R)$.
In this case, the pair $(I, \alpha)$ is a Cartier-Witt divisor of $R$ (Definition \ref{definition:Cartier-Witt-divisor}) if and only if $f$ belongs to the subset $\WCart_{0}(R) \subseteq W(R)$.

Let $\WCart_{\mathrm{P}}(R)$ denote the full subcategory of $\WCart(R)$ spanned by the principalized Cartier-Witt divisors. This category can be described more concretely as follows:
\begin{itemize}
\item The objects of $\WCart_{\mathrm{P}}(R)$ are given by elements $f \in \WCart_0(R) \subseteq W(R)$ (which we identify with the Cartier-Witt divisor given by the multiplication map $W(R) \xrightarrow{f} W(R)$).
\item If $f$ and $f'$ are elements of $\WCart_0(R)$, then a morphism from $f$ to $f'$ in the category $\WCart_{\mathrm{P}}(R)$ is an invertible element $u \in W(R)$ satisfying $f' = u f$.
\end{itemize}
Consequently, the construction $f \mapsto (W(R), \alpha_{f} )$ determines a functor from the set $\WCart_{0}(R)$ (regarded as a category having only identity morphisms)
to the groupoid $\WCart(R)$ of Cartier-Witt divisors. This construction depends functorially on $R$, and therefore determines a morphism from the
affine formal scheme $\WCart_{0} = \Spf( A^{0} )$ to the Cartier-Witt stack $\WCart$.
\end{remark}

\begin{proposition}\label{proposition:presentation-as-quotient}
The morphism $\WCart_{0} \rightarrow \WCart$ of Remark \ref{remark:covering-of-stack} exhibits
$\WCart$ as the quotient $[\WCart_{0} / W^{\times}]$, formed in the $2$-category of stacks with respect to the Zariski topology on $\Spec(R)$.
\end{proposition}

\begin{proof}
By virtue of Remark \ref{remark:covering-of-stack} it will suffice to show that the map $\WCart_0 \rightarrow \WCart$ is Zariski-locally essentially surjective.
Let $R$ be a commutative ring and let $\alpha: I \rightarrow W(R)$ be a Cartier-Witt divisor in $\Spec(R)$. 
Unwinding the definitions, we see that the pair $(I, \alpha)$ belongs to the essential image of the map $\WCart_0(R) \rightarrow \WCart(R)$ if and only if
the invertible $W(R)$-module $I$ is isomorphic to $W(R)$. Since $p$ is nilpotent in $R$, this is equivalent to the requirement that
the invertible $R$-module $R \otimes_{W(R)} I$ is isomorphic to $R$, which is true Zariski-locally on $\Spec(R)$.
\end{proof}

It will be convenient to reformulate Proposition \ref{proposition:presentation-as-quotient} using the language of prisms.
We begin by observing that every prism $(A,I)$ determines a morphism of stacks $\rho_{A}: \Spf(A) \rightarrow \WCart$.

\begin{construction}[From Prisms to Cartier-Witt Divisors]\label{construction:point-of-prismatic-stack}
Let $(A,I)$ be a prism and let $f: A \rightarrow R$ be a homomorphism of commutative rings, so that $f$ lifts uniquely to a morphism of $\delta$-rings
$\widetilde{f}: A \rightarrow W(R)$. Let $\widetilde{f}^{\ast}(I)$ denote the tensor product $W(R) \otimes_{A} I$, which we regard as an invertible $W(R)$-module.
Then the inclusion $I \hookrightarrow A$ induces a $W(R)$-linear map $\alpha: \widetilde{f}^{\ast}(I) \rightarrow W(R)$, and the pair
$( \widetilde{f}^{\ast}(I), \alpha)$ can be regarded as a generalized Cartier divisor of $\Spec(W(R))$. The pair
$( \widetilde{f}^{\ast}(I), \alpha)$ is a Cartier-Witt divisor of $R$ precisely when the image of the ideal $I + (p)$ is nilpotent in $R$.
Consequently, the construction $(f: A \rightarrow R ) \mapsto ( \widetilde{f}^{\ast}(I), \alpha)$ determines a morphism of stacks
$\rho_{A}: \Spf(A) \rightarrow \WCart$, where $\Spf(A)$ denotes the formal scheme determined by the $(p,I)$-adic topology on $A$.
\end{construction}

\begin{example}\label{example:rho-A-special-case}
Let $A^{0} = \Z[ a_0, a_1^{\pm 1}, a_2, a_3, \cdots ]^{\wedge}_{(p,a_0)}$ denote the coordinate ring of the formal scheme $\WCart_{0}$ described in
Construction \ref{construction:Cartier-Witt-atlas}. Note that, for every commutative ring $R$, the Witt vector Frobenius $F: W(R) \rightarrow W(R)$
carries the subset $\WCart_{0}(R)$ into itself. The resulting map $F: \WCart_{0}(R) \rightarrow \WCart_{0}(R)$ depends functorially on $R$
and therefore determines an endomorphism of the formal scheme $\WCart_{0}$ (reducing the Frobenius modulo $p$), which we can identify with a $\delta$-structure
on the commutative ring $A^{0}$. Let $I^{0} \subseteq A^{0}$ be the principal ideal generated by the element $a_0$. Then the pair $(A^0, I^0 )$ is a prism,
and the map $\rho_{A^{0}}: \Spf( A^{0} ) \rightarrow \WCart$ of Construction \ref{construction:point-of-prismatic-stack} can be identified with the morphism
$\WCart_{0} \rightarrow \WCart$ of Remark \ref{remark:covering-of-stack}.
\end{example}

\begin{example}\label{example:de-Rham-point}
Applying Construction \ref{construction:point-of-prismatic-stack} to the crystalline prism $(\Z_p, (p) )$, we obtain a $\Spf(\Z_p)$-valued point of the Cartier-Witt stack.
We will refer to this point as the {\it de Rham point} and denote it by $\rho_{\dR}: \Spf(\Z_p) \rightarrow \WCart$. We will later justify this terminology by showing
that pullback along the de Rham point $\rho_{\Z_p}: \Spf(\Z_p) \rightarrow \WCart$ implements the comparison between (absolute) prismatic cohomology
and de Rham cohomology (see Proposition \ref{proposition:refined-de-Rham-comparison}).
\end{example}

\begin{remark}\label{remark:universal-of-A0}
Let $(A^{0}, I^{0})$ be the prism of Example \ref{example:rho-A-special-case}. For every prism $(B,J)$, the construction $f \mapsto f(a_0)$ induces a bijection
$$ \{ \text{Prism homomorphisms $(A^{0}, I^{0} ) \rightarrow (B,J)$} \} \rightarrow \{ \text{Generators of $J$} \}.$$
\end{remark}

\begin{proposition}\label{proposition:pullback-diagram-prisms}
Let $(A,I)$ be a bounded prism and let $(B,J)$ be a transveral prism, so that $(A,I)$ and $(B,J)$ admit a coproduct $(C,K)$ in the category of prisms
(Proposition \ref{proposition:transversal-coproduct}). Then the diagram of stacks
$$ \xymatrix@R=50pt@C=50pt{ \Spf(C) \ar[r] \ar[d] & \Spf(B) \ar[d]^{ \rho_{B} } \\
\Spf(A) \ar[r]^-{ \rho_{A} } & \WCart }$$
is a pullback square.
\end{proposition}

\begin{proof}
The conclusion of Proposition \ref{proposition:pullback-diagram-prisms} can be tested Zariski-locally on $\Spf(A)$ and $\Spf(B)$.
We may therefore assume without loss of generality that the ideals $I$ and $J$ are generated by distinguished elements
$d_{A} \in A$ and $d_{B} \in B$, respectively. Let $(A^{0}, I^{0} )$ be the prism of Example \ref{example:rho-A-special-case}, so that there is a unique morphism of prisms
$(A^{0}, I^{0}) \rightarrow (A,I)$ carrying $a_0$ to $d_A$. The prism $(A^{0}, I^{0})$ is also bounded (in fact, it is transversal), so that
the prisms $(A^0, I^{0} )$ and $(B,J)$ admit a coproduct $(C^0, K^0)$ in the category of prisms (Proposition \ref{proposition:transversal-coproduct}).
We have a commutative diagram
$$ \xymatrix@R=50pt@C=50pt{ \Spf( C ) \ar[r] \ar[d] & \Spf( C^{0} ) \ar[r] \ar[d] & \Spf(B) \ar[d]^{ \rho_{B} } \\
\Spf( A ) \ar[r] & \Spf( A^{0} ) \ar[r]^-{ \rho_{A^0}} & \WCart, }$$
where the square on the left is a pullback (note that the vertical maps are both flat, by virtue of Proposition \ref{proposition:transversal-coproduct}).
It will therefore suffice to show that the square on the right is also a pullback square. 

By virtue of Proposition \ref{proposition:presentation-as-quotient}, we can identify the pullback $\Spf(A^{0}) \times_{ \WCart} \Spf(A^{0} )$ with the affine formal scheme
$\Spf(B) \times W^{\times}$. Concretely, this affine formal scheme is given by the formal spectrum of a $\delta$-ring $B'$, given by the $(p,J)$-completion
of the free $\delta$-$B$-algebra $B\{u^{\pm 1} \}_{\delta}$ generated by an invertible element $u$. To complete the proof, it will suffice to show that
the diagram above exhibits $(B', JB')$ as a coproduct of $(B,J)$ with $(A^{0}, I^{0})$ in the category of prisms. Equivalently, we wish to show that if
$(B'', JB'')$ is any prism over $(B,J)$, then the canonical map
$$ \{ \text{Prism maps $f: (A^{0}, I^{0} ) \rightarrow (B'',JB'')$} \} \rightarrow \{ \text{Invertible elements of $B'$} \}
\quad \quad f \mapsto \frac{ f(a_0)}{ d_B }$$
is a bijection, which follows from Remark \ref{remark:universal-of-A0}.
\end{proof}

\begin{corollary}
Let $(B,J)$ be a prism. Then the map $\rho_{B}: \Spf(B) \rightarrow \WCart$ of Construction \ref{construction:point-of-prismatic-stack} is affine.
That is, for every for every commutative ring $R$ and every Cartier-Witt divisor of $R$,
the fiber product $\Spf(B) \times_{\WCart} \Spec(R)$ is an affine $R$-scheme.
\end{corollary}

\begin{proof}
By virtue of Proposition \ref{proposition:transversal-approximation}, we may assume without loss of generality that the prism $(B,J)$ is transversal.
The assertion is local on $\Spec(R)$, so we may also assume without loss of generality that the Cartier-Witt divisor of $R$ is principalized: that is,
it is classified by a morphism of formal schemes $\Spec(R) \rightarrow \WCart_{0} \simeq \Spf(A^{0})$, where $(A^{0}, I^{0})$ is the prism of
Example \ref{example:rho-A-special-case}. Let $(C,K)$ denote the coproduct of $(A^{0}, I^{0})$ with $(B,J)$ in the category of prisms.
Using Proposition \ref{proposition:pullback-diagram-prisms}, we can identify $\Spf(B) \times_{ \WCart} \Spec(R)$ with the affine $R$-scheme
$\Spf(C) \times_{ \Spf(A^{0})} \Spec(R)$.
\end{proof}

\begin{corollary}\label{corollary:flatness-v-transversality}
Let $(B,J)$ be a bounded prism. Then $(B,J)$ is transversal if and only if the affine map $\rho_{B}: \Spf(B) \rightarrow \WCart$ is flat.
Moreover, if $(B,J)$ is transversal and nonzero, then $\rho_{B}$ is faithfully flat.
\end{corollary}

\begin{proof}
Assume first that $(B,J)$ is transversal; we will show that $\rho_{B}$ is flat. Without loss of generality, we may assume $B \neq 0$; in this case, we will show that $\rho_{B}$ is faithfully flat. Let $(A^{0}, I^{0})$ be the prism of Example \ref{example:rho-A-special-case}. By virtue of Proposition \ref{proposition:presentation-as-quotient}, it will suffice to show that the $\pi: \Spf(B) \times_{ \WCart } \Spf(A^{0}) \rightarrow \Spf(A^{0} )$ is faithfully flat. Using Proposition \ref{proposition:pullback-diagram-prisms}, we can identify the source of $\pi$ with the formal spectrum
$\Spf(C)$, where $(C,K)$ is a coproduct of $(A^{0}, I^{0})$ with $(B, J)$ in the category of prisms. We are therefore reduced to checking that the commutative ring
$C$ is $(p,I^0)$-completely faithfully flat over $A^{0}$, which follows from Proposition \ref{proposition:characterize-transversal}.

For the converse, assume that $\rho_{B}$ is flat. Then $(C,K)$ is flat over the transversal prism $(A^{0}, I^{0})$, and is therefore transversal (Remark \ref{remark:flatness-and-transversality}).
It follows from Proposition \ref{proposition:characterize-transversal} that the homomorphism of prisms $(B,J) \rightarrow (C,K)$ is faithfully flat, so that $(B,J)$ is also
transversal (Remark \ref{remark:flatness-and-transversality}).
\end{proof}

\subsection{Complexes on the Cartier-Witt Stack}\label{subsection:complexes-on-WCart}

For our purposes, the Cartier-Witt stack is primarily a bookkeeping device for carrying information about sheaves.

\begin{definition}\label{definition:DWCart}
For every commutative ring $R$, let $\calD(R)$ denote the derived $\infty$-category of $R$-modules.
We let $\calD( \WCart )$ denote the inverse limit
$$ \varprojlim_{\Spec(R) \rightarrow \WCart} \calD(R),$$
indexed by the category of points of the Cartier-Witt stack $\WCart$, which we regard as a symmetric monoidal stable $\infty$-category.
We will refer to the objects of $\calD( \WCart )$ as {\it quasi-coherent complexes on $\WCart$},
and to $\calD(\WCart)$ as {\it the $\infty$-category of quasi-coherent complexes on $\WCart$}.
\end{definition}

\begin{remark}
Stated more informally, a quasi-coherent complex $\mathscr{F} \in \calD(\WCart)$ can be viewed
as a rule which associates to each commutative ring $R$ and each Cartier-Witt divisor $\alpha: I \rightarrow W(R)$
a complex of $R$-modules $\mathscr{F}_{\alpha}$, depending functorially on $R$ (up to quasi-isomorphism).
\end{remark}

\begin{example}\label{example:structure-sheaf-of-WCart}
We let $\calO_{\WCart}$ denote the unit object of the symmetric monoidal $\infty$-category $\calD( \WCart )$, which we
refer to as the {\it structure sheaf} of the Cartier-Witt stack. Concretely, $\calO_{\WCart}$ associates to
every Cartier-Witt divisor $\alpha: I \rightarrow W(R)$ the underlying commutative ring $R \in \calD(R)$.
\end{example}

\begin{example}[The Hodge-Tate Ideal Sheaf]\label{example:HTI}
For every commutative ring $R$ and every Cartier-Witt divisor $\alpha: I \rightarrow W(R)$, extension of scalars
along the restriction map $W(R) \twoheadrightarrow R$ determines an invertible $R$-module $R \otimes_{W(R)} I$.
The construction $(R, \alpha) \mapsto R \otimes_{W(R)} I$ determines an invertible object of the $\infty$-category
$\calD( \WCart )$, which we will denote by $\mathcal{I}$ and refer to as the {\it Hodge-Tate ideal sheaf}.
Note that the construction $(R, \alpha) \mapsto \alpha$ determines a morphism of quasi-coherent complexes $\calI \rightarrow \calO_{\WCart}$.
\end{example}

Let $(A,I)$ be a bounded prism and let $\rho_{A}: \Spf(A) \rightarrow \WCart$ be the morphism of Construction \ref{construction:point-of-prismatic-stack}.
Then pullback along $\rho_{A}$ induces a functor $\rho_{A}^{\ast}: \calD( \WCart ) \rightarrow \calD( \Spf(A) ) \simeq \widehat{\calD}(A)$, where
$\widehat{\calD}(A)$ denotes the full subcategory of $\calD(A)$ spanned by those complexes of $A$-modules which are $(p,I)$-complete.

\begin{proposition}\label{proposition:DWCart-prism-description}
The preceding construction induces an equivalence of $\infty$-categories
$$ \calD( \WCart ) \rightarrow \varprojlim_{ (A,I) } \widehat{\calD}(A),$$
where the limit is indexed by the category of all bounded prisms $(A,I)$.
\end{proposition}

\begin{warning}
\label{warn:boundedqcoh}
In the statement of Proposition \ref{proposition:DWCart-prism-description}, the boundedness hypothesis on $(A,I)$
guarantees that we can identify $\widehat{\calD}(A)$ with the $\infty$-category of quasi-coherent complexes on the formal scheme $\Spf(A)$.
However, the limit $\varprojlim_{(A,I)} \widehat{\calD}(A)$ does not change if we enlarge the index category to include non-bounded prisms,
by virtue of Corollary \ref{siftedness}.
\end{warning}

\begin{remark}
Thanks to Proposition~\ref{proposition:DWCart-prism-description} and Warning~\ref{warn:boundedqcoh}, we can also regard $\mathcal{D}(\WCart)$ as the $\infty$-category of crystals of $(p,I_\Prism)$-complete complexes on the absolute prismatic site $(\mathrm{Spf}(\mathbf{Z}_p)_\Prism, \mathcal{O}_\Prism)$. 
\end{remark}

\begin{example}[Breuil-Kisin Twists]\label{example:BK-twist-as-line}
By virtue of Proposition \ref{proposition:DWCart-prism-description} and Remark \ref{remark:BK-functoriality}, the construction
$(A,I) \mapsto A\{1\}$ determines an invertible object of the $\infty$-category $\calD( \WCart )$, which we will denote by
$\calO_{\WCart}\{1\}$ and refer to as the {\it Breuil-Kisin line bundle} on the Cartier-Witt stack $\WCart$. More generally, for every integer $n$ and
every object  $\mathscr{E} \in \calD(\WCart)$, we let  $\mathscr{E}\{n\}$ denote the tensor product of $\mathscr{E}$ with the $n$th power of $\mathscr{O}_{\WCart}\{1\}$.
\end{example}

Our proof of Proposition \ref{proposition:DWCart-prism-description} will make use of an auxiliary construction.

\begin{notation}\label{notation:simplicial-prism}
Let $(A^{0}, I^{0} )$ denote the prism of Example \ref{example:rho-A-special-case}. For each $n \geq 0$, we let $(A^{n}, I^{n})$ denote the coproduct
of $(n+1)$ copies of $(A^{0}, I^{0})$ in the category of prisms (which exists by virtue of Proposition \ref{proposition:transversal-coproduct}),
so that the construction $n \mapsto (A^{n}, I^{n} )$ determines a cosimplicial prism $( A^{\bullet}, I^{\bullet} )$.
\end{notation}

\begin{lemma}\label{lemma:tot-description}
Let $(A^{\bullet}, I^{\bullet} )$ be the cosimplicial prism of Notation \ref{notation:simplicial-prism}. Then pullback along the maps
$\rho_{A^{\bullet}}: \Spf( A^{\bullet} ) \rightarrow \WCart$ of Construction \ref{construction:point-of-prismatic-stack}
determines an equivalence of $\infty$-categories
$$ \calD( \WCart ) \rightarrow \Tot(  \calD( \Spf( A^{\bullet} ) ) ) \simeq \Tot( \widehat{\calD}(A^{\bullet} ) ).$$
\end{lemma}

\begin{proof}
It follows from Proposition \ref{proposition:presentation-as-quotient} that the map $\Spf( A^{0} ) \rightarrow \WCart$
is surjective with respect to the Zariski topology and from Proposition \ref{proposition:pullback-diagram-prisms} that, for each
$n \geq 0$, we can identify $\Spf(A^{n} )$ with the $(n+1)$-fold fiber power of $\Spf( A^{0} )$ over $\WCart$.
Lemma \ref{lemma:tot-description} now follows from cohomological descent (with respect to the Zariski topology).
\end{proof}

\begin{proof}[Proof of Proposition \ref{proposition:DWCart-prism-description}]
By virtue of Lemma \ref{lemma:tot-description}, it will suffice to show that the restriction functor
$$ \varprojlim_{ (A,I) } \widehat{\calD}( A ) \rightarrow \Tot( \widehat{\calD}( A^{n} ) )$$
is an equivalence of $\infty$-categories. Using flat cohomological descent, we can assume that the limit on the left hand side is
indexed by the category $\calC$ of bounded prisms $(A,I)$ for which the ideal $I = (d)$ is principal. In this case, the desired result follows
from the cofinality of the functor $$\mathbf{\Delta} \rightarrow \calC \quad \quad [n] \mapsto A^{n}.$$
\end{proof}

Lemma \ref{lemma:tot-description} provides another useful description of the $\infty$-category $\calD( \WCart )$. Stated informally,
it asserts that quasi-coherent complexes on the Cartier-Witt stack can be identified with cosimplicial complexes $M^{\bullet}$ of modules over the cosimplicial
commutative ring $A^{\bullet}$ of Notation \ref{notation:simplicial-prism}, where each $M^{n}$ is $(p, I^{n} )$-complete
(and each of the transition maps $A^{n'} \widehat{\otimes}_{A^{n} }^{L} M^{n} \rightarrow M^{n'}$ is invertible). This description is useful for computing
the global sections of the quasi-coherent complex on $\WCart$:

\begin{corollary}\label{corollary:global-sections-on-WCart}
The functor
$$ \calD(\Z) \rightarrow \calD( \WCart ) \quad \quad M \mapsto M \otimes \calO_{\WCart}$$
admits a right adjoint $\RGamma( \WCart, \bullet ): \calD( \WCart ) \rightarrow \calD(\Z)$,
which we will refer to as the {\it global sections functor}. Under the identification $\calD( \WCart ) \simeq \Tot( \widehat{\calD}( A^{\bullet} ) )$ of
Lemma \ref{lemma:tot-description}, this functor is given by the construction
$$M^{\bullet} \in \Tot( \widehat{\calD}( A^{\bullet} ) ) \mapsto \Tot( M^{\bullet} ) \in \widehat{\calD}(\Z_p) \subset \calD(\Z).$$
\end{corollary}

\subsection{The Hodge-Tate Divisor}\label{subsection:HT-divisor}

Let $R$ be a commutative ring. Every Cartier-Witt divisor $(I, \alpha)$ for $R$ determines a generalized Cartier divisor in the
affine scheme $\Spec(R)$, given by extending scalars along the restriction map $W(R) \twoheadrightarrow R$.
Allowing $R$ to vary, we obtain a (generalized) Cartier divisor in the Cartier-Witt stack $\WCart$.

\begin{definition}\label{definition:Hodge-Tate-locus}
Let $R$ be a commutative ring. We let $\WCart^{\mathrm{HT}}(R)$ denote the full subcategory of $\WCart(R)$
spanned by those Cartier-Witt divisors $\alpha: I \rightarrow W(R)$ for which the composite map $I \xrightarrow{\alpha} W(R) \twoheadrightarrow R$
is equal to zero. The construction $R \mapsto \WCart^{\mathrm{HT}}(R)$ determines a closed substack of the Cartier-Witt stack $\WCart$.
We denote this closed substack by $\WCart^{\mathrm{HT}}$ and refer to it as the {\it Hodge-Tate divisor}.
\end{definition}

\begin{remark}\label{remark:HT-point-of-prismatic-stack}
Let $(A,I)$ be a prism and regard the commutative ring $A$ as equipped with the $(p,I)$-adic topology. 
Let $\rho_{A}: \Spf(A) \rightarrow \WCart$ be the the morphism described in Construction \ref{construction:point-of-prismatic-stack}.
Then $\rho_{A}$ carries the formal subscheme $\Spf(A/I) \subset \Spf(A)$ to the Hodge-Tate divisor $\WCart^{\mathrm{HT}}$,
and therefore restricts to a morphism $\rho_{A}^{\mathrm{HT}}: \Spf(A/I) \rightarrow \WCart^{\mathrm{HT}}$. Moreover,
the diagram
$$ \xymatrix@R=50pt@C=50pt{ \Spf(A/I) \ar@{^{(}->}[d] \ar[r]^-{ \rho_{A}^{\mathrm{HT}}} & \WCart^{\mathrm{HT}} \ar@{^{(}->}[d] \\
\Spf(A) \ar[r]^-{\rho_{A} } & \WCart }$$
is a pullback square.
\end{remark}

\begin{example}\label{example:points-from-perfectoid-rings}
Let $\overline{A}$ be a perfectoid ring. Then we identify $\overline{A}$ with a quotient $A/I$, where $(A,I)$ is a perfect prism (which is uniquely determined up to unique isomorphism). Applying Remark \ref{remark:HT-point-of-prismatic-stack}, we obtain a morphism $\rho_{A}^{\mathrm{HT}}: \Spf(\overline{A}) \rightarrow \WCart^{\mathrm{HT}}$,
which depends functorially on the perfectoid ring $\overline{A}$.
\end{example}

Note that the Hodge-Tate divisor $\WCart^{\mathrm{HT}}$ can be described as the central fiber of the morphism
$\mu: \WCart \rightarrow [\widehat{\mathbf{A}}^{1} / \mathbf{G}_{m}]$ described in Remark \ref{remark:arithmetic-lambda-line}; that is,
we have a pullback diagram of stacks
$$ \xymatrix@R=50pt@C=50pt{ \WCart^{\mathrm{HT}} \ar[r]^-{ \mu^{\mathrm{HT}} } \ar@{^{(}->}[d] & [\{0\} / \mathbf{G}_{m}] \ar@{^{(}->}[d] \\
\WCart \ar[r]^-{\mu} & [\widehat{\mathbf{A}}^{1} / \mathbf{G}_{m}]. }$$
Our next goal is to show that the morphism
$\mu^{\mathrm{HT}}: \WCart^{\mathrm{HT}} \rightarrow [\{0\} / \mathbf{G}_m] = B\mathbf{G}_{m}$ is not far from being an isomorphism (see Theorem \ref{theorem:describe-HT} below).

\begin{construction}\label{construction:fiber-at-eta}
For every commutative ring $R$, let us write $V: W(R) \rightarrow W(R)$ for the Verschiebung operator. If $p$ is nilpotent in $R$, then multiplication by the element
$V(1) \in W(R)$ determines a (principalized) Cartier-Witt divisor for $R$. Since $V(1)$ is annihilated by the restriction map $W(R) \twoheadrightarrow R$,
this Cartier-Witt divisor can be regarded as an object of the category $\WCart^{\mathrm{HT}}(R)$. Allowing $R$ to vary, we obtain a morphism of stacks
$\eta: \Spf(\Z_p) \rightarrow \WCart^{\mathrm{HT}}$.
\end{construction}

Let $R$ be a commutative ring in which $p$ is nilpotent. We let $\Aut(\eta)(R)$ denote the automorphism group of the composite map
$$ \Spec(R) \rightarrow \Spf(\Z_p) \xrightarrow{\eta} \WCart^{\mathrm{HT}}.$$
Unwinding the definitions, we see that $\Aut(\eta)(R)$ can be identified with the subgroup of $W(R)^{\times}$ consisting of those units $u$
which satisfy the equation $u \cdot V(1) = V(1)$ in $W(R)$. Using the identity $u \cdot V(1) = V( F(u) \cdot 1) = V( F(u) )$ (and the injectivity of the Verschiebung operator $V$), we obtain an isomorphism
$$ \Aut(\eta) \simeq \Spf(\Z_p) \times W^{\times}[F],$$
where $W^{\times}[F]$ denotes the affine group scheme over $\Z$ given by the kernel of the Frobenius $F: W^{\times} \rightarrow W^{\times}$.

\begin{remark}\label{remark:comparison-W-times-F}
The map $\mu^{\mathrm{HT}}: \WCart^{\mathrm{HT}} \rightarrow B\mathbf{G}_{m}$ induces a homomorphism
$$ \Spf( \Z_p ) \times W^{\times}[F] \simeq \Aut(\eta) \rightarrow \Aut( \mu^{\mathrm{HT}} \circ \eta ) = \Spf( \Z_p ) \times \mathbf{G}_{m}$$
of formal group schemes over $\Z_p$. Concretely, this homomorphism is induced by the map of affine group schemes
$$W^{\times}[F] \hookrightarrow W^{\times} \twoheadrightarrow \mathbf{G}_{m}$$
obtained by composing the inclusion of $W^{\times}[F]$ into $W^{\times}$ with the epimorphism
$W^{\times} \twoheadrightarrow \mathbf{G}_{m}$ given by restriction.
\end{remark}

\begin{remark}\label{remark:easy-locus}
After extending scalars to $\Z[1/p]$, the group scheme $W^{\times}$ splits as a product of infinitely many copies of
$\mathbf{G}_{m}$ (via the ghost components), with the Witt vector Frobenius acting by a shift. Consequently,
the comparison map $W^{\times}[F] \rightarrow \mathbf{G}_{m}$ of Remark \ref{remark:comparison-W-times-F}
becomes an isomorphism after extending scalars to $\Z[1/p]$.
\end{remark}

Our next goal is to describe the structure of the group scheme $W^{\times}[F]$ at the prime $p$.

\begin{proposition}\label{proposition:flat-Frobenius-on-W}
Let $W \simeq \Spec( \Z[a_0, a_1, \cdots, ])$ be the affine scheme representing the Witt vector functor
$R \mapsto W(R)$. Then the Witt vector Frobenius map $F: W \rightarrow W$ is faithfully flat.
\end{proposition}

\begin{proof}
Let $A = \Z[ a_0, a_1, \cdots ]$ denote the coordinate ring of the affine scheme $W$, so that the Frobenius
determines a ring homomorphism $\varphi_{A}: A \rightarrow A$. We wish to show that $\varphi_{A}$ is faithfully flat.
Since $A$ is torsion-free, it suffices to show establish faithful flatness after inverting $p$ and dividing out by the ideal $pA$.
For the first case, we observe that the localization $A[1/p]$ can be identified with the polynomial ring
$\Z[1/p][\gamma_0, \gamma_1, \gamma_2, \cdots ]$, where $\gamma_i$ denotes the $i$th ghost component,
with the Frobenius acting by $\varphi(\gamma_{i} ) = \gamma_{i+1}$. In particular, the map
$\varphi_{A}$ exhibits $A[1/p]$ as a polynomial algebra over itself on a single generator $\gamma_0$.
In the second case, we observe that the quotient ring $A/pA$ is a polynomial ring $\F_p[ a_0, a_1, \cdots ]$
on which the endomorphism $\varphi_A$ acts by $a_i \mapsto a_i^{p}$.
\end{proof}

\begin{corollary}\label{corollary:W-F-flat}
The Witt vector Frobenius map $F: W^{\times} \rightarrow W^{\times}$ is faithfully flat. In particular,
the group scheme $W^{\times}[F]$ is flat over $\Z$.
\end{corollary}

\begin{notation}\label{notation:G_m-sharp}
Let $\mathbf{G}_{m}$ denote the multiplicative group (regarded as an affine group scheme over $\Z$) and let $\calO_{ \mathbf{G}_{m} } = \Z[ t^{\pm 1}]$
denote its coordinate ring. We let $\calO_{\mathbf{G}_{m}}^{\sharp}$ denote the subring of $\Q \otimes \calO_{\mathbf{G}_{m} } \simeq \Q[ t^{\pm 1} ]$
generated by $t^{-1}$ together with the divided powers $\frac{ (t-1)^{n}}{n!}$ (equivalently, the commutative ring $\calO^{\sharp}_{\mathbf{G}_{m}}$ can be described as the {\em divided power envelope}
of the ideal $(t-1) \subseteq \Z[ t^{\pm 1} ]$). We let $\mathbf{G}_{m}^{\sharp}$ denote the affine scheme $\Spec(\calO_{\mathbf{G}_{m}}^{\sharp} )$.
The inclusion $\calO_{\mathbf{G}_m} \hookrightarrow \calO^{\sharp}_{\mathbf{G}_{m}}$ induces a morphism of schemes $\beta: \mathbf{G}_{m}^{\sharp} \rightarrow \mathbf{G}_{m}$, which restricts to an isomorphism $$ \Spec(\Q) \times \mathbf{G}_{m}^{\sharp} \xrightarrow{\sim} \Spec(\Q) \times \mathbf{G}_{m}.$$
Note that there is a unique group structure on the affine scheme $\mathbf{G}_{m}^{\sharp}$ for which $\beta$ is a homomorphism of
group schemes.
\end{notation}

\begin{remark}\label{remark:Z-p-points}
For every commutative ring $R$, the morphism $\beta: \mathbf{G}_m^{\sharp} \rightarrow \mathbf{G}_m$ of Notation \ref{notation:G_m-sharp}
induces a group homomorphism $\beta(R): \mathbf{G}_{m}^{\sharp}(R) \rightarrow \mathbf{G}_{m}(R) = R^{\times}$. If the commutative ring $R$ is torsion-free, 
then this homomorphism is an injection, whose image can be identified with the subgroup consisting of those elements $u \in R^{\times}$ for which 
$u-1$ has divided powers (that is, each power $(u-1)^{n}$ is divisible by $n!$). In particular, we have a canonical isomorphism $\mathbf{G}_{m}^{\sharp}( \Z_p ) \simeq ( 1 + p \Z_p)^{\times}$.
\end{remark}

\begin{lemma}\label{lemma:sharp-SES}
The comparison map $W^{\times}[F] \rightarrow \mathbf{G}_{m}$ of Remark \ref{remark:comparison-W-times-F} lifts uniquely to
an isomorphism
$$\Spec( \Z_{(p)} ) \times W^{\times}[F] \xrightarrow{\sim} \Spec(\Z_{(p)}) \times \mathbf{G}_{m}^{\sharp}$$
of group schemes over $\Z_{(p)}$.
\end{lemma}

In the language of $\delta$-rings, this lemma asserts that the free $\delta$-$\mathbf{Z}_{(p)}$-algebra $A$ on an element $x \in A$ with $\phi(x)=1$ identifies with $\mathcal{O}(\mathbf{G}_m^\sharp)$. Similarly, Variant~\ref{Gasharpdef} asserts that the free $\delta$-$\mathbf{Z}_{(p)}$-algebra $B$ on an element $x \in A$ with $\phi(x)=0$ identifies with $\mathcal{O}(\mathbf{G}_a^\sharp)$. 

\begin{proof}
Let $S$ denote the coordinate ring of the affine scheme $\Spec( \Z_{(p)} ) \times W^{\times}[F]$. Then
$S$ is a flat $\Z_{(p)}$-algebra (Corollary \ref{corollary:W-F-flat}), and the comparison map
$$ \Spec( \Z_{(p)} ) \times W^{\times}[F] \rightarrow \mathbf{G}_{m}$$
is classified by an invertible element $t \in S$. The Witt vector Frobenius determines a
$\delta$-structure on the ring $S$ satisfying $\varphi_{S}(t) = 1$. In particular, the element
$t-1$ is annihilated by $\varphi_{S}$ and therefore admits a system of divided powers in $S$
(see Lemma~2.35 of \cite{prisms}). We can therefore promote $t$ to a ring homomorphism
$\nu: (\calO_{\mathbf{G}_m}^{\sharp})_{(p)} \rightarrow S$, which we can identify with a morphism
of affine schemes 
$$ \Spec( \Z_{(p)} ) \times W^{\times}[F] \to \Spec(\Z_{(p)}) \times \mathbf{G}_{m}^{\sharp}.$$
By construction, the restriction of $\nu$ to generic fibers is an isomorphism of group schemes over $\Q$
(see Remark \ref{remark:easy-locus}). Since $\Spec( \Z_{(p)} ) \times W^{\times}[F]$ is flat over $\Z_{(p)}$, it follows that $\nu$
is a morphism of group schemes over $\Z_{(p)}$. To complete the proof, it will suffice to show that it is also an isomorphism of special fibers: that is,
that the map $\nu$ induces an isomorphism of $\F_p$-algebras
$$ \overline{\nu}: \calO_{\mathbf{G}_m}^{\sharp} / p \calO_{\mathbf{G}_m}^{\sharp} \rightarrow S/pS.$$

Set $a = \nu( t-1) \in S$, and let $\overline{a}$ denote the image of $a$ in the quotient ring $S/pS$.
For each element $x \in (t-1) \calO_{\mathbf{G}_m}^{\sharp}$, let us write $\alpha(x)$ for the element $- \frac{ x^{p} }{p}$. 
Note that the element $\nu(x) \in S$ is annihilated by the Frobenius map $\varphi_{S}$, and therefore satisfies
$$ \delta( \nu(x) ) = \frac{ \varphi_{ S}( \nu(x) ) - \nu(x)^p}{p} = \nu( \alpha(x) ).$$
The quotient ring $\calO_{\mathbf{G}_m}^{\sharp} / p \calO_{\mathbf{G}_m}^{\sharp}$ has an $\F_p$-basis consisting of the images of monomials
$$(t-1)^{e_0} \alpha(t-1)^{e_1} \alpha^2( t-1)^{e_2} \cdots \alpha^{n}( t-1)^{ e_n},$$ where each exponent $e_i$ is strictly less than $p$. Moreover, the homomorphism $\overline{\rho}$ carries every such monomial
to the image of the product $a^{e_0} \delta( a )^{e_1} \cdots \delta^n( a )^{e_n}$. We conclude by observing that the latter monomials form an $\F_p$-basis for the
quotient ring $$S/pS \simeq \F_p[ a, \delta( a ), \delta^2( a ), \cdots ] / ( a^{p}, \delta( a)^{p}, \cdots ).$$
\end{proof}

\begin{variant}[The group scheme $\mathbf{G}_a^\sharp$]
\label{Gasharpdef}
Let $\mathbf{G}_a = \mathrm{Spec}(\mathbf{Z}[x])$ denote the additive group, and let $\mathbf{G}_a^\sharp = \mathrm{Spec}(\mathbf{Z}[\{\frac{x^n}{n!}\}_{n \geq 1}])$ denote the group scheme obtained by taking the divided power envelope at $0$. Write $W[F] = \ker(W \xrightarrow{F} W)$ for the kernel of Frobenius on $W$, again regarded as a group scheme over $\mathbf{Z}$. A similar argument to Lemma~\ref{lemma:sharp-SES} shows the following: the natural map $W[F] \subset W \twoheadrightarrow \mathbf{G}_a$ lifts uniquely to an isomorphism 
\[ \Spec( \Z_{(p)} ) \times W[F] \xrightarrow{\sim} \Spec(\Z_{(p)}) \times \mathbf{G}_{a}^{\sharp}\]
of group schemes over $\mathbf{Z}_{(p)}$.
\end{variant}

\begin{theorem}\label{theorem:describe-HT}
The morphism $\eta: \Spf(\Z_p) \rightarrow \WCart^{\mathrm{HT}}$ of Construction \ref{construction:fiber-at-eta} extends to an isomorphism of stacks
$\Spf(\Z_p) \times B \mathbf{G}_{m}^{\sharp} \simeq \WCart^{\mathrm{HT}}$; here
$B \mathbf{G}_{m}^{\sharp}$ denotes the classifying stack for the group scheme $\mathbf{G}_{m}^{\sharp}$
(with respect to the fpqc topology).
\end{theorem}

\begin{proof}
Lemma \ref{lemma:sharp-SES} supplies an isomorphism $\Aut( \eta) \simeq \Spf(\Z_p) \times \mathbf{G}_{m}^{\sharp}$.
It will therefore suffice to show that the map $\eta: \Spf(\Z_p) \rightarrow \WCart^{\mathrm{HT}}$ is locally surjective with respect to the fpqc topology.
Let $R$ be a commutative ring and suppose we are given a morphism $f: \Spec(R) \rightarrow \WCart^{\mathrm{HT}}$,
corresponding to a Cartier-Witt divisor $\alpha: I \rightarrow W(R)$ of $R$ which factors through
the submodule $VW(R)$. We wish to show that, after passing to a faithfully flat extension of $R$, the morphism $f$ factors through $\eta$. 
By virtue of Proposition \ref{proposition:presentation-as-quotient}, we may assume without loss of generality that
$I = W(R)$, so that $\alpha$ is given by multiplication by the element $\alpha(1) \in VW(R)$. Write $\alpha(1) = V(u)$,
where $u \in W(R)$ is a unit. Since the Witt vector Frobenius $F: W^{\times} \rightarrow W^{\times}$ is faithfully flat (Corollary \ref{corollary:W-F-flat}),
we can assume (after passing to a faithfully flat extension of $R$) that $u = F(u')$ for some $u' \in W(R)^{\times}$. Then $\alpha(1) = V(u) = u' \cdot V(1)$,
so that $(W(R), \alpha)$ is isomorphic to the Cartier-Witt divisor $(W(R), V(1))$ in the category $\WCart(R)$.
\end{proof}

\begin{remark}
The homomorphism of group schemes $\mathbf{G}_{m}^{\sharp} \rightarrow \mathbf{G}_{m}$ determines a map of classifying
stacks $B \mathbf{G}_{m}^{\sharp} \rightarrow B \mathbf{G}_{m}$. Composing with the isomorphism of Theorem \ref{theorem:describe-HT},
we recover the comparison map $\mu^{\mathrm{HT}}: \WCart^{\mathrm{HT}} \rightarrow B \mathbf{G}_{m}$ of
Remark \ref{remark:comparison-W-times-F}.
\end{remark}

\begin{example}\label{example:torsor-over-Hodge-Tate}
Let $(A,I)$ be a prism and let $\overline{A}$ denote the quotient $A/I$, so that Remark \ref{remark:HT-point-of-prismatic-stack} determines a morphism
$\rho_{A}^{\mathrm{HT}}: \Spf( \overline{A} ) \rightarrow \WCart^{\mathrm{HT}}$. By virtue of Theorem \ref{theorem:describe-HT}, we can identify $\rho$ with a
$W^{\times}[F]$-torsor $\mathcal{P}$ over the formal scheme $\Spf( \overline{A} )$. Note that the quotient map $A \twoheadrightarrow \overline{A}$
lifts uniquely to a $\delta$-homomorphism $\psi: A \rightarrow W( \overline{A} )$. Unwinding the definitions, we see that
trivializations of the torsor $\mathcal{P}$ can be identified with $A$-linear maps $\xi: I \rightarrow W( \overline{A})$
for which the diagram of $A$-modules $$ \xymatrix@R=50pt@C=50pt{ I \ar[r]^-{ \xi } \ar[d] & W( \overline{A} ) \ar[d]^{ V(1)  } \\
A \ar[r]^-{\psi} & W( \overline{A} ) }$$
is commutative.
\end{example}

For later use, we note the following consequence of Theorem \ref{theorem:describe-HT}:

\begin{corollary}\label{corollary:automorphism-of-CW}
Let $R$ be a commutative ring and let $(I, \alpha)$ be a Cartier-Witt divisor for $R$. Then
the automorphism group $\Aut( I, \alpha)$ is commutative and is annihilated by a power of $p$.
\end{corollary}

\begin{proof}
By definition, $\Aut(I, \alpha)$ is a subgroup of the automorphism group of $I$ as a $W(R)$-module, which is isomorphic to the abelian group $W(R)^{\times}$
(since $I$ is an invertible $W(R)$-module).Let $J \subseteq R$ denote the ideal generated by $p$ together with the image of the map $I \xrightarrow{\alpha} W(R) \twoheadrightarrow R$. For each $n \geq 0$, let $(I_{n}, \alpha_n)$ denote the Cartier-Witt divisor for $R/J^n$ obtained from $(I,\alpha)$ by extension of scalars.
Note that the ideal $J$ is nilpotent, so that we have $\Aut(I,\alpha) \simeq \Aut(I_n, \alpha_n)$ for $n \gg 0$. We will complete the proof by showing that
each of the automorphism groups $\Aut( I_n, \alpha_n)$ is annihilated by $p^n$. The proof proceeds by induction on $n$. A standard deformation-theoretic argument shows that, for $n > 1$, the kernel of the restriction map $\Aut( I_{n}, \alpha_n) \rightarrow \Aut( I_{n-1}, \alpha_{n-1} )$ admits the structure of an $(R/J)$-module, and is therefore annihilated by $p$. It will therefore suffice to treat the case $n=1$. Replacing $R$ by the quotient ring $R/J$, we may assume
that $R$ is an $\F_p$-algebra and that the Cartier-Witt divisor $(I,\alpha)$ belongs to the Hodge-Tate divisor $\WCart^{\mathrm{HT}}(R)$.
In this case, Theorem \ref{theorem:describe-HT} supplies an isomorphism $\Aut(I,\alpha) \simeq \mathbf{G}_{m}^{\sharp}(R)$.
We conclude by observing that the $\F_p$-group scheme $\mathbf{G}_{m}^{\sharp} \times \Spec(\F_p)$ is annihilated by the Frobenius
(since it is isomorphic to the kernel of the Frobenius endomorphism of $W^{\times} \times \Spec(\F_p)$; see Lemma \ref{lemma:sharp-SES}),
and is therefore also annihilated by multiplication by $p$.
\end{proof}

\subsection{Complexes on the Hodge-Tate Divisor}\label{subsection:HT-locus}

Let $\mathbf{G}_{m} = \Spec( \Z[t^{\pm 1}] )$ be the multiplicative group. Representations of the group scheme $\mathbf{G}_{m}$ are easy to describe:
endowing an abelian group $M$ with an algebraic action of $\mathbf{G}_{m}$ is equivalent to giving a direct sum decomposition $M \simeq \bigoplus_{n \in \Z} M(n)$.
In this case, the abelian group $M$ is equipped with a {\it weight endomorphism} $\Theta: M \rightarrow M$, characterized by the identity
$\Theta(x) = nx$ for $x \in M(n)$. Concretely, the morphism $\Theta$ can be realized as the composition
$$ M \xrightarrow{c} \mathcal{O}_{\mathbf{G}_{m}} \otimes M \xrightarrow{ \psi \otimes \id_M } \Z \otimes M \simeq M,$$
where $c: M \rightarrow \mathcal{O}_{\mathbf{G}_{m}} \otimes M$ describes the action of $\mathbf{G}_{m}$ on $M$
and $\psi: \mathcal{O}_{\mathbf{G}_{m}} \rightarrow \Z$ is the homomorphism of abelian groups given by $f(t) \mapsto \frac{ \partial f(t)}{\partial t}|_{t=1}$.

Let $\mathbf{G}_{m}^{\sharp}$ be the group scheme described in Notation \ref{notation:G_m-sharp}.
Note that, if $M$ is an abelian group equipped with an algebraic action of $\mathbf{G}^{\sharp}_{m}$, then the rational vector space
$M_{\Q} = \Q \otimes M$ inherits action of $\mathbf{G}_{m}$ (since the map $\mathbf{G}_m^{\sharp} \rightarrow \mathbf{G}_{m}$ becomes
an isomorphism after extending scalars to $\Q$), which we can identify with a grading $M_{\Q} \simeq \bigoplus_{n \in \Q} M_{\Q}(n)$.
Beware that this grading generally cannot be obtained from a grading of the abelian group $M$ itself. However, the weight endomorphism
of the rational vector space $M_{\Q}$ {\em does} always arise from an endomorphism $\Theta$ of the underlying abelian group $M$, given by the composition
$$ M \xrightarrow{c} \mathcal{O}^{\sharp}_{\mathbf{G}_{m}} \otimes M \xrightarrow{ \psi \otimes \id_M } \Z \otimes M \simeq M,$$
where $\psi$ is again defined by the construction $f(t) \mapsto \frac{ \partial f(t)}{\partial t}|_{t=1}$. Moreover,
the action of $\mathbf{G}_{m}^{\sharp}$ on $M$ can be completely recovered from the endomorphism $\Theta$. In this section,
we prove a version of this result in the $p$-complete setting (Theorem \ref{theorem:compute-with-HT}), where representations of $\mathbf{G}_{m}^{\sharp}$ can be identified with quasi-coherent sheaves on the Hodge-Tate divisor $\WCart^{\mathrm{HT}}$.

\begin{definition}
Let $\calD( \WCart^{\mathrm{HT}} )$ denote the $\infty$-category of quasi-coherent complexes on the Hodge-Tate divisor $\WCart^{\mathrm{HT}} \subseteq \WCart$. 
Formally, $\calD( \WCart^{\mathrm{HT}} )$ is defined as the inverse limit
$$\varprojlim_{ \Spec(R) \rightarrow \WCart^{\mathrm{HT}} } \calD(R),$$
indexed by the category of points of the stack $\WCart^{\mathrm{HT}}$.
\end{definition}

\begin{example}\label{example:BK-twist-on-HT}
Let $\mathscr{E}$ be a quasi-coherent complex on the Cartier-Witt stack $\WCart$. Then $\mathscr{E}$ determines a quasi-coherent complex $\mathscr{E}|_{ \WCart^{\mathrm{HT}} }$ on the Hodge-Tate divisor $\WCart^{\mathrm{HT}}$. 

For every integer $n$, we let $\calO_{ \WCart^{\mathrm{HT}}}\{n\}$ denote the restriction of $\calO_{\WCart}\{n\}$ to the Hodge-Tate divisor
$\WCart^{\mathrm{HT}}$. By virtue of Remark \ref{remark:twist-modulo-I}, we have canonical isomorphisms
$\calO_{ \WCart^{\mathrm{HT}}}\{n\} \simeq \mathcal{I}^{n}|_{ \WCart^{\mathrm{HT}} }$, where $\mathcal{I}$ denotes the Hodge-Tate ideal
sheaf of Example \ref{example:HTI}.
\end{example}

\begin{remark}\label{remark:flashcard}
For every prism $(A,I)$, let $\rho_{A}^{\mathrm{HT}}: \Spf(A/I) \rightarrow \WCart^{\mathrm{HT}}$ be the map described in Remark \ref{remark:HT-point-of-prismatic-stack}.
Pullback along $\rho_A^{\mathrm{HT}}$ determines a functor from $\calD( \WCart^{\mathrm{HT}} )$ to the $p$-complete derived $\infty$-category $\widehat{\calD}(A/I)$. Arguing as in Proposition \ref{proposition:DWCart-prism-description}, we see that these maps furnish an equivalence of $\infty$-categories
$\calD( \WCart^{\mathrm{HT}} ) \rightarrow \varprojlim_{(A,I)} \widehat{\calD}(A/I)$, where the limit is taken over the category of bounded prisms (or the larger category of all prisms, or the smaller category of transversal prisms). Alternatively, we can describe $\calD( \WCart^{\mathrm{HT} })$ as the totalization of the cosimplicial
$\infty$-category $\widehat{\calD}( A^{\bullet} / I^{\bullet} )$, where $(A^{\bullet}, I^{\bullet} )$ is the cosimplicial prism of Notation \ref{notation:simplicial-prism}
(see Lemma \ref{lemma:tot-description}).
\end{remark}

\begin{construction}[The Sen Operator]\label{construction:fiber-at-eta2}
Let $\Z[ \epsilon ] / (\epsilon^2)$ denote the ring of dual numbers and let $R$ be a $\Z[ \epsilon ] / (\epsilon^2)$-algebra.
Note that the Teichm\"{u}ller representative $[\epsilon] \in W(R)$ is annihilated by the Witt vector Frobenius, so that multiplication
by $1 + [ \epsilon ]$ acts by the identity on the ideal $VW(R) \subseteq W(R)$. Let $(I, \alpha)$ be a Cartier-Witt divisor for $R$. Suppose that $(I, \alpha)$ determines an $R$-valued point of the Hodge-Tate divisor $\WCart^{\mathrm{HT}}$: that is, the map $\alpha: I \rightarrow W(R)$ factors through the ideal $VW(R) \subseteq W(R)$. It follows that multiplication by $1 + [\epsilon]$ determines an automorphism of the Cartier-Witt divisor $(I,\alpha)$. This construction depends functorially on
$R$ (and on the Cartier-Witt divisor $(I, \alpha)$), and can therefore be regarded as an automorphism $\gamma$ of the projection map
$$ \pi: \WCart^{\mathrm{HT}} \times \Spec( \Z[ \epsilon ] / (\epsilon^2) ) \rightarrow \WCart^{\mathrm{HT}}.$$

For every quasi-coherent complex $\mathscr{E} \in \QCoh( \WCart^{\mathrm{HT}} )$, $\gamma$ induces an automorphism of the pullback
$\pi^{\ast} \mathscr{E}$, which we can identify with a morphism
$$ u: \mathscr{E} \rightarrow \Z[ \epsilon ] / (\epsilon^2) \otimes_{\Z} \mathscr{E}$$
in the $\infty$-category $\calD( \WCart^{\mathrm{HT} } )$. This morphism reduces to the identity modulo $\epsilon$, and can therefore be written as
$\id + \epsilon \Theta_{\mathscr{E}}$ for some endomorphism $\Theta_{\mathscr{E}}: \mathscr{E} \rightarrow \mathscr{E}$, which we will
refer to as the {\it Sen operator} on the complex $\mathscr{E}$.
\end{construction}

\begin{remark}[The Leibniz Rule]\label{remark:leibniz-for-Sen}
Let $\mathscr{E}$ and $\mathscr{E}'$ be quasi-coherent complexes on $\WCart^{\mathrm{HT}}$, and let $\mathscr{E} \otimes \mathscr{E}'$ denote their (derived) tensor product. Then the Sen operator $\Theta_{ \mathscr{E} \otimes \mathscr{E}'}$ can be identified with the sum
$(\Theta_{\mathscr{E}} \otimes \id_{ \mathscr{E}'}) + ( \id_{\mathscr{E}} \otimes \Theta_{\mathscr{E}'} )$. In particular,
when $\mathscr{E} = \calO_{ \WCart^{\mathrm{HT}} }$ is the structure sheaf, the Sen operator $\Theta_{\mathscr{E}}$ vanishes.
\end{remark}

\begin{example}\label{example:Sen-operator-on-I}
Let $\mathscr{I} \in \QCoh( \WCart )$ be the Hodge-Tate ideal sheaf (Example \ref{example:HTI}). To every commutative $\Z[ \epsilon ] / (\epsilon^2)$-algebra $R$ and every Cartier-Witt divisor $(I, \alpha)$ for $R$, the sheaf $\mathscr{I}$ associates the $R$-module $R \otimes_{ W(R)} I$ obtained from $I$
by extending scalars along the restriction map $W(R) \twoheadrightarrow R$. Since this restriction map carries $1 + [ \epsilon ] \in W(R)$ to
the element $1 + \epsilon \in R$, it follows that the Sen operator $\Theta_{ \mathscr{I}|_{ \WCart^{\mathrm{HT}} } }$ is equal to the identity map.
Combining this observation with Remark \ref{remark:leibniz-for-Sen} and Example \ref{example:BK-twist-on-HT}, we see that
the Sen operator on $\calO_{ \WCart^{\mathrm{HT}}}\{n\}$ is given by multiplication by $n$ (for every integer $n$).
\end{example}

\begin{notation}\label{notation:Sen-operator-second}
Let $\eta: \Spf(\Z_p) \rightarrow \WCart^{\mathrm{HT}}$ be the faithfully flat surjection of Construction \ref{construction:fiber-at-eta},
corresponding to the Cartier-Witt divisor $W(\Z) \xrightarrow{V(1)} W(\Z)$. For every quasi-coherent complex $\mathscr{E} \in \calD( \WCart^{\mathrm{HT}} )$, we let $\mathscr{E}_{\eta}$ denote the pullback $\eta^{\ast}(\mathscr{E})$, which we regard as an object of the $p$-complete derived $\infty$-category $\widehat{\calD}(\Z_p) \simeq \calD( \Spf(\Z_p) )$. Note that the Sen operator $\Theta_{\mathscr{E}}: \mathscr{E} \rightarrow \mathscr{E}$ of Construction \ref{construction:fiber-at-eta2} induces an endomorphism of $\mathscr{E}_{\eta}$, which (by slight abuse of notation) we will also denote by $\Theta_{\mathscr{E}}$ and refer to
as the {\it Sen operator} of $\mathscr{E}$.
\end{notation}

Our goal in this section is to show that a quasi-coherent complex $\mathscr{E}$ on $\WCart^{\mathrm{HT}}$ is completely
determined by the fiber $\mathscr{E}_{\eta}$, together with the Sen operator $\Theta_{\mathscr{E}}: \mathscr{E}_{\eta} \rightarrow \mathscr{E}_{\eta}$
of Notation \ref{notation:Sen-operator-second}. More precisely, we have the following:

\begin{theorem}\label{theorem:compute-with-HT}
The functor
$$ \calD( \WCart^{\mathrm{HT}} ) \rightarrow \calD( \Z[ \Theta ] ) \quad \quad \mathscr{E} \mapsto (\mathscr{E}_{\eta}, \Theta_{\mathscr{E}})$$
is fully faithful. Moreover, its essential image consists of those objects $M \in \calD( \Z[ \Theta] )$ which satisfy the following pair of conditions:
\begin{itemize}
\item[$(a)$] The complex $M$ is $p$-complete.
\item[$(b)$] The action of $\Theta^{p} - \Theta$ on the cohomology $\mathrm{H}^{\ast}( \F_p \otimes^{L} M)$ is locally nilpotent.
\end{itemize}
\end{theorem}

The proof of Theorem \ref{theorem:compute-with-HT} will require some preliminaries.

\begin{example}\label{example:regular-representation}
Let $\eta: \Spf(\Z_p) \rightarrow \WCart^{\mathrm{HT}}$ be the faithfully flat surjection of Construction \ref{construction:fiber-at-eta}.
Let $\mathscr{E} = \eta_{\ast} \calO_{ \Spf(\Z_p)} \in \calD(\WCart^{\mathrm{HT}})$ denote the direct image of the structure sheaf of
$\Spf(\Z_p)$. Theorem \ref{theorem:describe-HT} then supplies an identification of $\eta^{\ast} \mathscr{E} \simeq \widehat{\calO}_{ \mathbf{G}_{m}}^{\sharp}$, where
$\widehat{\calO}_{ \mathbf{G}_{m}}^{\sharp}$ denotes the $p$-completion of the coordinate ring $\calO_{\mathbf{G}_{m}}^{\sharp}$ for the affine
group scheme $\mathbf{G}_{m}^{\sharp}$.

Note that the automorphism $\gamma$ of Construction \ref{construction:fiber-at-eta2} restricts to an automorphism of the composite map
$\Spf( \Z_p[ \epsilon ] / (\epsilon^2) ) \rightarrow \Spf( \Z_p) \xrightarrow{\eta} \WCart^{\mathrm{HT}}$. Under the identification $\Aut(\eta) \simeq \mathbf{G}_{m}^{\sharp}$ supplied by Theorem \ref{theorem:describe-HT},
this automorphism corresponds to the unit  $1 + \epsilon \in \mathbf{G}_m^{\sharp}( \Z_p[ \epsilon ] / (\epsilon^2) ) \subseteq \mathbf{G}_m( \Z_p[ \epsilon ] /  (\epsilon^2) )$.
It follows that, for each element $f(t) \in \widehat{\calO}_{\mathbf{G}_m}^{\sharp}$, we have an identity
$f( (1+\epsilon) t ) = f(t) + \epsilon \Theta_{ \mathscr{E}}( f(t) )$. In other words, the Sen operator on $\mathscr{E}_{\eta} \simeq \widehat{\calO}_{\mathbf{G}_m}^{\sharp}$ is given by the differential operator
$f(t) \mapsto t \frac{ \partial }{\partial t} f(t)$.
\end{example}

Let $u: \calO_{\WCart^{\mathrm{HT}}} \rightarrow \eta_{\ast} \calO_{ \Spf(\Z_p) }$ denote the unit map. Note that,
since the Sen operator vanishes on $\calO_{ \WCart^{\mathrm{HT}} }$, the morphism $u$ factors through the fiber
$$ (\eta_{\ast} \calO_{ \Spf(\Z_p) })^{\Theta = 0} = \fib( \Theta_{\eta_{\ast} \calO_{\Spf(\Z_p)} }: \eta_{\ast} \calO_{ \Spf(\Z_p) } \rightarrow \eta_{\ast} \calO_{ \Spf(\Z_p) }).$$

\begin{proposition}\label{proposition:O-inside-regular}
The unit map $\calO_{\WCart^{\mathrm{HT}}} \rightarrow (\eta_{\ast} \calO_{ \Spf(\Z_p) })^{\Theta = 0}$ is an isomorphism in 
the $\infty$-category $\calD(\WCart^{\mathrm{HT}} )$. 
\end{proposition}

\begin{proof}
Since the morphism $\eta: \Spf(\Z_p) \rightarrow \WCart^{\mathrm{HT}}$ of Construction \ref{construction:fiber-at-eta} is faithfully flat,
it will suffice to show that the induced map $\eta^{\ast}( \calO_{\WCart^{\mathrm{HT}}} ) \rightarrow \eta^{\ast}( (\eta_{\ast} \calO_{ \Spf(\Z_p) })^{\Theta = 0} )$
is an isomorphism in the $\infty$-category $\widehat{\calD}(\Z_p)$. By virtue of Example \ref{example:regular-representation}, this is equivalent to the exactness
of the sequence of $p$-complete abelian groups
$$ 0 \rightarrow \Z_p \rightarrow \widehat{\calO}_{\mathbf{G}_{m}}^{\sharp}
\xrightarrow{ t \frac{\partial}{\partial t} }  \widehat{\calO}_{\mathbf{G}_{m}}^{\sharp} \rightarrow 0.$$
Note that $\widehat{\calO}_{\mathbf{G}_{m}}^{\sharp}$ can be identified with the collection of all formal series
$f(t) = \sum_{n \geq 0} c_n \frac{ (t-1)^{n} }{n!}$ where the coefficients $c_{n}$ converge to zero in $\Z_p$. Under this identification,
the differential operator $t \frac{ \partial}{\partial t}$ is given by the construction
$$  \sum_{n \geq 0} c_n \frac{ (t-1)^{n} }{n!} \mapsto \sum_{ n \geq 0} (c_{n+1} + n c_n) \frac{ (t-1)^{n} }{n!},$$
from which we immediately deduce that $t \frac{ \partial}{\partial t}$ is surjective and that its kernel consists
of those series $\sum_{n \geq 0} c_n \frac{ (t-1)^{n} }{n!}$ where the coefficients $c_{n}$ vanish for $n > 0$.
\end{proof}

Let $\mathscr{E}$ be any quasi-coherent complex on $\WCart^{\mathrm{HT}}$. Tensoring the isomorphism of Proposition \ref{proposition:O-inside-regular}
with $\mathscr{E}$ and applying the projection formula, we obtain a fiber sequence
$$ \mathscr{E} \rightarrow \eta_{\ast} \mathscr{E}_{\eta} \rightarrow \eta_{\ast} \mathscr{E}_{\eta}.$$
Passing to global sections over $\WCart$, we obtain the following:

\begin{proposition}\label{proposition:easy-version}
For every object $\mathscr{E} \in \calD( \WCart^{\mathrm{HT}} )$, we have a canonical fiber sequence
$$ \RGamma( \WCart^{\mathrm{HT}}, \mathscr{E} ) \rightarrow \mathscr{E}_{\eta} \xrightarrow{ \Theta_{\mathscr{E} } } \mathscr{E}_{\eta},$$
where $\Theta_{\mathscr{E}}$ is the Sen operator of Notation \ref{notation:Sen-operator-second}.
\end{proposition}

\begin{example}\label{example:cohomology-of-HT}
Proposition \ref{proposition:easy-version} supplies canonical isomorphisms
$$ \mathrm{H}^{n}( \WCart^{\mathrm{HT}}, \calO_{ \WCart^{\mathrm{HT}} } ) \simeq \begin{cases} \Z_p & \text{ if } n=0,1 \\
0 & \text{ otherwise. } \end{cases}$$
When $n=1$, the inverse isomorphism $\Z_p \simeq \mathrm{H}^{1}( \WCart^{\mathrm{HT}}, \calO_{ \WCart^{\mathrm{HT}} } )$
can be described explicitly: it carries the generator $1 \in \Z_p$ to the cohomology class determined by the logarithm map
$$ \log: \Spf(\Z_p) \times \mathbf{G}_{m}^{\sharp} \rightarrow \Spf(\Z_p) \times \mathbf{G}_{a}.$$
\end{example}

\begin{corollary}\label{corollary:HT-computation}
The global sections functor $\RGamma( \WCart^{\mathrm{HT}}, \bullet): \calD( \WCart^{\mathrm{HT} } ) \rightarrow \widehat{\calD}(\Z_p)$
commutes with colimits.
\end{corollary}

\begin{corollary}\label{srong}
Let $n$ be an integer and let $\mathscr{E}$ be a quasi-coherent complex on $\WCart^{\mathrm{HT}}$. Then
$\mathscr{E}$ is isomorphic to the Breuil-Kisin twist $\calO_{\WCart^{\mathrm{HT}}} \{n\}$ if and only if the following pair of conditions is satisfied:
\begin{itemize}
\item[$(1)$] The fiber $\mathscr{E}_{\eta}$ is isomorphic to $\Z_p$.
\item[$(2)$] The Sen operator $\Theta_{\mathscr{E}}: \mathscr{E}_{\eta} \rightarrow \mathscr{E}_{\eta}$
is given by multiplication by $n$.
\end{itemize}
\end{corollary}

\begin{proof}
The necessity of conditions $(1)$ and $(2)$ follows from Example \ref{example:Sen-operator-on-I}. To prove the converse, we can replace
$\mathscr{E}$ by the twist $\mathscr{E}(-n)$ and thereby reduce to the case $n=0$ (see Remark \ref{remark:leibniz-for-Sen}).
In this case, assumption $(2)$ asserts that the Sen operator $\Theta_{\mathscr{E}}: \mathscr{E}_{\eta} \rightarrow \mathscr{E}_{\eta}$ vanishes.
Choose an isomorphism $u: \Z_p \simeq \mathscr{E}_{\eta}$. Applying Proposition \ref{proposition:easy-version}, we can realize $u$ as the
pullback along $\eta$ of a morphism of  quasi-coherent complexes $\calO_{ \WCart^{\mathrm{HT}} } \rightarrow \calE$, which is also an isomorphism (since $\eta$ is faithfully flat).
\end{proof}

\begin{proposition}\label{proposition:twists-generate}
The $\infty$-category $\calD( \WCart^{\mathrm{HT}} )$ is generated, under shifts and colimits,
by the Breuil-Kisin twists $\calO_{ \WCart^{\mathrm{HT}} } \{n \}$ for $n \geq 0$.
\end{proposition}

\begin{proof}
Let $\calC \subseteq \calD( \WCart^{\mathrm{HT} })$ be the full subcategory generated under
shifts and colimits by the Breuil-Kisin twists $\calO_{ \WCart^{\mathrm{HT}} } \{n \}$ for $n \geq 0$.
We wish to show that every object $\mathscr{E} \in  \calD( \WCart^{\mathrm{HT} })$ belongs to $\calC$.
Let $\eta: \Spf(\Z_p) \rightarrow \WCart^{\mathrm{HT}}$ be the faithfully flat morphism of Construction \ref{construction:fiber-at-eta}, so that
Proposition \ref{proposition:O-inside-regular} supplies a fiber sequence
$$\mathscr{E} \rightarrow \eta_{\ast}( \mathscr{E}_{\eta} ) \rightarrow \eta_{\ast}( \mathscr{E}_{\eta} ).$$
Consequently, to show that $\mathscr{E}$ belongs to $\calC$, it will suffice to show that $\calC$ contains
the pushforward $\eta_{\ast}( \mathscr{E}_0 )$ for each object $\mathscr{E}_0 \in \calD( \Spf(\Z_p) )$. Since $\calD( \Spf(\Z_p) )$ is generated
(under shifts and colimits) by the structure sheaf of $\Spf(\Z_p)$ it will suffice to show that $\mathscr{F} = \eta_{\ast}( \calO_{ \Spf(\Z_p)})$ belongs to $\calC$.

Using Theorem \ref{theorem:describe-HT}, we can identify $\calD( \WCart^{\mathrm{HT} } )$ with the $\infty$-category of $\calO_{\mathbf{G}_{m}}^{\sharp}$-comodule objects of $\widehat{\calD}(\Z_p)$. Under this identification, $\eta_{\ast}( \calO_{ \Spf(\Z_p) )}$ corresponds to the $p$-complete
regular representation $\widehat{\calO}_{\mathbf{G}_{m}}^{\sharp}$ of $\calO_{\mathbf{G}_{m}}^{\sharp}$. For each $n \geq 0$, let $V_{\leq n}$ denote the
$\Z_p$-submodule of $\widehat{\calO}_{\mathbf{G}_{m}}^{\sharp}$ generated by the divided powers $\frac{ (t-1)^{m} }{m!}$ for $m \leq n$. Note that the coaction
$c: \widehat{\calO}_{\mathbf{G}_{m}}^{\sharp} \rightarrow \widehat{\calO}_{\mathbf{G}_{m}}^{\sharp} \widehat{\otimes} \widehat{\calO}_{\mathbf{G}_{m}}^{\sharp}$ is determined by the formula $c(t) = t \otimes t$, and therefore satisfies
\begin{eqnarray*}
c( \frac{ (t-1)^{n} }{n!} ) & = & \frac{ ( t \otimes t - 1 \otimes 1 )^{n}}{n!} \\
& = & \frac{ (t \otimes (t-1)  + (t-1) \otimes 1 )^{n} }{n!} \\
& = & \sum_{a+b = n} \frac{ t^{a} (t-1)^{b} }{b!} \otimes \frac{ (t-1)^{a} }{a!} \\
& \in & t^{n} \otimes \frac{ (t-1)^{n}}{n!} + (\calO_{\mathbf{G}_{m}}^{\sharp} \widehat{\otimes} V_{\leq n-1} ).
\end{eqnarray*}
It follows that each $V_{\leq n}$ inherits the structure of a $\calO_{\mathbf{G}_{m}}^{\sharp}$-comodule object of $\widehat{\calD}(\Z_p)$, and
therefore corresponds to a quasi-coherent complex $\mathscr{F}_{\leq n} \in \calD( \WCart^{\mathrm{HT} })$. Note that the identity 
$$t \frac{ \partial}{\partial t} \frac{(t-1)^{n}}{n!} = \frac{(t-1)^{n}}{(n-1)!} + \frac{(t-1)^{n-1}}{(n-1)!} \equiv n \frac{(t-1)^{n}}{n!} \pmod{ V_{n-1} }$$
guarantees that the Sen operator acts by multiplication by $n$ on each quotient $(\mathscr{F}_{\leq n} / \mathscr{F}_{\leq n-1})_{\eta}$.
Invoking Corollary \ref{srong}, we obtain fiber sequences
$$ \mathscr{F}_{\leq n-1} \rightarrow \mathscr{F}_{\leq n} \rightarrow \calO_{ \WCart^{\mathrm{HT}}}\{n\}.$$
It follows by induction on $n$ that each $\mathscr{F}_{\leq n}$ belongs to the category $\calC$, so that
$\mathscr{F} \simeq \varinjlim_{n} \mathscr{F}_{\leq n}$ also belongs to $\calC$.
\end{proof}

\begin{corollary}\label{corollary:HT-ideals-generate}
The $\infty$-category $\calD( \WCart )$ is generated, under shifts and colimits, by the invertible sheaves $\calI^{n}$ for $n \in \Z$
(here $\calI$ denotes the Hodge-Tate ideal sheaf of Example \ref{example:HTI}).
\end{corollary}

\begin{proof}
Let $\mathscr{E}$ be a quasi-coherent complex on the Cartier-Witt stack satisfying $\RHom( \calI^{n}, \mathscr{E} )
\simeq \RGamma( \WCart, \calI^{-n} \mathscr{E} )$ vanishes for every integer $n$. Using the fiber sequence
$$ \RGamma( \WCart, \calI^{1-n} \mathscr{E} ) \rightarrow \RGamma( \WCart, \calI^{-n} \mathscr{E} ) 
\rightarrow \RGamma( \WCart^{\mathrm{HT}}, (\calI^{-n} \mathscr{E})|_{ \WCart^{\mathrm{HT}} } ),$$
we deduce that each of the complexes $$\RGamma( \WCart^{\mathrm{HT}}, (\calI^{-n} \mathscr{E})|_{ \WCart^{\mathrm{HT}} }  )
\simeq \RHom( \calI^{n}|_{ \WCart^{\mathrm{HT}}}, \mathscr{E}|_{ \WCart^{\mathrm{HT}}} )$$ vanishes.
It follows from Proposition \ref{proposition:twists-generate} (and Example \ref{example:BK-twist-on-HT}) that the restriction $\mathscr{E}|_{ \WCart^{\mathrm{HT}}}$ vanishes, so that $\mathscr{E} \simeq 0$.
\end{proof}

\begin{proof}[Proof of Theorem \ref{theorem:compute-with-HT}]
We first show that the functor
$$ \calD( \WCart^{\mathrm{HT}} ) \rightarrow \calD( \Z[ \Theta ] ) \quad \quad \mathscr{E} \mapsto (\mathscr{E}_{\eta}, \Theta_{\mathscr{E}})$$
is fully faithful. Let $\mathscr{E}$ and $\mathscr{F}$ be quasi-coherent complexes on $\WCart^{\mathrm{HT} }$; we
wish to show that the natural map
$$ \Hom_{  \calD( \WCart^{\mathrm{HT}} )}( \mathscr{E}, \mathscr{F} ) \rightarrow \Hom_{ \calD( \Z[ \Theta ] ) }( \mathscr{E}_{\eta}, \mathscr{F}_{\eta} )$$
is a homotopy equivalence. By virtue of Proposition \ref{proposition:twists-generate}, we may assume that $\mathscr{E} = \mathcal{O}_{ \WCart^{\mathrm{HT} } }\{n\}$ for some integer $n$.
Replacing $\mathscr{F}$ by the Breuil-Kisin twist $\mathscr{F}\{-n\}$, we can reduce to the case $n = 0$. In this case, the desired result is a reformulation of Proposition \ref{proposition:easy-version}.

We next claim that, for every object $\mathscr{E} \in \calD( \WCart^{\mathrm{HT}} )$, the fiber $\mathscr{E}_{\eta}$ satisfies condition $(b)$ 
of Theorem \ref{theorem:compute-with-HT} (note that condition $(a)$ is automatic from the definition).
By virtue of Proposition \ref{proposition:twists-generate}, we may again assume without loss of generality that $\mathscr{E} = \mathcal{O}_{ \WCart^{\mathrm{HT} } }\{n\}$ for some integer $n$, so that the desired result follows from Example \ref{example:Sen-operator-on-I} together with the congruence $n^{p} \equiv n \pmod{p}$.

Let $\calC \subseteq \calD( \Z[ \Theta ] )$ denote the full subcategory spanned by those objects satisfying conditions $(a)$ and $(b)$. To complete the proof, it will suffice to show that $\calC$ is generated by the objects $\calO_{\WCart^{\mathrm{HT}} }\{n\}_{\eta}$ for $n \in \Z$. Equivalently, we must show that for every nonzero object $M \in \calC$, there exists a nonzero morphism $\calO_{\WCart^{\mathrm{HT}} }\{n\}_{\eta}[m] \rightarrow M$, for some pair of integers $m$ and $n$. Replacing $M$ by the tensor product
$\F_p \otimes^{L} M$, we may assume without loss of generality that some cohomology group $\mathrm{H}^{-m}(M)$ contains a nonzero element which is annihilated by $\Theta^{p} - \Theta = \prod_{0 \leq n < p} (\Theta-n)$, and therefore also a nonzero element which is annihilated by $\Theta - n$ for some integer $n$. It then follows that there exists
a morphism from $\mathscr{O}_{\WCart^{\mathrm{HT}}}\{n\}_{\eta}[m] \simeq \Z_p[\Theta] / (\Theta - n )[m]$ to $M$ which is nonzero on cohomology in degree $-m$.
\end{proof}

\begin{remark}[The Cartier dual of $\mathbf{G}_m^\sharp$]
\label{rmk:LogGmsharp}
Theorem \ref{theorem:compute-with-HT} can be regarded as form of Cartier duality. Let $\mathbf{G}_a^\sharp$ be the group scheme from Variant~\ref{Gasharpdef}; it can be  identified with the Cartier dual of the formal additive group $\widehat{\mathbf{G}}_{a}$ by inspection. To describe the Cartier dual of $\mathbf{G}_m^\sharp$, we shall use the following:

\begin{lemma}
\label{LogPullbackDesc}
After extending scalars to $R = \Z / p^{k} \Z$ for any $k \geq 0$, there is a pullback diagram of flat group schemes $$ \xymatrix@R=50pt@C=50pt{ \mathbf{G}_{m}^{\sharp} \ar[r]^-{x \mapsto \log(t)} \ar[d]^-{t \mapsto t} & \mathbf{G}_{a}^{\sharp} \ar[d]^-{ t \mapsto \exp(px)} \\
\mathbf{G}_{m} \ar[r]^-{ t \mapsto t^{p} } & \mathbf{G}_{m}. }$$
In particular, there is a short exact sequence
\[ 0 \to \mu_p \xrightarrow{[\cdot]} \mathbf{G}_m^\sharp \xrightarrow{\log(-)} \mathbf{G}_a^\sharp \to 0\]
of flat group schemes over $R$. 
\end{lemma}
\begin{proof}
The existence of the commutative diagram is clear. To prove the pullback property, as all group schemes are flat over $R$, we may assume $R=\mathbf{Z}/p$. Now recall that we have natural identifications
\[ \mathbf{G}_a^\sharp = W[F] \quad \text{and} \quad \mathbf{G}_m^\sharp = W^\times[F].\]
Moreover, as we work in characteristic $p$, there is an isomorphism of group schemes
\[ \mu_p \times W[F] \simeq W^\times[F]\]
determined by $(a,x) \mapsto [a]+Vx$. Using these descriptions, it suffices to show that the composition
\[ W[F] \xrightarrow{x \mapsto 1+Vx} W^\times[F] \xrightarrow{\log(-)} W[F]\]
is an isomorphism. In fact, we claim the stronger statement that this composition identifies with the map $V-\mathrm{id}$ (which is clearly an isomorphism: $V$ is a topologically nilpotent endomorphism of $W[F]$ as we work in characteristic $p$). Thus, we must check that
\[ \log(1+Vx) = Vx-x\]
for $x \in W[F](S)$ for any $\mathbf{F}_p$-algebra $S$. Consider the power series expansion
\[ \log(1+Vx) = \sum_{n \geq 1} (-1)^{n-1} \frac{(Vx)^n}{n}.\]
Note that $Vx \cdot Vy = V(Fx \cdot x) = 0$ for $x \in W[F](S)$; this kills the terms indexed by $1 < n < p$. Moreover, by the existence of divided powers of $Vx$, the terms for $n > p$ vanish as $W[F](S) = \mathbf{G}_a^\sharp(S)$ is $p$-torsion. Thus, we obtain
\[ \log(1+Vx) = Vx + \frac{(Vx)^p}{p},\]
so it suffices to show that $\frac{(Vx)^p}{p} = -x$. Now recall that $\frac{ (Vx)^p}{p}$ is defined by observing that $1+Vx \in W^\times[F](S) = \mathbf{G}_m^\sharp(S)$. Concretely, since $F(1+Vx) = 1$, we can write $F(Vx) = p \cdot 0$, whence standard arguments with $\delta$-rings give a well-defined divided power:
\[ \frac{(Vx)^p}{p} := \frac{ F(Vx) - p\delta(Vx)}{p} := 0 - \delta(Vx) = -\delta(Vx).\]
We are therefore reduced to checking that $\delta(Vx) = x$ for $x \in W[F](S) \subset W(S)$. But for any $y \in W(S)$, one has the formula
\[ \delta(Vy) = y - (p-1)! p^{p-2}V(y^p),\]
which one checks by reducing to a universal case. Specializing to our context, it is then enough to show that $x^p = 0$ for $x \in W[F](S) \subset W(S)$. But $F(x) = 0$, so $x^p = -p\delta(x)$ in $W(S)$. As $W[F](S) \subset W(S)$ is $p$-torsion, it suffices to show that $\delta(x) \in W[F](S)$; this follows because $F\delta = \delta F$.
\end{proof}

Passing to Cartier duals from the pullback square in Lemma~\ref{LogPullbackDesc}, we obtain a pushout diagram of formal group schemes
$$ \xymatrix@R=50pt@C=50pt{ \mathbf{D}( \mathbf{G}_{m}^{\sharp} ) & \widehat{\mathbf{G}}_{a} \ar[l] \\
\Z \ar[u] & p\Z, \ar[l] \ar[u] }$$
over $R=\mathbf{Z}/p^k$, which induces an isomorphism of the Cartier dual $\mathbf{D}( \mathbf{G}_{m}^{\sharp} )$ with the formal scheme
$X$ obtained by completing the affine line $\mathbf{G}_{a} = \Spec( R[ \Theta ] )$ along the closed subset consisting of its $\F_p$-rational points
(that is, the vanishing locus of the element $\Theta^{p} - \Theta$). 
Cartier duality then identifies $\calD( \Spec(R) \times B\mathbf{G}_{m}^{\sharp})$ with the full subcategory of
$\calD( R[ \Theta] )$ spanned by those objects which are complete with respect to the ideal $(\Theta^{p} - \Theta)$.
We recover the statement of Theorem \ref{theorem:compute-with-HT} by taking the inverse limit over $k$. In \cite[Appendix B]{DrinfeldFormalGroup}, one can find a different perspective on the same calculation.
\end{remark}

\subsection{The Frobenius Endomorphism of \texorpdfstring{$\WCart$}{WCart}}\label{subsection:Frobenius-on-stack}

We now observe that the Cartier-Witt stack $\WCart$ is equipped with a lift of Frobenius.

\begin{construction}[The Frobenius Endomorphism of $\WCart$]\label{construction:Frobenius-on-stack}
Let $R$ be a commutative ring. Then the Witt vector Frobenius $W(R) \rightarrow W(R)$ determines a functor
$\Cart( W(R) ) \rightarrow \Cart( W(R) )$, which carries Cartier-Witt divisors of $R$ to Cartier-Witt divisors of $R$.
This construction depends functorially on $R$, and therefore determines an endomorphism of the Cartier-Witt stack
which we denote by $F: \WCart \rightarrow \WCart$ and refer to as the {\it Frobenius morphism}.
\end{construction}

\begin{remark}\label{remark:WCart-lift-Frobenius}
Let $R$ be an $\F_p$-algebra. Then the Witt vector Frobenius $W(R) \rightarrow W(R)$ 
can be identified with $W( \varphi_{R} )$, where $\varphi_{R}: R \rightarrow R$ is the Frobenius endomorphism of $R$.
It follows that the induced map $F: \WCart(R) \rightarrow \WCart(R)$ coincides (up to canonical isomorphism)
with the map $\WCart( \varphi_{R} )$. In other words, the Frobenius endomorphism $F: \WCart \rightarrow \WCart$
of Construction \ref{construction:Frobenius-on-stack} restricts to the usual Frobenius endomorphism of the
$\F_p$-stack $\WCart \times \Spec(\F_p)$.
\end{remark}

\begin{remark}\label{remark:rho-compatible-with-Frobenius}
Let $(A,I)$ be a prism and let $\rho_{A}: \Spf(A) \rightarrow \WCart$ be the morphism of Construction \ref{construction:point-of-prismatic-stack}.
Then the diagram
$$ \xymatrix@R=50pt@C=50pt{ \Spf(A) \ar[r]^-{ \rho_{A} } \ar[d]^{ \varphi} & \WCart \ar[d]^{ \varphi } \\
\Spf(A) \ar[r]^-{ \rho_{A} } & \WCart }$$
commutes up to canonical isomorphism. Here the left vertical map is the Frobenius lift determined by the $\delta$-structure on the ring $A$, and the right vertical map
is the Frobenius of Construction \ref{construction:Frobenius-on-stack}.
\end{remark}

\begin{remark}
Let $(A^{\bullet}, I^{\bullet} )$ be the cosimplicial prism of Notation \ref{notation:simplicial-prism}, so that $\WCart$ can be identified with the geometric realization of the
simplicial formal scheme $\Spf( A^{\bullet} )$. Then the Frobenius endomorphism $F: \WCart \rightarrow \WCart$ of Construction \ref{construction:Frobenius-on-stack}
is induced by the endomorphism of simplicial formal scheme $\Spf( A^{\bullet} ) \rightarrow \Spf( A^{\bullet} )$, given at each level by the Frobenius lift
$\varphi_{A^{n}}: A^{n} \rightarrow A^{n}$ supplied by the $\delta$-structure on $A^{n}$.
\end{remark}

\begin{remark}\label{remark:Frobenius-pullback-on-WCart}
Let $F: \WCart \rightarrow \WCart$ be the Frobenius morphism of Construction \ref{construction:Frobenius-on-stack}. Then $F$ determines a pullback functor
$F^{\ast}: \calD( \WCart ) \rightarrow \calD( \WCart )$. For every integer $n$, Remark \ref{remark:Frobenius-on-twist} supplies a canonical isomorphism 
$F^{\ast} \calO_{\WCart}\{n\} \simeq \calI^{\otimes -n} \calO_{\WCart}\{n\}$ of invertible sheaves on $\WCart$.
\end{remark}

Beware that the Frobenius morphism $F: \WCart \rightarrow \WCart$ of Construction \ref{construction:Frobenius-on-stack} does not preserve the Hodge-Tate divisor
$\WCart^{\mathrm{HT}}$. Instead we have the following (see Lemma~4.5.3 of \cite{drinfeld-prismatic}):

\begin{proposition}\label{proposition:Frobenius-square}
The diagram of stacks
\begin{equation}
\begin{gathered}\label{equation:Frobenius-pushout-square} \xymatrix@R=50pt@C=50pt{ \WCart^{\mathrm{HT}} \ar[r] \ar[d] & \WCart \ar[d]^{F} \\
\Spf( \Z_p) \ar[r]^-{ \rho_{\dR} } & \WCart }
\end{gathered}
\end{equation}
commutes up to unique isomorphism. Here $\rho_{\dR}$ is the de Rham point of Example \ref{example:de-Rham-point}.
\end{proposition}

\begin{proof}
Let $R$ be a commutative ring and let $\alpha: I \rightarrow W(R)$ be a Cartier-Witt divisor for $\Spec(R)$ which determines an $R$-valued point of the Hodge-Tate divisor
$\WCart^{\mathrm{HT}}$. Write $F: W(R) \rightarrow W(R)$ for the Witt vector Frobenius and $V: W(R) \rightarrow W(R)$ for the Verschiebung, so that $\alpha$
factors uniquely as a composition $V \circ \beta$ for some $F$-semilinear map $\beta: I \rightarrow W(R)$. Let $F^{\ast}(I)$ denote the invertible $W(R)$-module
obtained from $I$ by extending scalars along $F$, so that $\beta$ induces a $W(R)$-linear morphism $\gamma: F^{\ast}(I) \rightarrow W(R)$.
Our assumption that $(I, \alpha)$ is a Cartier-Witt divisor guarantees that $\gamma$ is an isomorphism. It follows from the identity
$\alpha = V \circ \beta$ that the pullback $F^{\ast}(\alpha): F^{\ast}(I) \rightarrow W(R)$ coincides with the composition
$$ F^{\ast}(I) \xrightarrow{\gamma} W(R) \xrightarrow{V} W(R) \xrightarrow{F} W(R).$$
Consequently, $\gamma$ determines an isomorphism of Cartier-Witt divisors $( F^{\ast}(I), F^{\ast}(\alpha) ) \simeq ( W(R), p )$.
This construction depends functorially on $R$ and therefore determines an isomorphism of $F|_{ \WCart^{\mathrm{HT}}}$ with the composite map
$\WCart^{\mathrm{HT}} \rightarrow \Spf(\Z_p) \xrightarrow{ \rho_{\dR} } \WCart$. To prove the uniqueness of this isomorphism, it suffices to observe
that the automorphism group of $\rho_{\dR}$ is trivial (as a $\Spf(\Z_p)$-valued point of the Cartier-Witt stack), since $p$ is not a zero-divisor in the ring $W(\Z_p)$.
\end{proof}

It follows from Proposition \ref{proposition:Frobenius-square} that the Frobenius morphism $F: \WCart \rightarrow \WCart$ is not an isomorphism: it collapses
the Hodge-Tate divisor $\WCart^{\mathrm{HT}}$ to a point. Our goal for the remainder of this section is to show that this is essentially the {\em only} reason that
$F$ fails to be an isomorphism: that is, the diagram (\ref{equation:Frobenius-pushout-square}) behaves in certain respects like a pushout square.

\begin{theorem}\label{theorem:Frobenius-pushout-square}
Let $\mathscr{E}$ be a quasi-coherent complex on the Cartier-Witt stack $\WCart^{\mathrm{HT}}$. Then the diagram (\ref{equation:Frobenius-pushout-square})
determines a pullback square
$$ \xymatrix@R=50pt@C=50pt{ \RGamma( \WCart^{\mathrm{HT}}, (F^{\ast} \mathscr{E} )|_{ \WCart^{\mathrm{HT} }} ) & \RGamma( \WCart, F^{\ast} \mathscr{E} ) \ar[l] \\
\RGamma( \Spf(\Z_p), \rho_{\dR}^{\ast} \mathscr{E} ) \ar[u] & \RGamma( \WCart, \mathscr{E} ) \ar[u] \ar[l] }$$
in the $\infty$-category $\widehat{\calD}( \Z_p )$.
\end{theorem}

\begin{remark}\label{remark:Frobenius-exact-sequence}
Let $\mathscr{E}$ be a quasi-coherent complex on the Cartier-Witt stack $\WCart$, and let $F^{\ast} \mathscr{E}$ be its Frobenius pullback,
and let $(F^{\ast} \mathscr{E})_{\eta} \in \widehat{\calD}(\Z_p)$ be the fiber of $F^{\ast} \mathscr{E}$ at the point $\eta$ of Construction \ref{construction:fiber-at-eta}.
It follows from Proposition \ref{proposition:Frobenius-square} that we can identify $(F^{\ast} \mathscr{E})_{\eta}$ with the pullback
$\rho_{\dR}^{\ast} \mathscr{E}$, where $\rho_{\dR}$ is the de Rham point of Example \ref{example:de-Rham-point}.
Moreover, the Sen operator on $(F^{\ast} \mathscr{E} )_{\eta}$ is canonically trivial. Applying Proposition \ref{proposition:easy-version}, we obtain
canonical isomorphisms
\begin{eqnarray*}
\RGamma( \WCart^{\mathrm{HT}}, (F^{\ast} \mathscr{E} ))|_{ \WCart^{\mathrm{HT} }} & \simeq & \fib( \Theta: (F^{\ast} \mathscr{E})_{\eta} \rightarrow (F^{\ast} \mathscr{E})_{\eta} ) \\
& \simeq & (F^{\ast} \mathscr{E})_{\eta} \oplus (F^{\ast} \mathscr{E})_{\eta}[-1] \\
& \simeq & (\rho_{\dR}^{\ast} \mathscr{E}) \oplus ( \rho_{\dR}^{\ast} \mathscr{E})[-1].
\end{eqnarray*}
Consequently, the pullback diagram of Theorem \ref{theorem:Frobenius-pushout-square} determines a fiber sequence
$$ \RGamma( \WCart, \mathscr{E} ) \rightarrow \RGamma( \WCart, F^{\ast} \mathscr{E} ) \rightarrow \rho_{\dR}^{\ast} \mathscr{E}[-1]$$
in the $\infty$-category $\widehat{\calD}(\Z_p)$.
\end{remark}

\begin{corollary}\label{corollary:Frob-fully-faithful}
The diagram ( \ref{equation:Frobenius-pushout-square} ) determines a fully faithful functor
of $\infty$-categories
$$ \calD( \WCart ) \xrightarrow{ ( F^{\ast}, \rho_{\dR}^{\ast} )} \calD( \WCart ) \times_{ \calD( \WCart^{\mathrm{HT}} ) } \calD( \Spf(\Z_p) ).$$
\end{corollary}

\begin{proof}
Let $\calC$ denote the fiber product 
$$ \calD( \WCart ) \times_{ \calD( \WCart^{\mathrm{HT}} ) } \calD( \Spf(\Z_p) )$$
and let $\mathscr{E}$ and $\mathscr{E}'$ be quasi-coherent complexes on the Cartier-Witt stack; we wish to show that the functor
$F$ induces a homotopy equivalence of mapping spaces
$$\Hom_{\calD( \WCart )}( \mathscr{E}', \mathscr{E} ) \rightarrow \Hom_{\calC}( (F^{\ast}, \rho_{\dR}^{\ast})( \mathscr{E}'), (F^{\ast}, \rho_{\dR}^{\ast})( \mathscr{E}) ).$$
By virtue of Corollary \ref{corollary:HT-ideals-generate}, we may assume without loss of generality that $\mathscr{E}' = \calI^{n}$
is a power of the Hodge-Tate ideal sheaf of Example \ref{example:HTI}. Replacing $\mathscr{E}$ by the twist
$\calI^{-n} \mathscr{E}$, we can reduce to the case $n=0$, in which case the desired result follows from
Theorem \ref{theorem:Frobenius-pushout-square}.
\end{proof}

\begin{remark}
Let $\mathscr{E}$ be a quasi-coherent complex on the Cartier-Witt stack $\WCart$, so that
the Sen operator vanishes on the fiber $(F^{\ast} \mathscr{E})_{\eta}$ (Remark \ref{remark:Frobenius-exact-sequence}).
We can regard Corollary \ref{corollary:Frob-fully-faithful} (together with Theorem \ref{theorem:compute-with-HT})
as providing a partial converse: the complex $\mathscr{E}$ can be recovered (up to canonical isomorphism)
by the Frobenius pullback $F^{\ast} \mathscr{E}$ together with a trivialization of the Sen operator $\Theta_{ F^{\ast} \mathscr{E}}$ 
(as an endomorphism of the complex $(F^{\ast} \mathscr{E})_{\eta} \in \widehat{\calD}(\Z_p)$). 
\end{remark}

\begin{warning}
The fully faithful functor of Corollary \ref{corollary:Frob-fully-faithful} is not an equivalence of $\infty$-categories.
\end{warning}

Our first step is to reduce Theorem \ref{theorem:Frobenius-pushout-square} to a more concrete statement describing the Hodge-Tate divisor $\WCart^{\mathrm{HT}}$.

\begin{notation}
Let $F^{-1} \WCart^{\mathrm{HT}}$ denote the closed substack of $\WCart$ given by the inverse image of the Hodge-Tate divisor along the Frobenius map
$F: \WCart \rightarrow \WCart$. Concretely, for every commutative ring $R$, $(F^{-1} \WCart^{\mathrm{HT}})(R)$ is the full subcategory of $\WCart(R)$
consisting of those Cartier-Witt divisors $(I, \alpha)$ for which the composite map $I \xrightarrow{\alpha} W(R) \xrightarrow{F} W(R) \twoheadrightarrow R$ vanishes.
Note that $F$ restricts to a morphism of stacks $F^{-1} \WCart^{\mathrm{HT}} \rightarrow \WCart^{\mathrm{HT}}$, which we will also denote by $F$.
\end{notation}

\begin{remark}
Let $(A^{\bullet}, I^{\bullet})$ be the cosimplicial prism of Notation \ref{notation:simplicial-prism}. For each $n \geq 0$, the prism $(A^{n}, I^{n} )$ is transversal,
so that the Frobenius pullback $\varphi_{A^{n}}^{\ast}( I^{n} )$ can be identified with an invertible ideal in $A^{n}$ (Lemma \ref{lemma:easy-monomorphism}).
Unwinding the definitions, we can identify $F^{-1} \WCart^{\mathrm{HT}}$ with the geometric realization of the simplicial $p$-adic formal scheme
$\Spf( A^{\bullet} / \varphi_{A^{\bullet}}^{\ast} I^{\bullet} )$. 
\end{remark}

\begin{remark}\label{remark:slashop}
Let $R$ be a commutative ring and suppose we are given an $R$-valued point of the Hodge-Tate divisor $\WCart^{\mathrm{HT}}$, corresponding to a Cartier-Witt divisor
$(I,\alpha)$ for $\Spec(R)$. Then Proposition \ref{proposition:Frobenius-square} supplies an isomorphism of the Frobenius pullback of $(I,\alpha)$
with the Cartier-Witt divisor $(W(R), p)$. In particular, the following conditions are equivalent:
\begin{itemize}
\item[$(1)$] The commutative ring $R$ is an $\F_p$-algebra.
\item[$(2)$] The Cartier-Witt divisor $(I,\alpha)$ corresponds to an $R$-valued point of $F^{-1} \WCart^{\mathrm{HT}}$.
\end{itemize}
That is, the intersection $\WCart^{\mathrm{HT}} \cap F^{-1} \WCart^{\mathrm{HT}}$ coincides (as a closed substack of $\WCart$)
with the product $\WCart^{\mathrm{HT}}_{\F_p} = \WCart^{\mathrm{HT}} \times \Spec(\F_p)$.
\end{remark}

It follows from Remark \ref{remark:slashop} that the commutative diagram of (\ref{equation:Frobenius-pushout-square}) restricts to a diagram of closed substacks
\begin{equation}\label{equation:HT-pushout-square}
\xymatrix@R=50pt@C=50pt{ \WCart^{\mathrm{HT}}_{\F_p} \ar[r] \ar[d] & F^{-1} \WCart^{\mathrm{HT}} \ar[d]^{F} \\
\Spec(\F_p) \ar[r]^-{ \rho_{\dR}^{\mathrm{HT}} } & \WCart^{\mathrm{HT} }. }
\end{equation}

\begin{lemma}\label{lemma:sheaf-mod-p}
Let $\mathscr{E}$ be a quasi-coherent complex on $\WCart^{\mathrm{HT}}$, and let $F^{\ast} \WCart^{\mathrm{HT}}$ denote its
pullback to $F^{-1} \WCart^{\mathrm{HT}}$. Then the tautological map
$$ \rho: \F_p \otimes^{L} \RGamma( F^{-1} \WCart^{\mathrm{HT}}, F^{\ast} \mathscr{E} ) \rightarrow
\RGamma( \WCart^{\mathrm{HT}}_{\F_p}, (F^{\ast} \mathscr{E})|_{ \WCart^{\mathrm{HT}}_{\F_p}} )$$
is an isomorphism in the $\infty$-category $\calD( \F_p )$.
\end{lemma}

\begin{proof}
Let $\mathscr{F}$ denote the pullback $F^{\ast} \mathscr{E}$, and let $\mathscr{F}_0 = \mathscr{F}|_{ \WCart^{\mathrm{HT}}_{\F_p}}$.
Let $\calI \subseteq \calO_{ F^{-1} \WCart^{\mathrm{HT}} }$ denote the restriction of the Hodge-Tate ideal sheaf.
It follows from Remark \ref{remark:slashop} that the closed substack $\WCart^{\mathrm{HT}}_{ \F_p } \subseteq F^{-1} \WCart^{\mathrm{HT}}$ is
the vanishing locus of $\calI$. On the other hand, the closed substack $\Spec(\F_p) \times F^{-1} \WCart^{\mathrm{HT}}$ is the vanishing
locus of $\calI^{p}$. Unwinding the definitions, we can identify $\rho$ with the tautological map
$$ \RGamma( F^{-1} \WCart^{\mathrm{HT}}, \mathscr{F}/ \calI^{p} \mathscr{F} ) \rightarrow
\RGamma( F^{-1} \WCart^{\mathrm{HT}}, \mathscr{F} / \calI \mathscr{F} ).$$
To show that this map is an isomorphism, it will suffice to show that the complex
$$ \RGamma( F^{-1} \WCart^{\mathrm{HT}}, \calI^{n} \mathscr{F} / \calI^{n+1} \mathscr{F} )
\simeq \RGamma( \WCart^{\mathrm{HT}}_{\F_p}, \mathscr{F}_0\{n\} )$$ 
vanishes for $0 < n < p$. Let us abuse notation by identifying $\mathscr{F}$ with its direct image
in the Hodge-Tate divisor $\WCart^{\mathrm{HT}}$, and let $\Theta_{\mathscr{F}}$ be the Sen operator
of Construction \ref{construction:fiber-at-eta2}. It follows from the commutativity of the diagram
(\ref{equation:HT-pushout-square}) that $\Theta_{ \mathscr{F} }$ vanishes, so that
$\Theta_{ \mathscr{F}\{n\} }$ is given my multiplication by $n$. For $0 < n < p$, the
operator $\Theta_{ \mathscr{F}\{n\} }$ is invertible, so the desired vanishing follows from
Proposition \ref{proposition:easy-version}.
\end{proof}

\begin{lemma}\label{lemma:simple-structure}
Let $n$ be an integer. Then the complex $\RGamma( F^{-1} \WCart^{\mathrm{HT}}, F^{\ast} ( \calO_{ \WCart^{\mathrm{HT}}}\{n\}) ) \in \widehat{\calD}(\Z_p)$
is quasi-isomorphic either to a direct sum $\Z_p \oplus \Z_p[-1]$ or to $(\Z/p^{k}\Z)[-1]$ for some integer $k > 0$.
\end{lemma}

\begin{proof}
Since $\RGamma( F^{-1} \WCart^{\mathrm{HT}}, F^{\ast} ( \calO_{ \WCart^{\mathrm{HT}}}\{n\})) )$ is $p$-complete,
it will suffice to show that the derived tensor product
$$ \F_p \otimes^{L} \RGamma( F^{-1} \WCart^{\mathrm{HT}}, F^{\ast} ( \calO_{ \WCart^{\mathrm{HT}}}\{n\}))$$
is isomorphic to $\F_p \oplus \F_p[-1]$ as an object of $\widehat{\calD}(\F_p)$.
By virtue of Lemma \ref{lemma:sheaf-mod-p}, we can rewrite this tensor product as
$$ \RGamma( \WCart^{\mathrm{HT}}_{\F_p}, \calO_{\WCart^{\mathrm{HT}}_{\F_p}}(pn) ) \simeq
\F_p \otimes^{L} \RGamma( \WCart^{\mathrm{HT}}, \calO_{\WCart^{\mathrm{HT}}_{\F_p}}(pn) ),$$
so the desired result follows from Proposition \ref{proposition:easy-version}.
\end{proof}

\begin{lemma}\label{rusk}
The coherent cohomology groups of $F^{-1} \WCart^{\mathrm{HT}}$ are given by
$$ \mathrm{H}^{n}( F^{-1} \WCart^{\mathrm{HT}}, \calO_{ F^{-1} \WCart^{\mathrm{HT}} } ) \simeq \begin{cases} \Z_p & \text{ if $n = 0,1$} \\
0 & \text{ otherwise.} \end{cases}$$
\end{lemma}

\begin{proof}
Since $F^{-1} \WCart^{\mathrm{HT}}$ is not empty, the coordinate ring $\mathrm{H}^{0}( F^{-1} \WCart^{\mathrm{HT}}, \calO_{ F^{-1} \WCart^{\mathrm{HT}} } )$
is nonzero. The desired result now follows from Lemma \ref{lemma:simple-structure}.
\end{proof}

The main ingredient in our proof of Theorem \ref{theorem:Frobenius-pushout-square} is the following calculation:

\begin{proposition}\label{proposition:main-ingredient}
Pullback along the Frobenius map $F: F^{-1} \WCart^{\mathrm{HT}} \rightarrow \WCart^{\mathrm{HT}}$ induces a monomorphism
$$ \Z_p \simeq \mathrm{H}^{1}( \WCart^{\mathrm{HT}}, \calO_{\WCart^{\mathrm{HT}} } ) 
\rightarrow \mathrm{H}^{1}( F^{-1} \WCart^{\mathrm{HT}}, \calO_{ F^{-1} \WCart^{\mathrm{HT}} } ) \simeq \Z_p,$$
whose image is the subgroup $p\Z_p$.
\end{proposition}

\begin{proof}
Let $R$ be a commutative ring which is $p$-adically separated and complete, and let $u$ be an invertible element of $W(R)$. Writing $u = \sum V^{n} [a_n]$, we let
$\ell(u)$ denote the element of $R$ given by the power series expansion
$$ \frac{1}{p} \log( 1 + p \frac{a_1}{a_0^{p}} ) = \sum_{n > 0} \frac{ (-p)^{n-1} }{ n } \frac{ a_1^n}{ a_0^{pn} }.$$
Writing $\gamma_n(u) = a_0^{p^{n}} + p a_{1}^{p^{n-1}} + \cdots + p^{n} a_n$ for the $n$th ghost component of $u$,
we can write $\ell(u)$ more informally as $\frac{ 1}{p} \log( \frac{ \gamma_1(u) }{ \gamma_0(u)^{p}} )$. It follows that the construction $u \mapsto \ell(u)$
determines a group homomorphism $W(R)^{\times} \rightarrow R = \mathbf{G}_{a}(R)$, depending functorially on $R$.

In what follows, let us abuse notation by identifying the group schemes $W^{\times}$, $\mathbf{G}_{a}$, and $\mathbf{G}_{m}$
with their products with $\Spf(\Z_p)$ (that is, viewing them as functors defined only on the category of commutative rings in which $p$ is nilpotent).
The preceding construction determines a homomorphism $\ell: W^{\times} \rightarrow \mathbf{G}_{a}$ of commutative formal group schemes over $\Z_p$.
Let $\beta \in \mathrm{H}^{1}( \WCart, \calO_{ \WCart} )$ denote the cohomology class determined by the composite map
$$ \WCart \simeq [\WCart_{0} / W^{\times}] \rightarrow BW^{\times} \xrightarrow{\ell} B \mathbf{G}_{a}.$$
Note that the composite map $$\mathbf{G}_{m}^{\sharp} \simeq W^{\times}[F] \hookrightarrow W^{\times} \xrightarrow{\ell} \mathbf{G}_{a}$$
is given by $t \mapsto - \log(t)$. Invoking Example \ref{example:cohomology-of-HT}, we see that $x|_{\WCart^{\mathrm{HT}}}$ is a generator of the cohomology group
$\mathrm{H}^{1}( \WCart^{\mathrm{HT}}, \calO_{ \WCart^{\mathrm{HT}} }) \simeq \Z_p$. In particular, the restriction
$\beta|_{ \WCart^{\mathrm{HT}}_{\F_p} }$ is a generator of the group $\mathrm{H}^{1}( \WCart^{\mathrm{HT}}_{\F_p}, \calO_{ \WCart^{\mathrm{HT}}_{\F_p} }) \simeq \F_p$.
Using Lemmas \ref{rusk} and \ref{lemma:sheaf-mod-p}, we deduce that $\beta|_{ F^{-1} \WCart^{\mathrm{HT}} }$ is a generator of the group
$\mathrm{H}^{1}( F^{-1} \WCart^{\mathrm{HT}}, \calO_{ F^{-1} \WCart^{\mathrm{HT}} }) \simeq \Z_p$.
Let $F^{\ast}(\beta)$ denote the pullback of $\beta$ along the Frobenius map $F: \WCart \rightarrow \WCart$, which we regard as an element of
$\mathrm{H}^{1}( \WCart, \calO_{ \WCart} )$. We will complete the proof by showing that the cohomology classes $F^{\ast}(\beta)$ and 
$-p\beta$ have the same restriction to $F^{-1} \WCart$. 

Let $X \subseteq W$ be the $p$-adic formal scheme whose $R$-valued points are given by Witt vectors $x = \sum V^{n} [a_n] \in W(R)$
where $a_1$ is invertible and the ghost component $\gamma_1(x) = a_0^{p} + pa_1$ vanishes.
The formal group scheme $W^{\times}$ acts on $X$ by multiplication, and $F^{-1} \WCart$ can be identified with the stack-theoretic quotient
$[X / W^{\times}]$. Unwinding the definitions, we see that the cohomology class $(F^{\ast}(\beta) + p\beta)|_{ F^{-1} \WCart}$ is represented
by the $1$-cocycle
$$ g: W^{\times} \times X \rightarrow \mathbf{G}_{a} \quad \quad g(u,x) = \ell( F(u) ) + p \ell(u) = \ell( F(u)u^{p} ).$$
We will complete the proof by showing that $g$ is the coboundary of the $0$-cochain
$$ f: X \rightarrow \mathbf{G}_{a} \quad \quad f(x) = \ell( \delta(x) );$$
that is, the morphism $f$ satisfies $f(ux) = f(x) + g(u,x)$ for every pair of $R$-valued points $u \in W(R)^{\times}$ and $x \in X(R)$,
where $R$ is a commutative ring which is $p$-adically separated and complete. Note that it suffices to check this in the {\em universal} case
where $R$ is the coordinate ring of the affine formal scheme $W^{\times} \times X$; in particular, we may assume without loss of generality
that $R$ is $p$-torsion-free. Unwinding the definitions, we wish to establish the identity
$\ell( \delta(ux) ) = \ell( \delta(x)F(u)u^{p} )$. Note that, for an invertible element $v = \sum V^{n} [a_n] \in W(R)$, the
element $\ell(v)$ depends only on $\frac{a_1}{ a_0^{p}} = \frac{1}{p} ( \frac{ \gamma_1(v) }{ \gamma_0(v)^{p}} - 1)$,
and therefore only on the ratio $\frac{ \gamma_1(v) }{ \gamma_0(v)^{p} }$. We now compute
\begin{eqnarray*}
\frac{ \gamma_1( \delta(ux) )}{ \gamma_0( \delta(ux) )^{p} } & = & \frac{ \frac{ \gamma_1( \varphi(ux) - (ux)^{p} }{p} }{ ( \frac{ \gamma_0( \varphi(ux) - (ux)^p ) }{p} )^{p} }\\
& = & \frac{ \frac{ \gamma_2(u) \gamma_2(x) - \gamma_1(u)^{p} \gamma_1(x)^{p} }{p} }{ ( \frac{ \gamma_1(u) \gamma_1(x) - \gamma_0(u)^{p} \gamma_0(x)^p}{p} )^{p} } \\
& = & (-1)^{p} p^{p-1} \frac{ \gamma_2(u) \gamma_2(x) }{ \gamma_0(u)^{p^2} \gamma_0(x)^{p^2} }
\end{eqnarray*}
where the last equality follows from the identity $\gamma_1(x) = 0$. Combining this with the analogous calculation for $u=1$, we obtain
\begin{eqnarray*}
\frac{ \gamma_1( \delta(ux) )}{ \gamma_0( \delta(ux) )^{p} }  & = & \frac{ \gamma_1( \delta(x) )}{ \gamma_0( \delta(x) )^{p} } \cdot
\frac{ \gamma_2(u) }{ \gamma_0(u)^{p^2} } \\
& = & \frac{ \gamma_1( \delta(x) )}{ \gamma_0( \delta(x) )^{p} } \cdot \frac{\gamma_2(u)}{ \gamma_1(u)^p} \cdot \frac{ \gamma_1(u)^p}{ \gamma_0(u)^{p^2}} \\
& = & \frac{ \gamma_1( \delta(x) \cdot F(u) \cdot u^p) }{ \gamma_0( \delta(x) \cdot F(u) \cdot u^{p})^{p} },
\end{eqnarray*}
as desired.
\end{proof}

\begin{lemma}\label{lemma:Frobenius-contraction-HT}
Let $\mathscr{E}$ be a quasi-coherent complex on $\WCart^{\mathrm{HT}}$. Then the diagram (\ref{equation:HT-pushout-square}) induces a pullback diagram
\begin{equation}\label{equation:Frobenius-contraction-HT}
 \xymatrix@R=50pt@C=50pt{ \RGamma( \WCart^{\mathrm{HT}}_{\F_p}, (F^{\ast} \mathscr{E})|_{ \WCart^{\mathrm{HT}}_{\F_p}} ) & \RGamma( F^{-1} \WCart^{\mathrm{HT}}, F^{\ast} \mathscr{E} ) \ar[l] \\
\RGamma( \Spec(\F_p), \mathscr{E}|_{\Spec(\F_p)}) \ar[u] & \RGamma( \WCart^{\mathrm{HT}}, \mathscr{E} ) \ar[l] \ar[u] } \end{equation}
in the $\infty$-category $\widehat{\calD}(\Z_p)$.
\end{lemma}

\begin{proof}
By virtue of Corollary \ref{corollary:HT-computation}
$$ \RGamma( \WCart^{\mathrm{HT}}, \bullet): \calD( \WCart^{\mathrm{HT}} ) \rightarrow \widehat{\calD}(\Z_p)$$
commutes with colimits. It follows immediately that the global sections functor
$$ \RGamma( \WCart^{\mathrm{HT}}_{\F_p}, \bullet): \calD( \WCart^{\mathrm{HT}} ) \rightarrow \calD(\F_p)$$
has the same property. Applying Lemma \ref{lemma:sheaf-mod-p}, we see that the functor
$$ \RGamma( F^{-1} \WCart^{\mathrm{HT}}, F^{\ast}(\bullet)): \calD( \WCart^{\mathrm{HT}} ) \rightarrow \widehat{\calD}(\Z_p)$$
also commutes with colimits. Consequently, the collection of those objects $\mathscr{E} \in \calD( \WCart^{\mathrm{HT}} )$ for
which the diagram (\ref{equation:Frobenius-contraction-HT}) is a pullback square is closed under colimits. Invoking
Proposition \ref{proposition:twists-generate}, we are reduced to proving Lemma \ref{lemma:Frobenius-contraction-HT} in the special
case where $\mathscr{E} = \calO_{ \WCart^{\mathrm{HT}} }\{n\}$ for some integer $n$. We now proceed in several steps:
\begin{itemize}
\item Suppose first that $n=0$. By virtue of Lemma \ref{lemma:sheaf-mod-p}, we have a commutative diagram of fiber
sequences
$$ \xymatrix@R=50pt@C=50pt{  \fib(\epsilon) \ar[r]^-{\alpha} \ar[d] &  \RGamma( F^{-1} \WCart^{\mathrm{HT}}, \calO_{ F^{-1} \WCart^{\mathrm{HT}} } ) \ar[d]^{p} \\
\RGamma( \WCart^{\mathrm{HT}}, \calO_{ \WCart^{\mathrm{HT}} } ) \ar[d]^{\epsilon} \ar[r] & \RGamma( F^{-1} \WCart^{\mathrm{HT}}, \calO_{ F^{-1} \WCart^{\mathrm{HT}} } ) \ar[d] \\
\F_p \ar[r] & \RGamma( \WCart^{\mathrm{HT}}_{\F_p}, \calO_{\WCart^{\mathrm{HT}}_{\F_p}} ) }$$
and we wish to show that $\alpha$ induces an isomorphism on cohomology. Note that the cohomology groups of both of the upper complexes
are concentrated in degrees $0$ and $1$, and that $\alpha$ clearly induces an isomorphism on cohomology in degree zero
(since the restriction map $\mathrm{H}^{0}( \WCart^{\mathrm{HT}}, \calO_{\WCart^{\mathrm{HT}}} ) \rightarrow \mathrm{H}^{0}( F^{-1} \WCart^{\mathrm{HT}}, \calO_{ F^{-1} \WCart^{\mathrm{HT}} } )$ can be identified with a ring homomorphism from $\Z_p$ to itself, and is therefore an isomorphism).
We are therefore reduced to checking that $\alpha$ induces an isomorphism on cohomology in degree $1$, which is a reformulation of Proposition \ref{proposition:main-ingredient}.

\item Suppose that $n$ is divisible by $p$. Since (\ref{equation:Frobenius-contraction-HT}) is a diagram of
$p$-complete complexes, it suffices to show that it is a pullback diagram after replacing $\mathscr{E} = \calO_{\WCart^{\mathrm{HT}}}\{n\}$ by the
derived tensor product $\F_p \otimes^{L} \mathscr{E}$. It follows from Theorem \ref{theorem:compute-with-HT} that the line bundles
$\calO_{\WCart^{\mathrm{HT}}}$ and $\calO_{\WCart^{\mathrm{HT}}}\{n\}$ become isomorphic after extending scalars to $\F_p$, so the desired result
follows from the previous step.

\item Suppose that $n$ is not divisible by $p$. In this case, the complex $\RGamma( \WCart^{\mathrm{HT}}, \calO_{\WCart^{\mathrm{HT}}}\{n\} )$ is acyclic.
Set $\mathscr{L} = F^{\ast} \calO_{\WCart^{\mathrm{HT}}}\{n\} \in \calD( F^{-1} \WCart^{\mathrm{HT}} )$
and let $\mathscr{L}_0 = \mathscr{L}|_{ \WCart^{\mathrm{HT}}_{\F_p}}$. 
Then $\mathscr{L}_0 \simeq \calO_{\WCart^{\mathrm{HT}}_{\F_p}}\{pn\}$ is a trivial line bundle on
$\WCart^{\mathrm{HT}}_{\F_p}$, so the complex $\RGamma( \WCart^{\mathrm{HT}}_{\F_p}, \mathscr{L}_0)$
is isomorphic to $\F_p \oplus \F_p[-1]$ (as an object of $\calD(\F_p)$) and $\RGamma( \Spec(\F_p), \mathscr{E}|_{\Spec(\F_p)})$ maps isomorphically onto the first factor.
We will show that the complex $\RGamma( F^{-1} \WCart^{\mathrm{HT}}, \mathscr{L} )$ is isomorphic to $\F_p[-1]$. It will then
follow from Lemma \ref{lemma:sheaf-mod-p} that the restriction map 
$$\RGamma( F^{-1} \WCart^{\mathrm{HT}}, \mathscr{L} ) \rightarrow \RGamma( \WCart^{\mathrm{HT}}_{\F_p}, \mathscr{L}_{0} )$$
induces an isomorphism on cohomology in degree $1$, which will complete the proof that (\ref{equation:Frobenius-contraction-HT}) is a pullback diagram.

Assume, for a contradiction, that $\RGamma( F^{-1} \WCart^{\mathrm{HT}}, \mathscr{L} )$ is not isomorphic to $\F_p[-1]$. 
Applying Lemma \ref{lemma:simple-structure}, we see that the derived tensor product 
$$(\Z / p^2 \Z) \otimes^{L} \RGamma( F^{-1} \WCart^{\mathrm{HT}}, \mathscr{L} ) \simeq \RGamma( F^{-1} \WCart^{\mathrm{HT}}, \mathscr{L} / p^2 \mathscr{L} )$$
is isomorphic to $(\Z / p^2\Z) \oplus ( \Z / p^2 \Z)[-1]$. Choose a generator $e$ of the cyclic group
$$\mathrm{H}^{0}( F^{-1} \WCart^{\mathrm{HT}}, \mathscr{L} / p^2 \mathscr{L} ) \simeq \Z / p^2 \Z.$$ Note that $e$ restricts to a generator of
the group $\mathrm{H}^{0}( \WCart^{\mathrm{HT}}_{\F_p}, \mathscr{L}_0)$, and can therefore be regarded as a trivialization of the line bundle $\mathscr{L}$
over the closed substack $\Spec( \Z / p^2\Z) \times F^{-1} \WCart^{\mathrm{HT}}$. Let $u$ denote the image of $e$ under the boundary map
$$ \partial: \mathrm{H}^{0}(  F^{-1} \WCart^{\mathrm{HT}}, \mathscr{L} / p^2 \mathscr{L} )
\rightarrow \mathrm{H}^{1}(  F^{-1} \WCart^{\mathrm{HT}}, p^2 \mathscr{L} / p^3 \mathscr{L} ).$$
Since the target of this boundary map is a $p$-torsion abelian group, we have $pu= 0$. It follows that
the $p$th power $e^{p}$ extends to a trivialization of $\mathscr{L}^{p}$ over the closed substack
$\Spec( \Z / p^3 \Z) \times F^{-1} \WCart^{\mathrm{HT}}$. In particular, the derived tensor product
$$ ( \Z / p^3 \Z) \otimes^{L} \RGamma( F^{-1} \WCart^{\mathrm{HT}}, \mathscr{L}^{p} )$$
is quasi-isomorphic to $(\Z / p^3 \Z) \oplus (\Z / p^3 \Z)[-1]$, so the cyclic group
$\mathrm{H}^{1}( F^{-1} \WCart^{\mathrm{HT}}, \mathscr{L}^{p})$ must have order $> p^2$. 
On the other hand, applying Lemma \ref{lemma:Frobenius-contraction-HT} to the line bundle
$\mathscr{E} = \calO_{ \WCart^{\mathrm{HT}}}\{pn\}$ yields a short exact sequence of abelian groups
$$ \mathrm{H}^{1}( \WCart^{\mathrm{HT}}, \calO_{\WCart^{\mathrm{HT}}}\{pn\} )
\rightarrow \mathrm{H}^{1}(  F^{-1} \WCart^{\mathrm{HT}}, \mathscr{L}^{p} )
\rightarrow \mathrm{H}^{1}( \WCart^{\mathrm{HT}}_{\F_p}, \mathscr{L}_{0}^{p}),$$
where the third term is cyclic of order $p$. It follows that the order of the cohomology group
$ \mathrm{H}^{1}( \WCart^{\mathrm{HT}}, \calO_{\WCart^{\mathrm{HT}}}\{pn\} )$ must be larger than $p$.
Proposition \ref{proposition:easy-version} supplies an isomorphism
$$\mathrm{H}^{1}( \WCart^{\mathrm{HT}}, \calO_{\WCart^{\mathrm{HT}}}\{pn\} ) \simeq \Z_p / (pn) \Z_p,$$ contradicting
our assumption that $n$ is not divisible by $p$.
\end{itemize}
\end{proof}

\begin{remark}[Transversality]\label{remark:transversality-stuff}
Let $\iota: \WCart^{\mathrm{HT}} \hookrightarrow \WCart$ and $\iota': F^{-1} \WCart^{\mathrm{HT}} \hookrightarrow \WCart$ be the inclusion maps,
so that we have a commutative diagram of stacks
\begin{equation}\label{equation:transversal-diagram}
\xymatrix@R=50pt@C=50pt{ \WCart^{\mathrm{HT}}_{\F_p} \ar[r] \ar[d] & \WCart^{\mathrm{HT}} \ar[d]^{\iota} \\
F^{-1} \WCart^{\mathrm{HT}} \ar[r]^-{\iota'} \ar[d]^{F^{\mathrm{HT}} } & \WCart \ar[d]^{F} \\
\WCart^{\mathrm{HT}} \ar[r]^-{\iota} & \WCart. }
\end{equation}
where both squares are pullbacks (Remark \ref{remark:slashop}). These squares are transversal in the following sense:
\begin{itemize}
\item For every quasi-coherent complex $\mathscr{E}$ on $\WCart^{\mathrm{HT}}$, the bottom square in (\ref{equation:transversal-diagram})
induces an isomorphism
$$ \RGamma( \WCart, F^{\ast} \iota_{\ast} \mathscr{E} ) \rightarrow \RGamma( F^{-1} \WCart^{\mathrm{HT}}, F^{\mathrm{HT} \ast} \mathscr{E} )$$
in the $\infty$-category $\widehat{\calD}(\Z_p)$.

\item For every quasi-coherent complex $\mathscr{F}$ on $F^{-1} \WCart^{\mathrm{HT}}$, the upper square in (\ref{equation:transversal-diagram})
induces an isomorphism
$$ \RGamma( \WCart^{\mathrm{HT}} , (\iota'_{\ast} \mathscr{F})|_{ \WCart^{\mathrm{HT}} }) \rightarrow
\RGamma( \WCart^{\mathrm{HT}}_{\F_p}, \mathscr{F}|_{ \WCart^{\mathrm{HT}}_{\F_p} }).$$ 
\end{itemize}
Writing $(A^{\bullet}, I^{\bullet})$ for the cosimplicial prism of Notation \ref{notation:simplicial-prism}, this assertion follows from the observation
that both squares in the diagram
$$ \xymatrix@R=50pt@C=50pt{ A^{\bullet} / ( I^{\bullet}, p) & A^{\bullet} / I^{\bullet} \ar[l] \\
A^{\bullet} / F^{\ast}(I^{\bullet} ) \ar[u] & A^{\bullet} \ar[u] \ar[l] \\
A^{\bullet} /I^{\bullet} \ar[u]  & A^{\bullet} \ar[u]^{\varphi} \ar[l] }$$
are pushout squares of (cosimplicial) {\em animated} commutative rings: that is, each of the inclusion maps $I^{n} \hookrightarrow A^{n}$
remains injective after base change along the Frobenius map $\varphi: A^{n} \rightarrow A^{n}$ (Lemma \ref{lemma:easy-monomorphism} and
and after base change along the quotient map $A^{n} \twoheadrightarrow A^{n} / F^{\ast}(I^{n} )$.
\end{remark}

\begin{proof}[Proof of Theorem \ref{theorem:Frobenius-pushout-square}]
Let $\mathscr{E}$ be a quasi-coherent complex on the Cartier-Witt stack $\WCart$; we wish to show that the diagram $\sigma_{\mathscr{E}}:$
$$ \xymatrix@R=50pt@C=50pt{ \RGamma( \WCart^{\mathrm{HT}}, (F^{\ast} \mathscr{E} )|_{ \WCart^{\mathrm{HT} }} ) & \RGamma( \WCart, F^{\ast} \mathscr{E} ) \ar[l] \\
\RGamma( \Spf(\Z_p), \rho_{\dR}^{\ast} \mathscr{E} ) \ar[u] & \RGamma( \WCart, \mathscr{E} ) \ar[u] \ar[l] }$$
is a pullback square in $\widehat{\calD}(\Z_p)$. Let $\calI$ denote the Hodge-Tate ideal sheaf on $\WCart$. Writing
$\mathscr{E}$ as a limit of sheaves of the form $\mathscr{E} / \calI^{n} \mathscr{E}$, we are reduced to proving that each of the diagrams
$\sigma_{ \mathscr{E} / \calI^{n} \mathscr{E} }$ is a pullback square. Proceeding by induction on $n$, we can reduce to showing that
$\sigma_{\mathscr{E}}$ is a pullback square in the special case where $\mathscr{E} = \iota_{\ast} \mathscr{F}$ is a quasi-coherent complex
on the Hodge-Tate divisor $\WCart^{\mathrm{HT}}$. By virtue of Remark \ref{remark:transversality-stuff}, we can identify $\sigma_{\mathscr{F}}$
with the diagram
$$ \xymatrix@R=50pt@C=50pt{ \RGamma( \WCart^{\mathrm{HT}}_{\F_p}, (F^{\ast} \mathscr{F} )|_{ \WCart^{\mathrm{HT}}_{\F_p} } ) & \RGamma( F^{-1} \WCart^{\mathrm{HT}}, F^{\ast} \mathscr{F} ) \ar[l] \\
\RGamma( \Spec(\F_p),\mathscr{F}|_{\Spec(\F_p)} ) \ar[u] & \RGamma( \WCart^{\mathrm{HT}}, \mathscr{F} ) \ar[u] \ar[l] }$$
determined by (\ref{equation:HT-pushout-square}), which is a pullback square by virtue of Lemma \ref{lemma:Frobenius-contraction-HT}.
\end{proof}

\subsection{Exponentiating the Sen Operator}\label{subsection:exp-theta}

Let $\mathscr{E}$ be a quasi-coherent complex on the Hodge-Tate divisor $\WCart^{\mathrm{HT}}$, and let
$\mathscr{E}_{\eta} \in \widehat{\calD}(\Z_p)$ denote its fiber at the point $\eta: \Spf(\Z_p) \rightarrow \WCart^{\mathrm{HT}}$
of Construction \ref{construction:fiber-at-eta}. Note that that $\mathscr{E}_{\eta}$ is equipped with an action of the group
$\Aut(\eta)(\Z_p) = \mathbf{G}_{m}^{\sharp}(\Z_p) = (1+p\Z_p)^{\times}$. In particular, each element $u \in (1+p\Z_p)^{\times}$
determines an automorphism of $\mathscr{E}_{\eta}$, which we will denote by $\gamma_{u}$. By virtue of
Theorem \ref{theorem:compute-with-HT}, the automorphism $\gamma_{u}$ is completely determined by the
Sen operator $\Theta_{\mathscr{E}}$ of Notation \ref{notation:Sen-operator-second}. Heuristically, the
automorphism $\gamma_{u}$ can be described as the exponential $u^{\Theta_{\mathscr{E}}}$. This heuristic
can be made precise as follows:

\begin{proposition}\label{proposition:formula-for-action-simplified}
Let $\mathscr{E}$ be a quasi-coherent complex on $\WCart^{\mathrm{HT}}$ and let
$u$ be an element of the group $(1+p\Z_p)^{\times}$. Assume that, if $p=2$, then $u \equiv 1 \pmod{4}$.
Then the automorphism $\gamma_{u}$ of $\mathscr{E}_{\eta}$ satisfies the congruence
$$ \gamma_{u} \equiv  \exp( \log(u) \cdot \Theta_{\mathscr{E}} )(x) = \sum_{m \geq 0} \frac{ \log(u)^{m} }{m!} \Theta_{\mathscr{E}}^{m} \pmod{p^d}$$
for every integer $d \geq 0$.
\end{proposition}

\begin{remark}
In the situation of Proposition \ref{proposition:formula-for-action-simplified}, the sequence of divided powers $\{ \frac{ \log(u)^m }{ m!} \}_{m \geq 0}$ converges $p$-adically to zero, so that the sum $\sum_{m \geq 0} \frac{ \log(u)^{m} }{m!} (\Theta_{\mathscr{E}} )^{m}$ has only finitely many nonzero term modulo $p^d$
for any fixed integer $d \geq 0$. Beware that when $p=2$, this convergence requires the assumption that $u \equiv 1 \pmod{4}$.
\end{remark}

\begin{proof}[Proof of Proposition \ref{proposition:formula-for-action-simplified}]
Note that the unit map
$$\mathscr{E} \rightarrow \eta_{\ast} \eta^{\ast}(\mathscr{E}) \simeq \mathscr{E}_{\eta} \otimes^{L} \eta_{\ast}( \calO_{ \Spf(\Z_p) }  )$$
is a (nonequivariantly split) monomorphism. We may therefore replace $\mathscr{E}$ by
$\eta_{\ast} \eta^{\ast}( \mathscr{E} )$, and thereby reduce to the case where $\mathscr{E}$ is a tensor product of some complex $M \in \widehat{\calD}( \Z_p)$
with $\eta_{\ast} \calO_{\Spf(\Z_p )}$. Without loss of generality we may assume that $M = \Z_p$, so that $\mathscr{E}_{\eta}$ can be identified
with the $p$-completed coordinate ring $\widehat{\calO}_{\mathbf{G}_{m}}^{\sharp}$ of the group scheme $\mathbf{G}_{m}^{\sharp}$.
Under this identification, $\gamma_{u}$ corresponds to the automorphism $f(t) \mapsto f(u \cdot t)$ of $\widehat{\calO}_{\mathbf{G}_{m}}^{\sharp}$,
and the Sen operator $\Theta_{\mathscr{E}}$ is given by the differential operator $t \frac{ \partial }{\partial t}$ (Example \ref{example:regular-representation}).
Using these identifications, we are reduced to verifying the identity
$$ f(u \cdot t ) = \sum_{ m \geq 0} \frac{ \log(u)^m}{m!} ( t \frac{ \partial}{ \partial t} )^{m} f(t).$$
for $f \in \widehat{\calO}_{\mathbf{G}_{m}}^{\sharp}$. Setting $x = \log(u) \in p \Z_p$, $y = \log(t) \in (t-1) \Q_p[[t-1]]$,
and $g(y) = f( \exp(y) )$, this reduces to Taylor's formula
$$ g(x+y) = \sum_{m \geq 0} \frac{ x^m}{m!} (\frac{ \partial}{ \partial y})^{m} g(y).$$
\end{proof}

Let us describe an application of Proposition \ref{proposition:formula-for-action-simplified}
which will be useful in \S\ref{subsection:computing-with-WCart}. Note that, for each element
$u \in (1+p\Z_p)^{\times}$, the divided powers $\frac{ (u-1)^{n} }{n!}$ are divisible by $p$ for each $n > 0$.
It follows that $u$ is annihilated by the group homomorphism
$$(1+p\Z_p)^{\times} \simeq \mathbf{G}_{m}^{\sharp}(\Z_p) \rightarrow \mathbf{G}_{m}^{\sharp}(\F_p).$$
In particular, if $\mathscr{E}$ is a quasi-coherent complex on $\WCart^{\mathrm{HT}}$, then the automorphism $\gamma_{u}$ induces the identity on the complex $\F_p \otimes^{L} \mathscr{E}_{\eta}$, so that the difference $\gamma_{u} - \id_{\mathscr{E}_{\eta}}$ is divisible by $p$. In fact, we can be more precise:

\begin{proposition}\label{proposition:fiber-sequence-global-sections}
Let $p$ be an odd prime and let $u$ be a topological generator of the group $(1+p\Z_p)^{\times}$. For every
quasi-coherent complex $\mathscr{E}$ on $\WCart^{\mathrm{HT}}$, there is a canonical fiber sequence
$$ \RGamma( \WCart^{\mathrm{HT}}, \mathscr{E} ) \rightarrow \mathscr{E}_{\eta} \xrightarrow{ \xi_{u} } \mathscr{E}_{\eta},$$
where $\xi_{u}$ is an endomorphism of $\mathscr{E}_{\eta}$ satisfying $p \xi_u = \gamma_u - \id_{\mathscr{E}_{\eta}}$.
\end{proposition}

\begin{remark}
Stated more informally, Proposition \ref{proposition:fiber-sequence-global-sections} asserts that the endomorphism
$\gamma_{u} - \id$ is {\it canonically} divisible by $p$, and that the fiber of the resulting map
$\frac{ \gamma_u - \id}{p}: \mathscr{E}_{\eta} \rightarrow \mathscr{E}_{\eta}$ can be identified with the complex of global sections $\RGamma( \WCart^{\mathrm{HT}}, \mathscr{E} )$.
\end{remark}

\begin{remark}
Let $\mathscr{E}$ be a quasi-coherent complex on the Hodge-Tate divisor $\WCart^{\mathrm{HT}}$, and let $K \in \widehat{\calD}(\Z_p)$ denote
the cofiber of the restriction map $\RGamma( \WCart^{\mathrm{HT}}, \mathscr{E} ) \rightarrow \mathscr{E}_{\eta}$ denote the restriction map.
It follows from Proposition \ref{proposition:easy-version} that there is a canonical isomorphism $\alpha: K \simeq \mathscr{E}_{\eta}$
for which the composite map $\mathscr{E}_{\eta} \rightarrow K \xrightarrow{\alpha} \mathscr{E}_{\eta}$ is equal to the Sen operator $\Theta_{\mathscr{E}}$ of Construction \ref{construction:fiber-at-eta2}. The content of Proposition \ref{proposition:fiber-sequence-global-sections} is that, if $p$ is an odd prime
and $u$ is a topological generator of $(1+p\Z_p)^{\times}$, then there is a {\em different} isomorphism $\alpha'_{u}: K \simeq \mathscr{E}_{\eta}$ for which the composite map $$ \mathscr{\calE}_{\eta} \rightarrow K \xrightarrow{\alpha'_{u}} \mathscr{E}_{\eta} \xrightarrow{p} \mathscr{E}_{\eta}$$
can be identified with $\gamma_{u} - \id$. Note that the composition $\alpha'_{u} \circ \alpha^{-1}$ is an automorphism
of the complex $\mathscr{E}_{\eta}$. Heuristically, this automorphism is given by the formal series
$$ \frac{ u^{\Theta_{\mathscr{E}}} - \id}{ p \Theta_{\mathscr{E}}} = \sum_{n \geq 1} \frac{ \log(u)^n}{p \cdot n! } \Theta_{\mathscr{E}}^{n-1}$$
(see Proposition \ref{proposition:formula-for-action-simplified}). Beware that this formal series is not $p$-adically convergent in the case $p=2$.
(see Remark \ref{remark:even-prime-case} below).
\end{remark}

For our applications, it will be convenient to have a more canonical formulation of Proposition \ref{proposition:fiber-sequence-global-sections}, 
which does not depend on a choice of topological generator $u \in (1+p\Z_p)^{\times}$.
We begin by observing that if $\mathscr{E}$ is then the action of $(1+p\Z_p)^{\times}$ on $\mathscr{E}_{\eta}$
{\em continuous} with respect to the $p$-adic topologies on $\mathscr{E}$ and $(1+p\Z_p)^{\times}$, respectively:

\begin{notation}\label{notation:continuity-of-action}
In what follows, we regard $(1+p\Z_p)^{\times}$ as a profinite group, and we write $[\Spf(\Z_p) / (1+p\Z_p)^{\times}]$ for its $p$-complete classifying stack.
Roughly speaking, we can think of quasi-coherent complexes on the quotient stack $[\Spf(\Z_p) / (1+p\Z_p)^{\times}]$ as $p$-complete complexes
$M \in \widehat{\calD}(\Z_p)$ equipped with a $p$-adically continuous action of the profinite group $(1+p\Z_p)^{\times}$.

The action of $(1+p\Z_p)^{\times}$ on the map $\eta: \Spf(\Z_p) \rightarrow \WCart^{\mathrm{HT}}$ determines a morphism of stacks $\widetilde{\eta}: [\Spf(\Z_p) / (1+p\Z_p)^{\times}] \rightarrow \WCart^{\mathrm{HT}}$. If $\mathscr{E}$ is a quasi-coherent complex on $\WCart^{\mathrm{HT}}$, we write $\widetilde{\eta}^{\ast}(\mathscr{E})$ for the pullback of $\mathscr{E}$ to $[\Spf(\Z_p) / (1+p\Z_p)^{\times}]$. 
Note that, if $p$ is an odd prime and $u$ is a topological generator for the group $(1+p\Z_p)^{\times}$, then the complex of global sections
$\RGamma( [\Spf(\Z_p) / (1+p\Z_p)^{\times}], \widetilde{\eta}^{\ast}(\mathscr{E}) )$ can be identified with the fiber
$$\mathscr{E}_{\eta}^{u = 1} = \fib( \gamma_u - \id: \mathscr{E}_{\eta} \rightarrow \mathscr{E}_{\eta}) $$ 
\end{notation}

Since the group homomorphism $\mathbf{G}_{m}^{\sharp}( \Z_p) \rightarrow \mathbf{G}_m^{\sharp}( \F_p)$ is trivial, the morphism
$\widetilde{\eta}$ of Notation \ref{notation:continuity-of-action} fits into a commutative diagram of stacks
\begin{equation}\label{equation:clean-statement-of-fiber} 
\xymatrix@R=50pt@C=50pt{ \WCart^{\mathrm{HT} } & [\Spf(\Z_p) / (1+p\Z_p)^{\times}] \ar[l]_-{\widetilde{\eta}}  \\
\Spec(\F_p) \ar[u] & [\Spec(\F_p) / (1+p\Z_p)^{\times}]. \ar[u] \ar[l] }
\end{equation}

\begin{proposition}\label{proposition:clean-statement-of-fiber}
Let $p$ be an odd prime and let $\mathscr{E}$ be a quasi-coherent complex on the Hodge-Tate divisor $\WCart^{\mathrm{HT}}$.
Then the diagram (\ref{equation:clean-statement-of-fiber}) determines a pullback square
\begin{equation}\label{equation:clean-statement-of-fiber2}
\xymatrix@R=50pt@C=50pt{ \RGamma( \WCart^{\mathrm{HT}}, \mathscr{E}) \ar[r] \ar[d] & \RGamma( [\Spf(\Z_p) / (1+p\Z_p)^{\times}], \widetilde{\eta}^{\ast} \mathscr{E}) \ar[d] \\
\F_p \otimes^{L} \mathscr{E}_{\eta} \ar[r] & \F_p \otimes^{L} \RGamma( [\Spf(\Z_p) / (1+p\Z_p)^{\times}], \widetilde{\eta}^{\ast} \mathscr{E}) }
\end{equation}
in the $p$-complete derived $\infty$-category $\widehat{\calD}(\Z_p)$.
\end{proposition}

\begin{warning}
Proposition \ref{proposition:clean-statement-of-fiber} is false in the case $p=2$. See Remark \ref{remark:even-prime-case} below.
\end{warning}

\begin{remark}\label{remark:clean-statement-concrete}
In the situation of Proposition \ref{proposition:clean-statement-of-fiber}, let $u$ be a topological generator of the group $1+p\Z_p$.
Then we can rewrite (\ref{equation:clean-statement-of-fiber2}) as a diagram of complexes
$$\xymatrix@R=50pt@C=50pt{ \RGamma( \WCart^{\mathrm{HT}}, \mathscr{E}) \ar[r] \ar[d] & \mathscr{E}_{\eta}^{u=1} \ar[d] \\
\F_p \otimes^{L} \mathscr{E}_{\eta} \ar[r] & \F_p \otimes^{L} \mathscr{E}_{\eta}^{u=1}. }$$
\end{remark}

\begin{proof}[Proof of Proposition \ref{proposition:fiber-sequence-global-sections} from Proposition \ref{proposition:clean-statement-of-fiber}]
Let $p$ be an odd prime, let $u$ be a topological generator of the group $(1+p\Z_p)^{\times}$, let 
$\mathscr{E}$ be a quasi-coherent complex on $\WCart^{\mathrm{HT}}$, and let $K$ denote the cofiber of the map
$\F_p \otimes^{L} \mathscr{E}_{\eta} \rightarrow \F_p \otimes^{L} \mathscr{E}_{\eta}^{u=1}$. Using Proposition \ref{proposition:clean-statement-of-fiber}
(and Remark \ref{remark:clean-statement-concrete}), we obtain a commutative diagram of fiber sequences
$$ \xymatrix@R=50pt@C=50pt{  \RGamma( \WCart^{\mathrm{HT}}, \mathscr{E}) \ar[r] \ar[d] & \mathscr{E}_{\eta}^{u=1} \ar[r] \ar[d] & K \ar[d] \\
\mathscr{E}_{\eta} \ar@{=}[r] \ar[d]^{\xi_u} & \mathscr{E}_{\eta} \ar[r] \ar[d]^{\gamma_u - \id} & 0 \ar[d] \\
\mathscr{E}_{\eta} \ar[r]^-{p} &  \mathscr{E}_{\eta} \ar[r] & \F_p \otimes^{L} \mathscr{E}_{\eta}. }$$
\end{proof}

\begin{proof}[Proof of Proposition \ref{proposition:clean-statement-of-fiber}]
For every quasi-coherent complex $\mathscr{E}$ on $\WCart^{\mathrm{HT}}$, 
Proposition \ref{proposition:O-inside-regular} supplies a fiber sequence
$$ \mathscr{E} \rightarrow \eta_{\ast}( \mathscr{E}_{\eta} ) \rightarrow \eta_{\ast}( \mathscr{E}_{\eta} ).$$
It will therefore suffice to show that the conclusion of Proposition \ref{proposition:clean-statement-of-fiber} holds in the special case
where $\mathscr{E}$ has the form $M \otimes \eta_{\ast}( \calO_{\Spf(\Z_p)} )$, for some object $M \in \widehat{\calD}(\Z_p)$.
Using Corollary \ref{corollary:HT-computation} we can reduce further to the case $M = \Z_p$, so that $\mathscr{E} \simeq \eta_{\ast}( \calO_{ \Spf(\Z_p)} )$.
In this case, we can identify $\mathscr{E}_{\eta}$ with the $p$-completed coordinate ring $\widehat{\calO}_{\mathbf{G}_{m}}^{\sharp}$ of the
group scheme $\mathbf{G}_{m}^{\sharp}$ (Theorem \ref{theorem:describe-HT}), and $\RGamma( \WCart^{\mathrm{HT}}, \mathscr{E} )$ with
the subalgebra $\Z_p \subset \widehat{\calO}_{\mathbf{G}_{m}}^{\sharp}$.
Let $u$ be a topological generator of the group $(1+p\Z_p)^{\times}$. Unwinding the definitions (as in the proof of Proposition \ref{proposition:fiber-sequence-global-sections}), we are reduced to showing that the sequence of abelian groups
$$ 0 \rightarrow \Z_p \rightarrow \widehat{\calO}_{\mathbf{G}_{m}}^{\sharp} \xrightarrow{ \frac{ \gamma_u - \id}{p} } \widehat{\calO}_{\mathbf{G}_{m}}^{\sharp} \rightarrow 0$$
is exact, where $\gamma_{u}$ denotes the automorphism of $\widehat{\calO}_{\mathbf{G}_{m}}^{\sharp}$ given by the construction $f(t) \mapsto f(ut)$.
Since $\widehat{\calO}_{\mathbf{G}_{m}}^{\sharp}$ is $p$-adically separated and complete, we can use Proposition
\ref{proposition:formula-for-action-simplified} to identify $\frac{ \gamma_{u} - \id}{p}$ with the formal power series
$$ \sum_{n \geq 1} \frac{ \log(u)^{n} }{p (n!)} \Theta^{n}_{\mathscr{E}} = \Theta_{\mathscr{E}} \cdot \sum_{n \geq 1} \frac{ \log(u)^{n} }{p (n!)} \Theta^{n-1}_{\mathscr{E}},$$
where $\Theta_{\mathscr{E}}$ is the Sen operator (given on $\widehat{\calO}_{\mathbf{G}_m}^{\sharp}$ by the differential operator $t \frac{ \partial}{ \partial t}$; see Example \ref{example:regular-representation}).
Note that the second factor is invertible (it is congruent modulo $p$ to the scalar $\log(u)/p$, which is invertible by virtue of our assumption that $u$ is a topological generator of $(1+p\Z_p)^{\times}$). We are therefore reduced to showing that the sequence
$$ 0 \rightarrow \Z_p \rightarrow  \widehat{\calO}_{\mathbf{G}_{m}}^{\sharp} \xrightarrow{ t \frac{ \partial}{\partial t} } \widehat{\calO}_{\mathbf{G}_{m}}^{\sharp} \rightarrow 0$$
is exact, which was established in the proof of Proposition \ref{proposition:O-inside-regular}.
\end{proof}

\begin{corollary}\label{corollary:presentation-of-HT}
Let $p$ be an odd prime. Then the commutative diagram (\ref{equation:clean-statement-of-fiber})
induces a fully faithful functor of $\infty$-categories
$$ F: \calD( \WCart^{\mathrm{HT}} ) \rightarrow \calD( [\Spf(\Z_p) / (1+p\Z_p)^{\times} )]
\times_{ \calD( [\Spec(\F_p) / (1+p\Z_p)^{\times}] ) } \calD( \Spec(\F_p) ).$$
\end{corollary}

\begin{remark}
Stated more informally, Corollary \ref{corollary:presentation-of-HT} asserts that a quasi-coherent complex $\mathscr{E}$ on
$\WCart^{\mathrm{HT}}$ can be recovered from the action of the profinite group $(1+ p \Z_p)^{\times}$ on the fiber
$\mathscr{E}_{\eta}$, together with the trivialization of this action modulo $p$. Alternatively, $\mathscr{E}$ can be recovered
from the complex $\mathscr{E}_{\eta}$ together with the endomorphism $\frac{ \gamma_{u} - \id}{p}$, where $u$ is a topological
generator of $(1+ p \Z_p)^{\times}$ (compare with Theorem \ref{theorem:describe-HT}).
\end{remark}

\begin{proof}[Proof of Corollary \ref{corollary:presentation-of-HT}]
Let $\calC$ denote the fiber product $$\calD( [\Spf(\Z_p) / (1+p\Z_p)^{\times}] ) \times_{ \calD( [\Spec(\F_p) / (1+p\Z_p)^{\times}] ) } \calD( \Spec(\F_p) )$$
and let $\mathscr{E}$ and $\mathscr{E}'$ be quasi-coherent complexes on $\calD( \WCart^{\mathrm{HT} })$; we wish to show that the functor
$F$ induces a homotopy equivalence of mapping spaces
$$\Hom_{\calD( \WCart^{\mathrm{HT}} )}( \mathscr{E}', \mathscr{E} ) \rightarrow \Hom_{\calC}( F( \mathscr{E}'), F( \mathscr{E}) ).$$
By virtue of Proposition \ref{proposition:twists-generate}, we can assume without loss of generality that
$\mathscr{E}' = \calO_{ \WCart^{\mathrm{HT}}}\{n\}$ for some integer $n$. Replacing $\mathscr{E}$ by the twist
$\mathscr{E}\{-n\}$, we can reduce to the case $n=0$, in which case the desired result follows from
Proposition \ref{proposition:clean-statement-of-fiber}.
\end{proof}

\begin{warning}\label{warning:not-essentially-surjective}
The functor $F$ of Corollary \ref{corollary:presentation-of-HT} is not essentially surjective. For example,
let $\mathscr{E}$ be an object of the category
$$ \Vect( [\Spf(\Z_p) / (1+p\Z_p)^{\times}] )
\times_{ \Vect( [\Spec(\F_p) / (1+p\Z_p)^{\times}] ) } \Vect( \Spec(\F_p) )$$
of {\em vector bundles} on the pushout 
$$ [\Spf( \Z_p ) / (1+p\Z_p)^{\times}] \coprod_{ [\Spec(\F_p) / (1+p\Z_p)^{\times}]} \Spec(\F_p).$$
Choosing a topological generator $u$ for the group $(1+p\Z_p)^{\times}$, we can identify $\mathscr{E}$ with
pair $(M, \gamma_{u} )$, where $M$ is a projective $\Z_p$-module of finite rank and $\gamma_{u}$ is an
automorphism of $M$ which is congruent to the identity modulo $p$. Let $\{ \lambda_i \}_{i \in I}$ be the collection of all eigenvalues for
the induced action of $\gamma_{u}$ on the vector space $\overline{\Q}_p \otimes_{\Z_p} M$; our assumption then guarantees that
each $\lambda_{i}$ is congruent to $1$ modulo $p$ as an element of the valuation ring $\calO_{ \overline{\Q}_p}$.
Using Theorem \ref{theorem:compute-with-HT} and Proposition \ref{proposition:formula-for-action-simplified}, it is not difficult to see that $\mathscr{E}$ belongs to the essential image of $F$ if and only if each $\lambda_{i}$ is congruent to an integer modulo $p^2 \calO_{ \overline{\Q}_p }$.
\end{warning}

\begin{remark}\label{remark:even-prime-case}
Suppose that $p=2$, so that the group $\mathbf{G}_{m}^{\sharp}( \Z_p ) = (1 + p \Z_p)^{\times}$ factors as a direct product
$\langle \pm 1 \rangle \times (1 + 4 \Z_p)^{\times}$. In this case, the action of the profinite group $(1+p \Z_p)^{\times}$ on
$\eta: \Spf(\Z_p) \rightarrow \WCart^{\mathrm{HT}}$ extends to an action of the group scheme
$\mu_{2} \times (1+4\Z_p)^{\times}$, and therefore determines a map of stacks
$$ \widetilde{\eta}: [\Spf(\Z_2) / ( \mu_2 \times (1+4\Z_2)^{\times} )] \rightarrow \WCart^{\mathrm{HT}}.$$
Since the restriction map $\mathbf{G}_{m}^{\sharp}( \Z_2 ) \rightarrow \mathbf{G}_m^{\sharp}( \Z / 4\Z)$ vanishes on
$(1+4\Z_2)^{\times}$, we can extend $\widetilde{\eta}$ to a commutative diagram of stacks
\begin{equation}\label{equation:clean-statement-of-fiber3} 
\xymatrix@R=50pt@C=50pt{ \WCart^{\mathrm{HT} } & [\Spf(\Z_2) / (\mu_2 \times (1+4\Z_2)^{\times})] \ar[l]_-{\widetilde{\eta}}  \\
[\Spec(\Z / 4 \Z) / \mu_2] \ar[u] & [\Spec( \Z / 4 \Z) / (\mu_2 \times (1+4\Z_2)^{\times})]. \ar[u] \ar[l] }
\end{equation}
Proposition \ref{proposition:fiber-sequence-global-sections}, Proposition \ref{proposition:clean-statement-of-fiber}, and  
Corollary \ref{corollary:presentation-of-HT} have analogues in the case $p=2$, where we use (\ref{equation:clean-statement-of-fiber3})
as a replacement for the diagram (\ref{equation:clean-statement-of-fiber}). For example, if $\mathscr{E}$ is a quasi-coherent complex
on the Hodge-Tate divisor $\WCart^{\mathrm{HT}}$, then there is a canonical fiber sequence
$$ \RGamma( \WCart^{\mathrm{HT}}, \mathscr{E} )
\rightarrow \mathscr{E}_{\eta}^{+} \xrightarrow{ \frac{ \gamma_{u} - \id}{4} } \mathscr{E}_{\eta}^{+},$$
where $\mathscr{E}_{\eta}^{+}$ denotes the direct summand of $\mathscr{E}_{\eta}$ on which the Sen operator $\Theta_{\mathscr{E}}$ is
topologically nilpotent, and $u$ is a topological generator of the group $(1+4\Z_2)^{\times}$.
\end{remark}

\subsection{Comparison with the \texorpdfstring{$q$}{q}-de Rham Prism}\label{subsection:computing-with-WCart}

Let $\mathscr{E}$ be a quasi-coherent complex on the Cartier-Witt stack. It follows from Corollary \ref{corollary:global-sections-on-WCart} that the complex $\RGamma( \WCart, \mathscr{E} )$ can be realized as the totalization of a cosimplicial complex defined over the cosimplicial $\delta$-ring
$A^{\bullet}$ of Notation \ref{notation:simplicial-prism}. While this description is explicit in some sense, it is not very practical: for example, each term of the cosimplicial $\delta$-ring $A^{\bullet}$ is non-Noetherian. In this section, we will provide a much more concrete description of the complex $\RGamma( \WCart, \mathscr{E})$, at least when $p$ is an odd prime, which can be regarded as a ``$q$-deformation'' of Proposition \ref{proposition:clean-statement-of-fiber} (Theorem \ref{theorem:qdR-WCart}). 

\begin{notation}\label{notation:qdR-presentation}
Let $( \Z_p[[q-1]], ( [p]_q ) )$ denote the $q$-de Rham prism of Example \ref{example:q-prism}. Applying Construction \ref{construction:point-of-prismatic-stack}, we obtain a morphism of stacks $\rho_{\qdR}: \Spf( \Z_p[[q-1]] ) \rightarrow \WCart$. Concretely, to each $\Z_p[[q-1]]$-algebra $R$ in which the elements $p$ and
$(q-1)$ are nilpotent, the map $\rho_{\qdR}$ associates the Cartier-Witt divisor for $R$ given by the map $$W(R) \xrightarrow{ 1 + [q] + \cdots + [q]^{p-1} } W(R),$$
where $[q]$ denotes the Teichm\"{u}ller representative of $q$ in $W(R)$.
\end{notation}

Note that, since the $q$-de Rham prism is transversal (and nonzero), the morphism $\rho_{\qdR}$ is faithfully flat (Corollary \ref{corollary:flatness-v-transversality}).
We would like to use this observation to present $\WCart$ as a quotient of the affine formal scheme $\Spf(\Z_p[[q-1]])$. Note that the
$q$-de Rham prism $( \Z_p[[q-1]], ( [p]_q ) )$ carries an action of the group of $p$-adic units $\Z_p^{\times}$, where an element $u \in \Z_p^{\times}$
acts by the automorphism $q \mapsto q^{u}$. It follows by naturality that we can promote $\rho_{\qdR}: \Spf( \Z_p[[q-1]] ) \rightarrow \WCart$
to a $\Z_p^{\times}$-equivariant morphism of stacks (where the group $\Z_p^{\times}$ acts trivially on $\WCart$). The resulting action of
$\Z_p^{\times}$ on $\rho_{\qdR}$ is continuous in the following sense:

\begin{proposition}\label{proposition:continuous-action}
The preceding construction determines a morphism of stacks
$$ [\Spf( \Z_p[[q-1]] ) / \Z_p^{\times}] \rightarrow \WCart,$$
where $[\Spf( \Z_p[[q-1]] ) / \Z_p^{\times}]$ denotes the stack-theoretic quotient of $\Spf( \Z_p[[q-1]] )$ by the profinite group $\Z_p^{\times}$ (regarded as a pro-(finite \'etale) group scheme). 
\end{proposition}

\begin{proof}
Let $R$ be a $\Z_p[[q-1]]$-algebra in which the element $p$ and $(q-1)$ are nilpotent, and let $G \subseteq \Z_p^{\times}$ be an open subgroup
which stabilizes the associated $R$-valued point of $\Spf( \Z_p[[q-1]] )$ (so that $q^{u} = q$ in $R$ for each $u \in G$). It follows that
the group $G$ acts on the Cartier-Witt divisor $W(R) \xrightarrow{ 1 + [q] + \cdots + [q]^{p-1} } W(R)$ described in Notation \ref{notation:qdR-presentation}.
To complete the proof, we must show that this action is trivial when restricted to an open subgroup of $G$. In fact, it is trivial
on image of the map
$$ G \xrightarrow{ x \mapsto x^{p^n}} G$$
for $n \gg 0$; see Corollary \ref{corollary:automorphism-of-CW}.
\end{proof}

In what follows, we will abuse notation by denoting the morphism $[\Spf( \Z_p[[q-1]] ) / \Z_p^{\times}] \rightarrow \WCart$ of
Proposition \ref{proposition:continuous-action} also by $\rho_{\qdR}$. This map is not quite an isomorphism. Note that the composite map
$$ \Spf(\Z_p) \xrightarrow{q \mapsto 1} \Spf( \Z_p[[q-1]] ) \xrightarrow{ \rho_{\qdR} } \WCart$$
can be identified with de Rham point $\rho_{\dR}$ of Example \ref{example:de-Rham-point}. Since the action of $\Z_p^{\times}$ on $( \Z_p[[q-1]], ( [p]_q ) )$ is trivial modulo $(q-1)$,
we obtain a commutative diagram
\begin{equation}\label{equation:qdR-WCart} 
\xymatrix@R=50pt@C=50pt{ \WCart & \Spf(\Z_p) \ar[l]_-{ \rho_{\dR}} \\
[ \Spf( \Z_p[[q-1]] ) / \Z_p^{\times}] \ar[u]_{ \rho_{\qdR }} & [\Spf(\Z_p) / \Z_p^{\times}]. \ar[u] \ar[l] }
\end{equation}
We can now formulate the main result of this section:

\begin{theorem}\label{theorem:qdR-WCart}
Let $p$ be an odd prime. Then the diagram (\ref{equation:qdR-WCart}) determines a fully faithful functor
of $\infty$-categories
$$ \calD( \WCart ) \rightarrow \calD( [\Spf( \Z_p[[q-1]] ) / \Z_p^{\times}] ) \times_{ \calD(  [\Spf(\Z_p) / \Z_p^{\times}] )} 
\calD( \Spf(\Z_p) ).$$
In particular, if $\mathscr{E}$ is a quasi-coherent complex on $\WCart$, then the diagram
\begin{equation}\label{equation:qdR-WCart-global} 
\xymatrix@R=50pt@C=50pt{ \RGamma(\WCart, \mathscr{E}) \ar[r] \ar[d] & \RGamma( \Spf(\Z_p), \rho_{\dR}^{\ast} \mathscr{E} ) \ar[d] \\
\RGamma( \Spf( [\Z_p[[q-1]] ) / \Z_p^{\times}], \rho_{\qdR}^{\ast} \mathscr{E} ) \ar[r] & \RGamma(  [\Spf(\Z_p) / \Z_p^{\times}], \rho_{\dR}^{\ast} \mathscr{E} ) }
\end{equation}
is a pullback square in the $\infty$-category $\widehat{\calD}(\Z_p)$.
\end{theorem}

\begin{warning}
Theorem \ref{theorem:qdR-WCart} is false in the case $p=2$ (see Remark \ref{remark:even-prime-case}). 
\end{warning}

\begin{warning}
The fully faithful functor
$$\calD( \WCart ) \rightarrow \calD( [\Spf( \Z_p[[q-1]] ) / \Z_p^{\times} ]) \times_{ \calD([ \Spf(\Z_p) / \Z_p^{\times}] )} 
\calD( \Spf(\Z_p) )$$
of Theorem \ref{theorem:qdR-WCart} is not essentially surjective (Warning \ref{warning:not-essentially-surjective}).
\end{warning}

It will be convenient to give a slightly more concrete formulation of Theorem \ref{theorem:qdR-WCart}.
Note that the short exact sequence
$$ 0 \rightarrow (1 + p \Z_p)^{\times} \rightarrow \Z_p^{\times} \rightarrow \F_p^{\times} \rightarrow 0$$
has a unique splitting, given by the homomorphism
$$ \F_p^{\times} \rightarrow \Z_p^{\times} \quad \quad e \mapsto [e].$$
This splitting determines an action of the group $\F_p^{\times}$ on $\Z_p[[q-1]]$, which carries each element
$e \in \F_p^{\times}$ to the automorphism of $\Z_p[[q-1]]$ given by $q \mapsto q^{[e]}$. We let $\slashp$ denote the sum
$$ \sum_{e \in \F_p} q^{[e]} = 1 + \sum_{e \in \F_p^{\times}} q^{[e]},$$
formed in the power series ring $\Z_p[[q-1]]$. By construction, the element $\slashp$ is invariant under the action of $\F_p^{\times}$.

\begin{proposition}\label{runnup}
The element $\slashp$ is a generator of the principal ideal $( [p]_q ) \subseteq \Z_p[[q-1]]$.
\end{proposition}

\begin{proof}
For each integer $0 \leq e < p$, let $[e] \in \Z_p$ denote the Teichm\"{u}ller representative of the image of $e$ in $\F_p^{\times}$.
Then $[e]$ is congruent to $e$ modulo $p$ (as an element of $\Z_p$), and therefore $q^{[e]}$ is congruent to $q^{e}$ modulo
$q^{p}-1$ (as an element of the ring $\Z_p[[q-1]]$). Summing over $e$, we obtain 
$\slashp = [p]_{q} + (q^{p}-1)f$ for some element $f \in \Z_p[[q-1]]$. It follows that $\slashp$ is the product of $[p]_{q}$ with the invertible
element $1 + (q-1)f \in \Z_p[[q-1]]$.
\end{proof}

\begin{example}\label{example:slashed-p}
For $p=2$, we have $\slashp = 1 + q = [p]_q$. For $p=3$, we have 
$$\slashp = q^{[0]} + q^{[1]} + q^{[2]} = 1 + q + q^{-1}  = q^{-1} [p]_{q}.$$
\end{example}

\begin{corollary}\label{corollary:reduced-q-de-Rham-structure}
There is a unique morphism of power series rings $f: \Z_p[[ x ]] \rightarrow \Z_p[[q-1]]$ satisfying $f(x) = \slashp$.
Moreover, $f$ is an isomorphism from $\Z_p[[x]]$ to the subring $\Z_p[[q-1]]^{\F_p^{\times}}$ of $\F_p^{\times}$-invariant elements of $\Z_p[[q-1]]$.
\end{corollary}

\begin{proof}
The existence and uniqueness of $f$ is clear (since $\slashp$ belongs to the maximal ideal of the complete local ring $\Z_p[[q-1]]$).
Let $R \subseteq \Z_p[[q-1]]$ denote the subring consisting of elements which are fixed by the action of $\F_p^{\times}$; we wish to show
that $f: \Z_p[[x]] \rightarrow R$ is an isomorphism. Since the source and target of $f$ are $(x)$-complete and $x$-torsion-free,
it will suffice to show that the induced map $\Z_p \rightarrow R/( \slashp )$ is an isomorphism. Note that the exact sequence
$$ 0 \rightarrow \slashp \Z_p[[q-1]] \rightarrow \Z_p[[q-1]] \rightarrow \Z_p[[ q-1]] / ( \slashp) \rightarrow 0$$
remains exact after taking $\F_p^{\times}$-fixed points (since the $\F_p^{\times}$ is a finite group of order relatively prime to $p$).
It follows that $R / (\slashp )$ can be identified with the ring of $\F_p^{\times}$-invariant elements of
the quotient $\Z_p[[q-1]] / ( \slashp ) \simeq \Z_p[[q-1]] / ( [p]_q )$ (see Proposition \ref{runnup}): that is, with the ring of
$\F_p^{\times}$-invariant elements of the discrete valuation ring $\Z_p[\zeta_p]$, which is evidently isomorphic to $\Z_p$.
\end{proof}

\begin{notation}\label{notation:reduced-q-de-Rham-prism}
Let $\Z_p[[ \slashp ]]$ denote the ring of $\F_p^{\times}$-invariant elements of $\Z_p[[q-1]]$
(by virtue of Corollary \ref{corollary:reduced-q-de-Rham-structure}, $\Z_p[[\slashp]]$ is isomorphic to a power series algebra over $\Z_p$ on one generator,
as suggested by the notation). Note that the $\delta$-ring structure on $\Z_p[[q-1]]$ restricts to a $\delta$-ring structure on the subring
$\Z_p[[ \slashp ]]$ (since $\F_p^{\times}$ acts by $\delta$-algebra automorphisms of $\Z_p[[q-1]]$). Moreover, it follows from
Proposition \ref{runnup} that the element $\slashp$ is a distinguished element of $\Z_p[[\slashp]]$: that is, the pair
$( \Z_p[[ \slashp ], (\slashp ) )$ is a prism. Note that the action of $\Z_p^{\times}$ on $\Z_p[[q-1]]$ restricts to an action of $\Z_p^{\times} \simeq \F_p^{\times} \times (1+p\Z_p)^{\times}$
on $\Z_p[[\slashp]]$, which is trivial on the first factor. In particular, the prism $( \Z_p[[\slashp]], (\slashp ) )$ inherits
an action of the $p$-profinite group $(1+p\Z_p)^{\times}$.
\end{notation}

Applying Construction \ref{construction:point-of-prismatic-stack} to the prism $(\Z_p[[ \slashp]], (\slashp ))$,
we obtain a morphism of stacks $\rho_{\Z_p[[\slashp]]}: \Spf( \Z_p[[ \slashp ]] ) \rightarrow \WCart$. Arguing as in
the proof of Proposition \ref{proposition:continuous-action}, we see that $\rho_{\Z_p[[\slashp]]}$ can be refined to a morphism
$$ [\Spf( \Z_p[[ \slashp ]] ) / (1+p\Z_p)^{\times}] \rightarrow \WCart,$$
which (by slight abuse of notation) we will also denote by $\rho_{\Z_p[[\slashp]]}$. Note that we can expand (\ref{equation:qdR-WCart}) to a diagram
\begin{equation}\label{equation:qdR-WCart-refined} 
\xymatrix@R=50pt@C=50pt{ \WCart & \Spf(\Z_p) \ar[l]_-{ \rho_{\dR}} \\
[\Spf( \Z_p[[\slashp]] ) / (1+p\Z_p)^{\times}] \ar[u]_{ \rho_{\Z_p[[\slashp]]} } & [\Spf(\Z_p) / (1+p\Z_p)^{\times}] \ar[u] \ar[l]_-{ \slashp \mapsto p} \\
[ \Spf( \Z_p[[q-1]] ) / \Z_p^{\times}] \ar[u] & [\Spf(\Z_p) / \Z_p^{\times}]. \ar[u] \ar[l] }
\end{equation}
Moreover, the bottom vertical maps induce fully faithful pullback functors
$$ \calD( [\Spf( \Z_p[[\slashp]] ) / (1+p\Z_p)^{\times}] ) \hookrightarrow
\calD(  [\Spf( \Z_p[[q-1]] ) / \Z_p^{\times}] )$$
$$\calD([ \Spf(\Z_p) / (1+p\Z_p)^{\times}] ) \hookrightarrow
\calD( [\Spf(\Z_p) / \Z_p^{\times}] ).$$
Consequently, Theorem \ref{theorem:qdR-WCart} can be reformulated as follows:

\begin{theorem}\label{theorem:qdR-WCart-refined}
Let $p$ be an odd prime. Then the diagram (\ref{equation:qdR-WCart-refined}) determines a fully faithful functor
of $\infty$-categories
$$ F: \calD( \WCart ) \rightarrow \calD( [\Spf( \Z_p[[\slashp]] ) / (1+p\Z_p)^{\times}] ) \times_{ \calD( [\Spf(\Z_p) / (1+p\Z_p)^{\times}] )} 
\calD( \Spf(\Z_p) ).$$
\end{theorem}

Theorem \ref{theorem:qdR-WCart-refined} is equivalent to the following {\it a priori} weaker result:

\begin{theorem}\label{theorem:qdR-WCart-refined-reformulation}
Let $p$ be an odd prime, and let $\mathscr{E}$ be a quasi-coherent complex on $\WCart$. Then the diagram
\begin{equation}\label{equation:qdR-WCart-refined-global} 
\xymatrix@R=50pt@C=50pt{ \RGamma(\WCart, \mathscr{E}) \ar[r] \ar[d] & \RGamma( \Spf(\Z_p), \rho_{\dR}^{\ast} \mathscr{E} ) \ar[d] \\
\RGamma( [\Spf( \Z_p[[\slashp ]] ) / (1+p\Z_p)^{\times}], \rho_{\Z_p[[\slashp]]}^{\ast} \mathscr{E} ) \ar[r] & \RGamma(  [\Spf(\Z_p) / (1+p\Z_p)^{\times}], \rho_{\dR}^{\ast} \mathscr{E} ) }
\end{equation}
is a pullback square in the $\infty$-category $\widehat{\calD}(\Z_p)$.
\end{theorem}

\begin{proof}[Proof of Theorem \ref{theorem:qdR-WCart-refined} from Theorem \ref{theorem:qdR-WCart-refined-reformulation}]
Let $\calC$ denote the fiber product 
$$\calD( [\Spf( \Z_p[[\slashp]] ) / (1+p\Z_p)^{\times} )] \times_{ \calD(  [\Spf(\Z_p) / (1+p\Z_p)^{\times} )]} 
\calD( \Spf(\Z_p) )$$
and let $\mathscr{E}$ and $\mathscr{E}'$ be quasi-coherent complexes on the Cartier-Witt stack; we wish to show that the functor
$F$ induces a homotopy equivalence of mapping spaces
$$\Hom_{\calD( \WCart )}( \mathscr{E}', \mathscr{E} ) \rightarrow \Hom_{\calC}( F( \mathscr{E}'), F( \mathscr{E}) ).$$
By virtue of Corollary \ref{corollary:HT-ideals-generate}, we may assume without loss of generality that $\mathscr{E}' = \calI^{n}$
is a power of the Hodge-Tate ideal sheaf of Example \ref{example:HTI}. Replacing $\mathscr{E}$ by the twist
$\calI^{-n} \mathscr{E}$, we can reduce to the case $n=0$, in which case the desired result follows from
Theorem \ref{theorem:qdR-WCart-refined-reformulation}.
\end{proof}

We will deduce Theorem \ref{theorem:qdR-WCart-refined-reformulation} by combining Proposition \ref{proposition:clean-statement-of-fiber} with
the following:

\begin{proposition}\label{proposition:old-diagram}
Let $\rho_{\Z_p[[\slashp]]}^{\mathrm{HT}}: \Spf(\Z_p) \rightarrow \WCart^{\mathrm{HT}}$ be the morphism obtained
by applying Remark \ref{remark:HT-point-of-prismatic-stack} to the prism $( \Z_p[[ \slashp]], ( \slashp ) )$ of Notation \ref{notation:reduced-q-de-Rham-prism}.
Then $\rho_{\Z_p[[\slashp]]}^{\mathrm{HT}}$ is isomorphic to the morphism $\eta: \Spf(\Z_p) \rightarrow \WCart^{\mathrm{HT}}$ of
Construction \ref{construction:fiber-at-eta}. Moreover, any choice of isomorphism is equivariant with respect to the action
of the group $(1+ p\Z_p)^{\times}$ (which acts on $\eta$ via the isomorphism $\Aut(\eta)(\Z_p) = \mathbf{G}_{m}^{\sharp}(\Z_p) \simeq (1+p\Z_p)^{\times}$
of Theorem \ref{theorem:describe-HT}).
\end{proposition}

\begin{proof}
Let us identify the quotient $\Z_p[[ q-1]] / ( \slashp )$ with the ring $\Z_p[ \zeta_p ]$ obtained from $\Z_p$ by adjoining a primitive $p$th roof of unity
$\zeta_p$ (given by the image of $q$). Note that the quotient map $\overline{\xi}: \Z_p[[ q-1]] \twoheadrightarrow \Z_p[ \zeta_p ]$ lifts uniquely to a morphism
of $\delta$-rings $\xi: \Z_p[[q-1]] \rightarrow W( \Z_p[ \zeta_p] )$, which carries $q$ to the Teichm\"{u}ller representative $[ \zeta_p ]$.
It follows that
$$ \xi( \slashp ) = f( \sum_{e \in \F_p } q^{[e]} ) = \sum_{e \in \F_p} [ \zeta_{p}^{e} ] \in W(\Z_p[ \zeta_p] ).$$
Note that $\xi(\slashp)$ belongs to the kernel of the restriction map $W(\Z_p[\zeta_p]) \twoheadrightarrow \Z_p[\zeta_p]$, 
and can therefore be written as $V(x)$ for some unique element  $x \in W(\Z_p[\zeta_p])$. We then compute
$$ px = F(V(x)) = F( \xi( \slashp) ) = \sum_{e \in \F_p} [ \zeta_{p}^{pe} ] = p.$$
Since $p$ is not a zero divisor in $W( \Z_p[ \zeta_p] )$, it follows that $x =1$: that is, the diagram of $\Z_p[[q-1]]$-modules
$$ \xymatrix@R=50pt@C=50pt{ {\slashp} \Z_p[[q-1]] \ar[r]^-{ \slashp \mapsto 1} \ar[d] & W( \Z_p[ \zeta_p] ) \ar[d]^{ V(1) } \\
\Z_p[[q-1]] \ar[r]^-{ \xi } & W( \Z_p[ \zeta_p] ) }$$
is commutative. Passing to $\F_p^{\times}$-fixed points, we obtain a commutative diagram of $\Z_p[[\slashp]]$-modules
$$ \xymatrix@R=50pt@C=50pt{ {\slashp} \Z_p[[\slashp]] \ar[d] \ar[r]^-{\slashp \mapsto 1} & W( \Z_p ) \ar[d]^{ V(1) } \\
\Z_p[[\slashp]] \ar[r]^-{ \xi_0 } & W( \Z_p), }$$
which determines an isomorphism of $\rho_{\Z_p[[\slashp]]}^{\mathrm{HT}}$ with $\eta$ (see Example \ref{example:torsor-over-Hodge-Tate}).
We will complete the proof by showing that this isomorphism is equivariant with respect to the group $(1+p\Z_p)^{\times}$
(it then follows formally that the same result holds for any other isomorphism of $\rho_{\Z_p[[\slashp]]}^{\mathrm{HT}}$ with $\eta$, 
since the automorphism group $\Aut(\eta) \simeq \mathbf{G}_{m}^{\sharp}$ is abelian). 

Let $u$ be an element of the group $(1+p\Z_p)^{\times}$, and let $\gamma_{u}$ denote the corresponding automorphism
of $\Z_p[[ q-1]]$ (given by $\gamma_{u}(q) = q^{u}$). Unwinding the definitions, we see that the image of $u$
in $\Aut(\eta)(\Z_p) \simeq \mathbf{G}_{m}^{\sharp}(\Z_p) \subseteq \Z_p^{\times}$ can be identified with the image
of $\frac{ \gamma_{u}( \slashp )}{ \slashp }$ under the composite map
$$ \Z_p[[\slashp]] \xrightarrow{\xi_0} W(\Z_p) \twoheadrightarrow \Z_p \quad \quad \slashp \mapsto 0,$$
given by the restriction of $\overline{\xi}$ to $\F_p^{\times}$-invariant elements of $\Z_p[[q-1]]$.
Identifying $\F_p^{\times}$ with the set of integers $\{ 0 < 1 < \cdots < p-1 \}$, we compute
\begin{eqnarray*} \gamma_{u}( \slashp) & = & 
\sum_{0 \leq e < p} q^{u [e] } \\
& = & [p]_{q} + \sum_{0 \leq e < p } q^{u[e]} - q^{e} \\
& = & [p]_{q} (1 + (q-1) \sum_{0 \leq e < p} q^{e} \frac{ q^{u[e] - e} }{ q^{p} - 1} ).
\end{eqnarray*}
Dividing by $[p]_{q}$ and applying the homomorphism $\overline{\xi}$, we obtain
\begin{eqnarray*}
\overline{\xi}( \frac{ \gamma_{u}( \slashp) }{ [p]_q } ) & = & \overline{\xi}( 1 + (q-1) \sum_{0 \leq e < p} q^{e} \frac{ q^{u[e] - e} }{ q^{p} - 1} ) \\
& = & 1 + (\zeta_p-1) \sum_{0 \leq e < p} \zeta_{p}^{e} \overline{\xi}( \frac{ q^{u[e] - e} }{ q^{p} - 1}) \\
& = & 1 + (\zeta_p-1) \sum_{0 \leq e < p} \zeta_p^{e} \frac{ u[e] - e}{p} \\
& = & u \frac{ \zeta_p-1}{p} \sum_{0 \leq e < p} \zeta_p^{e} [e];
\end{eqnarray*}
here the last equality follows from the identity $p = (\zeta_p -1) \sum_{0 \leq e < p} e \zeta_p^{e}$. We therefore compute
$$ \overline{\xi}(\frac{ \gamma_u( \slashp )}{ \slashp } )
 = \overline{\xi}( \frac{ \gamma_u(\slashp) }{ [p]_q } ) \overline{\xi}( \frac{ \gamma_1(\slashp)}{ [p]_q} )^{-1}
 = u,$$
as desired.
\end{proof}

\begin{warning}
The isomorphism of Proposition \ref{proposition:old-diagram} is not unique: it is ambiguous up to the action of the automorphism
group $\Aut(\eta)(\Z_p) \simeq \mathbf{G}_m^{\sharp}(\Z_p) = (1+p\Z_p)^{\times}$.
\end{warning}

\begin{proof}[Proof of Theorem \ref{theorem:qdR-WCart-refined-reformulation}]
Note that every quasi-coherent complex $\mathscr{E}$ on $\WCart$ is complete with respect to the Hodge-Tate ideal sheaf
$\mathscr{I}$: that is, $\mathscr{E}$ can be identified with the inverse limit of the tower $\varprojlim_{n \geq 0} \mathscr{E} / \calI^{n} \mathscr{E}$.
It will therefore suffice to prove that the diagram (\ref{equation:qdR-WCart-refined-global} ) is a pullback square after
replacing $\mathscr{E}$ by each quotient $\mathscr{E} / \calI^{n} \mathscr{E}$. Proceeding by induction on $n$, we are reduced to proving
Theorem \ref{theorem:qdR-WCart-refined-reformulation} in the special case where $\mathscr{E}$ is the direct image of a quasi-coherent complex $\mathscr{F}$ on
the Hodge-Tate divisor $\WCart^{\mathrm{HT}}$. Passing to the inverse image of the Hodge-Tate divisor, (\ref{equation:qdR-WCart-refined})
yields a diagram of stacks
\begin{equation}\label{equation:qdR-WCart-HT} 
\xymatrix@R=50pt@C=50pt{ \WCart^{\mathrm{HT}} & \Spec(\F_p) \ar[l]_-{ \rho^{\mathrm{HT}}_{\Z_p}} \\
[\Spf( \Z_p) / (1+p\Z_p)^{\times}] \ar[u]_{ \rho_{\Z_p[[\slashp]]}^{\mathrm{HT}} } & [\Spec(\F_p) / (1+p\Z_p)^{\times}], \ar[u] \ar[l] }
\end{equation}
and we wish to show that the induced diagram
$$\xymatrix@R=50pt@C=50pt{ \RGamma(\WCart^{\mathrm{HT}}, \mathscr{F}) \ar[r] \ar[d] & \RGamma( \Spec(\F_p), (\rho_{\dR}^{\mathrm{HT}})^{\ast} \mathscr{F} ) \ar[d] \\
\RGamma( [\Spf( \Z_p) / (1+p\Z_p)^{\times}], (\rho_{\Z_p[[\slashp]]}^{\mathrm{HT}})^{\ast} \mathscr{F} ) \ar[r] & \RGamma( [\Spec(\F_p) / (1+p\Z_p)^{\times}], 
(\rho_{\dR}^{\mathrm{HT}})^{\ast} \mathscr{F} ) }$$
is a pullback square in the $\infty$-category $\widehat{\calD}(\Z_p)$. By virtue of Proposition \ref{proposition:old-diagram}, we can identify
(\ref{equation:qdR-WCart-HT}) with the diagram (\ref{equation:clean-statement-of-fiber}), so that the desired result is a restatement of
Proposition \ref{proposition:clean-statement-of-fiber}.
\end{proof}

\subsection{Sen Theory}\label{subsection:Sen-theory}

Let $k$ be a perfect field of characteristic $p$, which we regard as fixed throughout this section. Let $K = W(k)[1/p]$ denote the fraction
field of $W(k)$, and choose an algebraic closure $\overline{K}$ of $K$. Let $C$ be the completion of $\overline{K}$ with respect to the unique valuation extending the $p$-adic valuation on $K$, and let $K_{\infty} = \bigcup_{n \geq 1} K( \zeta_{p^n} )$ denote the subfield of $\overline{K}$ generated by $K$ together with the $p$th power roots of unity. Let $\Gal( \overline{K}/K)$ denote the absolute Galois group of
$K$, so that $C$ carries an action of $\Gal( \overline{K}/K)$ which is continuous for the $p$-adic topology. Let 
$$\chi: \Gal( \overline{K}/K ) \twoheadrightarrow \Gal(K_{\infty}/K) \simeq \Z_p^{\times}$$
denote the cyclotomic character of $\Gal( \overline{K}/K)$. Our starting point is the following theorem of Sen (Theorem~4 of \cite{SenShankar1980CcaG}):

\begin{theorem}[Sen]\label{theorem:Sen-theory}
Let $V$ be a finite-dimensional vector space over $C$ equipped with a continuous semilinear action of $\Gal( K )$. Then there exists
a unique subgroup $V_{\infty} \subseteq V$ and endomorphism $\Theta: V_{\infty} \rightarrow V_{\infty}$ with the following properties:
\begin{itemize}
\item[$(1)$] The subgroup $V_{\infty}$ is a vector space over the subfield $K_{\infty} \subset C$, and the endomorphism $\Theta: V_{\infty} \rightarrow V_{\infty}$ is $K_{\infty}$-linear.

\item[$(2)$] The inclusion $V_{\infty} \hookrightarrow V$ induces a $C$-linear isomorphism $C \otimes_{ K_{\infty} } V_{\infty} \rightarrow V$; in particular, we have $\dim_{K_{\infty} }( V_{\infty} ) = \dim_{ C}(V) < \infty$.

\item[$(3)$] For each element $x \in V_{\infty}$, there exists an open subgroup $G \subseteq \Gal( \overline{K}/K)$ with the property that, for every element $g \in G$, we have $$ g(x) = \exp( \log( \chi(g) ) \cdot \Theta )(x).$$
\end{itemize}
\end{theorem}

In this section, we study a special class of semilinear Galois representations associated to perfect complexes on the product
$\WCart^{\mathrm{HT}}_{W(k)} = \Spf( W(k) ) \times \WCart^{\mathrm{HT}}$ (such complexes arise naturally when studying the Hodge-Tate cohomology of
formal schemes which are proper and smooth over $\Spf(W(k))$; see Remark \ref{remark:perfect-complex-absolute}). For such representations, we prove a
refinement of Theorem \ref{theorem:Sen-theory}, which describes the action of the entire Galois group $\Gal( \overline{K}/K)$ when $p$ is an odd prime (and subgroup of index $2$ in the case $p=2$).

\begin{construction}\label{construction:Galois-from-Hodge-Tate}
Let $\calO_{C}$ denote the valuation ring of $C$, and let $\Spf( \calO_{C} )$ denote its formal spectrum (with respect to the $p$-adic topology on $\calO_{C}$). Then $\calO_{C}$ is a perfectoid ring.
Combining the construction of Example \ref{example:points-from-perfectoid-rings} with the inclusion map $W(k) \hookrightarrow \calO_{C}$, we obtain a morphism of stacks
$\rho_{C}^{\mathrm{HT}}: \Spf( \calO_C ) \to \WCart^{\mathrm{HT}}_{W(k)}$. Let $\mathscr{E}$ be a perfect complex on $\WCart^{\mathrm{HT}}_{W(k)}$: that is, a dualizable object of
the derived $\infty$-category $\calD( \WCart^{\mathrm{HT}}_{W(k)} )$. For every integer $n$, we let $V^{n}( \mathscr{E} )$ denote the localization $\mathrm{H}^{n}( \rho_{\calO_{C}}^{\ast} \mathscr{E} )[1/p]$, which we regard as a finite-dimensional vector space over $C$. Note that the construction $\mathscr{E} \mapsto V^{n}( \mathscr{E} )$ depends functorially on the choice of algebraic closure
of $K$, and therefore carries semilinear action of the absolute Galois group $\Gal( K ) = \Gal( \overline{K} / K )$ (which is easily seen to be continuous with respect to the $p$-adic topology on $V^{n}( \mathscr{E} )$). 
\end{construction}

\begin{example}\label{example:representation-from-twist}
Let $(A,I)$ denote the perfect prism with quotient $A/I \simeq \calO_{C}$. Then there is a canonical $\Gal( \overline{K}/K)$-equivariant isomorphism
$$ V^{\ast}( \calO_{\WCart^{\mathrm{HT}}_{W(k)}}\{1\} ) \simeq \begin{cases} (I / I^2)[1/p] & \text{ if } \ast = 0 \\
0 & \text{ otherwise.} \end{cases}$$
More generally, for every integer $m$, the graded Galois representation $V^{\ast}( \calO_{\WCart^{\mathrm{HT}}_{W(k)} }\{m\} )$ can be identified with the Tate twist $C(m)$ (concentrated in degree zero).
\end{example}

\begin{remark}
\label{rmk:FCrysCrystalline}
Let $\mathscr{E}$ be a perfect complex on the stack $\WCart_{W(k)} = \Spf(W(k)) \times \WCart$ which is equipped with the structure of an {\it $F$-crystal}: that is, a pair of morphisms
$$ \mathscr{I}^{n} \otimes F^{\ast}( \mathscr{E} ) \rightarrow \mathscr{E} \quad \quad \mathscr{I}^{n} \otimes \mathscr{E} \rightarrow F^{\ast}(\mathscr{E})$$
for some $n \gg 0$ which are mutually inverse to one another (up to powers of the Hodge-Tate ideal sheaf $\mathscr{I}$). In this case, there is a canonical isomorphism
$$ V^{\ast}( \mathscr{E}|_{ \WCart_{W(k)}^{\mathrm{HT}}} ) \simeq C \otimes_{\Q_p} V_0^{\ast}( \mathscr{E} ),$$
where $V_0^{\ast}( \mathscr{E} )$ is a graded $\Q_p$-vector space equipped with a crystalline action of $\Gal(\overline{K}/K)$; we refer to \cite{bhatt2021prismatic} for more on the relation of $F$-crystals over $\WCart_{W(k)}$ and lattices in crystalline representations of $\Gal(\overline{K}/K)$.
\end{remark}

Let $\eta: \Spf(\Z_p) \rightarrow \WCart^{\mathrm{HT}}$ be the morphism of Construction \ref{construction:fiber-at-eta} and let $\eta_{W(k)}: \Spf(W(k)) \rightarrow \WCart^{\mathrm{HT}}_{W(k)}$ be the product of $\eta$ with the identity map on $\Spf(W(k))$. For every quasi-coherent complex $\mathscr{E}$ on $\WCart^{\mathrm{HT}}_{W(k)}$, we write $\mathscr{E}_{\eta}$ for the
pullback $\eta_{W(k)}^{\ast} \mathscr{E}$, which we identify with a $p$-complete complex of $W(k)$-modules. Our main result can now be stated as follows:

\begin{theorem}\label{theorem:Sen-theory-refined}
Let $\mathscr{E}$ be a perfect complex on $\WCart^{\mathrm{HT}}_{W(k)}$, and let $V^{\ast}(\mathscr{E})$ be the graded $C$-vector space of
Construction \ref{construction:Galois-from-Hodge-Tate}. Then there is a canonical isomorphism of graded $C$-vector spaces
$$\alpha: C \otimes_{W(k)} \mathrm{H}^{\ast}( \mathscr{E}_{\eta} ) \xrightarrow{\sim} V^{\ast}( \mathscr{E} )$$
with the following property:
\begin{itemize}
\item Let $g$ be an element of $\Gal( \overline{K} / K )$. Assume that, if $p=2$, then $\chi(g) \equiv 1 \pmod{4}$. Then
the diagram of graded abelian groups
\begin{equation}\label{equation:diagram-Sen-theory} \xymatrix@R=50pt@C=90pt{ \mathrm{H}^{\ast}( \mathscr{E}_{\eta} ) \ar[r]^-{ \exp( \log( \chi(g) ) \cdot \Theta_{\mathscr{E}} ) } \ar[d]^{\alpha}
& \mathrm{H}^{\ast}( \mathscr{E}_{\eta} ) \ar[d]^{\alpha} \\
V^{\ast}(\mathscr{E}) \ar[r]^-{g} & V^{\ast}( \mathscr{E} ) }
\end{equation}
is commutative, where $\Theta_{\mathscr{E}}$ is (induced by) the Sen operator of Notation \ref{notation:Sen-operator-second}.
\end{itemize}
\end{theorem}

\begin{remark}
Let $\mathscr{E}$ be a perfect complex on $\WCart^{\mathrm{HT}}_{W(k)}$, and let $V^{\ast}(\mathscr{E} )$ be the semilinear Galois representation of
Construction \ref{construction:Galois-from-Hodge-Tate}. It follows that the subspace $V^{\ast}( \mathscr{E} )_{\infty}$ appearing in the statement of Theorem \ref{theorem:Sen-theory}
coincides with the image of $K_{\infty} \otimes_{W(k)} \mathrm{H}^{\ast}( \mathscr{E}_{\eta} )$ under the isomorphism $\alpha$ of Theorem \ref{theorem:Sen-theory-refined}.
Moreover, the endomorphism $\Theta$ which appears in Theorem \ref{theorem:Sen-theory} can be identified with the $K_{\infty}$-linear extension of the
Sen operator $\Theta_{\mathscr{E}}$ on $\mathrm{H}^{\ast}( \mathscr{E}_{\eta} )$.
\end{remark}

Theorem \ref{theorem:Sen-theory-refined} is a consequence of a more refined statement which holds at the level of $\calO_{C}$-complexes and does not require the assumption that
$\mathscr{E}$ is perfect. Note that the morphism $\rho_{C}^{\mathrm{HT}}$ appearing in Construction \ref{construction:Galois-from-Hodge-Tate} descends to morphism of stacks
$[\Spf( \calO_{C} ) / \Gal( \overline{K} / K )] \rightarrow \WCart^{\mathrm{HT}}_{W(k)}$, which (by slight abuse of notation) we will also denote by $\rho_{C}^{\mathrm{HT}}$.
Let $\pi: \Z_p^{\times} \rightarrow (1+p\Z_p)^{\times}$ be the unique splitting of the exact sequence $0 \rightarrow  (1+p\Z_p)^{\times} \rightarrow
\Z_p^{\times} \rightarrow \F_p^{\times} \rightarrow 0$, so that the composite map
$$ \Gal( \overline{K} / K ) \xrightarrow{\chi} \Z_p^{\times} \xrightarrow{\pi} (1+p\Z_p)^{\times}$$
induces a map of quotient stacks $[\Spf( \calO_C) / \Gal( \overline{K} / K)] \rightarrow [\Spf( \Z_p ) / (1+p\Z_p)^{\times}]$. We will prove the following:

\begin{proposition}\label{proposition:commutative-diagram-of-stacks}
The diagram of stacks
$$ \xymatrix@R=50pt@C=50pt{ [\Spf( \calO_C ) / \Gal( \overline{K}/K)] \ar[r]^-{ \rho_{C}^{\mathrm{HT}} } \ar[d] & \WCart^{\mathrm{HT}}_{W(k)} \ar[d] \\
[\Spf( \Z_p ) / (1+p\Z_p)^{\times}] \ar[r]^-{\eta} & \WCart^{\mathrm{HT}} }$$
commutes (up to a canonical isomorphism).
\end{proposition}

\begin{proof}
Write the perfectoid ring $\calO_{C}$ as a quotient $A/I$, where $(A,I)$ is a perfect prism. Choosing a compatible system of primitive $p^{n}$th roots of unity in $\calO_{C}$,
we obtain a $\Gal( \overline{K} / K )$-equivariant morphism of prisms $( \Z_p[[q-1]], ( [p]_q ) ) \rightarrow (A,I)$, where $\Gal( \overline{K} / K )$ acts on
$\Z_p[[q-1]]$ via the cyclotomic character $q \mapsto q^{\chi(\bullet)}$. Restricting to the prism $( \Z_p[[ \slashp ]], (\slashp) )$ of
Notation \ref{notation:reduced-q-de-Rham-prism} and invoking the functoriality of Remark \ref{remark:HT-point-of-prismatic-stack}, we see that the diagram
$$ \xymatrix@R=50pt@C=50pt{ [\Spf( \calO_C ) / \Gal( \overline{K}/K)] \ar[r]^-{\rho_{C}^{\mathrm{HT}}} \ar[d] & \WCart^{\mathrm{HT}}_{W(k)} \ar[d] \\
[\Spf(\Z_p) / (1+p\Z_p)^{\times}] \ar[r]^-{ \rho_{ \Z_p[[\slashp]] }^{\mathrm{HT}} } & \WCart^{\mathrm{HT}} }$$
commutes up to canonical isomorphism. Proposition \ref{proposition:commutative-diagram-of-stacks} now follows from Proposition \ref{proposition:old-diagram}.
\end{proof}

\begin{proof}[Proof of Theorem \ref{theorem:Sen-theory-refined}]
Let $\mathscr{E}$ be a perfect complex on $\WCart_{W(k)}^{\mathrm{HT}}$. Proposition \ref{proposition:commutative-diagram-of-stacks}
then supplies a $\Gal( \overline{K}/K)$-equivariant isomorphism $\overline{\alpha}: \calO_{C} \otimes^{L}_{W(k)} \mathscr{E}_{\eta} \rightarrow (\rho_{C}^{\mathrm{HT}})^{\ast} \mathscr{E}$,
where $\Gal( \overline{K}/K)$ acts on $\mathscr{E}_{\eta}$ via the composite map 
$$\Gal( \overline{K} / K ) \xrightarrow{\chi} \Z_p^{\times} \xrightarrow{\pi} (1+p\Z_p)^{\times} = \mathbf{G}_{m}^{\sharp}(\Z_p) \simeq \Aut(\eta)(\Z_p).$$
Passing to cohomology and inverting $p$, we obtain an isomorphism of graded vector spaces
$\alpha: C \otimes_{W(k)} \mathrm{H}^{\ast}( \mathscr{E}_{\eta} ) \xrightarrow{\sim} V^{\ast}( \mathscr{E} )$. The commutativity
of the diagram (\ref{equation:diagram-Sen-theory}) follows by combining the Galois equivariance of $\overline{\alpha}$ with
Proposition \ref{proposition:formula-for-action-simplified}.
\end{proof}

\newpage \section{Absolute Prismatic Cohomology}\label{section:absolute-prismatic-cohomology}

Let $(A,I)$ be a bounded prism, and let $\mathfrak{X}$ be a $p$-adic formal scheme which is smooth over $\Spf(A/I)$. In \cite{prisms}, the first author
and Scholze introduced a complex of $A$-modules $\RGamma_{\Prism}(\mathfrak{X} / A)$, which we will refer to as the {\it prismatic complex of $\mathfrak{X}$ relative to $A$} and whose cohomology we will refer to as the {\it prismatic cohomology of $\mathfrak{X}$ relative to $A$}. As noted in \cite{prisms} (Remark~1.7), it is also possible to consider an {\em absolute} version of prismatic cohomology, which does not involve a choice of base prism $(A,I)$. Our goal in this section is to develop this absolute theory.

In the special case where $\mathfrak{X} = \Spf(R)$ is an affine formal $\overline{A}$-scheme, we denote the relative prismatic complex $\RGamma_{\Prism}( \mathfrak{X} / A)$ by $\Prism_{R/A}$. As in \cite{prisms}, it will be convenient to formally extend the definition of the complex $\Prism_{R/A}$ to the case where $R$ is
an arbitrary $p$-complete (animated) commutative $\overline{A}$-algebra (in fact, we will go slightly further and drop the assumption that $R$ is $p$-complete,
but this is just for notational convenience: see Corollary \ref{salmage}). In \S\ref{subsection:relative-prismatic-affine}, we review the construction of the complexes $\Prism_{R/A}$ and recall some of their essential properties. In \S\ref{subsection:relative-prismatic-global}, we extend this construction to arbitrary
formal $\overline{A}$-schemes $\mathfrak{X}$ by defining $\RGamma_{\Prism}( \mathfrak{X} / A)$ to be the homotopy inverse limit of the complexes
$\Prism_{R/A}$, indexed by the category of all maps $\Spec(R) \rightarrow \mathfrak{X}$ (Construction \ref{construction:prismatic-complex-of-scheme}).
In the case where $\mathfrak{X} = \Spf(R)$ is the formal spectrum of a $p$-complete $\overline{A}$-algebra $R$ with bounded $p$-power torsion, then this construction recovers the relative prismatic complex $\Prism_{R/A}$ (Corollary \ref{corollary:sashy}); this can be regarded as a ``$p$-adic continuity'' property of the functor $R \mapsto \Prism_{R/A}$. In \S\ref{subsection:relative-prismatic-site}, we show that in many cases,
the prismatic complex $\RGamma_{\Prism}( \mathfrak{X} / A)$ can be computed using the relative prismatic site $( \mathfrak{X} / A )_{\Prism}$ introduced in
\cite{prisms} (Theorem \ref{theorem:site-theoretic-equivalence}). In particular, this is true when $\mathfrak{X}$ is smooth over $\Spf( \overline{A} )$
(so that our notation is compatible with that of \cite{prisms}).

In \S\ref{subsection:absolute-prismatic}, we associate to each (animated) commutative ring $R$ a complex of abelian groups $\Prism_{R}$ (well-defined up to quasi-isomorphism), which we will refer to as the {\it absolute prismatic complex} of $R$ (Construction \ref{construction:absolute-prismatic-cohomology-general}). 
We will obtain $\Prism_{R}$ as the global sections of a quasi-coherent complex $\mathscr{H}_{\Prism}(R)$ on the Cartier-Witt stack $\WCart$,
which assembles the relative prismatic complexes $\Prism_{ (\overline{A} \otimes^{L} R) / A}$ as $(A,I)$ varies over all bounded prisms (Construction \ref{construction:prismatic-cohomology-sheaves}). As in the relative case, we can formally extend this construction to associate an absolute prismatic complex
$\RGamma_{\Prism}( \mathfrak{X} )$ to an arbitrary $p$-adic formal scheme $\mathfrak{X}$ (Construction \ref{construction:absolute-prismatic-complex-of-scheme}).
Under some mild assumptions, we show that the absolute prismatic complex $\RGamma_{\Prism}( \mathfrak{X} )$ can also be computed using an {\em absolute} version of the prismatic site introduced in \cite{prisms} (Theorem \ref{theorem:absolute-compute-with-stack}).
If $\mathfrak{X}$ is smooth over $\Spf(\Z_p)$ and the prime $p$ is odd, then $\RGamma_{\Prism}( \mathfrak{X} )$ has a simple description in terms of the
$q$-de Rham complex of \cite{scholzeq}, which we explain in \S\ref{subsection:concrete-qdR} (see Proposition \ref{proposition:concrete-qdR}
and Remark \ref{remark:concrete-qdR}).

Let $R$ be a commutative ring. In \S\ref{subsection:absolute-HT}, we introduce a complex of $R$-modules $\overline{\Prism}_{R}$ (well-defined up to quasi-isomorphism), which we will refer to as the {\it absolute Hodge-Tate complex of $R$} (Construction \ref{construction:absolute-HT}). There is a tautological map 
$\Prism_{R} \rightarrow \overline{\Prism}_{R}$, which behaves heuristically as if $\overline{\Prism}_{R}$ were obtained from $\Prism_{R}$
by reducing modulo a locally principal ideal (see Example \ref{example:filtered-derived-2}). Formally, we define $\overline{\Prism}_{R}$ as the global sections
of a quasi-coherent complex $\mathscr{H}_{\overline{\Prism}}(R) = \mathscr{H}_{\Prism}(R)|_{ \WCart^{\mathrm{HT}} }$ on the Hodge-Tate divisor
$\WCart^{\mathrm{HT}}$. By exploiting properties of the global section functor $\RGamma( \WCart^{\mathrm{HT}}, \bullet )$, we establish
cocontinuity properties of the functors $R \mapsto \overline{\Prism}_{R}$ and $R \mapsto \Prism_{R}$ (see Propositions \ref{proposition:HT-sifted-colimit}
and \ref{proposition:sokon}).

For some purposes, it is convenient to replace the absolute Hodge-Tate complex $\overline{\Prism}_{R}$ by a more fundamental invariant. 
Let $\eta: \Spf(\Z_p) \rightarrow \WCart^{\mathrm{HT}}$ be the point described in Construction \ref{construction:fiber-at-eta}.
The pullback $\eta^{\ast} \mathscr{H}_{\overline{\Prism}}(R)$ can be regarded as a complex of $R$-modules, which we will denote by
$\widehat{\Omega}_{R}^{\DHod}$ and refer to it as the {\it $p$-complete diffracted Hodge complex of $R$}. It is equipped with
an endomorphism $\Theta$ (the {\it Sen operator} of Notation \ref{notation:Sen-operator-second}), and the absolute Hodge-Tate complex
$\overline{\Prism}_{R}$ can be identified with the fiber $(\widehat{\Omega}_{R}^{\DHod})^{\Theta =0} = \fib( \Theta: \widehat{\Omega}_{R}^{\DHod} \rightarrow \widehat{\Omega}_{R}^{\DHod}$) (see Remark \ref{remark:diffracted-vs-prismatic2} for a slightly stronger statement). In \S\ref{subsection:diffracted-Hodge-integral}, we show that, as the prime $p$ varies, the $p$-complete diffracted Hodge complexes $\widehat{\Omega}_{R}^{\DHod}$ can be assembled into a single
complex $\Omega_{R}^{\DHod} \in \calD(R)$ which we refer to as the {\it diffracted Hodge complex of $R$} (Construction \ref{construction:diffracted-Hodge-integral}); we will see later that this complex arises naturally when when studying the ``decompleted'' analogues of the motivic filtration on topological Hochschild homology (see Proposition \ref{proposition:compute-gr-gr}).

\subsection{Relative Prismatic Cohomology: Affine Case}\label{subsection:relative-prismatic-affine}

Let $(A,I)$ be a bounded prism, which we regard as fixed throughout this section, and let $\overline{A}$ denote the quotient ring $A/I$.
Our goal is to provide a quick review of the theory of (derived) prismatic cohomology relative to $(A,I)$, following \cite{prisms}.

\begin{definition}[The Relative Prismatic Site]\label{definition:relative-prismatic-site}
Let $\mathfrak{X}$ be a bounded $p$-adic formal scheme over $\overline{A}$. We define a category $( \mathfrak{X} / A)_{\Prism}$ as follows:
\begin{itemize}
\item The objects of $( \mathfrak{X} / A)_{\Prism}$ are pairs $(B,v)$, where $B$ is a $\delta$-algebra over $A$
such that $(B,IB)$ is a prism, and $v: \Spf(B/IB) \rightarrow \mathfrak{X}$ is a morphism of formal $\overline{A}$-schemes.

\item A morphism from $(B,v)$ to $(B',v')$ in the category $( \mathfrak{X} /A)_{\Prism}$ is a homomorphism
$f: B \rightarrow B'$ of $\delta$-algebras over $A$ for which the composition
$\Spf(B'/IB') \rightarrow \Spf(B/IB) \xrightarrow{v} \mathfrak{X}$ is equal to $v'$.
\end{itemize}
We will refer to $( \mathfrak{X} / A)_{\Prism}$ as the {\it prismatic site of $\mathfrak{X}$ relative to the prism $(A,I)$}. In the special case where $\mathfrak{X}$ is the formal
completion of an affine scheme $\Spec(R)$, we denote the prismatic site $(\mathfrak{X}/A)_{\Prism}$ by $(R/A)_{\Prism}$.
\end{definition}

\begin{remark}\label{remark:flat-topology-on-relative-site}
Let $\mathfrak{X}$ be a $p$-adic formal scheme over $\overline{A}$. We define the {\it flat topology} on the relative prismatic site
$(\mathfrak{X}/A)^{\op}_{\Prism}$ to be the Grothendieck topology generated by those finite collections of morphisms
$$ \{ (B,v) \rightarrow (B_{s}, v_s) \}_{s \in S}$$
for which the induced ring homomorphism $B \rightarrow \prod_{s \in S} B_s$ is $(p,I)$-completely faithfully flat.
\end{remark}

\begin{construction}[Relative Prismatic Cohomology]\label{construction:relative-prismatic-cohomology}
Let $\CAlg^{\anim}_{\overline{A}}$ denote the $\infty$-category of animated $\overline{A}$-algebras, and let $\widehat{\calD}(A)$ denote the $(p,I)$-complete derived $\infty$-category of $A$. It follows from the universal property of Proposition \ref{proposition:universal-of-animated} that there is an essentially unique functor
$$ \Prism_{\bullet/A}: \CAlg^{\anim}_{\overline{A}} \rightarrow \widehat{\calD}(A) \quad \quad R \mapsto \Prism_{R/A}$$
with the following properties:
\begin{itemize}
\item[$(1)$] The functor $\Prism_{\bullet/A}$ commutes with sifted colimits, and is therefore a left Kan extension of its restriction
to the full subcategory $\Poly_{\overline{A}} \subseteq \CAlg^{\anim}_{ \overline{A}}$ spanned by the finitely generated polynomial rings over $\overline{A}$.

\item[$(2)$] Let $R$ be a finitely generated polynomial algebra over $\overline{A}$. Then $\Prism_{R/A}$ is given by the limit $\varprojlim_{ (B,v) \in (R/A)_{\Prism} } B$, formed in the $\infty$-category $\widehat{\calD}(A)$.
\end{itemize}
If $R$ is an animated commutative $\overline{A}$-algebra, we will refer to $\Prism_{R/A}$ as the {\it prismatic complex of $R$ relative to $A$}.
\end{construction}

\begin{remark}\label{remark:p-completion-irrelevant}
In \cite{prisms}, the relative prismatic complex $\Prism_{R/A}$ is defined under the assumption that the animated $\overline{A}$-algebra $R$
is $p$-complete. However, the additional generality of Construction \ref{construction:relative-prismatic-cohomology} is only illusory. If $R$ is any animated commutative
$\overline{A}$-algebra and $\widehat{R}$ is its $p$-completion, then the tautological map $R \mapsto \widehat{R}$ induces an isomorphism
$\Prism_{R/A} \rightarrow \Prism_{ \widehat{R} / A}$ in the $\infty$-category $\widehat{\calD}(A)$ (see Corollary \ref{salmage}).
\end{remark}

\begin{remark}[Change of Prism]\label{remark:change-of-prism}
Let $f: (A,I) \rightarrow (B,J)$ be a morphism of bounded prisms, and let $R$ be an animated commutative algebra over the quotient ring $\overline{A} = A/I$.
Set $\overline{B} = B/J$ and regard the derived tensor product $\overline{B} \otimes_{\overline{A}}^{L} R$ as an animated commutative $\overline{B}$-algebra.
Then there is a canonical isomorphism
$$ \rho: \Prism_{ (\overline{B} \otimes_{\overline{A}}^{L} R) / B } \simeq B \widehat{\otimes}^{L}_{A} \Prism_{R/A}$$
in the $\infty$-category $\widehat{\calD}(B)$, where $\widehat{\otimes}_{A}^{L}$ denotes the $(J,p)$-completed derived tensor product (see Construction~7.6 of \cite{prisms}). 
\end{remark}

\begin{variant}[The Relative Hodge-Tate Complex]\label{variant:relative-Hodge-Tate}
Let $R$ be an animated commutative $\overline{A}$-algebra. We let $\overline{\Prism}_{R/A}$ denote the derived tensor product $\overline{A} \otimes^{L}_{A} \Prism_{R/A}$,
which we regard as a commutative algebra object of the $p$-complete derived $\infty$-category $\widehat{\calD}( \overline{A} )$. We will refer to $\overline{\Prism}_{R/A}$ as the {\it Hodge-Tate complex of $R$ relative to $A$}.
\end{variant}

\begin{remark}[The Conjugate Filtration]\label{remark:derived-HT-filtration}
For every animated commutative $\overline{A}$-algebra $R$, the Hodge-Tate complex $\overline{\Prism}_{R/A}$ is equipped with an exhaustive increasing filtration
$$ \Fil_{0}^{\conj} \overline{\Prism}_{R/A} \rightarrow
\Fil_{1}^{\conj} \overline{\Prism}_{R/A} \rightarrow \Fil_{2}^{\conj} \overline{\Prism}_{R/A} \rightarrow \cdots,$$
which we will refer to as the {\it conjugate filtration}. This filtration is characterized by the requirement that it depends functorially on $R$ and satisfies the following conditions:
\begin{itemize}
\item[$(1)$] For each $n \geq 0$, the functor
$$ \Fil_{n}^{\conj} \overline{\Prism}_{\bullet/A}: \CAlg^{\anim}_{\overline{A}} \rightarrow \widehat{\calD}(A) \quad \quad R \mapsto \Fil_{n}^{\conj} \overline{\Prism}_{R/A}$$
commutes with sifted colimits, and is therefore a left Kan extension of its restriction
to the full subcategory $\Poly_{\overline{A}} \subseteq \CAlg^{\anim}_{ \overline{A}}$ spanned by the finitely generated polynomial rings over $\overline{A}$.

\item[$(2)$] When $R$ is a finitely generated polynomial ring over $\overline{A}$, the conjugate filtration coincides with the Postnikov filtration on the Hodge-Tate complex $\overline{\Prism}_{R/A}$.
That is, each $\Fil_{n}^{\conj} \overline{\Prism}_{R/A}$ can be identified with the cohomological truncation $\tau^{\leq n} \overline{\Prism}_{R/A}$.
\end{itemize}

Let $R$ be an animated commutative $\overline{A}$-algebra. By convention, we define $\Fil_{n}^{\conj} \overline{\Prism}_{R/A}$ for $n < 0$.
For each $n \in \Z$, we let $\gr_{n}^{\conj} \overline{\Prism}_{R/A}$ denote the cofiber of the map
$\Fil_{n-1}^{\conj} \overline{\Prism}_{R/A}  \rightarrow \Fil_{n}^{\conj} \overline{\Prism}_{R/A}$. Then Theorem~4.10 of \cite{prisms} supplies
a {\it Hodge-Tate comparison isomorphism}
$$ \xi_R: \gr_{n}^{\conj} \overline{\Prism}_{R/A} \simeq L \widehat{\Omega}^{n}_{ R / \overline{A} }[-n]\{-n\},$$
where $L \widehat{\Omega}^{n}_{R/\overline{A}}$ is the complex defined in Construction \ref{construction:exterior-of-cotangent}.
\end{remark}

\begin{remark}[K\"{u}nneth Formula]\label{remark:Kunneth-relative-prismatic}
For every animated commutative $\overline{A}$-algebra $R$, the prismatic complex $\Prism_{R/A}$ can be regarded as a commutative algebra object of the
$\infty$-category $\widehat{\calD}(A)$: that is, as an $E_{\infty}$-algebra over $A$ (which is complete with respect to the ideal $I + (p)$).
This commutative algebra structure depends functorially on $R$: that is, we can promote $\Prism_{\bullet/A}$ to a functor
$\CAlg^{\anim}_{\overline{A}} \rightarrow \CAlg( \widehat{\calD}(A) )$, which is characterized by the obvious analogue of properties $(1)$ and $(2)$ of Construction \ref{construction:relative-prismatic-cohomology}. Combining the Hodge-Tate comparison of Remark \ref{remark:derived-HT-filtration} with
Remark \ref{remark:kunneth-for-derived-Hodge}, we conclude that the functor
$$ \Prism_{\bullet/A}: \CAlg^{\anim}_{\overline{A}} \rightarrow \CAlg(\widehat{\calD}(A)) \quad \quad R \mapsto \Prism_{R/A}$$
commutes with {\em all} colimits. In particular, we have canonical isomorphisms
$$ \Prism_{ (R \otimes^{L}_{ \overline{A} } R') / A} \simeq \Prism_{R/A} \widehat{\otimes}^{L}_{A} \Prism_{R'/A},$$
where the (derived) tensor product on the right is completed with respect to the ideal $I + (p)$.
\end{remark}

\begin{remark}\label{remark:conjugate-equals-postnikov-conditions}
Let $R$ be an $\overline{A}$-algebra having bounded $p$-power torsion, and suppose that the completed cotangent complex
$L\widehat{\Omega}^{1}_{R/ \overline{A}}$ is $p$-completely flat over $R$ (this condition is satisfied, for example, if the $p$-completion $\widehat{R}$ is
$p$-completely smooth over $\overline{A}$). Then, for every integer $n$, the complex 
$\gr_{n}^{\conj} \overline{\Prism}_{R/A} \simeq L \widehat{\Omega}^{n}_{R/ \overline{A}}[-n]\{-n\}$ is concentrated in cohomological
degree $n$. It follows that the conjugate filtration on $\overline{\Prism}_{R/A}$ coincides with the Postnikov filtration
$\{ \tau^{\leq n} \overline{\Prism}_{R/A} \}_{n \in \Z}$.
\end{remark}

\begin{remark}\label{remark:HT-de-Rham}
Let $R$ be an animated commutative $\overline{A}$-algebra. For every integer $n$, the fiber sequence
$$ I^{n+1} / I^{n+2} \otimes^{L}_{A} \Prism_{R/A} \rightarrow I^{n} / I^{n+2} \otimes^{L}_{A} \Prism_{R/A} \rightarrow
I^{n} / I^{n+1} \otimes^{L}_{A} \Prism_{R/A}$$
determines a boundary map
$$ \beta: \mathrm{H}^{n}( \overline{\Prism}_{R/A}\{n\} ) \rightarrow \mathrm{H}^{n+1}( \overline{\Prism}_{R/A}\{n+1\}).$$
When $R$ is $p$-completely smooth over $\overline{A}$, the boundary map $\beta$ fits into a commutative diagram
$$ \xymatrix@R=50pt@C=50pt{ \mathrm{H}^{n}( \overline{\Prism}_{R/A}\{n\} ) \ar[r]^-{\beta} \ar[d]^{\sim} &  \mathrm{H}^{n+1}( \overline{\Prism}_{R/A}\{n+1\}) \ar[d]^{\sim} \\
 \widehat{\Omega}^{n}_{R/\overline{A}} \ar[r]^-{d} &  \widehat{\Omega}^{n+1}_{R/\overline{A}}, }$$
where the lower horizontal map is the de Rham differential and the vertical maps are induced by the Hodge-Tate comparison isomorphisms 
of Remark \ref{remark:derived-HT-filtration}. Moreover, the Hodge-Tate comparison isomorphism is essentially determined by this property
(together with the requirement that it depends functorially on $R$). 
\end{remark}

\begin{remark}\label{remark:algebra-structure-on-relative-prismatic-complex}
For every animated commutative $\overline{A}$-algebra $R$, the relative Hodge-Tate complex $\overline{\Prism}_{R/\overline{A}}$ admits the structure of a commutative algebra
object of $\widehat{\calD}( \overline{A} )$, which is compatible with the conjugate filtration and the Hodge-Tate comparison isomorphism of
Remark \ref{remark:derived-HT-filtration}. In particular, we have a canonical isomorphism
$$\Fil_{0}^{\conj} \overline{\Prism}_{R/A} \simeq L \widehat{\Omega}^{0}_{R/\overline{A}} = \widehat{R},$$
where $\widehat{R}$ denotes the $p$-completion of $R$, so that $\overline{\Prism}_{R/\overline{A}}$ inherits the structure of a $p$-complete
$E_{\infty}$-algebra over $R$ (that is, a commutative algebra object of the $\infty$-category $\widehat{\calD}(R)$). 
\end{remark}

\begin{remark}\label{remark:connectivity-of-conjugate}
Let $R$ be an animated commutative $\overline{A}$-algebra. For every integer $n \geq 0$, the cohomology groups of the complex $L \widehat{\Omega}^{n}_{R/ \overline{A}}$ are concentrated in degrees $\leq 0$. Proceeding by induction on $n$, we obtain the following:
\begin{itemize}
\item The cohomology groups of the complex $\Fil_{n}^{\conj} \overline{\Prism}_{R/A}$ are concentrated in degrees $\leq n$.
\item The natural map 
$$ \mathrm{H}^{n}( \Fil_{n}^{\conj} \overline{\Prism}_{R/A} )
\rightarrow \mathrm{H}^{n}( \gr_{n}^{\conj} \overline{\Prism}_{R/A} ) \simeq \mathrm{H}^{0}( L \widehat{\Omega}^{n}_{ R / \overline{A} }\{-n\} )
\simeq \widehat{\Omega}^{n}_{R/\overline{A}}\{-n\}$$
is an isomorphism.
\end{itemize}
\end{remark}

\begin{proposition}\label{proposition:formally-etale-relative-Hodge-Tate}
Let $f: R \rightarrow S$ be a morphism of animated commutative $\overline{A}$-algebras which is $p$-completely formally \'{e}tale: that is, the $p$-complete
relative cotangent complex $L \widehat{\Omega}^{1}_{S/R}$ vanishes. Then the tautological map $S \widehat{\otimes}^{L}_{R} \overline{\Prism}_{R/A} \rightarrow \overline{\Prism}_{S/A}$
is an isomorphism in the $\infty$-category $\widehat{\calD}( \overline{A} )$.
\end{proposition}

\begin{proof}
The vanishing assumption on $\widehat{L}_{S/R}$ guarantees that, for each $m \geq 0$, the induced map $S \otimes^{L}_{R} L\Omega^{n}_{R} \rightarrow L \Omega^{n}_{S}$
becomes an isomorphism after $p$-completion. Using the Hodge-Tate comparison (Remark \ref{remark:derived-HT-filtration}), it follows by induction on $n$ that
$f$ induces isomorphisms $S \widehat{\otimes}^{L}_{R} \Fil^{\conj}_{n} \overline{\Prism}_{R/A} \rightarrow \Fil^{\conj}_{n} \overline{\Prism}_{S/A}$
for each $n \geq 0$. Proposition \ref{proposition:formally-etale-relative-Hodge-Tate} follows by passing to the colimit over $n$.
\end{proof}

\begin{corollary}\label{salmage}
Let $R$ be an animated commutative $\overline{A}$-algebra and let $\widehat{R}$ denote the $p$-completion of $R$.
Then the tautological maps
$$ \overline{\Prism}_{R/A} \rightarrow \overline{\Prism}_{ \widehat{R} / A} \quad \quad \Prism_{R/A} \rightarrow \Prism_{ \widehat{R} / A}$$
are isomorphisms.
\end{corollary}

\begin{corollary}\label{corollary:etale-descent-relative-prismatic}
The functor
$$ \CAlg^{\anim}_{\overline{A}} \rightarrow \widehat{\calD}(A) \quad \quad R \mapsto \Prism_{R/A}$$
satisfies descent for the \'{e}tale topology.
\end{corollary}

\begin{proof}
Extension of scalars along the quotient map $A \twoheadrightarrow \overline{A}$ induces a functor
$\widehat{\calD}(A) \rightarrow \widehat{\calD}( \overline{A} )$ which is conservative and preserves small limits.
It will therefore suffice to show that the functor $R \mapsto \overline{\Prism}_{R/A}$ satisfies descent for the \'{e}tale topology,
which is an immediate consequence of Proposition \ref{proposition:formally-etale-relative-Hodge-Tate}.
\end{proof}

\begin{variant}\label{variant:conjugate-fpqc-descent}
For each integer $n$, the functor $R \mapsto \Fil_{n}^{\conj} \overline{\Prism}_{R/A}$ satisfies
$p$-complete faithfully flat descent. This follows by induction on $n$, using the fact
that the functor $R \mapsto \gr_{n}^{\conj} \overline{\Prism}_{R/A} \simeq L \widehat{\Omega}^{n}_{R/\overline{A}}\{-n\}[-n]$
satisfies $p$-complete faithfully flat descent (see \cite[Remark 2.8]{bhatt-completions}, \cite[Theorem 3.1]{BMS2}). 
\end{variant}

\begin{corollary}\label{corollary:etale-change-of-prism}
Let $(A,I) \rightarrow (B,IB)$ be a morphism of bounded prisms for which the induced map $\overline{A} \rightarrow B/IB$
is $p$-completely formally \'{e}tale, and let $R$ be an animated commutative $(B/IB)$-algebra. Then the tautological map
$\Prism_{R/A} \rightarrow \Prism_{R/B}$ is an isomorphism in the $\infty$-category $\widehat{\calD}(A)$.
\end{corollary}

\begin{proof}
As in the proof of Corollary \ref{corollary:etale-descent-relative-prismatic}, it will suffice to show that the natural map
$\theta: \overline{\Prism}_{R/A} \rightarrow \overline{\Prism}_{R/B}$ is an isomorphism in $\widehat{\calD}( \overline{A} )$. Setting $R' = (B/IB) \otimes^{L}_{\overline{A}} R$, we observe that
$\theta$ factors as a composition
$$ \overline{\Prism}_{R/A} \xrightarrow{\theta'} R \widehat{\otimes}^{L}_{R'} \overline{\Prism}_{R'/B}
\xrightarrow{\theta''} \overline{\Prism}_{R/B}$$
where $\theta'$ and $\theta''$ are invertible (by virtue of Remark \ref{remark:change-of-prism} and
Proposition \ref{proposition:formally-etale-relative-Hodge-Tate}, respectively).
\end{proof}

\begin{remark}\label{remark:hypothesis-for-discreteness}
Let $R$ be an $\overline{A}$-algebra which satisfies the following condition:
\begin{itemize}
\item[$(\ast)$] The commutative ring $R$ has bounded $p$-torsion and the derived tensor product $R/pR \otimes_{R}^{L} L_{R/\overline{A}} \in \calD(R/pR)$ has $\Tor$-amplitude contained in $[-1,0]$.
\end{itemize}
Then:
\begin{itemize}
\item For every integer $n \geq 0$, the derived tensor product $R/pR \otimes_{R}^{L} L \widehat{\Omega}^{n}_{R/\overline{A}} \in \calD(R/pR)$
has $\Tor$-amplitude contained in $[-n,0]$.

\item For every integer $n \geq 0$, the complex 
$$R/pR \otimes^{L}_{R} \gr_{n}^{\conj} \overline{\Prism}_{R/A}$$ 
has $\Tor$-amplitude contained in $[0,n]$ as an object of $\calD( R/ pR )$.

\item For every integer $n \geq 0$, the complex $(R/pR) \otimes_{R}^{L} \Fil_{n}^{\conj} \overline{\Prism}_{R/A}$
has $\Tor$-amplitude contained in $[0,n]$ as an object of $\calD(R/pR)$.

\item The complex $R/pR \otimes^{L}_{R} \overline{\Prism}_{R/A}$ has $\Tor$-amplitude contained in $[0,\infty]$
as an object of $\calD(R/pR)$.

\item The cohomology groups of the relative Hodge-Tate complex $\overline{\Prism}_{R/A}$ and the relative prismatic complex $\Prism_{R/A}$ are concentrated in degrees $\geq 0$.
\end{itemize}
\end{remark}

\begin{variant}\label{variant:prismatic-coh-discrete}
Let $R$ be an $\overline{A}$-algebra which satisfies the following condition:
\begin{itemize}
\item[$(\ast^{+})$] The commutative ring $R$ has bounded $p$-torsion and the complex $R/pR \otimes_{R}^{L} L_{R/\overline{A}}$ is isomorphic to a flat $R/pR$-module, concentrated
in cohomological degree $-1$.
\end{itemize}
Then:
\begin{itemize}
\item For every integer $n \geq 0$, the derived tensor product $R/pR \otimes_{R}^{L} L \widehat{\Omega}^{n}_{R/\overline{A}} \in \widehat{\calD}(R/pR)$
is isomorphic to a flat $R/pR$-module, concentrated in cohomological degree $-n$.

\item For every integer $n \geq 0$, the derived tensor product
$$R/pR \otimes^{L}_{R} \gr_{n}^{\conj} \overline{\Prism}_{R/A}$$
is isomorphic to a flat $R/pR$-module, concentrated in cohomological degree $0$.

\item For every integer $n \geq 0$, the derived tensor product
$$R/pR \otimes^{L}_{R} \Fil_{n}^{\conj} \overline{\Prism}_{R/A}$$
is isomorphic to a flat $R/pR$-module, concentrated in cohomological degree $0$.

\item The complex $R/pR \otimes^{L}_{R} \overline{\Prism}_{R/A}$ is isomorphic to a flat $R/pR$-module, concentrated in cohomological degree zero.

\item The relative Hodge-Tate complex $\overline{\Prism}_{R/A}$ is concentrated in cohomological degree zero and is $p$-completely flat as an $R$-module.
Consequently, the relative prismatic complex $\Prism_{R/A}$ is also concentrated in cohomological degree zero.
\end{itemize}
If, in addition, $R$ is $p$-completely flat over $\overline{A}$, it follows that $\overline{\Prism}_{R/A}$ is $p$-completely flat over $\overline{A}$,
so that $\Prism_{R/A}$ is $(p,I)$-completely flat over $A$. 
\end{variant}

\begin{remark}\label{remark:hypothesis-for-surjectivity}
Let $u: R' \twoheadrightarrow R$ be a surjective homomorphism of $\overline{A}$-algebras.
Assume that $R$ and $R'$ satisfy condition $(\ast^{+})$ of Variant \ref{variant:prismatic-coh-discrete} and that the $p$-complete
cotangent complex $L \widehat{\Omega}^{1}_{R/R'}$ has cohomology concentrated in degrees $\leq -2$. Then:
\begin{itemize}
\item For every integer $n \geq 0$, the induced map
$$ (R/pR) \otimes^{L}_{R'} L \widehat{\Omega}^{n}_{R'/\overline{A}} \rightarrow (R/pR) \otimes^{L}_{R}  L \widehat{\Omega}^{n}_{R'/\overline{A}}$$
is a surjection of flat $(R/pR)$-modules, concentrated in cohomological degree $-n$.

\item For every integer $n \geq 0$, the induced map
$$R/pR \otimes^{L}_{R'} \gr_{n}^{\conj} \overline{\Prism}_{R'/A} \rightarrow R/pR \otimes^{L}_{R} \gr_{n}^{\conj} \overline{\Prism}_{R/A}$$
is a surjection of flat $R/pR$-modules, concentrated in cohomological degree $0$.

\item For every integer $n \geq 0$, the induced map
$$ (R/pR) \otimes^{L}_{R'} \Fil_{n}^{\conj} \overline{\Prism}_{R'/A} \rightarrow (R/pR) \otimes^{L}_{R} \Fil_{n}^{\conj} \overline{\Prism}_{R/A}$$
is a surjection of flat $(R/pR)$-modules, concentrated in cohomological degree $0$.

\item The induced map
$$ (R/pR) \otimes_{R'}^{L} \overline{\Prism}_{R'/A} \rightarrow (R/pR) \otimes_{R}^{L} \overline{\Prism}_{R/A}$$
is a surjection of flat $(R/pR)$-modules, concentrated in cohomological degree zero.

\item The induced map $\overline{\Prism}_{R'/A} \rightarrow \overline{\Prism}_{R/A}$ is a surjection of $p$-adically separated and complete abelian groups.

\item The induced map $\Prism_{R'/A} \rightarrow \Prism_{R/A}$ is surjective.
\end{itemize}
\end{remark}

\subsection{Relative Prismatic Cohomology: Globalization}\label{subsection:relative-prismatic-global}

Let $(A,I)$ be a bounded prism, which we regard as fixed throughout this section, and let $\overline{A}$ denote the quotient ring $A/I$.

\begin{construction}\label{construction:prismatic-complex-of-scheme}
Let $X$ be scheme, formal scheme, or algebraic stack equipped with a morphism $\pi: X \rightarrow \Spec( \overline{A} )$. We let
$\RGamma_{\Prism}(X/A)$ denote the inverse limit $\varprojlim_{ \Spec(R) \rightarrow X} \Prism_{R/A}$, formed in the $\infty$-category $\widehat{\calD}(A)$. Here the limit is indexed by
the {\it category of points} of $X$: that is, the category of pairs $(R, f)$ where $R$ is a commutative ring and $f: \Spec(R) \rightarrow X$ is a morphism
(so that the composition $\pi \circ f$ determines an $\overline{A}$-algebra structure on $R$). We will refer to $\RGamma_{\Prism}(X/A)$ as the {\it prismatic complex of $X$ relative to $A$}.
We let $\RGamma_{\overline{\Prism}}( X / A)$ denote the derived tensor product $\overline{A} \otimes^{L}_{A} \RGamma_{\Prism}(X/A)$; we will refer to $\RGamma_{\overline{\Prism}}( X / A)$ as the {\it 
Hodge-Tate complex of $X$ relative to $A$}. 

For every integer $n$, we let $\mathrm{H}^{n}_{\Prism}( X/A )$ denote the $n$th cohomology group of the complex $\RGamma_{\Prism}(X/A)$, and we let
$\mathrm{H}^{n}_{\overline{\Prism} }( X/ A)$ denote the $n$th cohomology group of the complex $\RGamma_{ \overline{\Prism} }(X/A)$. We will refer to 
we refer to $\mathrm{H}^{n}_{\Prism}(X/A)$ as the {\it $n$th prismatic cohomology group of $X$ relative to $A$}, and to $\mathrm{H}^{n}_{ \overline{\Prism} }(X/A)$ as the
{\it $n$th Hodge-Tate cohomology group of $X$ relative to $A$}.
\end{construction}

\begin{remark}\label{remark:etale-descent-relative}
The construction $X \mapsto \RGamma_{\Prism}(X/A)$ satisfies descent with respect to the \'{e}tale topology (this is a formal consequence of Corollary \ref{corollary:etale-descent-relative-prismatic}).
\end{remark}

\begin{example}\label{example:sotto}
Let $X = \Spec(R)$ be an affine $\overline{A}$-scheme. Then the prismatic complex $\RGamma_{\Prism}(X/A)$ of Construction \ref{construction:prismatic-complex-of-scheme}
can be identified with the relative prismatic complex $\Prism_{R/A}$ of Construction \ref{construction:relative-prismatic-cohomology}.
\end{example}

\begin{example}\label{example:compute-for-formal-scheme}
Let $\mathfrak{X}$ be a bounded $p$-adic formal scheme over $\overline{A}$. For each $n \geq 0$, let $\mathfrak{X}_{n} \subseteq \mathfrak{X}$
denote the closed subscheme given by the vanishing locus of $p^{n}$, so that we have $\mathfrak{X} = \varinjlim_{n} \mathfrak{X}_{n}$ (as set-valued functors on the category of commutative rings).
It follows that $\RGamma_{\Prism}( \mathfrak{X} / A )$ can be identified with the limit $\varprojlim_{n} \RGamma_{\Prism}( \mathfrak{X}_{n} / A )$, formed in the $\infty$-category
$\widehat{\calD}(A)$. In particular, the construction $\mathfrak{X} \mapsto \RGamma_{\Prism}( \mathfrak{X} / A)$ satisfies descent for the \'{e}tale topology on
the category of bounded $p$-adic formal $\overline{A}$-schemes.
\end{example}

In practice, little generality is lost by restricting Construction \ref{construction:prismatic-complex-of-scheme} to the setting of formal schemes:

\begin{proposition}\label{proposition:scheme-vs-formal}
Let $X$ be an $\overline{A}$-scheme for which the structure sheaf $\calO_{X}$ has bounded $p$-power torsion (this condition is satisfied, for example, if $X$ is smooth over $\overline{A}$), and let $\mathfrak{X} = \Spf( \overline{A} ) \times_{\Spec(\overline{A} )} X$ denote its underlying $p$-adic formal scheme.
Then the restriction map $\RGamma_{\Prism}(X/A) \rightarrow \RGamma_{\Prism}(\mathfrak{X}/A)$ is an isomorphism in the $\infty$-category $\widehat{\calD}(A)$.
\end{proposition}

\begin{corollary}\label{corollary:sashy}
Let $R$ be an $\overline{A}$-algebra having bounded $p$-power torsion, and let $\widehat{R}$ denote the $p$-completion of $R$. Then the tautological map
$\Prism_{R/A} \rightarrow \RGamma_{\Prism}( \Spf( \widehat{R} ) / A)$ is an isomorphism in $\widehat{\calD}(A)$.
\end{corollary}

\begin{proof}
Combine Proposition \ref{proposition:scheme-vs-formal} with Example \ref{example:sotto}.
\end{proof}

Our proof of Proposition \ref{proposition:scheme-vs-formal} will require some preliminaries.

\begin{notation}\label{notation:F-p-bullet}
For each integer $n \geq 0$, let $\F_p^{\otimes n}$ denote the derived tensor product
$$ \F_p \otimes^{L}_{\Z} \F_p \otimes^{L}_{\Z} \cdots \otimes^{L}_{\Z} \F_p$$
of $n$ copies of the finite field $\F_p$, which we regard as an animated commutative ring. The construction $[n] \mapsto \F_p^{\otimes n+1}$ determines
a cosimplicial object of the $\infty$-category $\CAlg^{\anim}$, which we will denote by $\F_p^{\otimes \bullet +1 }$. 
\end{notation}

\begin{proposition}[Derived Descent for Relative Prismatic Cohomology]\label{proposition:derived-descent-relative}
Let $R$ be an animated commutative $\overline{A}$-algebra, and regard $R^{\bullet} = \F_p^{\otimes \bullet+1} \otimes^{L} R$ as a cosimplicial animated
commutative $R$-algebra. Then the tautological map $\Prism_{R/A} \rightarrow \Tot( \Prism_{ R^{\bullet} / A} )$ is an isomorphism in $\widehat{\calD}(A)$.
In other words, $\Prism_{R/A}$ can be identified with the inverse limit of the diagram $[n] \mapsto \Prism_{ ( \F_p^{\otimes n+1} \otimes^{L} R)/A}$.
\end{proposition}

\begin{proof}
By virtue of Remark \ref{remark:change-of-prism}, the assertion can be tested locally (with respect to the flat topology) on the category of prisms,
so we may assume without loss of generality that the ideal $I = (d)$ is principal. Let $S$ denote the derived tensor product
$\F_p \otimes^{L} \overline{A}$; it will then suffice to show that the comparison map
$$ S \otimes^{L}_{A} \Prism_{R/A} \rightarrow S \otimes^{L}_{A} \Tot( \Prism_{R^{\bullet}/A} ) \simeq
\Tot( S \otimes^{L}_{A} \Prism_{R^{\bullet}/A} )$$
is invertible. Let $S^{\bullet}$ denote the cosimplicial commutative algebra object of $\calD(S)$ given by
$$S^{\bullet} = S \otimes^{L}_{A} \Prism_{ (\overline{A} \otimes^{L} \F_p^{\otimes \bullet+1} ) / A}.$$
Using the K\"{u}nneth formula of Remark \ref{remark:Kunneth-relative-prismatic}, we are reduced to proving that
the tautological map $M \rightarrow \Tot( M \otimes^{L}_{S} S^{\bullet} )$ is an isomorphism in the $\infty$-category $\calD( S)$ for
$M = S \otimes^{L}_{A} \Prism_{ R / A}$. Moreover, the cosimplicial object $S^{\bullet}$ is given by the iterated (derived) tensor powers of $S^{0}$ relative to $S$. We will complete the proof by showing that the unit map $u: S \rightarrow S^{0}$ generates a covering sieve with respect to the faithfully flat topology: that is, there exists a faithfully flat
morphism of $S \rightarrow S'$ which factors through $u$. Note that this assertion does not depend on the animated commutative $\overline{A}$-algebra $R$.

Let $(A^{0}, I^{0} )$ be the prism of Example \ref{example:rho-A-special-case}. Since the ideal $I$ is principal, there exists a morphism of prisms $(A^{0}, I^{0} ) \rightarrow (A,I)$.
Using Remark \ref{remark:change-of-prism}, we can reduce to the case where $(A,I) = (A^{0}, I^{0} )$. In particular, $(A,I)$ is a transversal prism, so that
$S = \overline{A} / p \overline{A}$ is an ordinary commutative ring.
Let $W = \Spec( \Z[a_0, a_1, \cdots ])$ denote the $\Z_p$-scheme representing the Witt vector functor $B \mapsto W(B)$,
so that we can identify the $\Spf(A)$ with the formal subscheme $\WCart_0 \subseteq W$ described in
Construction \ref{construction:Cartier-Witt-atlas}. This restricts to an identification of
$\Spec(S)$ with the locally closed subscheme of $\Spec(\F_p) \times W$ given by the locus where $a_0$ vanishes and $a_1$ is invertible
(hence an isomorphism of commutative rings $S \simeq \F_p[ a_{1}^{\pm 1}, a_2, a_3, \cdots ]$). Similarly, we can identify
$\Prism_{S/A}$ with the $\delta$-ring
$A\{ \frac{p}{I^0} \}$, so that $( \Prism_{S/A}, I^{0} \Prism_{S/A } )$ is the universal crystalline prism over $(A, I )$ (see 
Example \ref{example:transversal-prism-modulo-p} below).
In particular, there is an identification of the formal spectrum $\Spf( \Prism_{S/A })$ with the product $\Spf(\Z_p) \times W^{\times} \subseteq W^{\times}$,
fitting into a commutative diagram
\begin{equation}
\begin{gathered}\label{equation:diagram-Frobenius-factorization}
\xymatrix@R=50pt@C=50pt{ \Spf( \Prism_{S/A }) \ar@{^{(}->}[r] \ar[d] & W^{\times} \ar[d]^{p} \\
\Spf( A ) \ar@{^{(}->}[r] & W. }
\end{gathered}
\end{equation}
When $B$ is an $\F_p$-algebra, the multiplication map $p: W(B) \rightarrow W(B)$ factors as a composition $W(B) \xrightarrow{F} W(B) \xrightarrow{V} W(B)$.
We can therefore expand (\ref{equation:diagram-Frobenius-factorization}) to a commutative diagram of formal schemes
$$ \xymatrix@R=50pt@C=50pt{ \Spec(\F_p) \times W^{\times} \ar@{^{(}->}[r] \ar[d]^{F} & \Spf( \Prism_{S/A }) \ar@{^{(}->}[r] \ar[dd] & W^{\times} \ar[dd]^{p} \\
\Spec(\F_p) \times W^{\times} \ar[d]^{V}_{\sim} & & \\
\Spec(S) \ar@{^{(}->}[r] & \Spf( A ) \ar@{^{(}->}[r] & W. }$$
We complete the proof by observing that the morphism of affine schemes $F: \Spec(\F_p) \times W^{\times} \rightarrow \Spec(\F_p) \times W^{\times}$
is faithfully flat (Corollary \ref{corollary:W-F-flat}).
\end{proof}

\begin{corollary}\label{corollary:mod-pn}
Let $R$ be an animated commutative $\overline{A}$-algebra. For each $n \geq 0$, let $R_{n}$ denote the derived tensor product $(\Z / p^{n} \Z) \otimes^{L} R$.
Then the tautological map $\Prism_{R/A} \rightarrow \varprojlim_{n} \Prism_{R_n/A}$ is an isomorphism in the $\infty$-category $\widehat{\calD}(A)$.
\end{corollary}

\begin{proof}
By virtue of Proposition \ref{proposition:derived-descent-relative}, we may assume without loss of generality that $R$ is an algebra over the derived tensor product
$\F_p \otimes^{L} \overline{A}$. In this case, the tower $\{ R_n \}_{n \geq 0}$ is isomorphic to $R$ as a pro-object of the $\infty$-category $\CAlg^{\anim}_{\overline{A}}$,
so the result is clear.
\end{proof}

\begin{proof}[Proof of Proposition \ref{proposition:scheme-vs-formal}]
Let $X$ be an $\overline{A}$-scheme for which the structure sheaf $\calO_{X}$ has bounded $p$-power torsion, and let $\mathfrak{X} = \Spf( \overline{A} ) \times_{\Spec( \overline{A} )} X$
denote its $p$-completion. We wish to show that the restriction map $\RGamma_{\Prism}( X/A) \rightarrow \RGamma_{\Prism}( \mathfrak{X} / A)$ is an isomorphism in
$\widehat{\calD}(A)$. Without loss of generality, we may assume that $X = \Spec(R)$ is affine (Remark \ref{remark:etale-descent-relative}), so that $R$ is an $\overline{A}$-algebra of bounded $p$-power torsion. Using Example \ref{example:compute-for-formal-scheme},
we are reduced to showing that the canonical map $\Prism_{R/A} \rightarrow  \varprojlim_{n} \Prism_{ (R/p^nR) / A} \simeq \RGamma_{\Prism}( \mathfrak{X} / A)$
is an isomorphism in $\widehat{\calD}(A)$. This follows from Corollary \ref{corollary:mod-pn}, since the tautological maps
$$ (\Z / p^{n} \Z) \otimes^{L} R \rightarrow (\Z / p^{n} \Z) \otimes R = R/p^{n}R$$
define an isomorphism $\{  (\Z / p^{n} \Z) \otimes^{L} R \}_{n \geq 0} \xrightarrow{\sim} \{ R/p^{n} R \}_{n \geq 0}$ of pro-objects of $\CAlg^{\anim}_{\overline{A}}$.
\end{proof}

For later use, we record a version of Proposition \ref{proposition:derived-descent-relative} for the conjugate filtration on relative Hodge-Tate cohomology.

\begin{variant}\label{variant:conjugate-derived-descent}
Let $R$ be an animated commutative $\overline{A}$-algebra and let $R^{\bullet} = R \otimes^{L} \F_p^{\otimes \bullet +1 }$ be as in the statement of
Proposition \ref{proposition:derived-descent-relative}. For every complex $M \in \widehat{\calD}(\Z_p)$ and every integer $m$, the tautological map
$$ \theta: M \widehat{\otimes}^{L} \Fil_{m}^{\conj} \overline{\Prism}_{ R/A} \rightarrow
\Tot( M \widehat{\otimes}^{L} \Fil_{m}^{\conj} \overline{\Prism}_{ R^{\bullet} /A} )$$
is an isomorphism in $\widehat{\calD}(R)$.
\end{variant}

\begin{proof}
Proceeding by induction on $m$ and using Remark \ref{remark:derived-HT-filtration}), we are reduced to proving the following assertion:
\begin{itemize}
\item[$(\ast)$] Let $M$ be an object of $\widehat{\calD}(\Z_p)$. Then the comparison map
$$ M \widehat{\otimes}^{L} L\widehat{\Omega}^{m}_{R/\overline{A}} \rightarrow 
\Tot( M \widehat{\otimes}^{L} L \widehat{\Omega}^{m}_{R^{\bullet}/ \overline{A}} )$$
is an isomorphism.
\end{itemize}
We have a fiber sequence of cosimplicial $R^{\bullet}$-modules
$R^{\bullet} \otimes_{R}^{L} L\widehat{\Omega}^{1}_{R/\overline{A}} \rightarrow L\widehat{\Omega}^{1}_{ R^{\bullet} / \overline{A} } \rightarrow L\widehat{\Omega}^{1}_{ R^{\bullet} / R }$. This determines a finite filtration of $L \widehat{\Omega}^{m}_{R^{\bullet}/ \overline{A}} )$ whose associated graded is given by
$(L\widehat{\Omega}^{m-d}_{R/\overline{A}}) \widehat{\otimes}_{R}^{L} (L \widehat{\Omega}^{d}_{R^{\bullet}/R} )$.
Moreover, the comparison map of $(\ast)$ can be promoted to a map of filtered objects, given at the associated graded level by
$$ M \widehat{\otimes}^{L} ( L \widehat{\Omega}^{m-d}_{R/\overline{A}}) \widehat{\otimes}_{R}^{L} ( L \widehat{\Omega}^{d}_{R/R} )
\rightarrow \Tot(M \widehat{\otimes} (L \widehat{\Omega}^{m-d}_{R/\overline{A}}) \widehat{\otimes}_{R}^{L} ( L\widehat{\Omega}^{d}_{R^{\bullet}/R} )).$$
Setting $N = M \widehat{\otimes} \widehat{\Omega}^{m-d}_{R/\overline{A}}$, we are reduced to proving the following:
\begin{itemize}
\item[$(\ast')$] Let $N$ be an object of $\widehat{\calD}(R)$ and let $d \geq 0$ be an integer. Then the comparison map
$$ N \widehat{\otimes}^L_{R} L \widehat{\Omega}^{d}_{ R/ R} \rightarrow \Tot( N \widehat{\otimes}^{L}_{R}  L \widehat{\Omega}^{d}_{ R^{\bullet} / R} )$$
is an isomorphism.
\end{itemize}
Since both sides are $p$-complete, this can be checked after extending scalars to $\F_p$. We may therefore assume without loss of generality that
$R$ has the structure of an animated $\F_p$-algebra, in which case the result is clear (since the augmented cosimplicial object $R \mapsto R^{\bullet}$ splits).
\end{proof}

\subsection{Calculation via the Relative Prismatic Site}\label{subsection:relative-prismatic-site}

We now compare Construction \ref{construction:prismatic-complex-of-scheme} with the definition of prismatic cohomology given in \cite{prisms}.

\begin{notation}\label{notation:comparison-with-site-theoretic}
Let $\mathfrak{X}$ be a bounded $p$-adic formal scheme over $\overline{A}$, and let $( \mathfrak{X} / A)_{\Prism}$ be the relative prismatic site of
Definition \ref{definition:relative-prismatic-site}. We let $\RGamma_{\Prism}^{\mathrm{site}}( \mathfrak{X} / A)$ denote the limit
$$ \varprojlim_{ (B,v) \in ( \mathfrak{X} / A )_{\Prism} } B,$$
formed in the $\infty$-category $\widehat{\calD}(A)$. We will refer to $\RGamma_{\Prism}^{\mathrm{site}}( \mathfrak{X} / A)$
as the {\it site-theoretic prismatic complex} of $\mathfrak{X}$. 

Let $R$ be an $\overline{A}$-algebra with bounded $p$-torsion, and let $\mathfrak{X} = \Spf( \overline{A} ) \times_{\Spec(A)} \Spec(R)$
denote the associated $p$-adic formal scheme. In this case, we will denote the complex $\RGamma_{\Prism}^{\mathrm{site}}( \mathfrak{X} / A)$ by
$\Prism_{R/A}^{\mathrm{site}}$.
\end{notation}

\begin{warning}
In the definition of $\Prism_{R/A}^{\mathrm{site}}$, we implicitly assume that the limit $\varprojlim_{(B,J,v) \in (R/A)_{\Prism} } B$ exists in the $\infty$-category $\widehat{\calD}(A)$. This is not completely obvious, since the prismatic site $(R/A)_{\Prism}$ is not a small category (or even equivalent to a small category). In what follows, we will ignore this technicality (it can be addressed, for example, by regarding $\Prism_{R/A}^{\mathrm{site}}$ as belonging to a larger universe).
\end{warning}

\begin{example}\label{example:isomorphism-for-polynomial}
By construction, the relative prismatic complex $\Prism_{R/A}$ coincides with $\Prism_{R/A}^{\mathrm{site}}$ in the special case where
$R$ is a finitely generated polynomial algebra over $\overline{A}$. In particular, Corollary \ref{corollary:sashy}
supplies a canonical isomorphism $\RGamma_{\Prism}( \Spf( \widehat{R} ) / A) \simeq \RGamma_{\Prism}^{\mathrm{site}}( \Spf( \widehat{R} ) / A)$.
\end{example}

\begin{proposition}\label{proposition:comparison-with-site}
Let $\mathfrak{X}$ be a bounded $p$-adic formal scheme over $\overline{A}$. Then there is a comparison map
$$ \xi_{\mathfrak{X}}: \RGamma_{\Prism}( \mathfrak{X} / A) \rightarrow \RGamma_{\Prism}^{\mathrm{site}}( \mathfrak{X} / A),$$
which is characterized (up to homotopy) by the requirement that it depends functorially on $\mathfrak{X}$
and agrees with the isomorphism of Example \ref{example:isomorphism-for-polynomial} when $\mathfrak{X}$ has the form
$\Spf( \widehat{R} )$, where $R$ is a finitely generated polynomial algebra over $\overline{A}$.
\end{proposition}

\begin{proof}
Let $\calC$ be the opposite of the category of bounded $p$-adic formal schemes over $\overline{A}$ and let $\calC' \subseteq \calC$
be the full subcategory spanned by the affine formal schemes. It is not difficult to see that the construction
$\mathfrak{X} \mapsto \RGamma_{\Prism}^{\mathrm{site}}( \mathfrak{X} / A)$ satisfies descent for the \'{e}tale topology, and is therefore
a right Kan extension of its restriction to the subcategory $\calC' \subseteq \calC$. It will therefore suffice to construct
the comparison map $\xi_{\mathfrak{X}}$ when $\mathfrak{X}$ is affine. Let $\calC'' \subseteq \calC'$ be the full subcategory spanned
by formal $\overline{A}$-schemes of the form $\Spf( \widehat{R} )$, where $R$ is a finitely generated polynomial algebra over $R$.
Using Corollary \ref{corollary:sashy} (and the definition of $\Prism_{R/A}$ as a right Kan extension), we see that the functor
$$(\mathfrak{X} \in \calC') \mapsto \RGamma_{\Prism}( \mathfrak{X} / A)$$
is a left Kan extension of its restriction to $\calC''$. It will therefore suffice to construct the comparison map $\xi_{\mathfrak{X}}$ for
$\mathfrak{X} \in \calC''$, where it is supplied by Example \ref{example:isomorphism-for-polynomial}.
\end{proof}

\begin{notation}
Let $R$ be an $\overline{A}$-algebra with bounded $p$-torsion. Then the comparison map $\xi_{\mathfrak{X}}$ of Proposition \ref{proposition:comparison-with-site} can be regarded as a morphism from
$\Prism_{R/A} \rightarrow \Prism_{R/A}^{\mathrm{site}}$, which we will denote by $\xi_{R}$.
\end{notation}

We can now state our main result:

\begin{theorem}\label{theorem:site-theoretic-equivalence}
Let $\mathfrak{X}$ be a bounded $p$-adic formal scheme over $\overline{A}$. Suppose that, for every affine open subset
$\Spf(R) \subseteq \mathfrak{X}$, the $\overline{A}$-algebra $R$ satisfies condition $(\ast)$ of Remark \ref{remark:hypothesis-for-discreteness}. Then
the comparison map 
$$ \xi_{\mathfrak{X}}: \RGamma_{\Prism}( \mathfrak{X} / A) \rightarrow \RGamma_{\Prism}^{\mathrm{site}}( \mathfrak{X} / A)$$
of Proposition \ref{proposition:comparison-with-site} is an isomorphism in $\widehat{\calD}(A)$.
\end{theorem}

\begin{remark}[Comparison with \cite{prisms}]
Let $\mathfrak{X}$ be a formal scheme which is smooth over $\Spf( \overline{A} )$. In \cite{prisms}, the prismatic complex
$\RGamma_{\Prism}( \mathfrak{X} / A)$ is defined to be the site-theoretic prismatic complex $\RGamma_{\Prism}^{\mathrm{site}}( \mathfrak{X} / A)$
of Notation \ref{notation:comparison-with-site-theoretic}. Theorem \ref{theorem:site-theoretic-equivalence} implies that this agrees with
with prismatic complex of Construction \ref{construction:prismatic-complex-of-scheme} (up to canonical isomorphism). Moreover, a similar result holds for
many formal schemes which are not smooth over $\overline{A}$.
\end{remark}

\begin{remark}\label{remark:perfect-relative-prismatic}
Let $\mathfrak{X}$ be a formal scheme which is smooth and proper over $\Spf( \overline{A} )$. Then the relative prismatic complex $\RGamma_{\Prism}( \mathfrak{X} / A)$
is perfect (see Theorem~1.8 of \cite{prisms}).
\end{remark}

We begin with an easy special case of Theorem \ref{theorem:site-theoretic-equivalence}.

\begin{remark}\label{remark:make-idempotent}
Let $R$ be a commutative algebra over $\overline{A}$ which satisfies condition $(\ast^{+} )$ of Variant \ref{variant:prismatic-coh-discrete}, so that the prismatic complex $\Prism_{R/A}$
is concentrated in degree zero, and let $\mathfrak{X}$ denote the $p$-completion of the affine scheme $\Spec(R)$. Applying Lemma~7.7 of \cite{prisms}, we deduce:
\begin{itemize}
\item The pair $( \Prism_{R/A}, I \Prism_{R/A} )$ can be regarded as an object of the relative prismatic site $(R/A)_{\Prism}$.
\item The relative prismatic site $(R/A)_{\Prism}$ has an initial object, which can then be identified with 
$( \Prism_{R/A}^{\mathrm{site}}, I \Prism_{R/A}^{\mathrm{site}} )$. 
\item The comparison map $\xi_{R}: \Prism_{R/A} \rightarrow \Prism_{R/A}^{\mathrm{site}}$
can be regarded a morphism in the prismatic site $( \mathfrak{X} / A)_{\Prism}$.
\end{itemize}
It follows that the composition
$$ ( \Prism_{R/A}, I \Prism_{R/A} ) \xrightarrow{ \xi_{R} } ( \Prism_{R/A}^{\mathrm{site}}, I \Prism_{R/A}^{\mathrm{site}} )
\rightarrow ( \Prism_{R/A}, I \Prism_{R/A} )$$
determines an idempotent endomorphism $\iota_{R/A}$ of the $\delta$-ring $\Prism_{R/A}$; this endomorphism is surjective exactly when it is an isomorphism, which happens if and only if $\xi_{R}$ is an isomorphism.

Note that $\iota_{R/A}$ depends functorially on both $R$ and $A$: if $(A,I) \rightarrow (A', I')$ is a morphism of bounded prisms
and $R'$ is a commutative algebra over the quotient ring $\overline{A}' = A' /I'$ which also satisfies condition $(\ast^{+})$ of Variant \ref{variant:prismatic-coh-discrete}, then we have
a commutative diagram of prisms
$$ \xymatrix@R=50pt@C=50pt{ (\Prism_{R/A}, I \Prism_{R/A} ) \ar[r]^-{\xi_{R}} \ar[d] & ( \Prism_{R/A}^{\mathrm{site}}, I \Prism_{R/A}^{\mathrm{site}}) \ar[r] \ar[d] 
& ( \Prism_{R/A}, I \Prism_{R/A} ) \ar[d] \\
 (\Prism_{R'/A'}, I' \Prism_{R'/A'} ) \ar[r]^-{\xi_{R'}} & ( \Prism_{R'/A'}^{\mathrm{site}}, I' \Prism_{R'/A'}^{\mathrm{site}}) \ar[r] & ( \Prism_{R'/A'}, I' \Prism_{R'/A'} ) }$$
where the top horizontal composition is $\iota_{R/A}$ and the bottom horizontal composition is $\iota_{R'/A'}$. Here the left square commutes by the construction of $\xi$, and the right square by virtue of the initiality of $( \Prism^{\mathrm{site}}_{R/A}, I \Prism^{\mathrm{site}}_{R/A} )$ as an object of the relative prismatic site
$( \mathfrak{X} /A)_{\Prism}$.
\end{remark}

\begin{lemma}\label{lemma:really-easy-case}
Let $R$ be a quotient of $\overline{A}$ which satisfies condition $(\ast)$ of 
Remark \ref{remark:hypothesis-for-discreteness}. Then the comparison map
$\xi_{R}: \Prism_{R/A} \rightarrow \Prism^{\mathrm{site}}_{R/A}$ is an isomorphism.
\end{lemma}

\begin{proof}
It follows from the surjectivity of the unit map $u: \overline{A} \twoheadrightarrow R$ that the module of K\"{a}hler differentials $\Omega_{R/\overline{A}}$ vanishes,
so that $R$ satisfies the stronger condition $(\ast^{+} )$ of $(\ast^{+})$ of Variant \ref{variant:prismatic-coh-discrete}. Let $\iota_{R/A}$ be the idempotent endomorphism of Remark \ref{remark:make-idempotent}. We wish to show that $\iota_{R/A}$ is an isomorphism, or equivalently that is surjective.
Let $\{ f_{s} \}_{s \in S}$ be a collection of generators for the kernel $\ker(u)$, and let $R'$ denote the quotient of
$\overline{A}[ Y_{t}^{1/p^{\infty}} ]$ by the ideal generated by the elements $\{ Y_{t} - f_{t} \}_{t \in T}$. Then there is a surjection of $\overline{A}$-algebras
$v: R' \twoheadrightarrow R$ which is uniquely determined by the requirement that $v( Y_{t}^{1/p^{n}} )$ vanishes for each $t \in T$ and each integer $n \geq 0$.
By construction, the induced map $(R/p) \otimes^{L}_{R'} L_{ R' / \overline{A} }[-1] \rightarrow (R/p) \otimes (R/p) \otimes^{L}_{R} L_{R/ \overline{A}}[-1]$ is a surjection of flat
$R/pR$-algebras. By virtue of Remark \ref{remark:hypothesis-for-surjectivity}, the induced map of $\delta$-rings $\Prism_{R'/A} \rightarrow \Prism_{R/A}$ is surjective.
Consequently, to show that the idempotent $\iota_{R/A}$ is surjective, it will suffice to show that $\iota_{R'/A}$ is surjective. Since the functor
$$ \Prism_{\bullet/A}: \CAlg^{\anim}_{\overline{A}} \rightarrow \CAlg( \widehat{\calD}(A) )$$
preserves colimits (Remark \ref{remark:Kunneth-relative-prismatic}), it will suffice to show that the idempotent $\iota_{S/A}$ is an isomorphism
in the special case where $S$ is an $\overline{A}$-algebra of the form $\overline{A}[ Y^{1/p^{\infty}}] / ( Y - f )$, for some element $f \in \overline{A}$. In this case,
the desired result follows from Example~7.9 of \cite{prisms}.
\end{proof}

The following result is a slight generalization of Proposition~7.10 of \cite{prisms}:

\begin{lemma}\label{lemma:easy-qrsp-case}
Let $(A,I)$ be a bounded prism and let $R$ be a commutative algebra over the quotient ring $\overline{A} = A/I$.
Assume that $R$ satisfies condition $(\ast)$ of Remark \ref{remark:hypothesis-for-discreteness} and that there exists
a collection of elements $\{ x_s \}_{s \in S}$ satisfying the following:
\begin{itemize}
\item[$(a)$] Each $x_{s}$ admits a compatible system of $p$th-power roots $\{ x_{s}^{1/p^{n}} \}_{n \geq 0}$.
\item[$(b)$] The images of the elements $x_{s}$ generate the quotient ring $R/pR$.
\end{itemize}
Then the comparison map $\xi_{R}: \Prism_{R/A} \rightarrow \Prism^{\mathrm{site}}_{R/A}$ is an isomorphism in the $\infty$-category $\widehat{\calD}(A)$.
\end{lemma}

\begin{proof}
It follows from conditions $(a)$ and $(b)$ that the quotient ring $R/pR$ is semiperfect, so that the module of relative K\"{a}hler differentials $\Omega_{R/\overline{A}}$ is 
$p$-divisible. It follows that $R$ satisfies the stronger condition $(\ast^{+})$ of Variant \ref{variant:prismatic-coh-discrete}. As in the proof of Lemma \ref{lemma:really-easy-case}, we
are reduced to showing that the idempotent endomorphism $\iota_{R/A}$ of Remark \ref{remark:make-idempotent} is an isomorphism.
Let $A'$ be the $(p,I)$-completion of $A[ X_{s}^{1/p^{\infty}} ]$, endowed with the $\delta$-structure determined by the requirement that each $X_{s}^{1/p^n}$ is rank one.
Let $\widehat{R}$ denote the $p$-completion of $R$ (which agrees with the separated $p$-completion, since $R$ has bounded $p$-torsion).
It follows from condition $(a)$ that there is a commutative diagram of rings
$$ \xymatrix@R=50pt@C=50pt{ A/I \ar[d] \ar[r] & R \ar[d] \\
A'/IA' \ar[r]^-{u} & \widehat{R}, }$$
where $u$ satisfies $u(X_s) = x_s$. It follows from $(b)$ that $u$ is surjective modulo $p$ and therefore surjective (since both $A'/IA'$ and $\widehat{R}$ are $p$-complete).
Since $A'/IA'$ is $p$-completely formally \'{e}tale over $A/IA$, it follows from Corollaries \ref{corollary:etale-change-of-prism} and \ref{salmage}
that the induced map of prismatic complexes $\Prism_{R/A} \rightarrow \Prism_{\widehat{R}/A'}$ is invertible. It will therefore suffice to show that the idempotent $\iota_{ \widehat{R}/A'}$ is an isomorphism,
which follows from Lemma \ref{lemma:really-easy-case}.
\end{proof}

To apply Lemma \ref{lemma:easy-qrsp-case}, we will appeal to the fact that both $\Prism_{\bullet/A}$ and $\Prism^{\mathrm{site}}_{\bullet/A}$ have good descent properties.

\begin{lemma}\label{lemma:qs-site-descent}
Let $R$ be an $\overline{A}$-algebra which is $p$-complete
and has bounded $p$-torsion, and let $f: R \rightarrow R^{0}$ be a $p$-quasisyntomic cover. 
Let $R^{\bullet}$ denote the cosimplicial $\overline{A}$-algebra given by the $p$-completed tensor powers of $R^{0}$ over $R$.
Then the canonical map
$$ \Prism^{\mathrm{site}}_{R/A} \rightarrow \Tot( \Prism^{\mathrm{site}}_{ R^{\bullet}/A} )$$
is an isomorphism in the $\infty$-category $\widehat{\calD}(A)$.
\end{lemma}

\begin{proof}
Unwinding the definitions, we see that $\Tot( \Prism^{\mathrm{site}}_{ R^{\bullet}/A} )$ can be identified with the inverse limit
$\varprojlim_{(B,v) \in \calC} B$, where $\calC \subseteq (R/A)_{\Prism}$ is the full subcategory spanned by those
objects $(B,v)$ for which the unit map $R \rightarrow B/IB$ factors through $R^{0}$. It follows from Proposition~7.11 of \cite{prisms}
that (the opposite of) $\calC$ is a sieve with respect to the flat topology of Remark \ref{remark:flat-topology-on-relative-site}.
The desired result now follows from the observation that the functor
$$ (R/A)_{\Prism} \rightarrow \widehat{\calD}(A) \quad \quad (B,v) \mapsto B$$
is a sheaf for the flat topology.
\end{proof}

\begin{lemma}\label{lemma:qs-descent-relative}
Let $R$ be a $p$-complete $\overline{A}$-algebra which satisfies condition $(\ast)$ of Remark \ref{remark:hypothesis-for-discreteness}, and let $f: R \rightarrow R^{0}$ be a $p$-quasisyntomic cover.
Let $R^{\bullet}$ denote the cosimplicial $\overline{A}$-algebra given by the $p$-completed tensor powers of $R^{0}$ over $R$.
Then the canonical maps
$$ \Prism_{R/A} \rightarrow \Tot( \Prism_{ R^{\bullet}/A} ) \quad \quad \overline{\Prism}_{R/A} \rightarrow \Tot( \overline{\Prism}_{R^{\bullet}/A} )$$
are isomorphisms (in the $\infty$-categories $\widehat{\calD}(A)$ and $\widehat{\calD}( \overline{A})$, respectively).
\end{lemma}


\begin{proof}
We will prove the assertion for the relative Hodge-Tate complex; the analogous statement for the prismatic complex is a formal consequence.
Using Remark \ref{remark:hypothesis-for-discreteness}, we see that each of the complexes $\Fil_{n}^{\conj} \overline{\Prism}_{R/A}$
and $\Fil_{n}^{\conj} \overline{\Prism}_{R^{m}/A}$ is coconnective. Since the formation of totalizations of coconnective complexes commutes with the formation
of filtered colimits, it will suffice to show that each of the natural maps
$$\Fil_{n}^{\conj} \overline{\Prism}_{R/A} \rightarrow \Tot( \Fil_{n}^{\conj} \overline{\Prism}_{R^{\bullet}/A} )$$
is an isomorphism in $\widehat{\calD}(\overline{A} )$. Proceeding by induction on $n$, we can reduce further to showing that
each of the maps
$$\gr_{n}^{\conj} \overline{\Prism}_{R/A} \rightarrow \Tot( \gr_{n}^{\conj} \overline{\Prism}_{R^{\bullet}/A} )$$
is an isomorphism. Invoking the Hodge-Tate comparison (Remark \ref{remark:derived-HT-filtration}), we are reduced to checking that
each of the maps
$$ L \widehat{\Omega}^{n}_{R/\overline{A}} \simeq
\Tot( R^{\bullet}) \widehat{\otimes}_{R} L \widehat{\Omega}^{n}_{R/\overline{A}} )
\rightarrow \Tot( R^{\bullet} \widehat{\otimes}_{R} L \widehat{\Omega}^{n}_{R/\overline{A}} )
\rightarrow \Tot( L \widehat{\Omega}^{n}_{ R^{\bullet}/ \overline{A}} )$$
is invertible (here the first identification follows from the $p$-complete faithful flatness of the map $R \rightarrow R^{0}$). 
The cofiber of this map has a finite filtration, where each successive quotient can be realized as the totalization of a cosimplicial object of the form
$$L \widehat{\Omega}^{m}_{R^{\bullet}/R} \widehat{\otimes}^{L}_{R} L \widehat{\Omega}^{n-m}_{R/\overline{A}},$$
and therefore vanishes (see Lemma~2.6 of \cite{bhatt-completions}).
\end{proof}

\begin{proof}[Proof of Theorem \ref{theorem:site-theoretic-equivalence}]
Without loss of generality, we may assume that $\mathfrak{X} = \Spf(R)$ is affine, so that $R$ is a $p$-complete $\overline{A}$-algebra
satisfying condition $(\ast)$ of Remark \ref{remark:hypothesis-for-discreteness}.
Choose a $p$-quasisyntomic cover $R \rightarrow R^{0}$, where every element of $R^{0}$ has a $p$th root,
and let $R^{\bullet}$ be as in Lemma \ref{lemma:qs-descent-relative}. By construction, each of the $\overline{A}$-algebras
$R^{n}$ satisfies the hypotheses of Lemma \ref{lemma:easy-qrsp-case}, so that the comparison map $\xi_{ R^{n} }: \Prism_{R^{n}/A} \rightarrow \Prism^{\mathrm{site}}_{R^{n}/A}$
is an isomorphism. By functoriality, we have a commutative diagram
$$ \xymatrix@R=50pt@C=50pt@R=50pt@C=50pt{ \Prism_{R/A} \ar[r]^-{\xi_{R}}_{\sim} \ar[d] & \Prism_{R/A}^{\mathrm{site}} \ar[d] \\
\Tot( \Prism_{R^{\bullet}/A} ) \ar[r]^-{ \Tot( \xi_{R^{\bullet}})} & \Tot( \Prism^{\mathrm{site}}_{ R^{\bullet}/A} ),}$$
where the vertical maps are isomorphisms by virtue of Lemmas \ref{lemma:qs-site-descent} and
\ref{lemma:qs-descent-relative}. It follows that $\xi_{R}$ is also an isomorphism.
\end{proof}

\begin{corollary}\label{corollary:relative-initial-object}
Let $R$ be an $\overline{A}$-algebra which satisfies condition $(\ast^{+})$ of Variant \ref{variant:prismatic-coh-discrete}. Then the pair
$( \Prism_{R/A}, I \Prism_{R/A} )$ can be regarded as an initial object of the relative prismatic site $(R/A)_{\Prism}$.
\end{corollary}

\begin{proof}
Combine Theorem \ref{theorem:site-theoretic-equivalence} with Remark \ref{remark:make-idempotent}.
\end{proof}

\begin{example}\label{example:pushout-property}
Suppose that $(A,I)$ is a tranversal prism and let $(B,J)$ be a perfect prism with quotient $\overline{B} = B/J$. Then the tensor product
$R = \overline{A} \otimes^{L} \overline{B}$ satisfies condition $(\ast^{+})$ of Variant \ref{variant:prismatic-coh-discrete}. It follows that we can regard
$(\Prism_{R/A}, I \Prism_{R/A})$ as an initial object of the relative prismatic site $(R/A)_{\Prism}$. Note that if $(C,v)$ is any object of $(R/A)_{\Prism}$, the ring
homomorphism $\overline{B} \rightarrow C/IC$ determined by $v$ lifts {\em uniquely} to a $\delta$-ring homomorphism $B \rightarrow C$ (since the prism $(B,J)$ is perfect).
It follows that $(\Prism_{R/A}, I \Prism_{R/A} )$ can be identified with the coproduct of $(A,I)$ and $(B,J)$ in the category of bounded prisms (see Proposition \ref{proposition:transversal-coproduct}).
\end{example}

\begin{example}\label{example:transversal-prism-modulo-p}
Suppose that $(A,I)$ is a transversal prism (Definition \ref{definition:prism}) and let $R$ denote the quotient ring $\overline{A} / p \overline{A}$.
Then $R$ satisfies condition $(\ast^{+} )$ of Variant \ref{variant:prismatic-coh-discrete}, so we can
regard $( \Prism_{R/A}, I \Prism_{R/A} )$ as an initial object of the relative prismatic site $(R/A)_{\Prism}$. It follows from Example \ref{example:pushout-property}
that $( \Prism_{R/A}, I \Prism_{R/A} )$ is the universal {\em crystalline} prism over $(A,I)$ (that is, the coproduct of $(A,I)$ with $(\Z_p, (p) )$ in the category of prisms). In particular, we can identify $\Prism_{R/A}$ with the $\delta$-ring $A\{ \frac{p}{I} \}^{\wedge}$ described
in Corollary~3.14 of \cite{prisms}.
\end{example}

\begin{example}\label{example:q-divided-powers}
Let $(A,I) = ( \Z_p[[q-1]], [p]_q )$ denote the $q$-de Rham prism of Example \ref{example:q-prism}, let
$(B,IB)$ be the prism over $(A,I)$ described in Proposition \ref{proposition:universal-divided-power-element} (freely generated by
a rank $1$ unit $u$ satisfying $u^{p} \equiv 1 \pmod{IB}$). Let $C$ denote the $(p,I)$-completion of the
tensor product $B \otimes_{ \Z[u] } \Z[u^{1/p^{\infty}} ]$, so that $(C,IC)$ inherits the structure of a prism over $(A,I)$.
Set $\overline{A} = A/I$ and $R = \overline{A}[ u^{1/p^{\infty}} ] / (u^{p} - 1)$. Then:
\begin{itemize}
\item Let $(D,ID)$ be any prism over $(A,I)$. Then $\overline{A}$-algebra homomorphisms $R \rightarrow D/ID$ can be identified
with elements $w \in (D/ID)^{\flat \times}$ satisfying $(w^p)^{\sharp} = 1$ in $(D/ID)^{\times}$.

\item Every element $w \in (D/ID)^{\flat \times}$ lifts uniquely to an element $\overline{w} \in D^{\flat \times}$ (Proposition \ref{proposition:lift-flat}),
and the image $w^{\sharp} \in D^{\times}$ is automatically of rank $1$ (Proposition \ref{proposition:pth-roots-rank-1}).
\end{itemize}
It follows that $(C,IC)$ can be regarded as an initial object of the relative prismatic site $( R/A)_{\Prism}$, so that
we have a canonical isomorphism $\Prism_{R/A} \simeq C$. Applying Proposition \ref{proposition:q-divided-power-structure}, we see
that $C$ can be identified with the $(p,q-1)$-completion of the free $A$-module generated by elements of the form
$$ u^{\alpha} \gamma_{n,q}(u-1) = u^{\alpha} \frac{ (u-1) (u-q) \cdots (u-q^{n-1})}{ [1]_q [2]_q \cdots [n]_{q} },$$
where $n$ ranges over nonnegative integers and $\alpha$ ranges over the set $\{ \alpha \in \Z[1/p]: 0 \leq \alpha < 1 \}$.
\end{example}

\subsection{Absolute Prismatic Cohomology}\label{subsection:absolute-prismatic}

We now introduce an absolute variant of prismatic cohomology, which does not depend on a choice of base prism.

\begin{construction}[Prismatic Cohomology Sheaves]\label{construction:prismatic-cohomology-sheaves}
Let $R$ be an animated commutative ring. For every bounded prism $(A,I)$ with quotient ring $\overline{A} = A/I$, let
$\Prism_{ (\overline{A} \otimes^{L} R) / A}$ denote the prismatic complex of the derived tensor product $\overline{A} \otimes^{L} R$ relative to $A$
(Construction \ref{construction:relative-prismatic-cohomology}), which we regard as an object of the $(p,I)$-complete derived $\infty$-category $\widehat{\calD}(A)$.
For every morphism of bounded prisms $(A,I) \rightarrow (B,IB)$, Remark \ref{remark:change-of-prism} supplies a canonical isomorphism
$$ B \widehat{\otimes}^{L}_{A} \Prism_{ ( \overline{A} \otimes^{L} R) / A} \simeq \Prism_{ (\overline{B} \otimes^{L} R) / B}.$$
where $\overline{B}$ denotes the quotient $B/IB$. By virtue of Proposition \ref{proposition:DWCart-prism-description}, the construction
$(A,I) \mapsto \Prism_{ (\overline{A} \otimes^{L} R) / A}$ determines an object of the $\infty$-category $\calD(\WCart)$.
We will denote this object by $\mathscr{H}_{\Prism}(R)$ and refer to it as the {\it prismatic cohomology sheaf} of $R$.
\end{construction}

\begin{remark}\label{remark:sifted-colimit-prismatic-sheaf}
The construction $R \mapsto \mathscr{H}_{\Prism}(R)$ determines a functor of $\infty$-categories $\mathscr{H}_{\Prism}(\bullet): \CAlg^{\anim} \rightarrow \calD( \WCart )$ which
commutes with sifted colimits (this follows immediately from the analogous property of relative prismatic complexes; see Construction \ref{construction:relative-prismatic-cohomology}).
\end{remark}

\begin{remark}[\'{E}tale Descent]\label{remark:etale-descent-sheaf}
The functor $(R \in \CAlg^{\anim} ) \mapsto ( \mathscr{H}_{\Prism}(R) \in \calD( \WCart) )$ satisfies descent for the \'{e}tale topology (this follows immediately
from Corollary \ref{corollary:etale-descent-relative-prismatic}).
\end{remark}

\begin{remark}\label{remark:derived-descent-sheafy}
Let $\F_p^{\otimes \bullet+1}$ be the cosimplicial animated commutative ring introduced in Notation \ref{notation:F-p-bullet}, let
$R$ be an animated commutative ring, and let $R^{\bullet}$ denote the derived tensor product $R \otimes^{L} \F_p^{\otimes \bullet+1}$.
Then the canonical map $\mathscr{H}_{\Prism}(R) \rightarrow \Tot( \mathscr{H}_{\Prism}( R^{\bullet} ) )$ is an isomorphism in the $\infty$-category $\calD( \WCart )$.
This follows immediately from Proposition \ref{proposition:derived-descent-relative}
\end{remark}

\begin{remark}
Let $R$ be an animated commutative ring. Then the prismatic cohomology sheaf $\mathscr{H}_{\Prism}(R)$ can be identified with the inverse limit
$\varprojlim_{n} \mathscr{H}_{\Prism}( (\Z / p^{n} \Z) \otimes^{L} R)$, formed in the $\infty$-category $\calD( \WCart )$. This follows immediately
from Corollary \ref{corollary:mod-pn} (alternatively, it can be deduced formally from Remark \ref{remark:derived-descent-sheafy} by repeating the proof of
Corollary \ref{corollary:mod-pn}).
\end{remark}

\begin{variant}[Globalization to Formal Schemes]\label{variant:prismatic-sheaf-globalized}
Let $\mathfrak{X}$ be a bounded $p$-adic formal scheme which is quasi-compact and quasi-separated. Then the construction
$$ (A,I) \mapsto \RGamma_{ \Prism}( \Spf(A/I) \times \mathfrak{X} / A),$$
defined on the category of {\em transversal} prisms $(A,I)$, determines a quasi-coherent complex on the Cartier-Witt stack $\WCart$, which we will denote by $\mathscr{H}_{\Prism}( \mathfrak{X} )$ and
refer to as the {\it prismatic cohomology sheaf of $\mathfrak{X}$}. This construction is characterized by the following properties:
\begin{itemize}
\item The construction $\mathfrak{X} \mapsto \mathscr{H}_{\Prism}( \mathfrak{X} )$ satisfies descent for the \'{e}tale topology on the category of bounded $p$-adic formal schemes.
\item If $\mathfrak{X} = \Spf(R)$ is the spectrum of a $p$-complete commutative ring $R$ with bounded $p$-torsion, then
$\mathscr{H}_{\Prism}( \mathfrak{X} )$ agrees with the prismatic cohomology sheaf $\mathscr{H}_{\Prism}(R)$ of Construction \ref{construction:prismatic-cohomology-sheaves}
(see Corollary \ref{corollary:sashy}).
\end{itemize}
\end{variant}

\begin{remark}\label{remark:perfect-complex-absolute}
Let $\mathfrak{X}$ be a formal scheme which is smooth and proper over $\Spf(\Z_p)$. Then the prismatic cohomology sheaf
$\mathscr{H}_{\Prism}( \mathfrak{X} )$ of Variant \ref{variant:prismatic-sheaf-globalized} is a perfect complex on the Cartier-Witt stack $\WCart$
(see Remark \ref{remark:perfect-relative-prismatic}). More generally, if $k$ is a perfect field and $\mathfrak{X}$ is proper and smooth over $\Spf( W(k))$,
then $\mathscr{H}_{\Prism}( \mathfrak{X} )$ can be regarded as the direct image of a perfect complex on the product
$\WCart \times \Spf(W(k))$. 
\end{remark}

Let $(B,J)$ be a bounded prism with quotient $\overline{B} = B/J$, and let $\rho_{B}: \Spf(B) \rightarrow \WCart$ be the morphism of Construction \ref{construction:point-of-prismatic-stack}.
For every animated commutative ring $R$, we will identify the pullback $\rho_{B}^{\ast} \mathscr{H}_{\Prism}(R) \in \calD( \Spf(B) ) \simeq \widehat{\calD}(B)$ with
the relative prismatic complex $\Prism_{ ( \overline{B} \otimes^{L} R) / B}$. In particular, if $R$ is an animated commutative $\overline{B}$-algebra, then
the multiplication map $\overline{B} \otimes^{L} R \rightarrow R$ induces an $B$-linear map $\rho_{B}^{\ast} \mathscr{H}_{\Prism}(R) \rightarrow \Prism_{R/B}$, which
we can identify with a map $\xi: \mathscr{H}_{\Prism}(R) \rightarrow \rho_{B \ast} \Prism_{R/B}$ in the derived $\infty$-category $\WCart$.

\begin{proposition}\label{proposition:prismatic-sheaf-perfect-case}
Let $(B,J)$ be a perfect prism and let $R$ be an animated commutative algebra over the quotient ring $\overline{B} = B/J$. Then the preceding construction
determines an isomorphism $\xi: \mathscr{H}_{\Prism}(R) \xrightarrow{\sim} \rho_{B \ast} \Prism_{R/B}$ in the $\infty$-category $\calD( \WCart)$.
\end{proposition}

\begin{proof}
Fix a nonzero transversal prism $(A,I)$, and let $\rho_{A}: \Spf(A) \rightarrow \WCart$ be the morphism of Construction \ref{construction:point-of-prismatic-stack}.
Then $\rho_{A}$ is faithfully flat (Corollary \ref{corollary:flatness-v-transversality}). It will therefore suffice to show that $\rho_{A}^{\ast}(\xi)$ is an isomorphism
in $\widehat{\calD}(A)$. Let $(C,K)$ be a coproduct of $(A,I)$ with $(B,J)$ in the category of prisms, so that Proposition \ref{proposition:pullback-diagram-prisms}
supplies a pullback diagram
$$ \xymatrix@R=50pt@C=50pt{ \Spf(C) \ar[r] \ar[d] & \Spf(B) \ar[d]^{ \rho_B } \\
\Spf(A) \ar[r]^-{\rho_{A}} & \WCart. }$$
Setting $\overline{A} = A/I$ and $\overline{C} = C/K$ and invoking Remark \ref{remark:change-of-prism}, we can identify $\rho_{A}^{\ast}(\xi)$ with the composition of the bottom horizontal
maps in the diagram
$$ \xymatrix@R=50pt@C=50pt{ \Prism_{ \overline{A} \otimes^{L} \overline{B} / A } \ar[r] \ar[d] & \Prism_{ \overline{C} \otimes^{L} \overline{B} / C} \ar[r] \ar[d] & \Prism_{ \overline{C} / C } \ar[d] \\
\Prism_{ \overline{A} \otimes^{L} R / A } \ar[r] & \Prism_{ \overline{C} \otimes^{L} R / C } \ar[r] & \Prism_{ \overline{C} \otimes^{L}_{\overline{B} } R / C }}$$
of commutative algebra objects of $\widehat{\calD}(A)$. Note that the squares on the left and right are pushout diagrams by virtue of
Remarks \ref{remark:change-of-prism} and \ref{remark:Kunneth-relative-prismatic}, respectively. It will therefore suffice to show that the upper horizontal composition
$\Prism_{ \overline{A} \otimes^{L} \overline{B} / A } \rightarrow \Prism_{\overline{C} /C } \simeq C$ is an isomorphism, which follows from Example \ref{example:pushout-property}.
\end{proof}

\begin{corollary}
Let $(A,I)$ be a perfect prism, let $\overline{A}$ denote the quotient ring $A/I$, and let $\mathfrak{X}$ be a bounded $p$-adic formal $\overline{A}$-scheme which 
is quasi-compact and quasi-separated. Then the prismatic cohomology sheaf $\mathscr{H}_{\Prism}( \mathfrak{X} )$ of Variant \ref{variant:prismatic-sheaf-globalized}
can be identified with the pushforward $\rho_{A \ast} \RGamma_{\Prism}( \mathfrak{X} / A)$, where $\rho_{A}: \Spf(A) \rightarrow \WCart$ is the morphism of
Construction \ref{construction:point-of-prismatic-stack}.
\end{corollary}

\begin{construction}[Absolute Prismatic Complexes: Affine Case]\label{construction:absolute-prismatic-cohomology-general}
Let $R$ be an animated commutative ring. For every pair of integers $m$ and $n$, we let $\Prism_{R}^{[m]}\{n\} \in \widehat{\calD}(\Z_p)$ denote complex
$\RGamma( \WCart, \mathcal{I}^{\otimes m} \otimes \mathscr{H}_{\Prism}(R)\{n\} )$, where $\mathscr{I}$ is the Hodge-Tate ideal sheaf of
Example \ref{example:HTI} and $\mathscr{H}_{\Prism}(R)$ is the prismatic cohomology sheaf of Construction \ref{construction:prismatic-cohomology-sheaves}.
In the special case $m=0$, we denote $\Prism_{R}^{[m]}\{n\}$ by $\Prism_{R}\{n\}$.  In the special case $n=0$, we denote $\Prism_{R}\{n\}$ simply by $\Prism_{R}$, and refer to it
as the {\it absolute prismatic complex of $R$}.
\end{construction}

\begin{remark}\label{remark:tautological-calculation-absolute-prismatic}
Let $(A^{\bullet}, I^{\bullet} )$ be the cosimplicial prism of Notation \ref{notation:simplicial-prism}, and let $\overline{A}^{\bullet}$ be the cosimplicial commutative ring $A^{\bullet} / I^{\bullet}$.
For every animated commutative ring $R$, the absolute prismatic complex $\Prism_{R}$ can be realized as the totalization of the cosimplicial complex
$\Prism_{ ( \overline{A}^{\bullet} \otimes^{L} R) / A^{\bullet} }$ (see Corollary \ref{corollary:global-sections-on-WCart}). In \S\ref{subsection:concrete-qdR}, we will
give a more efficient description of the $\Prism_{R}$, at least when $p$ is an odd prime (Proposition \ref{proposition:concrete-qdR}).
\end{remark}

Let $(A,I)$ be a bounded prism and let $R$ be an animated commutative algebra over the quotient ring $A/I$. For every pair of integers $m$ and $n$, we have
an evident comparison map
$$ \Prism^{[m]}_{R}\{n\} = \RGamma( \WCart, \mathcal{I}^{\otimes m} \otimes \mathscr{H}_{\Prism}(R)\{n\} )
\rightarrow \rho_{A}^{\ast}( \mathcal{I}^{\otimes m} \otimes \mathscr{H}_{\Prism}(R)\{n\} ) =
 I^{m} \Prism_{R/A}\{n\}.$$

\begin{proposition}\label{proposition:absolute-vs-relative}
Let $(A,I)$ be a perfect prism and let $R$ be an animated commutative algebra over the quotient ring $A/I$.
For every pair of integers $m$ and $n$, the comparison map
$$ \Prism^{[m]}_{R}\{n\} \rightarrow I^{m} \Prism_{ R / A}\{n\}$$
is an isomorphism in the $\infty$-category $\widehat{\calD}(A)$.
\end{proposition}

\begin{proof}
Apply Proposition \ref{proposition:prismatic-sheaf-perfect-case}.
\end{proof}

\begin{example}\label{example:prismatic-sheaf-qrsp}
Let $R$ be a quasiregular semiperfectoid ring. Then there exists a perfect prism $(A,I)$ and a ring homomorphism $\overline{A} = A/I \rightarrow R$.
It follows from Corollary \ref{corollary:relative-initial-object} that the relative prismatic complex $\Prism_{R/A}$ is concentrated in cohomological degree zero, and the pair
$( \Prism_{R/A}, I \Prism_{R/A} )$ can be viewed as prism over $(A,I)$ (which is an initial object of the relative prismatic site $(R/A)_{\Prism}$). 
Applying Proposition \ref{proposition:absolute-vs-relative}, we obtain isomorphisms
$$ \Prism_{R} \simeq \Prism_{R/A}  \quad \quad \Prism_{R}^{[1]} \simeq I \Prism_{R/A};$$
in particular, we can regard $( \Prism_{R}, \Prism_{R}^{[1]} )$ as a prism. Moreover, the prismatic sheaf $\mathscr{H}_{\Prism}(R)$ can be identified
with the direct image of the structure sheaf of the formal scheme $\Spf( \Prism_{R} )$ along the map
$\Spf( \Prism_R ) \rightarrow \WCart$ given by Construction \ref{construction:point-of-prismatic-stack}.
\end{example}

\begin{proposition}\label{proposition:absolute-etale-affine}
For every pair of integers $m$ and $n$, the functor
$$ \CAlg^{\anim} \rightarrow \widehat{\calD}(\Z_p) \quad \quad R \mapsto \Prism_{R}^{[m]}\{n\}$$
satisfies descent for the \'{e}tale topology.
\end{proposition}

\begin{proof}
Apply Remark \ref{remark:etale-descent-sheaf}.
\end{proof}

\begin{proposition}\label{proposition:derived-descent-absolute}
Let $\F_p^{\otimes \bullet+1}$ be the cosimplicial animated commutative ring introduced in Notation \ref{notation:F-p-bullet}.
Let $R$ be an animated commutative ring, and regard $R^{\bullet} = R \otimes^{L} \F_p^{\otimes \bullet+1}$ as a cosimplicial object of
$\CAlg^{\anim}$. For every pair of integers $m$ and $n$, the tautological map 
$$ \Prism^{[m]}_{R}\{n\} \rightarrow \Tot( \Prism^{[m]}_{ R^{\bullet} }\{n\} )$$
is an isomorphism in the $\infty$-category $\widehat{\calD}(\Z_p)$.
\end{proposition}

\begin{proof}
Apply Remark \ref{remark:derived-descent-sheafy}.
\end{proof}

\begin{remark}
In the situation of Proposition \ref{proposition:derived-descent-absolute}, each individual term $R^{k}$ of the cosimplicial animated commutative ring $R^{\bullet}$
admits the structure of animated commutative $\F_p$-algebra, so that the absolute prismatic complex $\Prism_{ R^{k} }$ computes the (derived) crystalline cohomology of
$R^{k}$ (see Remark \ref{remark:absolute-as-derived-crystalline}). Consequently, Proposition \ref{proposition:derived-descent-absolute} provides a mechanism for using crystalline cohomology
to prove results about prismatic cohomology, which we will exploit in \S\ref{section:first-chern}.
\end{remark}

\begin{corollary}\label{corollary:reduce-mod-n}
Let $R$ be an animated commutative ring. For every pair of integers $m$ and $n$, the tautological map
$$ \Prism^{[m]}_{R}\{n\} \rightarrow \varprojlim_{k}  \Prism^{[m]}_{(\Z / p^{k} \Z) \otimes^{L} R}\{n\}$$
is an isomorphism in the $\infty$-category $\widehat{\calD}(\Z_p)$.
\end{corollary}

\begin{corollary}\label{corollary:absolute-completion-invariance}
Let $R$ be an animated commutate ring and let $\widehat{R}$ denote the $p$-completion of $R$.
For every pair of integers $m$ and $n$, the tautological map $\Prism^{[m]}_{R}\{n\} \rightarrow \Prism^{[m]}_{ \widehat{R}} \{n\}$ is an isomorphism in the
$\infty$-category $\widehat{\calD}(\Z_p)$.
\end{corollary}

We now introduce a global variant of Construction \ref{construction:absolute-prismatic-cohomology-general}.

\begin{construction}\label{construction:absolute-prismatic-complex-of-scheme}
Let $X$ be scheme, formal scheme, or algebraic stack. For every pair of integers $m$ and $n$, we let $\RGamma_{\Prism}^{[m]}(X)\{n\}$ denote the inverse limit $\varprojlim_{ \Spec(R) \rightarrow X} \Prism_{R/A}^{[m]}\{n\}$, formed in the $\infty$-category $\widehat{\calD}(\Z_p)$. Here the limit is indexed by the {\it category of points} of $X$: that is, the category of pairs $(R, f)$ where $R$ is a commutative ring and $f: \Spec(R) \rightarrow X$ is a morphism. In the special case $m =0$, we denote $\RGamma_{\Prism}^{[m]}(X)\{n\}$ by $\RGamma_{\Prism}(X)\{n\}$. 
In the special case $n=0$, we denote $\RGamma_{\Prism}(X)\{n\}$ simply by $\RGamma_{\Prism}(X)$, and refer to it as the {\it absolute prismatic complex of $X$}.
We denote the cohomology groups of $\RGamma_{\Prism}(X)$ by $\mathrm{H}^{\ast}_{\Prism}(X)$ and refer to them as the {\it absolute prismatic cohomology groups} of $X$.
\end{construction}

\begin{example}
Let $\mathfrak{X}$ be a bounded $p$-adic formal scheme which is quasi-compact and quasi-separated. For every pair of integers $m$ and $n$, we have a canonical isomorphism
$$ \RGamma_{\Prism}^{[m]}( \mathfrak{X} )\{n\} \simeq \RGamma( \WCart, \calI^{m} \otimes \mathscr{H}_{\Prism}(\mathfrak{X})\{n\}),$$
where $\mathscr{H}_{\Prism}( \mathfrak{X} )$ is the prismatic cohomology sheaf of $\mathfrak{X}$ defined in Variant \ref{variant:prismatic-sheaf-globalized}.
\end{example}

\begin{remark}
Let $(A,I)$ be a bounded prism and let $X$ be a scheme, formal scheme, or algebraic stack which is defined over the quotient ring $\overline{A} = A/I$.
For every pair of integers $m$ and $n$, there is a tautological comparison map
$$ \RGamma_{\Prism}^{[m]}( X )\{n\} \rightarrow I^{m}\{n\} \otimes^{L}_{A} \RGamma_{\Prism}(X/A),$$
where $\RGamma_{\Prism}(X/A)$ is the prismatic complex of $X$ relative to $A$ (Construction \ref{construction:prismatic-complex-of-scheme}).
If the prism $(A,I)$ is perfect, then this comparison map is an isomorphism in the $\infty$-category $\widehat{\calD}(\Z_p)$ (see Proposition \ref{proposition:absolute-vs-relative}).
\end{remark}

\begin{remark}\label{remark:etale-descent-absolute}
For every pair of integers $m$ and $n$, the construction $X \mapsto \RGamma_{\Prism}^{[m]}(X)\{n\}$ satisfies descent with respect to the \'{e}tale topology 
(see Proposition \ref{proposition:absolute-etale-affine}).
\end{remark}

\begin{example}\label{example:sotto-absolute}
Let $X = \Spec(R)$ be an affine scheme. For every pair of integers $m$ and $n$, the complex $\RGamma_{\Prism}^{[m]}(X)\{n\}$ can be identified
with the absolute prismatic complex $\Prism_{R}^{[m]}\{n\}$ of Construction \ref{construction:absolute-prismatic-cohomology-general}.
\end{example}

\begin{example}\label{example:absolute-compute-for-formal-scheme}
Let $\mathfrak{X}$ be a bounded $p$-adic formal scheme. For each $k \geq 0$, let $\mathfrak{X}_{k} \subseteq \mathfrak{X}$
denote the closed subscheme given by the vanishing locus of $p^{k}$, so that we have $\mathfrak{X} = \varinjlim_{k} \mathfrak{X}_{k}$ (as set-valued functors on the category of commutative rings).
It follows that for every pair of integers $m$ and $n$, the canonical map $\RGamma_{\Prism}^{[m]}( \mathfrak{X} )\{n\} \rightarrow \varprojlim_{k} \RGamma_{\Prism}^{[m]}( \mathfrak{X}_{k} )\{n\}$
is an isomorphism in $\widehat{\calD}(\Z_p)$.
\end{example}

\begin{proposition}\label{proposition:scheme-vs-formal-absolute}
Let $X$ be a scheme for which the structure sheaf $\calO_{X}$ has bounded $p$-power torsion, and let $\mathfrak{X} = \Spf(\Z_p) \times X$ be the associated $p$-adic formal scheme.
For every pair of integers $m$ and $n$, the restriction map
$\RGamma_{\Prism}^{[m]}(X)\{n\} \rightarrow \RGamma_{\Prism}^{[m]}(\mathfrak{X})\{n\}$ is an isomorphism in the $\infty$-category $\widehat{\calD}(\Z_p)$.
\end{proposition}

\begin{proof}
Using Remark \ref{remark:etale-descent-absolute}, we can reduce to the case where $X = \Spec(R)$ is the spectrum of a commutative ring $R$ with bounded $p$-torsion.
In this case, the desired result follows by combining Example \ref{example:absolute-compute-for-formal-scheme}, Example \ref{example:sotto-absolute}, and Corollary \ref{corollary:reduce-mod-n}.
\end{proof}

\begin{corollary}\label{corollary:sashy-absolute}
Let $R$ be a commutative ring having bounded $p$-power torsion, and let $\widehat{R}$ be the $p$-completion of $R$. For every pair of integers $m$ and $n$,
the tautological map $\Prism_{R}^{[m]}\{n\} \rightarrow \RGamma_{\Prism}^{[m]}( \Spf(\widehat{R} )  )\{n\}$ is an isomorphism in $\widehat{\calD}(\Z_p)$.
\end{corollary}

In good cases, the absolute prismatic cohomology of a $p$-adic formal scheme $\mathfrak{X}$ can be computed using a variant of
Definition \ref{definition:relative-prismatic-site}.

\begin{definition}[Remark~4.6 of \cite{prisms}]\label{definition:absolute-prismatic-site}
Let $\mathfrak{X}$ be a bounded $p$-adic formal scheme. We define a category $(\mathfrak{X})_{\Prism}$ as follows:
\begin{itemize}
\item The objects of $( \mathfrak{X} )_{\Prism}$ are triples $(A,I,u)$, where $(A,I)$ is a bounded prism and $u: \Spf(A/I) \rightarrow \mathfrak{X}$ is a morphism of $p$-adic formal schemes.
\item A morphism from $(A,I,u)$ to $(B,J,v)$ in $(\mathfrak{X})_{\Prism}$ is a morphism of prisms $(A,I) \rightarrow (B,J)$ for which the composite map
$$ \Spf(B/J) \rightarrow \Spf(A/I) \xrightarrow{u} \mathfrak{X}$$
is equal to $v$.
\end{itemize}
We refer to $(\mathfrak{X})_{\Prism}$ as the {\it absolute prismatic site} of $\mathfrak{X}$. For every pair of integers $m$ and $n$, we let
$\RGamma_{\Prism}^{\mathrm{site},[m]}( \mathfrak{X} )\{n\}$ denote the inverse limit
$$ \varprojlim_{ (A,I,u) \in (\mathfrak{X})_{\Prism} } I^{m}\{n\},$$
formed in the $\infty$-category $\widehat{\calD}(\Z_p)$.
\end{definition}

\begin{remark}\label{remark:flat-topology-on-absolute}
Let $\mathfrak{X}$ be a bounded $p$-adic formal scheme. We define the {\it flat topology} on the category
$(\mathfrak{X})^{\op}_{\Prism}$ to be the Grothendieck topology generated by those finite collections of morphisms
$$ \{ (A,I,u) \rightarrow (A_{s}, I_s, u_s) \}_{s \in S}$$
for which the induced ring homomorphism $A \rightarrow \prod_{s \in S} A_s$ is $(p,I)$-completely faithfully flat.
Note that, for every pair of integers $m$ and $n$, the functor $(A,I,u) \mapsto I^{m}\{n\}$ is a $\widehat{\calD}(\Z_p)$-valued sheaf with respect to the flat topology.
\end{remark}

Let $\mathfrak{X}$ be a $p$-adic formal scheme. For every object $(A,I,u)$ of the absolute prismatic site $(\mathfrak{X})_{\Prism}$, restriction along $u$ determines comparison maps
$$ \RGamma^{[m]}_{\Prism}(\mathfrak{X})\{n\} \rightarrow \RGamma^{[m]}_{\Prism}( \Spf(A/I) )\{n\}
\simeq \Prism_{\overline{A}}^{[m]}\{n\} \rightarrow I^{m}\{n\} \otimes^{L}_{A} \Prism_{\overline{A}/A} \simeq I^{m}\{n\},$$
where $\overline{A}$ denotes the quotient ring $A/I$. Passing to the inverse limit over $(\mathfrak{X})_{\Prism}$, we obtain a comparison map
$$ \xi: \RGamma^{[m]}_{\Prism}(\mathfrak{X})\{n\} \rightarrow \RGamma_{\Prism}^{\mathrm{site},[m]}( \mathfrak{X} )\{n\}.$$

\begin{example}\label{example:prismatic-cohomological-qrsp}
Let $R$ be a quasiregular semiperfectoid ring, and let $\mathfrak{X} = \Spf(R)$ be the associated $p$-adic formal scheme.
Then the category $(\mathfrak{X})_{\Prism}$ has an initial object $(A,I,u)$, where
$(A,I)$ is given by the absolute prismatic cohomology $( \Prism_{R}, \Prism_{R}^{[1]} )$ (see Example \ref{example:prismatic-sheaf-qrsp}).
It follows that preceding construction supplies an isomorphism
$$ \xi: I^m\{n\} \simeq \RGamma^{[m]}_{\Prism}(\mathfrak{X})\{n\} \rightarrow \RGamma_{\Prism}^{\mathrm{site},[m]}( \mathfrak{X} )\{n\}$$
for every pair of integers $m$ and $n$.
\end{example}

We will prove the following generalization of Example \ref{example:prismatic-cohomological-qrsp}:

\begin{theorem}\label{theorem:absolute-compute-with-stack}
Let $\mathfrak{X}$ be a bounded $p$-adic formal scheme. Assume that, for every affine open subset $\Spf(R) \subseteq \mathfrak{X}$, the
coordinate ring $R$ is $p$-quasisyntomic. Then, for every pair of integers $m$ and $n$, the comparison map
$$ \xi: \RGamma^{[m]}_{\Prism}(\mathfrak{X})\{n\} \rightarrow \RGamma_{\Prism}^{\mathrm{site},[m]}( \mathfrak{X} )\{n\}.$$
is an isomorphism in the $\infty$-category $\widehat{\calD}(\Z_p)$.
\end{theorem}

The proof of Theorem \ref{theorem:absolute-compute-with-stack} will require some preliminaries.

\begin{proposition}\label{proposition:site-theoretic-descent}
Let $\calC$ be the category of bounded $p$-adic formal schemes. For every pair of integers $m$ and $n$, the functor
$$ \calC^{\op} \rightarrow \widehat{\calD}(\Z_p) \quad \quad \mathfrak{X} \mapsto \RGamma_{\Prism}^{[m], \mathrm{site}}( \mathfrak{X} )\{n\}$$
satisfies descent with respect to the $p$-quasisyntomic topology on $\calC$.
\end{proposition}

\begin{proof}
Let $\pi: \mathfrak{X}_{0} \rightarrow \mathfrak{X}$ be a $p$-quasisyntomic covering of $p$-adic formal schemes, and let
$\mathfrak{X}_{\bullet}$ be the simplicial formal scheme given by the fiber powers of $\mathfrak{X}_0$ relative to $\mathfrak{X}$.
We wish to show that the tautological map
$$ \rho: \RGamma_{\Prism}^{[m], \mathrm{site}}( \mathfrak{X} )\{n\} \rightarrow \Tot( 
\RGamma_{\Prism}^{[m], \mathrm{site}}( \mathfrak{X}_{\bullet} )\{n\})$$
is an isomorphism in the $\infty$-category $\widehat{\calD}(\Z_p)$. Unwinding the definitions, we see that the target of
$\rho$ can be identified with the limit $\varprojlim_{ (A,I,u) \in ( \mathfrak{X} )'_{\Prism} } I^{m}\{n\}$, where
$( \mathfrak{X} )'_{\Prism}$ is the full subcategory of $( \mathfrak{X} )_{\Prism}$ spanned by those triples 
$(A,I,u)$ for which the map $u: \Spf(A/I) \rightarrow \mathfrak{X}$ factors through $\pi$. It will therefore suffice to show that the subcategory
$( \mathfrak{X} )'_{\Prism}$ is a covering sieve with respect to the flat topology of Remark \ref{remark:flat-topology-on-absolute}, which follows
from Proposition~7.11 of \cite{prisms}.
\end{proof}

\begin{proposition}\label{proposition:absolute-quasisyntomic-descent}
Let $\calC'$ be the category of $p$-quasisyntomic $p$-adic formal schemes. For every pair of integers $m$ and $n$, the functor
$$ \calC'^{\op} \rightarrow \widehat{\calD}(\Z_p) \quad \quad \mathfrak{X} \mapsto \RGamma_{\Prism}^{[m]}( \mathfrak{X} )\{n\}$$
satisfies descent with respect to the $p$-quasisyntomic topology on $\calC'$.
\end{proposition}

\begin{proof}
By virtue of Remark \ref{remark:etale-descent-absolute}, the functor $\mathfrak{X} \mapsto \RGamma_{\Prism}^{[m]}( \mathfrak{X} )\{n\}$
satisfies \'{e}tale descent; it will therefore suffice to prove that it satisfies $p$-quasisyntomic descent when restricted to the full subcategory of
$\calC'$ spanned by the {\em affine} $p$-quasisyntomic formal schemes. Let $R$ be a $p$-quasisyntomic ring, let $f: R \rightarrow R^{0}$ be a $p$-quasisyntomic ring homomorphism which is $p$-completely faithfully flat, and let $R^{\bullet}$ be the cosimplicial commutative ring given by the $p$-completion of the iterated tensor powers of $R^{0}$ over $R$. We wish to show that the canonical map
$\theta: \Prism_{R}^{[m]}\{n\} \rightarrow \Tot( \Prism_{R^{\bullet} }^{[m]}\{n\} )$ is an isomorphism in $\widehat{\calD}(\Z_p)$. It follows from
Lemma \ref{lemma:qs-descent-relative} that the canonical map $\mathscr{H}_{\Prism}(R) \rightarrow \Tot( \mathscr{H}_{\Prism}(R^{\bullet} ))$ is an isomorphism in the
$\infty$-category $\calD( \WCart )$; the desired result now follows by tensoring both sides with the invertible sheaf $\calI^{m}\{n\}$ and passing to global sections.
\end{proof}

\begin{proof}[Proof of Theorem \ref{theorem:absolute-compute-with-stack}]
By virtue of Propositions \ref{proposition:site-theoretic-descent} and \ref{proposition:absolute-quasisyntomic-descent}, we may assume without loss of generality that $\mathfrak{X} = \Spf(R)$
is the formal spectrum of a quasiregular semiperfectoid ring $R$, in which case the desired result follows from Example \ref{example:prismatic-cohomological-qrsp}.
\end{proof}

\subsection{Absolute Hodge-Tate Cohomology}\label{subsection:absolute-HT}

For every animated commutative ring $R$, the absolute prismatic complex $\Prism_{R}$ is equipped with a decreasing filtration, given by the complexes
$\{ \Prism_{R}^{[n]} \}_{n \geq 0}$. We now study the associated graded of this filtration.

\begin{construction}[Hodge-Tate Cohomology Sheaves]\label{construction:absolute-Hodge-Tate-filtration}
Let $R$ be an animated commutative ring and let $n$ be an integer. For every bounded prism $(A,I)$ with quotient $\overline{A} = A/I$,
let $\Fil_{n}^{\conj} \overline{\Prism}_{ (\overline{A} \otimes^{L} R) / A}$ be the $n$th stage of the conjugate filtration of
the relative Hodge-Tate complex $\overline{\Prism}_{ (\overline{A} \otimes^{L} R) / A}$, which we regard as an object of the $p$-complete
derived $\infty$-category $\widehat{\calD}( \overline{A} )$. By virtue of Remark \ref{remark:flashcard}, the construction $(A,I) \mapsto \Fil_{n}^{\conj} \overline{\Prism}_{ (\overline{A} \otimes^{L} R) / A}$
can be identified with a quasi-coherent complex on the Hodge-Tate divisor $\WCart^{\mathrm{HT}}$, which we will denote by $\Fil_{n}^{\conj} \mathscr{H}_{\overline{\Prism} }(R)$.

Allowing $n$ to vary, we obtain a direct system
$$ \Fil_{0}^{\conj} \mathscr{H}_{ \overline{\Prism} }(R) \rightarrow
\Fil_{1}^{\conj} \mathscr{H}_{ \overline{\Prism} }(R) \rightarrow
\Fil_{2}^{\conj} \mathscr{H}_{ \overline{\Prism} }(R) \rightarrow \cdots$$
in the $\infty$-category $\calD( \WCart^{\mathrm{HT}} )$. We denote the colimit of this system by $\mathscr{H}_{ \overline{\Prism} }(R)$ and refer to it as the
{\it Hodge-Tate cohomology sheaf of $R$}, and we will refer to the direct system $\Fil_{\bullet}^{\conj}  \mathscr{H}_{ \overline{\Prism} }(R)$ as the {\it conjugate filtration} on
$\mathscr{H}_{\overline{\Prism}}(R)$.
\end{construction}

\begin{remark}
Let $R$ be an animated commutative ring. Then the Hodge-Tate cohomology sheaf $\mathscr{H}_{\overline{\Prism}}(R)$ can be identified with the restriction
$\mathscr{H}_{\Prism}(R)|_{ \WCart^{\mathrm{HT} } }$, where $\mathscr{H}_{\Prism}(R) \in \calD( \WCart )$ is the prismatic cohomology sheaf of Construction \ref{construction:prismatic-cohomology-sheaves}.
\end{remark}

\begin{remark}\label{remark:sheafy-HT-comparison}
Let $R$ be an animated commutative ring. For every integer $n$, let $\gr_{n}^{\conj} \mathscr{H}_{\overline{\Prism}}(R)$ denote the cofiber of the map
$\Fil_{n-1}^{\conj} \mathscr{H}_{\overline{\Prism}}(R) \rightarrow \Fil_{n}^{\conj} \mathscr{H}_{\overline{\Prism}}(R)$, formed in the $\infty$-category $\calD( \WCart^{\mathrm{HT}} )$. 
The Hodge-Tate comparison of Remark \ref{remark:derived-HT-filtration} then supplies a canonical isomorphism
$$\gr_{n}^{\conj} \mathscr{H}_{\overline{\Prism}}(R) \simeq L \Omega^{n}_{R} \otimes \mathcal{O}_{ \WCart^{\mathrm{HT}} }\{-n\}[-n],$$
where $L \Omega^{n}_{R}$ is the complex of Construction \ref{construction:exterior-of-cotangent}.
\end{remark}

\begin{remark}
Let $R$ be an animated commutative ring. Then $\Fil_{\bullet}^{\conj} \mathscr{H}_{ \overline{\Prism} }(R)$ can be regarded as a commutative algebra object in
the $\infty$-category of filtered objects of $\calD( \WCart^{\mathrm{HT}} )$. Moreover, Remark \ref{remark:sheafy-HT-comparison} supplies an isomorphism
$$ \Fil_{0}^{\conj} \mathscr{H}_{ \overline{\Prism}}(R) \simeq R \otimes \mathcal{O}_{ \WCart^{\mathrm{HT} }}.$$
It follows that each $\Fil_{n}^{\conj} \mathscr{H}_{\overline{\Prism}}(R)$ can be regarded as an $R$-module object of the $\infty$-category $\widehat{\calD}( \WCart^{\mathrm{HT}} )$.
\end{remark}

\begin{construction}\label{construction:absolute-HT}
Let $R$ be an animated commutative ring. For each integer $n$, we let $\overline{\Prism}_{R}\{n\}$ denote the complex $\RGamma( \WCart^{\mathrm{HT}}, \mathscr{H}_{\overline{\Prism}}(R)\{n\} )$,
which we regard as an object of the $p$-complete derived $\infty$-category $\widehat{\calD}(R)$. In the special case $n =0$, we will denote $\overline{\Prism}_{R}\{n\}$ by $\overline{\Prism}_{R}$
and refer to it as the {\it absolute Hodge-Tate complex of $R$}. 

For every pair of integers $m$ and $n$, we let $\Fil_{m}^{\conj} \overline{\Prism}_{R}\{n\}$ denote the complex 
$$\RGamma( \WCart^{\mathrm{HT}}, \Fil_{m}^{\conj} \mathscr{H}_{\overline{\Prism}}(R)\{n\} ).$$
Allowing $m$ to vary, we obtain a direct system
$$ \Fil_{0}^{\conj} \overline{\Prism}_{R}\{n\} \rightarrow \Fil_{1}^{\conj} \overline{\Prism}_{R}\{n\} \rightarrow \Fil_{2}^{\conj}  \overline{\Prism}_{R}\{n\} \rightarrow \cdots$$
whose colimit can be identified with $\overline{\Prism}_{R}\{n\}$ (Corollary \ref{corollary:HT-computation}). We will denote this direct system by 
$\Fil_{\bullet}^{\conj} \overline{\Prism}_{R}\{n\}$ and refer to it as the {\it conjugate filtration} on $\overline{\Prism}_{R}\{n\}$.
\end{construction}

\begin{example}\label{example:relative-HT-over-perfect}
Let $(A,I)$ be a perfect prism and let $R$ be an animated commutative algebra over the quotient ring $\overline{A} = A/I$,
and let $\rho_{A}^{\mathrm{HT}}: \Spf(A/I) \rightarrow \WCart^{\mathrm{HT}}$ be the morphism of Remark \ref{remark:HT-point-of-prismatic-stack}.
Then Proposition \ref{proposition:prismatic-sheaf-perfect-case} supplies an isomorphism 
$$\mathscr{H}_{\overline{\Prism}}(R) \simeq (\rho_{A}^{\mathrm{HT}})_{\ast} \overline{\Prism}_{R/A}$$
in the $\infty$-category $\calD( \WCart^{\mathrm{HT} } )$. Tensoring with $\mathcal{O}_{\WCart^{\mathrm{HT}}}\{n\}$ and passing to global sections, we obtain canonical isomorphisms
$\overline{\Prism}_{R}\{n\} \xrightarrow{\sim} \overline{\Prism}_{R/A}\{n\}$.
\end{example}

\begin{remark}[Relationship with Absolute Prismatic Cohomology]\label{remark:diffracted-vs-prismatic1}
Let $\iota: \WCart^{\mathrm{HT}} \hookrightarrow \WCart$ be the inclusion of the Hodge-Tate divisor and let $R$ be an animated commutative ring. For every pair of integers $m$ and $n$,
we have a canonical fiber sequence
$$ \mathcal{I}^{m+1} \mathscr{H}_{\Prism}(R)\{n\} \rightarrow \mathcal{I}^{m} \mathscr{H}_{\Prism}(R)\{n\} \rightarrow \iota_{\ast}(\mathscr{H}_{\overline{\Prism}}(R)\{m+n\})$$
in the $\infty$-category $\calD(\WCart)$. Passing to global sections, we obtain a fiber sequence
$$ \Prism_{R}^{[m+1]}\{n\} \rightarrow \Prism_{R}^{[m]}\{n\} \rightarrow \overline{\Prism}_{R}\{m+n\}$$
in the $\infty$-category $\widehat{\calD}(\Z_p)$. 
\end{remark}

\begin{proposition}\label{proposition:HT-sifted-colimit}
For every pair of integers $m$ and $n$, the functor
$$ \CAlg^{\anim} \rightarrow \widehat{\calD}(\Z_p) \quad \quad R \mapsto \Fil_{m}^{\conj} \overline{\Prism}_{R}\{n\}$$
commutes with sifted colimits. In particular, the functor $R \mapsto \overline{\Prism}_{R}\{n\}$ commutes with sifted colimits.
\end{proposition}

\begin{proof}
It follows immediately from the definition that the functor 
$$ \CAlg^{\anim} \rightarrow \calD( \WCart^{\mathrm{HT}} ) \quad \quad R \mapsto \Fil_{m}^{\conj} \mathscr{H}_{\overline{\Prism}}(R)\{n\}$$ commutes
with sifted colimits. Proposition \ref{proposition:HT-sifted-colimit} now follows from the fact that the global sections functor
$\RGamma: \calD( \WCart^{\mathrm{HT}} ) \rightarrow \widehat{\calD}(\Z_p)$ preserves colimits (Corollary \ref{corollary:HT-computation}).
\end{proof}

Beware that, when regarded as a functor from $\CAlg^{\anim}$ to $\widehat{\calD}(\Z_p)$, the functor $R \mapsto \Prism_{R}$ does {\em not} commute with filtered colimits. We can salvage the situation by regarding $\Prism_{R}$ as a filtered complex.

\begin{example}\label{example:filtered-derived-2}
Let $R$ be an animated commutative ring. For every integer $n$, the diagram
$$ \cdots \rightarrow \Prism_{R}^{[2]}\{n\} \rightarrow \Prism_{R}^{[1]}\{n\} \rightarrow \Prism_{R}^{[0]}\{n\} \rightarrow \Prism_{R}^{[-1]}\{n\} \rightarrow \Prism_{R}^{[-2]}\{n\} \rightarrow 0$$
determines an object of $\DFiltI(\Z_p)$, which we will denote by $\Prism_{R}^{[\bullet]}\{n\}$. By virtue of Remark \ref{remark:diffracted-vs-prismatic1}, the associated graded complex is given by $\gr^{m} \Prism_{R}^{[\bullet]}\{n\} \simeq \overline{\Prism}_{R}\{m+n\}$.
\end{example}

\begin{proposition}\label{proposition:sokon}
Let $n$ be an integer. Then:
\begin{itemize}
\item[$(1)$] For every animated commutative ring $R$, the object $\Prism_{R}^{[\bullet]}\{n\} \in \DFiltI(\Z_p)$ is filtration-complete.
\item[$(2)$] The construction $R \mapsto \Prism_{R}^{[\bullet]}\{n\}$ determines a functor of $\infty$-categories $\CAlg^{\anim} \rightarrow \DFiltIComp(\Z_p)$ which commutes with sifted colimits.
\end{itemize}
\end{proposition}

\begin{proof}
Let $R$ be an animated commutative ring. Then the limit of the tower
$$ \cdots \rightarrow \mathcal{I}^{\otimes 2} \otimes \mathscr{H}_{\Prism}(R) \rightarrow \mathcal{I} \otimes \mathscr{H}_{\Prism}(R) \rightarrow \mathscr{H}_{\Prism}(R)$$ 
is zero (as an object of the $\infty$-category $\calD( \WCart )$). Tensoring with $\calO_{ \WCart }\{n\}$ and taking global sections, we deduce that the limit of the tower
$$ \cdots \rightarrow \Prism_{R}^{[2]}\{n\} \rightarrow \Prism_{R}^{[1]}\{n\} \rightarrow \Prism_{R}\{n\}$$
is also zero (as an object of the $\infty$-category $\widehat{\calD}(\Z_p)$). This proves $(1)$. To prove $(2)$, we note that the functor
$\Fil^{\bullet}(M) \mapsto \gr^{\bullet}(M)$ determines a conservative and colimit-preserving functor from $\DFiltIComp(\Z_p)$ to the product $\prod_{m \in \Z} \widehat{\calD}(\Z_p)$. It will therefore suffice to show that, for every integer $m$, the functor
$$ R \mapsto \gr^{m} \Prism_{R}^{[\bullet]}\{n\} \simeq \overline{\Prism}_{R}\{m+n\}$$
commutes with sifted colimits, which follows from Proposition \ref{proposition:HT-sifted-colimit}.
\end{proof}


\subsection{Comparison with Crystalline Cohomology}\label{subsection:crystalline-comparison}

For every $\F_p$-scheme $X$, let $\RGamma_{\crys}(X, \F_p)$ denote the crystalline cochain complex of
$X$ (Construction \ref{construction:crystalline-cochain-complex}) and let $$\epsilon_{\crys}: \RGamma_{\crys}(X/\Z_p) \rightarrow \RGamma(X, \calO_X)$$
denote the crystalline augmentation map of Notation \ref{notation:crystalline-augmentation}. Our goal in this section is to prove the following comparison result:

\begin{theorem}\label{theorem:crystalline-comparison}
Let $X$ be a quasisyntomic $\F_p$-scheme. Then there is a canonical isomorphism
$$ \gamma_{\Prism}^{\crys}: \RGamma_{\Prism}(X) \simeq \RGamma_{\crys}(X / \Z_p)$$
of commutative algebra objects of $\widehat{\calD}(\Z_p)$, which is characterized (up to homotopy)
by the requirement that it depends functorially on $X$ and that the diagram
$$ \xymatrix@R=50pt@C=50pt{ \RGamma_{\crys}(X / \Z_p) \ar[r]^-{ ( \gamma_{\Prism}^{\crys} )^{-1} } \ar[d]^{\epsilon_{\crys}} & \RGamma_{\Prism}(X) \ar[r]^-{\varphi} & \RGamma_{\Prism}(X) \ar[d] \\
\RGamma(X, \calO_X) \ar[r]^-{\sim} & \Fil_{0}^{\conj} \RGamma_{\overline{\Prism}}(X/ \Z_p) \ar[r] & \RGamma_{ \overline{\Prism} }(X) }$$
commutes (up to a homotopy depending functorially on $X$).
\end{theorem}

\begin{warning}\label{warning:Frob-semilinearity}
Let $k$ be a perfect field of characteristic $p$ and let $X$ be a $p$-quasisyntomic $k$-scheme.
Then we can identify $\RGamma_{\Prism}(X)$ with the prismatic complex $\RGamma_{\Prism}(X / W(k) )$ of $X$ relative to the perfect prism
$( W(k), (p) )$ (Proposition \ref{proposition:absolute-vs-relative}). Similarly, the crystalline cochain complex $\RGamma_{\crys}(X / \Z_p )$
of $X$ relative to $\Z_p$ can be identified with its crystalline cochain complex $\RGamma_{\crys}( X / W(k) )$ relative to the ring $W(k)$.
Beware that the resulting identification of $\RGamma_{\Prism}(X / W(k) )$ with $\RGamma_{\crys}( X / W(k) )$ is {\em not} $W(k)$-linear.
Instead it is Frobenius semilinear: more precisely, Theorem \ref{theorem:crystalline-comparison} supplies an isomorphism 
$\RGamma_{\crys}( X / W(k) ) \simeq F^{\ast} \RGamma_{\Prism}( X / W(k) )$ in the $\infty$-category $\calD( W(k) )$
(see Example \ref{example:Frob-inverse-semilinearity} below).
\end{warning}

\begin{remark}
When $X$ is a smooth $\F_p$-scheme, the existence of an isomorphism
$\RGamma_{\Prism}(X) \simeq \RGamma_{\crys}( X/ \Z_p )$ is a special case of Theorem~5.2 of \cite{prisms}.
\end{remark}

\begin{notation}
Let $X$ be a quasisyntomic $\F_p$-scheme. We will refer to the isomorphism $\gamma_{\Prism}^{\crys}: \RGamma_{\Prism}(X) \simeq \RGamma_{\crys}(X / \Z_p)$
as the {\it crystalline comparison map}. 
\end{notation}

\begin{remark}\label{remark:absolute-as-derived-crystalline}
By virtue of Theorem \ref{theorem:crystalline-comparison}, the construction $R \mapsto \Prism_{R} \simeq \Prism_{R/\Z_p}$ determines a functor of $\infty$-categories
$$ \CAlg^{\anim}_{\F_p} \rightarrow \widehat{\calD}(\Z_p)$$
which commutes with sifted colimits and is isomorphic to $\RGamma_{\crys}( \bullet / \Z_p)$ when restricted to the category $\Poly_{\F_p}$ of finitely
generated polynomial algebras over $\F_p$. It follows that, when $R$ is an animated $\F_p$-algebra, the absolute prismatic complex $\Prism_{R}$
computes the {\em derived} crystalline cohomology of $R$: that is, $R \mapsto \Prism_{R}$ can be regarded as the nonabelian left derived functor
of the crystalline cochain complex functor $R \mapsto \RGamma_{\crys}( R / \Z_p)$, in the sense of Proposition \ref{proposition:universal-of-animated}.
\end{remark}

\begin{remark}\label{remark:crystalline-comparison-general}
The functor
$$ \CAlg_{\F_p}^{\anim} \rightarrow \widehat{\calD}(\Z_p) \quad \quad R \mapsto \RGamma_{\Prism}( \Spec(R) ) \simeq \Prism_{R/\Z_p}$$
commutes with sifted colimits, and is therefore a left Kan extension of its restriction to the category $\CAlg_{\F_p}^{\QSyn}$ of $p$-quasisyntomic
$\F_p$-algebras (or even the much smaller category of finitely generated polynomial algebras over $\F_p$). It follows that the crystalline comparison map
$\gamma_{\Prism}^{\crys}: \Prism_{R} \rightarrow \RGamma_{\crys}(R/\Z_p)$ admits an essentially unique functorial extension to the category of
{\em all} commutative $\F_p$-algebras $R$. By Zariski descent, we obtain a crystalline comparison map
$\gamma_{\Prism}^{\crys}: \RGamma_{\Prism}(X) \rightarrow \RGamma_{\crys}( X / \Z_p)$ for every $\F_p$-scheme $X$. Beware that this
map is generally not an isomorphism if the $\F_p$-scheme $X$ is not $p$-quasisyntomic.
\end{remark}

\begin{corollary}\label{corollary:qs-descent-for-crys}
The functor
$$ \{ \textnormal{Quasisyntomic $\F_p$-schemes} \}^{\op} \rightarrow \widehat{\calD}(\Z_p) \quad \quad X \mapsto \RGamma_{\crys}(X / \Z_p)$$
satisfies descent for the $p$-quasisyntomic topology.
\end{corollary}

\begin{proof}
Combine Theorem \ref{theorem:crystalline-comparison} with Proposition \ref{proposition:absolute-quasisyntomic-descent}.
\end{proof}



We begin our proof of Theorem \ref{theorem:crystalline-comparison} by considering the case where $X = \Spec(R)$ is the spectrum of quasiregular semiperfect $\F_p$-algebra $R$.
By virtue of Remark \ref{remark:senty}, we can identify the crystalline cochain complex $\RGamma_{\crys}( X / \Z_p )$ with the commutative ring $A_{\crys}(R)$ (see Construction \ref{construction:Acrys}), regarded as a cochain complex concentrated in cohomological degree zero. Similarly, the absolute prismatic complex $\RGamma_{\Prism}(X) \simeq \Prism_{R}$
is concentrated in cohomological degree zero by Example \ref{example:prismatic-cohomological-qrsp}. In this case, the inverse of the isomorphism
$\gamma_{\Prism}^{\crys}$ has a simple description.

\begin{lemma}\label{lemma:transformation-in-easy-case}
Let $R$ be a quasiregular semiperfect $\F_p$-algebra. Then there is a unique ring homomorphism $\beta_{R}: A_{\crys}(R) \rightarrow \Prism_{R}$
for which the diagram of ring homomorphisms
$$ \xymatrix@R=50pt@C=50pt{ A_{\crys}(R) \ar[r]^-{\beta_{R} } \ar[d]^{\epsilon_{\crys}} & \Prism_{R} \ar[r]^-{\varphi} & \Prism_{R} \ar[d] \\
R \ar[rr] & & \overline{\Prism}_{R} }$$
is commutative.
\end{lemma}

\begin{proof}
By virtue of Proposition~7.10 of \cite{prisms}, we can view the pair $( \Prism_{R}, \Prism_{R}^{[1]} )$ as a prism (Example \ref{example:prismatic-cohomological-qrsp}).
Since the quotient $\overline{\Prism}_{R} = \Prism_{R} / \Prism_{R}^{[1]}$ has the structure of an $R$-algebra, the ideal $\Prism_{R}^{[1]}$ contains the distinguished element
$p \in \Prism_{R}$ and is therefore generated by $p$ (in particular, $p$ is a regular element of $\Prism_{R}$). By Theorem~12.2 of \cite{prisms},
the composite map $$\Prism_{R} \xrightarrow{\varphi} \Prism_{R} \twoheadrightarrow \Prism_{R} / \Prism_{R}^{[1]} = \overline{\Prism}_{R}$$
has image $R \simeq \Fil_{0}^{\conj} \overline{\Prism}_{R/\Z_p}$, and therefore restricts to a surjection $\Prism_{R} \twoheadrightarrow R$
with kernel $J = \{ x \in \Prism_{R}: \varphi(x) \in p \Prism_{R} \}$. By virtue of Lemma~2.35 of \cite{prisms}, every element $x \in J$
has divided powers in $\Prism_{R}$. Moreover, for each $n > 0$, $\varphi( \frac{ x^n }{n!} ) = \frac{ \varphi(x)^{n} }{n!}$ is divisible by $\frac{ p^n }{n!}$ and therefore by $p$, so that
$\frac{ x^{n} }{n!}$ also belongs to $J$. The existence and uniqueness of $\beta_{R}$ now follow from the universal property of $A_{\crys}(R)$ (see Remark \ref{remark:senty}).
\end{proof}

\begin{example}\label{example:Frob-inverse-semilinearity}
Let $R$ be a perfect $\F_p$-algebra, so that both $A_{\crys}(R)$ and $\Prism_{R}$ can be identified with the ring of Witt vectors $W(R)$
(in the latter case, we normalize this identification to be a lift of the unit map $R \xrightarrow{\sim} \overline{\Prism}_{R}$). Under these identifications, the homomorphism $\beta_{R}$ of Lemma \ref{lemma:transformation-in-easy-case} corresponds to the automorphism $W( \varphi_{R}^{-1} ): W(R) \rightarrow W(R)$ (that is, the {\em inverse} of the Witt vector Frobenius).
\end{example}

It follows immediately from the proof of Lemma \ref{lemma:transformation-in-easy-case} that the construction $R \mapsto \beta_{R}$ is functorial: that is,
it can be viewed as a natural transformation of functors from the category $\CAlg_{\F_p}^{\qrsp}$ to the category of commutative rings. This immediately
yields the following weak version of Theorem \ref{theorem:crystalline-comparison}:

\begin{corollary}\label{corollary:transformation-Kan-extended}
Let $X$ be a quasisyntomic $\F_p$-scheme. There is a canonical morphism $\beta_{X}:
\RGamma_{\crys}( X / \Z_p) \rightarrow \RGamma_{\Prism}(X)$ in the $\infty$-category $\CAlg( \widehat{\calD}( \Z_p ) )$, which is characterized up to homotopy by the requirement that it depends functorially on $X$ and that the diagram
$$ \xymatrix@R=50pt@C=50pt{ \RGamma_{\crys}(X / \Z_p) \ar[r]^-{ \beta_{X} } \ar[d]^{\epsilon_{\crys}} & \RGamma_{\Prism}(X) \ar[r]^-{\varphi} & \RGamma_{\Prism}(X) \ar[d] \\
\RGamma(X, \calO_X) \ar[r]^-{\sim} & \Fil_{0}^{\conj} \RGamma_{\overline{\Prism}}(X/ \Z_p) \ar[r] & \RGamma_{ \overline{\Prism} }(X) }$$
commutes (up to a homotopy depending functorially on $X$).
\end{corollary}

\begin{proof}
Combine Lemma \ref{lemma:transformation-in-easy-case} with the observation that the functors $X \mapsto \RGamma_{\Prism}(X)$ and $X \mapsto \RGamma_{\overline{\Prism}}(X)$
satisfy descent for the $p$-quasisyntomic topology (Proposition \ref{proposition:absolute-quasisyntomic-descent}).
\end{proof}

We will prove Theorem \ref{theorem:crystalline-comparison} by showing that the map $\beta_{X}$ of Corollary \ref{corollary:transformation-Kan-extended} is an isomorphism
(we then take the crystalline comparison map $\beta_{\Prism}^{\crys}$ to be the inverse isomorphism $\beta_{X}^{-1}$). We begin by treating the case where
$X = \Spec(R)$ is the spectrum of a quasiregular semiperfect $\F_p$-algebra.

\begin{lemma}\label{lemma:nonzero-constant}
Let $R$ be a quasiregular semiperfect $\F_p$-algebra and let $I$ denote the kernel of the quotient map $R^{\flat} \twoheadrightarrow R$,
so that the relative cotangent complex $L \Omega^{1}_{ R/ \F_p}$ can be identified with $(I/I^2)[1]$. Let 
$$\overline{\beta}_{R}: A_{\crys}(R) / p A_{\crys}(R) \rightarrow \overline{\Prism}_{R} \simeq \overline{\Prism}_{R/\Z_p}$$
denote the mod $p$ reduction of the ring homomorphism $\beta_{R}$ of Lemma \ref{lemma:transformation-in-easy-case}.
Let $x$ be an element of $I$ having image $\overline{x}$ in $I/I^2$.
\begin{itemize}
\item[$(1)$] For each integer $n \geq 0$, the ring homomorphism
$\overline{\beta}_{R}$ carries the divided power $\gamma_{pn}(x) \in A_{\crys}(R) / p A_{\crys}(R)$ 
into the subgroup $\Fil_{n}^{\conj}( \overline{\Prism}_{R/\Z_p} ) \subseteq \overline{\Prism}_{R/\Z_p}$.

\item[$(2)$] There exists a nonzero constant $\lambda(n) \in \F_p$ (not depending on $R$ or $x$) for which the map
$$ \Fil_{n}^{\conj }( \overline{\Prism}_{R/\Z_p} ) \twoheadrightarrow \gr_{n}^{\conj}( \overline{\Prism}_{R/ \Z_p} )
\simeq \Gamma^{n}_{ R }( I/I^2)$$
carries $\overline{\beta}_{R}( \gamma_{pn}(x) )$ to $\lambda(n) \gamma_{n}( \overline{x} )$.
\end{itemize}
\end{lemma}

\begin{proof}
By functoriality, it suffices to treat the case where $R^{\flat}$ is the perfect polynomial ring $k[ x^{1/p^{\infty}} ]$ for some perfect field $k$
and $R$ is the quotient ring $R^{\flat} / (x)$ (in fact, we may assume that $k = \F_p$). Moreover, we are free to enlarge $k$
and may therefore assume without loss of generality that there exists an element $t \in k$ is transcendental over the ground field. 
Then $I/I^2$ is a free $R$-module of rank $1$ generated by $\overline{x}$. It follows from Example~7.9 of \cite{prisms}
that the morphism $\beta_{R}$ is an isomorphism in this case, so that $\overline{\beta}_{R}( \gamma_{pn}(x) )$ is a nonzero element of
the commutative ring $\overline{\Prism}_{R}$. Consequently, there exists a unique integer $m \geq 0$ for which
$\overline{\beta}_{R}( \gamma_{pn}(x) )$ belongs to the subgroup $\Fil_{m}^{\conj }( \overline{\Prism}_{R/\Z_p})$ and has nonzero image in
$\gr_{m}^{\conj}(  \overline{\Prism}_{R/\Z_p} ) \simeq \RGamma_{n}( \overline{x} )$; let us denote this image
by $f \gamma_m(x)$.

Let $\sigma$ denote the automorphism of $R$ given by the construction $x^{1/p^n} \mapsto (tx)^{1/p^n}$. Then
induces a $k$-linear automorphism of $A_{\crys}(R) / p A_{\crys}(R)$ for which the element $\gamma_{pn}(x)$ is an 
eigenvector with eigenvalue $t^{pn}$. Example \ref{example:Frob-inverse-semilinearity} implies that
$\overline{\beta}_{R}$ is $\varphi_{k}^{-1}$-semilinear, so that the $f \gamma_m(x)$ is also an eigenvector for the action of
$\sigma$ with eigenvalue $t^{n}$. It follows that the element $f \in R$ is a nonzero eigenvector for the action of
$\sigma$ on $R$ with eigenvalue $t^{n-m}$. We now observe that $R = k[ x^{1/p^{\infty}} ] / (x)$ has a $k$-basis
consisting of elements $\{ x^{\alpha} \}_{ \alpha \in \Z[1/p], 0 \leq \alpha < 1 }$, where each $x^{\alpha}$
is an eigenvector for $\sigma$ with eigenvalue $t^{\alpha}$. In particular, the only integral power of $t$
which appears as an eigenvalue for $\sigma$ on $R$ is the power $t^{0} = 1$, and the corresponding eigenspace
$R^{ \sigma = 1}$ is equal to the ground field $k$. This proves that $m = n$ and that $f \in k$ is a scalar.
By functoriality, this scalar does not depend on the field $k$ and is therefore a nonzero element $\lambda(n)$ of the prime field $\F_p$.
\end{proof}

\begin{remark}
In the situation of Lemma \ref{lemma:nonzero-constant}, it is not difficult to see that the scalars $\lambda(n) \in \F_p^{\times}$ must satisfy
$\lambda(n) = \lambda^n$ for a {\em single} nonzero constant $\lambda = \lambda(1) \in \F_p^{\times}$. The precise value of this constant is dependent on sign conventions, so
we will not attempt to evaluate it.
\end{remark}

For every semiperfect $\F_p$-algebra $R$, let us regard the quotient $A_{\crys}(R) / p A_{\crys}(R)$ as equipped with the conjugate
filtration $\Fil_{\bullet}^{\conj} A_{\crys}(R) / p A_{\crys}(R)$ of Construction \ref{construction:conjugate-filtration}.

\begin{lemma}\label{lemma:conjugate-filtration}
Let $R$ be a quasiregular semiperfect $\F_p$-algebra. Then:
\begin{itemize}
\item[$(1)$] For each $n \geq 0$, the homomorphism 
$$\overline{\beta}_{R}: A_{\crys}(R) / p A_{\crys}(R) \rightarrow \overline{\Prism}_{R} = \overline{\Prism}_{R/\Z_p}$$
carries $\Fil^{\conj}_{n}( A_{\crys}(R) / p A_{\crys}(R) )$ into $\Fil_{n}^{\conj}( \overline{\Prism}_{R/\Z_p} )$. 

\item[$(2)$] The induced map
$$ \gr( \overline{\beta}_{R} ): \gr_{\ast}^{\conj}( A_{\crys}(R) / p A_{\crys}(R) ) \rightarrow \gr_{\ast}^{\conj}( \overline{\Prism}_{R/\Z_p})$$
is an isomorphism of graded abelian groups.

\item[$(3)$] For each $n \geq 0$, the map $\overline{\beta}_{R}$ induces an isomorphism
$$\Fil_{n}^{\conj}( A_{\crys}(R) / p A_{\crys}(R) ) \simeq \Fil_{n}^{\conj}( \overline{\Prism}_{R / \Z_p} ).$$

\item[$(4)$] The map $\overline{\beta}_{R}: A_{\crys}(R) / p A_{\crys}(R) \rightarrow \overline{\Prism}_{R}$ is an isomorphism.

\item[$(5)$] The map $\beta_{R}: \RGamma_{\crys}(R/\Z_p) \rightarrow \Prism_{R}$ is an isomorphism.
\end{itemize}
\end{lemma}

\begin{proof}
Assertion $(1)$ follows immediately from Lemma \ref{lemma:nonzero-constant}, and the implications $(2) \Rightarrow (3) \Rightarrow (4) \Rightarrow (5)$
are formal. We will prove $(2)$. Let $I$ denote the kernel of the map $R^{\flat} \twoheadrightarrow R$, and let
$\xi: \Gamma^{\ast}_{R}(I / I^2) \twoheadrightarrow  \gr_{\ast}^{\conj}( A_{\crys}(R) / p A_{\crys}(R) )$ be the surjection of
Proposition \ref{proposition:describe-gr-conjugate}. It follows from Lemma \ref{lemma:nonzero-constant} that the composite map
$$ \Gamma^{\ast}_{R}(I / I^2) \xrightarrow{\xi}  \gr_{\ast}^{\conj}( A_{\crys}(R) / p A_{\crys}(R) ) \xrightarrow{ \gr( \overline{\beta}_{R} ) }
\gr_{\ast}^{\conj} ( \overline{\Prism}_{R/\Z_p} ) \simeq 
\Gamma^{\ast}_{R}( I/ I^2 )$$
is an isomorphism, which immediately implies that $\gr( \overline{\beta}_{R} )$ is also an isomorphism.
\end{proof}

\begin{remark}\label{remark:gr-conjugate-isomorphism}
Let $R$ be a semiperfect $\F_p$-algebra, let $I$ denote the kernel of the map $R^{\flat} \twoheadrightarrow R$,
and let $$\xi: \Gamma^{\ast}_{R}( I/ I^2 ) \twoheadrightarrow \gr_{\ast}^{\conj}( A_{\crys}(R) / p A_{\crys}(R) )$$ be the 
surjection of Proposition \ref{proposition:describe-gr-conjugate}. The proof of Lemma \ref{lemma:conjugate-filtration} shows that,
if $R$ is quasiregular semiperfect, then $\xi$ is an isomorphism.
\end{remark}

\begin{corollary}\label{corollary:A-crys-torsion}
Let $R$ be a quasiregular semiperfect $\F_p$-algebra. Then the commutative ring $A_{\crys}(R)$ is $p$-torsion-free.
\end{corollary}


\begin{proof}[Proof of Theorem \ref{theorem:crystalline-comparison}]
It will suffice to show that if $X$ is a quasisyntomic $\F_p$-scheme, then the map $\beta_{X}: \RGamma_{\crys}( X / \Z_p ) \rightarrow \RGamma_{\Prism}(X)$
of Corollary \ref{corollary:transformation-Kan-extended} is an isomorphism in the $\infty$-category $\widehat{\calD}(\Z_p)$. The assertion is Zariski local on $X$; we may therefore
assume without loss of generality that $X = \Spec(R)$ is the spectrum of a $p$-quasisyntomic commutative $\F_p$-algebra $R$.
Choose a surjection $P \twoheadrightarrow R$, where $P = \F_p[ \{ x_s \} ]$ is a polynomial ring over $\F_p$ (possibly on an infinite set of generators).
Let $P^{0} = \F_p[ \{ x_s^{1/p^{\infty}} \} ]$ be the perfection of $P^{0}$, let $P^{\bullet}$ be the cosimplicial $P$-algebra given by the iterated
tensor powers of $P^{0}$ relative to $P$, and set $R^{\bullet} = P^{\bullet} \otimes_{P} R$. The morphism $\beta_{X}$ fits into a commutative diagram
$$ \xymatrix@R=50pt@C=50pt{ \RGamma_{\crys}( X / \Z_p ) \ar[r]^-{ \beta_{X} } \ar[d] & \RGamma_{\Prism}(X) \ar[d] \\
\Tot( \RGamma_{\crys}( R^{\bullet} / \Z_p) ) \ar[r]^-{ \Tot( \beta_{R^{\bullet}} )} &  \Tot( \Prism_{R^{\bullet}} ), }$$
where the left vertical map is an isomorphism by virtue of Proposition \ref{proposition:descend-to-semiperfect} and the right vertical map
is an isomorphism by virtue of Proposition \ref{proposition:absolute-quasisyntomic-descent}. Consequently, to show that
$\beta_{X}$ is an isomorphism, it will suffice to show that $\beta_{R^{m}}: \RGamma_{\crys}( R^{m} / \Z_p) \rightarrow \Prism_{ R^{m} }$ is an isomorphism for each $m \geq 0$. This follows
from Lemma \ref{lemma:conjugate-filtration}, since $R^{m}$ is quasiregular semiperfect.
\end{proof}

\subsection{Diffracted Hodge Cohomology of Formal Schemes}\label{subsection:diffracted-Hodge-completed}

We now introduce a variant of the absolute Hodge-Tate complex $\overline{\Prism}_{R}$ which retains more information about the Hodge-Tate cohomology sheaf $\mathscr{H}_{\overline{\Prism}}(R)$ of Construction \ref{construction:absolute-Hodge-Tate-filtration}.

\begin{construction}[The $p$-Complete Diffracted Hodge Complex]\label{construction:complete-diffracted-Hodge}
Let $R$ be an animated commutative ring, and let $\mathscr{H}_{\Prism}(R) \in \calD( \WCart )$ be the prismatic cohomology sheaf of Construction \ref{construction:prismatic-cohomology-sheaves}. We let $\widehat{\Omega}^{\DHod}_{R} \in \widehat{\calD}(\Z_p)$ denote the pullback of $\mathscr{H}_{\overline{\Prism}}(R)$ along the map $\eta: \Spf(\Z_p) \rightarrow \WCart^{\mathrm{HT}} \subseteq \WCart$ described in Construction \ref{construction:fiber-at-eta}. We will refer to $\widehat{\Omega}^{\DHod}_{R}$ as the {\it $p$-complete diffracted Hodge complex of $R$}.

For every integer $n$, we let $\Fil_{n}^{\conj} \widehat{\Omega}^{\DHod}_{R}$ denote the pullback $\eta^{\ast} \Fil_{n}^{\conj} \mathscr{H}_{\overline{\Prism}}(R)$. Allowing $n$ to vary, we obtain a direct system
$$ \Fil_0^{\conj} \widehat{\Omega}^{\DHod}_{R} \rightarrow \Fil_1^{\conj} \widehat{\Omega}^{\DHod}_{R} \rightarrow
\Fil_2^{\conj} \widehat{\Omega}^{\DHod}_{R} \rightarrow \cdots$$
with having colimit $\widehat{\Omega}^{\DHod}_{R}$ in the $p$-complete derived $\infty$-category $\widehat{\calD}(\Z_p)$. We will denote this direct system by
$\Fil_{\bullet}^{\conj} \widehat{\Omega}^{\DHod}_{R}$ and refer to it as the {\it conjugate filtration} on $\widehat{\Omega}^{\DHod}_{R}$.

By virtue of Remark \ref{remark:sheafy-HT-comparison}, we have canonical isomorphisms $\gr_{n}^{\conj} \widehat{\Omega}^{\DHod}_{R} \simeq L \widehat{\Omega}^{n}_{R}$
for each integer $n$. In particular, we can identify $\Fil_{0}^{\conj} \widehat{\Omega}^{\DHod}_{R}$ with the $p$-completion of $R$, so that
$\Fil_{\bullet}^{\conj} \widehat{\Omega}^{\DHod}_{R}$ can be promoted to a filtered object of the $p$-complete derived $\infty$-category $\widehat{\calD}(R)$.
\end{construction}

\begin{notation}[The Sen Operator]\label{notation:Sen-operator-on-diffracted}
Let $R$ be an animated commutative ring. For every integer $n$, Notation \ref{notation:Sen-operator-second} determines
an endomorphism of the complex $\Fil_{n}^{\conj} \widehat{\Omega}^{\DHod}_{R} \in \calD(R)$, which we will denote by
$\Theta$ and refer to as the {\it Sen operator}. This construction depends functorially on $n$, and therefore also
determines an endomorphism $\Theta$ of the $p$-completed colimit $\widehat{\Omega}^{\DHod}_{R} = \varinjlim_{n} \Fil_{n}^{\conj} \widehat{\Omega}^{\DHod}_{R}$.

Using the isomorphism $$ \gr_{n}^{\conj} \mathcal{H}_{\overline{\Prism}}(R) \simeq L \widehat{\Omega}^{n}_{R} \otimes
\calO_{ \WCart^{\mathrm{HT}} }[-n]\{-n\}$$
supplied by  Remark \ref{remark:sheafy-HT-comparison}, we see that the Sen endomorphism of $\gr_{n}^{\conj} \widehat{\Omega}^{\DHod}_{R}$
is given by multiplication by $-n$.
\end{notation}

\begin{remark}\label{remark:connectivity-of-diffracted-Hodge}
Let $R$ be an animated commutative ring. For every integer $n$, the cohomology groups of the complex
$\Fil_{n}^{\conj} \widehat{\Omega}^{\DHod}_{R}$ are concentrated in degrees $\leq n$ (this follows by induction on $n$,
using the isomorphism $\gr_{n}^{\conj} \widehat{\Omega}^{\DHod}_{R} \simeq L \widehat{\Omega}^{n}_{R}[-n]$).
\end{remark}

\begin{remark}\label{remark:conjugate-equals-Postnikov1}
Let $R$ be a commutative ring which is $p$-torsion-free and suppose that the quotient ring $R/pR$ is a regular Noetherian $\F_p$-algebra.
Then, for every integer $n$, the complex $\gr_{n}^{\conj} \widehat{\Omega}^{\DHod}_{R} \simeq L \widehat{\Omega}^{n}_{R}[-n]$ is concentrated in cohomological degree $n$. It follows that the conjugate filtration on $\widehat{\Omega}^{\DHod}_{R}$ coincides with the Postnikov filtration:
that is, we have $\Fil_{n}^{\conj} \widehat{\Omega}^{\DHod}_{R} \simeq \tau^{\leq n} \widehat{\Omega}^{\DHod}_{R}$ for every integer $n$.
\end{remark}

\begin{remark}\label{remark:diffracted-vs-prismatic2}
Let $R$ be an animated commutative ring. For every integer $n$, we let $(\widehat{\Omega}^{\DHod}_{R})^{\Theta = n}$ denote the fiber of the map
$$ \Theta - n: \widehat{\Omega}^{\DHod}_{R} \rightarrow \widehat{\Omega}^{\DHod}_{R}.$$
Applying Proposition \ref{proposition:easy-version} to the complex $\mathscr{H}_{\overline{\Prism}}(R)\{-n\} \in \calD( \WCart^{\mathrm{HT} } )$, we obtain
a canonical isomorphism $\overline{\Prism}_{R}\{-n\} \simeq (\widehat{\Omega}^{\DHod}_{R})^{\Theta = n}$ in the $p$-complete derived $\infty$-category $\widehat{\calD}(R)$.
More generally, we have isomorphisms $\Fil_{m}^{\conj} \overline{\Prism}_{R}\{-n\} \simeq  (\Fil_{m}^{\conj} \widehat{\Omega}^{\DHod}_{R})^{\Theta = n}$
for every integer $m$.
\end{remark}

\begin{example}\label{example:diffracted-Hodge-of-perfectoid}
Let $R$ be a perfectoid ring. Then the $p$-complete diffracted Hodge complex $\widehat{\Omega}^{\DHod}_{R}$
can be identified with the coordinate ring of the $\mathbf{G}_{m}^{\sharp}$-torsor $\mathcal{P} \rightarrow \Spf(R)$
of Example \ref{example:torsor-over-Hodge-Tate} (this is an immediate consequence of Example \ref{example:prismatic-sheaf-qrsp}).
In particular, $\widehat{\Omega}^{\DHod}_{R}$ is concentrated in cohomological degree zero.
\end{example}

\begin{remark}[Relative Hodge-Tate Comparison]\label{remark:diffracted-relative-comparison-affine}
Let $(A,I)$ be a bounded prism with quotient ring $\overline{A} = A/I$, and suppose we are given a trivialization of 
the torsor $\mathcal{P}$ of Example \ref{example:torsor-over-Hodge-Tate}, which identifies the composite map
$\Spf( \overline{A} ) \rightarrow \Spf(\Z_p) \xrightarrow{\eta} \WCart$ with the map $\rho^{\mathrm{HT}}_{A}: \Spf( \overline{A} ) \rightarrow \WCart$ of Remark \ref{remark:HT-point-of-prismatic-stack}.
For every animated commutative ring $R$, we obtain an isomorphism
$\overline{\Prism}_{ (R \otimes^{L} \overline{A} )/A} \simeq  \overline{A} \widehat{\otimes}^{L} \widehat{\Omega}^{\DHod}_{R}$
in the derived $\infty$-category $\widehat{\calD}( \overline{A} )$. Beware that this isomorphism is not completely canonical: it depends on a choice of trivialization of the torsor $\mathcal{P}$.
\end{remark}

\begin{example}\label{example:diffracted-Hodge-as-relative-prismatic}
Let $( \Z_p[[ \slashp ]], (\slashp) )$ be the prism of Notation \ref{notation:reduced-q-de-Rham-prism}, and let $R$ be an animated commutative $\Z_p$-algebra.
Combining Remark \ref{remark:diffracted-relative-comparison-affine} with Proposition \ref{proposition:clean-statement-of-fiber}, we obtain an isomorphism
of the $p$-complete diffracted Hodge complex $\widehat{\Omega}^{\DHod}_{R}$ with the relative Hodge-Tate complex $\overline{\Prism}_{R/\Z_p[[\slashp]]}$. More precisely, we have an isomorphism of filtered complexes $$ \Fil_{\bullet}^{\conj} \widehat{\Omega}^{\DHod}_{R} \simeq \Fil_{\bullet}^{\conj} \overline{\Prism}_{R/\Z_p[[\slashp]]},$$
which is equivariant with respect to the action of the profinite group $\mathbf{G}_{m}^{\sharp}(\Z_p) \simeq (1+p\Z_p)^{\times}$.
\end{example}

\begin{remark}[Flat Descent]\label{remark:fpqc-descent-diffracted-hodge}
For every integer $n$, the functor $R \mapsto \Fil_{n}^{\conj} \widehat{\Omega}^{\DHod}_{R}$ satisfies
$p$-complete faithfully flat descent. By virtue of Example \ref{example:diffracted-Hodge-as-relative-prismatic}, this is
a special case of Variant \ref{variant:conjugate-fpqc-descent}.
\end{remark}

\begin{warning}
Let $(A,I)$ be a perfect prism and let $R$ be animated commutative algebra over the quotient ring $\overline{A} = A/I$,
so that Example \ref{example:relative-HT-over-perfect} supplies an isomorphism
$\overline{\Prism}_{R} \xrightarrow{\sim} \overline{\Prism}_{R/A}$. It is not difficult to see that this isomorphism can be refined to a morphism of filtered objects 
$$ \alpha: \Fil_{\bullet}^{\conj} \overline{\Prism}_{R} \rightarrow \Fil_{\bullet}^{\conj} \overline{\Prism}_{R/A},$$
where the left hand side is defined using the ``absolute'' conjugate filtration of Construction \ref{construction:absolute-HT} and
and the right hand side is defined using the relative conjugate filtration of Remark \ref{remark:derived-HT-filtration}.
Beware that $\alpha$ is {\em not} an isomorphism. For example, if $R$ is $p$-complete, then the complex $\Fil_{0}^{\conj} \overline{\Prism}_{R/A}$ can be identified with $R$,
while the complex $$\Fil_{0}^{\conj} \overline{\Prism}_{R} = ( \Fil_{0}^{\conj} \widehat{\Omega}^{\DHod}_{R} )^{\Theta = 0} \simeq
R^{\Theta=0}$$ can be identified with the direct sum $R \oplus R[-1]$.
\end{warning}

For every animated commutative ring $R$, the identification $\Fil_{0}^{\conj}  \widehat{\Omega}^{\DHod}_{R} \simeq \widehat{R}$ supplies each
$\Fil_{n}^{\conj} \widehat{\Omega}^{\DHod}_{R}$ with the structure of an $R$-module. These $R$-modules are compatible with \'{e}tale base change:

\begin{proposition}\label{proposition:etale-base-change-diffracted-complete}
Let $f: R \rightarrow S$ be a morphism of animated commutative rings which is $p$-completely formally \'{e}tale: that is, for which the $p$-complete relative cotangent complex
$L \widehat{\Omega}^{1}_{S/R}$ vanishes. Then, for every integer $n$, the natural map $S \widehat{\otimes}^{L}_{R} \Fil_{n}^{\conj} \widehat{\Omega}^{\DHod}_{R} \rightarrow \widehat{\Omega}^{\DHod}_{S}$
is an isomorphism in the $p$-complete derived $\infty$-category $\widehat{\calD}(S)$. In particular, the natural map $S \widehat{\otimes}^{L}_{R} \widehat{\Omega}^{\DHod}_{R} \rightarrow \widehat{\Omega}^{\DHod}_{S}$ is an isomorphism in $\widehat{\calD}(S)$.
\end{proposition}

\begin{proof}
Replacing $R$ and $S$ by their $p$-completions, we may assume without loss of generality that $R$ is a $\Z_p$-algebra. In this case, the result follows by
applying Proposition \ref{proposition:formally-etale-relative-Hodge-Tate} to the prism $(\Z_p[[\slashp]], (\slashp) )$ of Notation \ref{notation:reduced-q-de-Rham-prism}
(see Example \ref{example:diffracted-Hodge-as-relative-prismatic}).
\end{proof}

Using Proposition \ref{proposition:etale-base-change-diffracted-complete}, we can globalize Construction \ref{construction:complete-diffracted-Hodge} to the setting of $p$-adic formal schemes:

\begin{notation}\label{notation:diffracted-Hodge-formal-scheme}
Let $\mathfrak{X}$ be a bounded $p$-adic formal scheme, and let $\calU$ denote the collection of all affine open subsets of $X$. For every integer $n$, 
Proposition \ref{proposition:etale-base-change-diffracted-complete} guarantees that the construction
$$ (U \in \calU) \mapsto \Fil_{n}^{\conj} \widehat{\Omega}^{\DHod}_{ \calO_{\mathfrak{X}}(U) }$$
determines an object of the derived $\infty$-category $\calD( \mathfrak{X} )$, which we will denote by $\Fil_{n}^{\conj} \Omega^{\DHod}_{ \mathfrak{X} }$.
Allowing $n$ to vary, we obtain a direct system
$$\calO_{\mathfrak{X}} \simeq \Fil_{0}^{\conj} \Omega^{\DHod}_{\mathfrak{X} } \rightarrow \Fil_{1}^{\conj} \Omega^{\DHod}_{\mathfrak{X}} \rightarrow \Fil_{2}^{\conj} \Omega^{\DHod}_{\mathfrak{X}} \rightarrow \cdots$$
whose colimit we will denote by $\Omega^{\DHod}_{\mathfrak{X}}$ and refer to as the {\it diffracted Hodge complex} of $\mathfrak{X}$. We will refer to the (hypercohomology) groups
$\mathrm{H}^{\ast}( \mathfrak{X}, \Omega^{\DHod}_{\mathfrak{X}} )$ as the {\em diffracted Hodge cohomology groups} of $\mathfrak{X}$.
\end{notation}

\begin{example}
Let $R$ be commutative ring having bounded $p$-power torsion, let $\widehat{R}$ denote the $p$-completion of $R$, and let $\mathfrak{X} = \Spf( \widehat{R} )$ be the associated $p$-adic formal scheme. Then the diffracted Hodge complex $\Omega^{\DHod}_{\mathfrak{X}}$ can be identified with the image of $\widehat{\Omega}^{\DHod}_{R}$ under the equivalence of $\infty$-categories
$\widehat{\calD}( R ) \simeq \calD( \Spf( \widehat{R} ) )$. In particular, we have canonical isomorphisms
$\mathrm{H}^{\ast}( \mathfrak{X}, \Omega^{\DHod}_{\mathfrak{X} } ) \simeq \mathrm{H}^{\ast}( \widehat{\Omega}^{\DHod}_{R} )$.
\end{example}

\begin{remark}\label{remark:conjugate-filtration-in-general}
Let $\mathfrak{X}$ be a bounded $p$-adic formal scheme. Globalizing the Hodge-Tate comparison of Construction \ref{construction:complete-diffracted-Hodge}, we obtain canonical isomorphisms
$$ \gr_{n}^{\conj} \Omega^{\DHod}_{\mathfrak{X}} \simeq L \widehat{\Omega}^{n}_{\mathfrak{X}}[-n]$$
in the derived $\infty$-category $\calD( \mathfrak{X}$); here $L \widehat{\Omega}^{n}_{\mathfrak{X}}$ denotes the $n$th derived exterior power of the absolute cotangent complex of $\mathfrak{X}$.
\end{remark}

\begin{example}\label{example:diffracted-Hodge-of-smooth}
Let $k$ be a perfect field of characteristic $p$ and let $\mathfrak{X}$ be a $p$-adic formal scheme which is smooth over $\Spf( W(k) )$. Then Remark \ref{remark:conjugate-filtration-in-general}
supplies isomorphisms $$\gr_{n}^{\conj} \widehat{\Omega}^{\DHod}_{\mathfrak{X}}
\simeq \widehat{\Omega}^{n}_{\mathfrak{X}/W(k)}[-n],$$ where $\widehat{\Omega}^{n}_{\mathfrak{X}/W(k)}$
denotes the locally free $\calO_{\mathfrak{X}}$-module of differential forms of degree $n$ on $\mathfrak{X}$ (relative to $W(k)$).
It follows that the conjugate filtration on $\widehat{\Omega}^{\DHod}_{\mathfrak{X}}$ coincides with its Postnikov filtration, and its cohomology sheaves are given by $\mathcal{H}^{n}( \widehat{\Omega}^{\DHod}_{\mathfrak{X}} ) \simeq \widehat{\Omega}^{n}_{ \mathfrak{X} / W(k)}$; in particular, they are coherent $\calO_{ \mathfrak{X} }$-modules. 
We therefore obtain a convergent spectral sequence
\begin{equation}\label{equation:spectral-sequence-diffracted}
\mathrm{H}^{a}(\mathfrak{X}, \widehat{\Omega}^{b}_{ \mathfrak{X} / W(k) } ) \Rightarrow \mathrm{H}^{a+b}(\mathfrak{X}, \Omega^{\DHod}_{\mathfrak{X}} ) \end{equation}
relating the diffracted Hodge cohomology of Notation \ref{notation:diffracted-Hodge-formal-scheme} to the classical Hodge cohomology groups of $\mathfrak{X}$ relative to $W(k)$.
In particular, if $\mathfrak{X}$ is proper and smooth over $\Spf( W(k))$, then the diffracted Hodge cohomology groups $\mathrm{H}^{\ast}(\mathfrak{X}, \Omega^{\DHod}_{\mathfrak{X} })$ are finitely generated $W(k)$-modules.
\end{example}

\begin{remark}[\'{E}tale Comparison]\label{remark:diffracted-etale-comparison}
Let $k$ be a perfect field of characteristic $p$. Choose an embedding $W(k) \rightarrow \calO_{C}$, where
$C$ is an algebraically closed perfectoid field, and set $A = W( \calO_{C}^{\flat} )$.
Let $X$ be a scheme which is smooth and proper over $W(k)$, let $\mathfrak{X} = \Spf( W(k) ) \times_{ \Spec(W(k))} X$ denote its formal completion along the vanishing locus of $p$, and let $X_{C} = \Spec(C) \times_{ \Spec(W(k))} X$ be its geometric generic fiber. Combining Theorem~14.3 of \cite{BMS1} with Theorem~17.2 of \cite{prisms},
we obtain a canonical isomorphism
$$ C \otimes^{L}_{\Z_p} \RGamma_{\mathet}( X_{C}, \Z_p ) \simeq
\overline{\Prism}_{ ( X \times_{\Spec(W(k))} \Spf(\calO_C) ) / A }[1/p].$$
By virtue of Remark \ref{remark:diffracted-relative-comparison-affine}, any trivialization of the torsor $\mathcal{P}$ of
Example \ref{example:torsor-over-Hodge-Tate} determines an isomorphism
$$ \overline{\Prism}_{ ( X \times_{\Spec(W(k))} \Spf(\calO_C) ) / A } \simeq \calO_{C} \otimes^{L}_{W(k)}
\RGamma( \mathfrak{X}, \widehat{\Omega}^{\DHod}_{\mathfrak{X}} )$$
(note that there is no need to $p$-complete the tensor product on the right hand side,
since $\RGamma( \mathfrak{X}, \widehat{\Omega}^{\DHod}_{\mathfrak{X}} )$ is already a perfect complex of $W(k)$-modules).
Combining these isomorphisms, we obtain an isomorphism
$$ C \otimes^{L}_{\Z_p}  \RGamma_{\mathet}( X_{C}, \Z_p ) \simeq
C \otimes^{L}_{ W(k) } \RGamma( \mathfrak{X},  \widehat{\Omega}^{\DHod}_{\mathfrak{X}} )$$
in the derived $\infty$-category $\calD(C)$. Beware that this isomorphism depends on the choice
of trivialization of $\mathcal{P}$ (though it has a canonical trivialization if we fix a system of primitive $p^n$th roots of unity in
$C$: see Proposition \ref{proposition:commutative-diagram-of-stacks}).
\end{remark}

Using diffracted Hodge cohomology and the Sen operator, we can give a new perspective on (and a refinement of) the Deligne-Illusie decomposition in \cite{MR894379}; this connection was first observed by Drinfeld (in somewhat different language).

\begin{example}[Formality of the diffracted Hodge complex in small dimensions]
\label{example:dimension-below-p}
Let $\mathfrak{X}$ be a bounded $p$-adic formal scheme. The Sen operator of Notation \ref{notation:Sen-operator-on-diffracted} determines an endomorphism
$\Theta$ of the diffracted Hodge complex $\Omega^{\DHod}_{\mathfrak{X}}$, which preserves the conjugate filtration $\Fil_{\bullet}^{\conj} \Omega^{\DHod}_{\mathfrak{X}}$
and acts by multiplication by $-n$ on each $\gr_{n}^{\conj} \Omega^{\DHod}_{\mathfrak{X}}$. Suppose that $\mathfrak{X}$ is smooth of dimension $< p$ over $\Spf( W(k) )$, for some perfect field $k$.
In this case, the associated graded complexes $\gr_{n}^{\conj} \Omega^{\DHod}_{\mathfrak{X}}$ vanish for $n \geq p$. It follows that the conjugate filtration of $\Omega^{\DHod}_{\mathfrak{X}}$
admits a canonical splitting, given by the generalized eigenspaces of $\Theta$. We therefore obtain isomorphisms 
\begin{equation}
\label{eq:DeligneIllusie}
 \Omega^{\DHod}_{\mathfrak{X}} \simeq \bigoplus_{b \geq 0} \widehat{\Omega}^{b}_{\mathfrak{X}/W(k)}[-b] \quad \quad \mathrm{H}^{n}(\mathfrak{X}, \Omega^{\DHod}_{\mathfrak{X} }) \simeq \bigoplus_{a+b=n} \mathrm{H}^{a}(\mathfrak{X}, \widehat{\Omega}^{b}_{ \mathfrak{X} / W(k) } ),
 \end{equation}
and similarly after reduction modulo $p$. 
\end{example}

\begin{remark}[Connection to Deligne-Illusie]
\label{remark:Deligne-Illusie}
Keep notation as in Example~\ref{example:dimension-below-p}, so $\mathfrak{X}/W(k)$ is a smooth $p$-adic formal scheme of dimension $<p$. Let us explain why the decomposition \eqref{eq:DeligneIllusie} can be regarded as an integral refinement of the Deligne-Illusie decomposition \cite{MR894379}. Recall that the prismatic-crystalline comparison theorem gives a natural identification
\[ \Omega^{\DHod}_{\mathfrak{X}^{(1)}} \otimes_{W(k)} k \simeq F_{\mathfrak{X}_k/k,*} \Omega^\bullet_{\mathfrak{X}_k/k}\]
in the quasi-coherent derived category $\mathcal{D}(\mathfrak{X}_k^{(1)})$, where $\mathfrak{X}^{(1)} = \varphi^* \mathfrak{X}$ is the twist of $\mathfrak{X}/W$ by the Frobenius $\varphi$ on $W(k)$, and  $F_{\mathfrak{X}_k/k}: \mathfrak{X}_k \to \mathfrak{X}_k^{(1)}$ is the relative Frobenius for $\mathfrak{X}_k/k$. Combining this isomorphism with the mod $p$ reduction of the first decomposition in \eqref{eq:DeligneIllusie} applied to $\mathfrak{X}^{(1)}$ gives a decomposition
\[ F_{\mathfrak{X}_k/k,*} \Omega^\bullet_{\mathfrak{X}_k/k} \simeq  \bigoplus_{b \geq 0} \Omega^b_{\mathfrak{X}_k^{(1)}/k}[-b].\]
A similar decomposition was proven in \cite{MR894379}; in fact one can show that the two decompositions coincide, though we do not show that here  (see \cite{LiMondalEnd} for this compatibility as well as an alternate perspective on the above decomposition and the nilpotent operators of Remark~\ref{rmk:Sen-nilpotent}). Note that the result of \cite{MR894379} has weaker hypotheses than the ones above: \cite{MR894379} only requires that the smooth $k$-scheme $\mathfrak{X}_k$ is endowed with a flat lift to $W_2(k)$, while our argument ostensibly uses the $W(k)$-lift $\mathfrak{X}$ of $\mathfrak{X}_k$; this can be avoided by contemplating the prismatization of $\mathrm{Spec}(\mathbf{Z}/p^2)$, see \cite[Remark 5.16]{BhattLurieAPCsequel}.
\end{remark}

\begin{remark}[Splitting small truncations of the de Rham complex]
\label{remark:Achinger}
The analysis in Example~\ref{example:dimension-below-p} and Remark~\ref{remark:Deligne-Illusie} applies more generally to show the following statement: if $k$ is a perfect field of characteristic $p$,  $\mathfrak{X}/W(k)$ is a smooth $p$-adic formal scheme (of any dimension), and $a \in \mathbf{Z}$ is an integer, then one has a natural decomposition
\begin{equation}
\label{eq:AchingerDecomp}
 \tau^{[a,a+p-1]} \left( F_{\mathfrak{X}_k/k,*} \Omega^\bullet_{\mathfrak{X}_k/k} \right) \simeq  \bigoplus_{b = a}^{a+p-1} \Omega^b_{\mathfrak{X}_k^{(1)}/k}[-b].
 \end{equation}
A related statement was previously shown independently by Achinger \cite{AchingerDI}. More precisely,  \cite{AchingerDI} applies to the slightly smaller truncation $\tau^{[a,a+p-2]} \left( F_{\mathfrak{X}_k/k,*} \Omega^\bullet_{\mathfrak{X}_k/k} \right)$, but only assumes $W_2(k)$-liftability. 
\end{remark}

\begin{remark}[Nilpotent operators attached to the Sen operator]
\label{rmk:Sen-nilpotent}
The arguments in Example~\ref{example:dimension-below-p} and Remarks~\ref{remark:Deligne-Illusie} and \ref{remark:Achinger} only use the generalized eigenspace decomposition provided by the Sen operator $\Theta$ to recover known decompositions of the de Rham complex in characteristic $p$. However, the Sen operator itself provides additional structure in the form of certain nilpotent operators on pieces of the de Rham complex; this structure seems previously unexplored. To explain the construction of these operators, fix a smooth $p$-adic formal  scheme $\mathfrak{X}/W(k)$. The generalized eigenspace decomposition for the $\Theta$ action on $\Omega^{\DHod}_{\mathfrak{X}}$ gives a decomposition
\[ \Omega^{\DHod}_{\mathfrak{X}} \simeq \bigoplus_{i=0}^{p-1} \Omega^{\DHod}_{\mathfrak{X},i},\]
with the $i$-th summand on the right corresponding to generalized eigenvalue $-i$, i.e., the summand whose cohomology groups are annihilated by $(\Theta + i)^N$ for $N \gg 0$. Each of the summands $\Omega^{\DHod}_{\mathfrak{X},i}$ thus comes equipped with a residual nilpotent operator $\Theta_i := \Theta+i$.  It would be interesting to understand $\Theta_i$ 
more explicitly; for instance, is its order of nilpotence bounded independently of $\dim(\mathfrak{X}_k)$? Perhaps the first question is whether the operators $\Theta_i$ can be nonzero; this has been answered affirmatively by Petrov (in preparation), who has constructed a smooth projective variety where $\Theta_0$ already acts non-trivially on the truncated mod $p$ reduction $\tau^{\leq p} (\Omega^{\DHod}_{\mathfrak{X},0}) \otimes^{L}_{W(k)} k$.
\end{remark}

\begin{remark}
In the situation of Remark \ref{rmk:Sen-nilpotent}, one can regard the Sen operator $\Theta$ as encoding an action of the group scheme
$\mathbf{G}_m^\sharp$-action on the diffracted Hodge complex $\Omega^{\DHod}_{\mathfrak{X}}$ (see Remark~\ref{rmk:LogGmsharp}). The decomposition 
$\Omega^{\DHod}_{\mathfrak{X}} \simeq \bigoplus_{i=0}^{p-1} \Omega^{\DHod}_{\mathfrak{X},i}$ then corresponds to the $(\Z/p\Z)$-grading
obtained by restricting to the subgroup $\mu_p \subset \mathbf{G}_m^\sharp$, and the nilpotent operator $\Theta_i$ encodes the residual action of the quotient
group $\mathbf{G}_m^\sharp/\mu_p \simeq \mathbf{G}_a^\sharp$ on each summand $\Omega^{\DHod}_{\mathfrak{X},i}$.
\end{remark}

\begin{remark}[Splitting the conjugate filtration rationally]
Let $\mathfrak{X}$ be a bounded $p$-adic formal scheme. In general, the conjugate filtration $\Fil_{\bullet}^{\conj} \Omega^{\DHod}_{\mathfrak{X}}$ need not be split.
However, it has a canonical splitting after inverting $p$ (into eigenspaces for the Sen operator $\Theta$). In particular, the spectral sequence 
(\ref{equation:spectral-sequence-diffracted}) is always rationally degenerate for $\mathfrak{X}/W(k)$ smooth. Via the comparisons in Remark~\ref{rmk:FCrysCrystalline} and Theorem~\ref{theorem:Sen-theory-refined}, this yields the Hodge-Tate decomposition on the $p$-adic \'etale cohomology groups $\mathrm{H}^*(\mathfrak{X}_C; C)$ when $\mathfrak{X}$ is proper and smooth over $W(k)$. 
\end{remark}

\begin{remark}[The Sen operator and the Hodge filtration]
Let $k$ be a perfect field of characteristic $p$, and let $\mathfrak{X}/W(k)$ is a smooth $p$-adic formal scheme. The identification $\Omega_{\mathfrak{X}}^{\DHod} \otimes_{W(k)} k \simeq F_{\mathfrak{X}_k/k,*} \Omega^\bullet_{\mathfrak{X}_k/k}$ (as in Remark~\ref{remark:Deligne-Illusie}) endows $\Omega_{\mathfrak{X}^{(1)}}^{\DHod} \otimes_{W(k)} k$ with two filtrations in $\calD(\mathfrak{X}_k^{(1)})$: an increasing conjugate filtration $\Fil_*^{\conj}$, and a decreasing Hodge filtration $\Fil^*_{\Hodge}$. By construction, the Sen operator $\Theta$ preserves the conjugate filtration. However, we do not expect $\Theta$ to preserve the Hodge filtration; instead, using Drinfeld's refined prismatization $\Sigma'$ from \cite{drinfeld-prismatic},  we expect (in forthcoming work) to show a version of Griffiths transversality in this context: $\Theta$ refines to a filtered map $\Fil^\bullet_{\Hodge} (\Omega_{\mathfrak{X}^{(1)}}^{\DHod} \otimes_{W(k)} k) \to \Fil^{\bullet-p}_{\Hodge}(\Omega_{\mathfrak{X}^{(1)}}^{\DHod} \otimes_{W(k)} k)$.
\end{remark}


\subsection{Comparison with \texorpdfstring{$q$}{q}-de Rham Cohomology}\label{subsection:concrete-qdR}

If $p$ is an odd prime, then we can use Theorem \ref{theorem:qdR-WCart} to connect absolute prismatic cohomology with the $q$-de Rham cohomology of \cite{scholzeq}.
We begin by adopting a definition of the latter which will be convenient for our discussion.

\begin{definition}\label{definition:q-de-Rham-cohomology}
Let $(\Z_p[[q-1]], ( [p]_q ) )$ be the $q$-de Rham prism of Example \ref{example:q-prism}, and let us denote the quotient ring $\Z_p[[q-1]] / ( [p]_q )$ by $\Z_p[ \zeta_p ]$
(where $\zeta_p$ is the image of $q$). For every animated commutative $\Z_p$-algebra $R$ with $p$-completion $\widehat{R}$, we let 
$\widehat{R}[ \zeta_p ]$ denote the (derived) tensor product $\Z_p[\zeta_p] \widehat{\otimes}^{L}_{\Z_p} \widehat{R}$, and we let $\qOmega_{R} \in \widehat{\calD}( \Z_p[[q-1]] )$ denote the relative prismatic complex $\Prism_{ \widehat{R}[ \zeta_p ] / \Z_p[[q-1]] }$
of Construction \ref{construction:relative-prismatic-cohomology}. We will refer to $\qOmega_{R}$ as the {\it $q$-de Rham complex of $R$}.
\end{definition}

\begin{remark}
Let $R$ be a $p$-complete commutative ring which is $p$-completely smooth over $\Z_p$. By virtue of Theorem~16.17 of \cite{prisms}, the $q$-de Rham complex $\qOmega_{R}$ can be computed
using the {\it $q$-crystalline site} of $R$ (see Definition~16.12 of \cite{prisms}). In particular, when given an \'{e}tale coordinate system $\Spf(R) \rightarrow \mathbf{A}^{n} \times \Spf(\Z_p)$,
then $\qOmega_{R}$ is modeled by an explicit $q$-deformation of the $p$-complete de Rham complex $( \widehat{\Omega}^{\ast}_{R/\Z_p}, d)$ (see Theorem~16.21 of \cite{prisms}).
\end{remark}

For every animated commutative ring $R$, the action of the profinite group $\Z_p^{\times}$ on the $q$-de Rham prism $(\Z_p[[q-1]], ( [p]_q ) )$ induces an action of
$\Z_p^{\times}$ on the $q$-de Rham complex $\qOmega_{R}$. In what follows, it will be convenient to restrict our attention to the subcomplex given by the homotopy invariants
for the finite subgroup $\F_p^{\times} \subseteq \Z_p^{\times}$ of $(p-1)$st roots of unity. This subcomplex also has a description in terms of relative prismatic cohomology:

\begin{construction}[$\widetilde{p}$-de Rham Cohomology]\label{construction:tilde-p-deRham}
Let $(\Z_p[[ \slashp]], (\slashp) )$ be the prism of Notation \ref{notation:reduced-q-de-Rham-prism}. For every animated commutative ring $R$, we regard the
$p$-completion $\widehat{R}$ as an animated commutative algebra over the quotient ring $\Z_p[[ \slashp]] / (\slashp) \simeq \Z_p$, and we let
$\slashOmega_{R}$ denote the relative prismatic complex $\Prism_{ \widehat{R} / \Z_p[[ \slashp ]] }$ of Construction \ref{construction:relative-prismatic-cohomology},
which we regard as an object of $\widehat{\calD}( \Z_p[[ \slashp]] )$. We will refer to $\slashOmega_{R}$ as the {\it $\slashp$-de Rham complex of $R$}.
\end{construction}

\begin{remark}[Relationship with $q$-de Rham Cohomology]\label{remark:q-de-Rham-tilde}
For every animated commutative ring $R$, Remark \ref{remark:change-of-prism} supplies a canonical isomorphism 
$$ \Z_p[[q-1]] \widehat{\otimes}^{L}_{ \Z_p[[\slashp]] } \slashOmega_{R} \xrightarrow{\sim} \qOmega_{R}$$
in the $\infty$-category $\widehat{\calD}( \Z_p[[q-1]] )$, which is equivariant with respect to the action of the profinite group $\Z_p^{\times}$.
Passing to homotopy fixed points for the action of the finite subgroup $\F_p^{\times} \subset \Z_p^{\times}$ (which acts trivially on $\slashOmega_{R}$), we obtain an isomorphism
$$ \slashOmega_{R} \simeq \RGamma( \F_p^{\times}, \qOmega_R )$$
in the $\infty$-category $\widehat{\calD}( \Z_p[[\slashp]] )$ which is equivariant with respect to the action of the $p$-profinite group $(1+p\Z_p)^{\times}$.
Consequently, the $q$-de Rham complex $\qOmega_{R}$ and the $\widetilde{p}$-de Rham complex $\slashOmega_{R}$ contain essentially the same information.
\end{remark}

\begin{variant}[Globalization to Formal Schemes]
Let $\mathfrak{X}$ be a bounded $p$-adic formal scheme. We let $\RGamma_{\slashdR}( \mathfrak{X} )$ denote the relative prismatic complex 
$$ \RGamma_{\Prism}( \mathfrak{X} / \Z_p[[\slashp]] ) \in \widehat{\calD}( \Z_p[[ \slashp ]] ).$$
We will refer to $\RGamma_{\slashdR}( \mathfrak{X} )$ as the {\it $\slashp$-de Rham complex of $\mathfrak{X}$}. We denote cohomology of
$\RGamma_{\slashdR}( \mathfrak{X} )$ by $\mathrm{H}^{\ast}_{\slashdR}( \mathfrak{X} )$, which we refer to as the {\it $\slashp$-de Rham cohomology of $\mathfrak{X}$}.
Note that, in the special case where $\mathfrak{X} = \Spf(R)$ is an affine formal scheme, we can identify $\RGamma_{\slashdR}( \mathfrak{X} )$ with the $\widetilde{p}$-de Rham complex $\slashOmega_{R}$ of Construction \ref{construction:tilde-p-deRham}.
\end{variant}

\begin{remark}[Relationship with Diffracted Hodge Cohomology]\label{remark:diffracted-Hodge-extended-scalars}
Let $R$ be an animated commutative ring, and let $\widehat{\Omega}^{\DHod}_{R}$ denote the $p$-complete diffracted Hodge complex of $R$.
By virtue of Example \ref{example:diffracted-Hodge-as-relative-prismatic}, the complex $\widehat{\Omega}^{\DHod}_{R}$ can be obtained from
the $\widetilde{p}$-de Rham complex $\slashOmega_{R}$ by (derived) extension of scalars along the quotient map 
$\Z_p[[ \slashp ]] \twoheadrightarrow \Z_p[[ \slashp ]] / ( \slashp ) \simeq \Z_p$. Similarly, if $\mathfrak{X}$ is a bounded formal scheme,
then the complex $\RGamma(\mathfrak{X}, \Omega^{\DHod}_{\mathfrak{X}} )$ is obtained from $\RGamma_{\slashdR}( \mathfrak{X} )$ by
extending scalars along $\Z_p[[ \slashp]] \twoheadrightarrow \Z_p$. In particular, we have a long exact sequence of cohomology groups
$$ \cdots \rightarrow \mathrm{H}^{\ast-1}( \mathfrak{X}, \Omega_{\mathfrak{X}}^{\DHod})
\rightarrow \mathrm{H}^{\ast}_{\slashdR}( \mathfrak{X} ) \xrightarrow{\slashp}
\mathrm{H}^{\ast}_{\slashdR}( \mathfrak{X} )  \rightarrow
\mathrm{H}^{\ast}( \mathfrak{X}, \Omega^{\DHod}_{\mathfrak{X}} ) \rightarrow  \mathrm{H}^{\ast+1}_{\slashdR}( \mathfrak{X} ) \rightarrow \cdots$$
\end{remark}

\begin{remark}[Relationship with Derived de Rham Cohomology]
Let $R$ be an animated commutative ring, and let $\widehat{\dR}_{R}$ denote the {\it $p$-complete derived de Rham complex of $R$}
(see Construction \ref{construction:derived-de-Rham}). By virtue of the de Rham-crystalline and prismatic crystalline comparisons (see Theorem~\ref{theorem:deduce-comparison}), the complex $\widehat{\dR}_{R}$ can be computed as
the prismatic cohomology of $\F_p \otimes^{L} R$ relative to the crystalline prism $(\Z_p, (p) )$. Applying Remark \ref{remark:change-of-prism} to the map of prisms
$$ (\Z_p[[ \slashp]], ( \slashp ) ) \rightarrow ( \Z_p, (p) ) \quad \quad \slashp \mapsto p,$$
we see that $\widehat{\dR}_{R}$ can be obtained from $\slashOmega_{R}$ by (derived) extension of scalars along the quotient map
$\Z_p[[\slashp]] \twoheadrightarrow \Z_p[[ \slashp]] / ( \slashp - p) \simeq \Z_p$. Similarly, if $\mathfrak{X}$ is a bounded $p$-adic formal scheme, then
the derived de Rham complex $\RGamma_{\dR}( \mathfrak{X} / \Z_p )$ can obtained from $\RGamma_{\slashdR}( \mathfrak{X} )$ by extension of scalars
along $\Z_p[[ \slashp]] \twoheadrightarrow \Z_p$. In particular, we have a long exact sequence of cohomology groups 
$$ \cdots \rightarrow \mathrm{H}^{\ast-1}_{\dR}( \mathfrak{X})
\rightarrow \mathrm{H}^{\ast}_{\slashdR}( \mathfrak{X} ) \xrightarrow{\slashp-p}
\mathrm{H}^{\ast}_{\slashdR}( \mathfrak{X} )  \rightarrow
\mathrm{H}^{\ast}_{\dR}( \mathfrak{X} ) \rightarrow  \mathrm{H}^{\ast+1}_{\slashdR}( \mathfrak{X} ) \rightarrow \cdots$$
\end{remark}

Note that the action of $\Z_p^{\times}$ on the $q$-de Rham prism $(\Z_p[[q-1]], ( [p]_q ) )$ restricts to an action of the $p$-profinite group $(1+p\Z_p)^{\times}$ on
$\Z_p[[\slashp]]$. Moreover, the quotient map
$$ \Z_p[[ \slashp ]] \twoheadrightarrow \Z_p \quad \quad \slashp \mapsto p$$
is invariant under the action of $(1+p\Z_p)^{\times}$. It follows that, for every animated commutative ring $R$, the complex $\slashOmega_{R}$ is equipped with a semilinear
action of $(1+p\Z_p)^{\times}$. For each element $u \in (1+p\Z_p)^{\times}$, we write $\slashOmega_{R}^{u=1}$ for the fiber of the map
$(\gamma_u - \id): \slashOmega_{R} \rightarrow \slashOmega_{R}$, where $\gamma_u$ is the automorphism of $\slashOmega_{R}$ determined by $u$.
Applying Theorem \ref{theorem:qdR-WCart-refined-reformulation} to the prismatic cohomology sheaf $\mathscr{H}_{\Prism}(R)$ of Construction \ref{construction:prismatic-cohomology-sheaves}, we obtain the following:

\begin{proposition}\label{proposition:concrete-qdR}
Let $p$ be an odd prime, let $u$ be a topological generator for the profinite group $(1+p\Z_p)^{\times}$. For
every animated commutative ring $R$, the diagram of stacks (\ref{equation:qdR-WCart-refined})
induces a pullback diagram
$$ \xymatrix@R=50pt@C=50pt{ \Prism_{R} \ar[r] \ar[d] & \widehat{\dR}_{R} \ar[d] \\
\slashOmega_{R}^{u=1} \ar[r] & \widehat{\dR}_{R}^{u=1} }$$
in the $\infty$-category $\widehat{\calD}(\Z_p)$.
\end{proposition}

\begin{remark}\label{remark:concrete-qdR}
Let $R$ be an animated commutative ring, let $u$ be a topological generator of the profinite group $(1+p\Z_p)^{\times}$, and let
$\gamma_{u}$ denote the associated automorphism of the complex $\slashOmega_{R}$. We can state Proposition \ref{proposition:concrete-qdR} more informally as follows:
the endomorphism $\gamma_{u} - \id$ is canonically divisible by the difference $\widetilde{p} - p$, and we have a fiber sequence
$$ \Prism_{R} \rightarrow \slashOmega_{R} \xrightarrow{ \frac{ \gamma_u - \id}{ \widetilde{p} - p} } \slashOmega_{R}$$
in the $\infty$-category $\widehat{\calD}(\Z_p)$. Similarly, if $\mathfrak{X}$ is a bounded $p$-adic formal scheme, we have a fiber sequence
$$ \RGamma_{\Prism}( \mathfrak{X} ) \rightarrow \RGamma_{\slashdR}(\mathfrak{X}) \xrightarrow{ \frac{ \gamma_u - \id}{ \widetilde{p} - p}  } \RGamma_{\slashdR}(\mathfrak{X} ),$$
which yields a long exact sequence of cohomology groups
$$ \cdots \rightarrow \mathrm{H}^{\ast-1}_{\slashdR}( \mathfrak{X})
\rightarrow \mathrm{H}^{\ast}_{\Prism}( \mathfrak{X} ) \rightarrow
\mathrm{H}^{\ast}_{\slashdR}( \mathfrak{X} )  \xrightarrow{ \frac{ \gamma_u - \id}{ \slashp-p} } 
\mathrm{H}^{\ast}_{\slashdR}( \mathfrak{X} ) \rightarrow  \mathrm{H}^{\ast+1}_{\Prism}( \mathfrak{X} ) \rightarrow \cdots$$
\end{remark}

\subsection{Diffracted Hodge Cohomology of Schemes}\label{subsection:diffracted-Hodge-integral}

For the duration of this section, we suspend our convention that $p$ denotes a fixed prime number. Instead, we allow the prime number $p$ to vary and assemble the results of Construction \ref{construction:complete-diffracted-Hodge} to a single invariant, which we will refer to as the {\it diffracted Hodge complex}. Our construction will exploit Sullivan's arithmetic fracture square together with the observation that the Sen operator is rationally diagonalizable.

\begin{construction}[The Diffracted Hodge Complex]\label{construction:diffracted-Hodge-integral}
Let $R$ be an animated commutative ring. For every prime number $p$, we write $\widehat{\Omega}^{\DHod}_{R,p}$ for the $p$-complete diffracted Hodge complex of
$R$ and $\Fil_{\bullet}^{\conj} \widehat{\Omega}^{\DHod}_{R,p}$ for its conjugate filtration (Construction \ref{construction:complete-diffracted-Hodge}). 
For each integer $n$, Remark \ref{remark:diffracted-relative-comparison-affine} supplies an isomorphism of 
$\gr_{n}^{\conj} \widehat{\Omega}^{\DHod}_{R,p}$ with the $p$-completion of the complex $L \Omega^{n}_{R}[-n]$, which we (temporarily, only in this subsection) denote by $L\Omega^{n}_{R}[-n]^{\wedge}_{(p)}$ to emphasize its dependence on $p$. Let $\Theta: \Fil_{\bullet}^{\conj} \widehat{\Omega}^{\DHod}_{R,p} \rightarrow \Fil_{\bullet}^{\conj}  \widehat{\Omega}^{\DHod}_{R,p}$ be the Sen operator of Notation \ref{notation:Sen-operator-on-diffracted}, so that $\Theta$ acts on each $\gr_{n}^{\conj} \widehat{\Omega}^{\DHod}_{R,p}$ by multiplication by $-n$. In particular, regarding each $L\Omega^n_R[-n]$ (and thus its $p$-completion $L\Omega^{n}_{R}[-n]^{\wedge}_{(p)}$) as a $\mathbf{Z}[\Theta]$-complex where $\Theta$ acts by $-n$, we obtain unique  $\mathbf{Z}[\Theta]$-equivariant isomorphisms
\begin{equation}
\label{eq:DiffHodZ1}
\Fil_{n}^{\conj} \widehat{\Omega}^{\DHod}_{R,p} \simeq  \bigoplus_{m=0}^n L\Omega^{m}_{R}[-m]^{\wedge}_{(p)} \quad \text{for} \ 0 \leq n < p 
\end{equation}
and 
\begin{equation}
\label{eq:DiffHodZ2}
 \Fil_{n}^{\conj} \widehat{\Omega}^{\DHod}_{R,p}[1/p] \simeq  \bigoplus_{i=0}^n L\Omega^{n}_{R}[-n]^{\wedge}_{(p)}[1/p]  \quad \forall n \geq 0
 \end{equation}
 inducing the isomorphism $\gr_{n}^{\conj} \widehat{\Omega}^{\DHod}_{R,p} \simeq L\Omega^{n}_{R}[-n]^{\wedge}_{(p)}$ of Remark \ref{remark:diffracted-relative-comparison-affine} on the generalized eigenspaces of $\Theta$. This allows us to form a pullback diagram
$$ \xymatrix@R=50pt@C=50pt{ \Fil^{\conj}_n \Omega^{\DHod}_{R} \ar[r] \ar[d] & \prod_{p} \Fil_{n}^{\conj} \widehat{\Omega}^{\DHod}_{R,p} \ar[d] \\
\Q \otimes^{L} \left( \bigoplus_{m=0}^{n}  L\Omega^{m}_{R}[-m] \right) \ar[r] & \Q \otimes \left(\prod_{p} \Fil_{n}^{\conj} \widehat{\Omega}^{\DHod}_{R,p}\right)  }$$
where the products on the right are taken over the set of all prime numbers $p$, the bottom horizontal map comes from \eqref{eq:DiffHodZ1} and \eqref{eq:DiffHodZ2}.
This construction depends functorially on $n$, and therefore determines a diagram
$$ \Fil_{0}^{\conj} \Omega^{\DHod}_{R} \rightarrow \Fil_{1}^{\conj} \Omega^{\DHod}_{R} \rightarrow \Fil_{2}^{\conj} \Omega^{\DHod}_{R} \rightarrow \cdots$$
We will refer to the underlying complex $\varinjlim_{n} \Fil_{n}^{\conj} \Omega^{\DHod}_{R}$ as the {\it diffracted Hodge complex} of $R$ and denote it by $\Omega^{\DHod}_{R}$. The preceding diagram determines an increasing filtration of $\Omega^{\DHod}_{R}$, which we denote by $\Fil^{\conj}_{\bullet} \Omega^{\DHod}_{R}$ and refer to as the {\it conjugate filtration} on $\Omega^{\DHod}_{R}$.
\end{construction}

\begin{remark}
Let $R$ be an animated commutative ring and let $\Omega^{\DHod}_{R}$ be the diffracted Hodge complex of
Construction \ref{construction:diffracted-Hodge-integral}. For every prime number $p$, the $p$-completion of
$\Omega^{\DHod}_{R}$ can be identified with the $p$-complete diffracted Hodge complex $\widehat{\Omega}^{\DHod}_{R,p}$ of 
Construction \ref{construction:complete-diffracted-Hodge}.
\end{remark}

\begin{remark}[The Associated Graded for the Conjugate Filtration]\label{remark:conjugate-filtration-of-diffracted-Hodge}
For every integer $d \geq 0$, we let $\gr_{d}^{\conj} \Omega^{\DHod}_{R}$ denote the cofiber of the natural map $\Fil_{d-1}^{\conj} \Omega^{\DHod}_{R} \rightarrow \Fil_{d}^{\conj} \Omega^{\DHod}_{R}$.
We then have a pullback diagram
$$ \xymatrix@R=50pt@C=50pt{ \gr_{d}^{\conj} \Omega^{\DHod}_{R} \ar[r] \ar[d] & \prod_{p} \gr_{d}^{\conj} \widehat{\Omega}^{\DHod}_{R,p} \ar[d] \\
\Q \otimes^{L} L\Omega^{d}_{R}[-d] \ar[r] & \Q \otimes^{L} \prod_{p} L\Omega^{d}_{R}[-d]^{\wedge}_{(p)} }$$
which supplies an isomorphism $\gr_{d}^{\conj} \Omega^{\DHod}_{R} \simeq L\Omega^{d}_{R}[-d]$ (in the $\infty$-category $\calD(\Z)$).
\end{remark}

\begin{remark}\label{remark:conjugate-equals-Postnikov2}
Let $R$ be an animated commutative ring. For every integer $d \geq 0$, the complex
$L\Omega^{d}_{R}$ is concentrated in cohomological degrees $\leq d$, so that $\gr_{d}^{\conj} \Omega^{\DHod}_{R} \simeq L \Omega^{d}_{R}[-d]$ is concentrated in cohomological
degrees $\leq d$. It follows by induction that $\Fil_{d}^{\conj} \Omega^{\DHod}_{R}$ is also concentrated in cohomological degrees $\leq d$.

Let $R$ be a commutative ring for which the relative cotangent complex $L\Omega^{1}_{R}$ is a flat $R$-module (this condition is satisfied, for example, if $R$ is smooth over $\Z$). In this case, each $L\Omega^{d}_{R}$ can be identified with the $R$-module $\Omega^{d}_{R}$, concentrated in cohomological degree zero. It follows that each of the complexes
$\gr_{d}^{\conj} \Omega^{\DHod}_{R}$ is concentrated in cohomological degree $d$: that is, the conjugate filtration on $\Omega^{\DHod}_{R}$ coincides with 
its Postnikov filtration $d \mapsto \tau^{\leq d} \Omega^{\DHod}_{R}$.
\end{remark}

\begin{variant}\label{variant:coconnectivity-diffracted-hodge}
Let $R$ be a commutative ring for which the relative cotangent complex $L\Omega^{1}_{R}$ has $\Tor$-amplitude contained in the interval $[-1,0]$.
Then, for every integer $d$, the complex $L\Omega^{d}_{R}$ has $\Tor$-amplitude contained in the interval $[-d,0]$.
It follows by induction that each of the complexes $\Fil_{d}^{\conj} \Omega^{\DHod}_{R}$ has $\Tor$-amplitude contained in the interval $[0,d]$:
in particular, the cohomology groups of the complex $\Fil_{d}^{\conj} \Omega^{\DHod}_{R}$ are concentrated in degrees $\geq 0$.
Passing to the colimit over $d$, we deduce that the cohomology groups of the diffracted Hodge complex $\Omega^{\DHod}_{R}$ are concentrated in nonnegative degrees.
\end{variant}

\begin{proposition}\label{proposition:sifted-colimit-diffracted-Hodge}
For each integer $n$, the functor
$$ \CAlg^{\anim} \rightarrow \calD(\Z) \quad \quad R \mapsto \Fil_{n}^{\conj} \Omega^{\DHod}_{R}$$
commutes with sifted colimits and satisfies fpqc descent.
\end{proposition}

\begin{proof}
We proceed by induction on $n$, the case $n < 0$ being trivial. To carry out the inductive step, it suffices to observe that the functor
$R \mapsto \gr_{n}^{\conj} \Omega^{\DHod}_{R} \simeq L \Omega^{n}_{R}[-n]$ commutes with sifted colimits and satisfies fpqc descent
(see Remark \ref{remark:descent-for-LOmega}).
\end{proof}

\begin{corollary}
The functor 
$$ \CAlg^{\anim} \rightarrow \calD(\Z) \quad \quad R \mapsto \Omega^{\DHod}_{R}$$
which commutes with sifted colimits.
\end{corollary}

\begin{example}
Let $R$ be a smooth $\Z$-algebra. Then, for each integer $n$, the complex $\gr_{n}^{\conj} \Omega^{\DHod}_{R} \simeq \Omega^{n}_{R}[-n]$ is concentrated in cohomological degree $n$.
It follows that the conjugate filtration on $\Omega^{\DHod}_{R}$ coincides with the Postnikov filtration. Together with Proposition \ref{proposition:sifted-colimit-diffracted-Hodge},
this observation characterizes the conjugate filtration for a general animated commutative ring $R$.
\end{example}

\begin{remark}[The Sen Operator]\label{remark:Sen-operator-integral}
Let $R$ be an animated commutative ring. Then the diffracted Hodge complex $\Omega^{\DHod}_{R}$ of Construction \ref{construction:diffracted-Hodge-integral} comes equipped with
an endomorphism $\Theta: \Omega^{\DHod}_{R} \rightarrow \Omega^{\DHod}_{R}$; we will refer to $\Theta$ as the {\it Sen operator} on $\Omega^{\DHod}_{R}$.
By construction, it refines to an operator on the filtered complex $\Fil_{\bullet}^{\conj} \Omega^{\DHod}_{R}$, and the induced action on each $\gr_{n}^{\conj} \Omega^{\DHod}_{R} \simeq L\Omega^{n}_{R}[-n]$
is given by multiplication by $(-n)$. With a bit more effort, one can show that $\Theta$ underlies an action of the group scheme $\mathbf{G}_{m}^{\sharp}$ on the diffracted Hodge
complex $\Omega^{\DHod}_{R}$: that is, $\Omega^{\DHod}_{R}$ can be promoted to a quasi-coherent complex on the classifying stack $B\mathbf{G}_{m}^{\sharp}$.
\end{remark}

\begin{remark}\label{remark:factor-Theta-plus}
Let $R$ be an animated commutative ring. Since the Sen operator $\Theta$ acts by multiplication by $(-n)$ on the complex $\gr_{n}^{\conj} \Omega^{\DHod}_{R}$,
the composite map
$$ \Fil_{n}^{\conj} \Omega^{\DHod}_{R} \xrightarrow{\Theta+n} \Fil_{n}^{\conj} \Omega^{\DHod}_{R} \rightarrow \gr_{n}^{\conj} \Omega^{\DHod}_{R}$$
is nullhomotopic. In fact, there is an essentially unique nullhomotopy which depends functorially on $R$. To see this, we can reduce to the case
where $R$ is a polynomial algebra over $\Z$ (Proposition \ref{proposition:sifted-colimit-diffracted-Hodge}). In this case, the desired result follow
from the observation that for any morphism $f: R \rightarrow S$ in $\Poly_{\Z}$, the mapping space
 $\Hom_{ \calD(R) }(  \Fil_{n}^{\conj} \Omega^{\DHod}_{R}, \gr_{n}^{\conj} \Omega^{\DHod}_{S} )$ is discrete
(since $\Fil_{n}^{\conj} \Omega^{\DHod}_{R}$ is concentrated in cohomological degrees $\leq n$ and
$\gr_{n}^{\conj} \Omega^{\DHod}_{S}$ is concentrated in cohomological degree $n$; see Remark \ref{remark:conjugate-equals-Postnikov2}).
It follows that $\Theta+n$ factors through a map
$$ \Fil_{n}^{\conj} \Omega^{\DHod}_{R}  \rightarrow \Fil_{n-1}^{\conj} \Omega^{\DHod}_{R},$$
which (by slight abuse of notation) we will also denote by $\Theta+n$.
\end{remark}

\begin{remark}
For every animated commutative ring $R$, the filtered complex $\Fil_{\bullet}^{\conj} \Omega^{\DHod}_{R}$ can be regarded as a commutative algebra object in the $\infty$-category
of filtered objects of $\calD(\Z)$. Moreover, the isomorphism $\Fil_{0}^{\conj} \Omega^{\DHod}_{R} \simeq L \Omega^{0}_{R} \simeq R$ is an isomorphism of commutative
algebra objects of $\calD(\Z)$. It follows that each $\Fil_{n}^{\conj} \Omega^{\DHod}_{R}$ can be promoted to an object of the $\infty$-category $\calD(R)$, and
that $\Omega^{\DHod}_{R}$ can be regarded as a commutative algebra object of $\calD(R)$.
\end{remark}

We now introduce a global version of Construction \ref{construction:diffracted-Hodge-integral}. We will make use of the following observation:

\begin{proposition}\label{proposition:etale-base-change-diffracted}
Let $f: R \rightarrow S$ be a morphism of animated commutative rings which is formally \'{e}tale: that is, for which the relative cotangent complex
$L\Omega^{1}_{S/R}$ vanishes. Then, for every integer $n$, the natural map $S \otimes^{L}_{R} \Fil_{n}^{\conj} \Omega^{\DHod}_{R} \rightarrow  \Fil_{n}^{\conj} \Omega^{\DHod}_{S}$
is an isomorphism in the $\infty$-category $\calD(S)$. In particular, the natural map $S \otimes^{L}_{R} \Omega^{\DHod}_{R} \rightarrow \Omega^{\DHod}_{S}$
is an isomorphism in $\calD(S)$.
\end{proposition}

\begin{proof}
Proceeding by induction on $n$, we are reduced to showing that the natural map $\theta: S \otimes_{R}^{L} \gr_{n}^{\conj}  \Omega^{\DHod}_{R} \rightarrow
\gr_{n}^{\conj} \Omega^{\DHod}_{S}$ is an isomorphism in $\calD(S)$. We now observe that, up to a cohomological shift, $\theta$ can be identified
with the $n$th (derived) exterior power of the comparison map $S \otimes_{R}^{L} L\Omega^{1}_{R} \rightarrow L \Omega^{1}_{S}$, which is an isomorphism
by virtue of the vanishing of $L \Omega^{1}_{S/R}$.
\end{proof}

\begin{corollary}
The functor
$$ \CAlg^{\anim} \rightarrow \calD(\Z) \quad \quad R \mapsto \Omega^{\DHod}_{R}$$
satisfies descent for the \'{e}tale topology.
\end{corollary}

\begin{notation}
Let $X$ be a scheme, and let $\calU$ denote the collection of all affine open subsets of $X$. For every integer $n$, 
Proposition \ref{proposition:etale-base-change-diffracted} guarantees that the construction
$$ (U \in \calU) \mapsto \Fil_{n}^{\conj} \Omega^{\DHod}_{ \calO_{X}(U) }$$
determines an object of the derived $\infty$-category $\calD(X)$, which we will denote by $\Fil_{n}^{\conj} \Omega^{\DHod}_{X}$.
Allowing $n$ to vary, we obtain a direct system
$$\calO_{X} \simeq \Fil_{0}^{\conj} \Omega^{\DHod}_{X} \rightarrow \Fil_{1}^{\conj} \Omega^{\DHod}_{X} \rightarrow \Fil_{2}^{\conj} \Omega^{\DHod}_{X} \rightarrow \cdots$$
whose colimit we will denote by $\Omega^{\DHod}_{X}$ and refer to as the {\it diffracted Hodge complex} of $X$. We will refer to the (hypercohomology) groups
$\mathrm{H}^{\ast}( X, \Omega^{\DHod}_{X} )$ as the {\em diffracted Hodge cohomology groups} of $X$.
\end{notation}

\newpage \section{The Nygaard Filtration}\label{section:Nygaard-filtration}

Let $k$ be a perfect field of characteristic $p$ and let $X$ be a smooth $k$-scheme. Then the crystalline
cohomology of $X$ can be computed as the hypercohomology of $X$ with coefficients in the {\it de Rham-Witt complex}
$$ W\Omega^{\ast}_{X} = ( W \Omega^{0}_{X} \rightarrow W \Omega^{1}_{X} \rightarrow W \Omega^{2}_{X} \rightarrow \cdots )$$
of Bloch and Deligne-Illusie (see \cite{MR565469}, or \cite{DRW} for an alternative approach). In \cite{nygaard},
Nygaard introduces a filtration of $W\Omega^{\ast}_{X}$ by subcomplexes
$$ \cdots \rightarrow \calN^{2} W\Omega^{\ast}_{X} \subseteq \calN^{1} W\Omega^{\ast}_{X}
\subseteq \calN^{0} W \Omega^{\ast}_{X} = W \Omega^{\ast}_{X},$$
which we will refer to as the {\it Nygaard filtration} on $W \Omega^{\ast}_{X}$. This construction has the following feature:
\begin{itemize}
\item[$(\ast)$] For each integer $m$, multiplication by $p$ carries each $\calN^{m-1} W\Omega^{\ast}_{X}$ into $\calN^{m} W\Omega^{\ast}_{X}$.
Moreover, the quotient complex $\calN^{m} W\Omega^{\ast}_{X} / p \calN^{m-1} W\Omega^{\ast}_{X}$
is quasi-isomorphic to the truncated de Rham complex
$$ \Omega^{\geq m}_X := \left(\Omega^{m}_{X} \xrightarrow{d} \Omega^{m+1}_{X} \xrightarrow{d} \Omega^{m+2}_{X} \rightarrow \cdots\right)[-m].$$
\end{itemize}

Note that the crystalline cochain complex $\RGamma_{\crys}(X / \Z_p) = \RGamma(X, W\Omega^{\ast}_{X})$ 
can be identified with the Frobenius-twisted prismatic complex $F^{\ast} \RGamma_{\Prism}( X / W(k) )$
(see Warning \ref{warning:Frob-semilinearity}). In \cite{prisms}, the first author and Scholze introduced an analogue of the
Nygaard filtration for an arbitrary bounded prism $(A,I)$. To every $p$-adic formal scheme $\mathfrak{X}$ which is smooth
over $\overline{A} = A/I$, their construction supplies a filtered complex
$$ \cdots \rightarrow \Fil^{2}_{\Nyg} F^{\ast} \RGamma_{\Prism}( \mathfrak{X} / A)\rightarrow \Fil^{1}_{\Nyg} F^{\ast} \RGamma_{\Prism}( \mathfrak{X} / A) 
\rightarrow \Fil^{0}_{\Nyg} F^{\ast} \RGamma_{\Prism}( \mathfrak{X} / A) \simeq F^{\ast} \RGamma_{\Prism}( \mathfrak{X} / A)$$
satisfying the following variant of $(\ast)$: 
\begin{itemize}
\item[$(\ast')$] For every integer $m$, there is a canonical fiber sequence
$$ I \otimes_{A}^{L} \Fil^{m-1}_{\Nyg} F^{\ast} \RGamma_{\Prism}( \mathfrak{X} / A) \rightarrow
\Fil^{m}_{\Nyg} F^{\ast} \RGamma_{\Prism}( \mathfrak{X} / A) \rightarrow \RGamma( \mathfrak{X}, \widehat{\Omega}^{\geq m}_{\mathfrak{X} / \overline{A} } )$$
\end{itemize}
We review the construction of the filtered complex $\Fil^{\bullet}_{\Nyg} F^{\ast} \RGamma_{\Prism}( \mathfrak{X} / A)$ in \S\ref{subsection:relative-Nygaard}. To simplify the discussion, we will focus primarily on the case where $\mathfrak{X} = \Spf(R)$ is an affine formal scheme,
in which case we can identify $\RGamma_{\Prism}( \mathfrak{X} / A)$ with the relative prismatic complex $\Prism_{R/A}$ of 
Construction \ref{construction:relative-prismatic-cohomology}. However, we place no restrictions $R$, which we allow to
be an arbitrary (animated) commutative $\overline{A}$-algebra. 
In \S\ref{subsection:relative-dr-comparison}, we apply this result to prove a version of $(\ast')$ for an arbitrary bounded prism
(Corollary \ref{corollary:relative-Nygaard-deRham-refined}).
In \S\ref{subsection:Nygaard-filtration-comparison},
we specialize to the case where $(A,I)$ is the crystalline prism $(\Z_p, (p))$, and show that our construction agrees with the classical Nygaard filtration on the de Rham-Witt complex whenever $R$ is a regular $\F_p$-algebra (Proposition \ref{proposition:concrete-Nygaard-regular}).

Our primary goal in this section is to construct a counterpart of the Nygaard filtration in the setting of {\em absolute} prismatic cohomology.
Let $R$ be a commutative ring. Assume for simplicity that $p$ is not a zero-divisor in $R$, and let $\overline{R}$ denote the quotient
ring $R/pR$. Recall that the absolute prismatic complex $\Prism_{R}$ can be realized as the global sections of a quasi-coherent complex $\mathscr{H}_{\Prism}(R)$ on the Cartier-Witt stack $\WCart$ 
(Construction \ref{construction:absolute-prismatic-cohomology-general}). By virtue of Theorem \ref{theorem:Frobenius-pushout-square},
$\Prism_{R}$ can be realized as the pullback of a diagram
\begin{align}\label{equation:pullback-prism} \RGamma( \WCart, F^{\ast} \mathscr{H}_{\Prism}(R) ) & \rightarrow &
& \RGamma( \WCart^{\mathrm{HT}} , (F^{\ast} \mathscr{H}_{\Prism}(R))|_{\WCart^{\mathrm{HT}}} )
& \leftarrow & \RGamma( \Spf(\Z_p), \rho_{\dR}^{\ast} \mathscr{H}_{\Prism}(R) ). \end{align}
Note that the complex on the right can be identified with the relative prismatic complex $\Prism_{\overline{R} / \Z_p}$.
If $\overline{R}$ is a smooth $\F_p$-algebra, then Theorem \ref{theorem:crystalline-comparison} identifies $\Prism_{ \overline{R} / \Z_p }$
with the crystalline cochain complex $\RGamma_{\crys}( \overline{R} / \Z_p )$. By a classical result of Berthelot, $\RGamma_{\crys}( \overline{R} / \Z_p )$ has an explicit model at the level of cochain complexes, given by the $p$-complete de Rham complex $( \widehat{\Omega}^{\ast}_{R}, d)$.
More generally, if $R$ is an arbitrary $p$-torsion-free commutative ring, then there is a canonical isomorphism
$\Prism_{ \overline{R} / \Z_p } \simeq \widehat{\dR}_{R}$, where $\widehat{\dR}_{R}$ denotes the $p$-complete derived de Rham complex of $R$. In \S\ref{subsection:absolute-de-Rham-comparison}, we give a direct construction of this isomorphism, which does not make reference
to the theory of crystalline cohomology (Theorem \ref{theorem:deduce-comparison}). In \S\ref{subsection:absolute-nygaard-filtration}, we apply this
construction to refine (\ref{equation:pullback-prism}) to a diagram of filtered complexes
\begin{align}\label{equation:pullback-prism2} \RGamma( \WCart, \Fil^{\bullet}_{\Nyg} F^{\ast} \mathscr{H}_{\Prism}(R) ) & \rightarrow &
& \RGamma( \WCart^{\mathrm{HT}} , (\Fil^{\bullet}_{\Nyg} F^{\ast} \mathscr{H}_{\Prism}(R))|_{\WCart^{\mathrm{HT}}} )
& \leftarrow & \Fil^{\bullet}_{\Hodge} \widehat{\dR}_{R}, \end{align}
where $\Fil^{\bullet}_{\Nyg} F^{\ast} \mathscr{H}_{\Prism}(R)$ is a filtered complex on $\WCart$ obtained from the relative theory.
We denote the pullback of this diagram by $\Fil^{\bullet}_{\Nyg} \Prism_{R}$, which we refer to as the {\it absolute Nygaard filtration}
on the prismatic complex $\Prism_{R}$ (Construction \ref{construction:absolute-Nygaard-untwisted}).

Let $(A,I)$ be a bounded prism. It follows
immediately from our construction that, if $R$ is an algebra over the quotient ring $\overline{A} = A/I$, then the composite map
$\Prism_{R} \rightarrow \Prism_{R/A} \rightarrow F^{\ast} \Prism_{R/A}$ can be refined to a map of filtered complexes
$\Fil^{\bullet}_{\Nyg} \Prism_{R} \rightarrow \Fil^{\bullet}_{\Nyg} F^{\ast} \Prism_{R/A}$. In \S\ref{subsection:relative-nygaard-comparison},
we show that this map is an isomorphism in the case where the prism $(A,I)$ is perfect (Theorem \ref{theorem:compare-Nygaard-filtrations}).

\begin{remark}
Let $R$ be a quasiregular semiperfectoid ring, so that the absolute prismatic complex $\Prism_{R}$ can be identified with
a commutative ring, concentrated in cohomological degree zero. In \cite{prisms}, the Nygaard filtration of $\Prism_{R}$ is
defined as the system of ideals
$$ \cdots \subseteq \Fil^{2}_{\Nyg} \Prism_{R} \subseteq \Fil^{1}_{\Nyg} \Prism_{R} \subseteq \Fil^{0}_{\Nyg} \Prism_{R} = \Prism_{R},$$
given by the formula
\begin{eqnarray}\label{equation:absolute-nygaard-defined} \Fil^{n}_{\Nyg} \Prism_{R} & = & \{ x \in \Prism_{R}: \varphi(x) \in \Prism_{R}^{[n]} \}
\end{eqnarray}
We will show in \S\ref{subsection:relative-nygaard-comparison} that this agrees with the definition of the absolute Nygaard filtration
given in this paper (Corollary \ref{corollary:absolute-Nygaard-qrsp}). From the present perspective, this is not obvious from the definition: it depends on a nontrivial calculation.

It is possible to adopt another approach to the theory of the absolute Nygaard filtration. One could take (\ref{equation:absolute-nygaard-defined})
as a definition in the case where $R$ is quasiregular semiperfectoid, and extend to general $p$-quasisyntomic commutative rings using descent.
However, from this perspective, the relationship between the Nygaard filtration on $\Prism_{R}$ and the Hodge filtration on $\widehat{\dR}_{R}$
is somewhat more opaque. Moreover, the calculations of \S\ref{subsection:relative-nygaard-comparison} are not avoided: they
reappear in the analysis of the connectivity properties of the relative Frobenius map $\varphi: \Fil^{\bullet}_{\Nyg} \Prism_{R} \rightarrow
\Prism^{[\bullet]}_{R}$, which we will study in \S\ref{subsection:Frobenius}.
\end{remark}

Let $R$ be a commutative ring. We let $\widehat{\Prism}_{R}$ denote the {\em Nygaard-completed} prismatic complex of
$R$, defined as the limit of the diagram
$$ \cdots \rightarrow \Prism_{R} / \Fil^{3}_{\Nyg} \Prism_{R} \rightarrow \Prism_{R} / \Fil^{2}_{\Nyg} \Prism_{R} 
\rightarrow \Prism_{R} / \Fil^{1}_{\Nyg} \Prism_{R} \simeq R.$$
There is a tautological comparison map $\Prism_{R} \rightarrow \widehat{\Prism}_{R}$. In \S\ref{subsection:Nygaard-completion},
we show that this map is an isomorphism when $R$ is $p$-torsion-free and the quotient ring $R/pR$ is regular (Proposition \ref{proposition:smooth-absolute-Nygaard-complete}), but not in general (see Example \ref{example:non-nygaard-complete}).

\subsection{The Relative Nygaard Filtration}\label{subsection:relative-Nygaard}

Throughout this section, we fix a bounded prism $(A,I)$ and we let $\overline{A}$ denote the quotient ring $A/I$. If $R$ is an animated commutative $\overline{A}$-algebra,
we let $F^{\ast} \Prism_{R/A}$ denote the commutative algebra object of $\widehat{\calD}(A)$ obtained from the relative prismatic complex $\Prism_{R/A}$ by
(completed) extension of scalars along the Frobenius morphism $\varphi: A \rightarrow A$. In \S15 of \cite{prisms}, the complex $F^{\ast} \Prism_{R/A}$ is equipped with a decreasing filtration
$$ \cdots \rightarrow \Fil^{2}_{\Nyg} F^{\ast} \Prism_{R/A}
\rightarrow  \Fil^{1}_{\Nyg} F^{\ast} \Prism_{R/A} \rightarrow \Fil^{0}_{\Nyg} F^{\ast} \Prism_{R/A} = F^{\ast} \Prism_{R/A},$$
which we will refer to as the {\it relative Nygaard filtration}. Our goal in this section is to review the definition of this filtration in a form which will be convenient for our purposes. Note that the Frobenius endomorphism of $\Prism_{R/A}$ induces an $A$-linear morphism $\varphi: F^{\ast} \Prism_{R/A} \rightarrow \Prism_{R/A}$, which we will refer to as the {\em relative Frobenius}.

\begin{proposition}\label{proposition:relative-nygaard-characterization}
Let $R$ be an animated commutative $\overline{A}$-algebra. Then the relative Frobenius morphism $\varphi: F^{\ast} \Prism_{R/A} \rightarrow \Prism_{R/A}$
can be promoted to a morphism of filtered complexes $\Fil^{\bullet}(\varphi): \Fil^{\bullet}_{\Nyg} F^{\ast} \Prism_{R/A} \rightarrow I^{\bullet} \Prism_{R/A}$,
which is determined (up to essentially unique isomorphism) by the requirement that it depends functorially on $R$ and has the following additional properties:
\begin{itemize}
\item[$(1)$] For each integer $n$, the functor
$$ \CAlg^{\anim}_{\overline{A}} \rightarrow \widehat{\calD}( \overline{A} ) \quad \quad R \mapsto \Fil^{n}_{\Nyg} F^{\ast} \Prism_{R/A}$$
commutes with sifted colimits.

\item[$(2)$] Let $R$ be a finitely generated polynomial algebra over $\overline{A}$. Then, for every integer $n$, the map
 $$ \gr^{n}(\varphi): \gr^{n}_{\Nyg} F^{\ast} \Prism_{R/A} \rightarrow I^{n} / I^{n+1} \otimes^{L}_{A} \Prism_{R/A} = \overline{\Prism}_{R/A}\{n\}$$
 induces isomorphisms
 $$ \mathrm{H}^{m}( \gr^{n}_{\Nyg} F^{\ast} \Prism_{R/A} ) \xrightarrow{\sim} \begin{cases} \mathrm{H}^{m}_{\overline{\Prism}}(R/A)\{n\}
 \simeq \widehat{\Omega}^{m}_{R/\overline{A}}\{n-m\} & \text{ if } m \leq n \\
 0 & \text{ otherwise. } \end{cases}$$
\end{itemize}
\end{proposition}

\begin{remark}\label{remark:polynomial-to-smooth}
Let $R$ be an animated commutative $\overline{A}$-algebra. It follows from properties $(1)$ and $(2)$ of Proposition \ref{proposition:relative-nygaard-characterization} that, for every integer $n$, the relative Frobenius map $\gr^{n}(\varphi): \gr^{n}_{\Nyg} F^{\ast} \Prism_{R/A} \rightarrow \overline{\Prism}_{R/A}\{n\}$
factors through an isomorphism $\gr^{n}_{\Nyg} F^{\ast} \Prism_{R/A} \simeq \Fil_{n}^{\conj} \overline{\Prism}_{R/A}\{n\}$, depending
functorially on $R$. Combining this observation with Remark \ref{remark:conjugate-equals-postnikov-conditions}, we deduce
the following stronger version of $(2)$:
\begin{itemize}
\item[$(2')$] Let $R$ be an $\overline{A}$-algebra whose $p$-completion $\widehat{R}$ is $p$-completely smooth over $\overline{A}$.
For every integer $n$, the induced map $$ \gr^{n}_{\Nyg} F^{\ast} \Prism_{R/A} \rightarrow I^{n} / I^{n+1} \otimes^{L}_{A} \Prism_{R/A}
= \overline{\Prism}_{R/A}\{n\}$$
 induces isomorphisms
$$ \mathrm{H}^{m}( \gr^{n}_{\Nyg} F^{\ast} \Prism_{R/A} ) \xrightarrow{\sim} \begin{cases} \mathrm{H}^{m}_{\overline{\Prism}}(R/A)\{n\}
\simeq \widehat{\Omega}^{m}_{R/\overline{A}}\{n-m\} & \text{ if } m \leq n \\
0 & \text{ otherwise. } \end{cases}$$

\end{itemize}
\end{remark}

\begin{remark}[The relative Nygaard filtration is Beilinson connective]
\label{remark:Nygaard-Beilinson}
Let $R$ be an animated commutative $\overline{A}$-algebra. If the $p$-completion of $R$ is $p$-completely smooth over $\overline{A}$, then assertion (2') in Remark~\ref{remark:polynomial-to-smooth} together with the definition of the Beilinson $t$-structure (Definition~\ref{definition:Beilinson-t}) implies that the filtered object $\Fil^{\bullet}_{\Nyg} F^{\ast} \Prism_{R/A} \in \DFilt(A)$ is connective for the Beilinson $t$-structure. Using Proposition~\ref{proposition:relative-nygaard-characterization} (1), it then follows that for any animated $\overline{A}$-algebra $R$, the filtered object $\Fil^{\bullet}_{\Nyg} F^{\ast} \Prism_{R/A} \in \DFilt(A)$ is connective for the Beilinson $t$-structure.
\end{remark}

\begin{example}\label{example:relative-Nygaard-trivial}
In the special case $R = \overline{A}$, we can identify $F^{\ast} \Prism_{R/A}$ with the commutative ring $A$ (regarded as a chain complex
concentrated in degree zero). Under this identification, the Nygaard filtration is given by
$$ \Fil^{n}_{\Nyg} F^{\ast} \Prism_{R/A} = \begin{cases} I^{n} & \text{ if $n \geq 0$} \\
A & \text{ if $n \leq 0$. } \end{cases}$$
\end{example}

\begin{remark}\label{remark:Nygaard-degree-zero}
Let $R$ be an animated commutative $\overline{A}$-algebra. Combining the isomorphism $\gr^{n}_{\Nyg} F^{\ast} \Prism_{R/A} \simeq \Fil_{n}^{\conj} \overline{\Prism}_{R/A}\{n\}$
of Remark \ref{remark:polynomial-to-smooth} with the isomorphism $\Fil_0^{\conj} \overline{\Prism}_{R/A} \simeq \widehat{R}$ of Remark \ref{remark:algebra-structure-on-relative-prismatic-complex}, we obtain a canonical isomorphism $\gr^{0}_{\Nyg} F^{\ast} \Prism_{R/A} \simeq \widehat{R}$.
\end{remark}

\begin{remark}\label{remark:Nygaard-filtration-by-integers}
Let $R$ be an animated commutative $\overline{A}$-algebra. It follows from Remark \ref{remark:polynomial-to-smooth} that
the complex $\gr^{n}_{\Nyg} F^{\ast} \Prism_{R/A}$ vanishes for $n < 0$. Consequently, the maps
$$  F^{\ast} \Prism_{R/A} = \Fil^{0}_{\Nyg} F^{\ast} \Prism_{R/A} \rightarrow
 \Fil^{-1}_{\Nyg} F^{\ast} \Prism_{R/A} \rightarrow  \Fil^{-2}_{\Nyg} F^{\ast} \Prism_{R/A} \rightarrow \cdots$$
 are isomorphisms (in the $\infty$-category $\widehat{\calD}(A)$). It is therefore harmless to think of the relative Nygaard filtration of
Proposition \ref{proposition:relative-nygaard-characterization} as indexed by the nonnegative integers, rather than the collection of all integers.
\end{remark}

For the proof of Proposition \ref{proposition:relative-nygaard-characterization}, we need the following slight variant of Lemma \ref{lemma:qs-descent-relative}:

\begin{lemma}\label{lemma:qs-descent-for-pullback}
Let $\CAlg^{\QSyn}_{\overline{A}}$ denote the category of $p$-quasisyntomic $\overline{A}$-algebras.
The functor
$$ \CAlg^{\QSyn}_{\overline{A}} \rightarrow \widehat{\calD}(A) \quad \quad R \mapsto F^{\ast} \Prism_{R/A}$$
satisfies descent for the $p$-quasisyntomic topology.
\end{lemma}

\begin{proof}
We show more generally that for every object complex $M \in \widehat{\calD}(A)$ whose cohomology groups are
bounded below, the functor 
$$ \CAlg^{\QSyn}_{\overline{A}} \rightarrow \widehat{\calD}(A) \quad \quad R \mapsto M \widehat{\otimes}^{L}_{A} \Prism_{R/A}$$
satisfies descent for the $p$-quasisyntomic topology. Without loss of generality, we may assume that $M$ is the image of an object
of $\widehat{\calD}( \overline{A} )$, which we will also denote by $M$. Let $R$ be a $p$-quasisyntomic $\overline{A}$-algebra,
let $R \rightarrow R^{0}$ be a $p$-quasisyntomic cover, and let $R^{\bullet}$ be the cosimplicial $R$-algebra given by
the $p$-complete tensor powers of $R^{0}$ over $R$. We wish to show that the canonical map
$$ \theta: M \widehat{\otimes}^{L}_{\overline{A}} \overline{\Prism}_{R/A} \rightarrow 
\Tot^{\bullet}( M \widehat{\otimes}^{L} \overline{\Prism}_{R^{\bullet}/A} )$$
is an isomorphism in $\widehat{\calD}( \overline{A} )$. Since the degrees of the nonzero cohomology groups of the complexes
$M \widehat{\otimes}^{L} \overline{\Prism}_{R^{\bullet}/A}$ are uniformly bounded below, we can realize $\theta$ as a filtered colimit of
comparison maps
$$ \theta_{n}: M \widehat{\otimes}^{L}_{\overline{A}} \Fil_{n}^{\conj} \overline{\Prism}_{R/A} \to \Tot^{\bullet}( M \widehat{\otimes}^{L} \Fil_{n}^{\conj} \overline{\Prism}_{R^{\bullet}/A} ).$$
It will therefore suffice to show that each $\theta_n$ is an isomorphism. Let $M^{\bullet}$ denote the cosimplicial
object of $\widehat{\calD}(R)$ given by the completed tensor product $R^{\bullet} \widehat{\otimes}^{L}_{\overline{A}} M$
Proceeding by induction on $n$ and invoking the Hodge-Tate comparison (Remark \ref{remark:derived-HT-filtration}), we are reduced to checking that
each of the comparison maps
$$ M \widehat{\otimes}^{L}_{ \overline{A}} L\Omega^{n}_{R/ \overline{A}} \rightarrow
\Tot( M \widehat{\otimes}^{L}_{ \overline{A}}  L\Omega^{n}_{R^{\bullet}/ \overline{A}} )$$
is an isomorphism in $\widehat{\calD}( \overline{A} )$, which follows from Lemma~2.6 of \cite{bhatt-completions}.
\end{proof}

\begin{proof}[Proof of Proposition \ref{proposition:relative-nygaard-characterization}]
We proceed as in \S15 of \cite{prisms}. Let $\CAlg_{\overline{A}}^{\QSyn}$ denote the category of $p$-quasisyntomic
$\overline{A}$-algebras (Definition~\ref{definition:qsyn-topology}), and let $\calC \subseteq \CAlg_{\overline{A}}^{\QSyn}$ denote the full subcategory spanned by those
$p$-quasisyntomic $\overline{A}$-algebras $R$ for which the quotient ring $R/pR$ is generated by the images over $\overline{A}$ and $(R/pR)^{\flat}$
(that is, the collection of $\overline{A}$-algebras which are ``large'' in the sense of \cite{prisms}). 
For $R \in \calC$, the relative prismatic complex $\Prism_{R/A}$ can be identified with a $(p,I)$-completely flat $A$-algebra
(concentrated in cohomological degree zero), so that the Frobenius pullback $F^{\ast} \Prism_{R/A}$ can also be viewed as a commutative ring concentrated in cohomological degree zero. Let us write $\Fil^{n}_{\Nyg} F^{\ast} \Prism_{R/A}$ for the ideal of
$F^{\ast} \Prism_{R/A}$ consisting of those elements $x$ which satisfy $\varphi(x) \in I^{n} \Prism_{R/A}$, where $\varphi: F^{\ast} \Prism_{R/A} \rightarrow \Prism_{R/A}$ is the relative Frobenius map. By construction, the relative Frobenius map
$\varphi: F^{\ast} \Prism_{R/A} \rightarrow \Prism_{R/A}$ refines uniquely to a map of filtered abelian groups $\Fil^{\bullet}(\varphi): \Fil^{\bullet}_{\Nyg} F^{\ast} \Prism_{R/A} \rightarrow I^{\bullet} \Prism_{R/A}$, depending functorially on $R \in \calC$.

The construction $R \mapsto \{ \Fil^{n}_{\Nyg} F^{\ast} \Prism_{R/A} \}_{n \in \Z}$ then determines a functor from $\calC$ to the $\infty$-category $\DFiltI(A)$ of $I$-complete filtered complexes (Notation~\ref{notation:filtered-derived-I-complete}). This functor admits a right Kan extension $\CAlg_{\overline{A}}^{\QSyn} \rightarrow \DFiltI(A)$, which we will also denote by $R \mapsto \Fil^{\bullet}_{\Nyg} F^{\ast} \Prism_{R/A}$.
Note that $\calC$ forms a basis for the $p$-quasisyntomic topology on $\CAlg_{ \overline{A} }^{\QSyn}$. Consequently, to show that
a functor $U: \CAlg_{\overline{A}}^{\QSyn} \rightarrow \DFiltI(A)$ is a right Kan extension of its restriction to $\calC$,
it suffices to show that $U$ satisfies descent for the $p$-quasisyntomic topology.
Since these conditions are satisfied for the functor $R \mapsto F^{\ast} \Prism_{R/A}$ (Lemma \ref{lemma:qs-descent-for-pullback}), we have a canonical isomorphism $\Fil^{0}_{\Nyg} F^{\ast} \Prism_{R/A} \simeq F^{\ast} \Prism_{R/A}$ for $R \in \CAlg^{\QSyn}_{\overline{A}}$.
Similarly, Lemma \ref{lemma:qs-descent-relative} and Corollary \ref{salmage} guarantee that the functor
$$ \CAlg^{\QSyn}_{\overline{A}} \rightarrow \DFiltI(A) \quad \quad R \mapsto I^{\bullet} \Prism_{R/A}$$
is a right Kan extension of its restriction to $\calC$. It follows that the relative Frobenius map $\Fil^{\bullet}(\varphi): \Fil^{\bullet}_{\Nyg} F^{\ast} \Prism_{R/A} \rightarrow I^{\bullet} \Prism_{R/A}$ admits an essentially unique (functorial) extension to all $R \in \CAlg^{\QSyn}_{\overline{A}}$. 

Let $n$ be an integer. For each $R \in \CAlg^{\QSyn}_{\overline{A}}$, the relative Frobenius map $\Fil^{\bullet}(\varphi)$ induces a map
$\gr^{n}(\varphi): \gr^{n}_{\Nyg} F^{\ast} \Prism_{R/A} \rightarrow \overline{\Prism}_{R/A}\{n\}$. When $R \in \calC$, Theorem~15.2 of \cite{prisms} guarantees
that $\gr^{n}(\varphi)$ is a monomorphism of abelian groups whose image is the subgroup $\Fil_{n}^{\conj} \overline{\Prism}_{R/A}\{n\}$. Since
the functor
$$ \CAlg^{\QSyn}_{\overline{A}} \rightarrow \widehat{\calD}( \overline{A} ) \quad \quad R \mapsto \Fil_{n}^{\conj} \overline{\Prism}_{R/A}\{n\}$$
is also a right Kan extension of its restriction to $\calC$, it follows that there is an essentially unique isomorphism
$\beta_{R}: \gr^{n}_{\Nyg} F^{\ast} \Prism_{R/A} \simeq \Fil_{n}^{\conj} \overline{\Prism}_{R/A}\{n\}$ which depends functorially on
$R \in \CAlg^{\QSyn}_{\overline{A}}$ for which the composition $\gr^{n}_{\Nyg} F^{\ast} \Prism_{R/A}
\xrightarrow[\sim]{\beta_R} \Fil_{n}^{\conj} \overline{\Prism}_{R/A}\{n\} \rightarrow \overline{\Prism}_{R/A}\{n\}$ coincides
with $\gr^{n}(\varphi)$ (up to a homotopy depending functorially on $R$). From this observation, we immediately
deduce that $\Fil^{\bullet}(\varphi)$ satisfies requirement $(2)$ of Proposition \ref{proposition:relative-nygaard-characterization}.

We next claim that the functor
$$ \CAlg^{\QSyn}_{\overline{A}} \rightarrow \DFiltI(A) \quad \quad R \mapsto \Fil^{\bullet}_{\Nyg} F^{\ast} \Prism_{R/A}$$
is a left Kan extension of its restriction to the category $\Poly_{\overline{A}}$ of finitely generated polynomial algebras over $\overline{A}$.
Equivalently, for each $n \in \Z$, the functor $R \mapsto  \Fil^{n}_{\Nyg} F^{\ast} \Prism_{R/A}$ is a left Kan extension of its restriction to $\Poly_{\overline{A}}$.
To prove this, we proceed by induction on $n$, the case $n \leq 0$ being trivial (since the functor $R \mapsto F^{\ast} \Prism_{R/A}$ commutes with sifted colimits). 
To handle the inductive step, it suffices to observe that the functor $R \mapsto \gr^{n}_{\Nyg} F^{\ast} \Prism_{R/A}$ is a left Kan extension of its
restriction to $\Poly_{\overline{A}}$, which follows from the isomorphism $\gr^{n}_{\Nyg} F^{\ast} \Prism_{-/A} \simeq \Fil_{n}^{\conj} \overline{\Prism}_{-/A}\{n\}$
(since the functor $R \mapsto  \Fil_{n}^{\conj} \overline{\Prism}_{R/A}\{n\}$ also commutes with sifted colimits). 

Applying Proposition \ref{proposition:universal-of-animated}, we see that the functor $R \mapsto  \Fil^{\bullet}_{\Nyg} F^{\ast} \Prism_{R/A}$
admits an essentially unique extension to the $\infty$-category $\CAlg^{\anim}_{\overline{A}}$ of {\em all} animated commutative $\overline{A}$-algebras
which commutes with sifted colimits (and therefore satisfies requirement $(1)$ of Proposition \ref{proposition:relative-nygaard-characterization}). Moreover,
this functor is characterized by the requirement that it is a left Kan extension of its restriction to $\CAlg^{\QSyn}_{\overline{A}}$ (or equivalently to the subcategory
$\Poly_{\overline{A}}$), so the relative Frobenius map $\Fil^{\bullet}(\varphi):   \Fil^{\bullet}_{\Nyg} F^{\ast} \Prism_{R/A} \rightarrow I^{\bullet} \Prism_{R/A}$
admits an essentially unique extension to all animated commutative $\overline{A}$-algebras $R$. Since the functor
$R \mapsto F^{\ast} \Prism_{R/A}$ is a left Kan extension of its restriction to $\CAlg^{\QSyn}_{\overline{A}}$, we can identify
$\Fil^{0}_{\Nyg} F^{\ast} \Prism_{R/A}$ with the Frobenius pullback $F^{\ast} \Prism_{R/A}$ and
$\Fil^{0}(\varphi)$ with the relative Frobenius map $\varphi$. This completes the proof of the existence assertion of Proposition \ref{proposition:relative-nygaard-characterization}.

We now prove uniqueness. Suppose we are given another filtration $\Fil'^{\bullet} F^{\ast} \Prism_{R/A}$, depending functorially on $R$, together with
an extension of the relative Frobenius to a natural transformation $\Fil'^{\bullet}( \varphi): \Fil'^{\bullet} F^{\ast} \Prism_{R/A} \rightarrow I^{\bullet} \Prism_{R/A}$ 
satisfying conditions $(1)$ and $(2)$. We first show that the identity map  $F^{\ast} \Prism_{R/A}$ admits an essentially unique promotion to a morphism of
filtered complexes $\Fil^{\bullet}(\alpha_{R}): \Fil'^{\bullet} F^{\ast} \Prism_{R/A} \rightarrow \Fil^{\bullet}_{\Nyg} F^{\ast} \Prism_{R/A}$ for which the diagram
$$ \xymatrix@R=50pt@C=50pt{ \Fil'^{\bullet} F^{\ast} \Prism_{R/A} \ar[rr]^{ \Fil^{\bullet}(\alpha_R)} \ar[dr]^{ \Fil^{\bullet}( \varphi' ) } & & \Fil^{\bullet}_{\Nyg} F^{\ast} \Prism_{R/A} \ar[dl]_{ \Fil^{\bullet}( \varphi )} \\
& I^{\bullet} \Prism_{R/A} & }$$
commutes (up to a homotopy depending functorially on $R$). Since the functor $R \mapsto \Fil'^{\bullet} F^{\ast} \Prism_{R/A}$ satisfies condition $(1)$,
it is a left Kan extension of its restriction to $\CAlg_{\overline{A}}^{\QSyn}$. As noted above, the functors
$$ \CAlg^{\QSyn}_{\overline{A}} \rightarrow \DFiltI(\overline{A}) \quad \quad R \mapsto \Fil^{\bullet}_{\Nyg} F^{\ast} \Prism_{R/A}, I^{\bullet} \Prism_{R/A}$$
are right Kan extensions of their restrictions to $\calC$. Consequently, to construct the natural transformation $R \mapsto \Fil^{\bullet}( \alpha_{R} )$, we
may restrict our attention to $\overline{A}$-algebras $R$ which belong to the category $\calC$. In this case, we will prove the following assertion (for each integer $n \geq 0$):
\begin{itemize}
\item[$(\ast_n)$] For $R \in \calC$, the complex $\Fil'^{n} F^{\ast} \Prism_{R/A}$ is an abelian group (regarded as a chain complex concentrated in cohomological degree zero).
Moreover, the tautological map $\Fil'^{n} F^{\ast} \Prism_{R/A} \rightarrow F^{\ast} \Prism_{R/A}$ is a monomorphism whose image is the subgroup
$$ \Fil^{n}_{\Nyg} F^{\ast} \Prism_{R/A} = \{ x \in F^{\ast} \Prism_{R/A}: \varphi(x) \in I^{n} \Prism_{R/A} \}$$
defined above.
\end{itemize}
To prove $(\ast_n)$, we proceed by induction on $n$. The case $n=0$ is trivial. To handle the inductive step, we apply Remark \ref{remark:polynomial-to-smooth}
to the functor $R \mapsto \Fil'^{\bullet} F^{\ast} \Prism_{R/A}$ to obtain a fiber sequence
$$ \Fil'^{n+1} F^{\ast} \Prism_{R/A} \rightarrow \Fil'^{n} F^{\ast} \Prism_{R/A} \xrightarrow{ \gr^{n}(\varphi) } \Fil_{n}^{\conj} \overline{\Prism}_{R/A}\{n\}.$$
Using $(\ast_n)$, we can identify $\Fil'^{n} F^{\ast} \Prism_{R/A}$ with the ideal $\Fil^{n} F^{\ast} \Prism_{R/A}$. It follows that $\gr^{n}(\varphi)$ is
a surjection of abelian groups (regarded as cochain complexes concentrated in cohomological degree zero). Consequently,
we can identify $\Fil'^{n+1} F^{\ast} \Prism_{R/A}$ with the abelian group given by the kernel of the map
$$ \Fil^{n}_{\Nyg} F^{\ast} \Prism_{R/A} \xrightarrow{ \gr^{n}( \varphi) } \Fil_{n}^{\conj} \overline{\Prism}_{R/A}\{n\} \hookrightarrow \overline{\Prism}_{R/A}\{n\},$$
which immediately implies $(\ast_{n+1} )$. 

To complete the proof of Proposition \ref{proposition:relative-nygaard-characterization}, it will suffice to show that for every integer $n \geq 0$
and every animated commutative $\overline{A}$-algebra $R$, the map $\Fil^{n}(\alpha_R): \Fil'^{n} F^{\ast} \Prism_{R/A} \rightarrow \Fil^{n}_{\Nyg} F^{\ast} \Prism_{R/A}$
is an isomorphism. Since the construction $R \mapsto \Fil^{n}(\alpha_R)$ commutes with sifted colimits, we may assume without loss of generality
that $R$ is a finitely generated polynomial algebra over $\overline{A}$. Proceeding by induction on $n$, we are reduced to showing that the associated graded map
$\gr^{n}( \alpha_R): \gr'^{n} F^{\ast} \Prism_{R/A} \rightarrow \gr^{n}_{\Nyg} F^{\ast} \Prism_{R/A}$ induces an isomorphism on cohomology.
This follows from the commutativity of the diagram
$$ \xymatrix@R=50pt@C=50pt{ \gr'^{n} F^{\ast} \Prism_{R/A} \ar[rr]^{ \gr^{n}(\alpha_R)} \ar[dr]^{ \gr^{n}( \varphi' ) } & & \gr^{n}_{\Nyg} F^{\ast} \Prism_{R/A} \ar[dl]_{ \gr^{n}( \varphi )} \\
& \overline{\Prism}_{R/A}\{n\}, & }$$
noting that the vertical maps are isomorphisms in cohomological degrees $\leq n$ (and that the complexes $\gr'^{n} F^{\ast} \Prism_{R/A}$
and $\gr^{n}_{\Nyg} F^{\ast} \Prism_{R/A}$ are concentrated in cohomological degrees $\leq n$), by virtue of assumption $(2)$.
\end{proof}

\begin{remark}[Multiplicativity]\label{remark:multiplicative-relative-Nygaard}
Let $R$ be an animated commutative $\overline{A}$-algebra, and let us regard $F^{\ast} \Prism_{R/A}$ as a commutative algebra object of
the $\infty$-category $\widehat{\calD}(A)$ (that is, as an $E_{\infty}$-algebra over $A$). The commutative algebra structure on $F^{\ast} \Prism_{R/A}$
admits an essentially unique refinement to a commutative algebra structure on the filtered complex $\Fil^{\bullet}_{\Nyg} F^{\ast} \Prism_{R/A}$ which depends functorially
on $R$; in particular, the multiplication on $F^{\ast} \Prism_{R/A}$ admits a refinement
$$ (\Fil^{m}_{\Nyg} F^{\ast} \Prism_{R/A}) \otimes^{L}_{A} (\Fil^{n}_{\Nyg} F^{\ast} \Prism_{R/A}) \rightarrow \Fil^{m+n}_{\Nyg} F^{\ast} \Prism_{R/A},$$
depending functorially on $m$ and $n$. To prove this, we can use Kan extension arguments to reduce to the case where $R$ belongs
to the category $\calC$ appearing in the proof of Proposition \ref{proposition:relative-nygaard-characterization}. The desired assertion
then reduces to the observation that if $x,y \in F^{\ast} \Prism_{R/A}$ satisfy $\varphi(x) \in I^{m} \Prism_{R/A}$ and $\varphi(y) \in I^{n} \Prism_{R/A}$,
then $\varphi(xy) = \varphi(x) \varphi(y)$ belongs to $I^{m+n} \Prism_{R/A}$.
\end{remark}

\begin{remark}\label{remark:consequence-of-multiplicativity}
Let $R$ be an animated commutative $\overline{A}$-algebra. Combining Remark \ref{remark:multiplicative-relative-Nygaard} with Example \ref{example:relative-Nygaard-trivial},
we see that the multiplication on the filtered complex $\Fil^{\bullet}_{\Nyg} F^{\ast} \Prism_{R/A}$ induces a map of filtered complexes
$$ v: I \otimes^{L}_{A} \Fil^{\bullet}_{\Nyg} F^{\ast} \Prism_{R/A} \rightarrow \Fil^{\bullet+1}_{\Nyg} F^{\ast} \Prism_{R/A}.$$
Moreover, the composition of $v$ with the tautological map $\Fil^{\bullet+1}_{\Nyg} F^{\ast} \Prism_{R/A} \rightarrow \Fil^{\bullet}_{\Nyg} F^{\ast} \Prism_{R/A}$
agrees with the map $I \otimes^{L}_{A} \Fil^{\bullet}_{\Nyg} F^{\ast} \Prism_{R/A} \rightarrow \Fil^{\bullet}_{\Nyg} F^{\ast} \Prism_{R/A}$ induces by the inclusion $I \hookrightarrow A$.
\end{remark}

\begin{remark}[Change of Prism]\label{remark:Nygaard-change-of-prism}
Let $f: (A,I) \rightarrow (B,J)$ be a morphism of bounded prisms and let $R$ be an animated commutative algebra over the quotient ring $\overline{A} = A/I$.
Then there is a canonical isomorphism
$$ B \widehat{\otimes}^{L}_{A} \Fil^{\bullet}_{\Nyg} F^{\ast} \Prism_{R/A} \simeq \Fil^{\bullet}_{\Nyg} F^{\ast} \Prism_{(B \otimes^{L}_{A} R) / B }$$
in the $\infty$-category $\DFiltI(B)$ of filtered complexes over $B$, which is characterized (up to homotopy) by the fact that it depends functorially on $R$
and is given (in filtration degree zero) by the Frobenius pullback of the isomorphism $B \widehat{\otimes}_{A} \Prism_{R/A} \simeq
\Prism_{(B \otimes_{A}^{L} R)/B}$ of Remark \ref{remark:change-of-prism}.
\end{remark}

\subsection{The Relative de Rham Comparison}\label{subsection:relative-dr-comparison}

Throughout this section, we fix a bounded prism $(A,I)$ with quotient ring $\overline{A} = A/I$.

\begin{construction}[The Bockstein Operator]\label{construction:Bockstein-relative}
Let $R$ be an animated commutative $\overline{A}$-algebra. For every integer $n$, the relative Nygaard filtration of Proposition \ref{proposition:relative-nygaard-characterization} determines a fiber sequence $$ \gr^{n+1}_{\Nyg} F^{\ast} \Prism_{R/A} \rightarrow \Fil^{n}_{\Nyg} F^{\ast} \Prism_{R/A} / \Fil^{n+2}_{\Nyg} F^{\ast} \Prism_{R/A}
\rightarrow \gr^{n}_{\Nyg} F^{\ast} \Prism_{R/A}.$$
Passing to cohomology, we obtain a boundary map
$$ \beta: \mathrm{H}^{n}( \gr^{n}_{\Nyg} F^{\ast} \Prism_{R/A} ) \rightarrow \mathrm{H}^{n+1}( \gr^{n+1}_{\Nyg} F^{\ast} \Prism_{R/A} ),$$
which we will refer to as the {\it Bockstein operator}.
\end{construction}

\begin{proposition}\label{proposition:ralphus}
Let $R$ be an animated commutative $\overline{A}$-algebra. Then there is a unique homomorphism of commutative differential graded algebras
$$ \theta: \bigoplus_{n} \mathrm{H}^{n}( \gr^{n}_{\Nyg} F^{\ast} \Prism_{R/A}, \beta ) \rightarrow  ( \widehat{\Omega}^{\ast}_{R/\overline{A}}, d) $$
with the property that, in degree zero, it coincides with the ring isomorphism 
$$\mathrm{H}^{0}( \gr^{0}(\varphi)): \mathrm{H}^{0}( \gr^{0}_{\Nyg} F^{\ast} \Prism_{R/A} ) \simeq
\mathrm{H}^{0}( \widehat{R} ) = \widehat{\Omega}^{0}_{R/\overline{A}}$$ supplied by Remark \ref{remark:Nygaard-degree-zero}.
Here $( \widehat{\Omega}^{\ast}_{R/\overline{A}}, d)$ denotes the $p$-complete de Rham complex of $R$ relative to $\overline{A}$ (Notation \ref{notation:classical-de-Rham}).
Moreover, $\theta$ is an isomorphism.
\end{proposition}

\begin{proof}
For each integer $n$, Remarks \ref{remark:polynomial-to-smooth} and \ref{remark:connectivity-of-conjugate} supply isomorphisms
$$  \mathrm{H}^{n}( \gr^{n}_{\Nyg} F^{\ast} \Prism_{R/A}) \xrightarrow{ \mathrm{H}^{n}( \gr^{n}(\varphi) )}
\mathrm{H}^{n}( \Fil_{n}^{\conj} \overline{\Prism}_{R/A}\{n\} ) \simeq \widehat{\Omega}^{n}_{R/\overline{A}},$$
depending functorially on $R$. We claim that these isomorphisms identify the Bockstein operator of Construction \ref{construction:Bockstein-relative} with the de Rham differential on
$\widehat{\Omega}^{\ast}_{R/\overline{A}}$. To prove this, we can assume without loss of generality that $R$ is a polynomial algebra over $\overline{A}$, in which case it follows from Remark \ref{remark:HT-de-Rham}. This proves the existence of the isomorphism
$$\theta: \bigoplus_{n} \mathrm{H}^{n}( \gr^{n}_{\Nyg} F^{\ast} \Prism_{R/A}, \beta ) \simeq  ( \widehat{\Omega}^{\ast}_{R/\overline{A}}, d).$$
To prove uniqueness, let $\theta':  ( \bigoplus_{n} \mathrm{H}^{n}( \gr^{n}_{\Nyg} F^{\ast} \Prism_{R/A}, \beta ) \simeq  ( \widehat{\Omega}^{\ast}_{R/\overline{A}}, d)$
be any other homomorphism of commutative differential graded algebras which coincides with $\theta$ in degree zero. Then $\theta' \circ \theta^{-1}$
is an endomorphism of the $p$-complete de Rham complex $( \widehat{\Omega}^{\ast}_{R/\overline{A}}, d)$ which is the identity in degree zero.
We wish to show that $\theta' \circ \theta^{-1}$ is the identity. For this, it will suffice to show that $\theta' \circ \theta^{-1}$ restricts to the identity on
the usual de Rham complex $( \Omega^{\ast}_{R/\overline{A}}, d)$. This is clear, since $ \Omega^{\ast}_{R/\overline{A}}$ is generated (as a differential graded algebra)
by elements of degree zero.
\end{proof}

If $R$ is an animated commutative $\overline{A}$-algebra, we write $\widehat{\dR}_{R/ \overline{A}}$ for the $p$-complete derived de Rham complex of
$R$ relative to $\overline{A}$ (Construction \ref{construction:derived-de-Rham}), which we equip with the Hodge filtration $\Fil^{\bullet}_{\Hodge} \widehat{\dR}_{ R / \overline{A} }$.

\begin{proposition}\label{proposition:Nygaard-to-de-Rham-relative}
Let $R$ be an animated commutative $\overline{A}$-algebra. Then there is a canonical map
$$ \Fil^{\bullet}(\gamma_R): \Fil^{\bullet}_{\Nyg} F^{\ast} \Prism_{R/A} \rightarrow \Fil^{\bullet}_{\Hodge} \widehat{\dR}_{ R / \overline{A} }$$
of commutative algebra objects of $\DFiltI(A)$ which is determined (up to essentially unique homotopy) by the requirement that it depends functorially on $R$ and that 
$\gr^{0}( \gamma_{R} )$ agrees with the isomorphism
$$ \gr^{0}_{\Nyg} F^{\ast} \Prism_{R/A} \xrightarrow{ \gr^{0}(\varphi)} \Fil_{0}^{\conj} \overline{\Prism}_{R/A} \simeq \widehat{R} \simeq \gr^{0}_{\Hodge} \widehat{\dR}_{R/\overline{A}}$$
(up to a homotopy which depends functorially on $R$).
\end{proposition}

\begin{proof}
Since the functor $R \mapsto \Fil^{\bullet}_{\Nyg} F^{\ast} \Prism_{R/A}$ commutes with sifted colimits, it is a left Kan extension of its restriction to $\Poly_{\overline{A}}$.
It will therefore suffice to construct the morphism $\Fil^{\bullet}( \gamma_{R} )$ (and the homotopy $\gr^{0}( \gamma_{R} ) \simeq \gr^0( \varphi )$)
in the case where $R$ is a finitely generated polynomial algebra over $\overline{A}$. In this case, we can identify 
$\Fil^{\bullet}_{\Hodge} \widehat{\dR}_{ R / \overline{A} }$ with the $p$-complete de Rham complex $( \widehat{\Omega}^{\geq \bullet}_{R/\overline{A}}, d)$,
regarded as an object of the heart of the Beilinson t-structure on $\DFilt(A)$ (see Definition~\ref{definition:Beilinson-t} and Example~\ref{example:Beilinson-t-heart}). Since $\Fil^{\bullet}_{\Nyg} F^{\ast} \Prism_{R/A}$ is connective with respect to the Beilinson t-structure  (Remark \ref{remark:Nygaard-Beilinson}), we can identify morphisms $\Fil^{\bullet}(\gamma): \Fil^{\bullet}_{\Nyg} F^{\ast} \Prism_{R/A} \rightarrow
\Fil^{\bullet}_{\Hodge} \widehat{\dR}_{ R / \overline{A} }$ with morphisms of chain complexes
$$ ( \bigoplus_{n} \mathrm{H}^{n}( \gr^{n}_{\Nyg} F^{\ast} \Prism_{R/A}), \beta ) \rightarrow  ( \widehat{\Omega}^{\ast}_{R/\overline{A}}, d).$$
The desired result now follows from Proposition \ref{proposition:ralphus}.
\end{proof}

\begin{example}\label{example:comparison-map-dRW}
Let $(A,I)$ be the prism $(\Z_p, (p) )$ and let $R$ be a regular Noetherian $\F_p$-algebra. In the forthcoming Proposition \ref{proposition:concrete-Nygaard-regular}, we shall identify $\Fil^{\bullet}_{\Nyg} F^{\ast} \Prism_{R/A}$ with the Nygaard-filtered de Rham Witt complex $\calN^{\bullet} W\Omega^{\ast}_{R}$. Under this identification,
the comparison map $\Fil^{\bullet}(\gamma_R)$ of Proposition \ref{proposition:Nygaard-to-de-Rham-relative} has a concrete model at the level of (filtered) cochain complexes,
given by the quotient map $W \Omega^{\ast}_{R} \twoheadrightarrow \Omega^{\ast}_{R}$ which annihilates the Verschiebung operator on $W \Omega^{\ast}_{R}$.
\end{example}

Let $R$ be an animated commutative $\overline{A}$-algebra. Restricting the map of Proposition \ref{proposition:Nygaard-to-de-Rham-relative} to filtration degree zero, we obtain a morphism
$\gamma_{R}: F^{\ast} \Prism_{R/A} \rightarrow \widehat{\dR}_{R/ \overline{A} }$ in the $\infty$-category $\widehat{\calD}(A)$. Extending scalars along the quotient map $A \twoheadrightarrow \overline{A}$, we obtain a map $\overline{\gamma}_{R}: \overline{A} \otimes_{A}^{L} F^{\ast} \Prism_{R/A} \rightarrow \widehat{\dR}_{R/ \overline{A}}$.
The following result is Corollary~15.4 of \cite{prisms}:

\begin{proposition}[de Rham Comparison for Relative Prismatic Cohomology]\label{proposition:de-Rham-comparison-relative}
Let $R$ be an animated commutative $\overline{A}$-algebra. Then the map
$$ \overline{\gamma}_{R}: \overline{A} \otimes^{L}_{A} F^{\ast} \Prism_{R/A} \rightarrow \widehat{\dR}_{R/ \overline{A}}$$
is an isomorphism in the $\infty$-category $\widehat{\calD}( \overline{A} )$.
\end{proposition}

\begin{proof}
Since the construction $R \mapsto \overline{\gamma}_{R}$ commutes with sifted colimits, it will suffice to prove Proposition \ref{proposition:de-Rham-comparison-relative}
in the special case where $R = \overline{A}[x_1, \cdots, x_d]$ is a finitely generated polynomial algebra over $\overline{A}$.
It follows from Remark \ref{remark:Nygaard-change-of-prism} that the construction $$(A,I) \mapsto \fib( \overline{\gamma}_{\overline{A}[x_1, \ldots, x_d]} ) \in \widehat{\calD}( \overline{A} )$$
is compatible with base change, and therefore defines a quasi-coherent complex $\mathscr{E}$ on the Hodge-Tate divisor $\WCart^{\mathrm{HT}}$ (see Remark \ref{remark:flashcard}).
We wish to show that $\mathscr{E}$ vanishes. To prove this, it will suffice to show that $\mathscr{E}$ vanishes after pullback along the map
$\Spec(\F_p) \rightarrow \WCart^{\mathrm{HT}}$ obtained by applying Remark \ref{remark:HT-point-of-prismatic-stack} to the perfect prism $( \Z_p, (p) )$. In other words,
it suffices to prove Proposition \ref{proposition:de-Rham-comparison-relative} after replacing $(A,I)$ by the prism $(\Z_p, (p) )$. Since $R = \F_p[ x_1, \ldots, x_d]$
is a regular Noetherian $\F_p$-algebra, Example \ref{example:comparison-map-dRW} supplies an explicit model for $\overline{\gamma}_{R}$ at the level of cochain complexes, given by the quotient map $W\Omega^{\ast}_{R} / p W \Omega^{\ast}_{R} \twoheadrightarrow \Omega^{\ast}_{R}$. This map is an isomorphism at the level of cohomology (see Remark~4.3.6 of
\cite{DRW}).
\end{proof}

\begin{example}\label{example:qdR-vs-dR}
Let $\Z_p[[ \slashp]]$ be the $\delta$-ring of Notation \ref{notation:reduced-q-de-Rham-prism}. Note that the composite map
$$ \Z_p[[ \slashp ]] \xrightarrow{\varphi} \Z_p[[\slashp]] \twoheadrightarrow \Z_p[[ \slashp ]] / ( \slashp ) \simeq \Z_p$$
is a surjection which carries $\slashp$ to the element $p \in \Z_p$.

Let $R$ be an animated commutative $\Z_p$-algebra and let $\slashOmega_{R} = \Prism_{ R/ \Z_p[[ \slashp ]] }$ be the $\slashp$-de Rham complex 
of Construction \ref{construction:tilde-p-deRham}. Applying Proposition \ref{proposition:de-Rham-comparison-relative} to the
prism $( \Z_p[[\slashp]], (\slashp) )$, we obtain a canonical isomorphism
$$ \Z_p[[\slashp]]/(\slashp-p) \otimes^{L}_{ \Z_p[[\slashp]] } \slashOmega_{R} \simeq \widehat{\dR}_{R}.$$
of commutative algebra objects of $\widehat{\calD}(\Z_p)$. Equivalently, we have a canonical isomorphism
$$ \Z_p[[q-1]]/(q-1) \otimes^{L}_{ \Z_p[[q-1]] } \qOmega_{R} \simeq \widehat{\dR}_{R},$$
where $\qOmega_{R}$ is the $q$-de Rham complex of Definition \ref{definition:q-de-Rham-cohomology}
(see Remark \ref{remark:q-de-Rham-tilde}).

Stated more informally: the $\slashp$-de Rham complex $\slashOmega_{R}$ specializes to the usual $p$-complete derived de Rham complex $\widehat{\dR}_{R}$
by setting $\widetilde{p}$ equal to $p$, and the $q$-de Rham complex $\qOmega_{R}$ specializes to $\widehat{\dR}_{R}$ by setting $q$ equal to $1$.
\end{example}

\begin{remark}\label{remark:diffracted-Hodge-cohomology-mod-p}
Let $R$ be an animated commutative ring and let $\slashOmega_{R}$ be the $\slashp$-de Rham complex of $R$, which we regard
as a $(p, \slashp)$-complete object of the derived $\infty$-category $\calD( \Z_p[[\slashp]] )$. Then:
\begin{itemize}
\item The $p$-complete diffracted Hodge complex $\widehat{\Omega}^{\DHod}_{R}$ can be obtained from $\slashOmega_{R}$
by extending scalars along the map $\Z_p[[ \slashp ]] \twoheadrightarrow \Z_p$ carrying $\slashp$ to $0$ (Remark \ref{remark:diffracted-Hodge-extended-scalars}).
\item The $p$-complete derived de Rham complex $\widehat{\Omega}_{R}$ can be obtained from $\slashOmega_{R}$ by extending scalars along the map
$\Z_p[[ \slashp ]] \twoheadrightarrow \Z_p$ carrying $\slashp$ to $p$ (Example \ref{example:qdR-vs-dR}).
\end{itemize}
We therefore obtain a canonical isomorphism
$$ \F_p \otimes^{L} \widehat{\Omega}^{\DHod}_{R} \simeq \F_p \otimes^{L} \widehat{\dR}_{R} \simeq \dR_{ \overline{R} / \F_p },$$
where $\overline{R}$ denotes the derived tensor product $\F_p \otimes^{L} R$.
\end{remark}

Proposition \ref{proposition:de-Rham-comparison-relative} admits a refinement, which describes the Hodge filtration on derived de Rham cohomology:

\begin{corollary}\label{corollary:relative-Nygaard-deRham-refined}
For every animated commutative $\overline{A}$-algebra $R$, there is a fiber sequence
$$ I \otimes^{L}_{A} \Fil^{\bullet-1}_{\Nyg} F^{\ast} \Prism_{R/A} \xrightarrow{v} 
\Fil^{\bullet}_{\Nyg} F^{\ast} \Prism_{R/A} \xrightarrow{ \Fil^{\bullet}(\gamma_R) } \Fil^{\bullet}_{\Hodge} \widehat{\dR}_{ R/ \overline{A}}$$
in the $\infty$-category $\DFiltI(A)$, where $v$ is the morphism of Remark \ref{remark:consequence-of-multiplicativity} and $\Fil^{\bullet}(\gamma_R)$ is the morphism of 
Proposition \ref{proposition:Nygaard-to-de-Rham-relative}. This fiber sequence determined (up to isomorphism) by the requirement that it depends functorially on $R$.
\end{corollary}

\begin{proof}
We first claim that the composite map
$$ (\Fil^{\bullet}(\gamma_R) \circ v): I \otimes^{L}_{A} \Fil^{\bullet-1}_{\Nyg} F^{\ast} \Prism_{R/A} \rightarrow  \Fil^{\bullet}_{\Hodge} \widehat{\dR}_{ R/ \overline{A}}$$
admits an essentially unique nullhomotopy which depends functorially on $R$. Since the functor $R \mapsto I \otimes^{L}_{A} \Fil^{\bullet-1}_{\Nyg} F^{\ast} \Prism_{R/A}$
commutes with sifted colimits, it will suffice to construct this nullhomotopy in the case where $R$ is a finitely generated polynomial algebra over $\overline{A}$.
We now observe that the mapping space $\Hom_{ \DFiltI(A) }( I \otimes^{L}_{A} \Fil^{\bullet-1}_{\Nyg} F^{\ast} \Prism_{R/A}, \Fil^{\bullet}_{\Hodge} \widehat{\dR}_{ R'/ \overline{A}} )$
is contractible for $R, R' \in \Poly_{\overline{A}}$, since $\Fil^{\bullet}_{\Hodge} \widehat{\dR}_{ R'/ \overline{A}} \simeq ( \widehat{\Omega}^{\geq \bullet}_{R'/\overline{A}}, d)$ lies in the heart of the Beilinson t-structure on $\DFiltI(A)$, while the complex $I \otimes^{L}_{A} \Fil^{\bullet-1}_{\Nyg} F^{\ast} \Prism_{R/A}$ belongs to $\DFiltI(A)^{\leq -1}$
(Remark \ref{remark:Nygaard-Beilinson}).

Let us denote the cofiber of $v$ by $\Fil^{\bullet}_{\Nyg} F^{\ast} \Prism_{R/A} / I \Fil^{ \bullet-1}_{\Nyg} F^{\ast} \Prism_{R/A}$, so that
$\Fil^{\bullet}_{\gamma}$ descends to a morphism of filtered complexes
$$ \Fil^{\bullet}( \overline{\gamma}_{R} ): \Fil^{\bullet}_{\Nyg} F^{\ast} \Prism_{R/A} / I \Fil^{ \bullet-1}_{\Nyg} F^{\ast} \Prism_{R/A} \rightarrow \Fil^{\bullet}_{\Hodge} \widehat{\dR}_{ R/ \overline{A} },$$
depending functorially on $R$. We wish to show that, for every integer $n$, the induced map
$$ \Fil^{n}( \overline{\gamma}_{R} ): \Fil^{n}_{\Nyg} F^{\ast} \Prism_{R/A} / I \Fil^{n-1}_{\Nyg} F^{\ast} \Prism_{R/A} \rightarrow
\Fil^{n}_{\Hodge} \widehat{\dR}_{R/\overline{A}}$$
is an isomorphism. For $n \leq 0$, this is a restatement of Proposition \ref{proposition:de-Rham-comparison-relative}. The general case follows by induction on $n$, using the observation that the associated graded map
$$ \gr^{n}( \overline{\gamma}_{R} ): \gr^{n}_{\Nyg} F^{\ast} \Prism_{R/A} / I \gr^{n-1}_{\Nyg} F^{\ast} \Prism_{R/A} \rightarrow
\gr^{n}_{\Hodge} \widehat{\dR}_{R/\overline{A}}$$
is obtained by compositing the relative Frobenius isomorphism
$$ \varphi: \gr^{n}_{\Nyg} F^{\ast} \Prism_{R/A} / I \gr^{n-1}_{\Nyg} F^{\ast} \Prism_{R/A} \simeq \Fil_{n}^{\conj} \overline{\Prism}_{R/A}\{n\} / \Fil_{n-1}^{\conj} \overline{\Prism}_{R/A}\{n\}$$
of Remark \ref{remark:polynomial-to-smooth} with the Hodge-Tate comparsion isomorphism $\gr_{n}^{\conj} \overline{\Prism}_{R/A}\{n\} \simeq L\widehat{\Omega}^{n}_{R/\overline{A}}[n]$
of Remark \ref{remark:derived-HT-filtration}.
\end{proof}

\begin{corollary}
Let $R$ be an animated commutative $\overline{A}$-algebra. Then the comparison map of Proposition \ref{proposition:Nygaard-to-de-Rham-relative} induces an isomorphism
$$ \overline{A} \otimes^{L}_{A} (\varprojlim_{m} \Fil^{m}_{\Nyg} F^{\ast} \Prism_{R/A}) \rightarrow \varprojlim_{m} \Fil^{m}_{\Hodge} \widehat{\dR}_{R / \overline{A} }.$$
\end{corollary}

\begin{corollary}\label{corollary:Nygaard-vs-Hodge-complete}
Let $R$ be an animated commutative $\overline{A}$-algebra. Then the complex $F^{\ast} \Prism_{R/A}$ is Nygaard-complete
(that is, the limit $\varprojlim_{m} \Fil^{m}_{\Nyg} F^{\ast} \Prism_{R/A}$ vanishes) if and only if the derived de Rham complex $\widehat{\dR}_{R/ \overline{A} }$ is
Hodge-complete (that is, the limit $\varprojlim_{m} \Fil^{m}_{\Hodge} \widehat{\dR}_{R/ \overline{A}}$ vanishes).
\end{corollary}

\begin{corollary}\label{corollary:smooth-Nygaard-complete}
Let $R$ be a $p$-completely flat $\overline{A}$-algebra for which the quotient $R/pR$ can be written as a filtered colimit of smooth $\overline{A} / p \overline{A}$-algebras.
Then the complex $F^{\ast} \Prism_{R/A}$ is Nygaard-complete.
\end{corollary}

\begin{proof}
Proposition \ref{proposition:derived-to-classical-de-Rham} guarantees that the comparison map 
$$\Fil^{\bullet}_{\Hodge} \widehat{\dR}_{R/\overline{A}}\rightarrow ( \widehat{\Omega}^{\geq \bullet}_{R/ \overline{A}} , d)$$
is an isomorphism in the filtered derived $\infty$-category $\DFiltI(\overline{A})$, so that the derived de Rham complex $\widehat{\dR}_{R/\overline{A} }$ is Hodge-complete.
Applying Corollary \ref{corollary:Nygaard-vs-Hodge-complete}, we conclude that $F^{\ast} \Prism_{R/A}$ is Nygaard-complete.
\end{proof}

\begin{example}
Let $R$ be a regular Noetherian $\F_p$-algebra. It follows from Corollary \ref{corollary:smooth-Nygaard-complete} that the prismatic complex
$\Prism_{R} \simeq F^{\ast} \Prism_{R/\Z_p}$ is complete with respect to its Nygaard filtration. This can also be deduced from the forthcoming Proposition \ref{proposition:concrete-Nygaard-regular}: the Nygaard filtration on the prismatic complex $F^{\ast} \Prism_{R/\Z_p}$ has a concrete
cochain-level model given by the Nygaard filtration on the de Rham-Witt complex $W\Omega_{R}^{\ast}$, where we can witness completeness
at the level of cochain complexes: $W\Omega_{R}^{\ast}$ is the inverse limit of the system of quotient complexes $\{ W \Omega_{R}^{\ast} / \calN^{m} W \Omega_{R}^{\ast} \}_{m \geq 0}$
(see Remark~8.1.2 of \cite{DRW}). 
\end{example}

\begin{corollary}\label{corollary:Nygaard-completeness-regular}
Let $R$ be a regular Noetherian $\F_p$-algebra. Then relative Nygaard filtration $\Fil^{\bullet}_{\Nyg} F^{\ast} \Prism_{R/\Z_p}$ is complete:
that is, the limit $\varprojlim_{m} \Fil^{m}_{\Nyg} F^{\ast} \Prism_{R/\Z_p}$ vanishes in the $\infty$-category $\widehat{\calD}(\Z_p)$.
\end{corollary}

\begin{corollary}\label{corollary:Beilinson-connective-cover}
Let $R$ be a $p$-completely smooth $\overline{A}$-algebra. Then the relative Frobenius map $\Fil^{\bullet}(\varphi): \Fil^{\bullet}_{\Nyg} F^{\ast} \Prism_{R/A} \rightarrow I^{\bullet} \Prism_{R/A}$
exhibits $\Fil^{\bullet}_{\Nyg} F^{\ast} \Prism_{R/A}$ as the connective cover of $I^{\bullet} \Prism_{R/A}$ with respect to the
Beilinson t-structure on $\DFiltI(A)$.
\end{corollary}

\begin{proof}
Corollary \ref{corollary:smooth-Nygaard-complete} guarantees that the complex $F^{\ast} \Prism_{R/A}$ is complete with respect to its Nygaard filtration.
Using the criterion of Remark \ref{remark:characterize-Beilinson-connective}, we are reduced to showing that relative Frobenius map
$\gr^{n}( \varphi): \gr^{n}_{\Nyg}F^*  \Prism_{R/A} \rightarrow \overline{\Prism}_{R/A}\{n\}$ exhibits $\gr^{n}_{\Nyg} F^* \Prism_{R/A}$ as the truncation of $\overline{\Prism}_{R/A}\{n\}$ in cohomological degrees $\leq n$, for every integer $n$. This follows from Remark \ref{remark:polynomial-to-smooth}.
\end{proof}

\begin{corollary}\label{respit}
Let $R$ be an $\overline{A}$-algebra for which the $p$-completion $\widehat{R}$ is $p$-completely smooth of relative dimension $\leq d$ over $\overline{A}$.
Then the map $v: I \otimes^{L}_{A} \Fil^{n}_{\Nyg} F^{\ast} \Prism_{R/A} \rightarrow \Fil^{n+1}_{\Nyg} F^{\ast} \Prism_{R/A}$ of Remark \ref{remark:consequence-of-multiplicativity}
is an isomorphism for $n \geq d$.
\end{corollary}

\begin{proof}
By virtue of Corollary \ref{corollary:relative-Nygaard-deRham-refined}, it suffices to observe that the complex 
$$\Fil^{n}_{\Hodge} \widehat{\dR}_{R/ \overline{A}} \simeq ( \widehat{\Omega}^{\geq n}_{R /\overline{A}}, d)$$
vanishes for $n > d$. 
\end{proof}

\begin{corollary}\label{corollary:increasing-filtration-on-prism}
Let $R$ be an animated commutative $\overline{A}$-algebra. Then the relative Frobenius map exhibits $\Prism_{R/A}$ as the $(p,I)$-completed direct limit
of the diagram
$$ F^{\ast} \Prism_{R/A} \rightarrow I^{-1} \Fil^{1}_{\Nyg} F^{\ast} \Prism_{R/A} \rightarrow I^{-2} \Fil^{2}_{\Nyg} F^{\ast} \Prism_{R/A} \rightarrow \cdots$$
\end{corollary}

\begin{proof}
Let $T(R)$ denote the colimit $\varinjlim_{m} I^{-m} \Fil^{m}_{\Nyg} F^{\ast} \Prism_{R/A}$, formed in the $(p,I)$-complete derived $\infty$-category
$\widehat{\calD}(A)$. The relative Frobenius supplies a comparison map $u_{R}: T(R) \rightarrow \Prism_{R/A}$, depending functorially on $R$.
We wish to show that $u_R$ is an isomorphism for every animated commutative $\overline{A}$-algebra $R$. Since the construction $R \mapsto u_R$
commutes with filtered colimits, we may assume without loss of generality that $R = \overline{A}[x_1, \cdots, x_d]$ is a finitely generated
polynomial algebra over $\overline{A}$. In this case, Corollary \ref{corollary:relative-Nygaard-deRham-refined} guarantees that the diagram
$\{ I^{-m} \Fil^{m}_{\Nyg} F^{\ast} \Prism_{R/A} \}_{m \geq 0}$ is constant for $m \geq d$. Tensoring both sides with $I^{d}$, we are reduced to showing
that the relative Frobenius map
$$ \Fil^{d}( \varphi ): \Fil^{d}_{\Nyg} F^{\ast} \Prism_{R/A} \rightarrow I^{d} \Prism_{R/A}$$
is an isomorphism. To check this, we can extend scalars along the map $A \twoheadrightarrow \overline{A}$. Using Corollary \ref{respit}, we are
reduced to showing that the map $\gr^{d}( \varphi): \gr^{d}_{\Nyg} F^{\ast} \Prism_{R/A} \rightarrow \overline{\Prism}_{R/A}\{d\}$ is an isomorphism.
By virtue of Remark \ref{remark:polynomial-to-smooth}, we can identify $\gr^{d}(\varphi)$ with the tautological map
$\Fil_{d}^{\conj} \overline{\Prism}_{R/A}\{d\} \rightarrow \overline{\Prism}_{R/A}\{d\}$. Since the conjugate filtration on $\overline{\Prism}_{R/A}$ is exhaustive,
we are reduced to showing that the successive quotients $\gr_{n}^{\conj} \overline{\Prism}_{R/A}$ vanish for $n > d$. This follows from
the Hodge-Tate comparison isomorphism $\gr_{n}^{\conj} \overline{\Prism}_{R/A} \simeq L \widehat{\Omega}^{n}_{R/\overline{A}}\{-n\}[-n]$ of Remark
\ref{remark:derived-HT-filtration}, since $L \widehat{\Omega}^{1}_{R/\overline{A}}$ is a free $\widehat{R}$-module of rank $d < n$.
\end{proof}

\subsection{The Crystalline Nygaard Filtration}\label{subsection:Nygaard-filtration-comparison}

Let $R$ be an $\F_p$-algebra and let $\Prism_{R/ \Z_p}$ denote the prismatic complex of $R$ relative to the perfect prism
$( \Z_p, (p) )$. Note that the Frobenius is the identity on the $\delta$-ring $\Z_p$, so we can identify $\Prism_{R/\Z_p}$
with its Frobenius pullback $F^{\ast} \Prism_{R/\Z_p}$. We will write 
 $\Fil^{\bullet}_{\Nyg} \Prism_{R/\Z_p}$ for the
image of the relative Nygaard filtration $\Fil^{\bullet}_{\Nyg} F^{\ast} \Prism_{R/\Z_p}$ under this identification; we shall write 
\[\gamma_R:\Prism_{R/\Z_p} \to \dR_{R/\F_p}\]
for the map coming from Proposition~\ref{proposition:Nygaard-to-de-Rham-relative}. If the $\F_p$-algebra $R$ is $p$-quasisyntomic, then Theorem \ref{theorem:crystalline-comparison} supplies an isomorphism
$$ \gamma_{\Prism}^{\crys}: \Prism_{R/ \Z_p} \simeq \RGamma_{\crys}(R / \Z_p ).$$
of commutative algebra objects of $\widehat{\calD}(\Z_p)$. Our goal in this section is to describe the image of the
Nygaard filtration under this isomorphism in two special cases: 
when $R$ is a quasiregular semiperfect $\F_p$-algebra (Proposition \ref{proposition:concrete-Nygaard-qrsp}), and when $R$ is a regular Noetherian $\F_p$-algebra (Proposition \ref{proposition:concrete-Nygaard-regular}). We also show that
the image of $\Fil^{1}_{\Nyg} \Prism_{R/\Z_p}$ can always be characterized as the fiber of the crystalline
augmentation map $\epsilon_{\crys}: \RGamma_{\crys}(R/\Z_p) \rightarrow R$. This is an immediate consequence
of the following compatibility:

\begin{proposition}\label{proposition:augmentation-compatibility}
Let $R$ be an $\F_p$-algebra. Then the diagram
$$  \xymatrix@R=50pt@C=50pt{ \Prism_{R/\Z_p} \ar[d] \ar[r]^-{ \gamma_{\Prism}^{\crys} } & \RGamma_{\crys}(R/\Z_p) \ar[d]^{ \epsilon_{\crys} } \\
\dR_{R/\F_p} \ar[r]^-{ \epsilon_{\dR} } & R }$$
commutes (up to a homotopy which depends functorially on $R$). Here $\epsilon_{\crys}$
and $\epsilon_{\dR}$ denote the crystalline and de Rham augmentation maps (Notations \ref{notation:crystalline-augmentation}
and \ref{notation:dra}), and the right vertical map is the composition
$$ \Prism_{R/\Z_p} \xrightarrow{\sim} F^{\ast} \Prism_{R/\Z_p}
\rightarrow \F_p \otimes^{L} F^{\ast} \Prism_{R/\Z_p} \simeq \dR_{R/\F_p}.$$
\end{proposition}

\begin{proof}
The functor
$$ \CAlg_{\F_p}^{\anim} \rightarrow \widehat{\calD}(\Z_p) \quad \quad R \mapsto \Prism_{R/\Z_p}$$
commutes with sifted colimits, and is therefore a left Kan extension of its restriction to the category of finitely
generated polynomial algebras over $\F_p$. It will therefore suffice to prove Proposition \ref{proposition:augmentation-compatibility}
in the case where $R$ is $p$-quasisyntomic. Since the identity functor $R \mapsto R$ satisfies flat descent,
we can further reduce to the case where $R$ is quasiregular semiperfect. In this case, we can identify $\gamma_{\Prism}^{\crys}$ with the inverse of the ring homomorphism $\beta_{R}: A_{\crys}(R) \rightarrow \Prism_{R}$ appearing in
the statement of Lemma \ref{lemma:transformation-in-easy-case}. Consider the diagram of commutative rings
$$  \xymatrix@R=50pt@C=50pt{ \Prism_{R/\Z_p} \ar[d] \ar[r]^-{ \gamma_{\Prism}^{\crys} } & A_{\crys}(R) \ar[d]^{ \epsilon_{\crys} } \ar[r]^-{\beta_R} & \Prism_{R/\Z_p} \ar[d] \\
\dR_{R/\F_p} \ar[r]^-{ \epsilon_{\dR} } & R \ar@{^{(}->}[r] & \overline{\Prism}_{R/\Z_p}. }$$
Note that the right square commutes (by the construction of $\beta_{R}$), and that the bottom right horizontal map
$$R \simeq \Fil^{\conj}_{0} \overline{\Prism}_{R/\Z_p} \rightarrow \overline{\Prism}_{R/\Z_p} \simeq \overline{\Prism}_{R}$$
is a monomorphism. It will therefore suffice to prove the commutativity of the outer rectangle, which is implicit
in the construction of the left vertical map. 
\end{proof}

\begin{notation}\label{notation:Nygaard-crystalline}
Let $X$ be an $\F_p$-scheme. We write $\Fil^{1}_{\Nyg} \RGamma_{\crys}(X/\Z_p)$ denote the fiber of the crystalline augmentation morphism
$$\epsilon_{\crys}: \RGamma_{\crys}(X/\Z_p) \rightarrow \RGamma(X, \calO_X)$$ 
of Notation \ref{notation:crystalline-augmentation}. By virtue of Proposition \ref{proposition:augmentation-compatibility},
the comparison map 
$$ \gamma_{\Prism}^{\crys}: \RGamma_{\Prism}( X / \Z_p ) \rightarrow \RGamma_{\crys}( X / \Z_p )$$
of Remark \ref{remark:crystalline-comparison-general} restricts to a map
$$ \Fil^{1}( \gamma_{\Prism}^{\crys} ): \Fil^{1}_{\Nyg} \RGamma_{\Prism}(X/ \Z_p) \rightarrow \Fil^{1}_{\Nyg} \RGamma_{\crys}( X / \Z_p ),$$
which is an isomorphism if $X$ is $p$-quasisyntomic (Theorem \ref{theorem:crystalline-comparison}).
\end{notation}

\begin{remark}\label{remark:phi-over-p}
Let $R$ be a quasisyntomic $\F_p$-algebra and let $\Fil^{\bullet}(\varphi): \Fil^{\bullet}_{\Nyg} \Prism_{R/\Z_p} \rightarrow p^{\bullet} \Prism_{R/\Z_p}$ be the morphism of filtered complexes appearing in Proposition \ref{proposition:relative-nygaard-characterization}. Restricting to filtration degree $1$ and invoking the identifications
$$ \RGamma_{\Prism}(X/ \Z_p) \simeq \Fil^{1}_{\Nyg} \RGamma_{\crys}(X/\Z_p)
\quad \quad  \Fil^{1}_{\Nyg} \RGamma_{\Prism}(X/ \Z_p) \simeq \Fil^{1}_{\Nyg} \RGamma_{\crys}(X/\Z_p),$$
we obtain a {\it divided Frobenius morphism}
$$ \varphi/p: \Fil^{1}_{\Nyg} \RGamma_{\crys}(X/\Z_p) \rightarrow \RGamma_{\crys}(X/\Z_p),$$
which fits into a commutative diagram
$$ \xymatrix@R=50pt@C=50pt{ \Fil^{1}_{\Nyg} \RGamma_{\crys}(X/\Z_p) \ar[r]^-{\varphi/p} \ar[d] & \RGamma_{\crys}(X/\Z_p) \ar[d]^{p} \\
\RGamma_{\crys}(X/\Z_p) \ar[r]^-{\varphi} & \RGamma_{\crys}(X/\Z_p). }$$
Stated more informally, the Frobenius endomorphism of the crystalline cochain complex
$\RGamma_{\Prism}(X/\Z_p)$ is divisible by $p$ when restricted to the fiber 
$$\Fil^{1}_{\Nyg} \RGamma_{\crys}(X/\Z_p) = \fib( \RGamma_{\crys}(X/\Z_p) \rightarrow \RGamma(X, \calO_X) ).$$
\end{remark}

\begin{notation}
Let $R$ be a quasiregular semiperfect $\F_p$-algebra, and let $A_{\crys}(R)$ denote the commutative ring of Construction \ref{construction:Acrys}.
The Frobenius endomorphism $\varphi: R \rightarrow R$ induces an endomorphism of $A_{\crys}(R)$, which we will also denote by $\varphi$.
For each integer $n \geq 0$, we let $\Fil^{n}_{\Nyg} A_{\crys}(R) \subseteq A_{\crys}(R)$ consisting of those elements
$x \in A_{\crys}(R)$ for which $\varphi(x)$ is divisible by $p^n$. It follows from the proof of Proposition \ref{proposition:relative-nygaard-characterization} that we can identify
$\Fil^{n}_{\Nyg} A_{\crys}(R)$ with the image of $\Fil^{n}_{\Nyg} \Prism_{R/ \Z_p}$ under the crystalline comparison isomorphism
$\Prism_{R/\Z_p} \simeq \RGamma_{\crys}( R / \Z_p) \simeq A_{\crys}(R)$ of Theorem \ref{theorem:crystalline-comparison}.
\end{notation}

\begin{remark}\label{remark:short-exact-sequence-Nygaard-conjugate}
Let $R$ be a quasiregular semiperfect $\F_p$-algebra. For each integer $n \geq 0$, the construction $x \mapsto \varphi(x) / p^n$ determines
a map $\Fil^{n}_{\Nyg} A_{\crys}(R) \rightarrow A_{\crys}(R)$, which induces a monomorphism of abelian groups
$$ \Fil^{n}_{\Nyg} A_{\crys}(R) / \Fil^{n+1} A_{\crys}(R) \hookrightarrow A_{\crys}(R) / p A_{\crys}(R).$$
Combining Remark \ref{remark:polynomial-to-smooth} with Lemma \ref{lemma:conjugate-filtration}, we see that the image of this monomorphism
is the abelian group $\Fil_{n}^{\conj}(A_{\crys}(R) / p A_{\crys}(R))$ introduced in Construction \ref{construction:conjugate-filtration}.
That is, we have a short exact sequence of abelian groups
$$ 0 \rightarrow \Fil^{n+1}_{\Nyg} A_{\crys}(R) \rightarrow \Fil^{n}_{\Nyg} A_{\crys}(R) \xrightarrow{ \varphi / p^n } \Fil_{n}^{\conj}(A_{\crys}(R) / p A_{\crys}(R) ) \rightarrow 0$$
\end{remark}

\begin{proposition}\label{proposition:concrete-Nygaard-qrsp}
Let $R$ be a quasiregular semiperfect $\F_p$-algebra and let $I$ denote the kernel of the quotient map $R^{\flat} \twoheadrightarrow R$.
For each $n \geq 0$, $\Fil^{n}_{\Nyg} A_{\crys}(R)$ coincides with the ideal $J(n) \subseteq A_{\crys}(R)$ generated by elements of the form
$$ p^{m_0} \gamma_{m_1}([x_1]) \gamma_{m_2}( [x_2 ] ) \cdots \gamma_{m_k}( [x_k] ),$$
where $x_1, x_2, \cdots, x_k \in I$ and $m_0 + m_1 + \cdots + m_k \geq n$; here we write $[x_i]$ for the image of the Teichm\"{u}ller representative of
$x_i$ under the tautological map $W( R^{\flat} ) \rightarrow A_{\crys}(R)$.
\end{proposition}

\begin{proof}
To show that $J(n)$ is contained in $\Fil^{n}_{\Nyg} A_{\crys}(R)$, we observe that for every collection of elements $x_1, \cdots, x_k \in I$ and every collection of integers $m_0, m_1, \cdots, m_k \geq 0$ satisfying $m_0 + m_1 + \cdots + m_k \geq n$, the expression
$$ \varphi( p^{m_0} \gamma_{m_1}([x_1]) \cdots \gamma_{m_k}( [x_k] ) ) =
p^{m_0} \gamma_{m_1}( [x_1]^p )) \cdots \gamma_{m_k}( [ x_k]^p ) $$
is divisible by $p^{n}$ in $A_{\crys}(R)$. This follows from the observation that each $\gamma_{m_i}( [x_i]^{p} )$ is equal to
$\frac{ (pm_i)!}{m_i!} \gamma_{pm_i}( [x_i ] )$ as an element of $A_{\crys}(R)$, and is therefore divisible by $p^{m_i}$.

We now wish to show that the inclusion $J(n) \subseteq \Fil^{n}_{\Nyg} A_{\crys}(R)$ is an equality. We proceed by induction on $n$, the case $n=0$ being trivial.
To handle the inductive step, we observe that Remark \ref{remark:short-exact-sequence-Nygaard-conjugate} supplies a commutative diagram of short
exact sequences
$$ \xymatrix@R=50pt@C=30pt{ 0 \ar[r] & J(n+1) \ar[r] \ar[d] & J(n) \ar[r] \ar[d] & J(n) / J(n+1) \ar[d]^{\beta} \ar[r] & 0 \\
0 \ar[r] & \Fil^{n+1}_{\Nyg} A_{\crys}(R) \ar[r] & \Fil^{n}_{\Nyg} A_{\crys}(R) \ar[r]^-{ \varphi / p^{n} } & \Fil_{n}^{\conj} A_{\crys}(R) / p A_{\crys}(R) \ar[r] & 0. }$$
Combining the snake lemma with our inductive hypothesis, we are reduced to proving that the map $\beta: J(n) / J(n+1) \rightarrow \Fil_{n}^{\conj} A_{\crys}(R) / p A_{\crys}(R)$ is an isomorphism.

For $-1 \leq s \leq n$, let $J(n,s) \subseteq J(n)$ denote the ideal generated by $J(n+1)$ together with elements of the form
\begin{equation}\label{equation:define-y} y = p^{m_0} \gamma_{m_1}([x_1]) \gamma_{m_2}( [x_2 ] ) \cdots \gamma_{m_k}( [x_k] ), \end{equation}
where $x_1, x_2, \cdots, x_k \in I$, and $m_0 + \cdots + m_k = n$ and $m_1 + \cdots + m_k \leq s$. 
By construction, we have inclusion maps
$$ J(n+1) = J(n,-1) \hookrightarrow J(n,0) \hookrightarrow J(n,1) \hookrightarrow \cdots \hookrightarrow J(n,n) = J(n).$$
Note that, for $y$ as in (\ref{equation:define-y}), we have 
$$ \varphi(y) / p^{n} \prod_{i=1}^{k} \frac{ (pm_i)!}{ p^{m_i} m_i!} \gamma_{pm_i}( [x_i ] )
\equiv \prod_{i=1}^{k} (-1)^{m_i} \gamma_{pm_i}( [x_i ] ) \pmod{p},$$
so that $\beta$ restricts to a homomorphism $\beta_{s}: J(n,s) / J(n+1) \rightarrow \Fil_{s}^{\conj}( A_{\crys}(R) / p A_{\crys}(R) )$ for each $-1 \leq s \leq n$.
We will complete the proof by showing that each $\beta_{s}$ is an isomorphism. The proof proceeds by induction on $s$, the case $s = -1$ being trivial.
To carry out the inductive step, we use the commutative diagram of short exact sequences
$$ \xymatrix@R=50pt@C=50pt{ 0 \ar[d] & 0 \ar[d] \\
 J(n,s-1) / J(n+1) \ar[r]^-{\beta_{s-1}} \ar[d] & \Fil_{s-1}^{\conj}( A_{\crys}(R) / p A_{\crys}(R) ) \ar[d] \\
J(n,s) / J(n+1) \ar[d] \ar[r]^-{\beta_s} &  \Fil_{s}^{\conj}( A_{\crys}(R) / p A_{\crys}(R) ) \ar[d] \\
J(n,s) / J(n,s-1) \ar[d] \ar[r]^-{\overline{\beta}_s} & \gr_{s}^{\conj}( A_{\crys}(R) / p A_{\crys}(R) ) \ar[d] \\ 
0 & 0. }$$
By virtue of the snake lemma, it will suffice to show that the map $\overline{\beta}_{s}$ is an isomorphism for $0 \leq s \leq n$.

We now observe that the construction 
$$ \gamma_{m_1}(x_1) \cdots \gamma_{m_k}(x_k) \mapsto p^{n-s} \gamma_{m_1}([x_1]) \gamma_{m_2}( [x_2 ] ) \cdots \gamma_{m_k}( [x_k] )$$
determines an $R^{\flat}$-linear surjection $\Gamma^{s}_{R^{\flat}}(I) \twoheadrightarrow J(n,s) / J(n,s-1)$, which factors through an $R$-linear surjection
$\alpha: \Gamma^{s}_{R}(I/I^2) \twoheadrightarrow J(n,s) / J(n,s-1)$. Note that the composite map
$$ \Gamma^{s}_{R}( I/I^2) \xrightarrow{\alpha} J(n,s) / J(n,s-1) \xrightarrow{ \overline{\beta}_{s} } \gr_{s}^{\conj}( A_{\crys}(R) / p A_{\crys}(R) )$$
is given (up to sign) by the surjection $\xi$ of Proposition \ref{proposition:describe-gr-conjugate}, which is an isomorphism
by virtue of our assumption that $R$ is quasiregular semiperfect (Remark \ref{remark:gr-conjugate-isomorphism}). It follows
that $\overline{\beta}_{s}$ is also an isomorphism, as desired.
\end{proof}
  
\begin{remark}\label{remark:vshod}
Let $S$ be a perfect $\F_p$-algebra and let $R = S / (x_1, x_2, \cdots, x_k )$ be the quotient of $S$ by a regular sequence $(x_1, \cdots, x_k )$.
For $1 \leq i \leq d$, let us abuse notation by identifying $x_i$ with its image in $R^{\flat}$ and the Teichm\"{u}ller representative $[ x_i ] \in W( R^{\flat} )$ with
its image in $A_{\crys}(R)$. For each integer $n \geq 0$, let $J(n) \subseteq A_{\crys}(R)$ be the ideal defined in Proposition \ref{proposition:concrete-Nygaard-qrsp}.
Then the quotient $J(n) / J(n+1)$ is free as an $R$-module, with basis given by the collection of all products
$$ p^{m_0} \gamma_{m_1}([x_1]) \gamma_{m_2}( [x_2 ] ) \cdots \gamma_{m_k}( [x_k] ),$$
with $m_0 + m_1 + \cdots + m_k = n$.
\end{remark}

We now turn to the case where $R$ is a regular Noetherian $\F_p$-algebra. Let $W\Omega^{\ast}_{R}$ denote the de Rham-Witt complex of $R$. Following the convention of \cite{DRW}, we regard $W\Omega^{\ast}_{R}$ as an object of the ordinary category of cochain complexes over $\Z_p$, and we write
$W \Omega_{R}$ for its image in the $\infty$-category $\widehat{\calD}(\Z_p)$, so that we have a canonical isomorphism $\RGamma_{\crys}(R / \Z_p) \simeq W \Omega_{R}$
when $R$ is regular and Noetherian (see Theorem~10.1.1 of \cite{DRW} and its proof). For each integer $n \geq 0$, we let $\calN^{n} W\Omega^{\ast}_{R}$
denote the $m$th stage of the Nygaard filtration of the cochain complex $W\Omega^{\ast}_{R}$ (Construction~8.1.1 of \cite{DRW}), and let $\calN^{n} W\Omega_{R}$ denote its image
in $\widehat{\calD}(\Z_p)$. Concretely, $\calN^{n} W\Omega^{\ast}_{R}$ is the subcomplex of $W \Omega^{\ast}_{R}$ consisting of those elements $x$ for which
$\varphi(x)$ is divisible by $p^n$, where $\varphi: W\Omega^{\ast}_{R} \rightarrow W\Omega^{\ast}_{R}$ is the endomorphism of $W\Omega^{\ast}_{R}$ induced by the Frobenius endomorphism of
$R$. The following result is essentially Theorem~8.17 of \cite{BMS2}:

\begin{proposition}\label{proposition:concrete-Nygaard-regular}
Let $R$ be a regular Noetherian $\F_p$-algebra. Then the crystalline comparison isomorphism $\gamma_{\Prism}^{\crys}: \Prism_{R} \simeq \RGamma_{\crys}(R/\Z_p)$ of Theorem \ref{theorem:crystalline-comparison} can be promoted to an isomorphism
$\Fil^{\bullet}_{\Nyg} \Prism_{R/\Z_p} \simeq \calN^{\bullet} W\Omega_{R}$ in the filtered derived $\infty$-category $\DFiltI(\Z_p)$.
\end{proposition}

\begin{proof}
By virtue of Proposition \ref{proposition:universal-of-animated}, there is an essentially unique functor 
$$ \CAlg^{\anim}_{\F_p} \rightarrow \DFiltI(\Z_p) \quad \quad R \mapsto \calN^{n} LW\Omega_{R}$$
which commutes with sifted colimits and agrees with the functor $R \mapsto \calN^{n} W\Omega_{R}$ in the case where $R$ is a finitely
generated polynomial algebra over $\F_p$. By virtue of Theorem \ref{theorem:crystalline-comparison}, there is a functorial identification of
$\calN^{0} LW\Omega_{R}$ with the prismatic complex $\Prism_{R/\Z_p}$ (see Remark \ref{remark:absolute-as-derived-crystalline}),
so we can regard $\calN^{\bullet} W\Omega_{R}$ as a filtration on the complex $\Prism_{R/\Z_p} \simeq F^{\ast} \Prism_{R/\Z_p}$.

Note that, when $R$ is a polynomial algebra, the Frobenius on $W\Omega^{\ast}_{R}$ defines a map of filtered cochain complexes
$\calN^{\bullet} W\Omega^{\ast}_{R} \rightarrow p^{\bullet} W\Omega^{\ast}_{R}$. 
Passing to the filtered derived $\infty$-category $\DFiltI(\Z_p)$, we obtain a map 
$$\Fil^{\bullet}(\varphi): \calN^{\bullet} W \Omega_{R} \rightarrow p^{\bullet} W\Omega_{R} \simeq p^{\bullet} \Prism_{R/\Z_p}.$$
Since the functor $R \mapsto \calN^{\bullet} LW\Omega_{R}$ is a left Kan extension of its restriction to polynomial algebras over $\F_p$,
this construction admits an essentially unique extension to a natural transformation
$$\Fil^{\bullet}(\varphi): \calN^{\bullet} LW \Omega_{R} \rightarrow p^{\bullet} \Prism_{R/\Z_p}$$
for all animated commutative $\F_p$-algebras $R$. We claim that this map satisfies the requirements of Proposition \ref{proposition:relative-nygaard-characterization},
and therefore identifies $\calN^{\bullet} LW\Omega_{R}$ with the relative Nygaard filtration on the complex $F^{\ast} \Prism_{R/\Z_p}$.
To prove this, it will suffice to show that if $R$ is a finitely generated polynomial ring over $\F_p$ and $n \geq 0$ is an integer,
then $\gr^{n}(\varphi)$ induces an isomorphism from the cofiber $\calN^{n} W\Omega_{R} / \calN^{n+1} W\Omega_{R}$
to the cohomological truncation $\tau^{\leq n} (p^{n} W \Omega_{R} / p^{n+1} W\Omega_{R} )$ in the derived $\infty$-category $\widehat{\calD}(\Z_p)$.
In fact, an even stronger result is true: the divided Frobenius $\varphi$ induces an isomorphism
$$ \calN^{n} W\Omega^{\ast}_{R} / \calN^{n+1} W\Omega^{\ast}_{R} \rightarrow \tau^{\leq n} ( p^{n} W\Omega^{\ast}_{R} / p^{n+1} W \Omega^{\ast}_{R} )$$
in the abelian category of cochain complexes: see Proposition~8.2.1 of \cite{DRW}.

Note that, since the functor $R \mapsto \calN^{\bullet} LW\Omega_{R}$ is a left Kan extension of its restriction to polynomial algebras,
there is an essentially unique comparison map 
$$\Fil^{\bullet}(\alpha_R): \calN^{\bullet} LW\Omega_{R} \rightarrow \calN^{\bullet} W\Omega_{R}$$
which depends functorially on $R$ and agrees with the identity map when $R$ is a finitely generated polynomial algebra over $\F_p$.
To complete the proof of Proposition \ref{proposition:concrete-Nygaard-regular}, it will suffice to show that if $R$ is a regular Noetherian
$\F_p$-algebra and $n \geq 0$ is an integer, then the map $\Fil^{n}( \alpha_R): \calN^{n} LW\Omega_{R} \rightarrow \calN^{n} W\Omega_{R}$
is an isomorphism in the $\infty$-category $\widehat{\calD}(\Z_p)$. In the case $n=0$, this follows from Theorem \ref{theorem:crystalline-comparison}.
In particular, we can identify the quotient $p^{n} W\Omega_{R} /p^{n+1} W\Omega_{R}$ with the relative Hodge-Tate complex $\overline{\Prism}_{R/\Z_p}\{n\}$.
Using Proposition~8.2.1 of \cite{DRW}, we obtain a commutative diagram of fiber sequences
$$ \xymatrix@R=50pt@C=50pt{ \calN^{n+1} LW\Omega_{R} \ar[d]^{ \Fil^{n+1}(\alpha_R)} \ar[r] & \calN^{n} LW\Omega_{R} \ar[d]^{ \Fil^{n}(\alpha_R) } \ar[r] & \Fil_{n}^{\conj} \overline{\Prism}_{R/\Z_p}\{n\} \ar[d]^{\gr^{n}(\alpha_R)} \\
\calN^{n+1} W\Omega_{R} \ar[r] & \calN^{n} W\Omega_{R} \ar[r] & \tau^{\leq n} \overline{\Prism}_{R/\Z_p}\{n\}. }$$
We are therefore reduced to showing that the map $\gr^{n}(\alpha_R)$ is an isomorphism in $\widehat{\calD}(\Z_p)$: that is, that the
conjugate filtration on $\overline{\Prism}_{R/\Z_p}$ coincides with the Postnikov filtration. This is a special case of Remark \ref{remark:conjugate-equals-postnikov-conditions}:
our hypothesis that $R$ is regular guarantees that the cotangent complex $L\Omega^{1}_{R/\F_p}$ is a flat $R$-module (concentrated in cohomological degree zero).
\end{proof}

\subsection{The Absolute de Rham Comparison}\label{subsection:absolute-de-Rham-comparison}

Let $\mathfrak{X}$ be a $p$-adic formal scheme which is smooth over $\Spf(\Z_p)$, and let $\RGamma_{\dR}( \mathfrak{X} )$
denote its de Rham complex. The theory of crystalline cohomology originates from the observation that, up to quasi-isomorphism, the complex
$\RGamma_{\dR}( \mathfrak{X} )$ depends only on the special fiber of $\mathfrak{X}$. More precisely, we have the following:

\begin{theorem}[Berthelot, \cite{BerthelotCrysCoh}]\label{theorem:de-Rham-crys-comparison}
Let $\mathfrak{X}$ be a formal scheme which is smooth over $\Spf(\Z_p)$ and let $X = \Spec(\F_p) \times_{ \Spf(\Z_p) } \mathfrak{X}$ 
denote its special fiber. Then there is a canonical isomorphism
$$ \RGamma_{\crys}( X / \Z_p) \simeq \RGamma_{\dR}( \mathfrak{X} )$$
of commutative algebra objects of $\widehat{\calD}(\Z_p)$.
\end{theorem}

Theorem~\ref{theorem:de-Rham-crys-comparison} was conjectured by Grothendieck \cite{GrothendieckCrystals}, who was motivated by Monsky-Washnitzer's then-announced result \cite{MW} that $\RGamma_{\dR}(\mathfrak{X})$ only depended on $X$ in the above situation if one additionally also assumes that $X$ is proper. We refer to \cite{BhattdJCrysdR} for a relatively efficient recent proof of Theorem~\ref{theorem:de-Rham-crys-comparison}.

Our goal in this section is to sketch a different proof of Theorem \ref{theorem:de-Rham-crys-comparison}, which exploits the fact that both de Rham and
crystalline cohomology are closely related to prismatic cohomology. To simplify the discussion, let us suppose that the formal scheme
$\mathfrak{X} = \Spf(R)$ is affine. In this case, Theorem \ref{theorem:crystalline-comparison} identifies
$\RGamma_{\crys}( X / \Z_p )$ with the absolute prismatic complex of the quotient ring $R/pR$, and $\RGamma_{\dR}( \mathfrak{X} )$
can be identified with the $p$-complete derived de Rham complex $\widehat{\dR}_{R}$ (Proposition \ref{proposition:derived-to-classical-de-Rham}).
Theorem \ref{theorem:de-Rham-crys-comparison} can therefore be regarded as a special case of the following:

\begin{theorem}\label{theorem:deduce-comparison}
Let $R$ be an animated commutative ring and let $\overline{R}$ denote the derived tensor product $\F_p \otimes^{L} R$.
Then there is a canonical isomorphism $\Prism_{ \overline{R} } \simeq \widehat{\dR}_{R}$ of commutative
algebra objects of $\widehat{\calD}(\Z_p)$, depending functorially on $R$.
\end{theorem}

\begin{remark}
Let $R$ be a commutative ring which is $p$-torsion-free, and suppose that the quotient ring $\overline{R} = R/pR$ is regular and Noetherian.
Combining Theorems \ref{theorem:crystalline-comparison} and \ref{theorem:deduce-comparison} with Proposition \ref{proposition:derived-to-classical-de-Rham},
we obtain isomorphisms
$$ \RGamma_{\crys}( \overline{R} / \Z_p ) \simeq \Prism_{ \overline{R} } \simeq \widehat{\dR}_{R} \xrightarrow{\sim} ( \widehat{\Omega}^{\ast}_{R}, d).$$
These isomorphisms depend functorially on $R$, and therefore globalize to give an isomorphism
$$ \RGamma_{\crys}( X / \Z_p ) \simeq \RGamma_{\dR}( \mathfrak{X} )$$
when $\mathfrak{X}$ is a $p$-adic formal scheme which is flat over $\Spf(\Z_p)$ for which the special fiber $X = \Spec(\F_p) \times_{\Spf(\Z_p) } \mathfrak{X}$ is regular and Noetherian.
These hypotheses are satisfied, for example, if $\mathfrak{X}$ is smooth over $\Spf(W(k))$ for some perfect field $k$.
\end{remark}

Our proof of Theorem \ref{theorem:deduce-comparison} will exploit the geometry of the Cartier-Witt stack $\WCart$.

\begin{notation}
Let $\mathscr{E}$ be a quasi-coherent complex on the Cartier-Witt stack $\WCart$. We let $\mathscr{E}_{\dR}$ denote
the complex $\RGamma( \Spf(\Z_p), \rho_{\dR}^{\ast} \mathscr{E} )$, where $\rho_{\dR}: \Spf(\Z_p) \rightarrow \WCart$
is the de Rham point of Example \ref{example:de-Rham-point}. We regard $\mathscr{E}$ as an object of the derived $\infty$-category $\widehat{\calD}(\Z_p)$, which
we refer to as the {\it de Rham specialization} of $\mathscr{E}$.
\end{notation}

\begin{example}\label{example:dR-specialization-prismatic-sheaf}
Let $R$ be an animated commutative ring, let $\mathscr{H}_{\Prism}(R) \in \calD( \WCart )$ be the prismatic cohomology sheaf of
Construction \ref{construction:prismatic-cohomology-sheaves}, and set $\overline{R} = \F_p \otimes^{L} R$. By construction,
the de Rham specialization $\mathscr{H}_{\Prism}(R)_{\dR}$ can be identified with the prismatic complex
$\Prism_{ \overline{R} / \Z_p }$ of $\overline{R}$ relative to the crystalline prism $( \Z_p, (p) )$. Since the prism
$(\Z_p, (p) )$ is perfect, Proposition \ref{proposition:absolute-vs-relative} supplies an isomorphism
$\Prism_{ \overline{R} } \simeq \mathscr{H}_{\Prism}(R)_{\dR}$.
\end{example}

\begin{remark}\label{remark:de-Rham-specialization-pullback}
Let $\mathscr{E}$ be a quasi-coherent complex on the Cartier-Witt stack $\WCart$, and let $F^{\ast} \mathscr{E}$ denote the
pullback of $\mathscr{E}$ along the Frobenius map $F: \WCart \rightarrow \WCart$ of Construction \ref{construction:Frobenius-on-stack}
Then the commutative diagram
$$ \xymatrix@R=50pt@C=50pt{ \WCart^{\mathrm{HT}} \ar[r] \ar[d] & \WCart \ar[d]^{F} \\
\Spf( \Z_p) \ar[r]^-{ \rho_{\dR} } & \WCart }$$
of Proposition \ref{proposition:Frobenius-square} supplies a canonical isomorphism
$$ (F^{\ast} \mathscr{E})|_{ \WCart^{\mathrm{HT}} } \simeq \mathscr{E}_{\dR} \otimes \calO_{ \WCart^{\mathrm{HT}} }$$
in the derived $\infty$-category $\calD( \WCart^{\mathrm{HT}})$.
\end{remark}

\begin{construction}\label{construction:prismatic-sheaf-dr-pullback}
Let $R$ be an animated commutative ring, let $\mathscr{H}_{\Prism}(R) \in \calD( \WCart )$ denote the prismatic cohomology sheaf of
Construction \ref{construction:prismatic-cohomology-sheaves}, and let $F^{\ast} \mathscr{H}_{\Prism}(R) \in \calD( \WCart )$ denote the
pullback of $\mathscr{H}_{\Prism}(R)$ along the Frobenius morphism $F: \WCart \rightarrow \WCart$ of Construction \ref{construction:Frobenius-on-stack}.
For every bounded prism $(A,I)$, the pullback $\rho_{A}^{\ast} F^{\ast} \mathscr{H}_{\Prism}(R)$ can be identified with the complex
$F^{\ast} \Prism_{ (\overline{A} \otimes^{L} R) / A}$, where $\overline{A}$ denotes the quotient ring $A/I$ and
$\rho_{A}: \Spf(A) \rightarrow \WCart$ is the morphism described in Remark \ref{remark:HT-point-of-prismatic-stack}.
In particular, Proposition \ref{proposition:de-Rham-comparison-relative} supplies a canonical isomorphism
$$\overline{A} \otimes^{L}_{A} \rho_{A}^{\ast}(F^{\ast} \mathscr{H}_{\Prism}(R)) \simeq  \widehat{\dR}_{ ( \overline{A} \otimes^{L} R) / \overline{A} } 
\simeq \widehat{\dR}_{R} \widehat{\otimes}^{L} \overline{A}.$$
These isomorphisms depends functorially on the prism $(A,I)$, and therefore determine an isomorphism
$$ F^{\ast}\mathscr{H}_{\Prism}(R)|_{ \WCart^{\mathrm{HT}} } \simeq \widehat{\dR}_{R} \otimes \calO_{ \WCart^{\mathrm{HT}} }$$
of quasi-coherent complexes on the Hodge-Tate divisor $\WCart^{\mathrm{HT}}$. Combining this analysis with
Remark \ref{remark:de-Rham-specialization-pullback}, we obtain an canonical isomorphism
$$ \mathscr{H}_{\Prism}(R)_{\dR} \otimes \calO_{ \WCart^{\mathrm{HT}} } \simeq \widehat{\dR}_{R} \otimes \calO_{ \WCart^{\mathrm{HT}} }.$$
\end{construction}

 The terminology of Example \ref{example:de-Rham-point} is justified by the following observation:

\begin{proposition}\label{proposition:refined-de-Rham-comparison}
Let $R$ be an animated commutative ring. Then there is a canonical isomorphism
$$ \alpha_{R}: \mathscr{H}_{\Prism}(R)_{\dR} \simeq \widehat{\dR}_{R}$$
in the $\infty$-category $\widehat{\calD}(\Z_p)$, which is characterized by the requirement that it depends functorially on $R$
and that the induced map
$$ \mathscr{H}_{\Prism}(R)_{\dR} \otimes \calO_{ \WCart^{\mathrm{HT}} } \xrightarrow{ \alpha_{R} \otimes \id} \widehat{\dR}_{R} \otimes \calO_{ \WCart^{\mathrm{HT}} }$$
agrees with the isomorphism of Construction \ref{construction:prismatic-sheaf-dr-pullback} (up to a homotopy depending functorially on $R$).
\end{proposition}

\begin{proof}
Construction \ref{construction:prismatic-sheaf-dr-pullback} supplies a canonical map
$$ \alpha_{R}^{+}: \mathscr{H}_{\Prism}(R)_{\dR}
\rightarrow \RGamma( \WCart^{\mathrm{HT}}, F^{\ast} \mathscr{H}_{\Prism}(R)|_{ \WCart^{\mathrm{HT}} } )
\simeq \widehat{\dR}_{R} \widehat{\otimes} \RGamma( \WCart^{\mathrm{HT}}, \calO_{ \WCart^{\mathrm{HT}}} )$$
We will show that this map admits an essentially unique factorization 
$$ \mathscr{H}_{\Prism}(R)_{\dR} \xrightarrow{ \alpha_R } \widehat{\dR}_{R} \rightarrow \widehat{\dR}_{R} \widehat{\otimes} \RGamma( \WCart^{\mathrm{HT}}, \calO_{ \WCart^{\mathrm{HT}}} ) $$
which depends functorially on $R$. Note that in this case, the induced map
$$ \mathscr{H}_{\Prism}(R)_{\dR} \otimes \calO_{ \WCart^{\mathrm{HT}} } \rightarrow \widehat{\dR}_{R} \otimes \calO_{ \WCart^{\mathrm{HT}} }$$
of quasi-coherent complexes in $\WCart^{\mathrm{HT}}$ agrees with the isomorphism of Construction \ref{construction:prismatic-sheaf-dr-pullback},
so that $\alpha_{R}$ is automatically an isomorphism.

Using Example \ref{example:cohomology-of-HT}, we can identify the tensor product $\widehat{\dR}_{R} \widehat{\otimes} \RGamma( \WCart^{\mathrm{HT}}, \calO_{ \WCart^{\mathrm{HT}}} )$ with the direct sum $\widehat{\dR}_{R} \oplus \widehat{\dR}_{R}[-1]$. Composing $\alpha_{R}^{+}$ with the projection onto the second summand, we obtain a map
$$  \mathscr{H}_{\Prism}(R)_{\dR} \rightarrow \widehat{\dR}_{R}[-1].$$
We wish to show that this map has an essentially unique nullhomotopy which depends functorially on $R$.
Let $\calC$ denote the category of commutative rings by $R$ which are $p$-quasisyntomic and $p$-torsion-free,
and let $\calC_0 \subseteq \calC$ be the full subcategory spanned by those rings which are quasiregular semiperfectoid.
Since the functor $R \mapsto \mathscr{H}_{\Prism}(R)_{\dR} \in \widehat{\calD}(\Z_p)$ commutes with sifted colimits,
it is a left Kan extension of its restriction to $\calC$. The functor
$$ \calC \rightarrow \widehat{\calD}(\Z_p) \quad \quad R \mapsto \widehat{\dR}_{R}$$
is insensitive to $p$-completion and satisfies $p$-quasisyntomic descent (Variant \ref{variant:qs-de-Rham}), and is therefore a right Kan extension
of its restriction to $\calC_0$. We will complete the proof by showing that for $R,R' \in \calC_0$, the mapping space
$$\Hom_{ \widehat{\calD}(\Z_p)}(  \mathscr{H}_{\Prism}(R)_{|dR}, \widehat{\dR}_{R'}[-1] )$$
is contractible. This is clear, since the complexes $\mathscr{H}_{\Prism}(R)_{\dR} \simeq \Prism_{R/pR}$
and $\widehat{\dR}_{R'}$ are both concentrated in cohomological degree zero (see Example \ref{example:prismatic-sheaf-qrsp} and Proposition \ref{proposition:de-Rham-qrsp}).
\end{proof}

\begin{proof}[Proof of Theorem \ref{theorem:deduce-comparison}]
Combine Proposition \ref{proposition:refined-de-Rham-comparison} with Example \ref{example:dR-specialization-prismatic-sheaf}.
\end{proof}

\begin{construction}[The de Rham Comparison Map]\label{construction:de-Rham-comparison}
Let $R$ be an animated commutative ring. We let $\gamma_{\Prism}^{\dR}$ denote the composite map
$$ \Prism_{R} = \RGamma( \WCart, \mathscr{H}_{\Prism}(R) ) \rightarrow  \mathscr{H}_{\Prism}(R)_{\dR} \xrightarrow{\alpha_{R}} \widehat{\dR}_{R},$$
where $\alpha_{R}$ is the isomorphism of Proposition \ref{proposition:refined-de-Rham-comparison}. We will refer to
$\gamma_{\Prism}^{\dR}: \Prism_{R} \rightarrow \widehat{\dR}_{R}$ as the {\it de Rham comparison map} for the absolute prismatic complex $\Prism_{R}$.

Writing $\overline{R}$ for the derived tensor product $\F_p \otimes^{L} R$, we see that $\gamma_{\Prism}^{\dR}$ is obtained by composing the tautological map
$\Prism_{R} \rightarrow \Prism_{ \overline{R} }$ (given by the functoriality of absolute prismatic cohomology)
with the isomorphism $\Prism_{\overline{R}} \simeq \widehat{\dR}_{R}$ of Theorem \ref{theorem:deduce-comparison}.
\end{construction}

\begin{construction}[The Prismatic Augmentation]\label{construction:augmentation-on-prismatic-cohomology}
Let $R$ be an animated commutative ring and let $\widehat{R}$ denote its $p$-completion, which
we identify with $\gr^{0}_{\Hodge} \widehat{\dR}_{R}$. We let $\epsilon_{\Prism}: \Prism_{R} \rightarrow \widehat{R}$
denote the composite map
$$ \Prism_{R} \xrightarrow{ \gamma_{\Prism}^{\dR} } \widehat{\dR}_{R} \rightarrow \gr^{0}_{\Hodge} \widehat{\dR}_{R} \simeq \widehat{R}.$$
We will refer to $\epsilon_{\Prism}$ as the {\it prismatic augmentation}. Unwinding the definitions, we see that
$\epsilon_{\Prism}$ fits into a commutative diagram
\begin{equation}
\begin{gathered}\label{equation:augmentation-diagram}
\xymatrix@R=50pt@C=50pt{ \Prism_{R} \ar[r]^-{ \varphi } \ar[d]^{ \epsilon_{\Prism} } & \Prism_{R} \ar[d] \\
\widehat{R} \ar[r] & \overline{\Prism}_{R} }
\end{gathered}
\end{equation}
which depends functorially on $R$; here $\varphi$ denotes the Frobenius on the absolute prismatic complex $\Prism_{R}$, which we will study in \S\ref{subsection:Frobenius}.
\end{construction}

\begin{example}\label{example:semilinearity-warning}
Let $(A,I)$ be a perfect prism and let $\overline{A}$ denote the perfectoid ring $A/I$. It follows from
the commutativity of (\ref{equation:augmentation-diagram}) that, under the tautological identification $\Prism_{ \overline{A} } \simeq A$, the prismatic augmentation $\epsilon_{\Prism}: \Prism_{\overline{A}} \rightarrow \overline{A}$ corresponds to the composite map
$$ A \xrightarrow{ \varphi } A \twoheadrightarrow A/I = \overline{A}.$$
\end{example}

\begin{remark}
Let $R$ be an $\F_p$-algebra. Then the prismatic augmentation $\epsilon_{R}: \Prism_{R} \rightarrow R$
can be identified with the composition
$$ \Prism_{R} \xrightarrow{ \gamma_{\Prism}^{\crys} } \RGamma_{\crys}( R / \Z_p ) 
\xrightarrow{ \epsilon_{\crys} } R,$$
where $\gamma_{\Prism}^{\crys}$ is the comparison morphism of Remark \ref{remark:crystalline-comparison-general}.
This is essentially a restatement of Proposition \ref{proposition:augmentation-compatibility}.
\end{remark}


\begin{variant}[The de Rham Comparison Map: Twisted Version]\label{variant:de-Rham-comparison-twisted}
Let $R$ be an animated commutative ring. For every integer $n$, we let $\gamma_{\Prism}^{\dR}\{n\}$ denote the composite map
$$ \Prism_{R}\{n\} = \RGamma( \WCart, \mathscr{H}_{\Prism}(R)\{n\} ) \rightarrow \mathscr{H}_{\Prism}(R)_{\dR}\xrightarrow{\alpha_{R}} \widehat{\dR}_{R},$$
where we implicitly invoke the trivialization of the Breuil-Kisin twist $\calO_{\WCart}\{n\}_{\dR}$ supplied by generator
$e^{n}_{\Z_p} \in \Z_p\{n\}$ of Remark \ref{BKcrystalline}. We will refer to $\gamma_{\Prism}^{\dR}\{n\}$ as the {\it de Rham comparison map} for the twisted
absolute prismatic complex $\Prism_{R}\{n\}$.
\end{variant}

For every quasi-coherent complex $\mathscr{E}$ on the Cartier-Witt stack, Theorem \ref{theorem:Frobenius-pushout-square} guarantees that the diagram of restriction maps
$$ \xymatrix@R=50pt@C=50pt{ \RGamma( \WCart, \mathscr{E} ) \ar[r] \ar[d] & \mathscr{E}_{\dR} \ar[d] \\
\RGamma( \WCart, F^{\ast} \mathscr{E} ) \ar[r] & \RGamma( \WCart^{\mathrm{HT}}, (F^{\ast} \mathscr{E})|_{ \WCart^{\mathrm{HT}} } ) }$$
is a pullback square in the $\infty$-category $\widehat{\calD}(\Z_p)$. Applying this result in the case where the complex $\mathscr{E}$ has the form $\mathscr{H}_{\Prism}(R)\{n\}$ 
(and invoking the identifications of Construction \ref{construction:prismatic-sheaf-dr-pullback} and Proposition \ref{proposition:refined-de-Rham-comparison}), we obtain
the following:

\begin{corollary}\label{corollary:spanier}
Let $R$ be an animated commutative ring. For every integer $n$, the de Rham comparison map $\gamma_{\Prism}^{\dR}\{n\}: \Prism_{R} \rightarrow \widehat{\dR}_{R}$
fits into a pullback diagram
\begin{equation}
\begin{gathered}\label{equation:spanier}
\xymatrix@R=50pt@C=50pt{ \Prism_{R}\{n\} \ar[r]^-{\gamma_{\Prism}^{\dR}\{n\}} \ar[r] \ar[d] & \widehat{\dR}_{R} \ar[d] \\
\RGamma( \WCart, F^{\ast}( \mathscr{H}_{\Prism}(R)\{n\} ) ) \ar[r] & \widehat{\dR}_{R} \widehat{\otimes}^{L} \RGamma( \WCart^{\mathrm{HT}}, \mathscr{O}_{\WCart^{\mathrm{HT}}}). }
\end{gathered}
\end{equation}
\end{corollary}

\begin{remark}\label{remark:fiber-sequence-version-of-pullback}
In the situation of Corollary \ref{corollary:spanier}, we can identify $$\widehat{\dR}_{R} \widehat{\otimes}^{L} \RGamma( \WCart^{\mathrm{HT}}, \mathscr{O}_{\WCart^{\mathrm{HT}}})$$
with the direct sum $\widehat{\dR}_{R} \oplus \widehat{\dR}_{R}[-1]$ (see Example \ref{example:cohomology-of-HT}). Under this identification,
the right vertical map appearing in the diagram (\ref{equation:spanier}) corresponds to the inclusion of the first summand. It follows that
Corollary \ref{corollary:spanier} supplies a fiber sequence
$$ \Prism_{R}\{n\}  \rightarrow \RGamma( \WCart, F^{\ast}( \mathscr{H}_{\Prism}(R)\{n\} ) ) \rightarrow
\widehat{\dR}_{R}[-1]$$
in the $\infty$-category $\widehat{\calD}(\Z_p)$.
\end{remark}

\subsection{The Absolute Nygaard Filtration}\label{subsection:absolute-nygaard-filtration}

We now apply Corollary \ref{corollary:spanier} to construct a Nygaard filtration on the absolute prismatic complex of an animated commutative ring
(Construction \ref{construction:absolute-Nygaard-untwisted}).

\begin{notation}
Let $R$ be an animated commutative ring and let $n$ be an integer. For every bounded prism $(A,I)$ with quotient ring $\overline{A} = A/I$,
we can regard the tensor product $\overline{A} \otimes^{L} R$ as an animated commutative $\overline{A}$-algebra, and we write
$\Fil^{m}_{\Nyg} F^{\ast} \Prism_{  ( \overline{A} \otimes^{L} R)/A}$ for the $m$th stage of the relative Nygaard filtration
of Proposition \ref{proposition:relative-nygaard-characterization}. By virtue of Remark \ref{remark:Nygaard-change-of-prism},
the construction
$$ (A,I) \mapsto \Fil^{m}_{\Nyg} F^{\ast} \Prism_{  ( \overline{A} \otimes^{L} R)/A} \in \widehat{\calD}(A)$$
is compatible with (completed) extension of scalars, and therefore determines a quasi-coherent complex
on the Cartier-Witt stack which we will denote by $\Fil^{m}_{\Nyg} F^{\ast} \mathscr{H}_{\Prism}(R) \in \calD( \WCart )$
(see Proposition \ref{proposition:DWCart-prism-description}). This construction depends functorially on $R$ and $m$. For fixed $R$,
we obtain a diagram
$$ \cdots \Fil^{2}_{\Nyg} F^{\ast} \mathscr{H}_{\Prism}(R) \rightarrow
\Fil^{1}_{\Nyg} F^{\ast} \mathscr{H}_{\Prism}(R) \rightarrow \Fil^{0}_{\Nyg} F^{\ast} \mathscr{H}_{\Prism}(R) \simeq F^{\ast} \mathscr{H}_{\Prism}(R)$$
in the $\infty$-category $\calD( \WCart )$, which we denote by $\Fil^{\bullet}_{\Nyg} F^{\ast} \mathscr{H}_{\Prism}(R)$ and refer to
as the {\it Nygaard filtration} on the complex $\mathscr{H}_{\Prism}(R)$.
\end{notation}

\begin{remark}\label{remark:Nygaard-vs-hodge-globalized}
Let $R$ be an animated commutative ring, and let $\iota: \WCart^{\mathrm{HT}} \hookrightarrow \WCart$ denote the inclusion of the Hodge-Tate divisor.
Globalizing Corollary \ref{corollary:relative-Nygaard-deRham-refined}, we obtain a fiber sequence
$$ \mathscr{I} \Fil^{\bullet-1}_{\Nyg} F^{\ast} \mathscr{H}_{\Prism}(R) \rightarrow \Fil^{\bullet}_{\Nyg} F^{\ast} \mathscr{H}_{\Prism}(R)
\rightarrow \iota_{\ast} ( \Fil^{\bullet}_{\Hodge} \widehat{\dR}_{R} \otimes \calO_{ \WCart^{\mathrm{HT}} } )$$
of filtered complexes on $\WCart$, which can be obtained in filtration degree zero from the isomorphism
$(F^{\ast} \mathscr{H}_{\Prism}(R))|_{ \WCart^{\mathrm{HT}} } \simeq \widehat{\dR}_{R} \otimes \calO_{ \WCart^{\mathrm{HT}} }$
of Construction \ref{construction:prismatic-sheaf-dr-pullback}.
\end{remark}

\begin{construction}[The Absolute Nygaard Filtration]\label{construction:absolute-Nygaard-untwisted}
Let $n$ be an integer. For every animated commutative ring $R$, we let $\Fil^{\bullet}_{\Nyg} F^{\ast}( \mathscr{H}_{\Prism}(R)\{n\} )$
denote the filtered complex on $\WCart$ given by the tensor product of $\Fil^{\bullet}_{\Nyg} F^{\ast}( \mathscr{H}_{\Prism}(R) )$ with the line bundle $F^{\ast}( \calO_{\WCart}\{n\} )
\simeq \calI^{-n} \calO_{\WCart}\{n\}$. Since the line bundle $F^{\ast}( \calO_{\WCart}\{n\} )$ is canonically trivialized on the Hodge-Tate divisor $\WCart^{\mathrm{HT}}$
(Example \ref{example:BK-twist-on-HT}), Remark \ref{remark:Nygaard-vs-hodge-globalized} yields a fiber sequence
$$ \mathscr{I} \Fil^{\bullet-1}_{\Nyg} F^{\ast}(\mathscr{H}_{\Prism}(R)\{n\}) \rightarrow
\Fil^{\bullet}_{\Nyg} F^{\ast}(\mathscr{H}_{\Prism}(R)\{n\}) \rightarrow \iota_{\ast}( \Fil^{\bullet}_{\Hodge} \widehat{\dR}_{R} \otimes \calO_{ \WCart^{\mathrm{HT}} }).$$
Passing to global sections, we obtain a comparison map
\begin{eqnarray*} \RGamma( \WCart, \Fil^{\bullet}_{\Nyg} F^{\ast}(\mathscr{H}_{\Prism}(R)\{n\}) ) & \rightarrow & \RGamma( \WCart^{\mathrm{HT}}, \Fil^{\bullet}_{\Hodge} \widehat{\dR}_{R} \otimes \calO_{ \WCart^{\mathrm{HT}} }) \\
& \simeq & \Fil^{\bullet}_{\Hodge} \widehat{\dR}_{R} \widehat{\otimes} \RGamma( \WCart^{\mathrm{HT}}, \calO_{\WCart^{\mathrm{HT}} } ).
\end{eqnarray*}
Form a pullback diagram
$$\xymatrix@R=50pt@C=50pt{ \Fil^{\bullet}_{\Nyg} \Prism_{R}\{n\} \ar[r]^-{\Fil^{\bullet}(\gamma_{\Prism}^{\dR}\{n\})} \ar[d] & \Fil^{\bullet}_{\Hodge} \widehat{\dR}_{R} \ar[d] \\
\RGamma( \WCart, \Fil^{\bullet}_{\Nyg} F^{\ast}( \mathscr{H}_{\Prism}(R)\{n\} ) ) \ar[r] & \Fil^{\bullet}_{\Hodge} \widehat{\dR}_{R} \widehat{\otimes}^{L} \RGamma( \WCart^{\mathrm{HT}}, \mathscr{O}_{\WCart^{\mathrm{HT}}}) }$$
in the filtered derived $\infty$-category $\DFiltI(\Z_p)$. Note that Corollary \ref{corollary:spanier} supplies a canonical isomorphism $\Fil^{0}_{\Nyg} \Prism_{R}\{n\} \simeq \Prism_{R}\{n\}$,
which identifies $\Fil^{0}( \gamma_{\Prism}^{\dR}\{n\} )$ with the de Rham comparison map $\gamma_{\Prism}^{\dR}\{n\}: \Prism_{R}\{n\} \rightarrow \widehat{\dR}_{R}$ of Variant \ref{variant:de-Rham-comparison-twisted}. We will henceforth identify $\Fil^{0}_{\Nyg} \Prism_{R}\{n\}$ with $\Prism_{R}\{n\}$, and refer to $\Fil^{\bullet}_{\Nyg} \Prism_{R}\{n\}$ as the {\it Nygaard filtration} on the absolute prismatic complex $\Prism_{R}\{n\}$. In the special case $n =0$, we denote $\Fil^{\bullet}_{\Nyg} \Prism_{R}\{n\}$
simply by $\Fil^{\bullet}_{\Nyg} \Prism_{R}$.
\end{construction}

\begin{notation}[Prismatic-to-Hodge Comparison Map]\label{notation:prismatic-to-Hodge}
Let $R$ be an animated commutative ring, let $n$ be an integer, and let 
$$ \Fil^{\bullet}( \gamma_{\Prism}^{\dR}\{n\} ): \Fil^{\bullet}_{\Nyg} \Prism_{R}\{n\} \rightarrow \Fil^{\bullet}_{\Hodge} \widehat{\dR}_{R}$$
be the comparison map of Construction \ref{construction:absolute-Nygaard-untwisted}. Passing to the associated graded (and invoking
Remark \ref{remark:gr-derived-deRham}), we obtain maps
$$ \gamma_{\Prism}^{\Hodge}: \gr^{m}_{\Nyg} \Prism_{R}\{n\} \rightarrow L \widehat{\Omega}^{m}_{R}[-m],$$
which we will refer to as the {\it prismatic-to-Hodge comparison maps}. Note that $\gamma_{\Prism}^{\Hodge}$ can also be computed
as the composition
$$ \gr^{m}_{\Nyg} \Prism_{R}\{n\} \xrightarrow{\varphi\{n\}} \Fil_{m}^{\conj} \widehat{\Omega}^{\DHod}_{R}
\rightarrow \gr_{m}^{\conj} \widehat{\Omega}^{\DHod}_{R} \simeq L \widehat{\Omega}_{R}^{m}.$$
\end{notation}

\begin{remark}
Let $R$ be an animated commutative ring and let $n$ be an integer. The complex $\Fil^{m}_{\Nyg} \Prism_{R}\{n\}$ is defined for all integers $m$.
However, for $m \leq 0$, the transition map $\Fil^{0}_{\Nyg} \Prism_{R}\{n\} \rightarrow \Fil^{m}_{\Nyg} \Prism_{R}\{n\}$ is an isomorphism.
Consequently, we would lose no information by regarding the Nygaard filtration as indexed only by the nonnegative integers.
\end{remark}

\begin{remark}\label{remark:fiber-sequence-for-Nygaard}
Let $R$ be an animated commutative ring and let $n$ be an integer. Then the pullback diagram of Construction \ref{construction:absolute-Nygaard-untwisted} supplies
a fiber sequence of filtered complexes
$$ \Fil^{\bullet}_{\Nyg} \Prism_{R}\{n\} \rightarrow
\RGamma( \WCart, \Fil^{\bullet}_{\Nyg} F^{\ast}( \mathscr{H}_{\Prism}(R)\{n\} ))  \rightarrow \Fil^{\bullet}_{\Hodge} \widehat{\dR}_{R}[-1].$$
\end{remark}

\begin{remark}[Multiplicativity]\label{remark:multiplicativity-of-Nygaard}
Let $R$ be an animated commutative ring. Then the direct sum
$$ \bigoplus_{n \in \Z} \Fil^{\bullet}_{\Nyg} \Prism_R\{n\}$$
can be regarded as a graded commutative algebra object of the $\infty$-category $\DFilt(\Z_p)$. In particular:
\begin{itemize}
\item The filtered complex $\Fil^{\bullet}_{\Nyg} \Prism_{R}$ is a commutative algebra object of $\DFilt(\Z_p)$.
In particular, the multiplication on $\Prism_{R}$ refines to a map
$$ \Fil^{s}_{\Nyg} \Prism_{R} \otimes^{L} \Fil^{t}_{\Nyg} \Prism_{R} \rightarrow \Fil^{s+t}_{\Nyg} \Prism_{R}$$
for every pair of integers $(s,t)$.

\item For each $n \in \Z$, the filtered complex $\Fil^{\bullet}_{\Nyg} \Prism_{R}\{n\}$ has the structure of a module
over $\Fil^{\bullet}_{\Nyg} \Prism_{R}$.

\item The sum $\bigoplus_{m,n \in \Z} \gr^{m}_{\Nyg} \Prism_{R}\{n\}$ can be regarded as a bigraded commutative algebra object of $\calD(\Z_p)$.
\end{itemize}
\end{remark}

\begin{remark}[The Associated Graded of the Nygaard Filtration]\label{remark:Nygaard-associated-graded}
Let $R$ be an animated commutative ring. For every pair of integers $m$ and $n$, we let
$\gr^{m}_{\Nyg} \Prism_{R}\{n\}$ denote the cofiber of the transition map $\Fil^{m+1}_{\Nyg} \Prism_{R}\{n\} \rightarrow \Fil^{m}_{\Nyg} \Prism_{R}\{n\}$.
Combining Construction \ref{construction:absolute-Nygaard-untwisted} with \ref{remark:relative-Frobenius-sheaf-filtered}, we obtain a pullback diagram
$$ \xymatrix@R=50pt@C=50pt{ \gr^{m}_{\Nyg} \Prism_{R}\{n\} \ar[r]^-{ \gamma_{\Prism}^{\Hodge} } \ar[d] & L \widehat{\Omega}^{m}_{R}[-m] \ar[d] \\
\RGamma( \WCart^{\mathrm{HT}}, \Fil_{m}^{\conj} \mathscr{H}_{\overline{\Prism}}\{m\} ) \ar[r] & \RGamma( \WCart^{\mathrm{HT}}, 
\gr_{m}^{\conj} \mathscr{H}_{\overline{\Prism}}\{m\} ), }$$
where the right vertical map induces the isomorphism $$\gr_{m}^{\conj} \mathscr{H}_{\overline{\Prism}}\{m\} \simeq L \widehat{\Omega}^{m}_{R}[-m] \widehat{\otimes} \calO_{ \WCart^{\mathrm{HT}} }$$
of Remark \ref{remark:sheafy-HT-comparison}. Writing $\widehat{\Omega}^{\DHod}_{R}$ for the $p$-complete diffracted Hodge complex of $R$ of Construction \ref{construction:complete-diffracted-Hodge},
and $\Theta$ for its Sen operator, we obtain a fiber sequence
\begin{equation}
\label{eq:NygFib}
 \gr^{m}_{\Nyg} \Prism_{R}\{n\} \rightarrow \Fil_{m}^{\conj} \widehat{\Omega}^{\DHod}_{R} \xrightarrow{\Theta+m} \Fil_{m-1}^{\conj} \widehat{\Omega}^{\DHod}_{R},
 \end{equation}
where $\Theta+m$ denotes the $p$-completion of the map described in Remark \ref{remark:factor-Theta-plus}. 
\end{remark}

\begin{remark}
A fiber sequence similar to \eqref{eq:NygFib} is also constructed in \cite[Corollary 5.21]{AMMN}, using the ideas of \cite{KrauseNikolaus}. As observed in {\em loc.\ cit.}, it follows from \eqref{eq:NygFib} that $\mathrm{gr}^m_{\Nyg} \Prism_R\{n\} \in D^{\leq m}$ for all $m$ and $n$: indeed, we have $\Fil_{i}^{\conj} \widehat{\Omega}^{\DHod}_{R} \in D^{\leq i}$ for all $i$. 
\end{remark}

\begin{remark}\label{remark:graded-Nygaard-sifted-colimit}
Let $m$ and $n$ be integers. It follows from Remark \ref{remark:Nygaard-associated-graded} that the functor
$$ \CAlg^{\anim} \rightarrow \widehat{\calD}(\Z_p) \quad \quad R \mapsto \gr^{m}_{\Nyg} \Prism_{R}\{n\}$$
commutes with sifted colimits. Beware that the analogous statement is not true for the functor
$R \mapsto \Fil^{m}_{\Nyg} \Prism_{R}\{n\}$ (since it fails in the case $m=0$).
\end{remark}

\begin{remark}[Absolute Prismatic Gauges]
Let us define an {\it absolute prismatic gauge} to be a triple $( \Fil^{\bullet} \mathscr{E}_{\dR}, \Fil^{\bullet} \mathscr{E}, \alpha)$, where
$\Fil^{\bullet} \mathscr{E}_{\dR}$ is an object of the filtered derived $\infty$-category $\DFiltI(\Z_p)$, $\Fil^{\bullet} \mathscr{E}$
is a filtered quasi-coherent complex on the Cartier-Witt stack $\WCart$ equipped with an action of the filtered complex $\mathscr{I}^{\bullet}$
determined by the Hodge-Tate ideal sheaf, and 
$$\alpha: \Fil^{\bullet} \mathscr{E}_{\dR} \otimes \calO_{\WCart^{\mathrm{HT}} } \xrightarrow{\sim} \Fil^{\bullet} \mathscr{E} / \mathscr{I} \Fil^{\bullet-1} \mathscr{E}$$
is an isomorphism of filtered complexes on the Hodge-Tate divisor $\WCart^{\mathrm{HT}}$. 
Every animated commutative ring $R$ determines an absolute prismatic gauge $\mathscr{G}_{R} = ( \Fil^{\bullet}_{\Hodge} \widehat{\dR}_{R}, \Fil^{\bullet}_{\Nyg} F^{\ast} \mathscr{H}_{\Prism}(R), \alpha)$, where the isomorphism $\alpha$ is obtained by globalizing the construction of Corollary \ref{corollary:relative-Nygaard-deRham-refined}.
From any absolute prismatic gauge $( \Fil^{\bullet} \mathscr{E}_{\dR}, \Fil^{\bullet} \mathscr{E}, \alpha)$, one can extract a filtered complex
by forming the fiber product
$$ \RGamma( \WCart, \Fil^{\bullet} \mathscr{E} ) \times_{ \RGamma( \WCart^{\mathrm{HT}}, \Fil^{\bullet} \mathscr{E} / \mathscr{I} \Fil^{\bullet-1} \mathscr{E})}
\Fil^{\bullet} \mathscr{E}_{\dR};$$
applied to the prismatic gauge $\mathscr{G}_{R}$, this reproduces the absolute Nygaard filtration $\Fil^{\bullet}_{\Nyg} \Prism_{R}$ of
Construction \ref{construction:absolute-Nygaard-untwisted}.

The collection of absolute prismatic gauges can be organized into an $\infty$-category, which seems to be closely related to the $\infty$-category of quasi-coherent
complexes on a certain enlargement of the Cartier-Witt stack introduced by Drinfeld (and denoted by $\Sigma'$) in \cite{drinfeld-prismatic}.
\end{remark}

\begin{proposition}\label{proposition:low-gr-Nygaard}
Let $R$ be an animated commutative ring, let $n$ be an integer, and let $$\Fil^{\bullet}(\gamma_{\Prism}^{\dR}\{n\}): \Fil^{\bullet}_{\Nyg} \Prism_R\{n\} \rightarrow \Fil^{\bullet}_{\Hodge} \widehat{\dR}_{R}$$ denote the de Rham specialization map of Construction \ref{construction:absolute-Nygaard-untwisted}. Then:
\begin{itemize}
\item For every integer $m \geq 0$, the induced map
$$ \Prism_{R}\{n\} / \Fil^{m}_{\Nyg} \Prism_{R}\{n\} \rightarrow \widehat{\dR}_{R} / \Fil^{m}_{\Hodge} \widehat{\dR}_{R}$$
is an isogeny: that is, its fiber is annihilated by $p^{k}$ for $k \gg 0$.

\item For $0 \leq m < p$, the induced map
$$ \Prism_{R}\{n\} / \Fil^{m}_{\Nyg} \Prism_{R}\{n\} \rightarrow \widehat{\dR}_{R} / \Fil^{m}_{\Hodge} \widehat{\dR}_{R}$$
is an isomorphism.
\end{itemize}
\end{proposition}

\begin{proof}
We proceed by induction on $m$. To prove the first assertion, it will suffice to show that the map
$$ \gr^{m}( \gamma_{\Prism}^{\dR}\{n\}): \gr^{m}_{\Nyg} \Prism_{R}\{n\} \rightarrow \gr^{m}_{\Hodge} \widehat{\dR}_{R} \simeq L\widehat{\Omega}^{m}_{R}[-m]$$
is an isogeny. Using Remark \ref{remark:Nygaard-associated-graded}, we see that the fiber of $\gr^{m}( \gamma_{\Prism}^{\dR}\{n\} )$ can be identified with 
the complex $( \Fil_{m-1}^{\conj} \widehat{\Omega}^{\DHod}_{R} )^{\Theta = -m}$. We will prove more generally that for $0 \leq i < m$, the complex
$( \Fil_{i}^{\conj} \widehat{\Omega}^{\DHod}_{R} )^{\Theta = -m}$ is annihilated by a power of $p$. Proceeding by induction on $i$, we are reduced to showing
that $( \gr_{i}^{\conj} \widehat{\Omega}^{\DHod}_{R} )^{\Theta = -m}$ is annihilated by a power of $p$. This is clear, since $\Theta+m$ acts by
the nonzero scalar $m-i$ on the complex $\gr_{i}^{\conj} \widehat{\Omega}^{\DHod}_{R}$. The proof of the second assertion is similar, noting that
$m-i$ is a unit in $\Z_p$ when $0 \leq i < m < p$.
\end{proof}

\begin{example}\label{example:prismatic-augmentation-vs-Nygaard}
Let $R$ be an animated commutative ring and let $\widehat{R}$ denote its $p$-completion. Applying Proposition \ref{proposition:low-gr-Nygaard}
in the case $m=1$, we obtain a fiber sequence
$$ \Fil^{1}_{\Nyg} \Prism_{R} \rightarrow \Prism_{R} \xrightarrow{ \epsilon^{\Prism}_{R} } \widehat{R},$$
where $\epsilon^{\Prism}_{R}$ denotes the prismatic augmentation of Construction \ref{construction:augmentation-on-prismatic-cohomology}.
\end{example}

\begin{corollary}\label{corollary:gr0-Nygaard}
Let $R$ be an animated commutative ring. For every integer $n$, there is a unique element
$e^{n} \in \mathrm{H}^{0}( \gr^{0}_{\Nyg} \Prism_{R}\{n\} )$ satisfying
$\gr^{0}(\gamma_{\Prism}^{\dR}\{n\})(e^n) = 1 \in \mathrm{H}^{0}( \gr^{0}_{\Hodge} \widehat{\dR}_{R} )$.
\end{corollary}

\begin{remark}[Periodicity]\label{remark:laurent-in-e}
Let $R$ be an animated commutative ring. For every pair of integers $m$ and $n$, multiplication by the element $e^n$ of Corollary \ref{corollary:gr0-Nygaard}
induces an isomorphism
$$ \gr^{m}_{\Nyg} \Prism_{R} \xrightarrow{\sim} \gr^{m}_{\Nyg} \Prism_{R}\{n\},$$ 
whose inverse is given by multiplication by the element $e^{-n}$. These isomorphisms are compatible with the multiplicative structure described in
Remark \ref{remark:multiplicativity-of-Nygaard}: that is, they supply an isomorphism
$$ ( \bigoplus_{m \in \Z} \gr^{m}_{\Nyg} \Prism_{R} )[ e^{\pm 1} ] \simeq \bigoplus_{m,n \in \Z} \gr^{m}_{\Nyg} \Prism_{R}\{n\}$$
of bigraded commutative algebra objects of $\calD(\Z_p)$.
\end{remark}

\begin{example}\label{example:qdr-periodicity}
Let $(A,I)$ be the perfection of the $q$-de Rham prism $( \Z_p[[q-1]], [p]_q )$ (so that $A$ is the completion of
$\Z[ q^{1/p^{\infty}} ]$ with respect to the ideal $(p, q-1)$ and $I$ is generated by the element $[p]_q$),
and let $\Z_p^{\cyc}$ denote the quotient ring $A/I$. Then the element $e \in \gr^{0}_{\Nyg} \Prism_{ \Z_p^{\cyc} }\{1\}$ 
appearing in Remark \ref{remark:laurent-in-e} can be identified with the image of the element
$\widetilde{e} = \frac{ \log_{\Prism}(q^{p})}{q-1} \in A\{1\}$. Unwinding definitions and using Example~\ref{remark:preferred-generator-universal-case}, we must check that the canonical map 
\[ I^{-1} A\{1\} \to I^{-1} \otimes_A I/I^2 \simeq A/I\] 
carries the element 
\[ \varphi_{A\{1\}}(\widetilde{e}) = \frac{\log_{\Prism}(q^p)}{q^p-1} = [p]_q^{-1} \widetilde{e} \in I^{-1} A\{1\}\]
to $1 \in A/I$. But observe that the canonical surjection $A\{1\} \to I/I^2$ carries $\log_{\Prism}(q^{p} ) \in A\{1\}$ to the element $q^{p} - 1 \in I/I^2$ (Proposition \ref{proposition:normalization-of-logarithm}), and hence carries $\widetilde{e} \in A\{1\}$ to $[p]_q \in I/I^2$ (as $q-1$ is a nonzerodivisor in $A/I$). Multiplying by $[p]_q^{-1}$ then yields the claim. 
\end{example}

\begin{warning}\label{warning:qdr-periodicity-again}
Let $(A,I)$ be the perfection of the $q$-de Rham prism $( \Z_p[[q-1]], [p]_q )$. By virtue of Proposition \ref{proposition:twist-in-q-de-Rham-case}, the Breuil-Kisin twist $A\{1\}$ is a free $A$-module of rank $1$, generated by the element $\widetilde{e} = \log_{\Prism}(q^{p})/ (q-1)$.
If $R$ is any animated commutative algebra over the quotient ring $\Z_{p}^{\cyc} = A/I$, then multiplication by $\widetilde{e}^{n}$ induces an isomorphism of filtered complexes
$$ \Fil^{\bullet}_{\Nyg} \Prism_{R} \simeq \Fil^{\bullet}_{\Nyg} \Prism_{R}\{n\}.$$
Note that the induced isomorphism of associated graded complexes $\gr^{\bullet}_{\Nyg} \Prism_{R} \simeq \gr^{\bullet}_{\Nyg} \Prism_{R}\{n\}$
agrees with the isomorphism appearing in Remark \ref{remark:laurent-in-e}: this follows from the observation that the tautological map
$$ A\{1\} \simeq \Prism_{ \Z_p^{\cyc} }\{1\} \rightarrow \Prism_{R}\{1\} \rightarrow \gr^{0}_{\Nyg} \Prism_{R}\{1\}$$
carries $\widetilde{e}$ to $e$ (see Example \ref{example:qdr-periodicity}).
Beware that, while the element $e$ is canonical (and can be defined without using the $\Z_p^{\cyc}$-algebra structure on $R$),
the element $\widetilde{e}$ is not; for example, $\widetilde{e}$ is not invariant under the action of the automorphism group
$\Aut( \Z_p^{\cyc} ) \simeq \Z_p^{\times}$.
\end{warning}

\begin{corollary}
Let $R$ be an animated commutative ring. For $n < p$, there is a canonical isomorphism
$\gr^{n}_{\Nyg} \Prism_{R} \simeq L \widehat{\Omega}^{n}_{R}[-n]$ in the $\infty$-category $\widehat{\calD}(\Z_p)$.
\end{corollary}

\begin{proposition}\label{proposition:Nygaard-filtration-Beilinson-connective}
Let $R$ be an animated commutative ring, and let $n$ be an integer. Then, for every integer $m$, the
cohomology groups of the complex $\gr^{m}_{\Nyg} \Prism_{R}\{n\}$ are concentrated in degrees $\leq m$. In other words,
the filtered complex $\Fil^{\bullet}_{\Nyg} \Prism_{R}\{n\}$ is connective with respect to the Beilinson t-structure on $\DFiltI(\Z_p)$.
\end{proposition}

\begin{proof}
Combine the fiber sequence
$$ \gr^{m}_{\Nyg} \Prism_{R}\{n\} \rightarrow \Fil_{m}^{\conj} \widehat{\Omega}^{\DHod}_{R} \xrightarrow{\Theta+m} \Fil_{m-1}^{\conj} \widehat{\Omega}^{\DHod}_{R}$$
of Remark \ref{remark:Nygaard-associated-graded} with Remark \ref{remark:connectivity-of-diffracted-Hodge}.
\end{proof}

\begin{proposition}\label{proposition:gr-nygaard-flat-descent}
For every pair of integers $m$ and $n$, the functor
$$ \CAlg^{\anim} \rightarrow \widehat{\calD}(\Z_p) \quad \quad R \mapsto \gr^{m}_{\Nyg} \Prism_{R}\{n\}$$
satisfies descent for the $p$-completely faithfully flat topology.
\end{proposition}

\begin{proof}
Combine Remarks \ref{remark:fpqc-descent-diffracted-hodge} and \ref{remark:Nygaard-associated-graded}.
\end{proof}

\begin{corollary}\label{corollary:qs-descent-nygaard}
For every pair of integers $m$ and $n$, the functor
$$ \CAlg^{\QSyn} \rightarrow \widehat{\calD}(\Z_p) \quad \quad R \mapsto \Fil^{m}_{\Nyg} \Prism_{R}\{n\}$$
satisfies descent for the $p$-quasisyntomic topology.
\end{corollary}

\begin{proof}
The case $m \leq 0$ follows from Proposition \ref{proposition:absolute-quasisyntomic-descent}. The general case follows by induction on $m$,
using Proposition \ref{proposition:gr-nygaard-flat-descent} and Corollary \ref{corollary:absolute-completion-invariance-Nygaard}.
\end{proof}

\begin{corollary}\label{corollary:absolute-etale-affine-Nygaard}
For every pair of integers $m$ and $n$, the functor
$$ \CAlg^{\anim} \rightarrow \widehat{\calD}(\Z_p) \quad \quad R \mapsto \Fil^{m}_{\Nyg} \Prism_{R}\{n\}$$
satisfies descent for the \'{e}tale topology.
\end{corollary}

\begin{proof}
The case $m \leq 0$ follows from Proposition \ref{proposition:absolute-etale-affine}. The general case follows by induction on $m$, using 
Proposition \ref{proposition:gr-nygaard-flat-descent}.
\end{proof}

\begin{notation}\label{notation:absolute-Nygaard-filtration-of-scheme}
Let $X$ be scheme, formal scheme, or algebraic stack. For every pair of integers $m$ and $n$, we let $\Fil^{m}_{\Nyg} \RGamma_{\Prism}(X)\{n\}$ denote the inverse limit $\varprojlim_{ \Spec(R) \rightarrow X} \Fil^{m}_{\Nyg} \Prism_{R/A}\{n\}$, formed in the $\infty$-category $\widehat{\calD}(\Z_p)$; here the limit is taken over the category of
all commutative rings $R$ equipped with a map $\Spec(R) \rightarrow X$. Note that, for $m \leq 0$, we can identify
$\Fil^{m}_{\Nyg} \RGamma_{\Prism}(X)\{n\}$ with the complex $\RGamma_{\Prism}(X)\{n\}$ of Construction \ref{construction:absolute-prismatic-complex-of-scheme}.
It follows from Corollary \ref{corollary:absolute-etale-affine-Nygaard} that the construction $X \mapsto \Fil^{m}_{\Nyg} \RGamma_{\Prism}(X)\{n\}$ satisfies
descent for the \'{e}tale topology.
\end{notation}

\begin{proposition}\label{proposition:derived-descent-absolute-Nygaard}
Let $\F_p^{\otimes \bullet+1}$ be the cosimplicial animated commutative ring introduced in Notation \ref{notation:F-p-bullet}.
Let $R$ be an animated commutative ring, and regard $R^{\bullet} = R \otimes^{L} \F_p^{\otimes \bullet+1}$ as a cosimplicial object of
$\CAlg^{\anim}$. For every pair of integers $m$ and $n$, the tautological map 
$$ \Fil^{m}_{\Nyg} \Prism_{R}\{n\} \rightarrow \Tot( \Fil^{m}_{\Nyg} \Prism_{ R^{\bullet} }\{n\} )$$
is an isomorphism in the $\infty$-category $\widehat{\calD}(\Z_p)$.
\end{proposition}

\begin{proof}
We proceed by induction on $m$. The case $m \leq 0$ follows from Proposition \ref{proposition:derived-descent-absolute}.
To carry out the inductive step, it will suffice to show that the comparison map
$$ \gr^{m}_{\Nyg} \Prism_{R}\{n\} \rightarrow \Tot( \gr^{m}_{\Nyg} \Prism_{ R^{\bullet} }\{n\} )$$
is an isomorphism. Using the fiber sequence of Remark \ref{remark:Nygaard-associated-graded},
we are reduced to proving that for every integer $d$, the comparison map
$$ \Fil_{d}^{\conj} \widehat{\Omega}^{\DHod}_{R} \rightarrow \Tot(  \Fil_{d}^{\conj} \widehat{\Omega}^{\DHod}_{R^{\bullet}} )$$
is an isomorphism. This is a special case of Variant \ref{variant:conjugate-derived-descent} (see Example \ref{example:diffracted-Hodge-as-relative-prismatic}).
\end{proof}

\begin{corollary}\label{corollary:reduce-mod-n-Nygaard}
Let $R$ be an animated commutative ring. For every pair of integers $m$ and $n$, the tautological map
$$ \Fil^{m}_{\Nyg} \Prism_{R}\{n\} \rightarrow \varprojlim_{k}  \Fil^{m}_{\Nyg} \Prism_{(\Z / p^{k} \Z) \otimes^{L} R}\{n\}$$
is an isomorphism in the $\infty$-category $\widehat{\calD}(\Z_p)$.
\end{corollary}

\begin{corollary}\label{corollary:absolute-completion-invariance-Nygaard}
Let $R$ be an animated commutative ring and let $\widehat{R}$ denote the $p$-completion of $R$.
For every pair of integers $m$ and $n$, the tautological map $\Fil^{m}_{\Nyg} \Prism_{R}\{n\} \rightarrow \Fil^{m}_{\Nyg} \Prism_{ \widehat{R}} \{n\}$ is an isomorphism in the
$\infty$-category $\widehat{\calD}(\Z_p)$.
\end{corollary}

\begin{corollary}
Let $X$ be a scheme for which the structure sheaf $\calO_{X}$ has bounded $p$-power torsion, and let $\mathfrak{X} = \Spf(\Z_p) \times X$ be the associated $p$-adic formal scheme.
For every pair of integers $m$ and $n$, the restriction map $\Fil^{m}_{\Nyg} \RGamma_{\Prism}(X)\{n\} \rightarrow \Fil^{m}_{\Nyg} \RGamma_{\Prism}(\mathfrak{X})\{n\}$
is an isomorphism.
\end{corollary}

\begin{proof}
Without loss of generality, we may assume that $X = \Spec(R)$ is affine, in which case the desired result is a special case of Corollary \ref{corollary:reduce-mod-n-Nygaard}.
\end{proof}

\subsection{The Case of a Perfect Prism}\label{subsection:relative-nygaard-comparison}

Let $(A,I)$ be a bounded prism. For every animated commutative algebra $R$ over the quotient ring $\overline{A} := A/I$, we have a comparison map
$$ \Prism_{R} = \RGamma( \WCart, \mathscr{H}_{\Prism}(R) ) \rightarrow \RGamma( \Spf(A), \rho_{A}^{\ast} \mathscr{H}_{\Prism}(R) ) =
\Prism_{ (\overline{A} \otimes^{L} R) / A} \rightarrow \Prism_{R/A},$$
which is an isomorphism if the prism $(A,I)$ is perfect (Proposition \ref{proposition:absolute-vs-relative}). Our goal in this section is to show that, under the same assumption, this isomorphism carries the absolute Nygaard filtration on $\Prism_{R}$ to the relative Nygaard filtration on the Frobenius pullback $F^{\ast} \Prism_{R/A}$.
of Construction \ref{construction:absolute-Nygaard-untwisted} (see Theorem \ref{theorem:compare-Nygaard-filtrations}).

\begin{construction}\label{construction:compare-Nygaard-filtrations}
Let $(A,I)$ be a bounded prism and let $R$ be an animated commutative algebra over the quotient ring $\overline{A} = A/I$. For every integer $n$, we let
$\Fil^{\bullet}_{\Nyg} F^{\ast}( \Prism_{R/A}\{n\} )$ denote the tensor product of the filtered complex $\Fil^{\bullet}_{\Nyg} F^{\ast} \Prism_{R/A}$ with the invertible $A$-module
$F^{\ast}( A\{n\} ) \simeq I^{-n} A\{n\}$. We then have a comparison map of filtered complexes
\begin{eqnarray*} \Fil^{\bullet}_{\Nyg} \Prism_{R}\{n\} & \rightarrow & \RGamma( \WCart, \Fil^{\bullet}_{\Nyg} F^{\ast}(\mathscr{H}_{\Prism}(R)\{n\}) )  \\
& \rightarrow & \RGamma( \Spf(A),  \rho_{A}^{\ast} \Fil^{\bullet}_{\Nyg} F^{\ast}(\mathscr{H}_{\Prism}(R)\{n\}) ) \\
& \simeq & \Fil^{\bullet}_{\Nyg} F^{\ast}(\Prism_{ (\overline{A} \otimes^{L} R) / A}\{n\}) \\
& \rightarrow & \Fil^{\bullet}_{\Nyg} F^{\ast}(\Prism_{ R / A}\{n\}). \end{eqnarray*}
where the source involves the Nygaard filtration on the absolute prismatic cohomology (Construction \ref{construction:absolute-Nygaard-untwisted})
and the target involves the Nygaard filtration on relative prismatic cohomology (Proposition \ref{proposition:relative-nygaard-characterization}).
\end{construction}

\begin{theorem}\label{theorem:compare-Nygaard-filtrations}
Let $(A,I)$ be a perfect prism. For every animated commutative $\overline{A}$-algebra $R$ and every integer $n$, the comparison map
$\Fil^{\bullet}_{\Nyg} \Prism_{R}\{n\} \rightarrow\Fil^{\bullet}_{\Nyg} F^{\ast}(\Prism_{R/A}\{n\})$ of Construction \ref{construction:compare-Nygaard-filtrations}
is an isomorphism (in the filtered derived $\infty$-category $\DFiltI(\Z_p)$).
\end{theorem}

Theorem \ref{theorem:compare-Nygaard-filtrations} is essentially equivalent to the following consequence, which does not mention the relative Nygaard filtration:

\begin{corollary}\label{corollary:absolute-Nygaard-qrsp}
Let $R$ be a commutative ring which is quasiregular semiperfectoid, so that we can regard $( \Prism_{R}, \Prism_{R}^{[1]})$ as a prism
(see Example \ref{example:prismatic-sheaf-qrsp}). For every pair of integers $m$ and $n$, the canonical map
$\Fil^{m}_{\Nyg} \Prism_{R}\{n\} \rightarrow \Prism_{R}\{n\}$ identifies $\Fil^{m}_{\Nyg} \Prism_{R}\{n\}$ with the subgroup of
$\Prism_{R}\{n\}$ consisting of those elements $x$ satisfying $\varphi(x) \in \Prism_{R}^{[m-n]}\{n\}$. In particular, each of the complexes $\Fil^{m}_{\Nyg} \Prism_{R}\{n\}$ is concentrated in cohomological degree zero.
\end{corollary}

\begin{proof}
Choose a ring homomorphism $\overline{A} \rightarrow R$, where $\overline{A}$ is a perfectoid ring, and write $\overline{A} = A/I$ for a perfect prism $(A,I)$.
We then have a commutative diagram of filtered complexes
$$ \xymatrix@R=50pt@C=50pt{ \Fil^{\bullet}_{\Nyg} \Prism_{R}\{n\} \ar[r] \ar[d]^{\varphi} & F^{\ast}( A\{n\} ) \otimes^{L}_{A} \Fil^{\bullet}_{\Nyg} F^{\ast} \Prism_{R/A} \ar[d]^{\varphi} \\
\Prism_{R}^{[\bullet-n]}\{n\} \ar[r] & I^{\bullet-n}\{n\} \otimes_{A}^{L} F^{\ast}( \Prism_{R/A} )\{n\} }$$
where the horizontal maps are isomorphisms by virtue of Theorem \ref{theorem:compare-Nygaard-filtrations} and Proposition \ref{proposition:absolute-vs-relative}.
It will therefore suffice to show that $\Fil^{m}_{\Nyg} F^{\ast} \Prism_{R/A}$ identifies with the $A$-submodule of $F^{\ast} \Prism_{R/A}$
consisting of those elements $y$ satisfying $\varphi(y) \in I^{m} \Prism_{R/A}$, which was established in the proof of Proposition \ref{proposition:relative-nygaard-characterization}.
\end{proof}

\begin{corollary}\label{corollary:qrsp-filtered-discrete}
Let $R$ be a commutative ring which is quasiregular semiperfectoid. Then, for every pair of integers $m$ and $n$, the complex $\gr^{m}_{\Nyg} \Prism_{R}\{n\}$ is concentrated in cohomological degree zero.
\end{corollary}

\begin{remark}
If $R$ is a $p$-quasisyntomic commutative ring, then the Nygaard filtration on the absolute prismatic complex $\Prism_R\{n\}$ is completely determined by
Corollaries \ref{corollary:absolute-Nygaard-qrsp} and \ref{corollary:qs-descent-nygaard}. More precisely, choose a $p$-quasisyntomic cover $R \rightarrow R^{0}$, where $R^{0}$ is quasiregular semiperfectoid, and let $R^{\bullet}$ denote the cosimplicial $R$-algebra given by the $p$-complete tensor powers of $R^{0}$.
Then each $\Fil^{m}_{\Nyg} \Prism_{R^{\bullet}}\{n\}$ is a cosimplicial abelian group (with an explicit description given in Corollary \ref{corollary:absolute-Nygaard-qrsp}),
whose image under the Dold-Kan correspondence is an explicit model for $\Fil^{m}_{\Nyg} \Prism_{R}\{n\}$ at the level of cochain complexes
(by virtue of Corollary \ref{corollary:absolute-Nygaard-qrsp}). 
\end{remark}

\begin{remark}
Let $R$ be a quasiregular semiperfectoid ring. Then there exists a perfect prism $(A,I)$ and a ring homomorphism
$A/I \rightarrow R$. For every integer $m$, we have a canonical isomorphism $\overline{\Prism}_{R}\{m\} \simeq \overline{\Prism}_{R/A}\{m\}$ (Example \ref{example:relative-HT-over-perfect}). Tensoring the relative conjugate filtration of Remark \ref{remark:derived-HT-filtration} with $A\{m\}$, we obtain a filtration
$$  \Fil_{0}^{\conj} \overline{\Prism}_{R/A}\{m\} \rightarrow
\Fil_{1}^{\conj} \overline{\Prism}_{R/A}\{m\} \rightarrow \Fil_{2}^{\conj} \overline{\Prism}_{R/A}\{m\} \rightarrow \cdots$$
of the absolute Hodge-Tate complex $\overline{\Prism}_{R}\{m\}$. Beware that this filtration depends not only on $R$, but also on the choice of perfect prism $(A,I)$.
However, Remark~\ref{remark:polynomial-to-smooth} and Thorem~\ref{theorem:compare-Nygaard-filtrations} ensure that the $m$th stage of this filtration depends only on $R$ (since it can be identified with $\gr^{m}_{\Nyg} \Prism_R$).
\end{remark}

\begin{remark}\label{remark:local-qrsp-conclusion}
Suppose that the conclusion of Theorem \ref{theorem:compare-Nygaard-filtrations} holds for a {\em particular} perfect prism $(A,I)$. In that case,
the conclusion of Corollary \ref{corollary:absolute-Nygaard-qrsp} holds for every quasiregular semiperfectoid ring $R$ which admits the structure of an
$A/I$-algebra. 
\end{remark}

We will use Remark \ref{remark:local-qrsp-conclusion} to reduce the proof of Theorem \ref{theorem:compare-Nygaard-filtrations} to the following special case:

\begin{proposition}\label{proposition:compare-Nygaard-filtrations-special}
Let $R$ be an animated commutative algebra over the perfectoid ring $\Z_p^{\cyc}$, and write
$\Z_p^{\cyc} = A/I$ for a perfect prism $(A,I)$. For every integer $n$, the comparison map
$$\Fil^{\bullet}_{\Nyg} \Prism_{R}\{n\} \rightarrow\Fil^{\bullet}_{\Nyg} F^{\ast}(\Prism_{R/A}\{n\})$$ 
of Construction \ref{construction:compare-Nygaard-filtrations} is an isomorphism (in the $\infty$-category $\DFiltI(\Z_p)$).
\end{proposition}

\begin{proof}[Proof of Theorem \ref{theorem:compare-Nygaard-filtrations} from Proposition \ref{proposition:compare-Nygaard-filtrations-special}]
Let $(A,I)$ be a perfect prism and let $R$ be an animated commutative algebra over the quotient ring $\overline{A} = A/I$. We wish to show that, for
every pair of integers $m$ and $n$, the comparison map
$$\Fil^{m}(\theta_{R}): \Fil^{m}_{\Nyg} \Prism_{R}\{n\} \rightarrow \Fil^{m}_{\Nyg} F^{\ast}( \Prism_{R/A}\{n\} )$$
of Construction \ref{construction:compare-Nygaard-filtrations} is an isomorphism. Without loss of generality we may assume $m \geq 0$. The proof proceeds by induction on $m$.
In the case $m = 0$, we observe that $\Fil^{0}(\theta_{R})$ can be obtained by composing the isomorphism $\Prism_{R}\{n\} \simeq \Prism_{R/A}\{n\}$
of Proposition \ref{proposition:absolute-vs-relative} with the tautological Frobenius-semilinear map $\Prism_{R/A}\{n\} \rightarrow F^{\ast} \Prism_{R/A}\{n\}$,
which is also an isomorphism by virtue of our assumption that $(A,I)$ is perfect. To carry out the inductive step, it will suffice to show that the associated graded map
$$ \gr^{m}(\theta_{R}): \gr^{m}_{\Nyg} \Prism_{R}\{n\} \rightarrow \gr^{m}_{\Nyg} F^{\ast}( \Prism_{R/A}\{n\} ) \simeq \Fil_{m}^{\conj} \overline{\Prism}_{R/A}\{m\}$$
is an isomorphism for each $m \geq 0$. 

Let $\CAlg_{\overline{A}}^{\QSyn}$ denote the category of $p$-quasisyntomic $\overline{A}$-algebras (Definition~\ref{definition:qsyn-topology}). Note that $\CAlg_{ \overline{A} }^{\QSyn}$ contains the category $\Poly_{ \overline{A}}$ of finitely generated polynomial rings over $\overline{A}$. Since the functors $R \mapsto  \gr^{m}_{\Nyg} \Prism_{R}\{n\}$ and $R \mapsto \Fil_{m}^{\conj} \overline{\Prism}_{R/A}\{m\}$ commute with sifted colimits,
they are left Kan extensions of their restrictions to $\CAlg_{\overline{A}}^{\QSyn}$. It will therefore suffice to show that $\gr^{m}( \theta_{R} )$ is an isomorphism for
$R \in \CAlg_{ \overline{A} }^{\QSyn}$. 

Let $\calC \subseteq \CAlg_{ \overline{A} }^{\QSyn}$ be the category of $p$-quasisyntomic $\overline{A}$-algebras $R$ which are quasiregular semiperfectoid and which admit
a ring homomorphism $\Z_p^{\cyc} \rightarrow R$. Note that $\calC$ forms a basis for the $p$-quasisyntomic topology on $\CAlg^{\QSyn}_{\overline{A}}$.
Since the functors 
$$ \CAlg_{\overline{A}}^{\QSyn} \rightarrow \widehat{\calD}(\Z_p) \quad \quad R \mapsto \gr^{m}_{\Nyg} \Prism_{R}\{n\}, R \mapsto \Fil_{m}^{\conj} \overline{\Prism}_{R/A}\{m\}$$ 
satisfy $p$-complete faithfully flat descent (Proposition \ref{proposition:gr-nygaard-flat-descent} Variant \ref{variant:conjugate-fpqc-descent}), they are right Kan extensions
of their restrictions to $\calC$. It will therefore suffice to show that $\gr^{m}( \theta_{R} )$ is an isomorphism in the special case where $R \in \calC$
is quasiregular semiperfectoid and admits the structure of an $\Z_p^{\cyc}$-algebra. In this case, we show that each $\Fil^{m}( \theta_R )$ is an isomorphism.
We have a commutative diagram
$$ \xymatrix@R=50pt@C=50pt{ \Fil^{m}_{\Nyg} \Prism_{R}\{n\} \ar[r]^-{ \Fil^{m}( \theta_R ) } \ar[d] & \Fil^{m}_{\Nyg} F^{\ast}( \Prism_{R/A}\{n\} ) \ar[d] \\
\Prism_{R}\{n\} \ar[r]^-{\sim} & F^{\ast}( \Prism_{R/A}\{n\} ), }$$
where the bottom horizontal map is an isomorphism of abelian groups (regarded as chain complexes concentrated in cohomological degree zero).
Moreover, the right vertical map identifies $\Fil^{m}_{\Nyg} F^{\ast}( \Prism_{R/A}\{n\} )$ with the subgroup of
$F^{\ast}( \Prism_{R/A}\{n\} )$ consisting of those elements $x$ whose image under the relative Frobenius map $F^{\ast}(\Prism_{R/A}\{n\}) \rightarrow I^{-n} \Prism_{R/A}\{n\}$
belongs to $I^{m-n} \Prism_{R/A}\{n\}$. We are therefore reduced to showing that the left vertical map identifies $\Fil^{m}_{\Nyg} \Prism_{R}\{n\}$ with the
subgroup of $\Prism_{R}\{n\}$ consisting of those elements $y$ whose image under the absolute Frobenius map $\varphi: \Prism_{R}\{n\} \rightarrow \Prism_{R}^{[-n]}\{n\}$
belongs to $\Prism_{R}^{[m-n]}\{n\}$. This follows from Proposition \ref{proposition:compare-Nygaard-filtrations-special} and Remark \ref{remark:local-qrsp-conclusion},
since $R$ admits the structure of a $\Z_p^{\cyc}$-algebra.
\end{proof}

The proof of Proposition \ref{proposition:compare-Nygaard-filtrations-special} will require some preliminaries.

\begin{notation}
Let $R$ be an animated commutative ring and let $R[x]$ denote the free animated $R$-algebra on one variable. For every object $M \in \widehat{\calD}(R)$,
we write $M \langle x \rangle$ for the $p$-completion of the complex $M[x] = R[x] \otimes^{L}_{R} M$, which we regard as an object of the $p$-complete derived
$\infty$-category $\widehat{\calD}(R[x])$.
\end{notation}

\begin{example}\label{example:differential-form-polynomial-ring}
Let $\overline{A}$ be a commutative ring and let $R$ be an animated commutative $\overline{A}$-algebra. For every integer $m$, we have a canonical isomorphism
$$ L \widehat{\Omega}^{m}_{R / \overline{A}} \langle x \rangle \oplus L \widehat{\Omega}^{m-1}_{R / \overline{A} } \langle x \rangle \rightarrow L \widehat{\Omega}^{m}_{R[x] / \overline{A}},$$ given on the first factor by functoriality and on the second factor by multiplication by the differential form $dx \in \widehat{\Omega}^{1}_{ R[x] / \overline{A} }$.
Here we adopt the convention that $L \widehat{\Omega}^{m}_{R/ \overline{A} }$ vanishes for $m < 0$.
\end{example}

\begin{variant}\label{variant:differential-form-conjugate-filtration}
Let $(A,I)$ be a bounded prism with quotient ring $\overline{A} = A/I$. Let us identify the differential form $dx$ with its image
under the Hodge-Tate comparison isomorphism $\widehat{\Omega}^{1}_{ \overline{A}[x] / \overline{A} } \simeq \mathrm{H}^{1}( \Fil_{1}^{\conj} \overline{\Prism}_{ \overline{A}[x] / A}\{1\} )$.
For every animated commutative $\overline{A}$-algebra $R$, we also identify $dx$ with its image in the cohomology group
$\mathrm{H}^{1}( \Fil_{1}^{\conj} \overline{\Prism}_{R[x] / A}\{1\} )$. Then, for every integer $m$, we have a canonical isomorphism
$$ \Fil_{m}^{\conj} \overline{\Prism}_{ R / A } \langle x \rangle
\oplus (\Fil_{m-1}^{\conj} \overline{\Prism}_{R/A})[-1]\{-1\} \langle x \rangle \rightarrow \Fil_{m}^{\conj} \overline{\Prism}_{R[x]/A},$$
given on the first factor by functoriality and on the second factor by multiplication by $dx$. This follows by induction on $m$, using Example \ref{example:differential-form-polynomial-ring}
and the Hodge-Tate comparison isomorphism.
\end{variant}

\begin{variant}\label{variant:differential-form-diffracted-Hodge}
Let $R$ be an animated commutative ring, and let us abuse notation by identifying the differential form $dx$ with its image under
the isomorphism $\widehat{\Omega}^{1}_{R[x]} \simeq \mathrm{H}^{1}( \Fil_{1}^{\conj} \widehat{\Omega}^{\DHod}_{R} )$.
Then we have a canonical isomorphism
$$ \Fil_{m}^{\conj} \widehat{\Omega}^{\DHod}_{R} \langle x \rangle \oplus 
(\Fil_{m-1}^{\conj} \widehat{\Omega}^{\DHod}_{R})[-1] \langle x \rangle \rightarrow \Fil_{m}^{\conj} \widehat{\Omega}^{\DHod}_{R[x]},$$
given on the first factor by functoriality and on the second factor by multiplication by $dx$. To prove this, we may assume without loss of generality
that $R$ is a $\Z_p$-algebra (replacing $R$ by its $p$-completion if necessary), in which case the assertion is a special case of Variant \ref{variant:differential-form-conjugate-filtration}
(see Example \ref{example:diffracted-Hodge-as-relative-prismatic}).
\end{variant}

\begin{variant}\label{variant:differential-form-Nygaard}
Let $R$ be an animated commutative ring, and let us abuse notation by identifying the differential form $dx$ with an element of $\mathrm{H}^{1}( \gr^{1}_{\Nyg} \Prism_{R[x]} )$
using the isomorphism of $\gr^{1}_{\Nyg} \Prism_{R[x]} \simeq L \widehat{\Omega}^{1}_{R[x]}[-1]$ supplied by Proposition \ref{proposition:derived-descent-absolute-Nygaard}.
For every integer $n$, we have a canonical isomorphism
$$ (\gr^{m}_{\Nyg} \Prism_{R}\{n\}) \langle x \rangle \oplus (\gr^{m-1}_{\Nyg} \Prism_{R}\{n\})[-1] \langle x \rangle \rightarrow \gr^{m}_{\Nyg} \Prism_{R[x]}\{n\},$$
given on the first factor by functoriality and on the second factor by multiplication by $dx$. This follows by combining Variant \ref{variant:differential-form-diffracted-Hodge}
with the fiber sequence of Remark \ref{remark:Nygaard-associated-graded}.
\end{variant}

Let $( \Z_p[[ \slashp ]], ( \slashp) )$ be the prism introduced in Notation \ref{notation:reduced-q-de-Rham-prism}. Suppose we are given a morphism of prisms
$$ f: ( \Z_p[[ \slashp]], (\slashp) ) \rightarrow (A,I),$$
where $(A,I)$ is a bounded prism having quotient ring $\overline{A} = A/I$. Then $f$ induces a map
$$ \xi: \widehat{\Omega}^{\DHod}_{\overline{A}} \simeq \overline{\Prism}_{ \overline{A} / \Z_p[[\slashp]] } \rightarrow \overline{\Prism}_{ \overline{A} / A} \simeq \overline{A}.$$
The main ingredient in our proof of Proposition \ref{proposition:compare-Nygaard-filtrations-special} is the following calculation:

\begin{lemma}\label{kurbal}
Let $(A,I)$ be a perfect prism equipped with a morphism of prisms $f: ( \Z_p[[\slashp]], (\slashp) ) \rightarrow (A,I)$.
 For every nonnegative integer $m$, the map $$ \Fil_{m}^{\conj} \widehat{\Omega}^{\DHod}_{\overline{A}} \xrightarrow{ (\Theta+m, \xi) } (\Fil_{m-1}^{\conj} \widehat{\Omega}^{\DHod}_{\overline{A}}) \oplus
\overline{A}$$
is an isomorphism between free $\overline{A}$-modules of rank $(m+1)$.
\end{lemma}

\begin{proof}
Let $\mathbf{G}_{m}^{\sharp}$ be the affine group scheme introduced in Notation \ref{notation:G_m-sharp}, let $\calO_{\mathbf{G}_{m}}^{\sharp} \subseteq \Q[ t^{\pm 1} ]$ be its coordinate ring, and let $\mathcal{P} \rightarrow \Spf( \overline{A} )$ be the $\mathbf{G}_{m}^{\sharp}$-torsor described in Example \ref{example:torsor-over-Hodge-Tate}. Since $\overline{A}$ is a perfectoid ring, we can identify the diffracted Hodge complex $\widehat{\Omega}^{\DHod}_{ \overline{A} }$ with the coordinate ring of $\mathcal{P}$ (Example \ref{example:diffracted-Hodge-of-perfectoid}). The map $\xi: \widehat{\Omega}^{\DHod}_{\overline{A}} \rightarrow \overline{A}$ is a ring homomorphism which we can regard as a trivialization
of the torsor $\mathcal{P}$ (see Proposition \ref{proposition:old-diagram}), and therefore determines an isomorphism of commtutative rings
$$u: \widehat{\Omega}^{\DHod}_{\overline{A}} \simeq \overline{A} \widehat{\otimes} \calO_{\mathbf{G}_{m}}^{\sharp},$$
carrying the Sen operator $\Theta$ to the differential operator $t \frac{ \partial}{\partial t}$.
It follows that, for every integer $n \geq 0$, there is a unique element $\gamma_{n} \in \widehat{\Omega}^{\DHod}_{\overline{A}}$
satisfying $u( \gamma_n ) = \frac{ (1 - t^{-1} )^{n}}{n!}$, and that every element of $\widehat{\Omega}^{\DHod}_{\overline{A}}$
can be written uniquely as an infinite sum $\sum_{n \geq 0} c_{n} \gamma_{n}$ for some sequence of elements $\{ c_n \in \overline{A} \}_{n \geq 0}$ which converges $p$-adically to zero. In terms of this representation, the Sen operator $\Theta$ and the homomorphism $\xi$ are described concretely by the formulae
$$ \Theta( \sum_{n \geq 0} c_{n} \gamma_{n} ) = \sum_{n \geq 0} ( c_{n+1} - n c_n ) \gamma_n \quad \quad 
\xi( \sum_{n \geq 0} c_n \gamma_n) = c_0.$$ 
We will prove the following assertion for $n \geq 0$:
\begin{itemize}
\item[$(\ast_n)$] The complex $\Fil_{n}^{\conj} \widehat{\Omega}^{\DHod}_{\overline{A}}$ can be identified with
the $\overline{A}$-submodule of $\widehat{\Omega}^{\DHod}_{\overline{A}}$ generated by the elements $\{ \gamma_i \}_{0 \leq i \leq n}$.
\end{itemize}
Assuming $(\ast_{m})$ and $(\ast_{m-1})$, we can regard the sequences $(\gamma_{m}, \cdots, \gamma_{0})$ and $(\gamma_{m-1}, \gamma_{m-2}, \cdots, \gamma_0)$
as bases for the $\overline{A}$-modules $\Fil_{m}^{\conj} \widehat{\Omega}^{\DHod}_{\overline{A}}$ and $\Fil_{m-1}^{\conj} \widehat{\Omega}^{\DHod}_{\overline{A}}$,
respectively. Lemma \ref{kurbal} then follows by observing that, with respect to these bases, the map
$$ \Fil_{m}^{\conj} \widehat{\Omega}^{\DHod} \xrightarrow{ (\Theta+m, \xi) } (\Fil_{m-1}^{\conj} \widehat{\Omega}^{\DHod}_{R}) \oplus
\overline{A}$$
is represented by the lower-triangular matrix
$$ \begin{bmatrix}
    1 & 0 & 0 & 0 & \cdots & 0 & 0 \\
    1 & 1 & 0 & 0 & \cdots & 0 & 0 \\
    0 & 2 & 1 & 0 & \cdots & 0 & 0 \\
    0 & 0 & 3 & 1 & \cdots & 0 & 0 \\
    \cdots & \cdots & \cdots & \cdots & \cdots & \cdots & \cdots \\
    0 & 0 & 0 & 0 & \cdots & m & 1 
    \end{bmatrix}$$

Our proof of $(\ast_n)$ proceeds by induction on $n$. Assume that $n \geq 0$ and that, if $n > 0$, then
condition $(\ast_{n-1} )$ is satisfied. Our assumption that $\overline{A}$ is a perfectoid ring guarantees that,
for every integer $d \geq 0$, the complex $\gr_{d}^{\conj} \widehat{\Omega}^{\DHod}_{\overline{A}} \simeq L\widehat{\Omega}^{d}_{\overline{A}}[-d]$ is a free $\overline{A}$-module of rank $1$.
It follows that the conjugate filtration on $\widehat{\Omega}^{\DHod}_{\overline{A}}$ admits a (noncanonical) splitting.
In particular $\Fil_{n}^{\conj} \widehat{\Omega}^{\DHod}_{\overline{A}}$ can be identified with a $\overline{A}$-submodule
of $\widehat{\Omega}^{\DHod}_{\overline{A}}$, concentrated in cohomological degree zero, and the
quotient $\widehat{\Omega}^{\DHod}_{\overline{A}} / \Fil^{\conj}_{n} \widehat{\Omega}^{\DHod}_{\overline{A}}$ is the $p$-completion of
a free $\overline{A}$-module of infinite rank. Choose an element
$x \in \Fil_{n}^{\conj} \widehat{\Omega}^{\DHod}_{\overline{A}}$ whose image is a generator of
$\gr_{n}^{\mathrm{\conj}} \widehat{\Omega}^{\DHod}_{ \overline{A} }$, and write $x = \sum_{d \geq 0} c_d \gamma_d$.
By virtue of $(\ast_{n-1})$, we may assume (modifying $x$ if necessary) that the coefficients $c_{d}$ vanish for $d < n$.
Since the Sen operator $\Theta$ acts by $(-n)$ on $\gr_{n}^{\mathrm{\conj}} \widehat{\Omega}^{\DHod}_{\overline{A}}$, we have
$$ (\Theta+n)(x) = \sum_{d \geq 0} (c_{d+1} + (n-d)c_{d}) \gamma_{d} \in \Fil_{n-1}^{\conj} \widehat{\Omega}^{\DHod}_{\overline{A}}.$$
By virtue of $(\ast_{n-1})$, this is equivalent to the statement that the coefficients $c_{d+1} + (n-d)c_{d}$ vanish for $d \geq n$.
It follows by induction on $d$ that we have $c_{d} = 0$ for $d > n$, so that $x = c_{n} \gamma_{n}$. Since
multiplication by $x$ induces a monomorphism $\overline{A} \hookrightarrow \widehat{\Omega}^{\DHod}_{\overline{A}}$,
the element $c_n \in \overline{A}$ is not a zero-divisor. It follows that $c_n$ is also not a zero divisor on the
topologically free $A$-module $M = \widehat{\Omega}^{\DHod}_{\overline{A}} / \Fil_{n}^{\conj}  \widehat{\Omega}^{\DHod}_{\overline{A}}$.
By construction, the image of $\gamma_n$ in $M$ is annihilated by $c_n$: the element $c_n \gamma_n = x$ lies in $\Fil_n^{\conj}  \widehat{\Omega}^{\DHod}_{\overline{A}}$ by assumption. It follows that the image of $\gamma_n$ must vanish in $M$:
that is, we have $\gamma_n \in \Fil_{n}^{\conj}  \widehat{\Omega}^{\DHod}_{\overline{A}}$. This proves that $c_n$ is a unit in $\overline{A}$,
so that  $$ \Fil^{\conj}_{n} \widehat{\Omega}^{\DHod}_{\overline{A}} = \overline{A}x + \Fil_{n-1}^{\conj} \widehat{\Omega}^{\DHod}_{\overline{A}}$$
is the $\overline{A}$-submodule generated by the elements $\{ \gamma_{d} \}_{0 \leq d \leq n}$, as desired.
\end{proof}

\begin{proof}[Proof of Proposition \ref{proposition:compare-Nygaard-filtrations-special}]
Let $(A,I)$ be a perfect prism with quotient ring $A/I = \Z_p^{\cyc}$ and let $R$ be an animated commutative $\Z_p^{\cyc}$-algebra. 
We wish to show that, for every pair of integers $m$ and $n$, the comparison map
$$\Fil^{m}(\theta_{R}): \Fil^{m}_{\Nyg} \Prism_{R}\{n\} \rightarrow \Fil^{m}_{\Nyg} F^{\ast}( \Prism_{R/A}\{n\} )$$
of Construction \ref{construction:compare-Nygaard-filtrations} is an isomorphism. Arguing as in the proof of
Theorem \ref{theorem:compare-Nygaard-filtrations}, we are reduced to showing that the induced map
$$ \gr^{m}( \theta_R): \gr^{m}_{\Nyg} \Prism_{R}\{n\} \rightarrow  \Fil_{m}^{\conj} \overline{\Prism}_{R/A}\{m\}$$
is an isomorphism. Since the construction $R \mapsto \gr^{m}( \theta_{R} )$ commutes with sifted colimits (see Remark \ref{remark:graded-Nygaard-sifted-colimit}). It
will therefore suffice to show that $\gr^{m}( \theta_{R} )$ is an isomorphism in the case where $R = \Z_p^{\cyc}[ x_1, \ldots, x_d ]$ is a finitely generated polynomial ring
over $\Z_p^{\cyc}$. We now proceed by induction on $d$. To carry out the inductive step, we observe that
Variants \ref{variant:differential-form-conjugate-filtration} and \ref{variant:differential-form-Nygaard} allow us to identify $\gr^{m}(\theta_{R[x]})$ with a direct sum of maps
$$ (\gr^{m}_{\Nyg} \Prism_{R}\{n\}) \langle x \rangle \rightarrow \Fil_{m}^{\conj} \overline{\Prism}_{R/A}\{m\} \langle x \rangle$$
$$ (\gr^{m-1}_{\Nyg} \Prism_{R}\{n\})[-1] \langle x \rangle \rightarrow \Fil_{m-1}^{\conj} \overline{\Prism}_{R/A}\{m\}[-1] \langle x \rangle$$
induced by $\gr^{m}(\theta_R)$ and $\gr^{m-1}( \theta_{R} )$, respectively. Consequently, if $\gr^{m}(\theta_R)$ and $\gr^{m-1}( \theta_{R} )$ are isomorphisms,
then $\gr^{m}(\theta_{R[x]})$ is an isomorphism. We can therefore assume without loss of generality that $d=0$, so that $R = \Z_p^{\cyc}$ is a perfectoid ring.
We may also assume that $m \geq 0$ (otherwise there is nothing to prove), so that $\Fil_{m}^{\conj} \overline{\Prism}_{R/A} = \overline{\Prism}_{R/A}$ can be identified with $\Z_p^{\cyc}$.
Note that there exists a map of prisms $( \Z_p[[\slashp]], (\slashp) ) \rightarrow (A,I)$. Using Example \ref{example:diffracted-Hodge-as-relative-prismatic}, we see
that the map $\gr^{m}( \theta_{\Z_p^{\cyc}} )$ factors as a composition
$$ \xymatrix@R=20pt@C=50pt{ \gr^{m}_{\Nyg} \Prism_{\Z_p^{\cyc} }\{n\} \ar[d]^{\sim} \\
\fib( \Theta+m: \Fil_{m}^{\conj} \widehat{\Omega}^{\DHod}_{\Z_p^{\cyc} } \rightarrow \Fil_{m-1}^{\conj} \widehat{\Omega}^{\DHod}_{\Z_p^{\cyc}}  ) \ar[d] \\
 \Fil_{m}^{\conj} \widehat{\Omega}^{\DHod}_{\Z_p^{\cyc} } \ar[d]^{\xi} \\
 \Z_p^{\cyc}, }$$
where $\xi$ is the map appearing in the statement of Lemma \ref{kurbal}. It follows that $\gr^{m}( \theta_{\Z_p}^{\cyc} )$ is an isomorphism if and only if the map
$$ \Fil_{m}^{\conj} \widehat{\Omega}^{\DHod}_{\Z_p^{\cyc}} \xrightarrow{ (\Theta+m, \xi) } (\Fil_{m-1}^{\conj} \widehat{\Omega}^{\DHod}_{\Z_p^{\cyc}}) \oplus
\Z_p^{\cyc}$$
is an isomorphism, which is a special case of Lemma \ref{kurbal}.
\end{proof}

\subsection{The Frobenius on Absolute Prismatic Cohomology}\label{subsection:Frobenius}

For every animated commutative ring $R$, the absolute prismatic complex $\Prism_{R}$ is equipped with a canonical {\em Frobenius} endomorphism $\varphi: \Prism_{R} \rightarrow \Prism_{R}$. To analyze the properties of this morphism, it will be convenient to factor it as a composition
$$ \Prism_{R} = \RGamma( \WCart, \mathscr{H}_{\Prism}(R) ) \rightarrow
\RGamma( \WCart, F^{\ast} \mathscr{H}_{\Prism}(R) ) \xrightarrow{ \RGamma( \WCart, \Phi) } \RGamma( \WCart, \mathscr{H}_{\Prism}(R) ) = 
\Prism_{R},$$
where the first map is given by pullback along the Frobenius endomorphism of the Cartier-Witt stack itself (Construction \ref{construction:Frobenius-on-stack}), and the second is induced by a morphism of quasi-coherent complexes $\Phi: F^{\ast}( \mathscr{H}_{\Prism}(R) ) \rightarrow
\mathscr{H}_{\Prism}(R)$ which we will refer to as the {\em relative Frobenius}.

\begin{construction}[The Relative Frobenius]\label{construction:relative-Frobenius-on-sheaf}
Let $R$ be an animated commutative ring. For every bounded prism $(A,I)$ with quotient ring $\overline{A} = A/I$,
we can view the derived tensor product $\overline{A} \otimes^{L} R$ as an animated commutative $\overline{A}$-algebra.
Let 
$$ \Fil^{\bullet}( \varphi): \Fil^{\bullet}_{\Nyg} F^{\ast} \Prism_{R/A} \rightarrow I^{\bullet} \Prism_{R/A}$$
be the morphism of filtered complexes appearing in the statement of Proposition \ref{proposition:relative-nygaard-characterization}.
Note that this morphism depends functorially on the prism $(A,I)$, and therefore determines a morphism
of filtered complexes on the Cartier-Witt stack $\WCart$. Passing to filtration degree zero, we obtain a map
$\Phi: F^{\ast} \mathscr{H}_{\Prism}(R) \rightarrow \mathscr{H}_{\Prism}(R)$, which we will refer to as the
{\em relative Frobenius morphism}.
\end{construction}

\begin{remark}\label{remark:relative-Frobenius-sheaf-filtered}
Let $R$ be an animated commutative ring and let
$$\Fil^{\bullet}( \Phi): \Fil^{\bullet}_{\Nyg} F^{\ast} \mathscr{H}_{\Prism}(R) \simeq \mathscr{I}^{\bullet} \mathscr{H}_{\Prism}(R)$$
be as in Construction \ref{construction:relative-Frobenius-on-sheaf}. For every integer $m$, we obtain an induced map
$$ \gr^{m}(\Phi): \gr^{m}_{\Nyg} F^{\ast} \mathscr{H}_{\Prism}(R) \rightarrow \iota_{\ast} \mathscr{H}_{\overline{\Prism}}(R)\{m\},$$
where $\iota: \WCart^{\mathrm{HT}} \hookrightarrow \WCart$ is the inclusion of the Hodge-Tate divisor. Globalizing Remark \ref{remark:polynomial-to-smooth}, we see that this map factors through an isomorphism $\gr^{m}_{\Nyg} F^{\ast} \mathscr{H}_{\Prism}(R) \simeq \iota_{\ast} \Fil_{m}^{\conj} \mathscr{H}_{\overline{\Prism}}(R)\{m\}$ which (by a slight abuse of notation) we will also denote by $\gr^{m}(\Phi)$.
\end{remark}

\begin{remark}\label{ablos}
Let $R$ be an animated commutative ring. Globalizing Corollary \ref{corollary:increasing-filtration-on-prism}, we see that the relative Frobenius map of Construction \ref{construction:relative-Frobenius-on-sheaf} exhibits $\mathscr{H}_{\Prism}(R)$ as a colimit of the diagram
$$ F^{\ast} \mathscr{H}_{\Prism}(R) \rightarrow \mathscr{I}^{-1} \Fil^{1}_{\Nyg} F^{\ast} \mathscr{H}_{\Prism}(R) \rightarrow
\mathscr{I}^{-2} \Fil^{2}_{\Nyg} F^{\ast} \mathscr{H}_{\Prism}(R) \rightarrow \cdots$$
in the $\infty$-category $\calD( \WCart )$. Here each of the successive quotients
$$ \mathscr{I}^{-n} \Fil^{n}_{\Nyg} F^{\ast} \mathscr{H}_{\Prism}(R) / \mathscr{I}^{1-n} \Fil^{n-1}_{\Nyg} F^{\ast} \mathscr{H}_{\Prism}(R)$$
can be identified with $\iota_{\ast}( \Fil^{n}_{\Hodge} \widehat{\dR}_{R} \otimes \calO_{ \WCart^{\mathrm{HT}} }\{-n\} )$,
where $\iota: \WCart^{\mathrm{HT}} \hookrightarrow \WCart$ denotes the inclusion of the Hodge-Tate divisor
(see Remark \ref{remark:Nygaard-vs-hodge-globalized}).
\end{remark}

\begin{example}\label{example:reverse-Nygaard-filtration-finite}
Let $k$ be a perfect field of characteristic $p$ and let $R$ be a commutative ring which is $p$-completely smooth of relative dimension $\leq d$ over $W(k)$. Combining Remark \ref{ablos} with Corollary \ref{respit}, we deduce that the map
$$\Fil^{n}( \Phi): \Fil^{n}_{\Nyg} F^{\ast} \mathscr{H}_{\Prism}(R) \simeq \mathscr{I}^{n} \mathscr{H}_{\Prism}(R)$$
is an isomorphism for $n \geq d$. In particular, the relative Frobenius map $\Phi: F^{\ast} \mathscr{H}_{\Prism}(R) \rightarrow \mathscr{H}_{\Prism}(R)$ factors as a finite composition
$$ F^{\ast} \mathscr{H}_{\Prism}(R)  \rightarrow  \mathscr{I}^{-1} \Fil^{1}_{\Nyg} F^{\ast} \mathscr{H}_{\Prism}(R)
\rightarrow \cdots \rightarrow
 \mathscr{I}^{-d} \Fil^{d}_{\Nyg} F^{\ast} \mathscr{H}_{\Prism}(R) \xrightarrow{\sim} \mathscr{H}_{\Prism}$$
of morphisms whose cofibers are described in Remark \ref{ablos}.
\end{example}

\begin{notation}\label{notation:Frobenius-on-absolute}
Let $R$ be an animated commutative ring and let $n$ be an integer. Tensoring the relative Frobenius
$\Fil^{n}( \Phi): \Fil^{n}_{\Nyg} F^{\ast} \mathscr{H}_{\Prism}(R) \rightarrow \mathscr{I}^{n} \otimes \mathscr{H}_{\Prism}(R)$
with the isomorphism $F^{\ast}(\calO_{\WCart}\{n\}) \simeq \mathscr{I}^{-n}\{n\}$ of Remark \ref{remark:Frobenius-pullback-on-WCart},
we obtain a morphism 
$$ \Phi\{n\}: \Fil^{n}_{\Nyg} F^{\ast}(\mathscr{H}_{\Prism}(R)\{n\}) \rightarrow \mathscr{H}_{\Prism}(R)\{n\}$$
of quasi-coherent sheaves on $\WCart$. Passing to global sections over the Cartier-Witt stack, we obtain a map
$$ \RGamma( \WCart, \Phi\{n\} ): \RGamma( \WCart, \Fil^{n}_{\Nyg} F^{\ast}(\mathscr{H}_{\Prism}(R)\{n\}) \rightarrow \Prism_{R}\{n\}.$$
We let $\varphi\{n\}: \Fil^{n}_{\Nyg} \Prism_{R}\{n\} \rightarrow \Prism_{R}\{n\}$ denote the composition of
$\RGamma( \WCart, \Phi|\{n\} )$ with the tautological map 
$\Fil^{n}_{\Nyg} \Prism_{R}\{n\} \rightarrow \RGamma( \WCart, , \Fil^{n}_{\Nyg} F^{\ast}(\mathscr{H}_{\Prism}(R)\{n\})$.
In the special case $n=0$, we can regard $\varphi\{n\}$ as a morphism from the absolute prismatic complex $\Prism_{R}$
to itself, which we will denote by $\varphi: \Prism_{R} \rightarrow \Prism_{R}$.
\end{notation}

\begin{example}\label{example:describe-Frobenius-smooth}
Let $k$ be a perfect field of characteristic $p$ and let $R$ be a commutative ring which is $p$-completely smooth of relative dimension $\leq d$ over $W(k)$. Fix an integer $n$. For $n > d$, it follows from Example \ref{example:reverse-Nygaard-filtration-finite} and Remark \ref{remark:fiber-sequence-for-Nygaard} that the map $\varphi\{n\}: \Fil^{n}_{\Nyg} \Prism_{R}\{n\} \rightarrow \Prism_{R}\{n\}$ is an isomorphism.
For $n \leq d$, it admits a finite factorization
\begin{eqnarray*}
\Fil^{n}_{\Nyg} \Prism_{R}\{n\} & \xrightarrow{ \varphi_n } & \RGamma( \WCart, \Fil^{n}_{\Nyg} F^{\ast}( \mathscr{H}_{\Prism}(R)\{n\}) ) \\
& \xrightarrow{ \varphi_{n+1} } &  \RGamma( \WCart, \mathscr{I}^{-1} \Fil^{n+1}_{\Nyg} F^{\ast}( \mathscr{H}_{\Prism}(R)\{n\}) ) \\
& \rightarrow & \cdots \\
& \xrightarrow{\varphi_d} & \RGamma( \WCart, \mathscr{I}^{n-d} \Fil^{d}_{\Nyg} F^{\ast}( \mathscr{H}_{\Prism}(R)\{n\} ) ) \\
& \simeq & \RGamma( \WCart, \mathscr{H}_{\Prism}(R)\{n\} ) \\
& = & \Prism_{R}\{n\}.
\end{eqnarray*}
Here the cofiber of $\varphi_{n}$ can be identified with $\Fil^{n}_{\Hodge} \widehat{\dR}_{R}[-1]$ (Remark \ref{remark:fiber-sequence-for-Nygaard}),
and for $n < m \leq d$ we have isomorphisms
\begin{eqnarray*}
\cofib( \varphi_m) & \simeq & \RGamma( \WCart, \mathscr{I}^{n-m} \Fil^{m}_{\Nyg} F^{\ast}(\mathscr{H}_{\Prism}(R)\{n\} ) /
\mathscr{I}^{n+1-m} \Fil^{m-1}_{\Nyg} F^{\ast}(\mathscr{H}_{\Prism}(R)\{n\} ) \\
& \simeq & \RGamma( \WCart^{\mathrm{HT}}, \Fil^{m}_{\Hodge} \widehat{\dR}_R \otimes \mathcal{O}_{\WCart^{\mathrm{HT}} }\{n-m\}) \\
& \simeq & (\Fil^{m}_{\Hodge} \widehat{\dR}_{R} ) \otimes^{L} (\Z / (n-m)\Z)[-1]
\end{eqnarray*}
supplied by Remark \ref{ablos} and Proposition \ref{proposition:easy-version}.
\end{example}

\begin{example}
Let $k$ be a perfect field of characteristic $p$. Then Example \ref{example:describe-Frobenius-smooth} supplies
a fiber sequence
$\Prism_{ W(k) } \xrightarrow{\varphi} \Prism_{ W(k) } \rightarrow W(k)[-1]$ in the $\infty$-category $\widehat{\calD}(\Z_p)$. Moreover, the map $\varphi\{n\}: \Fil^{n}_{\Nyg} \Prism_{W(k)}\{n\} \rightarrow \Prism_{W(k)}\{n\}$
is an isomorphism for $n > 0$.
\end{example}

\begin{remark}\label{remark:filtered-Frobenius-absolute}
Let $R$ be an animated commutative ring. For every integer $n$, the morphism $\varphi\{n\}: \Fil^{n}_{\Nyg} \Prism_{R}\{n\}
\rightarrow \Prism_{R}\{n\}$ can be refined to a morphism of filtered complexes
$$ \Fil^{\bullet}( \varphi\{n\} ): \Fil^{\bullet}_{\Nyg} \Prism_{R}\{n\} \rightarrow \Prism^{[\bullet+n]}_{R}\{n\}.$$
For every integer $m$, the associated graded map $\gr^{m}( \varphi\{n\} ): \gr^{m}_{\Nyg} \Prism_{R}\{n\} \rightarrow \overline{\Prism}_{R}\{m\}$ identifies (using Remark \ref{remark:Nygaard-associated-graded}) with the map
$$\fib( \Theta+m: \Fil_{m}^{\conj} \widehat{\Omega}^{\DHod}_{R} \rightarrow \Fil_{m-1}^{\conj} \widehat{\Omega}^{\DHod}_{R} ) 
\rightarrow( \widehat{\Omega}^{\DHod}_{R} )^{\Theta = -m}$$
\end{remark}

The following result is an ``absolute'' counterpart of Corollary \ref{corollary:Beilinson-connective-cover}:

\begin{proposition}\label{proposition:absolute-nygaard-connective-cover}
Let $R$ be a commutative ring which is $p$-torsion-free and suppose that $R/pR$ is a regular Noetherian $\F_p$-algebra.
Then, for every integer $n$, the morphism 
$$ \Fil^{\bullet}(\varphi\{n\} ): \Fil^{\bullet}_{\Nyg} \Prism_{R}\{n\} \rightarrow \Prism_{R}^{[\bullet+n]}\{n\}$$
of Remark \ref{remark:filtered-Frobenius-absolute} exhibits $\Fil^{\bullet}_{\Nyg} \Prism_{R}\{n\}$ as the connective cover of
$\Prism_{R}^{[\bullet+n]}\{n\}$ with respect to the Beilinson t-structure.
\end{proposition}

\begin{proof}
By virtue of Proposition \ref{proposition:smooth-absolute-Nygaard-complete}, the absolute prismatic complex
$\Prism_{R}\{n\}$ is Nygaard-complete. By virtue of Remark \ref{remark:characterize-Beilinson-connective}, it will suffice to show that for every integer $m$, the associated graded map
$$ \gr^{m}(\varphi\{n\}): \gr^{m}_{\Nyg} \Prism_{R}\{n\} \rightarrow \overline{\Prism}_{R}\{m\}$$
exhibits $\gr^{m}_{\Nyg} \Prism_{R}\{n\}$ as the truncation of $\overline{\Prism}_{R}\{m\}$ to cohomological degrees $\leq m$.
Since the cohomology groups of $\gr^{m}_{\Nyg} \Prism_{R}\{n\}$ are concentrated in degrees $\leq m$ (Proposition \ref{proposition:Nygaard-filtration-Beilinson-connective}), this is equivalent to the assertion that the cohomology groups of the cofiber $\cofib( \gr^m(\varphi\{n\}) )$ are concentrated in degrees $> m$. By virtue of Remark \ref{remark:filtered-Frobenius-absolute}, we have a fiber sequence
$$ \cofib( \gr^{m}(\varphi\{n\}) ) \rightarrow
\widehat{\Omega}^{\DHod}_{R} / \Fil_{m}^{\conj} \widehat{\Omega}^{\DHod}_{R} 
\xrightarrow{\Theta+m} \widehat{\Omega}^{\DHod}_{R} / \Fil_{m-1}^{\conj} \widehat{\Omega}^{\DHod}_{R}.$$
It will therefore suffice to show that the complexes $\widehat{\Omega}^{\DHod}_{R} / \Fil_{m}^{\conj} \widehat{\Omega}^{\DHod}_{R}$
and $\widehat{\Omega}^{\DHod}_{R} / \Fil_{m-1}^{\conj} \widehat{\Omega}^{\DHod}_{R}$ are concentrated in cohomological degrees $> m$
and $\geq m$, respectively, which follows from Remark \ref{remark:conjugate-equals-Postnikov1}.
\end{proof}

\begin{remark}\label{remark:Frobenius-on-Nygaard-completion}
Let $R$ be an animated commutative ring. For every integer $n$, the filtered complex
$\Prism_{R}^{[\bullet+n]}\{n\}$ is filtration-complete. It follows that the map
$$\Fil^{\bullet}(\varphi): \Fil^{\bullet}_{\Nyg} \Prism_{R}\{n\} \rightarrow \Prism_{R}^{[\bullet+n]}\{n\}$$
of Remark \ref{remark:filtered-Frobenius-absolute} admits an essentially unique factorization
through a map $\Fil^{\bullet}( \widetilde{\varphi}\{n\}): \Fil^{\bullet}_{\Nyg} \widehat{\Prism}_{R}\{n\} \rightarrow \Prism_{R}^{[\bullet+n]}\{n\}$,
where $\widehat{\Prism}_{R}\{n\}$ denotes the Nygaard-completed prismatic complex of $R$
(Definition \ref{definition:Nygaard-completed-prismatic-complex}). By virtue of Proposition \ref{proposition:universal-of-animated}, the Nygaard-completed absolute prismatic complex $\widehat{\Prism}_{R}\{n\}$ and its Nygaard filtration are determined (up to isomorphism) by the following
properties:
\begin{itemize}
\item The functor 
$$ \CAlg^{\anim} \rightarrow \DFiltIComp(\Z_p) \quad \quad R \mapsto \Fil^{\bullet}_{\Nyg} \widehat{\Prism}_{R}\{n\}$$
commutes with sifted colimits (this is an immediate consequence of Remark \ref{remark:graded-Nygaard-sifted-colimit};
beware that we must compute the colimit in the $\infty$-category of {\em filtration-complete} filtered complexes).

\item When $R$ is a finitely generated polynomial algebra over $\Z$, the filtered complex $\Fil^{\bullet}_{\Nyg} \widehat{\Prism}_{R}\{n\}$ 
identifies (via the Frobenius morphism $\Fil^{\bullet}( \varphi\{n\} )$) with the connective cover of $\Prism_{R}^{[\bullet+n]}\{n\}$, taken
with respect to the Beilinson t-structure.
\end{itemize}
\end{remark}

\subsection{Nygaard Completion}\label{subsection:Nygaard-completion}

Let $R$ be an animated commutative ring and let $n$ be an integer. We will say that the absolute prismatic complex $\Prism_{R}\{n\}$ is
{\it Nygaard-complete} if the limit $\varprojlim_{m} \Fil^{m}_{\Nyg} \Prism_{R}$ vanishes (as an object of the $\infty$-category $\widehat{\calD}(\Z_p)$).

\begin{example}\label{example:non-nygaard-complete}
Let $S$ be the perfect $\F_p$-algebra $k[ x^{1/p^{\infty}} ]$ and let $R$ be the quotient ring $S / (x)$. By virtue of Lemma \ref{lemma:conjugate-filtration}, 
we can identify the absolute prismatic complex $\Prism_{R}$ with the ring $A_{\crys}(R)$ of Construction \ref{construction:Acrys}, defined
as the separated $p$-completion of the ring obtained from $W(S)$ by adjoining divided powers on $x$. Under this identification, 
$\Fil^{m}_{\Nyg} \Prism_{R}$ corresponds to the ideal of $A_{\crys}(R)$ generated by all expressions of the form $p^{a} \frac{ x^{b} }{b!}$ for $a+b \geq m$
(Proposition \ref{proposition:concrete-Nygaard-qrsp}). The formal series $\exp(x) = \sum_{n \geq 0} \frac{ x^{n} }{n!}$
does not converge in $A_{\crys}$ (with respect to the topology defined by the systems of ideals $\{ \Fil^{m}_{\Nyg} A_{\crys}(R) \}_{m \geq 0}$ ),
and therefore represents an element of $\widehat{\Prism}_{R}$ which does not belong to the image of the map $\Prism_{R} \rightarrow \widehat{\Prism}_{R}$.
In particular, the absolute prismatic complex $\Prism_{R}$ is not Nygaard-complete.
\end{example}

\begin{proposition}\label{proposition:smooth-absolute-Nygaard-complete}
Let $R$ be a commutative ring which is $p$-torsion-free, and suppose that the quotient $R/pR$ is a regular Noetherian ring.
Then, for every integer $n$, the prismatic complex $\Prism_{R}\{n\}$ is Nygaard-complete.
\end{proposition}

We begin by proving a ``local'' version of Proposition \ref{proposition:smooth-absolute-Nygaard-complete}.

\begin{lemma}\label{lemma:local-Nygaard-completeness}
Let $R$ be a commutative ring which is $p$-torsion-free, and suppose that the quotient $R/pR$ is a regular Noetherian ring.
Then the complex $F^{\ast} \mathscr{H}_{\Prism}(R)$ is complete with respect to its Nygaard filtration: that is, the limit
$\varprojlim_{n} \Fil^{n}_{\Nyg} F^{\ast} \mathscr{H}_{\Prism}(R)$ vanishes in the $\infty$-category $\calD( \WCart )$.
\end{lemma}

\begin{proof}
Let $(A^{\bullet}, I^{\bullet} )$ be the cosimplicial prism of Notation \ref{notation:simplicial-prism}. It follows
from Lemma \ref{lemma:tot-description} that every quasi-coherent complex $\mathscr{E}$ on $\WCart$
can be recovered as the totalization of the cosimplicial complex $\rho_{ A^{\bullet} \ast } \rho_{A^{\bullet}}^{\ast} \mathscr{E}$.
It will therefore suffice to show that, for every integer $k \geq 0$, the pullback
$\rho^{\ast}_{A^{k}}( F^{\ast} \mathscr{H}_{\Prism}(R) ) \in \calD( \Spf( A^{m} ) )$ is Nygaard-complete. Writing
$\overline{A}^{m}$ for the quotient ring $A^{m} / I^{m}$, we are reduced to showing that the complex $F^{\ast} \Prism_{ (\overline{A}^{m} \otimes R) / A^m }$
is complete with respect to the relative Nygaard filtration of Proposition \ref{proposition:relative-nygaard-characterization}, which follows from
Corollary \ref{corollary:smooth-Nygaard-complete}.
\end{proof}

\begin{remark}
Let $k$ be a perfect field of characteristic $p$, and let $R$ be a $p$-completely smooth $W(k)$-algebra. In this case,
we can give a simpler proof of Lemma \ref{lemma:local-Nygaard-completeness}. If $R$ has relative dimension $\leq d$ over $W(k)$,
then the complex $\Fil^{n}_{\Hodge} \widehat{\dR}_{R} \simeq \Fil^{n}_{\Hodge} \widehat{\dR}_{R/W(k)}$ vanishes for $n > d$.
Applying Remark \ref{remark:Nygaard-vs-hodge-globalized}, we see that the multiplication map
$\mathscr{I} \otimes \Fil^{n}_{\Nyg} F^{\ast} \mathscr{H}_{\Prism}(R) \rightarrow \Fil^{n+1}_{\Nyg} F^{\ast} \mathscr{H}_{\Prism}(R)$ is an
isomorphism for $n \geq d$. It follows that, in degrees $\geq d$, the Nygaard filtration on on $F^{\ast} \mathscr{H}_{\Prism}(R)$ agrees with the $\mathscr{I}$-adic filtration on
$\Fil^{d}_{\Nyg} F^{\ast} \mathscr{H}_{\Prism}(R)$: that is, we have isomorphisms $\Fil^{n}_{\Nyg} F^{\ast}  \mathscr{H}_{\Prism}(R) \simeq \mathscr{I}^{n-d} \Fil^{n-d}_{\Nyg} F^{\ast} \mathscr{H}_{\Prism}(R)$.
\end{remark}

\begin{proof}[Proof of Proposition \ref{proposition:smooth-absolute-Nygaard-complete}]
Let $R$ be a commutative ring which is $p$-torsion-free and for which $R/pR$ is a regular Noetherian $\F_p$-algebra.
Remark \ref{remark:fiber-sequence-for-Nygaard} supplies a fiber sequence of filtered complexes
$$ \Fil^{\bullet}_{\Nyg} \Prism_{R}\{n\} \rightarrow \RGamma( \WCart, \Fil^{\bullet}_{\Nyg} F^{\ast}( \mathscr{H}_{\Prism}(R)\{n\} ) ) 
\rightarrow \Fil^{\bullet}_{\Hodge} \widehat{\dR}_{R}[-1].$$
The middle term is filtration-complete by virtue of Lemma \ref{lemma:local-Nygaard-completeness}. It will therefore suffice to show
that the derived de Rham complex $\widehat{\dR}_{R}$ is complete with respect to its Hodge filtration. This follows from
the isomorphism $\Fil^{\bullet}_{\Hodge} \widehat{\dR}_{R} \rightarrow ( \widehat{\Omega}^{\geq \bullet}_{R}, d)$ of 
Proposition \ref{proposition:derived-to-classical-de-Rham}.
\end{proof}

In cases where the absolute prismatic complex $\Prism_{R}$ is not Nygaard-complete, it will be convenient to consider the following variant of
Construction \ref{construction:absolute-Nygaard-untwisted}:

\begin{definition}\label{definition:Nygaard-completed-prismatic-complex}
Let $R$ be an animated commutative ring. For every pair of integers $m$ and $n$, we let $\Fil^{m}_{\Nyg} \widehat{\Prism}_{R}\{n\}$ denote the limit of the tower
$$ \cdots \rightarrow
 \Fil^{m}_{\Nyg} \Prism_{R}\{n\} / \Fil^{m+2}_{\Nyg} \Prism_{R}\{n\} \rightarrow  \Fil^{m}_{\Nyg} \Prism_{R}\{n\} / \Fil^{m+1}_{\Nyg} \Prism_{R}\{n\} = \gr^{m}_{\Nyg} \Prism_{R}\{n\},$$
 formed in the $\infty$-category $\widehat{\calD}(\Z_p)$. We will denote $\Fil^{0}_{\Nyg} \widehat{\Prism}_{R}\{n\}$ by $\widehat{\Prism}_{R}\{n\}$, which
 we refer to as the {\it Nygaard-completed absolute prismatic complex} of $R$. By construction, there is a comparison map $\Prism_{R}\{n\} \rightarrow \widehat{\Prism}_R\{n\}$, which
 is an isomorphism if and only if the absolute prismatic complex $\Prism_{R}\{n\}$ is Nygaard-complete.

 Allowing $m$ to vary, we obtain an inverse system
 $$ \cdots \rightarrow \Fil^{3}_{\Nyg} \widehat{\Prism}_{R}\{n\} \rightarrow
 \Fil^{2}_{\Nyg} \widehat{\Prism}_{R}\{n\} \rightarrow \Fil^{1}_{\Nyg} \widehat{\Prism}_{R}\{n\} \rightarrow
 \Fil^{0}_{\Nyg} \widehat{\Prism}_{R}\{n\} = \widehat{\Prism}_{R}\{n\},$$
 which we will denote by $\Fil^{\bullet}_{\Nyg} \widehat{\Prism}_{R}$ and refer to as the {\it Nygaard filtration} on $\widehat{\Prism}_{R}$.
The comparison map $\Prism_{R}\{n\} \rightarrow \widehat{\Prism}_{R}\{n\}$ refines to a map of filtered complexes $\Fil^{\bullet}_{\Nyg} \Prism_{R}\{n\} \rightarrow \Fil^{\bullet}_{\Nyg} \widehat{\Prism}_{R}\{n\}$, which induces an isomorphism $\gr^{\bullet}_{\Nyg} \Prism_{R}\{n\} \simeq \gr^{\bullet}_{\Nyg} \widehat{\Prism}_{R}\{n\}$.
\end{definition}

\begin{remark}\label{remark:p-complete-fpqc-Nygaard}
For every integer $n$, the construction $R \mapsto \Fil^{\bullet}_{\Nyg} \widehat{\Prism}_{R}\{n\}$ determines a functor from
the $\infty$-category $\CAlg^{\anim}$ of animated commutative rings to the $\infty$-category $\DFiltIComp(\Z_p)$ of filtered complexes which are $p$-complete and
filtration-complete (see Notation \ref{notation:filtered-derived-I-complete}). This functor commutes with sifted colimits and satisfies $p$-complete fpqc descent
(in particular, it satisfies descent for the \'{e}tale topology). This reduces immediately to the analogous statements for the functors
$R \mapsto \gr^{m}_{\Nyg} \widehat{\Prism}_{R}\{n\} \simeq \gr^{m}_{\Nyg} \Prism_{R}\{n\}$; see Proposition \ref{proposition:gr-nygaard-flat-descent}
and Remark \ref{remark:graded-Nygaard-sifted-colimit}.
\end{remark}

\begin{remark}\label{remark:Nygaard-completed-deRham-comparison}
Let $R$ be an animated commutative ring, let $n$ be an integer, and let
$\Fil^{\bullet}( \gamma_{\Prism}^{\dR}\{n\} ): \Fil^{\bullet}_{\Nyg} \Prism_{R}\{n\} \rightarrow \Fil^{\bullet}_{\Hodge} \widehat{\dR}_{R}$
be the de Rham comparison map (see Construction \ref{construction:absolute-Nygaard-untwisted}). Passing to completions with
respect to the Nygaard and Hodge filtrations, we obtain a map of filtration-complete complexes
$\Fil^{\bullet}( \widehat{\gamma}_{\Prism}^{\dR}\{n\} ): \Fil^{\bullet}_{\Nyg} \widehat{\Prism}_{R}\{n\} \rightarrow \Fil^{\bullet}_{\Hodge} \widehat{\dR}_{R}^{\hc}$. Here $\widehat{\dR}_{R}^{\hc}$ denotes the Hodge completion of the $p$-complete derived de Rham complex
$\widehat{\dR}_{R}$ (see Construction \ref{construction:Hodge-complete-dR}).
\end{remark}

\begin{remark}\label{remark:qrsp-Nygaard-complete}
Let $R$ be a quasiregular semiperfectoid ring. Then, for every pair of integers $m$ and $n$, the complex $\Fil^{m}_{\Nyg} \widehat{\Prism}_{R}\{n\}$ 
can be viewed as an abelian group, concentrated in cohomological degree zero. This follows formally from the analogous statement for the complexes
$\gr^{m}_{\Nyg} \Prism_{R}\{n\}$ (see Corollary \ref{corollary:qrsp-filtered-discrete}).
\end{remark}

\newpage
\section{Applications to Periodic Cyclic Homology}\label{section:periodic-cyclic-homology}

The Nygaard-completed prismatic complex arises naturally in the theory of topological cyclic homology. Let 
$R$ be a $p$-quasisyntomic ring, let $\TC^{-}(R) = \THH(R)^{hS^1}$ denote the topological negative cyclic homology spectrum of $R$, and let $\TC^{-}(R)^{\wedge}_{p}$ denote its $p$-completion. In \cite{BMS2}, it is shown that the spectrum
$\TC^{-}(R)^{\wedge}_{p}$ admits a complete and exhaustive filtration $\Fil^{\bullet}_{\mot} \TC^{-}(R)^{\wedge}_{p}$,
whose successive quotients are described by the formula
$$ \gr^{n}_{\mot} \TC^{-}(R)^{\wedge}_{p} \simeq \Fil^{n}_{\Nyg} \widehat{\Prism}_{R}\{n\}[2n].$$
In \S\ref{subsection:p-complete-motivic}, we review the construction of this filtration and its extension to arbitrary animated commutative rings.
In \S\ref{subsection:non-complete-motivic}, we describe an analogous filtration on the spectrum $\TC^{-}(R)$ itself (Construction \ref{construction:integral-motivic-filtration}), whose associated graded invariants are closely related to the diffracted Hodge complex studied in
\S\ref{subsection:diffracted-Hodge-integral} (Proposition \ref{proposition:compute-gr-gr}).

\subsection{The Tate Filtration}\label{subsection:Tate-filtration}

We begin by recalling the construction of topological periodic cyclic homology and establishing some terminology.

\begin{notation}[Filtered Spectra]
Let $\calD( \Sphere )$ denote the $\infty$-category of spectra. We say that a spectrum $X$ is {\it $p$-complete} if the limit of the diagram
$$ \cdots \xrightarrow{p} X \xrightarrow{p} X \xrightarrow{p} X \xrightarrow{p} X$$
vanishes in the $\infty$-category $\calD(\Sphere)$. Equivalently, $X$ is $p$-complete if each homotopy group $\pi_{n}(X)$ is a $p$-complete
abelian group. We let $\widehat{\calD}(\Sphere)$ denote the full subcategory of $\calD( \Sphere)$ spanned by the $p$-complete spectra.
The inclusion functor $\widehat{\calD}(\Sphere) \hookrightarrow \calD(\Sphere)$ has a left adjoint which carries each spectrum $X$ to its
{\it $p$-completion}, which we will denote by $X^{\wedge}_{p}$. 

A {\it filtered spectrum} is a diagram
$$ \cdots \rightarrow \Fil^{2} X \rightarrow \Fil^{1} X \rightarrow \Fil^{0} X \rightarrow \Fil^{-1} X \rightarrow \Fil^{-2} X \rightarrow \cdots$$
in the $\infty$-category $\calD(\Sphere)$, which we will denote by $\Fil^{\bullet} X$. In this case, we will refer to the colimit
$X = \varinjlim_{m} \Fil^{-m} X$ as the {\it underlying spectrum} of $\Fil^{\bullet} X$, and refer to $\Fil^{\bullet} X$ as a
{\it filtration on $X$}. We also write $\gr^{n} X$ for the spectrum given by the cofiber of the map $\Fil^{n+1} X \rightarrow \Fil^{n} X$.

Filtered spectra can be regarded as objects of the functor $\infty$-category $\Fun( ( \Z, \geq)^{\op}, \calD(\Sphere ) )$, which we
will denote by $\DFilt(\Sphere)$ and refer to as {\it the $\infty$-category of filtered spectra}. This $\infty$-category has several subcategories which will be of interest to us:
\begin{itemize}
\item We say that a filtered spectrum $\Fil^{\bullet} X$ is {\it $p$-complete} if the spectrum $\Fil^{n} X$ is $p$-complete for every integer $n$. We let
$\DFiltI(\Sphere)$ denote the full subcategory of $\DFilt(\Sphere)$ spanned by the $p$-complete filtered spectra (that is, the functor
$\infty$-category $\Fun( ( \Z, \geq)^{\op}, \widehat{\calD}( \Sphere ) )$.

\item We say that a filtered spectrum $\Fil^{\bullet} X$ is {\it filtration-complete} if the limit $\varprojlim_{n} \Fil^{n} X$ vanishes in the $\infty$-category $\calD( \Sphere )$ of
spectra. We let $\DFiltComp(\Sphere)$ denote the full subcategory of $\DFilt(\Sphere)$ spanned by the filtration-complete filtered spectra.

\item We let $\DFiltIComp(\Sphere) = \DFiltComp(\Sphere) \cap \DFiltI(\Sphere)$ denote the full subcategory of $\DFilt(\Sphere)$ spanned by those
filtered spectra which are both $p$-complete and filtration-complete.
\end{itemize}
\end{notation}

\begin{example}
Let $M$ be a complex of abelian groups. In what follows, we will abuse notation by identifying $M$ with the generalized Eilenberg-MacLane spectrum associated to $M$.
This construction determines a functor $\calD(\Z) \rightarrow \calD( \Sphere)$, given restriction of scalars along a map of commutative ring
spectra $\Sphere \rightarrow \pi_0(\Sphere) = \Z$. In particular, if $R$ is an animated commutative ring, then we can view
$\Fil^{\bullet}_{\Nyg} \widehat{\Prism}_{R}\{n\}$ as a filtered spectrum (which is $p$-complete and filtration-complete).
\end{example}

\begin{example}[The Tate Construction]\label{example:Tate-spectrum}
Let $V = \mathbf{C}$ be the complex numbers, viewed as a representation of the circle $S^1 = \{ z \in \mathbf{C}: |z| = 1 \}$. We write $\Sphere(V)$ for the suspension spectrum
of the one-point compactification of $V$. The spectrum $\Sphere(V)$ is equipped with an action of $S^1$ (in the Borel-equivariant sense: that is,
it can be regarded as the fiber of a local system of spectra over the classifying space $BS^1 \simeq \CP^{\infty}$), which is nonequivariantly isomorphic to the $2$-sphere
$\Sphere[2]$. In particular, the spectrum $\Sphere(V)$ is invertible with respect to the smash product $\otimes_{\Sphere}$. For every integer $n$, we write $\Sphere(V)^{n}$ for the $n$-fold smash product of $\Sphere(V)$ with itself. 

Let $X$ be any spectrum equipped with an action of $S^1$. For every integer $n$, we let $\Fil^{n}_{\Tate} X^{tS^1}$ denote the spectrum
$( \Sphere(V)^{-n} \otimes_{\Sphere} X)^{hS^1}$ denote the homotopy fixed points for the induced action of $S^1$ on the smash product of
$X$ with $\Sphere(V)^{-n}$. Note that the inclusion map $\{0\} \rightarrow V$ induces an $S^1$-equivariant morphism $e: \Sphere \rightarrow \Sphere(V)$, which
determines a diagram
\begin{equation}\label{equation:Tate-filtration} \cdots \rightarrow \Fil^{2}_{\Tate} X^{tS^1} \rightarrow \Fil^{1}_{\Tate} X^{tS^1} \rightarrow
 \Fil^{0}_{\Tate} X^{tS^1} \rightarrow  \Fil^{-1}_{\Tate} X^{tS^1} \rightarrow \Fil^{-2}_{\Tate} X^{tS^1} \rightarrow \cdots \end{equation}
We denote the colimit of this diagram by $X^{tS^1}$, which we refer to as the {\it Tate construction} on $X$. The diagram
(\ref{equation:Tate-filtration}) itself can be viewed as a filtered spectrum whose underlying spectrum is $X^{tS^1}$. We denote
this filtered spectrum by $\Fil^{\bullet}_{\Tate} X^{tS^1}$ and refer to it as the {\it Tate filtration} on $X^{tS^1}$. Note that $\Fil^{0}_{\Tate} X^{tS^1}$
is the homotopy fixed point spectrum $X^{hS^1}$.
\end{example}

\begin{remark}
The morphism $e: \Sphere \rightarrow \Sphere(V)$ appearing in Example \ref{example:Tate-spectrum} is nullhomotopic as a morphism of spectra
(though not as a morphism of $S^1$-equivariant spectra). For any spectrum $X$ equipped with an action of $S^1$, it follows that the limit
$$ \varprojlim_{n} \Fil^{n}_{\Tate} X^{tS^1} = \varprojlim_{n} ( \Sphere(V)^{-n} \otimes_{\Sphere} X)^{hS^1} 
\simeq ( \varprojlim_{n} ( \Sphere(V)^{-n} \otimes_{\Sphere} X ) )^{hS^1}$$
vanishes in the $\infty$-category of spectra. In other words, the Tate filtration $\Fil^{\bullet}_{\Tate} X^{tS^1}$ is complete.
\end{remark}

\begin{remark}\label{remark:connectivity-estimate}
Let $X$ be a spectrum equipped with an $S^1$-action. For every integer $n$, there is a canonical isomorphism $\gr^{n}_{\Tate} X^{tS^1} \simeq X[-2n]$.
It follows that, if $X$ is a connective spectrum, then the tautological map $\Fil^{-n}_{\Tate} X^{tS^1} \rightarrow X^{tS^1}$ is 
$(2n+1)$-connective: that is, the map $\pi_{\ast}( \Fil^{-n}_{\Tate} X^{tS^1}  ) \rightarrow \pi_{\ast}( X^{tS^1} )$ is an isomorphism
for $\ast \leq 2n$ and a surjection for $\ast = 2n+1$.
\end{remark}

\begin{remark}
Let $X$ be a spectrum equipped with an action of $S^1$, and let $X^{\wedge}_{p}$ denote its $p$-completion. Then there is a canonical isomorphism of filtered spectra
$$ \Fil^{\bullet}_{\Tate} (X^{\wedge}_{p})^{tS^1} \simeq ( \Fil^{\bullet}_{\Tate} X^{tS^1} )^{\wedge}_{p}.$$
If $X$ is connective, then we can pass to the colimit over the Tate filtration (using the connectivity estimates of Remark \ref{remark:connectivity-estimate}) to obtain
an isomorphism $(X^{\wedge}_{p})^{tS^1} \simeq (X^{tS^1})^{\wedge}_{p}$.

In particular, if $X$ is $p$-complete, then the filtered spectrum $\Fil^{\bullet}_{\Tate} X^{tS^1}$ is also $p$-complete. If $X$ is $p$-complete and connective, then $X^{tS^1}$ is $p$-complete.
\end{remark}

\begin{remark}\label{remark:periodicity}
Let $A$ be a commutative ring spectrum equipped with an action of $S^1$. Suppose that the homotopy groups $\pi_{n}(A)$ vanish when $n$ is odd. 
Then the canonical map
$$ \rho: \pi_{-2} ( \Sphere(V)^{-1} \otimes_{\Sphere} A)^{hS^1} \rightarrow \pi_{-2}( \Sphere(V)^{-1} \otimes_{\Sphere} A) \simeq \pi_0(A)$$
is surjective. We can therefore choose an element $b \in \pi_{-2}( \Sphere(V)^{-1} \otimes_{\Sphere} A)^{hS^1}$ satisfying $\rho(b) = 1$.
The element $b$ can then be viewed as an $S^1$-equivariant $A$-module isomorphism $\Sphere(V) \otimes_{\Sphere} A \simeq A[2]$.
It follows that, for every integer $n$, multiplication by $b^{-n}$ induces an isomorphism
$$ \Fil^{n}_{\Tate} A^{tS^1} = ( \Sphere(V)^{-n} \otimes_{\Sphere} A)^{hS^1} \simeq A^{hS^1}[-2n].$$
Writing $\beta$ for the image of $b$ under the map $\pi_{-2}( \Fil^{1}_{\Tate} A^{tS^1} ) \rightarrow \pi_{-2}( \Fil^0_{\Tate} A^{tS^1} ) = \pi_{-2} A^{hS^1}$,
obtain an identification of $\Fil^{\bullet}_{\Tate} A^{tS^1}$ with the filtered spectrum
$$ \cdots \rightarrow A^{hS^1}[-4] \xrightarrow{\beta} A^{hS^1}[-2] \xrightarrow{\beta} A^{hS^1} \xrightarrow{\beta} A^{hS^1}[2] \xrightarrow{\beta}
A^{hS^1}[4] \rightarrow \cdots,$$
so that the Tate construction $A^{tS^1}$ is given by the localization $A^{hS^1}[ \beta^{-1} ]$.
\end{remark}

\begin{warning}
If $A$ is an associative ring spectrum equipped with an action of $S^1$, then the filtered spectrum $\Fil^{\bullet}_{\Tate} A^{tS^1}$ has the structure
of an associative algebra object of the $\infty$-category $\DFilt(\Sphere)$, which endows $A^{tS^1}$ with the structure of an associative ring spectrum.
If $A$ is a commutative ring spectrum, then $A^{tS^1}$ also inherits the structure of a commutative ring spectrum. Beware that the commutativity is not visible
at the level of the filtered spectrum $\Fil^{\bullet}_{\Tate} A^{tS^1}$: that is, $\Fil^{\bullet}_{\Tate} A^{tS^1}$ need not have the structure of a commutative
algebra object of the $\infty$-category $\DFilt(\Sphere)$ (even in the special case where $A = \Sphere$ is the sphere spectrum, equipped with the trivial action of $S^1$).
\end{warning}

\begin{example}\label{example:definition-of-TP}
Let $R$ be an animated commutative ring, which we identify with its underlying commutative ring spectrum.
We write $\THH(R)$ for the topological Hochschild homology of $R$, which we view as a commutative ring spectrum equipped with an action of $S^1$.
Let $\TP(R)$ denote the Tate construction $\THH(R)^{tS^1}$; we will refer to $\TP(R)$ as the {\it topological periodic cyclic homology of $R$}.
For every integer $n$, we write $\Fil^{n}_{\Tate} \TP(R)$ for the spectrum 
$\Fil^{n}_{\Tate} \THH(R)^{tS^1} = ( \Sphere(V)^{-n} \otimes_{\Sphere} \THH(R) )^{hS^1}$ introduced in Example \ref{example:Tate-spectrum}.
We regard the diagram
$$ \cdots \rightarrow \Fil^{2}_{\Tate} \TP(R) \rightarrow \Fil^{1}_{\Tate} \TP(R) \rightarrow \Fil^{0}_{\Tate} \TP(R) \rightarrow
\Fil^{-1}_{\Tate} \TP(R) \rightarrow \Fil^{-2}( \TP(R) ) \rightarrow \cdots$$
as a filtered spectrum, which we denote by $\Fil^{\bullet}_{\Tate} \TP(R)$ and refer to as the {\it Tate filtration} on $\TP(R)$.
We denote the spectrum $\Fil^{0}_{\Tate} \TP(R) = \THH(R)^{hS^1}$ by $\TC^{-}(R)$, which we refer to as the
{\it topological negative cyclic homology of $R$}.
\end{example}

\subsection{The Motivic Filtration: \texorpdfstring{$p$}{p}-Complete Case}
\label{subsection:p-complete-motivic}

In this section, we review \cite{BMS2}, giving a relationship between Nygaard-completed prismatic cohomology and topological Hochschild homology.
One of the central results of \cite{BMS2} is the following:

\begin{theorem}[\cite{BMS2}]\label{theorem:BMS2-calculation}
Let $R$ be a quasiregular semiperfectoid ring. Then the homotopy groups of the spectrum $\THH(R)^{\wedge}_{p}$ are concentrated in even degrees.
Moreover, there is canonical isomorphism of graded rings
$$ \pi_{\ast}( \TC^{-}(R)^{\wedge}_{p} ) \simeq \begin{cases} \Fil^{n}_{\Nyg} \widehat{\Prism}_R\{n\} & \text{ if } \ast = 2n \\
0 & \text{ otherwise, } \end{cases}$$
depending functorially on $R$.
\end{theorem}

Combining Theorem \ref{theorem:BMS2-calculation} with Remark \ref{remark:periodicity}, we obtain the following:

\begin{corollary}[\cite{BMS2}]\label{corollary:BMS2-calculation-TP}
Let $R$ be a quasiregular semiperfectoid ring. Then there is a canonical isomorphism of graded rings
$$ \pi_{\ast}( \TP(R)^{\wedge}_{p} ) \simeq \begin{cases} \widehat{\Prism}_R\{n\} & \text{ if } \ast = 2n \\
0 & \text{ otherwise, } \end{cases}$$
depending functorially on $R$.
\end{corollary}

We will need the following refinement of Corollary \ref{corollary:BMS2-calculation-TP}:

\begin{corollary}\label{corollary:BMS-shifted}
Let $R$ be a quasiregular semiperfectoid ring. For every pair of integers $m$ and $n$, the map of homotopy groups
$$\pi_{2n} (\Fil^{m}_{\Tate} \TP(R)^{\wedge}_{p} ) \rightarrow \pi_{2n}( \TP(R)^{\wedge}_{p} )$$
is a monomorphism, whose image can be identified with the subgroup 
$\Fil^{n+m}_{\Nyg} \widehat{\Prism}_{R}\{n\} \subseteq \widehat{\Prism}_R\{n\} \simeq \pi_{2n}(\TP(R)^{\wedge}_{p})$. 
Moreover, the homotopy groups of the spectrum $\Fil^{m}_{\Tate} \TP(R)^{\wedge}_{p}$ vanish in odd degrees.
\end{corollary}

\begin{proof}
Using the fiber sequence
$$ \Fil^{1}_{\Tate} \TP(R)^{\wedge}_{p} \rightarrow \TC^{-}(R)^{\wedge}_{p} \rightarrow \THH(R)^{\wedge}_{p}$$
and the connectivity of $\THH(R)$, we deduce that the natural map
$$\pi_{-2}( \Fil^{1}_{\Tate} \TP(R)^{\wedge}_{p} ) \xrightarrow{\sim} \pi_{-2}( \TC^{-}(R)^{\wedge}_{p} )$$
is an isomorphism. Let us henceforth use this isomorphism (and Theorem \ref{theorem:BMS2-calculation})
to identify $\pi_{-2}( \Fil^{1}_{\Tate} \TP(R)^{\wedge}_{p} )$ with $\widehat{\Prism}_{R}\{-1\}$, regarded
as an invertible module over the commutative ring $\widehat{\Prism}_{R} \simeq \pi_0 \TC^{-}(R)^{\wedge}_{p}$.
Note that this module is free (since $R$ is defined over a perfectoid ring); let
$b \in \widehat{\Prism}_{R}\{-1\}$ be a generator. Since the homotopy groups of $\THH(R)^{\wedge}_{p}$
are concentrated in even degrees, the map
$$ \rho: \widehat{\Prism}_{R}\{-1\} = \pi_{-2}( \Fil^{1}_{\Tate} \TP(R)^{\wedge}_{p} ) \rightarrow
\pi_{-2}( \gr^{1}_{\Tate} \TP(R)^{\wedge}_{p} ) \simeq R$$
is surjective. We may therefore assume without loss of generality that $\rho(b) = 1$. We then have a commutative diagram of spectra
$$ \xymatrix@R=50pt@C=50pt{ \Fil^{m}_{\Tate} \TP(R)^{\wedge}_{p} \ar[d]^{b^{-m}}_{\sim} \ar[r] & \TP(R)^{\wedge}_{p} \ar[d]^{b^{-m}}_{\sim} \\
\TC^{-}(R)^{\wedge}_{p}[-2m] \ar[r] & \TP(R)^{\wedge}_{p}[-2m], }$$
where the vertical maps are isomorphisms (see Remark \ref{remark:periodicity}). We are therefore reduced to proving
Corollary \ref{corollary:BMS-shifted} in the case $m=0$, which is clear.
\end{proof}

In \cite{BMS2}, Theorem \ref{theorem:BMS2-calculation} and Corollary \ref{corollary:BMS2-calculation-TP} were applied to construct
{\em motivic filtrations} on the spectra $\TP(R)^{\wedge}_{p}$, $\TC^{-}(R)^{\wedge}_{p}$, and $\THH(R)^{\wedge}_{p}$
in the case where $R$ is a $p$-quasisyntomic commutative ring. These constructions extend formally to the case where $R$ is an arbitrary animated commutative
ring, albeit with worse convergence properties (see Warning \ref{warning:only-filtration-colimit}).

\begin{theorem}[The motivic filtration following \cite{BMS2}]
\label{theorem:BMS2-main}
Let $R$ be an animated commutative ring. Then $\Fil^{\bullet}_{\Tate} \TP(R)^{\wedge}_{p}$ can be realized as the colimit of a diagram
$$ \cdots \rightarrow \Fil^{1}_{\mot} \Fil^{\bullet}_{\Tate} \TP(R)^{\wedge}_{p} \rightarrow \Fil^{0}_{\mot} \Fil^{\bullet}_{\Tate} \TP(R)^{\wedge}_{p} \rightarrow 
 \Fil^{-1}_{\mot} \Fil^{\bullet}_{\Tate} \TP(R)^{\wedge}_{p} \rightarrow \cdots$$
 in the $\infty$-category $\DFiltIComp{\Sphere}$. This realization is determined (up to isomorphism) by the requirement that it depends functorially on $R$
 and satisfies the following additional conditions:
\begin{itemize}
\item[$(1)$] For each integer $n$, the functor
$$ \CAlg^{\anim} \rightarrow \DFiltIComp(\Sphere) \quad \quad R \mapsto \Fil^{n}_{\mot} \Fil^{\bullet}_{\Tate} \TP(R)^{\wedge}_{p}$$
commutes with sifted colimits.

\item[$(2)$] For each integer $n$, the functor
$$ \CAlg^{\anim} \rightarrow \DFiltIComp(\Sphere) \quad \quad R \mapsto \Fil^{n}_{\mot} \Fil^{\bullet}_{\Tate} \TP(R)^{\wedge}_{p}$$
satisfies $p$-complete fpqc descent.

\item[$(3)$] When $R$ is a quasiregular semiperfectoid ring, each of the maps
$$ \Fil^{n}_{\mot} \Fil^{m}_{\Tate} \TP(R)^{\wedge}_{p} \rightarrow \Fil^{m}_{\Tate} \TP(R)^{\wedge}_{p}$$
exhibits $\Fil^{n}_{\mot} \Fil^{m}_{\Tate} \TP(R)^{\wedge}_{p}$ as the Postnikov truncation $\tau^{\leq -2n} \Fil^{m}_{\Tate} \TP(R)^{\wedge}_{p}$.
That is, the homotopy groups of the spectrum $\Fil^{n}_{\mot} \Fil^{m}_{\Tate} \TP(R)^{\wedge}_{p}$ are concentrated in degrees $\geq 2n$,
and the induced map $\pi_{\ast}(\Fil^{n}_{\mot} \Fil^{m}_{\Tate} \TP(R)^{\wedge}_{p}) \rightarrow \pi_{\ast}( \Fil^{m}_{\Tate} \TP(R)^{\wedge}_{p} )$
is an isomorphism for $\ast \geq 2n$.
\end{itemize}
\end{theorem}

We will sketch the proof of Theorem \ref{theorem:BMS2-main} at the end of this section.

\begin{definition}[The Motivic Filtration]
Let $R$ be an animated commutative ring. We will denote the diagram
$$ \cdots \rightarrow \Fil^{1}_{\mot} \Fil^{\bullet}_{\Tate} \TP(R)^{\wedge}_{p} \rightarrow
\Fil^{0}_{\mot} \Fil^{\bullet}_{\Tate} \TP(R)^{\wedge}_{p} \rightarrow \Fil^{-1}_{\mot} \Fil^{\bullet}_{\Tate} \TP(R)^{\wedge}_{p} \rightarrow \cdots$$
of Theorem \ref{theorem:BMS2-main} by $\Fil^{\bullet}_{\mot} \Fil^{\bullet}_{\Tate} \TP(R)^{\wedge}_{p}$ and refer to it as the
{\it motivic filtration} on the filtered spectrum $\Fil^{\bullet}_{\Tate} \TP(R)^{\wedge}_{p}$. 
\end{definition}

\begin{remark}[Uniqueness of the Motivic Filtration]\label{remark:motivic-uniqueness}
The uniqueness assertion of Theorem \ref{theorem:BMS2-main} is straightforward. Condition $(1)$ of Theorem \ref{theorem:BMS2-main}
is equivalent to the assertion that, for each integer $n$, the functor
$$ \CAlg^{\anim} \rightarrow \DFiltIComp(\Sphere) \quad \quad R \mapsto \Fil^{n}_{\mot} \Fil^{\bullet}_{\Tate} \TP(R)^{\wedge}_{p}$$
is a left Kan extension of its restriction to the category $\Poly_{\Z}$ of finitely generated polynomial algebras over $\Z$. In particular,
it is a left Kan extension of its restriction to the category $\CAlg^{\QSyn}$ of $p$-quasisyntomic commutative rings. It follows from condition $(2)$
of Theorem \ref{theorem:BMS2-main} that the functor
$$ \CAlg^{\QSyn} \rightarrow \DFiltIComp(\Sphere) \quad \quad R \mapsto \Fil^{n}_{\mot} \Fil^{\bullet}_{\Tate} \TP(R)^{\wedge}_{p}$$
is a right Kan extension of its restriction to the full subcategory $\CAlg^{\qrsp} \subset \CAlg^{\QSyn}$ of quasiregular semiperfectoid rings,
which is then determined (up to isomorphism) by condition $(3)$.
\end{remark}

\begin{warning}\label{warning:only-filtration-colimit}
Let $R$ be an animated commutative ring. Theorem \ref{theorem:BMS2-main} asserts only that $\Fil^{\bullet}_{\Tate} \TP(R)^{\wedge}_{p}$ can
be recovered as the colimit of the diagram
$$ \cdots \rightarrow \Fil^{0}_{\mot} \Fil^{\bullet}_{\Tate} \TP(R)^{\wedge}_{p} \rightarrow \Fil^{-1}_{\mot} \Fil^{\bullet}_{\Tate} \TP(R)^{\wedge}_{p} \rightarrow 
 \Fil^{-2}_{\mot} \Fil^{\bullet}_{\Tate} \TP(R)^{\wedge}_{p} \rightarrow \cdots$$
in the $\infty$-category $\DFiltIComp(\Sphere)$ of {\em filtration-complete} (and $p$-complete) filtered spectra. More concretely, this is
equivalent to the assertion that for every integer $m$, the canonical map
$$ \varinjlim_{n} \Fil^{-n}_{\mot} \gr^{m}_{\Tate} \TP(R)^{\wedge}_{p} \rightarrow \gr^{m}_{\Tate} \TP(R)^{\wedge}_{p}$$
is an isomorphism (where the colimit is formed in the $\infty$-category $\widehat{\calD}( \Sphere )$ of $p$-complete spectra).
Beware that that the analogous map
$$ \varinjlim_{n} \Fil^{-n}_{\mot} \Fil^{m}_{\Tate} \TP(R)^{\wedge}_{p} \rightarrow \Fil^{m}_{\Tate} \TP(R)^{\wedge}_{p}$$
need not be an isomorphism in general. However, it is an isomorphism in the case where $R$ is $p$-quasisyntomic: see 
Proposition~7.13 of \cite{BMS2}, and Corollary \ref{corollary:exhaustive-level} for an integral counterpart.
\end{warning}

For every animated commutative ring $R$ and every integer $n$, we write $\gr^{n}_{\mot} \Fil^{\bullet}_{\Tate} \TC(R)^{\wedge}_{p}$ for the cofiber of the tautological map
$\Fil^{n+1}_{\mot} \Fil^{\bullet}_{\Tate} \TC(R)^{\wedge}_{p} \rightarrow \Fil^{n}_{\mot} \Fil^{\bullet}_{\Tate} \TC(R)^{\wedge}_{p}$, which we regard as a filtered spectrum.

\begin{theorem}[Associated graded for the motivic filtration]\label{theorem:motivic-gr}
Let $n$ be an integer. For every animated commutative ring, there is a canonical isomorphism of filtered spectra
$$ \gr^{n}_{\mot} \Fil^{\bullet}_{\Tate} \TC(R)^{\wedge}_{p} \simeq \Fil^{\bullet+n}_{\Nyg} \widehat{\Prism}_{R}\{n\}[2n].$$
This isomorphism is determined by the requirement that it depends functorially on $R$ and that, when $R$ is a quasiregular semiperfectoid ring,
it reduces to the isomorphism of filtered abelian groups
$$\pi_{2n} ( \Fil^{\bullet}_{\Tate} \TC(R)^{\wedge}_{p} ) \simeq \Fil^{\bullet+n} \widehat{\Prism}_{R}\{n\}$$
described in Corollary \ref{corollary:BMS-shifted}.
\end{theorem}

\begin{proof}[Proof of Theorem \ref{theorem:motivic-gr} from Theorem \ref{theorem:BMS2-main}]
As in Remark \ref{remark:motivic-uniqueness}, we observe that the functor
$$ \CAlg^{\anim} \rightarrow \DFiltIComp(\Sphere) \quad \quad R \mapsto \gr^{n}_{\mot} \Fil^{\bullet}_{\Tate} \TC(R)^{\wedge}_{p}$$
is a left Kan of its restriction to the category $\CAlg^{\QSyn}$, which is a right Kan extension of its restriction to the subcategory $\CAlg^{\qrsp} \subseteq \CAlg^{\QSyn}$
of quasiregular semiperfectoid rings. It will therefore suffice to show that the functor $R \mapsto \Fil^{\bullet}_{\Nyg} \widehat{\Prism}_{R}\{n\}$
has the same properties, which follows from Remark \ref{remark:p-complete-fpqc-Nygaard}.
\end{proof}

\begin{example}[The Motivic Filtration on $\THH$]\label{example:THH-filtration}
Let $R$ be an animated commutative ring. For every integer $n$, let us write $\Fil^{n}_{\mot} \THH(R)^{\wedge}_{p}$ for the spectrum $\Fil^{n}_{\mot} \gr^{0}_{\Tate} \TP(R)^{\wedge}_{p}$.
We then obtain a diagram
$$ \cdots \rightarrow \Fil^{1}_{\mot} \THH(R)^{\wedge}_{p} \rightarrow
 \Fil^{0}_{\mot} \THH(R)^{\wedge}_{p} \rightarrow \Fil^{-1}_{\mot} \THH(R)^{\wedge}_{p} \rightarrow \cdots,$$
which we will denote by $\Fil^{\bullet}_{\mot} \THH(R)^{\wedge}_{p}$ and refer to it as the {\it motivic filtration} on the spectrum $\THH(R)^{\wedge}_{p}$.
Note that Theorem \ref{theorem:motivic-gr} supplies isomorphisms $$\gr^{n}_{\mot} \THH(R)^{\wedge}_{p} \simeq \gr^{n}_{\Nyg} \Prism_{R}\{n\}[2n].$$ 
In particular, $\gr^{n}_{\mot} \THH(R)^{\wedge}_{p}$ vanishes for $n < 0$, so that we have $$\Fil^{n}_{\mot} \THH(R)^{\wedge}_{p} \simeq \THH(R)^{\wedge}_{p}$$ for all $n \leq 0$. More generally, for every integer $m$, the natural map
$\Fil^{n}_{\mot} \gr^{m}_{\Tate} \TP(R) \rightarrow \gr^{m}_{\Tate} \TP(R)$ is an isomorphism for $n \leq -m$.
\end{example}

\begin{example}[The Motivic Filtration on $\TC^{-}$]\label{example:motivic-TC-minus}
Let $R$ be an animated commutative ring. For every integer $n$, let us write $\Fil^{n}_{\mot} \TC^{-}(R)^{\wedge}_{p}$ for the spectrum $\Fil^{n}_{\mot} \Fil^{0}_{\Tate} \TP(R)^{\wedge}_{p}$.
We then obtain a diagram
$$ \cdots \rightarrow \Fil^{1}_{\mot} \TC^{-}(R)^{\wedge}_{p} \rightarrow
 \Fil^{0}_{\mot} \TC^{-}(R)^{\wedge}_{p} \rightarrow \Fil^{-1}_{\mot} \TC^{-}(R)^{\wedge}_{p} \rightarrow \cdots,$$
which we will denote by $\Fil^{\bullet}_{\mot} \TC^{-}(R)^{\wedge}_{p}$ and refer to it as the {\it motivic filtration} on the spectrum $\TC^{-}(R)^{\wedge}_{p}$.
Beware that this terminology is slightly misleading: there is a canonical map from the $p$-completed colimit $\varinjlim_{n} \Fil^{-n}_{\mot} \TC^{-}(R)^{\wedge}_{p}$ to $\TC^{-}(R)^{\wedge}_{p}$ which is an isomorphism in the case where $R$ is $p$-quasisyntomic (Proposition~7.13 of \cite{BMS2}), but not in general (Warning \ref{warning:only-filtration-colimit}).
In other words, the motivic filtration of $\TC^{-}(R)^{\wedge}_{p}$ need not be exhaustive.
\end{example}

\begin{example}[The Motivic Filtration on $\TP$]
Let $R$ be an animated commutative ring. For every integer $n$, For every integer $n$, let us write $\Fil^{n}_{\mot} \TP(R)^{\wedge}_{p}$ for the 
colimit $\varinjlim_{m} \Fil^{n}_{\mot} \Fil^{-m}_{\Tate} \TP(R)^{\wedge}_{p}$. We then obtain a diagram
$$ \cdots \rightarrow \Fil^{1}_{\mot} \TP(R)^{\wedge}_{p} \rightarrow
 \Fil^{0}_{\mot} \TP(R)^{\wedge}_{p} \rightarrow \Fil^{-1}_{\mot} \TP(R)^{\wedge}_{p} \rightarrow \cdots,$$
which we will denote by $\Fil^{\bullet}_{\mot} \TP(R)^{\wedge}_{p}$ and refer to it as the {\it motivic filtration} on the spectrum $\TP(R)^{\wedge}_{p}$.
As in Example \ref{example:motivic-TC-minus}, this terminology is slightly misleading in general: there is a canonical map from $p$-completed
colimit $\varinjlim_{n} \Fil^{-n}_{\mot} \TP(R)^{\wedge}_{p}$ to $\TP(R)^{\wedge}_{p}$ which is an isomorphism in the case where $R$ is $p$-quasisyntomic,
but not in general. 
\end{example}

Our proof of Theorem \ref{theorem:BMS2-main} will also show the following:

\begin{proposition}\label{proposition:motivic-completeness}
Let $R$ be an animated commutative ring and let $n$ be an integer. 
Then the filtered spectrum $\Fil^{n}_{\mot} \Fil^{\bullet}_{\Tate} \TP(R)^{\wedge}_{p}$ is
$n$-connective with respect to the Beilinson t-structure. In other words, for every integer $m$,
the homotopy groups of the spectrum $\Fil^{n}_{\mot} \gr^{m}_{\Tate} \TP(R)^{\wedge}_{p}$ are concentrated
in degrees $\geq n-m$. 
\end{proposition}

\begin{corollary}\label{corollary:motivic-completeness}
For every animated commutative ring $R$, the motivic filtration on $\Fil^{\bullet}_{\Tate} \TP(R)^{\wedge}_{p}$
is complete. That is, for every integer $m$, the limit $\varprojlim_{n} \Fil^{n}_{\mot} \Fil^{m}_{\Tate} \TP(R)^{\wedge}_{p}$
vanishes.
\end{corollary}

\begin{proof}
Since each $\Fil^{n}_{\mot} \Fil^{\bullet}_{\Tate} \TP(R)^{\wedge}_{p}$ is complete with respect to the Tate filtration,
it suffices to prove the vanishing of the limit $\varprojlim_{n} \Fil^{n}_{\mot} \gr^{m}_{\Tate} \TP(R)^{\wedge}_{p}$,
which is immediate from the connectivity estimate supplied by Proposition \ref{proposition:motivic-completeness}
\end{proof}

\begin{proof}[Proof of Theorem \ref{theorem:BMS2-main} and Proposition \ref{proposition:motivic-completeness}]
We follow the outline supplied by Remark \ref{remark:motivic-uniqueness}.
Let $\CAlg^{\qrsp}$ denote the category of quasiregular semiperfectoid commutative rings. For $R \in \CAlg^{\qrsp}$, we define
$\Fil^{\bullet}_{\mot} \Fil^{\bullet}_{\Tate} \TP(R)^{\wedge}_{p}$ by the formula
$$\Fil^{n}_{\mot} \Fil^{m}_{\Tate} \TP(R)^{\wedge}_{p} = \tau^{\leq -2n}  \Fil^{m}_{\Tate} \TP(R)^{\wedge}_{p}.$$
This construction tautologically satisfies condition $(3)$ of Theorem \ref{theorem:BMS2-main}, and guarantees the vanishing of the limit
the limit $\varprojlim_{n} \Fil^{n}_{\mot}  \Fil^{\bullet}_{\Tate} \TP(R)^{\wedge}_{p}$ when $R$ is quasiregular semiperectoid.
Moreover, Corollary \ref{corollary:BMS-shifted} supplies a canonical isomorphism of filtered spectra
$$\gr^{n}_{\mot} \Fil^{\bullet}_{\Tate} \TP(R)^{\wedge}_{p} \simeq \Fil^{\bullet+n}_{\Nyg} \widehat{\Prism}_{R}\{n\}[2n].$$
In particular, the functor
$$ \CAlg^{\qrsp} \rightarrow \DFilt(\Sphere) \quad \quad R \mapsto \gr^{n}_{\mot} \Fil^{\bullet}_{\Tate} \TP(R)^{\wedge}_{p}$$
takes values in the subcategory $\DFiltIComp(\Sphere)$ and satisfies descent for the $p$-quasisyntomic topology. It follows by induction that for every integer $k \geq 0$, the
functor $$R \mapsto \Fil^{n}_{\mot} \Fil^{\bullet}_{\Tate} \TP(R)^{\wedge}_{p} / \Fil^{n+k}_{\mot} \Fil^{\bullet}_{\Tate} \TP(R)^{\wedge}_{p}$$
has the same properties. Passing to the limit over $k$, we conclude that the functor
$$ \CAlg^{\qrsp} \rightarrow \DFilt(\Sphere) \quad \quad R \mapsto \Fil^{n}_{\mot} \Fil^{m}_{\Tate} \TP(R)^{\wedge}_{p}$$
takes values in the subcategory $\DFiltIComp(\Sphere)$ and satisfies descent for the $p$-quasisyntomic topology. 

Let $\CAlg^{\QSyn}$ denote the category of $p$-quasisyntomic commutative rings. For each integer $n$, we define
$$ \Fil^{n}_{\mot} \Fil^{\bullet}_{\Tate} \TP(R)^{\wedge}_{p}: \CAlg^{\QSyn} \rightarrow \DFiltIComp(\Sphere)$$
to be the right Kan extension of its restriction to $\CAlg^{\qrsp}$. Since the formation of right Kan extensions commutes with limits,
it follows that the limit $\varprojlim_{n} \Fil^{n}_{\mot} \Fil^{\bullet}_{\Tate} \TP(R)^{\wedge}_{p}$ vanishes for $R \in \CAlg^{\QSyn}$.
For every integer $m$, the proof of Theorem \ref{theorem:motivic-gr} supplies isomorphisms
$$ \gr^{n}_{\mot} \gr^{m}_{\Tate} \TP(R)^{\wedge}_{p} \simeq \gr^{m+n}_{\Nyg} \widehat{\Prism}_{R}\{n\}[2n],$$
so that the homotopy groups of $\gr^{n}_{\mot} \gr^{m}_{\Tate} \TP(R)^{\wedge}_{p}$ are concentrated in degrees $\geq n-m$ (see Proposition \ref{proposition:Nygaard-filtration-Beilinson-connective}). Combining this with the completeness of the motivic filtration, we deduce that the homotopy groups of the spectrum $\Fil^{n}_{\mot} \gr^{m}_{\Tate} \TP(R)^{\wedge}_{p}$ are also concentrated in degrees $\geq n-m$:
that is, Proposition \ref{proposition:motivic-completeness} holds in the case where $R$ is $p$-quasisyntomic.

We now claim that, for each integer $n$, the functor
$$ \CAlg^{\QSyn} \rightarrow \DFiltIComp(\Sphere) \quad \quad R \mapsto  \Fil^{n}_{\mot} \Fil^{\bullet}_{\Tate} \TP(R)^{\wedge}_{p}$$ 
is a left Kan extension of its restriction to the category $\Poly_{\Z}$ of finitely
generated polynomial algebras over $\Z$. To prove this, it will suffice to show that for every sifted diagram $\{ R_{\alpha} \}$ in $\CAlg^{\QSyn}$ whose colimit
$R$ (in the $\infty$-category $\CAlg^{\anim}$) also belongs to $\CAlg^{\QSyn}$, then $\Fil^{n}_{\mot} \Fil^{\bullet}_{\Tate} \TP(R)^{\wedge}_{p}$
is the colimit of the diagram $\{ \Fil^{n}_{\mot} \Fil^{\bullet}_{\Tate} \TP(R_{\alpha})^{\wedge}_{p} \}$ in the $\infty$-category $\DFiltIComp(\Sphere)$.
Equivalently, we must show that for every integer $m$, the map
$$ \gamma_{n}: \varinjlim_{\alpha} \Fil^{n}_{\mot} \gr^{m}_{\Tate} \TP( R_{\alpha} )^{\wedge}_{p} \rightarrow
\Fil^{n}_{\mot} \gr^{m}_{\Tate} \TP(R)^{\wedge}_{p}$$
is an equivalence in the $\infty$-category $\widehat{\calD}(\Sphere)$ of $p$-complete spectra. The preceding argument shows that
the source and target of $\gamma_n$ are $(n-m)$-connective spectra, so the cofiber $\cofib( \gamma_n )$ is $(n-m)$-connective.
Since the functor $R \mapsto \gr^{m+n}_{\Nyg} \widehat{\Prism}_{R}\{n\}$ commutes with sifted colimits (Remark \ref{remark:graded-Nygaard-sifted-colimit}), the spectrum
$\cofib( \gamma_n )$ is independent of $n$. It follows that $\cofib(\gamma_n) \simeq \cofib( \gamma_{n'} )$ is $(n'-m)$-connective for any $n' \geq n$, and therefore vanishes.

For each integer $n$, Proposition \ref{proposition:universal-of-animated} shows that
 $$ \CAlg^{\QSyn} \rightarrow \DFiltIComp(\Sphere) \quad \quad R \mapsto \Fil^{n}_{\mot} \Fil^{\bullet}_{\Tate} \TP(R)^{\wedge}_{p}$$ admits an essentially unique extension to the $\infty$-category $\CAlg^{\anim}$ of animated commutative rings which preserves sifted colimits (and therefore
satisfies condition $(1)$ of Theorem \ref{theorem:BMS2-main}). Note that, since the collection of connective spectra is closed under the
formation of colimits, this extension automatically satisfies the conclusion of Proposition \ref{proposition:motivic-completeness}.
In particular, the motivic filtration is complete (Corollary \ref{corollary:motivic-completeness}).
As in the proof of Theorem \ref{theorem:motivic-gr}, we obtain canonical isomorphisms
$$\gr^{n}_{\mot} \gr^{m}_{\Tate} \TP(R)^{\wedge}_{p} \simeq \gr^{n+m}_{\Nyg} \Prism_{R}\{n\}[2n],$$
so that the functor $R \mapsto \gr^{n}_{\mot} \gr^{m}_{\Tate} \TP(R)^{\wedge}_{p}$ satisfies descent for the $p$-complete faithfully flat
topology (Proposition \ref{proposition:gr-nygaard-flat-descent}). Combining this observation Corollary \ref{corollary:motivic-completeness},
we conclude that the functor $R \mapsto \Fil^{\bullet}_{\mot} \Fil^{\bullet}_{\Tate} \TP(R)^{\wedge}_{p}$ also satisfies $p$-complete fpqc descent.

We now complete the proof of Theorem \ref{theorem:motivic-gr} by showing that, for every animated commutative ring $R$, we
can identify $\Fil^{\bullet}_{\Tate} \TP(R)^{\wedge}_{p}$ with the colimit of the diagram $\{ \Fil^{-n}_{\mot} \Fil^{\bullet}_{\Tate} \TP(R)^{\wedge}_{p} \}$
is an isomorphism in the $\infty$-category $\DFiltIComp(\Sphere)$. Equivlaently, we show that for every integer $m$,
the natural map 
$$ \rho_{R}:  \varinjlim_{n} \Fil^{-n}_{\mot} \gr^{m}_{\Tate} \TP(R)^{\wedge}_{p} \rightarrow \gr^{m}_{\Tate} \TP(R)^{\wedge}_{p} \simeq \THH(R)^{\wedge}_{p}[-2m]$$
is an isomorphism, where the colimit is formed in the $\infty$-category of $p$-complete spectra. Note that the formation $R \mapsto \rho_{R}$ commutes with sifted colimits; we may
therefore assume without loss of generality that $R$ is $p$-quasisyntomic (or even a finitely generated polynomial algebra over $\Z$). Note that
the spectra $\gr^{-n}_{\mot} \gr^{m}_{\Tate} \TP(R)^{\wedge}_{p} \simeq \gr^{m-n}_{\Nyg} \widehat{\Prism}_{R}\{-n\}[-2n]$ vanish for $n > m$.
We can therefore identify the domain of $\rho_{R}$ with the spectrum $\Fil^{-m}_{\mot} \gr^{m}_{\Tate} \TP(R)^{\wedge}_{p}$. 
Using $(i)$ and fpqc descent $\THH$ (\cite[Corollary 3.4]{BMS2}), we see that the construction $R \mapsto \rho_{R}$ satisfies $p$-complete fpqc descent. 
We are therefore reduced showing that $\rho_R$ is an isomorphism in the special case where 
$R$ is a quasiregular semiperfectoid ring. In this case, we have a stronger result: the canonical map
$$ \varinjlim_{n} \Fil^{-n}_{\mot} \Fil^{m}_{\Tate} \TP(R)^{\wedge}_{p} \rightarrow \Fil^{m}_{\Tate} \TP(R)^{\wedge}_{p}$$
is an isomorphism for every integer $m$ (since every spectrum $X$ can be recovered as the colimit 
$\varinjlim_{n} \tau^{\leq n} X$ of its Postnikov approxmations).
\end{proof}

\begin{remark}\label{remark:connective-TC-to-TP}
Let $R$ be an animated commutative ring. For every pair of integers $m$ and $n$, the cofiber of the map
$$ \Fil^{n}_{\mot} \Fil^{m}_{\Tate} \TP(R)^{\wedge}_{p} \rightarrow 
\varinjlim_{m' \leq m} \Fil^{n}_{\mot} \Fil^{m'}_{\Tate} \TP(R)^{\wedge}_{p}$$
is $(n-m)$-connective. This follows from the observation that, for every integer $m' \leq m$, the complex $\Fil^{n}_{\mot} \gr^{m'}_{\Tate} \TP(R)^{\wedge}_{p}$
is $(n-m')$-connective (Proposition \ref{proposition:motivic-completeness}), and is therefore also $(n-m)$-connective. 
\end{remark}

\begin{corollary}\label{corollary:TP-filtration-complete}
Let $R$ be an animated commutative ring. Then the motivic filtrations on $\THH(R)^{\wedge}_{p}$, $\TC^{-}(R)^{\wedge}_{p}$, and $\TP(R)^{\wedge}_{p}$
are complete: that is, the limits
$$ \varprojlim_{n} \Fil^{n}_{\mot} \THH(R)^{\wedge}_{p} \quad \quad \varprojlim_{n} \Fil^{n}_{\mot} \TC^{-}(R)^{\wedge}_{p} \quad \quad
\varprojlim_{n} \Fil^{n}_{\mot} \TP(R)^{\wedge}_{p}$$
vanish (in the $\infty$-category $\calD(\Sphere)$).
\end{corollary}

\begin{proof}
For $\THH(R)^{\wedge}_{p}$ and $\TC^{-}(R)^{\wedge}_{p}$, this is an immediate consequence of Corollary \ref{corollary:motivic-completeness}.
To handle the case of $\TP(R)^{\wedge}_{p}$, it will suffice to show that the limit
$$ \varprojlim_{n} \cofib( \Fil^{n}_{\mot} \TC^{-}(R)^{\wedge}_{p}  \rightarrow \Fil^{n}_{\mot} \TP(R)^{\wedge}_{p} )$$
vanishes, which is immediate from the connectivity estimate of Remark \ref{remark:connective-TC-to-TP}.
\end{proof}

\subsection{The HKR Filtration}\label{subsection:HKR-filtration}

Let $R$ be an animated commutative ring. We write $\HH(R)$ denote the Hochschild complex of $R$, which we view as an animated commutative ring which represents
the functor $A \mapsto \mathscr{L} \Hom_{ \CAlg^{\anim} }( R, A)$; here we write $\mathscr{L}X$ for the free loop space $\Hom( S^1, X)$ of a space $X$.
We will generally abuse notation by identifying $\HH(R)$ its underlying commutative ring spectrum, which carries an action of the circle group $S^1$.
We let $\HP(R)$ denote the Tate construction $\HH(R)^{tS^1}$ and $\Fil^{\bullet}_{\Tate} \HP(R)$ for its Tate filtration (Example \ref{example:Tate-spectrum}), which
we view as an object of the filtered derived $\infty$-category $\DFiltComp(\Z)$. We write $\HC^{-}(R)$ for the negative cyclic complex $\Fil^{0}_{\Tate} \HP(R) = \HH(R)^{hS^1}$.

A classical result Hochschild-Kostant-Rosenberg asserts that, if $R$ is a smooth $\Z$-algebra, then the homology 
of the Hochschild chain complex $\HH(R)$ can be identified with the algebra of differential forms $\Omega^{\ast}_{R} = \Omega^{\ast}_{R/\Z}$
(see Theorem~5.2 of \cite{HKR}). One can use this calculation to endow the periodic cyclic complex $\HP(R)$ of an arbitrary animated commutative ring $R$ with an analogue of the motivic filtration of Theorem \ref{theorem:BMS2-main}. This filtration was described in \cite{BMS2} in the $p$-complete setting. We review the construction in an integral setting, following work of Antieau (see \cite{antieau}). 

\begin{construction}[The Hochschild-Kostant-Rosenberg Filtration]\label{construction:HKR-filtration}
Let $R$ be an animated commutative ring. For every integer $n$, we let $\tau^{\leq -n}_{\mathrm{Dec}} \Fil^{\bullet}_{\Tate} \HP(R)$
denote the $n$-connective cover of $\Fil^{\bullet}_{\Tate} \HP(R)$ with respect to the Beilinson t-structure on the filtered derived
$\infty$-category $\DFilt(\Z)$ (Proposition \ref{proposition:Beilinson-exists}). Allowing $n$ to vary, we obtain a diagram
$$ \cdots \rightarrow \tau^{\leq -1}_{\mathrm{Dec}} \Fil^{\bullet}_{\Tate} \HP(R) \rightarrow
\tau^{\leq 0}_{\mathrm{Dec}}  \Fil^{\bullet}_{\Tate} \HP(R) \rightarrow \tau^{\leq 1}_{\mathrm{Dec}} \Fil^{\bullet}_{\Tate} \HP(R) \rightarrow \cdots,
$$
which we will refer to as the {\it decalage filtration} on $\Fil^{\bullet}_{\Tate} \HP(R)$ (see Remark \ref{remark:decalage}). Note that, since $\Fil^{\bullet}_{\Tate} \HP(R)$ is
filtration-complete, each truncation $\tau^{\leq -n} \Fil^{\bullet}_{\Tate} \HP(R)$ is also filtration-complete (Remark \ref{remark:characterize-Beilinson-connective}).

Fix an integer $n$. By virtue of Proposition \ref{proposition:universal-of-animated}, there is an essentially unique functor
$$ \CAlg^{\anim} \rightarrow \DFiltComp(\Z) \quad \quad R \mapsto \Fil^{n}_{\HKR} \Fil^{\bullet}_{\Tate} \HP(R)$$
which commutes with sifted colimits and coincides with $R \mapsto \tau^{\leq -n}_{\mathrm{Dec}} \Fil^{\bullet}_{\Tate} \HP(R)$ in
the case where $R$ is a finitely generated polynomial algebra over $\Z$. This construction depends functorially on $n$, and
determines a diagram
$$ \cdots \rightarrow \Fil^{1}_{\HKR} \Fil^{\bullet}_{\Tate} \HP(R) \rightarrow
\Fil^{0}_{\HKR} \Fil^{\bullet}_{\Tate} \HP(R) \rightarrow \Fil^{-1}_{\HKR} \Fil^{\bullet}_{\Tate} \HP(R) \rightarrow \cdots$$
in the $\infty$-category $\DFiltComp(\Z)$. We denote this diagram by $\Fil^{\bullet}_{\HKR} \Fil^{\bullet}_{\Tate} \HP(R)$
and refer to it as the {\it Hochschild-Kostant-Rosenberg} filtration on $\Fil^{\bullet}_{\Tate} \HP(R)$.
\end{construction}

\begin{remark}\label{remark:HKR-calculation}
Let $R$ be a commutative ring. Then the quotient 
$$\tau^{\leq 0}_{\mathrm{Dec}} \Fil^{\bullet}_{\Tate} \HP(R) / \tau^{< 0}_{\mathrm{Dec}} \Fil^{\bullet}_{\Tate} \HP(R)$$
is a commutative algebra object of $\DFilt(\Z)$ which belongs to the heart of the Beilinson t-structure, and can therefore be identified with a differential graded ring
$$ \cdots \rightarrow A^{-2} \xrightarrow{\partial} A^{-1} \xrightarrow{\partial} A^0 \xrightarrow{\partial} A^{1} \xrightarrow{\partial} A^{2} \rightarrow \cdots$$
Concretely, each $A^{m}$ can be identified with the Hoschschild homology group 
$$\HH_m(R) \simeq \mathrm{H}^{m}( \HH(R)[-2m] ) = \mathrm{H}^{m}( \gr^{m}_{\Tate} \HP(R) ),$$
and the differential $\partial$ is obtained from the circle action on the Hochschild complex. In particular, the cochain complex
$(A^{\ast}, \partial)$ is concentrated in nonnegative cohomological degrees, and we have a canonical isomorphism
$\alpha: R \simeq \HH_0(R) = A^{0}$. A classical calculation of Hochschild-Kostant-Rosenberg (\cite{HKR}) guarantees that,
if $R$ is a smooth $\Z$-algebra, then $\alpha$ extends uniquely to an isomorphism of differential graded algebras $( \Omega^{\ast}_{R/\Z}, d) \rightarrow (A^{\ast}, \partial)$. More generally, for any integer $n$, we can identify
$$\tau^{\leq -n}_{\mathrm{Dec}} \Fil^{\bullet}_{\Tate} \HP(R) / \tau^{< -n}_{\mathrm{Dec}} \Fil^{\bullet}_{\Tate} \HP(R)$$
with the shifted de Rham complex $( \Omega^{\ast+n}_{R/\Z}, d)$.
\end{remark}

\begin{remark}[Successive Quotients]\label{remark:gr-of-HKR}
Let $R$ be a finitely generated polynomial ring over $\mathbf{Z}$. Then Remark \ref{remark:HKR-calculation} supplies isomorphisms
$$ \gr^{n}_{\HKR} \Fil^{\bullet}_{\Tate} \HP(R) \simeq ( \Omega^{\geq \bullet + n}_{R/\Z}, d)[n].$$
By construction, the functor
$$ \CAlg^{\anim} \rightarrow \DFiltComp(\Z) \quad \quad R \mapsto \gr^{n}_{\HKR} \Fil^{\bullet}_{\Tate} \HP(R)$$
commutes with sifted colimits. It follows that, in general, we have canonical isomorphisms
$$ \gr^{n}_{\HKR} \Fil^{\bullet}_{\Tate} \HP(R) \simeq \Fil^{\bullet+n}_{\Hodge} \dR_{R}^{\hc}[n],$$
where $\dR_{R}^{\hc}$ denotes the {\em Hodge-complete} derived de Rham complex of $R$ (see Construction
\ref{construction:Hodge-complete-dR}).
\end{remark}

\begin{remark}\label{silash}
Let $R$ be an animated commutative ring. For every integer $n$, the filtered complex
$\Fil^{n}_{\HKR} \Fil^{\bullet}_{\Tate} \HP(R)$ is $n$-connective with respect to the Beilinson
t-structure on $\DFiltComp(\Z)$. In other words, for every integer $m$, the complex
$\Fil^{n}_{\HKR} \gr^{m}_{\Tate} \HP(R)$ is $n-m$-connective. This is tautological in the case where $R$ is a finitely generated polynomial ring,
and follows in general since Beilinson-connectivity is preserved by the formation of colimits in the $\infty$-category $\DFiltComp(\Z)$).
\end{remark}

\begin{remark}\label{remark:HKR-complete}
Let $R$ be an animated commutative ring. Then the Hochschild-Kostant-Rosenberg on $\Fil^{\bullet}_{\Tate} \HP(R)$ is complete:
that is, the limit $\varprojlim_{n} \Fil^{n}_{\HKR} \Fil^{\bullet}_{\Tate} \HP(R)$ vanishes in the $\infty$-category $\DFiltComp(\Z)$.
To prove this, it suffices to show that for every integer $m$, the limit $\varprojlim_{n} \Fil^{n}_{\HKR} \gr^{m}_{\Tate} \HP(R)$
vanishes in $\calD(\Z)$, which is immediate from the connectivity estimate of Remark \ref{silash}
\end{remark}

\begin{remark}
Let $R$ be an animated commutative ring. Then the Hochschild-Kostant-Rosenberg filtration is exhaustive in the following
weak sense: the filtered complex $\Fil^{\bullet}_{\Tate} \HP(R)$ can be identified with the colimit of the diagram
$$ \cdots \rightarrow \Fil^{0}_{\HKR} \Fil^{\bullet}_{\Tate} \HP(R) \rightarrow  \Fil^{-1}_{\HKR} \Fil^{\bullet}_{\Tate} \HP(R) \rightarrow
 \Fil^{-2}_{\HKR} \Fil^{\bullet}_{\Tate} \HP(R) \rightarrow \cdots$$
 in the $\infty$-category $\DFiltComp(\Z)$ of filtered complexes which are {\em filtration-complete}. Equivalently, for every integer $m$,
we can identify $\gr^{m}_{\Tate} \HP(R)$ with the colimit of the diagram
$$ \cdots \rightarrow \Fil^{0}_{\HKR} \gr^{m}_{\Tate} \HP(R) \rightarrow  \Fil^{-1}_{\HKR} \gr^{m}_{\Tate} \HP(R) \rightarrow
 \Fil^{-2}_{\HKR} \gr^{m}_{\Tate} \HP(R) \rightarrow \cdots$$
 in $\calD(\Z)$. Beware that the analogous statement for $\Fil^{m}_{\Tate} \HP(R)$ is false in general,
 but holds under some mild additional assumptions on $R$: see Corollary~4.11 of \cite{antieau}, and
 Corollary \ref{corollary:exhaustive-level} for a related statement).
\end{remark}

\begin{remark}\label{remark:gr-of-HKR-HH}
Let $R$ be an animated commutative ring. 
Passing to the associated graded with respect to the Tate filtration,
Construction \ref{construction:HKR-filtration} filtration yields an exhaustive filtration
 $$ \cdots \rightarrow \Fil^{1}_{\HKR} \HH(R) \rightarrow
\Fil^{0}_{\HKR} \HH(R) \rightarrow \Fil^{-1}_{\HKR} \HH(R) \rightarrow \cdots$$
of the Hochschild complex $\HH(R) \simeq \gr^{0}_{\Tate} \HP(R)$.
Note that each $\Fil^{n}_{\HKR} \HH(R)$ is concentrated in cohomological degrees $\leq -n$,
and Remark \ref{remark:gr-of-HKR} supplies isomorphisms
$$ \gr^{n}_{\HKR} \HH(R) \simeq \gr^{n}_{\Hodge} \dR_{R}^{\hc}[n] \simeq L\Omega^{n}_{R}[n].$$
In particular, the complex $\gr^{n}_{\HKR} \HH(R)$ vanishes for $n<0$, so we have
$\Fil^{n}_{\HKR} \HH(R) \simeq \HH(R)$ for $n \leq 0$.
More generally, the natural map $\Fil^{n}_{\HKR} \gr^{m}_{\Tate} \HP(R) \rightarrow \gr^{m}_{\Tate} \HP(R)$
is an isomorphism for $n \leq -m$.
\end{remark}

\begin{proposition}\label{proposition:gr-HKR_HH}
Let $R$ be a commutative ring, and suppose that the absolute cotangent complex $L\Omega^{1}_{R}$
is a flat $R$-module (regarded as a complex concentrated in cohomological degree zero). Then
the Hochschild-Kostant-Rosenberg filtration of $\HH(R)$ coincides with the Postnikov filtration:
that is, for every integer $m$, the map $\Fil^{n}_{\HKR} \HH(R) \rightarrow \HH(R)$ is an isomorphism on cohomology
in degrees $\leq -n$ (and the cohomology groups of $\Fil^{n}_{\HKR} \HH(R)$ vanish in degrees $\leq n$). 
\end{proposition}

\begin{proof}
It follows from Remark \ref{remark:gr-of-HKR-HH} that for every integer $n$, the complex $\gr^{n}{\HKR} \HH(R) \simeq L\Omega^{n}_{R}[n]$
is concentrated in cohomological degree $-n$.
\end{proof}
 
For smooth $\Z$-algebras, the Hochschild-Kostant-Rosenberg filtration coincides with the decalage filtration:

\begin{corollary}\label{corollary:sixpence}
Let $R$ be a commutative ring, and suppose that the absolute cotangent complex $L\Omega^{1}_{R}$ is a flat $R$-module (for instance, $R$ could be smooth over $\Z$). Then, for every integer $n$, the canonical map
$$ \Fil^{n}_{\HKR} \Fil^{\bullet}_{\Tate} \TP(R) \rightarrow \Fil^{\bullet}_{\Tate} \TP(R)$$
exhibits $\Fil^{n}_{\HKR} \Fil^{\bullet}_{\Tate} \TP(R)$ as the $n$-connective cover of
$\Fil^{\bullet}_{\Tate} \TP(R)$ with respect to the Beilinson t-structure on $\DFilt(\Z)$.
\end{corollary}

\begin{proof}
Since $\Fil^{n}_{\HKR} \Fil^{\bullet}_{\Tate} \TP(R)$ and $\Fil^{\bullet}_{\Tate} \TP(R)$
are filtration-complete, this can be checked at the associated graded level:
that is, it suffices to show that each of the maps $$\Fil^{n}_{\HKR} \gr^{m}_{\Tate} \TP(R) \rightarrow
\gr^{m}_{\Tate} \TP(R)$$ exhibits $\Fil^{n}_{\HKR} \gr^{m}_{\Tate} \TP(R)$
as the $(n-m)$-connective cover of $\gr^{m}_{\Tate} \TP(R)$. To prove this, we may assume
without loss of generality that $m=0$, in which case it follows from Proposition \ref{proposition:gr-HKR_HH}.
\end{proof}

We will need a slight variant of Corollary \ref{corollary:sixpence}. Let $R$ be an animated commutative ring. For every pair of integers $m$ and $n$,
we let $\Fil^{n}_{\HKR} \Fil^{m}_{\Tate} \HP(R)^{\wedge}_{p}$ denote the $p$-completion
of the complex $\Fil^{n}_{\HKR} \Fil^{m}_{\Tate} \HP(R)$.

\begin{variant}\label{variant:elpence}
Let $R$ be a commutative ring having bounded $p$-power torsion, and suppose that the absolute
cotangent complex $L\Omega^{1}_{R}$ is $p$-completely flat over $R$. Then,
for every integer $n$, the canonical map
$$ \Fil^{n}_{\HKR} \Fil^{\bullet}_{\Tate} \TP(R)^{\wedge}_{p} \rightarrow \Fil^{\bullet}_{\Tate} \TP(R)^{\wedge}_{p}$$
exhibits $\Fil^{n}_{\HKR} \Fil^{\bullet}_{\Tate} \TP(R)^{\wedge}_p$ as the $n$-connective cover of $\Fil^{\bullet}_{\Tate} \TP(R)^{\wedge}_p$ with respect to the Beilinson t-structure on $\widehat{\DFilt}(\Z_p)$.
\end{variant}

\begin{proof}
Arguing as in the proof of Corollary \ref{corollary:sixpence}, we are reduced to showing that the
Hochschild-Kostant-Rosenberg filtration of $\HH(R)^{\wedge}_{p}$ coincides with its Postnikov filtration:
that is, that the associated graded complexes $\gr^{n}_{\HKR} \HH(R)^{\wedge}_{p} \simeq L \widehat{\Omega}^{n}_{R}[-n]$
are concentrated in cohomological degree $-n$. Our hypothesis that $L\Omega^{1}_{R}$ is $p$-completely
flat over $R$ guarantees that $L\widehat{\Omega}^{n}_{R}$ is also $p$-completely flat over $R$, and is
therefore concentrated in cohomological degree zero since $R$ has bounded $p$-power torsion.
\end{proof}

\subsection{The Global Motivic Filtration}\label{subsection:non-complete-motivic}

Throughout this section, we temporarily suspend our convention that the prime number $p$ is fixed. Our goal is to construct a decompleted
analogue of the motivic filtration of Theorem \ref{theorem:BMS2-main}. 

Let $R$ be an animated commutative ring. There is a comparison map $\THH(R) \rightarrow \HH(R)$,
given by extension of scalars along the map of commutative ring spectra $\THH(\Z) \rightarrow \Z$.
This map is $S^1$-equivariant, and therefore induces a map of Tate filtrations
$$ \Fil^{\bullet}_{\Tate} \TP(R) \rightarrow \Fil^{\bullet}_{\Tate} \HH(R).$$
Under this comparison, there is a close relationship between the motivic filtration
of Theorem \ref{theorem:BMS2-main} and the Hochschild-Kostant-Rosenberg filtration of Construction \ref{construction:HKR-filtration}.

\begin{proposition}\label{proposition:compare-filtrations}
Let $p$ be a prime number, let $n$ be an integer, and let $R$ be an animated commutative ring.
Then the map of filtered spectra
$$ \Fil^{n}_{\mot} \Fil^{\bullet}_{\Tate} \TP(R)^{\wedge}_{p} \rightarrow 
\Fil^{\bullet}_{\Tate} \TP(R)^{\wedge}_{p}
\rightarrow \Fil^{\bullet}_{\Tate} \HP(R)^{\wedge}_{p}$$ 
factors canonically through the filtered spectrum $\Fil^{n}_{\HKR} \Fil^{\bullet}_{\Tate} \HP(R)^{\wedge}_{p}$.
This factorization is determined (up to homotopy) by the requirement that it depends functorially on $R$.
\end{proposition}

\begin{proof}
By construction, the functor
$$ \CAlg^{\anim} \rightarrow \DFiltIComp(\Sphere) \quad \quad R \mapsto  \Fil^{n}_{\mot} \Fil^{\bullet}_{\Tate} \TP(R)^{\wedge}_{p}$$
commutes with sifted colimits and is therefore a left Kan extension of its restriction to the category of finitely generated
polynomial algebras over $\Z$ (see Theorem \ref{theorem:BMS2-main}). It will therefore suffice to construct
the desired factorization in the case where $R$ is a finitely generated polynomial algebra over $\Z$.
In this case, Variant \ref{variant:elpence} identifies $\Fil^{n}_{\HKR} \Fil^{\bullet}_{\Tate} \HP(R)^{\wedge}_{p}$
with the $n$-connective cover of $\Fil^{\bullet}_{\Tate} \HP(R)^{\wedge}_{p}$ with respect to the Beilinson
t-structure on $\DFiltIComp(\Z)$. The existence (and uniqueness up to contractible ambiguity) of the desired factorization
now follow from the observation that, for every animated commutative ring $R'$, the filtered spectrum $\Fil^{n}_{\mot} \Fil^{\bullet}_{\Tate} \TP(R')^{\wedge}_{p}$ is $n$-connective with respect to the Beilinson t-structure (Proposition \ref{proposition:motivic-completeness}).
\end{proof}

We now use Proposition \ref{proposition:compare-filtrations} to define an integral counterpart of the motivic filtration.

\begin{construction}[The Motivic Filtration]\label{construction:integral-motivic-filtration}
Let $R$ be an animated commutative ring and let $n$ be an integer. We let
$\Fil^{n}_{\mot} \Fil^{\bullet}_{\Tate} \TP(R)$ denote the pullback of the diagram of filtered spectra
$$ \Fil^{n}_{\HKR} \Fil^{\bullet}_{\Tate} \HP(R) \rightarrow
\prod_{p} \Fil^{n}_{\HKR} \Fil^{\bullet}_{\Tate} \HP(R)^{\wedge}_{p} \leftarrow
\prod_{p} \Fil^{n}_{\mot} \Fil^{\bullet}_{\Tate} \TP(R)^{\wedge}_{p}.$$
Here the product is taken over all prime numbers $p$, and the map on the right is given by
Proposition \ref{proposition:compare-filtrations}. This construction depends functorially on $n$ (and on $R$),
and therefore determines a diagram of filtered spectra
$$ \cdots \rightarrow \Fil^{1}_{\mot} \Fil^{\bullet}_{\Tate} \TP(R) \rightarrow
\Fil^{0}_{\mot} \Fil^{\bullet}_{\Tate} \TP(R) \rightarrow \Fil^{-1}_{\mot} \Fil^{\bullet}_{\Tate} \TP(R) \rightarrow \cdots$$
We will denote this diagram by $\Fil^{\bullet}_{\mot} \Fil^{\bullet}_{\Tate} \TP(R)$ and refer to it as
the {\it motivic filtration} on $\Fil^{\bullet}_{\Tate} \TP(R)$ (see Proposition \ref{proposition:motivic-filtration-exhaustive} below).
\end{construction}

\begin{remark}
Let $R$ be an animated commutative ring. For every prime number $p$,
the bifiltered spectrum $\Fil^{\bullet}_{\mot} \Fil^{\bullet}_{\Tate} \TP(R)^{\wedge}_{p}$ of
Theorem \ref{theorem:BMS2-main} can be recovered from the bifiltered spectrum
$\Fil^{\bullet}_{\mot} \Fil^{\bullet}_{\Tate} \TP(R)$ of Construction \ref{construction:integral-motivic-filtration} by termwise $p$-completion.
\end{remark}

\begin{remark}\label{remark:motivic-completeness-again}
Let $R$ be an animated commutative ring. Then the motivic filtration of Construction \ref{construction:integral-motivic-filtration}
is complete: that is, the filtered spectrum $\varprojlim_{n}  \Fil^{n}_{\mot} \Fil^{\bullet}_{\Tate} \TP(R)$ vanishes.
This follows by combining Remark \ref{remark:HKR-complete} with Corollary \ref{corollary:motivic-completeness}.
\end{remark}

\begin{proposition}\label{proposition:motivic-filtration-exhaustive}
Let $R$ be an animated commutative ring. Then $\Fil^{\bullet}_{\Tate} \TP(R)$ can be identified with the colimit of the diagram
$$ \cdots \rightarrow \Fil^{0}_{\mot} \Fil^{\bullet}_{\Tate} \TP(R) \rightarrow
\Fil^{-1}_{\mot} \Fil^{\bullet}_{\Tate} \TP(R) \rightarrow \Fil^{-2}_{\mot} \Fil^{\bullet}_{\Tate} \TP(R) \rightarrow \cdots$$
in the $\infty$-category $\DFiltComp(\Sphere)$.
\end{proposition}

\begin{proof}
The comparison map $\THH(R) \rightarrow \HH(R)$ is a rational isomorphism, and therefore determines a pullback diagram of
$S^1$-spectra
$$ \xymatrix@R=50pt@C=50pt{ \THH(R) \ar[r] \ar[d] & \prod_{p} \THH(R)^{\wedge}_{p} \ar[d] \\
\HH(R) \ar[r] & \prod_{p} \HH(R)^{\wedge}_{p}. }$$
Applying the Tate construction, we obtain a pullback diagram of filtered spectra
$$ \xymatrix@R=50pt@C=50pt{ \Fil^{\bullet}_{\Tate} \TP(R) \ar[r] \ar[d] & \prod_{p}  \Fil^{\bullet}_{\Tate} \TP(R)^{\wedge}_{p} \ar[d] \\
 \Fil^{\bullet}_{\Tate} \HP(R) \ar[r] & \prod_{p} \Fil^{\bullet}_{\Tate} \HP(R)^{\wedge}_{p}. }$$
It will therefore suffice to show that the vertical maps in the diagram
$$ \xymatrix@R=50pt@C=35pt{ \varinjlim_{n} \Fil^{-n}_{\HKR} \Fil^{\bullet}_{\Tate} \HP(R) \ar[r] \ar[d] &
\varinjlim_{n} \prod_{p} \Fil^{-n}_{\HKR} \Fil^{\bullet}_{\Tate} \HP(R)^{\wedge}_{p} \ar[d] & \varinjlim_{n} \prod_{p} \Fil^{-n}_{\mot} \Fil^{\bullet}_{\Tate} \TP(R)^{\wedge}_{p} \ar[d] \ar[l] \\
\Fil^{\bullet}_{\Tate} \HP(R) \ar[r] & \prod_{p} \Fil^{\bullet}_{\Tate} \HP(R)^{\wedge}_{p} & \prod_{p} \Fil^{\bullet}_{\Tate} \TP(R)^{\wedge}_{p} \ar[l] }$$
are isomorphisms, where the colimits are computed in the $\infty$-category $\DFiltComp(\Sphere)$. This is equivalent
to the assertion that for every integer $m$, the vertical maps in the diagram
$$ \xymatrix@R=50pt@C=35pt{ \varinjlim_{n} \Fil^{-n}_{\HKR} \gr^{m}_{\Tate} \HP(R) \ar[r] \ar[d] &
\varinjlim_{n} \prod_{p} \Fil^{-n}_{\HKR} \gr^{m}_{\Tate} \HP(R)^{\wedge}_{p} \ar[d] & \varinjlim_{n} \prod_{p} \Fil^{-n}_{\mot} \gr^{m}_{\Tate} \TP(R)^{\wedge}_{p} \ar[d] \ar[l] \\
\gr^{m}_{\Tate} \HP(R) \ar[r] & \prod_{p} \gr^{m}_{\Tate} \HP(R)^{\wedge}_{p} & \prod_{p} \gr^{m}_{\Tate} \TP(R)^{\wedge}_{p} \ar[l]. }$$
This is immediate from Remark \ref{remark:gr-of-HKR-HH} and Example \ref{example:THH-filtration} (each of the colimits stabilizes
at $n = -m$).
\end{proof}

\begin{construction}[Global Prismatic Complexes]
Let $R$ be an animated commutative ring. For every integer $n$, we let
$\widehat{\Prism}_{R}^{\glo}\{n\}$ denote the spectrum $\gr^{n}_{\mot} \Fil^{-n}_{\Tate} \TP(R)[-2n]$.
We will refer to $\widehat{\Prism}_{R}^{\glo}\{n\}$ as the {\it global prismatic complex of $R$} (with Breuil-Kisin twist $n$).
More generally, for every pair of integers $m$ and $n$, we let $\Fil^{m}_{\Nyg} \widehat{\Prism}_{R}^{\glo}\{n\}$ denote the spectrum
$\gr^{n}_{\mot} \Fil^{m-n}_{\Tate} \TP(R)[-2n]$. This construction determines a diagram
$$ \cdots \rightarrow \Fil^{2}_{\Nyg} \widehat{\Prism}_{R}^{\glo}\{n\} \rightarrow
\Fil^{1}_{\Nyg} \widehat{\Prism}_{R}^{\glo}\{n\} \rightarrow \Fil^{0}_{\Nyg} \widehat{\Prism}_{R}^{\glo}\{n\} = \widehat{\Prism}_{R}^{\glo}\{n\},$$
which we will denote by $\Fil^{\bullet}_{\Nyg} \widehat{\Prism}_{R}^{\glo}\{n\}$ and refer to as the {\it Nygaard filtration} on
$\widehat{\Prism}_{R}^{\glo}\{n\}$.
\end{construction}

\begin{remark}\label{smasy}
Let $R$ be an animated commutative ring. For every prime number $p$, let us denote the Nygaard-completed prismatic complex $\widehat{\Prism}_{R}\{n\}$ by $\widehat{\Prism}_{R,p}\{n\}$ to emphasize its dependence on the prime number $p$. Combining
Theorem \ref{theorem:motivic-gr} with Remark \ref{remark:gr-of-HKR}, we obtain a pullback diagram of filtered spectra
$$ \xymatrix@R=50pt@C=50pt{ \Fil^{\bullet}_{\Nyg} \widehat{\Prism}_{R}^{\glo}\{n\} \ar[r] \ar[d] & \prod_{p} \Fil^{\bullet}_{\Nyg} \widehat{\Prism}_{R,p}\{n\} \ar[d]^{\prod_{p} \Fil^{\bullet} \widehat{\gamma}_{\Prism}^{\dR}\{n\} }\\
\Fil^{\bullet}_{\Hodge} \dR_{R}^{\hc} \ar[r] & \prod_{p} \Fil^{\bullet}_{\Hodge} ( \dR_R^{\hc} )^{\wedge}_{p}; }$$
here the right vertical map is obtained concretely as the product over all prime numbers $p$ of the de Rham comparison maps appearing in Construction~\ref{construction:absolute-Nygaard-untwisted}. From this description, we see that $\Fil^{\bullet}_{\Nyg} \widehat{\Prism}_{R}^{\glo}\{n\}$ admits a canonical $\Z$-module structure: that is, it can be promoted to an object of the filtered derived $\infty$-category $\DFilt(\Z)$.
\end{remark}

\begin{proposition}\label{proposition:compute-gr-gr}
Let $m$ and $n$ be integers. For every animated commutative ring $R$, there is a canonical fiber sequence
$$ \gr^{m}_{\Nyg} \widehat{\Prism}^{\glo}_{R}\{n\} \rightarrow \Fil_m^{\conj} \Omega^{\DHod}_{R} \xrightarrow{\Theta+m} \Fil_{m-1}^{\conj} \Omega^{\DHod}_{R},$$
where $\Omega^{\DHod}_{R}$ denotes the diffracted Hodge complex of $R$ (Definition \ref{construction:diffracted-Hodge-integral}) and
$\Theta+m$ is the morphism described in Remark \ref{remark:factor-Theta-plus}.
\end{proposition}

\begin{proof}
Let $M$ denote the fiber of the map $\Theta+m: \Fil_m^{\conj} \Omega^{\DHod}_{R} \rightarrow \Fil_{m-1}^{\conj} \Omega^{\DHod}_{R}$. 
Writing $\Omega^{\DHod}_{R,p}$ for the $p$-complete diffracted Hodge complex, we see that the $p$-completion of $M$
can be identified with the fiber
$$ \fib( \Theta+m:  \Fil_m^{\conj} \Omega^{\DHod}_{R,p} \rightarrow \Fil_{m-1}^{\conj} \Omega^{\DHod}_{R,p} ) \simeq \gr^{m}_{\Nyg} \widehat{\Prism}_{R,p}\{n\}$$
(see Remark \ref{remark:Nygaard-associated-graded}). We therefore obtain a commutative diagram of fiber sequences
$$ \xymatrix@R=50pt@C=50pt{ ( \Fil_{m-1}^{\conj} \Omega^{\DHod}_{R} )^{\Theta=-m} \ar[r] \ar[d] & \prod_{p} (\Fil_{m-1}^{\conj} \Omega^{\DHod}_{R,p} )^{\Theta=-m} \ar[d] \\
M \ar[r] \ar[d] & \prod_{p}  \gr^{m}_{\Nyg} \widehat{\Prism}_{R,p}\{n\} \ar[d]^{\alpha} \\
L\Omega^{m}_{R}[-m] \ar[r] & \prod_{p} L\Omega^{m}_{R}[-m]^{\wedge}_{p}. }$$
Since $\Fil_{m-1}^{\conj} \Omega^{\DHod}_{R}$ admits a finite filtration by complexes on which the Sen operator acts by integers different from $-m$,
the upper horizontal map in this diagram is an isomorphism. Unwinding the definitions, we see that $\alpha$ can be identified with the product over all primes $p$
of the map
$$ \gr^{m}_{\Nyg} \widehat{\Prism}_{R,p}\{n\} \rightarrow \gr^{m}_{\Hodge} ( \dR_{R} )^{\wedge}_{p}$$
given by the associated graded of the de Rham comparison. It follows that we can identify $M$
with the pullback of the diagram
$$ \gr^{m}_{\Hodge} \dR_{R} \rightarrow \prod_{p} (\gr^{m}_{\Hodge} \dR_{R})^{\wedge}_{p} \leftarrow \prod_{p} \gr^{m}_{\Nyg} \widehat{\Prism}_{R,p}\{n\},$$
which agrees with the complex $\gr^{m}_{\Nyg} \widehat{\Prism}^{\glo}_{R}\{n\}$ by virtue of Remark \ref{smasy}.
\end{proof}

\begin{example}
For every animated commutative ring $R$, there is a canonical isomorphism $$\gr^{0}_{\Nyg} \widehat{\Prism}_R^{\glo} \simeq R.$$
\end{example}

\begin{corollary}\label{corollary:global-good-behavior}
For every pair of integers $m$ and $n$, the functor $R \mapsto \gr^{m}_{\Nyg} \widehat{\Prism}_{R}^{\glo}\{n\}$
commutes with sifted colimits and satisfies descent for the fpqc topology.
\end{corollary}

\begin{proof}
Combine Proposition \ref{proposition:compute-gr-gr} with Proposition \ref{proposition:sifted-colimit-diffracted-Hodge}.
\end{proof}

\begin{corollary}\label{corollary:global-fpqc}
For every pair of integers $m$ and $n$, the functor $\Fil^{m}_{\Nyg} \widehat{\Prism}_{R}^{\glo}\{n\}$
satisfies descent for the fpqc topology.
\end{corollary}

\begin{proof}
By construction, the global prismatic complex $\widehat{\Prism}_{R}^{\glo}\{n\}$ is complete with respect to its Nygaard filtration,
so the desired result follows from Corollary \ref{corollary:global-good-behavior}.
\end{proof}

\begin{corollary}
For every pair of integers $m$ and $n$, the functor
$$ R \mapsto \Fil^{n}_{\mot} \Fil^{m}_{\Tate} \TP(R)$$
satisfies descent with respect to the fpqc topology.
\end{corollary}

\begin{proof}
By virtue of Remark \ref{remark:motivic-completeness-again}, it suffices to prove the analogous result for
the functor $R \mapsto \gr^{n}_{\mot} \Fil^{m}_{\Tate} \TP(R)$, which is a restatement of Corollary \ref{corollary:global-fpqc}.
\end{proof}

\begin{corollary}
Let $R$ be an animated commutative ring and let $n$ be an integer. Then the complexes $\gr^{m}_{\Nyg} \widehat{\Prism}_{R}^{\glo}\{n\}$
vanish for $m < 0$. Consequently, the natural map $$\widehat{\Prism}_{R}^{\glo}\{n\} = \Fil^{0}_{\Nyg} \widehat{\Prism}_{R}^{\glo}\{n\}
\rightarrow  \Fil^{m}_{\Nyg} \widehat{\Prism}_{R}^{\glo}\{n\}$$ is an isomorphism for $m \leq 0$.
\end{corollary}

\begin{corollary}\label{corollary:coconnectivity-gr-global}
Let $R$ be a commutative ring for which the cotangent complex $L\Omega^{1}_{R}$ has $\Tor$-amplitude contained in $[0,1]$.
Let $m$ and $n$ be integers. Then:
\begin{itemize}
\item[$(a)$] The cohomology groups of the complex $\gr^{m}_{\Nyg} \widehat{\Prism}^{\glo}_{R}\{n\}$ are concentrated in degrees $\geq 0$.
\item[$(b)$] The cohomology groups of the complex $\Fil^{m}_{\Nyg} \widehat{\Prism}^{\glo}_{R}\{n\}$ are concentrated in degrees $\geq 0$.
\end{itemize}
\end{corollary}

\begin{proof}
Assertion $(a)$ follows by combining Proposition \ref{proposition:compute-gr-gr} with Variant \ref{variant:coconnectivity-diffracted-hodge}.
Assertion $(b)$ follows from $(a)$, using the completeness of the Nygaard filtration on $\widehat{\Prism}^{\glo}_{R}\{n\}$.
\end{proof}

\begin{corollary}\label{corollary:sixwell}
Let $R$ be a commutative ring for which the cotangent complex $L\Omega^{1}_{R}$ has $\Tor$-amplitude contained in $[0,1]$.
Let $m$ and $n$ be integers. Then:
\begin{itemize}
\item[$(a)$] The homotopy groups of the spectrum $\gr^{n}_{\mot} \gr^{m}_{\Tate} \TP(R)$ are concentrated in
degrees $\leq 2n$.
\item[$(b)$] The homotopy groups of the spectrum $\gr^m_{\Tate} \TP(R) / \Fil^{n+1}_{\mot} \gr^{m}_{\Tate} \TP(R)$
are concentrated in degrees $\leq 2n$.
\item[$(c)$] The homotopy groups of the spectrum $\Fil^{m}_{\Tate} \TP(R) / \Fil^{n+1}_{\mot} \Fil^{m}_{\Tate} \TP(R)$
are concentrated in degrees $\leq 2n$.
\end{itemize}
\end{corollary}

\begin{proof}
Assertion $(a)$ is a restatement of the corresponding assertion of Corollary \ref{corollary:coconnectivity-gr-global}.
It then follows by induction that, for every integer $k \geq 0$, the cofiber
$$\Fil^{n-k} \gr^m_{\Tate} \TP(R) / \Fil^{n+1}_{\mot} \gr^{m}_{\Tate} \TP(R)$$ has homotopy groups concentrated in degrees $\leq 2n$.
Assertion $(b)$ follows by passing to the colimit over $k$ (and using Proposition \ref{proposition:motivic-filtration-exhaustive}).
Assertion $(c)$ follows from $(b)$ using a similar argument, since both
$\Fil^{\bullet}_{\Tate} \TP(R)$ and $\Fil^{n+1}_{\mot} \Fil^{\bullet}_{\Tate} \TP(R)$ are filtration-complete.
\end{proof}

\begin{corollary}[Exhaustivity]\label{corollary:exhaustive-level}
Let $R$ be a commutative ring for which the cotangent complex $L\Omega^{1}_{R}$ has $\Tor$-amplitude contained in $[0,1]$.
Then, for every integer $m$, the spectrum $\Fil^{m}_{\Tate} \TP(R)$ is the colimit of the diagram
$$ \cdots \rightarrow \Fil^{0}_{\mot} \Fil^{\bullet}_{\Tate} \TP(R) \rightarrow
\Fil^{-1}_{\mot} \Fil^{\bullet}_{\Tate} \TP(R) \rightarrow \Fil^{-2}_{\mot} \Fil^{\bullet}_{\Tate} \TP(R) \rightarrow \cdots$$
in the $\infty$-category $\DFilt(\Sphere)$.
\end{corollary}

\begin{proof}
It follows from Corollary \ref{corollary:sixwell} that, for every integer $n$, the map
$$ \Fil^{-n}_{\mot} \Fil^{m}_{\Tate} \TP(R) \rightarrow \Fil^{m}_{\Tate} \TP(R)$$
induces an isomorphism on homotopy groups in degrees $\geq -2n-1$.
\end{proof}

\begin{example}[The Motivic Filtration on $\THH$]\label{example:THH-filtration-integral}
Let $R$ be an animated commutative ring. For every integer $n$, let us write $\Fil^{n}_{\mot} \THH(R)$ for the spectrum $\Fil^{n}_{\mot} \gr^{0}_{\Tate} \TP(R)$.
We then obtain a diagram
$$ \cdots \rightarrow \rightarrow \Fil^{1}_{\mot} \THH(R) \rightarrow
 \Fil^{0}_{\mot} \THH(R) \rightarrow \Fil^{-1}_{\mot} \THH(R) \rightarrow \cdots,$$
which we denote by $\Fil^{\bullet}_{\mot} \THH(R)$ and refer to it as the {\it motivic filtration} on the spectrum $\THH(R)$.
The associated graded of this filtration is given by
$$ \gr^{n}_{\mot} \THH(R) \simeq \gr^{n}_{\Nyg} \widehat{\Prism}^{\glo}_{R}\{n\}[2n]
\simeq \fib( \Theta+n: \Fil_{n}^{\conj} \Omega^{\DHod}_{R} \rightarrow
\Fil_{n-1}^{\conj} \Omega^{\DHod}_{R} )[2n].$$
\end{example}

\begin{remark}
One can use Example~\ref{example:THH-filtration-integral} in conjunction with the tools developed in this paper to calculate $\pi_* \mathrm{THH}(R)$ in some cases. For instance, if $R=\mathbf{Z}$, then the unit map
\[ \mathbf{Z} = \Fil^0_{\conj} \Omega_{\Z}^{\DHod} \xrightarrow{\sim} \Omega_{\mathbf{Z}}^{\DHod}\]
is an isomorphism, so the Sen operator on $\Omega_{\Z}^{\DHod}$ vanishes (see Remark~\ref{remark:conjugate-equals-Postnikov1}). For $n > 0$, the fiber sequence of Proposition \ref{proposition:compute-gr-gr} supplies isomorphisms
\[ \gr^{n}_{\Nyg} \Prism_{\mathbf{Z}}^{\mathrm{glo}}\{n\} \simeq \mathrm{fib}(\mathbf{Z} \xrightarrow{x \mapsto nx} \mathbf{Z}) \simeq (\mathbf{Z}/n \Z)[-1].\]
Combining this with the isomorphism $\gr^{0}_{\Nyg} \Prism_{\Z} \simeq \Z$, the 
filtration of Example~\ref{example:THH-filtration-integral} recovers B\"{o}kstedt's calculation \cite{BokstedtTHHZ} 
$$ \pi_{m} \mathrm{THH}(R) = \begin{cases} \Z_p & \text{ if } m=0 \\
\Z/ n \Z & \text{ if $m=2n-1 > 0$ } \\
0 & \text{ otherwise.} \end{cases}$$
\end{remark}

\begin{example}[The Motivic Filtration on $\TC^{-}$]\label{example:motivic-TC-minus-integral}
Let $R$ be an animated commutative ring. For every integer $n$, let us write $\Fil^{n}_{\mot} \TC^{-}(R)$ for the spectrum $\Fil^{n}_{\mot} \Fil^{0}_{\Tate} \TP(R)$.
We then obtain a diagram
$$ \cdots \rightarrow \rightarrow \Fil^{1}_{\mot} \TC^{-}(R) \rightarrow
 \Fil^{0}_{\mot} \TC^{-}(R) \rightarrow \Fil^{-1}_{\mot} \TC^{-}(R) \rightarrow \cdots,$$
which we will denote by $\Fil^{\bullet}_{\mot} \TC^{-}(R)$ and refer to it as the {\it motivic filtration} on the spectrum $\TC^{-}(R)$.
The associated graded of this filtration is given by
$$ \gr^{n}_{\mot} \TC^{-}(R) \simeq \Fil^{n}_{\Nyg} \widehat{\Prism}_{R}^{\glo}\{n\}[2n].$$
If $R$ is a commutative ring for which the cotangent complex $L\Omega^{1}_{R}$ has $\Tor$-amplitude concentrated in degrees $[0,1]$,
then Corollary \ref{corollary:exhaustive-level} supplies an isomorphism $\varinjlim_{n} \Fil^{-n}_{\mot} \TC^{-1}(R) \rightarrow \TC^{-1}(R)$.
\end{example}

\begin{example}[The Motivic Filtration on $\TP$]\label{example:motivic-TP-integral}
Let $R$ be an animated commutative ring. For every integer $n$, let us write $\Fil^{n}_{\mot} \TP(R)$ for the colimit
$\varinjlim_{m} \Fil^{n}_{\mot} \Fil^{-m}_{\Tate} \TP(R)$. We then obtain a diagram
$$ \cdots \rightarrow \rightarrow \Fil^{1}_{\mot} \TP(R) \rightarrow
 \Fil^{0}_{\mot} \TP(R) \rightarrow \Fil^{-1}_{\mot} \TP(R) \rightarrow \cdots,$$
which we will denote by $\Fil^{\bullet}_{\mot} \TP(R)$ and refer to it as the {\it motivic filtration} on the spectrum $\TP(R)$.
The associated graded of this filtration is given by
$$ \gr^{n}_{\mot} \TP(R) \simeq \widehat{\Prism}_{R}^{\glo}\{n\}[2n].$$
If $R$ is a commutative ring for which the cotangent complex $L\Omega^{1}_{R}$ has $\Tor$-amplitude concentrated in degrees $[0,1]$,
then Corollary \ref{corollary:exhaustive-level} supplies an isomorphism $\varinjlim_{n} \Fil^{-n}_{\mot} \TP(R) \rightarrow \TP(R)$.
\end{example}

\newpage \section{The First Chern Class}\label{section:first-chern}

Let $X$ be an $\F_p$-scheme. To every line bundle $\mathscr{L}$ on $X$, one can associate a class
$$ c_{1}^{\crys}( \mathscr{L} ) \in \mathrm{H}^{2}_{\crys}( X / \Z_p )$$
in the crystalline cohomology of $X$, which we will refer to as the {\it crystalline first Chern class} of $\mathcal{L}$.
This construction has the following features:
\begin{itemize}
\item The crystalline first Chern class $c_{1}^{\crys}( \mathscr{L} )$ is annihilated by the crystalline augmentation map
$$ \epsilon_{\crys}: \mathrm{H}^{2}_{\crys}( X / \Z_p ) \rightarrow \mathrm{H}^{2}( X, \calO_X ).$$

\item Writing $\varphi$ for the endomorphism of $\mathrm{H}^{2}_{\crys}(X/\Z_p)$ given by pullback along the absolute Frobenius morphism $\varphi_{X}: X \rightarrow X$, we have
$$\varphi( c_{1}^{\crys}(\mathscr{L} ) ) = c_{1}^{\crys}( \varphi_{X}^{\ast} \mathscr{L} ) =
c_{1}^{\crys}( \mathscr{L}^{\otimes p} ) = p \cdot c_{1}^{\crys}( \mathscr{L} ).$$
\end{itemize}

One can combine these observation to produce a refinement of the crystalline first Chern class.
Let $\Fil^{1}_{\Nyg} \RGamma_{\crys}(X/\Z_p)$ denote the fiber of the crystalline
augmentation map $\epsilon_{\crys}: \RGamma_{\crys}(X/ \Z_p ) \rightarrow \RGamma(X, \calO_X)$ of
Notation \ref{notation:crystalline-augmentation}. Assume that the $\F_p$-scheme $X$ is $p$-quasisyntomic, and let
$\varphi/p: \Fil^{1}_{\Nyg} \RGamma_{\crys}(X/\Z_p) \rightarrow \RGamma_{\crys}(X/\Z_p)$ denote the
divided Frobenius morphism of Remark \ref{remark:phi-over-p}. We let $\RGamma_{\Syn}( X, \Z_p(1) )$
denote the fiber of the map 
$$ \frac{ \varphi}{p} - 1: \Fil^{1}_{\Nyg} \RGamma_{\crys}(X/\Z_p) \rightarrow \RGamma_{\crys}(X/\Z_p).$$
We denote the cohomology groups of $\RGamma_{\Syn}(X, \Z_p(1) )$ by $\mathrm{H}^{\ast}_{\Syn}( X, \Z_p(1) )$
and refer to them as {\it syntomic cohomology groups of $X$}; this definition essentially goes back to
Fontaine and Messing (see \cite{syntomic-fm} and \cite{Kato}). The construction $\mathscr{L} \mapsto c_{1}^{\crys}( \mathscr{L} )$ then factors as a composition
$$ \Pic(X) \xrightarrow{ c_{1}^{\Syn} } \mathrm{H}^{2}_{\Syn}( X, \Z_p(1) ) \rightarrow \mathrm{H}^{2}_{\crys}( X, \Z_p(1)),$$
where $c_{1}^{\Syn}$ is a homomorphism of abelian groups which we refer to as the {\it syntomic first Chern class}.

In \cite{BMS2}, the first author, Morrow, and Scholze studied a mixed-characteristic generalization of the preceding construction.
Let $\mathfrak{X}$ be a $p$-quasisyntomic $p$-adic formal scheme. To every integer $n$, \cite{BMS2} associates a chain complex $\RGamma_{\Syn}( \mathfrak{X}, \Z_p(n) )$, which arises naturally as an associated graded piece
for the motivic filtration on the topological cyclic homology spectrum $\TC( \mathfrak{X} )$. Moreover, using the
cyclotomic trace, \cite{BMS2} constructs a map of chain complexes
$$ c_{1}^{\Syn}: \RGamma_{ \textnormal{\'{e}t} }( \mathfrak{X}, \mathbf{G}_{m} )[-1]  \rightarrow 
\RGamma_{\Syn}( \mathfrak{X}, \Z_p(1) ) $$
and shows that it becomes an isomorphism after $p$-completion (Proposition~7.17); passing to
cohomology in degree $2$, we obtain a group homomorphism $\Pic( \mathfrak{X} ) \rightarrow \mathrm{H}^{2}_{\Syn}( \mathfrak{X}, \Z_p(1) )$.

Our goal in this section is to give a purely algebraic exposition of the first Chern class in syntomic cohomology. 
Let $\mathfrak{X}$ be a bounded $p$-adic formal scheme. For every integer $n \geq 0$, let 
$\RGamma_{\Prism}( \mathfrak{X} )\{n\}$ denote the absolute prismatic complex of $\mathfrak{X}$
(Construction \ref{construction:absolute-prismatic-complex-of-scheme}) and let
$\Fil^{n}_{\Nyg} \RGamma_{\Prism}( \mathfrak{X} )\{n\}$ denote the $n$th stage of its Nygaard filtration
(Notation \ref{notation:absolute-Nygaard-filtration-of-scheme}). Writing $\iota$ for the tautological map
$$\Fil^{n}_{\Nyg} \RGamma_{\Prism}( \mathfrak{X})\{n\} \rightarrow \Fil^{0}_{\Nyg} \RGamma_{\Prism}( \mathfrak{X})\{n\} = \Prism_{R}\{n\},$$ we define
$\RGamma_{\Syn}( \mathfrak{X}, \Z_p(n) )$ to be the fiber of the map
$$ ( \varphi\{n\} - \iota): \Fil^{n}_{\Nyg} \RGamma_{\Prism}( \mathfrak{X})\{n\} \rightarrow  \RGamma_{\Prism}( \mathfrak{X})\{n\}$$
where $\varphi\{n\}$ denotes the Frobenius morphism studied in \S \ref{subsection:Frobenius}.
In \S\ref{subsection:syntomic-complex}, we study the properties of the complexes $\{ \RGamma_{\Syn}( \mathfrak{X}, \Z_p(n) ) \}_{n \geq 0}$, which we refer to as {\it syntomic complexes}. In \S\ref{subsection:syntomic-c1}, we use
the prismatic logarithm of \S\ref{section:twist-and-log} to construct a map
\begin{equation}\label{equation:main-theorem} c_{1}^{\Syn}: \RGamma_{ \textnormal{\'{e}t} }( \mathfrak{X}, \mathbf{G}_{m} )[-1]  \rightarrow
\RGamma_{ \Syn}( \mathfrak{X}, \Z_p(1) ), \end{equation}
and show that it induces an isomorphism after $p$-completion (Theorem \ref{theorem:syntomic-chern-class-isomorphism}).
When $\mathfrak{X}$ is a quasisyntomic $\F_p$-scheme, this construction refines the usual crystalline
first Chern class (Proposition \ref{proposition:crystalline-agreement}), whose construction we review in \S\ref{subsection:c1-crystalline}. In this case, we give a direct argument that $c_{1}^{\Syn}$ becomes an isomorphism after $p$-completion (Theorem \ref{theorem:syntomic-chern-class-isomorphism}); our proof ultimately relies on a concrete
calculation with crystalline period rings (Theorem \ref{theorem:logarithm-sequence}) which we explain in \S
\ref{subsection:logarithm-sequence}. Our proof of Theorem \ref{theorem:syntomic-chern-class-isomorphism}
in general will use formal arguments to reduce to the case where $\mathfrak{X}$ is an $\F_p$-scheme
exploiting a certain $p$-adic continuity property of the construction $R \mapsto \RGamma_{ \textnormal{\'{e}t} }( \Spf(R), \mathbf{G}_{m} )$ which we explain in \S\ref{subsection:G-m-cohomology} (see Proposition \ref{proposition:derived-descent-Gm}).

\subsection{The Logarithm Sequence}\label{subsection:logarithm-sequence}

Let $R$ be a semiperfect $\F_p$-algebra, let $R^{\flat}$ denote the perfect $\F_p$-algebra given by the inverse limit of the tower
$$ \cdots \xrightarrow{ \varphi_{R} } R \xrightarrow{ \varphi_{R} } R \xrightarrow{ \varphi_{R} } R,$$
and let $A_{\crys}(R)$ be as in Construction \ref{construction:Acrys}. For every element $x \in R^{\flat}$, we abuse notation
by writing $[x]$ for the image of the Teichm\"{u}ller representative of $x$ under the tautological map $W(R^{\flat} ) \rightarrow A_{\crys}(R)$.
Note that, if $x^{\sharp} = 1$ in $R$, then $[x] - 1$ belongs to the divided power ideal $I_{\crys}(R) = \ker( A_{\crys}(R) \twoheadrightarrow R)$,
so the logarithm $\log[x]$ is well-defined as an element of $A_{\crys}(R)$. Moreover, the Frobenius on $R$ induces a ring homomorphism
$\varphi$ from $A_{\crys}(R)$ to itself which satisfies the identities
$$ \varphi( [x] ) = [x]^{p} \quad \quad \varphi( \log[x] ) = p \log[x].$$
It follows that $\log[x]$ belongs to the ideal $$\Fil^{1}_{\Nyg} A_{\crys}(R) = \{ y \in A_{\crys}(R): \varphi(y) \in p A_{\crys}(R) \}.$$
The goal of this section is to prove the following:

\begin{theorem}\label{theorem:logarithm-sequence}
Let $R$ be a quasiregular semiperfect $\F_p$-algebra, and let $J$ denote the kernel of the homomorphism $R^{\flat} \twoheadrightarrow R$.
Then the sequence of abelian groups
\begin{equation}\label{equation:first-exact}
\xymatrix@R=50pt@C=50pt{ 0 \ar[r] &  (1 + J)^{\times} \ar[r]^-{\log[\bullet]}&  \Fil^{1}_{\Nyg} A_{\crys}(R) \ar[r]^-{ \varphi/p - 1} & A_{\crys}(R) \ar[r] & 0 }
\end{equation}
is exact.
\end{theorem}

\begin{proof}
To simplify the notation, we write $A$ for the commutative ring $A_{\crys}(R)$ and $\overline{A}$ for the quotient ring $A/pA$.
For $n \geq 0$, let $\Fil_{n}^{\conj} \overline{A}$ denote the $n$th stage of the conjugate filtration of $\overline{A}$ (Construction \ref{construction:conjugate-filtration}),
and let $\gr_{n}^{\conj} \overline{A}$ denote the quotient $\Fil_{n}^{\conj} \overline{A} / \Fil_{n-1}^{\conj} \overline{A}$. Note that the composite map
$\Fil^{1}_{\Nyg} A \xrightarrow{ \varphi / p } A \twoheadrightarrow \overline{A}$
factors through the subgroup $\Fil_{1}^{\conj} \overline{A} \subseteq \overline{A}$, and therefore induces a homomorphism
$$F: \Fil^{1}_{\Nyg} \overline{A} / p \Fil^{1}_{\Nyg} \overline{A} \rightarrow \Fil_{1}^{\conj} \overline{A}.$$
Similarly, the homomorphism $x \mapsto \log[x]$ reduces modulo $p$ to a homomorphism
$$\overline{\log}: (1+J) / (1+\varphi(J)) \rightarrow (\Fil^{1}_{\Nyg} A) / p(\Fil^{1}_{\Nyg} A).$$ Since the abelian groups
$(1+J)^{\times}$, $\Fil^{1}_{\Nyg} A$, and $A$ are $p$-complete and $p$-torsion-free, the exactness of (\ref{equation:first-exact}) is equivalent
to the exactness of the sequence
\begin{equation}\label{equation:second-exact}
\xymatrix@R=50pt@C=30pt{ 0 \ar[r] & (1+J) / (1+ \varphi(J) ) \ar[r]^-{ \overline{\log} } & (\Fil^{1}_{\Nyg} A) / p( \Fil^{1}_{\Nyg} A) \ar[r]^-{M-1} &  \overline{A} \ar[r] & 0. }
\end{equation}

Since the ring $R$ has characteristic $p$, the tautological quotient map $A \twoheadrightarrow R$ factors through a surjection $\epsilon: \overline{A} \twoheadrightarrow R$.
Moreover, the kernel $\ker(\epsilon)$ can be identified with the quotient $(\Fil^{1}_{\Nyg} A) / pA$, so we have a short exact sequence of abelian groups
$$ 0 \rightarrow A / \Fil^{1}_{\Nyg} A \xrightarrow{p} (\Fil^{1}_{\Nyg} A) / p (\Fil^{1}_{\Nyg} A) \rightarrow \ker(\epsilon: \overline{A} \rightarrow R) \rightarrow 0.$$
By virtue of Remark \ref{remark:short-exact-sequence-Nygaard-conjugate}, the composition
$$ A / \Fil^{1}_{\Nyg} A  \xrightarrow{p}  (\Fil^{1}_{\Nyg} A) / p (\Fil^{1}_{\Nyg} A) \xrightarrow{M} \Fil_{1}^{\conj} \overline{A}$$
is a monomorphism, whose image is equal to the subgroup $\Fil_{0}^{\conj} \overline{A}$. In particular, $M$ induces a homomorphism
$\overline{M}: \ker(\epsilon) \rightarrow \gr_{1}^{\conj} \overline{A}$. Moreover, the exactness of (\ref{equation:second-exact}) is equivalent to the exactness of
the sequence
\begin{equation}\label{equation:third-exact}
\xymatrix@R=50pt@C=30pt{ 0 \ar[r] & (1+J) / (1+ \varphi(J) ) \ar[r]^-{ \overline{\log} } & \ker(\epsilon: \overline{A} \rightarrow R) \ar[r]^-{ \overline{M}-1} & \overline{A} / \Fil_{0}^{\conj} \overline{A}  \ar[r] & 0. }
\end{equation}
Here we abuse notation by identifying $\overline{\log}$ with the composite homomorphism $$(1+J) / (1+\varphi(J) ) \xrightarrow{ \overline{\log} }
(\Fil^{1}_{\Nyg} A) / p (\Fil^{1}_{\Nyg} A) \twoheadrightarrow \ker(\epsilon: \overline{A} \rightarrow R).$$

Let $x$ be an element of the ideal $J \subseteq R^{\flat}$. It follows from an elementary calculation that the logarithm $\log[1-x] \in A$ is given by the $p$-adically convergent
sum
$$ \sum_{ \alpha > 0} \frac{ -[ x^{\alpha} ]}{\alpha},$$
where $\alpha$ ranges over all positive elements of the ring $\Z[1/p]$. Note that, when $\alpha$ is not an integer, the expression $[ x^{\alpha} ] / \alpha$ belongs to the ideal
$pA$. Moreover, if $\alpha$ is an integer $> p$, then the expression $[ x^{\alpha} ] / \alpha$ belongs to the ideal $pA$ by virtue of the fact that $[x]$ has divided powers in
$A$. We therefore have a congruence
\begin{equation}\label{equation:congruence-for-log} \log[1-x] \equiv \sum_{d=1}^{p} \frac{-[ x^{d} ]}{d} \pmod{p} \end{equation}
It follows that the homomorphism $\overline{\log}: (1+J) / (1+ \varphi(J) ) \rightarrow \ker( \epsilon: \overline{A} \rightarrow R)$ takes values in the intersection
$$ \ker(\epsilon: \overline{A} \rightarrow R) \cap \Fil_{1}^{\conj} \overline{A} = \ker( \epsilon: \Fil_{1}^{\conj} \overline{A} \rightarrow R).$$
Moreover, we have a monomorphism of short exact sequences
$$ \xymatrix@R=50pt@C=30pt{ 0 \ar[r] & \ker( \epsilon: \Fil_{1}^{\conj} \overline{A} \rightarrow R ) \ar[r] \ar[d]^{\overline{M}-1} & \ker(\epsilon: \overline{A} \rightarrow R) \ar[d]^{\overline{M}-1} 
\ar[r] & \overline{A} / \Fil_{1}^{\conj} \overline{A} \ar[r] \ar[d]^{-1} & 0 \\
0 \ar[r] & \gr_{1}^{\conj} \overline{A} \ar[r] & \overline{A} / \Fil_{0}^{\conj} \overline{A} \ar[r] & \overline{A} / \Fil_{1}^{\conj} \overline{A}, \ar[r] & 0 }$$
where the right vertical map is given by multiplication by $-1$ and is therefore an isomorphism. Consequently, the exactness of (\ref{equation:third-exact}) is 
equivalent to the exactness of the sequence
\begin{equation}\label{equation:fourth-exact}
\xymatrix@R=50pt@C=30pt{ 0 \ar[r] & (1+J) / (1+ \varphi(J) ) \ar[r]^-{ \overline{\log} } & \ker(\epsilon: \Fil_{1}^{\conj} \overline{A} \rightarrow R) \ar[r]^-{ \overline{M}-1} & \gr^{1}_{\conj}(A) \ar[r] & 0. }
\end{equation}

We next claim that the homomorphism $\overline{\log}$ carries $(1+J^2) / (1 + \varphi(J) )$ into the subgroup
$\Fil_{0}^{\conj} \overline{A}$. Using the congruence (\ref{equation:congruence-for-log}), we are reduced to proving that for every element $x \in J^2$,
the map $W(R^{\flat} ) \rightarrow A \twoheadrightarrow \overline{A}$
carries $[x]^{p}/p$ to an element of the group $\Fil_{0}^{\conj} \overline{A} = \im( W(R^{\flat} ) \rightarrow \overline{A} )$. In fact, we prove a slightly more precise claim (which will be useful below): for any ideal $I \subseteq J$, if $x$ belongs to $I^{2}$, then there exists an element $y \in I^{2p}$ such that
$[x]^{p}/[$ and $[y]$ have the same image in $\overline{A}$. Note that if $x$ decomposes as a sum $x' + x''$ for $x',x'' \in I^2$,
then the congruence $[x] \equiv [x'] + [x''] \pmod{p}$ implies
$$\frac{[x]^{p}}{p} \equiv \frac{[x']^{p}}{p} + \frac{ [x'']^{p} }{p} + [ \sum_{i=1}^{p-1}   \frac{ (p-1)!}{i! (p-i)!} x'^{i} x''^{p-i} ]. \pmod{p}$$
We may therefore reduce to the case where $x = uv$ for $u,v \in I \subseteq J$. In this case, we have $\frac{ [x]^{p} }{p} = p \frac{ [u]^{p} }{p} \frac{[v]^{p} }{p} \equiv 0 \pmod{p}$,
since the Teichm\"{u}ller representatives $[u]$ and $[v]$ both admit divided powers in $A$.

The preceding argument supplies a commutative diagram of short exact sequences
$$ \xymatrix@R=50pt@C=25pt{ 0 \ar[r] & (1+J^2) / (1+\varphi(J) ) \ar[d]^{\overline{\log} } \ar[r] & (1+J) / (1+ \varphi(J) ) \ar[d]^{ \overline{\log} } \ar[r] & (1+J) / (1+J^2) \ar[d]^{\rho} \ar[r] & 0 \\
0 \ar[r] & \ker( \epsilon: \Fil_{0}^{\conj} A \rightarrow R) \ar[r] & \ker(\epsilon: \Fil_{1}^{\conj} A \rightarrow R) \ar[r] & \gr_{1}^{\conj} \overline{A} \ar[r] & 0. }$$
Moreover, the congruence (\ref{equation:congruence-for-log}) shows that, for each element $x \in J$, the homomorphism $\rho$ carries the residue class of
$1-x$ to the residue class of $- [x]^{p}/p$. Using Lemma \ref{lemma:conjugate-filtration}, we deduce that $\rho$ is an isomorphism.
Consequently, the exactness of (\ref{equation:fourth-exact}) is equivalent to the exactness of the sequence
\begin{equation}\label{equation:fifth-exact}
\xymatrix@R=50pt@C=30pt{ 0 \ar[r] & (1+J^2) / (1+ \varphi(J) ) \ar[r]^-{ \overline{\log} } & \ker(\epsilon: \gr_{0}^{\conj} \overline{A} \rightarrow R) \ar[r]^-{ \overline{M} } & \gr^{1}_{\conj} \overline{A} \ar[r] & 0. }
\end{equation}

By construction, $\gr_{0}^{\conj} \overline{A}$ is equal to the image of the tautological map $R^{\flat} \rightarrow \overline{A}$. 
It follows from Remark \ref{remark:short-exact-sequence-Nygaard-conjugate} that the kernel of this map is the ideal $\varphi(J) \subseteq R^{\flat}$,
so that $\ker(\epsilon: \gr_{0}^{\conj} \overline{A} \rightarrow R)$ can be identified with the quotient $J/ \varphi(J)$ (regarded as an abelian group under addition).
Applying Lemma \ref{lemma:conjugate-filtration}, we see that the map $\overline{M}: J/\varphi(J) \rightarrow \gr^{1}_{\conj} \overline{A}$ is a surjection
whose kernel is the submodule $J^2 / \varphi(J)$. Consequently, to complete the proof, it will suffice to show that $\overline{\log}$ determines
an isomorphism of abelian groups
$$ (1+J^2) / (1+ \varphi(J) ) \rightarrow J^2 / \varphi(J)$$
(where the group structure on the left is given by multiplication and the group structure on the right is given by addition). We first prove injectivity.
Suppose that $x$ is an element of $J^2$ such that $\log[1-x]$ belongs to $pA$; we wish to show that $x$ belongs to $\varphi(J)$.
Choose a finitely generated ideal $I \subseteq J$ such that $x \in I^2$. By the argument given above, we can choose an element $y \in I^{2p}$ such
that $\frac{ [x]^{p}}{p}$ and $[y]$ have the same image in $\overline{A}$. Using the congruence (\ref{equation:congruence-for-log}), we conclude that the sum
$\sum_{d=1}^{p-1} \frac{ x^{d} }{d}$ belongs to $I^{2p} + \varphi(J)$, so that $x \in I^{4} + \varphi(J)$. Subtracting an element of $\varphi(J)$ from $x$,
we may assume that $x \in I^{4}$. Repeating the above argument, we can assume that $x \in I^{n}$ for $n$ arbitrarily large. We conclude by observing that,
if the ideal $I$ is generated by $m$ elements, then we have $I^{pm} \subseteq \varphi(I) \subseteq \varphi(J)$.

We now prove that the map $\overline{\log}: (1+J^2) / (1 + \varphi(J) ) \rightarrow J^2 / \varphi(J)$ is surjective. Since $\overline{\log}$ is a group homomorphism,
it will suffice to show that its image contains the residue class of every product $uv$ for $u,v \in J$. Let $J'$ denote the principal ideal $(uv) \subseteq R^{\flat}$
and let $\overline{J}$ denote its image in the quotient ring $R^{\flat} / \varphi(J)$. Since every element $x \in J'$ satisfies $\frac{ [x]^{p} }{p} \in pA$, the congruence (\ref{equation:congruence-for-log}) implies that the composite map
$$ \overline{J} \xrightarrow{ x \mapsto 1-x } ( 1 + J^2 ) / (1 + \varphi(J) ) \xrightarrow{-\overline{\log}} J^2 / \varphi(J)$$
is given by the power series $f(x) = \sum_{d=1}^{p-1} \frac{ x^d }{d}$. Since $f(x) \equiv x \pmod{x^2}$, the power
series $f$ is invertible with respect to composition. Let $g(y)$ denote the inverse power series. Since $\overline{J}$ is a nilpotent ideal of $R^{\flat} / \varphi(J)$,
it follows that the function $x \mapsto f(x)$ determines a bijection from $\overline{J}$ to itself (with inverse given by $y \mapsto g(y)$). In particular, the image of
$f$ contains the residue class of the product $uv$.
\end{proof}

\subsection{Cohomology with \texorpdfstring{$\mathbf{G}_m$}{Gm}-Coefficients}\label{subsection:G-m-cohomology}

Let $X$ be a scheme or formal scheme. We write $\RGamma_{ \textnormal{\'{e}t} }( X, \mathbf{G}_{m} )$ for the derived global sections of the functor $U \mapsto \calO_{X}(U)^{\times}$,
which we view as a sheaf of abelian groups on the \'{e}tale site of $X$. For each integer $n$, let $\mathrm{H}^{n}_{ \textnormal{\'{e}t} }( X, \mathbf{G}_{m} )$ denote the $n$th cohomology
group of the complex $\RGamma_{ \textnormal{\'{e}t} }( X, \mathbf{G}_{m} )$. In particular, we have
$$\mathrm{H}^{n}_{ \textnormal{\'{e}t} }( X, \mathbf{G}_{m} ) = \begin{cases} 0 & \text{ if $n < 0$} \\
\Gamma(X, \calO_{X})^{\times} & \text{ if $n=0$ } \\
\Pic(X) & \text{ if $n=1$.} \end{cases}$$

\begin{remark}\label{remark:Frobenius-G-m-cohomology}
Let $X$ be an $\F_p$-scheme, and let $\varphi_{X}: X \rightarrow X$ denote the absolute Frobenius morphism. Then the pullback map
$\varphi^{\ast}_{X}: \RGamma_{ \textnormal{\'{e}t} }( X, \mathbf{G}_{m} ) \rightarrow \RGamma_{ \textnormal{\'{e}t} }( X, \mathbf{G}_{m} )$ is
given by multiplication by $p$.
\end{remark}

\begin{remark}\label{remark:Tate-module-appearance}
Let $X$ be a scheme. A classical result of Grothendieck (Theorem~11.7 of \cite{MR244271}) asserts that the \'{e}tale cohomology of $X$ with coefficients in $\mathbf{G}_m$
coincides with the fppf cohomology of $X$ with coefficients in $\mathbf{G}_m$. Consequently, for each integer $n$, the short exact sequence of fppf sheaves
$$ 0 \rightarrow \mu_{p^{n}} \rightarrow \mathbf{G}_m \xrightarrow{p^{n}} \mathbf{G}_m \rightarrow 0$$
induces a fiber sequence
$$ \RGamma_{\textnormal{fppf}}( X, \mu_{p^{n}} ) \rightarrow \RGamma_{\textnormal{\'{e}t}}(X, \mathbf{G}_m) \xrightarrow{p^{n}} \RGamma_{\textnormal{\'{e}t}}(X, \mathbf{G}_m).$$
Passing to the homotopy limit over $n$, we obtain an isomoprhism
\begin{equation}
\label{eq:GmpCompInverseLimit}
 \RGamma_{\textnormal{\'{e}t}}( X, \mathbf{G}_m)^{\wedge} \simeq \varprojlim_{n} \RGamma_{\textnormal{fppf}}( X, \mu_{p^{n}} )[1],
 \end{equation}
where $\RGamma_{\textnormal{\'{e}t}}( X, \mathbf{G}_m)^{\wedge}$ denotes the $p$-completion of the complex
$\RGamma_{\textnormal{\'{e}t}}( X, \mathbf{G}_m)$. In particular:
\begin{itemize}
\item The cohomology groups of the complex $\RGamma_{\textnormal{\'{e}t}}( X, \mathbf{G}_m)[-1]^{\wedge}$ are concentrated in degrees $\geq 0$.
\item The $0$th cohomology group of the complex $\RGamma_{\textnormal{\'{e}t}}( X, \mathbf{G}_m)[-1]^{\wedge}$ is the
Tate module $T_p( A^{\times} )$, where $A = \Gamma( X, \calO_X)$ is the ring of global functions on $X$.
\item The $1$st cohomology group of the complex $\RGamma_{\textnormal{\'{e}t}}( X, \mathbf{G}_m)^{\wedge}$ can be identified
with the set of isomorphism classes of torsors for the group scheme $\varprojlim \mu_{p^{n}}$ (with respect to the fpqc topology on $X$).
\end{itemize}
\end{remark}

\begin{remark}[$\RGamma_{\mathet}(-,\mathbf{G}_m)$ for $p$-adic formal schemes]
\label{remark:compute-RGamma-G-m-formal}
Let $\mathfrak{X}$ be a bounded $p$-adic formal scheme, which we can write as a colimit of closed
subschemes $\mathfrak{X}_{n} = \Spec( \Z / p^{n} \Z) \times \mathfrak{X}$. Then
the restriction map
$$ \RGamma_{\mathet}( \mathfrak{X}, \mathbf{G}_m) \rightarrow \varprojlim_{n} \RGamma_{\mathet}( \mathfrak{X}_n, \mathbf{G}_m )$$
is an isomorphism in the $\infty$-category $\calD(\Z)$.
\end{remark}

\begin{remark}[Pro-fppf descent for $\RGamma_{\mathet}(-,\mathbf{G}_m)$]
\label{remark:profppf-Gm}
A theorem of Grothendieck (see Remark~\ref{remark:Tate-module-appearance}) guarantees that the $\calD^{\geq 0}(\Z)$-valued functor $\RGamma_{\mathet}(-,\mathbf{G}_m)$ is a sheaf for the fppf topology on the category of schemes. Moreover, given a scheme $X = \varprojlim_i X_i$ presented as a cofiltered limit of qcqs schemes $X_i$ along affine transition maps, the natural map 
\[ \colim_i \RGamma_{\mathet}(X_i,\mathbf{G}_m) \to \RGamma(X,\mathbf{G}_m)\]
is an equivalence in $\calD^{\geq 0}(\Z)$. As filtered colimits commute with totalizations for diagrams in $\calD^{\geq 0}(\Z)$, it follows that the $\calD^{\geq 0}(\Z)$-valued functor $\RGamma_{\mathet}(-,\mathbf{G}_m)$ is even a sheaf for the pro-fppf topology (that is, the Grothendieck topology on schemes generated by the fppf topology
together with maps $\Spec(B) \rightarrow \Spec(A)$, where $B$ can be realized as a filtered colimit of faithfully flat $A$-algebras of finite presentation).
\end{remark}

\begin{proposition}\label{proposition:degree-zero-qrsp}
Let $R$ be a quasiregular semiperfect $\F_p$-algebra. Then
$\RGamma_{\textnormal{\'{e}t}}( X, \mathbf{G}_m)[-1]^{\wedge}$ is concentrated in cohomological degree $0$, and
can therefore be identified with the Tate module $T_p(R^{\times})$.
\end{proposition}

\begin{proof}
We wish to show that the tautological map
$$ \theta: T_p(R^{\times} ) \rightarrow \RGamma_{\textnormal{\'{e}t}}( X, \mathbf{G}_m)[-1]^{\wedge}$$
is an isomorphism in $\widehat{\calD}(\Z_p)$. It suffices to prove this after reduction modulo $p$.
Using the exactness of the sequence
\begin{equation} \label{equation:T-p-exact}
\xymatrix@R=50pt@C=50pt{ 0 \ar[r] & T_p(R^{\times}) \ar[r]^-{p} &  T_p(R^{\times}) \ar[r]r &  \mu_p(R) \ar[r] & 0 }
\end{equation}
and Remark \ref{remark:Tate-module-appearance}, we are reduced to showing that the tautological map
$$ \mu_p(R) \rightarrow \RGamma_{\fppf}( \Spec(R), \mu_p )$$
is an isomorphism: that is, that the cohomology groups $\mathrm{H}^{n}_{\fppf}( \Spec(R), \mu_p )$ vanish for $n > 0$.
Since every \'{e}tale $R$ algebra is semiperfect, the sequence of abelian sheaves
$$ 0 \rightarrow \mu_p \rightarrow \mathbf{G}_m \xrightarrow{p} \mathbf{G}_m \rightarrow 0$$
is exact on the \'{e}tale site of $\Spec(R)$; we can therefore identify $\mathrm{H}^{n}_{\fppf}( \Spec(R), \mu_p )$
with the \'{e}tale cohomology group $\mathrm{H}^{n}_{\textnormal{\'{e}t}}( \Spec(R), \mu_p )$. To prove that
these groups vanish, it will suffice to show that the functor $R \mapsto \mu_p(R)$ satisfies \'{e}tale
descent when regarded as a functor from $\CAlg^{\qrsp}_{\F_p}$ to the derived $\infty$-category $\calD(\F_p)$.
Using the sequence (\ref{equation:T-p-exact}) again, we are reduced to proving the analogous property
for the functor $R \mapsto T_p(R)$. This follows from the exact sequence
$$ 0 \rightarrow T_p(R^{\times}) \xrightarrow{ \log[\bullet]} \Fil^{1}_{\Nyg} A_{\crys}(R) \xrightarrow{ \varphi/p - \id} A_{\crys}(R) \rightarrow 0.$$
of Theorem \ref{theorem:logarithm-sequence}, since the functor $R \mapsto A_{\crys}(R) \simeq \RGamma_{\crys}(\Spec(R)/ \Z_p)$ satisfies
$p$-quasisyntomic descent (Corollary \ref{corollary:qs-descent-for-crys}).
\end{proof}

\begin{variant}\label{variant:vanishing-H}
Let $X$ be a semiperfect $\F_p$-scheme; fix an integer $n$. Then the sequence of abelian sheaves
$$ 0 \rightarrow \mu_{p^n} \rightarrow \mathbf{G}_m \xrightarrow{p^n} \mathbf{G}_m \rightarrow 0$$
is exact on the \'{e}tale site of $X$: for any \'etale map $\mathrm{Spec}(R) \to X$ the ring $R$ is semiperfect as $X$ is so, whence the map $\mathbf{G}_m(R) \xrightarrow{p^n} \mathbf{G}_m(R)$ is surjective. Consequently, for every integer $n$, the canonical map
$$ ( \Z / p^{n} \Z) \otimes^{L} \RGamma_{\textnormal{\'{e}t}}(X, \mathbf{G}_m)[-1] \rightarrow \RGamma_{\textnormal{\'{e}t}}(X, \mu_{p^{n}} )$$
is an isomorphism. 
\end{variant}

%
%

It will be convenient to enlarge the domain of the functor $R \mapsto \RGamma_{\textnormal{\'{e}t}}( \Spec(R), \mathbf{G}_m)$.

\begin{notation}
Let $R$ be an animated commutative ring. We let $\mathbf{G}_m(R)$ denote the {\it multiplicative group} of $R$, which we regard as an object of the $\infty$-category $\calD(\Z)^{\leq 0}$
(that is, as an animated abelian group). It is characterized (up to canonical isomorphism) by the requirement that for every free abelian group $\Lambda$ of finite rank, we have
a homotopy equivalence $$ \Hom_{ \calD(\Z) }( M, \mathbf{G}_m(R) ) \simeq \Hom_{ \CAlg^{\anim} }( \Z[M], R ),$$
where $\Z[M]$ denotes the group algebra of $M$ over $\Z$.
\end{notation}

\begin{remark}\label{remark:Kan-extension-multiplicative}
Let $k$ be a commutative ring and let $\CAlg^{\anim}_{k}$ denote the $\infty$-category of animated $k$-algebras.
Then the functor $\mathbf{G}_m: \CAlg^{\anim}_{k}  \rightarrow \calD(\Z)$ is a left Kan extension of its restriction to full subcategory of
$\CAlg^{\anim}_{k}$ spanned by Laurent polynomial algebras $k[ x_1^{\pm 1}, \cdots, x_n^{\pm 1} ]$. Beware that it is
{\em not} a left Kan extension of its restriction to the category of polynomial rings over $k$.
\end{remark}

\begin{remark}\label{remark:homotopy-of-Gm}
Let $R$ be an animated commutative ring. Then the cohomology groups of the complex $\mathbf{G}_m(R)$ are related to the
homotopy groups of $R$ by the formula
$$ \mathrm{H}^{n}( \mathbf{G}_m(R) ) = \begin{cases} \pi_0(R)^{\times} & \textnormal{ if $n=0$ } \\
\pi_{-n}(R) & \textnormal{ if } n < 0 \\
0 & \textnormal{ if } n > 0.\end{cases}$$
In particular, if $R$ is an ordinary commutative ring, then $\mathbf{G}_m(R)$ can be identified with the abelian group $R^{\times}$ of units in $R$
(which we identify with an object of the derived $\infty$-category $\calD(\Z)$). 
\end{remark}

\begin{notation}\label{notation:etale-sheafification}
If we regard the construction $R \mapsto \mathbf{G}_m(R)$ as a functor from $\CAlg^{\anim}$ to the $\infty$-category $\calD(\Z)^{\leq 0}$, then it
satisfies descent with respect to the \'{e}tale topology. However, if we regard $\mathbf{G}_m$ as a functor taking values in the larger $\infty$-category $\calD(\Z)$, 
then it does {\em not} satisfy descent for the \'{e}tale topology. We will denote its \'{e}tale sheafification by $R \mapsto \RGamma_{\textnormal{\'{e}t}}( \Spec(R), \mathbf{G}_m)$.
We let $\RGamma_{\textnormal{\'{e}t}}( \Spec(R), \mathbf{G}_m)^{\wedge}$ denote the $p$-completion of the complex $\RGamma_{\textnormal{\'{e}t}}( \Spec(R), \mathbf{G}_m)$.
\end{notation}

\begin{proposition}\label{proposition:Hodge-c1-affine}
For every animated commutative ring $R$, there is a canonical map
$$ c_{1}^{\Hodge}: \RGamma_{\mathet}( \Spec(R), \mathbf{G}_m) \rightarrow L \Omega^{1}_{R},$$
characterized by the fact that it depends functorially on $R$ and that it is given on cohomology in degree zero by the construction
$$ (u \in \GL_1(R)) \mapsto (\dlog(u) = \frac{ du}{u} \in \Omega^{1}_{R} )$$
when $R$ is an ordinary commutative ring.
\end{proposition}

\begin{proof}
Since the functor $R \mapsto L \Omega^{1}_{R}$ satisfies descent for the \'{e}tale topology, it will suffice to show
that the construction $u \mapsto \dlog(u)$ extends uniquely to a natural transformation
$\mathbf{G}_m(R) \rightarrow L \Omega^{1}_{R}$. By virtue of Remark \ref{remark:Kan-extension-multiplicative}, we
may assume that $R$ is a smooth $\Z$-algebra, in which case $\mathbf{G}_m(R)$ and $L \Omega^{1}_{R}$
are concentrated in cohomological degree zero and there is nothing to prove.
\end{proof}

\begin{notation}\label{notation:hodge-c1}
Let $X$ be a scheme, formal scheme, or algebraic stack. Globalizing Proposition \ref{proposition:Hodge-c1-affine}, we obtain a map
$$ c_1^{\Hodge}: \RGamma_{\mathet}(X, \mathbf{G}_m) \rightarrow \RGamma(X, L \Omega^{1}_{X} ),$$
which we will refer to as the {\it Hodge first Chern class}. Passing to cohomology in degree $1$, we obtain a homomorphism of abelian groups
$\Pic(X) \rightarrow \mathrm{H}^{1}( X, L \Omega^{1}_{X} )$, which we also denote by $c_{1}^{\Hodge}$.
\end{notation}

\begin{proposition}[$p$-adic continuity, \v{C}esnavi\v{c}ius-Scholze {\cite[Theorem~5.3.4]{cesnavicius2019purity}}]
\label{proposition:p-adic-continuity}
Let $R$ be an animated commutative ring, and write $R^h$ for the $p$-henselization of $R$. For each integer $n \geq 0$, let
$R_{n}$ denote the derived tensor product $(\Z/p^{n}\Z) \otimes^{L} R$.  Let $G$ be a finite locally free commutative $R$-group scheme $G$ obtained via base change from a discrete ring mapping to $R$.  The restriction map
$$ \RGamma_{\fppf}( \Spec(R^h), G) \rightarrow \varprojlim_{n} \RGamma_{\fppf}(  \Spec(R_n), G) $$
is an isomorphism in the $\infty$-category $\calD(\Z)$.
\end{proposition}

We give a slightly different proof of this result than the one in \cite{cesnavicius2019purity}: our proof involves only \'etale localization (in contrast with the the fppf localization strategy used in \cite{cesnavicius2019purity}), but requires a stronger local input.

\begin{proof}[Proof of Proposition \ref{proposition:p-adic-continuity}]
Using animated deformation theory as in the first paragraph of the proof of \cite[Theorem 5.3.5]{cesnavicius2019purity}, we may assume $R$ is discrete (whence $R^h$ is also discrete). 

We begin by reducing the theorem to special $R$'s where most of the higher cohomology vanishes. For this, note that the functor $R \mapsto \RGamma_{\fppf}(\Spec(R), G)$ is an ind-fppf (and hence a pro-\'etale) sheaf. In particular, pro-Nisnevich descent gives a pullback square
\[ \xymatrix@R=50pt@C=50pt{ \RGamma_{\fppf}(\Spec(R), G) \ar[r] \ar[d] & \RGamma_{\fppf}(\Spec(R^h), G) \ar[d] \\
	\RGamma_{\fppf}(\Spec(R[1/p]), G) \ar[d]^{\sim} 
) \ar[r] & \RGamma_{\fppf}(\Spec(R^h[1/p]), G) \ar[d]^{\sim}  \\
	\RGamma_{\mathet}(\Spec(\pi_0(R[1/p])), G) \ar[r] & \RGamma_{\mathet}(\Spec(\pi_0(R^h)[1/p]), G)}
\]
where the lower vertical isomorphisms are obtained by observing that $G$ becomes a finite \'etale group scheme after inverting $p$ and using the identification of fppf and \'etale cohomology with \'etale coefficients, together with the identification of the \'etale site of an animated ring with that of its $\pi_0$. Regarding this as a diagram of $\calD(\Z)$-valued functors of the  ring $R$, the terms on left are pro-\'etale sheaves by the first sentence of this paragraph. Moreover, by $\arc_p$-descent, the same holds for the bottom right term (see \cite[Theorem 6.11, Corollary 6.17]{arcs}). Consequently, by the above pullback square, the functor $R \mapsto \RGamma_{\fppf}( \Spec(R^h), G)$ is a pro-\'etale sheaf. As the same is also true for the functors $R \mapsto  \RGamma_{\fppf}(  \Spec(R_n), G)$ for all $n$ as well as  their limit, we may pass to a pro-\'etale cover of $R$ to assume $R$ is \'etale-local, i.e., every faithfully flat \'etale map $R \to S$ admits a section.  In this case, we claim that each $R_n$ as well as $R^h$ are also \'etale local. The \'etale-locality of $R_n$ follows by observing that any faithfully flat \'etale map $R_n \to T$ has a lift to an \'etale map $R \to S$ \cite[Tag 04D1]{stacks-project} and hence to a faithfully flat \'etale map $R \to S \times R[1/p]$. This also implies that each $\pi_0(R_n)$ is \'etale-local as the \'etale site is insensitive to replacing an animated ring with its $\pi_0$. The \'etale-locality of $R^h$ then follows as $R^h$ is henselian along the kernel of $R \to \pi_0(R_1)$.

Thanks to the previous paragraph, we have reduced to checking the theorem when $R$, $R^h$ and each $R_n$ are \'etale-local. Using the Begueri resolution, one then learns that $\tau^{\leq 1} \RGamma_{\fppf}(\Spec(R^h),G) \to \RGamma_{\fppf}(\Spec(R^h),G)$ is an equivalence, and similarly for each $R_n$. Interpreting low degree cohomology via torsors, it is enough to show the following concrete statement: writing $BG(-)$ for the functor carrying an animated $R$-algebra to its $\infty$-groupoid of fppf $G$-torsors, the natural maps
\[ BG(R^h) \xrightarrow{\alpha} BG(\widehat{R}) \xrightarrow{\beta} \varprojlim BG(R_n)\]
are equivalences (where $\widehat{R} = \varprojlim_n R_n$ is the $p$-completion of $R$ or equivalently of $R^h$). The claim for $\beta$ follows by interpreting $G$-torsors in terms of their (locally free) $R$-algebra of functions and using the symmetric-monoidal equivalence $\mathrm{Vect}(\widehat{R}) \simeq \varprojlim_n \mathrm{Vect}(R_n)$ coming from base change. For $\alpha$, one argues similarly using Beauville-Laszlo gluing  (see \cite[Theorem 7.4.0.1]{LurieSAG} and \cite[Theorem 1.4]{MR3572635}) as well as the equivalence
\[ BG(R^h[1/p]) \simeq BG(\widehat{R}[1/p])\]
coming from the Fujiwara-Gabber theorem (see \cite[Theorem 6.11]{arcs}). 
\end{proof}

\begin{corollary}\label{snifflet}
Let $R$ be a $p$-complete commutative ring with bounded $p$-power torsion. 
Then the restriction map
$$ \RGamma_{\mathet}( \Spec(R), \mathbf{G}_m)^{\wedge}  \rightarrow
\RGamma_{\mathet}( \Spf(R), \mathbf{G}_m )^{\wedge}$$
is an isomorphism.
\end{corollary}

\begin{proof}
Using Remark \ref{remark:compute-RGamma-G-m-formal}, it is enough to show that
\[ \RGamma_{\mathet}( \Spec(R), \mathbf{G}_m)^{\wedge} \to \varprojlim_n \RGamma_{\mathet}( \Spec(R/p^n), \mathbf{G}_m)^{\wedge} \]
is an equivalence in $\calD(\Z)$. As both sides are $p$-complete, it suffices to check this modulo $p$. Using the formula $\mathbf{G}_m/p \simeq \mu_p[1]$ as fppf sheaves coming from the Kummer sequence (see Remark~\ref{remark:Tate-module-appearance}), we are then reduced to checking that 
\[ \RGamma_{\mathet}( \Spec(R), \mu_p)  \to \varprojlim_n \RGamma_{\mathet}( \Spec(R/p^n), \mu_p) \]
is an equivalence. As $R$ has bounded $p$-power torsion, we may replace $R/p^n$ with $R_n := R \otimes_{\mathbf{Z}}^L \mathbf{Z}/p^n$ without affecting the limit on the right. In this case, the claim follows by Proposition \ref{proposition:p-adic-continuity} applied to $G=\mu_p$.
\end{proof}

For our applications, we will need a slightly stronger version of Proposition \ref{proposition:p-adic-continuity}.

\begin{proposition}[Derived Descent for Reduction Modulo $p$]\label{proposition:derived-descent-Gm}
Let $R$ be a $p$-complete animated commutative ring, let $\F_p^{\otimes \bullet +1 }$ be the cosimplicial animated commutative ring
introduced in Notation \ref{notation:F-p-bullet}, and let $R^{\bullet}$ denote the cosimplicial $R$-algebra given by 
$R^{\bullet} = R \otimes^{L}_{\Z} \F_p^{\otimes \bullet +1 }$. Then the tautological map
$$ \RGamma_{\textnormal{\'{e}t}}( \Spec(R), \mathbf{G}_m)^{\wedge} \rightarrow \Tot( \RGamma_{\textnormal{\'{e}t}}( \Spec(R^{\bullet}), \mathbf{G}_m)^{\wedge} )$$
is an isomorphism in the $\infty$-category $\widehat{\calD}(\Z_p)$.
\end{proposition}

The proof of Proposition \ref{proposition:derived-descent-Gm} will require some preliminaries. For the remainder of this section,
we fix an animated commutative ring $R$ and we let $R^{\bullet}$ denote the cosimplicial animated $R$-algebra
appearing in the statement of Proposition \ref{proposition:derived-descent-Gm}.

\begin{lemma}\label{lemma:connective-partial-totalization}
For every integer $n \geq 0$, the partial totalization
$\Tot^{n}( \mathbf{G}_m( R^{\bullet} ) ) \in \calD(\Z)$ is connective: that is, its cohomology groups vanish in positive degrees.
\end{lemma}

\begin{proof}
We proceed by induction on $n$. In the case $n = 0$, the totalization $\Tot^{n}( \mathbf{G}_m( R^{\bullet} ) )$ is given by $\mathbf{G}_m(R^{0})$, which
is connective by construction (see Remark \ref{remark:homotopy-of-Gm}). Suppose that $n > 0$, and form a fiber sequence
$$ M \rightarrow \Tot^{n}( \mathbf{G}_m(R^{\bullet} ) ) \rightarrow \Tot^{n-1}( \mathbf{G}_m( R^{\bullet} ) ).$$
By virtue of our inductive hypothesis, it will suffice to show that the cohomology groups of the complex $M$ are concentrated in degrees $\leq 0$.

Note that the identification $R^{n} \simeq R^{0} \otimes^{L}_{\F_p} ( \F_p^{\otimes n+1} )$
determines an isomorphism of the homotopy ring $\pi_{\ast}( R^{n} )$ with exterior algebra over $\pi_{\ast}( R^{0} )$ generated by
elements $e_1, e_2, \ldots, e_n \in \pi_{1}(  \F_p^{\otimes n+1})$. For every integer $m \geq 0$, let $[m]$ denote the linearly ordered set
$\{ 0 < 1 < \cdots < m \}$. Then the complex $M[n]$ can be identified with the fiber of the tautological map
$$ \mathbf{G}_m( R^{n} ) \rightarrow \varprojlim_{ [n] \twoheadrightarrow [m] } \mathbf{G}_m( R^{m} ),$$
where the limit is taken over all monotone surjections $[n] \twoheadrightarrow [m]$ which are not bijective.
An elementary calculation shows that the map $M[n] \rightarrow \mathbf{G}_m( R^{n} )$ induces an injection on cohomology,
whose image is the summand of $\mathrm{H}^{< 0}( \mathbf{G}_m(R^{n} ) ) \simeq \pi_{> 0}( R^{n} )$ consisting of elements
which are divisible by the product $e_1 e_2 \cdots e_n$. In particular, the cohomology groups of $M[n]$ are concentrated in degrees $\leq -n$,
so the cohomology groups of $M$ are concentrated in degrees $\leq 0$.
\end{proof}

\begin{lemma}\label{lemma:tot-n-iso}
For every integer $n \geq 0$, the tautological map
$$ \theta: \mathbf{G}_m( \Tot^{n}( R^{\bullet} ) ) \rightarrow \Tot^{n}( \mathbf{G}_m( R^{\bullet} ) )$$
is an isomorphism in the $\infty$-category $\calD(\Z)$.
\end{lemma}

\begin{proof}
Note that, when regarded as a functor from the $\infty$-category $\CAlg^{\anim}$ of animated commutative rings to the
$\infty$-category $\calD(\Z)^{\leq 0}$ of {\em connective} chain complexes, the functor $\mathbf{G}_m$ preserves all inverse limits
(since it is right adjoint to the group algebra functor $M \mapsto \Z[M]$). Consequently, the invertibility of the map $\theta$
is equivalent to the connectivity of the partial totalization $\Tot^{n}( \mathbf{G}_m( R^{\bullet} ) )$, which follows from
Lemma \ref{lemma:connective-partial-totalization}.
\end{proof}

\begin{lemma}\label{lemma:tot-n-iso-sheaf}
For every integer $n \geq 0$, the tautological map
$$ \theta': \RGamma_{\textnormal{\'{e}t}}( \Spec( \Tot^{n}(R^{\bullet}) ), \mathbf{G}_m) \rightarrow \Tot^{n}(  \RGamma_{\textnormal{\'{e}t}}( \Spec(R^{\bullet}), \mathbf{G}_m)$$
is an isomorphism in the $\infty$-category $\calD(\Z)$.
\end{lemma}

\begin{proof}
For every \'{e}tale $S$-algebra $R$, let $S^{\bullet}$ denote the cosimplicial animated $S$-algebra given by $R^{\bullet} \otimes^{L}_{R} S$. Since sheafification
for the \'{e}tale topology commutes with finite limits and with the formation of direct images along the closed immersions 
$$\Spec(R^{m} ) \hookrightarrow \Spec(R) \hookleftarrow \Spec( \Tot^{n}(R) ),$$
the morphism $\theta'$ is induced by a map of $\calD(\Z)$-valued sheaves on the \'{e}tale site of $\Spec(R)$, given by \'{e}tale sheafification of the construction
$$ S \mapsto (  \mathbf{G}_m( \Tot^{n}( S^{\bullet} ) ) \rightarrow \Tot^{n}( \mathbf{G}_m( S^{\bullet} ) ) ).$$
By virtue of Lemma \ref{lemma:tot-n-iso}, this map is already an isomorphism at the level of $\calD(\Z)$-valued presheaves.
\end{proof}

\begin{proof}[Proof of Proposition \ref{proposition:derived-descent-Gm}]
Let $R$ be an animated commutative ring and let $R^{\bullet}$ be the cosimplicial $R$-algebra appearing in the statement of
Proposition \ref{proposition:derived-descent-Gm}. Note that, for every integer $n \geq 0$, we can identify the partial totalization
$\Tot^{n}(R^{\bullet} )$ with the derived tensor product $R_n = R \otimes^{L}_{\Z} (\Z / p^{n+1} \Z)$. We wish to show that,
if $R$ is $p$-complete, then the composite map
\begin{eqnarray*}
\RGamma_{\textnormal{\'{e}t}}( \Spec(R), \mathbf{G}_m)^{\wedge} & \xrightarrow{\rho} & \varprojlim_{n} \RGamma_{\textnormal{\'{e}t}}( \Spec(R_n), \mathbf{G}_m)^{\wedge} \\
& = & \varprojlim_{n} \RGamma_{\textnormal{\'{e}t}}( \Spec( \Tot^{n}( R^{\bullet} ) ), \mathbf{G}_m)^{\wedge} \\
& \xrightarrow{\rho'} &  \varprojlim_{n} \Tot^{n} \RGamma_{\textnormal{\'{e}t}}( \Spec(R^{\bullet}), \mathbf{G}_m)^{\wedge}  \\
& = & \Tot \RGamma_{\textnormal{\'{e}t}}( \Spec( R^{\bullet} ), \mathbf{G}_m)^{\wedge}.
\end{eqnarray*}
is an isomorphism. Note that $\rho'$ is given by an inverse limit of the $p$-completions of comparison maps
$$\RGamma_{\textnormal{\'{e}t}}( \Spec( \Tot^{n}(R^{\bullet}) ), \mathbf{G}_m) \rightarrow \Tot^{n}(  \RGamma_{\textnormal{\'{e}t}}( \Spec(R^{\bullet}), \mathbf{G}_m),$$
each of which is an isomorphism by virtue of Lemma \ref{lemma:tot-n-iso-sheaf}. It will therefore suffice to show that $\rho$ is an isomorphism, which follows from Proposition \ref{proposition:p-adic-continuity}.
\end{proof}

\subsection{The Crystalline First Chern Class}\label{subsection:c1-crystalline}

In this section, we review the definition of Chern classes of line bundles in the setting of crystalline cohomology.

\begin{construction}\label{construction:crystalline-c1}
Let $X$ be an $\F_p$-scheme, and let $\Crys(X/\Z_p)$ denote the category of triples $(A,I,v)$, where $(A,I)$ is a $p$-complete and separated
divided power algebra over $(\Z_p, (p) )$, and $v: \Spec(A/I) \rightarrow X$ is a morphism of $\F_p$-schemes (Notation \ref{notation:big-crystalline-site}).
Let $\mathcal{O}_{\Crys(X/\Z_p)}$ denote the sheaf of commutative rings on $\Crys(X/\Z_p)^{\op}$ given by $\mathscr{O}_{\Crys(X/\Z_p)}(A,I,v) = A$,
and let $\mathscr{I} \subseteq \mathcal{O}_{\Crys(X/\Z_p)}$ denote the ideal sheaf given by $\mathscr{I}(A,I,v) = I$.
Recall that the crystalline cochain complex $\RGamma_{\crys}(X/\Z_p)$ can be defined as the limit of the diagram
$$ \calO_{\Crys(X/\Z_p)}: \Crys(X/\Z_p) \rightarrow \widehat{\calD}(\Z_p) \quad \quad (A,I,v) \mapsto A.$$
Similarly, the complex $\Fil^{1}_{\Nyg} \RGamma_{\crys}(X/\Z_p)$ of Notation \ref{notation:Nygaard-crystalline} can be identified with the limit of the diagram
$$\mathscr{I}: \Crys(X/\Z_p) \rightarrow \widehat{\calD}(\Z_p) \quad \quad (A,I,v) \mapsto I.$$

For each object $(A,I,v)$ of the category $\Crys(X/\Z_p)$, the system of divided powers on the ideal $I$ determines a logarithm map
$(1+I)^{\times} \rightarrow I$ which depends functorially $(A,I,v)$, and therefore determines a map of abelian sheaves
$\log: (1+ \mathscr{I})^{\times} \rightarrow \mathscr{I}$. Form a pushout diagram
$$ \xymatrix@R=50pt@C=50pt{ (1+\mathscr{I})^{\times} \ar[r] \ar[d]^{\log} & \mathcal{O}^{\times}_{ \Crys(X/\Z_p)}  \ar[d] \\
\mathscr{I} \ar[r] & \mathscr{F} }$$
in the category of abelian presheaves on $\Crys(X/\Z_p)$, so that we have a short exact sequence of abelian presheaves
$$ 0 \rightarrow \mathscr{I} \rightarrow \mathscr{F} \rightarrow (\mathcal{O}_{ \Crys(X/\Z_p)} / \mathscr{I} )^{\times} \rightarrow 0,$$
which determines a boundary map $\delta: (\mathcal{O}_{ \Crys(X/\Z_p)} / \mathscr{I} )^{\times}[-1] \rightarrow \mathscr{I}$
of $\calD(\Z)$-valued presheaves on $\Crys(X/\Z_p)^{\op}$. Passing to the limit over the subcategory of $\Crys(X/\Z_p)$ spanned by the \'{e}tale maps $v: \Spec(A/I) \rightarrow X$,
we obtain a map
$$ c_{1}^{\crys}: \RGamma_{\textnormal{\'{e}t}}( X, \mathbf{G}_m)[-1] \rightarrow \Fil^{1}_{\Nyg} \RGamma_{\crys}(X/\Z_p),$$
which we will refer to as the {\it crystalline first Chern class}.
\end{construction}

\begin{notation}
Let $X$ be an $\F_p$-scheme. We will abuse notation by identifying the crystalline first Chern class $c_{1}^{\crys}$ of Construction \ref{construction:crystalline-c1}
with the composition
$$  \RGamma_{\textnormal{\'{e}t}}( X, \mathbf{G}_m)[-1] \xrightarrow{c_1^{\crys}} \Fil^{1}_{\Nyg} \RGamma_{\crys}(X/\Z_p) \RGamma_{\crys}(X/\Z_p).$$
Passing to cohomology in degree $2$, we obtain a homomorphism of abelian groups
$$ c_{1}^{\crys}: \Pic(X) = \mathrm{H}^{1}_{\textnormal{\'{e}t}}( X, \mathbf{G}_m) \rightarrow \mathrm{H}^{2}_{\crys}( X/ \Z_p).$$
\end{notation}

\begin{variant}\label{variant:crystalline-factor-through-completion}
Let $X$ be an $\F_p$-scheme. Since the complex $\Fil^{1}_{\Nyg} \RGamma_{\crys}(X/\Z_p)$ is $p$-complete, the morphism
$c_{1}^{\crys}$ of Construction \ref{construction:crystalline-c1} admits an essentially unique factorization as a composition
$$\RGamma_{\textnormal{\'{e}t}}( X, \mathbf{G}_m)[-1] \rightarrow \RGamma_{\textnormal{\'{e}t}}( X, \mathbf{G}_m)^{\wedge}[-1] \xrightarrow{ \widehat{c}_1^{\crys}} \Fil^{1}_{\Nyg} \RGamma_{\crys}(X/\Z_p).$$
\end{variant}

\begin{example}[The Semiperfect Case]\label{example:c1crys-semiperfect}
Let $R$ be a semiperfect $\F_p$-algebra, and let us identify $\Fil^{1}_{\Nyg} \RGamma_{\crys}(R/\Z_p)$ with the kernel of the augmentation map
$\epsilon^{\crys}_{R}: A_{\crys}(R) \rightarrow R$. Passing to cohomology in degree zero, the map
$$\widehat{c}_1^{\crys}: \RGamma_{\textnormal{\'{e}t}}( \Spec(R), \mathbf{G}_m)^{\wedge}[-1] \rightarrow \Fil^{1}_{\Nyg} A_{\crys}(R)$$
of Variant \ref{variant:crystalline-factor-through-completion} induces a group homomorphism $\theta: T_p(R^{\times}) \rightarrow \ker(\epsilon^{\crys}_{R} )$.
Writing $J$ for the kernel of the tautological map $R^{\flat} \twoheadrightarrow R$, we see that $\theta$ is given concretely by the composition
$$ T_p( R^{\times} ) \simeq (1+J)^{\times} \xrightarrow{ u \mapsto \log[u] } \Fil^{1}_{\Nyg} A_{\crys}(R)$$
here we write $[x]$ for the image of the Teichm\"{u}ller representative $[u]$ under the tautological map $W(R^{\flat} ) \rightarrow A_{\crys}(R)$
(so that $[u]-1$ belongs to the divided power ideal $\Fil^{1}_{\Nyg} A_{\crys}(R) = \ker( \epsilon_{\crys} )$).
\end{example}

Let $X$ be an $\F_p$-scheme. By functoriality, the absolute Frobenius map $\varphi_{X}: X \rightarrow X$
determines a commutative diagram
$$ \xymatrix@R=50pt@C=50pt{ \RGamma_{\textnormal{\'{e}t}}( X, \mathbf{G}_m)[-1]  \ar[r]^-{ c_{1}^{\crys} } \ar[d]^{ \varphi_{X}^{\ast} } & \Fil^{1}_{\Nyg} \RGamma_{\crys}(X / \Z_p) \ar[d]^{ \varphi_{X}^{\ast}}  \\
\RGamma_{\textnormal{\'{e}t}}( X, \mathbf{G}_m)[-1]  \ar[r]^-{ c_{1}^{\crys} } & \Fil^{1}_{\Nyg} \RGamma_{\crys}(X / \Z_p), }$$
where the left vertical map is given by multiplication by $p$ (Remark \ref{remark:Frobenius-G-m-cohomology}). If $X$ is $p$-quasisyntomic, this
observation admits a converse:

\begin{theorem}\label{theorem:syntomic-c1-crystalline}
Let $X$ be a quasisyntomic $\F_p$-scheme. Then the morphism $\widehat{c}_{1}^{\crys}$ of Variant \ref{variant:crystalline-factor-through-completion}
fits into a fiber sequence
$$ \RGamma_{\textnormal{\'{e}t}}( X, \mathbf{G}_m)^{\wedge}[-1]
\xrightarrow{ \widehat{c}_{1}^{\crys} } \Fil^{1}_{\Nyg} \RGamma_{\crys}(X/\Z_p) \xrightarrow{ \varphi/p - 1} 
\RGamma_{\crys}(X/\Z_p),$$
where $\varphi/p:  \Fil^{1}_{\Nyg} \RGamma_{\crys}(X/\Z_p) \rightarrow \RGamma_{\crys}(X/\Z_p)$ is the morphism described
in Remark \ref{remark:phi-over-p}. This fiber sequence is characterized (up to homotopy) by the requirement that it depends functorially on $X$.
\end{theorem}

\begin{proof}
For every $\F_p$-scheme $X$, let $M(X)$ denote the fiber of the morphism
$$ \varphi/p - 1: \Fil^{1}_{\Nyg} \RGamma_{\crys}(X/\Z_p) \rightarrow 
\RGamma_{\crys}(X/\Z_p).$$
It follows from Corollary \ref{corollary:qs-descent-for-crys} that, when restricted to quasisyntomic $\F_p$-schemes,
the construction $X \mapsto M(X)$ satisfies descent for the $p$-quasisyntomic topology. By virtue of Remark~\ref{remark:profppf-Gm}, the functor  $X \mapsto  \RGamma_{\textnormal{\'{e}t}}( X, \mathbf{G}_m)^{\wedge}[-1]$ is a sheaf for the pro-fppf topology. Consequently, to prove Theorem \ref{theorem:syntomic-c1-crystalline}, we are free to replace $X$ by a suitable pro-syntomic cover (which is automatically both a $p$-quasisyntomic cover and a pro-fppf cover) to reduce to the case where $X = \Spec(R)$ is the spectrum of a quasiregular semiperfect $\F_p$-algebra $R$.
In this case, we can use Proposition \ref{proposition:degree-zero-qrsp} to identify $\RGamma_{\textnormal{\'{e}t}}( X, \mathbf{G}_m)^{\wedge}[-1]$
with the Tate module $T_p(R^{\times})$, so that the desired fiber sequence is supplied by Theorem \ref{theorem:logarithm-sequence}.
\end{proof}

\subsection{Syntomic Cohomology of Formal Schemes}\label{subsection:syntomic-complex}

In this section, we review the syntomic complexes introduced in \cite{BMS2} and \cite{prisms}. Fix an integer $n$.
For every animated commutative ring $R$, we write $\iota$ for the canonical map $\Fil^{n}_{\Nyg} \Prism_{R}\{n\} \rightarrow
\Prism_{R}\{n\}$, and $\varphi\{n\}$ for the Frobenius map $\Fil^{n}_{\Nyg} \Prism_{R}\{n\} \rightarrow
\Prism_{R}\{n\}$ given by Notation \ref{notation:Frobenius-on-absolute}

\begin{construction}[Syntomic Complexes]\label{construction:syntomic-complex}
Let $R$ be an animated commutative ring. For every integer $n$, we let $\RGamma_{\Syn}( \Spf(R), \Z_p(n))$ denote the
fiber of the morphism $(\varphi\{n\} - \iota): \Fil^{n}_{\Nyg} \Prism_{R}\{n\} \rightarrow
\Prism_{R}\{n\}$, formed in the derived $\infty$-category $\widehat{\calD}(\Z_p)$. We will refer to $\RGamma_{\Syn}( \Spf(R), \Z_p(n))$
as the {\it $n$th syntomic complex} of $\Spf(R)$. We will denote the cohomology groups of these complexes
by $\mathrm{H}^{\ast}_{\Syn}( \Spf(R), \Z_p(n) )$ and refer to them as the {\it syntomic cohomology groups of $\Spf(R)$}.
In the special case $n=0$, we denote the complex $\RGamma_{\Syn}( \Spf(R), \Z_p(n) )$ by
$\RGamma_{\Syn}( \Spf(R), \Z_p)$ and its cohomology groups by $\mathrm{H}^{\ast}_{\Syn}( \Spf(R), \Z_p)$.
\end{construction}

\begin{warning}
In the situation of Construction \ref{construction:syntomic-complex}, we do not assume that the animated
commutative ring $R$ is $p$-complete. However, it would be harmless to add this assumption:
the tautological map from $R$ to its $p$-completion $\widehat{R}$ induces an isomorphism
$\RGamma_{\Syn}( \Spf(R), \Z_p(n)) \rightarrow \RGamma_{\Syn}( \Spf( \widehat{R} ), \Z_p(n))$.
In \S\ref{subsection:syntomic-complexes-general}, we will consider a variant of Construction \ref{construction:syntomic-complex}
which does not share this property (see Construction \ref{construction:syntomic-complexes-general}).
\end{warning}

\begin{example}\label{example:syntomic-of-qrsp}
Let $R$ be a quasiregular semiperfectoid ring. For every integer $n$, the complexes $\Fil^{n}_{\Nyg} \Prism_{R}\{n\}$
and $\Prism_{R}\{n\}$ are concentrated in cohomological degree zero (Corollary \ref{corollary:absolute-Nygaard-qrsp}). It follows that the cohomology groups $\mathrm{H}^{d}_{\Syn}( \Spf(R), \Z_p(n) )$ vanish for $d \notin \{0,1\}$, and we have an exact sequence of abelian groups
$$ 0 \rightarrow \mathrm{H}^{0}_{\Syn}( \Spf(R), \Z_p(n) ) \rightarrow \Fil^{n}_{\Nyg} \Prism_{R}\{n\} \xrightarrow{ \varphi\{n\} - \id}
\Prism_{R}\{n\} \rightarrow \mathrm{H}^{1}_{\Syn}( \Spf(R), \Z_p(n) ) \rightarrow 0.$$
In particular, we can identify $\mathrm{H}^{0}_{\Syn}( \Spf(R), \Z_p(n) ) $ with the abelian group
$\{ x \in \Prism_{R}\{n\}: \varphi\{n\}(x) = x \}$.
\end{example}

\begin{example}\label{example:syntomic-of-qrsp-char-p}
Let $R$ be a quasiregular semiperfect $\F_p$-algebra and let $n$ be an integer. Using Remark \ref{BKcrystalline} and
Lemma \ref{lemma:conjugate-filtration}, we can identify absolute prismatic complex $\Prism_{R}\{n\}$ with
the commutative ring $A_{\crys}(R)$ of Construction \ref{construction:Acrys}. Under this identification,
the Frobenius morphism $\varphi\{n\}: \Fil^{n}_{\Nyg} \Prism_{R}\{n\} \rightarrow \Prism_{R}\{n\}$ corresponds to the map
$$ \Fil^{n}_{\Nyg} A_{\crys}(R) \rightarrow A_{\crys}(R) \quad \quad x \mapsto \frac{ \varphi(x) }{p^n},$$
where $\Fil^{n}_{\Nyg} A_{\crys}(R) = \{ x \in A_{\crys}(R): \varphi(x) \in p^{n} A_{\crys}(R) \}$ denotes the ideal
described in Proposition \ref{proposition:concrete-Nygaard-qrsp}. It follows that 
$\RGamma_{\Syn}( \Spf(R), \Z_p(n))$ is represented concretely by the two-term chain complex of abelian groups
$$ \Fil^{n}_{\Nyg} A_{\crys}(R) \rightarrow A_{\crys}(R) \quad \quad x \mapsto \frac{ \varphi(x)}{p^n} - x.$$
\end{example}

\begin{example}\label{example:syntomic-1-qrsp}
Let $R$ be a quasiregular semiperfect $\F_p$-algebra. Then the exact sequence
$$ 0 \rightarrow T_p(R^{\times} ) \rightarrow  \Fil^{1}_{\Nyg} A_{\crys}(R) \xrightarrow{  \varphi/p - 1 } A_{\crys}(R) \rightarrow 0$$
of Theorem \ref{theorem:logarithm-sequence} supplies an isomorphism of the syntomic complex $\RGamma_{\Syn}( \Spf(R), \Z_p(1) )$ with the Tate module $T_p(R^{\times})$, regarded as a chain complex concentrated in cohomological degree zero.
\end{example}

Let $R$ be an animated commutative ring and let $n$ be an integer. By virtue of
Remark \ref{remark:Frobenius-on-Nygaard-completion}, the Frobenius map $\varphi\{n\}$
admits a canonical factorization
$$ \Fil^{n}_{\Nyg} \Prism_{R}\{n\} \rightarrow \Fil^{n}_{\Nyg} \widehat{\Prism}_{R}\{n\} \xrightarrow{\widetilde{\varphi}\{n\}} \Prism_{R}\{n\},$$
where $\widehat{\Prism}_{R}\{n\}$ denotes the Nygaard-completed absolute prismatic complex of Definition \ref{definition:Nygaard-completed-prismatic-complex}. Let us write $\widehat{\varphi}\{n\}: \Fil^{n}_{\Nyg} \widehat{\Prism}_{R}\{n\} \rightarrow \widehat{\Prism}_{R}\{n\}$
for the composition of $\widetilde{\varphi}\{n\}$ with the tautological map $\Prism_{R}\{n\} \rightarrow \widehat{\Prism}_{R}\{n\}$,
and let us write $\widetilde{\iota}: \Fil^{n}_{\Nyg} \widehat{\Prism}_{R}\{n\} \rightarrow \widehat{\Prism}_{R}\{n\}$
for the map supplied by the Nygaard filtraton on $\widehat{\Prism}_{R}$.

\begin{proposition}\label{proposition:syntomic-Nygaard-complete}
Let $R$ be an animated commutative ring. For every integer $n$, the tautological map
\begin{eqnarray*}
\RGamma_{\Syn}( \Spf(R), \Z_p(n) ) & = & \fib( \varphi\{n\} - \iota: \Fil^{n}_{\Nyg} \Prism_{R}\{n\} \rightarrow \Prism_{R}\{n\}) \\
& \simeq & \fib( \widehat{\varphi}\{n\} - \widehat{\iota}: \Fil^{n}_{\Nyg} \widehat{\Prism}_{R}\{n\} \rightarrow \widehat{\Prism}_{R}\{n\} )
\end{eqnarray*}
is an isomorphism in the $\infty$-category $\widehat{\calD}(\Z_p)$.
\end{proposition}

\begin{proof}
We have a commutative diagram
$$ \xymatrix@R=50pt@C=50pt{ \Fil^{n}_{\Nyg} \Prism_{R}\{n\} \ar[r]^-{ (\varphi, \iota) } \ar[d] & \Prism_{R}\{n\} \times \Prism_{R}\{n\} \ar[d] & \Prism_{R}\{n\} \ar@{=}[d] \ar[l]  \\
\Fil^{n}_{\Nyg} \widehat{\Prism}_{R}\{n\} \ar[r]^-{ (\widetilde{\varphi}\{n\}, \widehat{\iota}) } \ar@{=}[d] & \Prism_{R}\{n\} \times \widehat{\Prism}_{R}\{n\} \ar[d] & \Prism_{R}\{n\} \ar[d] \ar[l] \\
\Fil^{n}_{\Nyg} \widehat{\Prism}_R\{n\} \ar[r]^-{ ( \widehat{\varphi}\{n\}, \widehat{\iota}) } & \widehat{\Prism}_{R}\{n\} \times \widehat{\Prism}_{R}\{n\} & \widehat{\Prism}_{R}\{n\}. \ar[l] }$$
To prove Proposition \ref{proposition:syntomic-Nygaard-complete}, it will suffice to show that the vertical maps in this diagram
induce isomorphisms between the limits of the rows. For the upper rectangle, this follows since the upper left square is a pullback and
the vertical map on the upper right is an isomorphism. For the lower rectangle, it follows since the lower right square is a pullback and the 
vertical map on the lower left is an isomorphism.
\end{proof}

\begin{proposition}\label{proposition:syntomic-sheaf-descent}
For each integer $n$, the functor 
$$ \CAlg^{\anim} \rightarrow \widehat{\calD}(\Z_p) \quad \quad R \mapsto \RGamma_{\Syn}( \Spf(R), \Z_p(n) )$$
satisfies descent for the $p$-complete fpqc topology.
\end{proposition}

\begin{proof}
By virtue of Proposition \ref{proposition:syntomic-Nygaard-complete}, it will suffice to show that
the functors $R \mapsto \Fil^{n}_{\Nyg} \widehat{\Prism}_{\widehat{R}}\{n\}$ and $R \mapsto \widehat{\Prism}_{\widehat{R}}\{n\}$
satisfy descent for the $p$-complete fpqc topology, which follows from Remark \ref{remark:p-complete-fpqc-Nygaard}.
\end{proof}

\begin{proposition}\label{proposition:sifted-colimit-syntomic}
For each integer $n$, the functor
$$ \CAlg^{\anim} \rightarrow \widehat{\calD}(\Z_p) \quad \quad R \mapsto \RGamma_{\Syn}( \Spf(R), \Z_p(n) )$$
commutes with sifted colimits (and is therefore a left Kan extension of its restriction to the category $\Poly_{\Z}$ of finitely generated polynomial rings).
\end{proposition}

\begin{proof}
Set $n' = \max \{ 0, n \}$. Since $n' \geq 0$, we have a tautological map 
$\Prism_R\{n\} \rightarrow \Prism_{R}^{[-n']}\{n\}$, which we will denote by $\id$.
Since $n' \geq n$, the construction of Notation \ref{notation:Frobenius-on-absolute}
determines a Frobenius map $\varphi\{n\}: \Prism_{R}\{n\} \rightarrow \Prism_{R}^{[-n']}\{n\}$.
Let $F(R)$ denote the fiber of the difference
$$ (\varphi\{n\} - 1):  \Prism_{R}\{n\} \rightarrow \Prism_{R}^{[-n']}\{n\}.$$
We have a commutative diagram of fiber sequences
$$ \xymatrix@R=50pt@C=50pt{ \Fil^{n}_{\Nyg} \Prism_{R}\{n\} \ar[r] \ar[d]^{\varphi\{n\} - \iota} & \Prism_{R}\{n\} \ar[d]^{\varphi\{n\} - 1} \ar[r] & \Prism_{R}\{n\} / \Fil^{n}_{\Nyg} \Prism_{R}\{n\} \ar[d] \\
\Prism_{R}\{n\} \ar[r] & \Prism_{R}^{[-n']}\{n\} \ar[r] &  \Prism^{[-n']}_{R}\{n\} / \Prism_{R}\{n\} }$$
depending functorially on $R$. Moreover, the functors
$$R \mapsto \Prism_{R}\{n\} /  \Fil^{n}_{\Nyg} \Prism_{R}\{n\} \quad \quad R \mapsto \Prism_{R}^{[-n']}\{n\} / \Prism_{R}\{n\}$$
commute with sifted colimits (Remark \ref{remark:graded-Nygaard-sifted-colimit} and Proposition \ref{proposition:sokon}). 
It will therefore suffice to show that the functor $R \mapsto F(R)$ commutes with sifted colimits. Since
the extension of scalars functor
$$ \widehat{\calD}(\Z_p) \rightarrow \calD(\F_p) \quad \quad M \mapsto \F_p \otimes^{L} M$$
is conservative and preserves colimits, this is equivalent to the assertion that the functor $R \mapsto \F_p \otimes^{L} F(R)$
commutes with sifted colimits. For each integer $m > n$, the restriction of the Frobenius map
$$ \varphi\{n\}: \F_p \otimes^{L} \Prism_{R}\{n\} \rightarrow \F_p \otimes^{L} \Prism_{R}^{[-n']}\{n\},$$
to the complex $\F_p \otimes^{L} \Prism_{R}^{[m]}\{n\}$ factors naturally through $\F_p \otimes^{L} \Prism_{R}^{[pm-n]}\{n\}$,
and therefore also through $\F_p \otimes^{L} \Prism_{R}^{[m+1]}\{n\}$. We therefore obtain a map of filtered objects
$$ (\varphi\{n\} - 1): \F_p \otimes^{L} \Prism_{R}^{[\bullet]}\{n\} \rightarrow  \F_p \otimes^{L} \Prism_{R}^{[\bullet]}\{n\}$$
for $\bullet > n'$ which is multiplication by $-1$ at the associated graded level, and is therefore an isomorphism (since the filtration on both sides is complete: see Proposition \ref{proposition:sokon}). Choosing $m > \max \{n, -n' \}$, the commutative diagram of fiber sequences
$$ \xymatrix@R=50pt@C=50pt{ \F_p \otimes^{L} \Prism_{R}^{[m]}\{n\} \ar[d]^{\varphi\{n\}-1}_{\sim} \ar[r] & \F_p \otimes^{L} \Prism_{R}\{n\} \ar[d]^{\varphi\{n\}-1} \ar[r] & 
\F_p \otimes^{L} \Prism_{R}\{n\} / \Prism_{R}^{[m]}\{n\} \ar[d]^{\rho_{R}} \\
 \F_p \otimes^{L} \Prism_{R}^{[m]}\{n\}  \ar[r] &  \F_p \otimes^{L} \Prism_{R}^{[-n']}\{n\} \ar[r] &  \F_p \otimes^{L} (\Prism_{R}^{[-n']}\{n\}  / \Prism_{R}^{[m]}\{n\} ) }$$
determines a functorial identification $\F_p \otimes^{L} F(R) \simeq \fib( \rho_R )$. We conclude by observing that the source and target of $\rho_R$
commute with sifted colimits by virtue of Proposition \ref{proposition:sokon}.
\end{proof}

\begin{corollary}\label{corollary:vanishing-negative}
Let $R$ be an animated commutative ring. For every positive integer $n$, the complex
$\RGamma_{\Syn}( \Spf(R), \Z_p(-n) )$ vanishes.
\end{corollary}

\begin{proof}
By virtue of Proposition \ref{proposition:sifted-colimit-syntomic}, it will suffice to prove this when $R$ is $p$-quasisyntomic. Using Proposition \ref{proposition:syntomic-sheaf-descent}, we can further reduce to the case where $R$ is a quasiregular semiperfectoid ring, so that we can identify $( \Prism_{R}, \Prism_{R}^{[1]} )$ with a transversal prism $(A,I)$. Since $R$ is defined over a perfectoid ring, $A\{-n\}$ is a free
$A$-module of rank $1$; let $e \in A\{-n\}$ be a generator. We can then write $\varphi\{-n\}(e) = \lambda e$ for some
unique element $\lambda \in I^{n}$. Multiplication by $e^{-1}$ identifies $A\{-n\}$ with $A$, and therefore identifies
$\RGamma_{\Syn}( \Spf(R), \Z_p(-n) )$ with the two-term complex
$$ A \rightarrow A \quad \quad x \mapsto x - \lambda \varphi(x).$$
The differential on this complex has an inverse which is given explicitly by the construction
$$ x \mapsto x + \lambda \varphi(x) + \lambda \varphi(\lambda) \varphi^2(x) + \lambda \varphi(\lambda) \varphi^2(\lambda) \varphi^3(x) + \cdots;$$
note that this infinite sum is convergent because each product $\lambda \varphi(\lambda) \cdots \varphi^{r}(\lambda)$
belongs to the ideal $I_r$ of Notation \ref{notation:I-r} (see Remark \ref{remark:inverse-system-compute-A}).
\end{proof}

\begin{proposition}[Derived Descent for Reduction Modulo $p$]\label{proposition:derived-descent-syntomic}
Let $R$ be an animated commutative ring and let $R^{\bullet}$ denote the derived tensor product of $R$ with the
cosimplicial animated commutative ring $\F_p^{\otimes \bullet+1}$ of Notation \ref{notation:F-p-bullet}. Then, for every integer $n \geq 0$,
the tautological map $\RGamma_{\Syn}( \Spf(R), \Z_p(n)) \rightarrow \Tot \RGamma_{\Syn}( \Spf(R^{\bullet}), \Z_p(n) )$ is an isomorphism in the $\infty$-category
$\widehat{\calD}(\Z_p)$.
\end{proposition}

\begin{proof}
This is an immediate consequence of the corresponding statement for the functors $R \mapsto \Fil^{m}_{\Nyg} \Prism_{R}\{n\}$
(see Proposition \ref{proposition:derived-descent-absolute-Nygaard}).
\end{proof}

\begin{corollary}\label{corollary:reduce-mod-n-syntomic}
Let $R$ be an animated commutative ring. Then the tautological map
$$ \RGamma_{\Syn}( \Spf(R), \Z_p(n)) \rightarrow \varprojlim_{m}
\RGamma_{\Syn}( \Spf(R \otimes^{L} \Z/p^m \Z), \Z_p(n))$$
is an isomorphism.
\end{corollary}

\begin{variant}[Globalization to Formal Schemes]\label{variant:syntomic-sheaf-globalized}
Let $\mathfrak{X}$ be a bounded $p$-adic formal scheme. For every integer $n$, we let
$\RGamma_{\Syn}( \mathfrak{X}, \Z_p(n) )$ denote the limit
$$\varprojlim_{ \eta: \Spec(A) \rightarrow \mathfrak{X} } \RGamma_{\Syn}( \Spf(A), \Z_p(n)),$$
formed in the $\infty$-category$ \widehat{\calD}(\Z_p)$. This construction is characterized by the following properties:
\begin{itemize}
\item In the case where $\mathfrak{X} = \Spf(R)$ is affine (so that $R$ is a $p$-complete commutative ring with bounded $p$-power torsion), we can identify $\RGamma_{\Syn}( \mathfrak{X}, \Z_p(n) )$ with the syntomic complex
$\RGamma_{\Syn}( \Spf(R), \Z_p(n) )$ introduced in Construction \ref{construction:syntomic-complex}
(this follows from Corollary \ref{corollary:reduce-mod-n-syntomic}).

\item The construction $\mathfrak{X} \mapsto \RGamma_{\Syn}( \mathfrak{X}, \Z_p(n) )$
satisfies descent for the \'{e}tale topology on the category of bounded $p$-adic formal schemes
(see Proposition \ref{proposition:syntomic-sheaf-descent} for a stronger statement).
\end{itemize}

We will denote the cohomology groups of the complex $\RGamma_{\Syn}( \mathfrak{X}, \Z_p(n) )$ by
$\mathrm{H}^{\ast}_{\Syn}( \mathfrak{X}, \Z_p(n) )$ and refer to them as the {\it syntomic cohomology groups} of the formal scheme $\mathfrak{X}$. 
Note that we have a fiber sequence 
$$ \RGamma_{\Syn}( \mathfrak{X}, \Z_p(n) ) \rightarrow \Fil^{n}_{\Nyg} \RGamma_{\Prism}(\mathfrak{X})\{n\} \xrightarrow{ \varphi\{n\} - \iota} \RGamma_{\Prism}(\mathfrak{X})\{n\}.$$
\end{variant}

\subsection{The Syntomic First Chern Class}\label{subsection:syntomic-c1}

Let $X$ be a quasisyntomic $\F_p$-scheme. Combining Theorem \ref{theorem:syntomic-c1-crystalline}
with Theorem \ref{theorem:crystalline-comparison}, we see that the crystalline first Chern class
of Construction \ref{construction:crystalline-c1} can be promoted to an isomorphism
$$ \RGamma_{ \textnormal{\'{e}t}}( X, \mathbf{G}_m)^{\wedge}[-1]
\simeq \RGamma_{\Syn}(X, \Z_p(1) )$$
in the $\infty$-category $\widehat{\calD}(\Z_p)$. In this section, we use the prismatic logarithm of \S\ref{section:twist-and-log}
to extend this construction to all bounded $p$-adic formal schemes.

\begin{notation}
Let $R$ be a quasiregular semiperfectoid ring. Using Example \ref{example:syntomic-of-qrsp}, we can identify
the syntomic cohomology group cohomology group $\mathrm{H}^{0}_{\Syn}( \Spec(R), \Z_p(1) )$ can be identified with the abelian group $\{ x \in \Fil^{1}_{\Nyg} \Prism_{R}\{1\}: \varphi(x) = x \}$. Regarding $( \Prism_{R}, \Prism^{[1]}_{R})$ as a prism,
it follows from Proposition \ref{proposition:log-Nygaard} that the prismatic logarithm of Construction \ref{construction:Tate-logarithm} determines a group homomorphism
$$ \log_{\Prism}: T_p( \overline{\Prism}_{R}^{\times} ) \rightarrow \mathrm{H}^{0}_{\Syn}( \Spec(R), \Z_p(1) ).$$
\end{notation}

\begin{proposition}\label{proposition:prismatic-chern-class-construction}
Let $R$ be an animated commutative ring.
There is a canonical map
$$ c_{1}^{\Syn}: \RGamma_{\textnormal{\'{e}t}}( \Spec(R), \mathbf{G}_m)[-1] \rightarrow \RGamma_{\Syn}( \Spf( R), \Z_p(1) )$$
in the $\infty$-category $\widehat{\calD}(\Z_p)$, characterized (up to homotopy) by the requirement
that it depends functorially on $R$ and satisfies the following condition:
\begin{itemize}
\item[$(\ast)$] If $R$ is a quasiregular semiperfectoid ring, then the induced map
$$ \RGamma_{\textnormal{\'{e}t}}( \Spec(R), \mathbf{G}_m)^{\wedge}[-1] \rightarrow \RGamma_{\Syn}( \Spf(R), \Z_p(1) )$$
is given on cohomology in degree zero by homomorphism
\begin{equation}\label{equation:composition-defining-c1} T_p( R^{\times} ) \rightarrow T_p( \overline{\Prism}_{R}^{\times} ) \xrightarrow{ \log_{\Prism} } \mathrm{H}^{0}_{\Syn}( \Spf(R), \Z_p(1) ). \end{equation}
\end{itemize}
\end{proposition}

\begin{proof}
Let $\calC \subseteq \CAlg^{\qrsp}$ denote the category of quasiregular semiperfectoid rings $R$ for which
every element $x \in R$ admits a $p$th root. Note that, if this condition is satisfied, then we can the Tate
module $T_p(R^{\times})$ with the $p$-completion of the complex $\mathbf{G}_m(R)[-1]$ (as an object of the derived $\infty$-category $\calD(\Z)$). Let $\alpha_{R}$ denote the composite map
$$ \mathbf{G}_m(R)[-1] \rightarrow T_p(R) \rightarrow T_p( \overline{\Prism}_{R} ) \xrightarrow{ \log_{\Prism} }
\mathrm{H}^{0}_{\Syn}( \Spec(R), \Z_p(1) ) \rightarrow \RGamma_{\Syn}( \Spf(R), \Z_p(1) ).$$
The construction $R \mapsto \alpha_{R}$ determines a morphism from $\mathbf{G}_{m}[-1]$ to
$\RGamma_{\Syn}( \Spf(\bullet), \Z_p(1) )$ in the
$\infty$-category $\Fun( \calC, \widehat{\calD}(\Z_p) )$. Since the category $\calC$ forms a basis for the $p$-quasisyntomic topology
on $\CAlg^{\QSyn}$, the functor
$$ \CAlg^{\QSyn} \rightarrow \widehat{\calD}(\Z_p) \quad \quad R \mapsto \RGamma_{\Syn}( \Spf(R), \Z_p(1) )$$
is a right Kan extension of its restriction to $\calC$ (see Proposition \ref{proposition:syntomic-sheaf-descent}).
It follows that the construction $R \mapsto \alpha_{R}$ admits an essentially unique extension to the
category $\CAlg^{\QSyn}$ of $p$-quasisyntomic commutative rings. Since $\CAlg^{\QSyn}$ contains the Laurent polynomial rings $\Z[ x_1^{\pm 1}, \cdots, x_n^{\pm 1} ]$, the functor $\mathbf{G}_m[-1]$ is a left Kan extension of its restriction
to $\CAlg^{\QSyn}$ (Remark \ref{remark:Kan-extension-multiplicative}). Consequently, the extension
$R \mapsto \alpha_{R}$ admits an essentially unique extension to the $\infty$-category of {\em all}
animated commutative rings. Because the functor
$R \mapsto \RGamma_{\Syn}( \Spf(R), \Z_p(1) )$ satisfies \'{e}tale descent (Proposition \ref{proposition:syntomic-sheaf-descent}) the natural transformation $\alpha$ admits an essentially unique factorization as a composition
$$ \mathbf{G}_m(R)[-1] \rightarrow \RGamma_{\textnormal{\'{e}t} }( \Spf(R), \mathbf{G}_m)[-1] \xrightarrow{ c_{1}^{\Syn} } 
\RGamma_{\Syn}( \Spf(R), \Z_p(1)).$$
To complete the proof, it will suffice to show that this construction satisfies condition $(\ast)$ for
every quasiregular semiperfectoid ring $R$. Note that we can choose a $p$-quasisyntomic cover
$R \rightarrow R'$, where every element of $R'$ admits a $p$th root. Since the induced map
$\mathrm{H}^{0}_{\Syn}( \Spf(R), \Z_p(1) ) \rightarrow \mathrm{H}^{0}_{\Syn}( \Spf(R'), \Z_p(1) )$
is injective, we can replace $R$ by $R'$ and thereby reduce to verifying $(\ast)$ for quasiregular semiperfectoid rings which belong to $\calC$, in which case the desired result is immediate from the construction.
\end{proof}

\begin{notation}[The Syntomic First Chern Class]\label{notation:syntomic-c1}
Let $R$ be an animated commutative ring. We will refer to the morphism
$$ c_{1}^{\Syn}:  \RGamma_{\textnormal{\'{e}t}}( \Spec(R), \mathbf{G}_m)[-1] \rightarrow  \RGamma_{\Syn}( \Spf(R), \Z_p(1) )$$
as the {\it syntomic first Chern class}. We will refer to $c_{1}^{\Syn}$ as the {\it syntomic first Chern class}. We write
$c_{1}^{\Prism}$ for the composition 
$$  \RGamma_{\textnormal{\'{e}t}}( \Spec(R), \mathbf{G}_m)[-1] 
\rightarrow \RGamma_{\Syn}( \Spf(R), \Z_p(1) ) \rightarrow
\Fil^{1}_{\Nyg} \Prism_{R}\{1\},$$
which we refer to as the {\it prismatic first Chern class}. We will abuse notation by identifying
$c_{1}^{\Prism}$ with the composite map
$$\RGamma_{\textnormal{\'{e}t}}( \Spec(R), \mathbf{G}_m)[-1] \xrightarrow{ c_{1}^{\Prism} }
\Fil^{1}_{\Nyg} \Prism_{R}\{1\} \rightarrow \Prism_{R}\{1\},$$
which we also refer to as the {\it prismatic first Chern class}.
\end{notation}

\begin{variant}[Globalization to Formal Schemes]\label{variant:globalized-syntomic-c1}
Let $\mathfrak{X}$ be a bounded $p$-adic formal scheme. Then there is a canonical map
$$ c_{1}^{\Syn}: \RGamma_{\mathet}( \mathfrak{X}, \mathbf{G}_m)[-1]
\rightarrow \RGamma_{\Syn}(\mathfrak{X}, \Z_p(1) ),$$
which is characterized by the requirement that it depends functorially on $\mathfrak{X}$
and, when $\mathfrak{X} = \Spf(R)$ is affine, the composite map
$$   \RGamma_{\textnormal{\'{e}t}}( \Spec(R), \mathbf{G}_m)[-1] 
\rightarrow   \RGamma_{\textnormal{\'{e}t}}( \Spf(R), \mathbf{G}_m)[-1] 
\xrightarrow{ c_{1}^{\Syn} }  \RGamma_{\Syn}( \Spf(R), \Z_p(1) )$$
agrees with the map of Notation \ref{notation:syntomic-c1}. This follows 
by combining Corollary \ref{snifflet} with the observation that the
functor $\mathfrak{X} \mapsto \RGamma_{\Syn}( \mathfrak{X}, \Z_p(1) )$
satisfies \'{e}tale descent.

Passing to cohomology in degree $2$, we obtain a group homomorphism
$$ \Pic( \mathfrak{X} ) \rightarrow \mathrm{H}^{2}_{\Syn}( \mathfrak{X}, \Z_p(1) ),$$
which we will also refer to as the {\it syntomic first Chern class}.
\end{variant}

The prismatic first Chern class can be considered as a refinement of the crystalline first Chern class
of Construction \ref{construction:crystalline-c1}. For every $\F_p$-scheme $X$, let view the comparison map $\gamma_{\Prism}^{\crys}$ of Remark \ref{remark:crystalline-comparison-general}
as a morphism from $\RGamma_{\Prism}(X)\{1\}$ to $\RGamma_{\crys}( X/ \Z_p )$; here we implicitly invoke the identification
$$  \RGamma_{\Prism}(X)\{1\} \simeq \RGamma_{\Prism}(X/\Z_p)\{1\} \simeq \RGamma_{\Prism}(X/\Z_p) \simeq \RGamma_{\Prism}(X)$$
given by the trivialization of the Breuil-Kisin twist of the crystalline prism $( \Z_p, (p) )$.
By virtue of Proposition \ref{proposition:augmentation-compatibility}, $\gamma_{\Prism}^{\crys}$ induces a morphism
from $\Fil^{1}_{\Nyg} \RGamma_{\Prism}(X)\{\}$ to $\Fil^{1}_{\Nyg} \RGamma_{\crys}( X/ \Z_p)$, which we will also denote by
$\gamma_{\Prism}^{\crys}$.

\begin{proposition}\label{proposition:crystalline-agreement}
Let $X$ be an $\F_p$-scheme. Then the crystalline first Chern class
$$c_{1}^{\crys}:
\RGamma_{\textnormal{\'{e}t}}( X, \mathbf{G}_m)[-1] 
\rightarrow \Fil^{1}_{\Nyg} \RGamma_{\crys}( X / \Z_p)$$
agrees with the composition
\begin{eqnarray*}
\RGamma_{\textnormal{\'{e}t}}( X, \mathbf{G}_m)[-1] & \xrightarrow{c_{1}^{\Prism}} &
\Fil^{1}_{\Nyg} \RGamma_{\Prism}(X)\{1\} \\
& \xrightarrow{ \gamma^{\crys}_{\Prism} } & \Fil^{1}_{\Nyg} \RGamma_{\crys}( X / \Z_p ),
\end{eqnarray*}
up to a homotopy which depends functorially on $X$.
\end{proposition}

\begin{proof}
Without loss of generality we may assume that $X = \Spec(R)$ is affine. Let $u^{\crys}$ denote the composition
$$ \mathbf{G}_m(R)[-1] \rightarrow \RGamma_{\textnormal{\'{e}t}}(\Spec(R), \mathbf{G}_m)[-1] \xrightarrow{ c_{1}^{\crys} } \Fil^{1}_{\Nyg} \RGamma_{\crys}(R/\Z_p),$$
and define $u^{\Prism}: \mathbf{G}_m(R)[-1] \rightarrow \Fil^{1}_{\Nyg} \Prism_{R}$ similarly. Since the functor $R \mapsto \Fil_{1}^{\Nyg} \RGamma_{\crys}(R/\Z_p)$ satisfies descent for the \'{e}tale topology, it will suffice to show that $u^{\crys}$ agrees with the composition $\gamma_{\Prism}^{\crys} \circ u^{\Prism}$ (up to a homotopy which depends functorially on $R$). By virtue of Remark \ref{remark:Kan-extension-multiplicative}, the functor
$$ \CAlg_{\F_p}^{\anim} \rightarrow \calD(\Z) \quad \quad R \mapsto \mathbf{G}_m(R)$$
is a left Kan extension of its restriction to the category $\CAlg_{\F_p}^{\QSyn}$ of quasisyntomic $\F_p$-algebras; we may therefore restrict our attention to the case where $R$ is $p$-quasisyntomic. Since the functor 
$$ \CAlg_{\F_p}^{\QSyn} \rightarrow \calD(\Z) \quad \quad R \mapsto \Fil^{1}_{\Nyg} \RGamma_{\crys}(R/\Z_p)$$
satisfies descent for the $p$-quasisyntomic topology (Corollary \ref{corollary:qs-descent-for-crys}), we may further reduce to the case where $R$ is quasiregular semiperfect. In this case, the $p$-completion of the complex $\mathbf{G}_m(R)[-1]$ can be identified with the Tate module $T_p( R^{\times} )$. Since $\Fil^{1}_{\Nyg} \RGamma_{\crys}(R/\Z_p)$ is $p$-complete complete, the map $u^{\crys}$ determines a homomorphism of abelian groups
$\widehat{u}^{\crys}: T_p(R^{\times} ) \rightarrow \Fil^{1}_{\Nyg} \RGamma_{\crys}(R/\Z_p) \simeq \Fil^{1}_{\Nyg} A_{\crys}(R)$,
and $u^{\Prism}$ a homomorphism $\widehat{u}^{\Prism}: T_p( R^{\times} ) \rightarrow \Fil^{1}_{\Nyg} \Prism_{R}$.
We will complete the proof by showing that $\widehat{u}^{\crys}$ coincides with $\gamma^{\crys}_{\Prism} \circ \widehat{u}^{\Prism}$ (as a homomorphism of abelian groups).

Since $R$ is quasiregular semiperfect, the comparison map $\gamma_{\crys}^{\Prism}$ is an isomorphism, whose inverse is the map $\beta_{R}: A_{\crys}(R) \rightarrow \Prism_{R}$ which appears in the statement of Lemma \ref{lemma:transformation-in-easy-case} (after restriction
to the first stage of the Nygaard filtration). By construction, we have a commutative diagram
\begin{eqnarray}\label{equation:crystalline-agreement} \xymatrix@R=50pt@C=50pt{ A_{\crys}(R) \ar[d]^{ \epsilon_{\crys} } \ar[r]^-{\beta_{R}}_{\sim} & \Prism_R \ar[r]^-{\varphi} & \Prism_{R} \ar[d] \\
R \ar@{=}[rr] & & \overline{\Prism}_{R}.  }\end{eqnarray}
Let $J$ denote the kernel of the quotient map $R^{\flat} \twoheadrightarrow R$, let $x$ be an element of the abelian group $T_p( R^{\times} ) \simeq (1+J)^{\times}$, and let $[x]$ denote the image of the Teichm\"{u}ller representative of $x$ in $A_{\crys}(R)$. We then compute
\begin{eqnarray*}
\widehat{u}^{\Prism}(x) & = & \log_{\Prism}( \varphi( \beta_{R}([x]) ) ) \\
& = & \beta_{R}( \log[x] ) \\
& = & \beta_{R}( \widehat{u}^{\crys}(x) )
\end{eqnarray*}
where the first equality follows from the defining property of $c_{1}^{\Prism}$ together with the commutativity of (\ref{equation:crystalline-agreement}), the
second equality from Corollary \ref{corollary:crystalline-formula-for-log}, and the third equality from Example \ref{example:c1crys-semiperfect}. The desired result now follows by applying the isomorphism $\gamma_{\crys}^{\Prism}$ to both sides
\end{proof}

Let $R$ be an animated commutative ring. Since the syntomic complex $\RGamma_{\Syn}( \Spf(R), \Z_p(1))$
is $p$-complete, the morphism $c_{1}^{\Syn}$ of Proposition \ref{proposition:prismatic-chern-class-construction}
admits an essentially unique factorization as a composition 
$$ \RGamma_{\textnormal{\'{e}t}}( \Spec(R), \mathbf{G}_m)[-1] \rightarrow 
\RGamma_{\textnormal{\'{e}t}}( \Spec(R), \mathbf{G}_m)^{\wedge}[-1] 
\xrightarrow{ \widehat{c}_{1}^{\Syn} } \RGamma_{\Syn}( \Spf(R), \Z_p(1) ).$$
We can now formulate our main result:

\begin{theorem}\label{theorem:syntomic-chern-class-isomorphism}
Let $R$ be a $p$-complete animated commutative ring. Then the morphism
$$ \widehat{c}_{1}^{\Syn}: \RGamma_{\textnormal{\'{e}t}}( \Spec(R), \mathbf{G}_m)^{\wedge}[-1] \rightarrow \RGamma_{\Syn}( \Spf(R), \Z_p(1) )$$
is an isomorphism (in the derived $\infty$-category $\widehat{\calD}(\Z_p)$).
\end{theorem}

Theorem \ref{theorem:syntomic-chern-class-isomorphism} is essentially proven in \cite{BMS2}, at least for $p$-quasisyntomic commutative rings (see Proposition~7.17 of \cite{BMS2}). We present here a different proof, which
avoids the use of algebraic $K$-theory and instead proceeds by reduction to the analogous statement
for crystalline cohomology (Theorem \ref{theorem:syntomic-c1-crystalline}), which was proved using
the exact sequence of Theorem \ref{theorem:logarithm-sequence}. We proceed in several steps.

\begin{lemma}\label{lemma:main-theorem-crys}
Let $R$ be a quasisyntomic $\F_p$-algebra. Then the morphism
$$ \widehat{c}_{1}^{\Syn}: \RGamma_{\textnormal{\'{e}t}}( \Spec(R), \mathbf{G}_m)^{\wedge}[-1] \rightarrow \RGamma_{\Syn}( \Spf(R), \Z_p(1) )$$
is an isomorphism.
\end{lemma}

\begin{proof}
By virtue of Proposition \ref{proposition:crystalline-agreement}, this is equivalent to the assertion
that the map
$$ \widehat{c}_{1}^{\crys}: \RGamma_{\textnormal{\'{e}t}}( \Spec(R), \mathbf{G}_m)^{\wedge}[-1]
\rightarrow \RGamma_{\crys}( R / \Z_p )$$
is an isomorphism, which follows from Theorem \ref{theorem:syntomic-c1-crystalline}.
\end{proof}

\begin{lemma}\label{lemma:make-alpha-R}
Let $R$ be an animated $\F_p$-algebra. Then there is a canonical map
$$ \rho_{R}: \RGamma_{\Syn}( \Spf(R), \Z_p(1) ) \rightarrow \RGamma_{\textnormal{\'{e}t}}( \Spec(R), \mathbf{G}_m)^{\wedge}[-1],$$
which is characterized by the requirement that it depends functorially on $R$ (as an object of the $\infty$-category $\CAlg_{\F_p}^{\anim}$)
and the composition
$$ \RGamma_{\Syn}( \Spf(R), \Z_p(1) ) \xrightarrow{ \rho_{R} } \RGamma_{\textnormal{\'{e}t}}( \Spec(R), \mathbf{G}_m)^{\wedge}[-1] \xrightarrow{ \widehat{c}_{1}^{\Syn} }
\RGamma_{\Syn}( \Spf(R), \Z_p(1) )$$
is the identity (up to a homotopy which depends functorially on $R$). Here $\widehat{c}_{1}^{\Syn}$
denotes the $p$-completion of the morphism $c_{1}^{\Syn}$ 
\end{lemma}

\begin{proof}
By virtue of Proposition \ref{proposition:sifted-colimit-syntomic}, the functor
$$ \CAlg^{\anim}_{\F_p} \rightarrow \widehat{\calD}(\Z_p) \quad \quad R \mapsto \RGamma_{\Syn}( \Spf(R), \Z_p(1) )$$
commutes with sifted colimits, and is therefore a left Kan extension of its restriction to the category $\Poly_{\F_p}$ of finitely
generated polynomial algebras over $\F_p$. It will therefore suffice to construction the morphism $\rho_{R}$
(and the homotopy $\widehat{c}_{1}^{\Syn} \circ \rho_{R} \simeq \id$) in the case where $R$ is a polynomial algebra over $\F_p$.
In this case, the morphism $\widehat{c}_{1}^{\Syn}$ is an isomorphism (Lemma \ref{lemma:main-theorem-crys}).
\end{proof}

\begin{lemma}\label{lemma:main-theorem-Fp-case}
Let $R$ be an animated $\F_p$-algebra. Then the morphism
$$ \widehat{c}_{1}^{\Syn}: \RGamma_{\textnormal{\'{e}t}}( \Spec(R), \mathbf{G}_m)^{\wedge}[-1] \rightarrow \RGamma_{\Syn}( \Spf(R), \Z_p(1) )$$
is an isomorphism.
\end{lemma}

\begin{proof}
Let $\rho_{R}: \RGamma_{\Syn}( \Spf(R), \Z_p(1) ) \rightarrow \RGamma_{\textnormal{\'{e}t}}( \Spec(R), \mathbf{G}_m)^{\wedge}[-1]$ be the morphism of Lemma \ref{lemma:make-alpha-R}. By construction, $\rho_{R}$ is a right homotopy inverse to $\widehat{c}_{1}^{\Syn}$.
We claim that it is also left homotopy inverse: that is, that the composition
$$  \RGamma_{\textnormal{\'{e}t}}( \Spec(R), \mathbf{G}_m)^{\wedge}[-1] \xrightarrow{ \widehat{c}_{1}^{\Syn} }  \RGamma_{\Syn}( \Spf(R), \Z_p(1) )
\xrightarrow{\rho_{R}} \RGamma_{\textnormal{\'{e}t}}( \Spec(R), \mathbf{G}_m)^{\wedge}[-1]$$
is homotopic to the identity. Moreover, we will show that the homotopy can be chosen functorially in $R \in \CAlg^{\anim}_{\F_p}$.

For every animated $\F_p$-algebra $R$, let $u_{R}: \mathbf{G}_{m}(R)[-1] \rightarrow \RGamma_{\textnormal{\'{e}t}}( \Spec(R), \mathbf{G}_m)^{\wedge}[-1]$
be the tautological map. Since the functor $R \mapsto \RGamma_{\textnormal{\'{e}t}}( \Spec(R), \mathbf{G}_m)^{\wedge}[-1]$ takes $p$-complete values and satisfies \'{e}tale descent, it will suffice to show $u_{R}$ agrees with $u_{R} \circ \rho_{R} \circ \widehat{c}_{1}^{\Syn}$ (up to a homotopy which can be chosen functorially in $R$).
By virtue of Remark \ref{remark:Kan-extension-multiplicative}, the functor $R \mapsto \mathbf{G}_{m}(R)[-1]$ is a left Kan extension of its restriction to smooth $\F_p$-algebras. It will therefore suffice to show that $\rho_{R}$ is a left homotopy inverse
to $\widehat{c}_{1}^{\Syn}$ when restricted to smooth $\F_p$-algebras, which follows from Lemma \ref{lemma:main-theorem-crys}.

\end{proof}

\begin{proof}[Proof of Theorem \ref{theorem:syntomic-chern-class-isomorphism}]
Let $R$ be a $p$-complete animated commutative ring, let $\F_p^{\otimes \bullet +1 }$ be the cosimplicial animated commutative ring
introduced in Notation \ref{notation:F-p-bullet}, and let $R^{\bullet}$ denote the cosimplicial $R$-algebra given by 
$R^{\bullet} = R \otimes^{L}_{\Z} \F_p^{\otimes \bullet +1 }$. We then have a commutative diagram
$$ \xymatrix@R=50pt@C=50pt{ \RGamma_{\textnormal{\'{e}t}}( \Spec(R), \mathbf{G}_m)^{\wedge}[-1] \ar[r]^-{ \widehat{c}^{\Syn}_{1}} \ar[d] &  \RGamma_{\Syn}( \Spf(R), \Z_p(1) ) \ar[d] \\
\Tot( \RGamma_{\textnormal{\'{e}t}}( \Spec(R^{\bullet}), \mathbf{G}_m)^{\wedge}[-1] ) \ar[r] & \Tot(\RGamma_{\Syn}( \Spf(R), \Z_p(1)) ), }$$
where the vertical maps are isomorphisms by virtue of Propositions \ref{proposition:derived-descent-syntomic} and \ref{proposition:derived-descent-Gm}.
Consequently, to show that the map $\widehat{c}_{1}^{\Syn}: \RGamma_{\textnormal{\'{e}t}}( \Spec(R), \mathbf{G}_m)^{\wedge}[-1] \rightarrow
\RGamma_{\Syn}( \Spf(R), \Z_p(1))$ is an isomorphism, it suffices to prove the analogous statement for the animated
commutative rings $R^{n}$ for each $n \geq 0$. This follows from Lemma \ref{lemma:main-theorem-Fp-case}, since $R^{n}$ admits the structure of an animated $\F_p$-algebra.
\end{proof}

\subsection{The de Rham Specialization}\label{subsection:de-Rham-c1}

Our goal in this section is to compare the syntomic first Chern class of Proposition \ref{proposition:prismatic-chern-class-construction} with the classical construction of Chern classes in (derived) de Rham cohomology.
We begin by reviewing how the latter are defined.

\begin{proposition}\label{proposition:construction-of-c1dr}
For every animated commutative ring $R$, there is a canonical map
$$ c_{1}^{\dR}: \RGamma( \Spec(R)_{\textnormal{\'{e}t}}, \mathbf{G}_m)[-1] \rightarrow \Fil^{1}_{\Hodge} \widehat{\dR}_{R}$$
in the $\infty$-category $\widehat{\calD}(\Z)$, which is determined (up to essentially unique homotopy) by the requirement that it
depends functorially on $R$ and that, when $R$ is an ordinary commutative ring, the composite map
$$ \mathbf{G}_{m}(R)[-1] \rightarrow \RGamma_{\textnormal{\'{e}t}}( \Spec(R), \mathbf{G}_m)[-1]
\xrightarrow{ c_{1}^{\dR} }  \Fil^{1}_{\Hodge} \widehat{\dR}_{R}
\rightarrow (\widehat{\Omega}^{\geq 1}_{R},d) \rightarrow \widehat{\Omega}^{1}_{R}[-1]$$
coincides with the homomorphism of abelian groups
$$ R^{\times} \rightarrow \widehat{\Omega}^{1}_{R} \quad \quad u \mapsto \dlog(u) = \frac{du}{u}.$$
\end{proposition}

\begin{proof}
Since the functor $R \mapsto \Fil^{1}_{\Hodge} \widehat{\dR}_{R}$ satisfies \'{e}tale descent and
the functor $R \mapsto \mathbf{G}_m(R)$ is a left Kan extension of its restriction to smooth $\Z$-algebras (Remark \ref{remark:Kan-extension-multiplicative}),
we may restrict our attention to the case where $R$ is a smooth $\Z$-algebra. In this case, the comparison map
$\Fil^{1}_{\Hodge} \widehat{\dR}_{R} \rightarrow (\widehat{\Omega}^{\geq 1}_{R},d)$ is an isomorphism in $\calD(\Z)$ (see Proposition \ref{proposition:derived-to-classical-de-Rham}). We are therefore reduced to showing that the map $\dlog: \mathbf{G}_m(R)[-1] \rightarrow \widehat{\Omega}^{1}_{R}[-1]$ admits an essentially unique
factorization through the truncated de Rham complex $(\widehat{\Omega}^{\geq 1}_{R},d)$. This follows from the assertion that
for every invertible element $u \in R$, the differential form $\dlog(u) = du/u$ is closed.
\end{proof}

We will refer to the morphism $c_{1}^{\dR}: \RGamma_{\textnormal{\'{e}t}}( \Spec(R), \mathbf{G}_m)[-1] \rightarrow \Fil^{1}_{\Hodge} \widehat{\dR}_{R}$
as the {\it de Rham first Chern class}. Our main result can now be stated as follows:

\begin{theorem}\label{theorem:dR-comparison-c1}
Let $R$ be an animated commutative ring. Then the de Rham first Chern class $c_{1}^{\dR}:  \RGamma( \Spec(R)_{\textnormal{\'{e}t}}, \mathbf{G}_m)[-1] \rightarrow \Fil^{1}_{\Hodge} \widehat{\dR}_{R}$
of Proposition \ref{proposition:construction-of-c1dr} agrees with the composition
$$ \RGamma_{\textnormal{\'{e}t}}( \Spec(R), \mathbf{G}_m)[-1]
\xrightarrow{ c_{1}^{\Prism} }  \Fil^{1}_{\Nyg} \Prism_{R}\{1\} \xrightarrow{ \gamma_{\Prism}^{\dR}\{1\} } \Fil^{1}_{\Hodge} \widehat{\dR}_{R}$$
(up to a homotopy which depends functorially on $R$); here $\gamma_{\Prism}^{\dR}\{1\}$ is the de
Rham comparison morphism (see Construction \ref{construction:absolute-Nygaard-untwisted}).
\end{theorem}

\begin{proof}
For every animated commutative ring $R$ with underlying commutative ring $\pi_0(R)$, let $u_{R}$ denote the composite map
\begin{eqnarray*}
\mathbf{G}_m(R)[-1] & \rightarrow &  \RGamma_{\textnormal{\'{e}t}}(\Spec(R), \mathbf{G}_m)[-1] \\
& \xrightarrow{ c_{1}^{\Prism} } &  \Fil^{1}_{\Nyg} \Prism_{R}\{1\} \\
& \xrightarrow{ \gamma_{\Prism}^{\dR}\{1\} } & \Fil^{1}_{\Hodge} \widehat{\dR}_{R} \\
& \rightarrow & \gr^{1}_{\Hodge} \widehat{\dR}_{R} \\
& \simeq & L\widehat{\Omega}^{1}_{R}[-1];
\end{eqnarray*}
here $L\widehat{\Omega}^{1}_{R}$ denotes the $p$-complete absolute cotangent complex of $R$.
Specializing to the case where $R$ is an ordinary commutative ring and passing to cohomology in
degree $1$, we obtain a homomorphism of abelian groups $\rho_{R}: R^{\times} \rightarrow \widehat{\Omega}^{1}_{R}$.
We will complete the proof by showing that the homomorphism $\rho_{R}$ is given by $u \mapsto \dlog(u)$.

Let $R_0$ denote the Laurent polynomial ring $\Z[ t^{\pm 1} ]$ and let $\widehat{R}_0$ denote its $p$-completion, so that
$\widehat{\Omega}^{1}_{R_0}$ is freely generated as a $\widehat{R}_0$-module by the element $\dlog(t) = \frac{dt}{t}$.
We can therefore write $\rho_{R_0}(t) = f(t) \dlog(t)$ for some unique element $f(t) \in \widehat{R}_0$. It then
follows by functoriality that, for every invertible element $u$ of any commutative ring $R$, we have
$\rho_{R}(u) = f(u) \dlog(u)$ in $\widehat{\Omega}^{1}_{R}$. Specializing to the case $R = \Z[ u^{\pm 1}, v^{\pm 1} ]$, we obtain an identity
\begin{eqnarray*}
f(u) \dlog(u) + f(v) \dlog(v) & = & \rho_R(u) + \rho_R(v) \\
& = & \rho_{R}(uv) \\
& = & f(uv) \dlog(uv) \\
& = & f(uv) \dlog(u) + f(uv) \dlog(v).
\end{eqnarray*}
Since $\widehat{\Omega}^{1}_{R}$ is freely generated by $\dlog(u)$ and $\dlog(v)$ as a module over $\widehat{R}$,
it follows that we have an equality $f(u) = f(uv) = f(v)$ in $\widehat{R}$, so that $f(t) = c$ for some scalar $c \in \Z_p$.
We will complete the proof by showing that $c = 1$.

Let $(A,I)$ be the perfection of the $q$-de Rham prism (so that $A$ is the completion of $\Z[ q^{1/p^{\infty}} ]$ with respect to the ideal $(p,q-1)$
and $I$ is the principal ideal generated by the element $[p]_{q} = \frac{ q^{p} - 1}{q-1}$) let $\Z_p^{\cyc}$ denote the quotient ring
$A/I$. To avoid confusion, let us write $R$ for same commutative ring $\Z_p^{\cyc}$, regarded as an $A$-algebra
via the prismatic augmentation $\epsilon_{\Prism}: A \rightarrow \Z_p^{\cyc}$ (that is, for the quotient ring $A/J$, where $J = \varphi^{-1}(I)$ is the ideal generated by the element $[p]_{q^{1/p}}$). Note that the $p$-complete complex $L\widehat{\Omega}_{R}[-1]$ can be identified with the quotient
$J/J^2$, so that $\rho_{R}$ induces a morphism of Tate modules $\widehat{\rho}_{R}: T_p( R^{\times} ) \rightarrow J/J^2$.
Let $\epsilon \in T_p( R^{\times} )$ denote the element given by
the compatible system of $p$-power roots of unity $(q, q^{1/p}, q^{1/p^2}, \cdots )$.
Using the commutativity of the diagram
$$ \xymatrix@R=50pt@C=50pt{ 0 \ar[r] & \Z_p \ar[r] \ar[d] & \Q_p \ar[d]^{\gamma \mapsto \dlog(q^{\gamma}) } \ar[r] & \Q_p / \Z_p \ar[d] \ar[r] & 0 \\
 & J/J^2 \ar[r]^-{d} & R \otimes_{A} \Omega^{1}_{A} \ar[r]  & \Omega^{1}_{R} \ar[r] & 0 }$$
and the congruence $\dlog(q) \equiv d(q-1) \pmod{q-1}$, we obtain $\widehat{\rho}_{R}(\epsilon) = c(q-1)$
as an element of $J/J^2$. On the other hand, the prismatic logarithm $\log_{\Prism}$ carries $\epsilon$ to the element
$\log_{\Prism}(q^{p}) \in A\{1\}$. Using the definition of the comparison map $\gamma_{\Prism}^{\dR}\{1\}$ and
Corollary \ref{corollary:crystalline-formula-for-log}, we see that $\widehat{\rho}_{R}(\epsilon)$ is the image of $\log(q)$ under the quotient map $A_{\crys}(R/pR) \twoheadrightarrow A/J^2$. It follows that the scalar $c$ can be identified with the coefficient of
$q-1$ in the power series expansion 
$$ \log(q) = \log( 1 + (q-1) ) = \sum_{m > 0} (-1)^{m+1} \frac{ (q-1)^{m} }{m},$$
which is equal to $1$.
\end{proof}

\newpage \section{Comparison with \'{E}tale Cohomology}\label{section:etale-comparison}

Let $X$ be a scheme, and assume for simplicity that the structure sheaf $\calO_{X}$ has bounded $p$-power torsion. We can then consider two different $p$-adic cohomological invariants of $X$:
\begin{itemize}
\item[$(1)$] The $p$-adic \'{e}tale cohomology groups $\mathrm{H}^{\ast}_{\mathet}( U, \Z_p(n) )$, where
$U = \Spec( \Z[1/p]) \times X$ is the open subscheme of $X$ where $p$ is invertible.

\item[$(2)$] The syntomic cohomology groups $\mathrm{H}^{\ast}_{\Syn}( \mathfrak{X}, \Z_p(n) )$, where
$\mathfrak{X} = \Spf(\Z_p) \times X$ is the formal completion of $X$ along the vanishing locus of $p$.
\end{itemize}
Our goal in this section is to study the relationship between these invariants. Let us begin by considering the case
$n=1$. In this case, both $\mathrm{H}^{\ast}_{\Syn}( \mathfrak{X}, \Z_p(1) )$ and $\mathrm{H}^{\ast}_{\mathet}(U, \Z_p(1))$ are closely related to the cohomology groups $\mathrm{H}^{\ast}_{\mathet}( X, \mathbf{G}_m )$.
More precisely, there are canonical maps
$$ \xymatrix@R=50pt@C=50pt{ & \mathrm{H}^{\ast}_{\mathet}( X, \mathbf{G}_m) \ar[dl] \ar[dr] & \\
\mathrm{H}^{\ast}_{\mathet}( \mathfrak{X}, \mathbf{G}_m) \ar[d]^{ c_{1}^{\Syn} }  & & \mathrm{H}^{\ast}_{\mathet}(U, \mathbf{G}_m) \ar[d]^{ c_{1}^{\mathet}} \\
 \mathrm{H}^{\ast+1}_{\Syn}( \mathfrak{X}, \Z_p(1) )  & & \mathrm{H}^{\ast+1}_{\mathet}(U, \Z_p(1) ) }$$
where the left vertical map is (induced by) the syntomic first Chern class of Variant \ref{variant:globalized-syntomic-c1},
and the right vertical map is its counterpart in $p$-adic \'{e}tale cohomology (whose definition we recall in
\S\ref{subsection:etale-chern-class}). These can be obtained from maps of cochain complexes
\begin{equation}
\begin{gathered}\label{diagram:two-chern-classes}
 \xymatrix@R=50pt@C=50pt{ & \RGamma_{\mathet}(X, \mathbf{G}_m)[-1] \ar[dl]_{c_{1}^{\Syn}} \ar[dr]^{ c_{1}^{\mathet} } & \\
\RGamma_{\Syn}( \mathfrak{X}, \Z_p(1) ) & & \RGamma_{\mathet}(U, \Z_p(1) ). }
\end{gathered}
\end{equation}

Suppose now that the scheme $X = \Spec(R)$ is affine, and that the commutative ring $R$ is $p$-complete.
In this case, Theorem \ref{theorem:syntomic-chern-class-isomorphism} guarantees that the map
$c_{1}^{\Syn}$ exhibits $\RGamma_{\Syn}( \mathfrak{X}, \Z_p(1) )$ as the $p$-completion of 
$ \RGamma_{\mathet}(X, \mathbf{G}_m)[-1]$ (as an object of the derived $\infty$-category $\calD(\Z)$). 
Since the complex $\RGamma_{\mathet}(U, \Z_p(1) )$ is $p$-complete, we can expand
(\ref{diagram:two-chern-classes}) to a commutative diagram
$$ \xymatrix@R=50pt@C=50pt{ & \RGamma_{\mathet}( \Spec(R), \mathbf{G}_m)[-1] \ar[dl]_{c_{1}^{\Syn}} \ar[dr]^{ c_{1}^{\mathet} } & \\
\RGamma_{\Syn}( \Spf(R), \Z_p(1) ) \ar[rr]^{ \gamma_{\Syn}^{\mathet}\{1\} } & & \RGamma_{\mathet}( \Spec(R[1/p]), \Z_p(1) ). }$$
In \S\ref{subsection:etale-comparison-morphism}, we show that $\gamma_{\Syn}^{\mathet}\{1\}$ admits an essentially
unique extension to a multiplicative map
$$ \bigoplus_{n \in \Z} \RGamma_{\Syn}( \Spf(R), \Z_p(n)) \rightarrow \bigoplus_{n \in \Z} \RGamma_{\mathet}( \Spec( R[1/p] ), \Z_p(n) )$$
which depends functorially on $R$ (Theorem \ref{subsection:etale-comparison-morphism}). Restricting to graded degree
$n$, we obtain a map $\gamma_{\Syn}^{\mathet}\{n\}: \RGamma_{\Syn}( \Spf(R), \Z_p(n)) \rightarrow
\RGamma_{\mathet}( \Spec( R[1/p] ), \Z_p(n) )$ which we will refer to as the {\it \'{e}tale comparison morphism} for syntomic cohomology. 

In \S\ref{subsection:syntomic-complexes-general}, we use the \'{e}tale comparison morphism of
\S\ref{subsection:etale-comparison-morphism} to amalgamate $p$-adic \'{e}tale cohomology (regarded as an invariant of $\Z[1/p]$-schemes) and syntomic cohomology (regarded as an invariant of $p$-adic formal schemes) into a single invariant. Let $R$ be an animated commutative ring and let
$\widehat{R}$ denote its $p$-completion. We let $\RGamma_{\Syn}( \Spec(R), \Z_p(n) )$ denote the pullback 
of the diagram $$ \RGamma_{\Syn}( \Spf(\widehat{R}), \Z_p(n) )
\xrightarrow{ \gamma_{\Syn}^{\mathet}\{n\} } \RGamma_{\mathet}( \Spec( \widehat{R}[1/p]), \Z_p(n) )
\leftarrow \RGamma_{\mathet}( \Spec(R[1/p]), \Z_p(n)).$$
We denote the cohomology groups of $\RGamma_{\Syn}( \Spec(R), \Z_p(n) )$ by $\mathrm{H}^{\ast}_{\Syn}( \Spec(R), \Z_p(n) )$, which we will refer to as the {\it syntomic cohomology groups} of $\Spec(R)$. This construction has the following features:

\begin{itemize}
\item The construction $R \mapsto \RGamma_{\Syn}( \Spec(R), \Z_p(n) )$ satisfies descent for the fpqc topology (Proposition \ref{rhoax}). In particular, it admits a canonical extension $X \mapsto \RGamma_{\Syn}( X, \Z_p(n) )$ to
the category of schemes (Variant \ref{variant:syntomic-sheaf-integrally-globalized}).

\item If $X$ is a scheme with the property that $p$ is invertible in $\calO_{X}$, then
$\RGamma_{\Syn}(X, \Z_p(n) )$ can be identified with the usual \'{e}tale cochain complex $\RGamma_{\mathet}( X, \Z_p(n) )$
(Remark \ref{remark:arithmetic-square-for-syntomic}).

\item Let $X$ be a scheme which is proper over a commutative ring $R$,
and assume that the structure sheaf $\calO_{X}$ has bounded $p$-power torsion. If
$R$ is $p$-complete, then $\RGamma_{\Syn}( X, \Z_p(n) )$ can
be identified with the syntomic complex of the formal scheme $\mathfrak{X} = \Spf(\Z_p) \times X$.

\item For every scheme $X$, there is a canonical isomorphism $\RGamma_{\Syn}( X, \Z_p(0) )
\simeq \RGamma_{\mathet}(X, \Z_p)$ (Proposition \ref{proposition:keron}). Our proof will use
an analysis of the untwisted syntomic complex $\RGamma_{\Syn}( \Spf(R), \Z_p )$,
which we carry out in \S\ref{subsection:etale-constant} (see Theorem \ref{theorem:syntomic-complex-degree-zero}).

\item The syntomic first Chern class of Variant \ref{variant:globalized-syntomic-c1} has a schematic counterpart.
To every scheme $X$, we can associate a canonical map
$$ c_{1}^{\Syn}: \RGamma_{\mathet}( X, \mathbf{G}_m)[-1] \rightarrow \RGamma_{\Syn}(X, \Z_p(1) ),$$
which exhibits $\RGamma_{\Syn}(X, \Z_p(1) )$ as a $p$-completion of $\RGamma_{\mathet}( X, \mathbf{G}_m)[-1]$
(Proposition \ref{vrox}). Passing to cohomology in degree $2$, we obtain a homomorphism of abelian groups
$$ \Pic(X) \rightarrow \mathrm{H}^{2}_{\Syn}( X, \Z_p(1) ) \quad \quad \mathscr{L} \mapsto c_{1}^{\Syn}( \mathscr{L} ).$$

\item Let $\Z[ \zeta_{p^{\infty}} ]$ denote the cyclotomic number ring $\bigcup_{n \geq 1} \Z[ \zeta_{p^{n}} ]$, so
that we can regard the sequence $\epsilon = (1, \zeta_p, \zeta_p^2, \cdots )$ as an element
of the Tate module $$T_p( \Z[ \zeta_{p^{\infty}} ]^{\times} ) \simeq \mathrm{H}^{0}_{\Syn}( \Spec( \Z[ \zeta_{p^{\infty}} ]), \Z_p(1) ).$$
Let $X$ be a $\Z[ \zeta_{p^{\infty}} ]$-scheme which is quasi-compact and quasi-separated, and let 
$$U = \Spec( \Q[ \zeta_{p^{\infty}} ]) \times_{ \Spec( \Z[ \zeta_{p^{\infty} }] )} X$$ denote its generic fiber.
In \S\ref{subsection:localization-generic-fiber}, we show that the \'{e}tale cochain complex $\RGamma_{\mathet}(U, \Z_p)$
can be identified with the $p$-completed direct limit of the diagram
$$ \RGamma_{\Syn}(X, \Z_p) \xrightarrow{\epsilon} \RGamma_{\Syn}(X, \Z_p(1) ) \xrightarrow{ \epsilon} 
\RGamma_{\Syn}(X, \Z_p(2) ) \xrightarrow{\epsilon} \cdots$$
(Theorem \ref{theorem:localization-to-generic-fiber}). Stated more informally, passage from syntomic cohomology to
$p$-adic \'{e}tale cohomology can be implemented by inverting the element $\epsilon$.
\end{itemize}

\subsection{\'{E}tale Cohomology with Constant Coefficients}\label{subsection:etale-constant}

Let $R$ be an animated commutative ring. For every integer $k \geq 0$, we let $\RGamma_{\mathet}( \Spec(R), \Z / p^{k} \Z )$ denote the derived global sections
of the constant sheaf $\underline{ \Z / p^{k} \Z}$ on the \'{e}tale site of $\Spec(R)$. We let $\RGamma_{\mathet}( \Spec(R), \Z_p )$ denote the inverse limit
$$ \varprojlim_{k} \RGamma_{\mathet}( \Spec(R), \Z / p^{k} \Z ),$$
formed in the derived $\infty$-category $\calD(\Z)$. 

\begin{remark}\label{remark:ignore-higher-homotopy}
For every animated commutative ring $R$, the \'{e}tale site of $\Spec(R)$ agrees with the \'{e}tale site of the
underlying commutative ring $\pi_0(R)$. It follows that the cochain complexes $\RGamma_{\mathet}( \Spec(R), \Z / p^{n} \Z )$
and $\RGamma_{\mathet}( \Spec(R), \Z_p )$ depend only on $\pi_0(R)$. Nevertheless, it will be convenient in what follows to
regard both constructions as defined on the $\infty$-category $\CAlg^{\anim}$ of animated commutative rings.
\end{remark}

\begin{example}\label{example:etale-cohomology-strictly-Henselian}
For every strictly Henselian local ring $R$, the unit map
$\Z_p \rightarrow \RGamma_{\mathet}( \Spec(R), \Z_p )$ is an isomorphism (in the $\infty$-category $\widehat{\calD}(\Z_p)$).
\end{example}

\begin{example}\label{example:kroun}
Let $R$ be an animated $\F_p$-algebra. Then the constant sheaf
$\underline{\Z / p \Z}$ on $\Spec(R)_{\mathet}$ can be identified with the fiber of the Artin-Schreier map
$$ \varphi - 1: \calO_{\Spec(R)} \rightarrow \calO_{\Spec(R)}.$$
Passing to global sections, we obtain a fiber sequence
$$ \RGamma_{\mathet}( \Spec(R), \Z / p \Z ) \rightarrow R \xrightarrow{ \varphi -1 } R$$
in the derived $\infty$-category $\calD(\F_p)$.
\end{example}

\begin{remark}\label{remark:etale-coh-sifted-colim}
The functor
$$ \CAlg_{\F_p}^{\anim} \rightarrow \widehat{\calD}(\Z_p) \quad \quad R \mapsto \RGamma_{\mathet}( \Spec(R), \Z_p )$$
commutes with sifted colimits. To prove this, it suffices to show that the functor
$$ \CAlg_{\F_p}^{\anim} \rightarrow \calD(\F_p) \quad \quad R \mapsto \RGamma_{\mathet}( \Spec(R), \Z/p\Z )$$
commutes with sifted colimits, which follows immediately from Example \ref{example:kroun}.
\end{remark}

\begin{remark}\label{remark:etale-cohomology-filtered-colimit}
The functor
$$ \CAlg^{\anim} \rightarrow \widehat{\calD}(\Z_p) \quad \quad R \mapsto \RGamma_{\mathet}( \Spec(R), \Z_p )$$
commutes with filtered colimits.
\end{remark}

\begin{remark}\label{remark:Gabber-affine-proper}
Let $R$ be a commutative ring and let $I \subseteq R$ be an ideal such that the pair $(R,I)$ is Henselian.
Then the quotient map $R \rightarrow R/IR$ induces an isomorphism
$\RGamma_{\mathet}( \Spec(R), \Z/p^n \Z ) \rightarrow \RGamma_{\mathet}( \Spec(R/I), \Z / p^{n} \Z )$ for all $n$.
This is a special case of Gabber's affine proper base change theorem (see \cite{gabberaffine} or \cite{arcs}).
\end{remark}

\begin{remark}[Derived Descent]\label{remark:derived-descent-for-etale}
Let $R$ be an animated commutative ring, and let $R^{\bullet}$ denote the derived tensor product of $R$ with the
cosimplicial animated commutative ring $\F_p^{\otimes \bullet+1}$ of Notation \ref{notation:F-p-bullet}. If
$R$ is $p$-complete, then Remarks \ref{remark:ignore-higher-homotopy} and Remark \ref{remark:Gabber-affine-proper}
guarantee that the tautological map $\RGamma_{\mathet}( \Spec(R), \Z_p ) \rightarrow \RGamma_{\mathet}( \Spec(R^{n}), \Z_p )$
is an isomorphism (in the $\infty$-category $\widehat{\calD}(\Z_p)$) for every nonnegative integer $n$. In particular,
we can recover $\RGamma_{\mathet}( \Spec(R), \Z_p )$ as the totalization of the (constant) cosimplicial 
object $\RGamma_{\mathet}( \Spec(R^{\bullet}), \Z_p )$.
\end{remark}

The construction $R \mapsto \RGamma_{\mathet}( \Spec(R), \Z_p )$ determines
a functor of $\infty$-categories $\CAlg^{\anim} \rightarrow \CAlg( \widehat{\calD}(\Z_p) )$ which is characterized by the following universal property:
for every functor $\mathscr{F}: \CAlg^{\anim} \rightarrow \CAlg( \widehat{\calD}(\Z_p) )$ which satisfies descent for the \'{e}tale topology,
there is an essentially unique natural transformation of functors $\RGamma_{\mathet}( \Spec(\bullet), \Z_p ) \rightarrow \mathscr{F}(\bullet)$.
Applying this universal property to the functor $\mathscr{F}(R) = \RGamma_{\Syn}( \Spec( \widehat{R}), \Z_p)$ of Construction \ref{construction:syntomic-complex} (and invoking Proposition \ref{proposition:syntomic-sheaf-descent}), we obtain the following:

\begin{proposition}\label{proposition:etale-to-syntomic}
For every animated commutative ring $R$, there is a canonical map
$$ \alpha_{R}: \RGamma_{\mathet}( \Spec(R), \Z_p) \rightarrow \RGamma_{\Syn}( \Spf(R), \Z_p),$$
of commutative algebra objects of the $\infty$-category $\widehat{\calD}(\Z_p)$, which is characterized (up to homotopy)
by the requirement that it depends functorially on $R$.
\end{proposition}

\begin{theorem}\label{theorem:syntomic-complex-degree-zero}
Let $R$ be a $p$-complete animated commutative ring. Then the comparison map
$$ \alpha_{R}: \RGamma_{\mathet}( \Spec(R), \Z_p) \rightarrow \RGamma_{\Syn}( \Spf(R), \Z_p)$$
of Proposition \ref{proposition:etale-to-syntomic} is an isomorphism (in the $\infty$-category $\widehat{\calD}(\Z_p)$).
\end{theorem}

Theorem \ref{theorem:syntomic-chern-class-isomorphism} is essentially proven in \cite{BMS2}, at least for $p$-quasisyntomic commutative rings (see Proposition~7.16 of \cite{BMS2}). We will give a different proof here, which avoids the use of algebraic $K$-theory. 

\begin{proof}[Proof of Theorem \ref{theorem:syntomic-complex-degree-zero}]
By virtue of Proposition \ref{proposition:derived-descent-syntomic} and Remark \ref{remark:derived-descent-for-etale}, we can assume without loss of generality that $R$ has the structure of an animated $\F_p$-algebra. Since the functors
$$ \CAlg^{\anim}_{\F_p} \rightarrow \widehat{\calD}(\Z_p) \quad \quad R \mapsto \RGamma_{\mathet}( \Spec(R), \Z_p), \RGamma_{\Syn}( \Spf(R), \Z_p)$$
both commute with sifted colimits (Proposition \ref{proposition:sifted-colimit-syntomic} and Remark \ref{remark:etale-coh-sifted-colim}), we can
further reduce to the case where $R$ is a smooth $\F_p$-algebra (or even a polynomial algebra over $\F_p$).
Using Propositions \ref{proposition:absolute-vs-relative} and \ref{proposition:de-Rham-comparison-relative},
we obtain isomorphisms 
$$\F_p \otimes^{L} \Prism_{R} \simeq \F_p \otimes^{L} \Prism_{R/\Z_p} \simeq \dR_{R/ \F_p}.$$
We are therefore reduced to showing that $\alpha_{R}$ induces an isomorphism
$$ \overline{\alpha}_{R}: \RGamma_{\mathet}( \Spec(R), \Z/ p\Z) \rightarrow \fib( 
\varphi - 1: \dR_{R/\F_p} \rightarrow \dR_{R/\F_p} )$$
in the $\infty$-category $\calD(\F_p)$. We now observe that the composition
$$ \RGamma_{\mathet}( \Spec(R), \Z/ p\Z) \xrightarrow{\alpha_{R}} \fib( \varphi-1: \dR_{R/\F_p} \rightarrow \dR_{R/\F_p} )
\xrightarrow{ \epsilon_{\dR} } \fib( \varphi - 1: R \rightarrow R )$$
is the isomorphism supplied by the Artin-Schreier sequence of Example \ref{example:kroun}. It will therefore suffice to show that the
second map is an isomorphism in $\calD(\F_p)$: that is, that the map $\varphi -1$ acts by an automorphism on the complex
$\Fil^{1}_{\Hodge} \dR_{R/ \F_p }$. This is clear: since $R$ is smooth over $\F_p$, we can identify
$\Fil^{1}_{\Hodge} \dR_{R/\F_p}$ with the truncated de Rham complex $(\Omega^{\geq 1}_{R/\F_p}, d)$, on which
the Frobenius map vanishes.
\end{proof}

\begin{corollary}\label{corollary:0th-syntomic}
Let $R$ be a $p$-complete commutative ring, and suppose that the Artin-Schreier map
$$R/pR \xrightarrow{ x \mapsto x^{p} - x} R/pR$$ is a surjection. Then the syntomic complex
$\RGamma_{\Syn}( \Spf(R), \Z_p)$ is concentrated in cohomological degree zero (and is therefore
isomorphic to the abelian group $\mathrm{H}^{0}_{\mathet}( \Spec(R), \Z_p)$ of 
continuous $\Z_p$-valued functions on $\Spec(R)$).
\end{corollary}

\begin{proof}
By virtue of Theorem \ref{theorem:syntomic-complex-degree-zero}, it will suffice to show
that the complex $\RGamma_{\mathet}( \Spec(R), \Z_p)$ is concentrated in cohomological degree zero.
Since this can be tested after reduction modulo $p$, we are reduced to checking that the cohomology groups
$\mathrm{H}^{i}_{\mathet}( \Spec(R), \F_p )$ vanish for $i > 0$. By virtue of Remark \ref{remark:Gabber-affine-proper},
and Example \ref{example:kroun}, this is equivalent to the surjectivity of the Artin-Schreier map
$R/pR \xrightarrow{ x \mapsto x^{p} - x} R/pR$.
\end{proof}

\subsection{The \'{E}tale First Chern Class}\label{subsection:etale-chern-class}

Let $R$ be an animated commutative ring in which $p$ is invertible. For every integer $k \geq 0$, the finite flat group scheme
$$ \mu_{p^{k}} = \ker( \mathbf{G}_{m} \xrightarrow{ p^{k} } \mathbf{G}_m )$$
determines an invertible sheaf of $(\Z / p^{k} \Z)$-modules on the \'{e}tale site of $\Spec(R)$. For every integer $n$, we denote the $n$th tensor power of this invertible
sheaf by $(\Z / p^{k}\Z )(n)$, and we denote its derived global sections by $\RGamma_{\mathet}( \Spec(R), (\Z/p^{k}\Z)(n) )$.
Let $\RGamma_{\mathet}( \Spec(R), \Z_p(n) )$ denote the inverse limit
$$ \varprojlim_{k} \RGamma_{\mathet}( \Spec(R), (\Z / p^{k} \Z)(n) ),$$
formed in the derived $\infty$-category $\calD(\Z)$. We denote the cohomology groups of this complex by
$\mathrm{H}^{\ast}_{\mathet}( \Spec(R), \Z_p(n) )$. 

\begin{remark}\label{remark:ignore-higher-homotopy-twisted}
Let $R$ be an animated commutative ring in which $p$ is invertible, and let $\pi_0(R)$ denote the underlying commutative ring of $R$.
For every integer $n$, the tautological map $$\RGamma_{\mathet}( \Spec(R), \Z_p(n) ) \rightarrow \RGamma_{\mathet}( \Spec( \pi_0(R) ), \Z_p(n) )$$
is an isomorphism. In other words, the complex $\RGamma_{\mathet}( \Spec(R), \Z_p(n) )$ depends only on the commutative ring $\pi_0(R)$. Nevertheless, it will be convenient for us to view $R \mapsto \RGamma_{\mathet}( \Spec(R), \Z_p(n) )$ as a functor defined on animated commutative rings.
\end{remark}

\begin{example}\label{example:Zp1-etale}
Let $R$ be an animated commutative ring. If $p$ is invertible in $\pi_0(R)$, then $\RGamma_{\mathet}( \Spec(R), \Z_p(1) )$ agrees with the complex $\RGamma_{\mathet}( \Spec(R), \mathbf{G}_m)[-1]^{\wedge}$ introduced in Notation \ref{notation:etale-sheafification}. In particular, if $R$ is a commutative ring, then the $0$th cohomology group $\mathrm{H}^{0}( \Spec(R), \Z_p(1) )$ is isomorphic to the Tate module $T_p( R^{\times} )$ (see Remark \ref{remark:Tate-module-appearance}).
\end{example}

\begin{remark}
Let $R$ be an animated commutative ring in which $p$ is invertible. Then we can regard the direct sum
$$ \bigoplus_{n \in \Z} \RGamma_{\mathet}( \Spec(R), \Z_p(n) )$$
as a graded commutative ring object of the $\infty$-category $\widehat{\calD}(\Z_p)$. In particular, each
of the complexes $\RGamma_{\mathet}( \Spec(R), \Z_p(n) )$ can be regarded as a module over $\RGamma_{\mathet}( \Spec(R), \Z_p )$.
\end{remark}

\begin{example}\label{example:trivialize-roots}
For each $n \geq 1$, let $\Q( \zeta_{p^{n}} )$ denote the cyclotomic field obtained from $\Q$ by adjoining a primitive $p^n$th root of unity
$\zeta_{p^{n}}$, and let $\Q( \zeta_{p^{\infty} })$ denote the direct limit $\varinjlim_{n} \Q( \zeta_{p^{n}} )$. Then
the system $\{ \zeta_{p^{n}} \}_{n \geq 0}$ determines an element $\epsilon$ of the Tate module $$T_p( \Q( \zeta_{p^{\infty} })^{\times}) \simeq
\mathrm{H}^{0}_{\mathet}( \Spec( \Q( \zeta_{p^{\infty} }) ), \Z_p(1)).$$
If $R$ is any animated commutative algebra over $\Q( \zeta_{p^{\infty}})$, then multiplication by $\epsilon^{n}$ induces an isomorphism
$\RGamma_{\mathet}( \Spec(R), \Z_p ) \simeq \RGamma_{\mathet}( \Spec(R), \Z_p(n) )$ in the $\infty$-category $\calD(\Z)$.
Stated more informally, we have an isomorphism 
$$\bigoplus_{n \in \Z} \RGamma_{\mathet}( \Spec(R), \Z_p(n) ) \simeq \RGamma_{\mathet}( \Spec(R), \Z_p)[ \epsilon^{\pm 1} ]$$ 
of graded commutative algebras objects in the $\infty$-category $\widehat{\calD}(\Z_p)$.
\end{example}

\begin{construction}[The \'{E}tale First Chern Class]\label{construction:etale-chern}
Let $R$ be an animated commutative ring. Then Example \ref{example:Zp1-etale} determines a canonical map
$$\RGamma_{\mathet}( \Spec(R), \mathbf{G}_m )[-1] \rightarrow \RGamma_{\mathet}( \Spec(R[1/p]), \mathbf{G}_m)[-1]^{\wedge}
\simeq \RGamma_{\mathet}( \Spec(R[1/p]), \Z_p(1) ),$$
which we will denote by $c_{1}^{\mathet}$ and refer to as the {\it \'{e}tale first Chern class}.
\end{construction}

\subsection{The \'{E}tale Comparison Morphism}\label{subsection:etale-comparison-morphism}

Let $R$ be an animated commutative ring which is $p$-complete. Using Theorem \ref{theorem:syntomic-chern-class-isomorphism}, we can identify the \'{e}tale first Chern class of Construction \ref{construction:etale-chern} with a morphism of complexes
$$ \gamma_{\Syn}^{\mathet}\{1\}: \RGamma_{\Syn}( \Spf(R), \Z_p(1) ) \rightarrow
\RGamma_{\mathet}( \Spec(R[1/p]), \Z_p(1) ).$$
Our goal in this section is to show that this construction admits an essentially unique multiplicative extension to general twists.

\begin{theorem}\label{theorem:etale-comparison}
Let $R$ be a $p$-complete animated commutative ring. Then there is a canonical map
$$ (\bigoplus_{n \in \Z} \gamma_{\Syn}^{\mathet}\{n\}): \bigoplus_{n \in \Z} \RGamma_{\Syn}( \Spf(R), \Z_p(n)) \rightarrow
\bigoplus_{n \in Z} \RGamma_{\mathet}( \Spec(R[1/p]), \Z_p(n) )$$
of graded commutative ring objects of $\calD(\Z)$. This map is characterized (up to homotopy) by the
requirement that it depends functorially on $R$ and that the composition
$$ \RGamma_{\mathet}( \Spec(R), \mathbf{G}_m )[-1]
\xrightarrow{ c_{1}^{\Syn} } \RGamma( \Spf(R), \Z_p(1) ) \xrightarrow{ \gamma_{\Syn}^{\mathet}\{1\} } \RGamma_{\mathet}( \Spec(R[1/p]), \Z_p(1) )$$
coincides with the \'{e}tale first Chern class $c_{1}^{\mathet}$ of Construction \ref{construction:etale-chern} (up to a homotopy which depends functorially on $R$).
\end{theorem}

To prove Theorem \ref{theorem:etale-comparison}, we will exploit the fact that there is a large class of commutative rings $R$ for which the syntomic complexes $\RGamma_{\Syn}( \Spf(R), \Z_p(n))$ can be described explicitly. In what follows, we write $\Z_{p}^{\cyc}$ for the perfectoid ring
given by the $p$-completion of $\Z[ q^{1/p^{\infty}} ] / ( [p]_q )$, and $\epsilon \in T_p( \Z_p^{\cyc \times} )$ for the compatible
system of $p^n$th power roots of unity given by the sequence $( q^{p}, q, q^{1/p}, \cdots )$.
Note that $\Z_p^{\cyc}$ is a local ring, whose residue field can be identified (uniquely) with the finite field $\F_p$.

\begin{lemma}\label{lemma:compute-syntomic-easy}
Let $R$ be a perfectoid $\Z_p^{\cyc}$-algebra. Then:
\begin{itemize}
\item[$(a)$] For each $n \geq 2$, multiplication by the element $\log_{\Prism}(\epsilon)$ induces an isomorphism 
$$ \RGamma_{\Syn}( \Spf(R), \Z_p(n-1) ) \rightarrow \RGamma_{\Syn}( \Spf(R), \Z_p(n) )$$
in the $\infty$-category $\widehat{\calD}(\Z_p)$.

\item[$(b)$] Suppose that the groups $\Ext^{i}_{\Z_p^{\cyc}}( \F_p, R)$ vanish for every integer $i$. Then
multiplication by $\log_{\Prism}(\epsilon)$ induces an isomorphism 
$$ \RGamma_{\Syn}( \Spf(R), \Z_p) \rightarrow \RGamma_{\Syn}( \Spf(R), \Z_p(1) ).$$
\end{itemize}
\end{lemma}

\begin{proof}
To simplify the notation, we will assume that $R$ is $p$-torsion-free (this is not actually essential to the argument, but will be satisfied in the cases we apply Lemma \ref{lemma:compute-syntomic-easy}). Write $R$ as a quotient
$A/I$, where $(A,I)$ is a perfect prism. The $\Z_{p}^{\cyc}$-algebra structure on $R$ then determines
a map of prisms $(A_0, ( [p]_{q} ) ) \rightarrow (A,I)$, where $A_0$ denotes the $(p,q-1)$-completion of the ring
$\Z[ q^{1/p^{\infty}} ]$. Since $R$ is $p$-torsion-free, the ring $A$ is $(q-1)$-torsion-free.

For each $n \geq 0$, the Frobenius on $A$ induces an isomorphism $\varphi\{n\}: \Fil^{n}_{\Nyg} A\{n\} \simeq A\{n\}$. We will denote its inverse $\varphi\{n\}^{-1}$, which we regard as an endomorphism of $A\{n\}$,
so that the syntomic complex $\RGamma_{\Syn}( \Spf(R), \Z_p(n))$ is given by the fiber of the map $(1 - \varphi\{n\}^{-1}): A\{n\} \rightarrow A\{n\}$.
By virtue of Proposition \ref{proposition:twist-in-q-de-Rham-case}, multiplication by $\log_{\Prism}(\epsilon)^n$ determines
a Frobenius-equivariant isomorphism $\frac{1}{ (q-1)^{n} } A \simeq A\{n\}$. 
For $n \geq 1$, we conclude that the cofiber of the map $\log_{\Prism}(\epsilon): \RGamma_{\Syn}( \Spf(R), \Z_p(n-1)) \rightarrow \RGamma_{\Syn}( \Spf(R), \Z_p(n))$
can be identified with the two-term complex
$$ (1 - \varphi^{-1}): (\frac{ 1}{(q-1)^{n}} A) / (\frac{ 1}{(q-1)^{n-1}} A) \rightarrow (\frac{ 1}{(q-1)^{n}} A) / (\frac{ 1}{(q-1)^{n-1}} A).$$
To prove $(a)$, we must show that the differential of this complex is an isomorphism for $n \geq 2$.
For this, it suffices to observe that the endomorphism 
$$\varphi^{-1}: \frac{ 1}{(q-1)^{n}} A / \frac{ 1}{(q-1)^{n-1}} A \rightarrow \frac{ 1}{(q-1)^{n}} A / \frac{ 1}{(q-1)^{n-1}} A$$
is divisible by $p^{n-1}$ (by virtue of the congruence $[p]_{q^{1/p}}^{n} \equiv p^{n-1} [p]_{q^{1/p}} \pmod{q-1}$).

We now prove $(b)$. Assume that the groups $\Ext^{i}_{\Z_p^{\cyc}}( \F_p, R)$ vanish for every integer $i$. Writing
$\mathfrak{m}$ for the maximal ideal of $\Z_p^{\cyc}$, we conclude that the unit map
$R \rightarrow \RHom_{\Z_p^{\cyc}}( \mathfrak{m}, R )$ is invertible: that is, $R$ can be identified
with the homotopy limit of the the diagram of fractional ideals $\{ \frac{ 1}{ \zeta_{p^{m+1}} - 1} R \}_{m \geq 0}$ in the $\infty$-category $\calD(R)$. It follows that $A$ can be identified with the homotopy limit of the diagram of fractional ideals $\{ \frac{1}{ q^{1/p^m} -1 } A \}_{m \geq 0}$ in the $\infty$-category $\calD(A)$.
Consequently, to show that $1 - \varphi^{-1}$ is an automorphism of the abelian group $(\frac{1}{q-1} A) / A$, it will suffice to show that it induces an automorphism
of the quotient $(\frac{ 1}{q-1} A) / (\frac{1}{q^{1/p^m}-1} A)$ for every integer $m \geq 0$. This follows by induction on $m$, since $\varphi^{-1}$
acts by zero on each of the successive quotients $(\frac{ 1}{q^{1/p^{m-1}}-1} A) / (\frac{1}{q^{1/p^{m}}-1} A)$.
\end{proof}

\begin{proposition}\label{proposition:compute-syntomic}
Let $R$ be a perfectoid $\Z_p^{\cyc}$-algebra satisfying the following conditions:
\begin{itemize}
\item[$(a)$] The abelian groups $\Ext^{i}_{\Z_p^{\cyc}}( \F_p, R)$ vanish for all integers $i$.
\item[$(b)$] The Artin-Schreier map 
$$ R/pR \xrightarrow{ x \mapsto x^{p}-x } R/pR$$
is surjective.
\end{itemize}
Then the tautological map
$$\mathrm{H}^{0}_{\mathet}( \Spf(R), \Z_p)[ \log_{\Prism}(\epsilon) ] \rightarrow \bigoplus_{n \in \Z} \RGamma_{\Syn}( \Spf(R), \Z_p(n))$$
is an isomorphism of graded commutative ring objects of $\calD(\Z)$. In particular, each of the syntomic complexes $\RGamma_{\Syn}( \Spf(R), \Z_p(n))$ has cohomology concentrated in degree zero.
\end{proposition}

\begin{proof}
Combine Lemma \ref{lemma:compute-syntomic-easy}, Corollary \ref{corollary:0th-syntomic}, and Corollary \ref{corollary:vanishing-negative}.
\end{proof}

\begin{example}\label{example:spherical-complete}
Let $C$ be a field of characteristic zero which is complete with respect to a non-archimedean absolute value, and whose residue field has characteristic $p$.
Assume that $C$ is algebraically closed and spherically complete. Any choice of a compatible system of primitive $p^n$th roots of unity determines a morphism
from $\Z_p^{\cyc}$ to the valuation ring $\calO_{C}$ which satisfies hypotheses $(a)$ and $(b)$ of Proposition \ref{proposition:compute-syntomic}.
\end{example}

\begin{remark}\label{remark:product-closure}
The collection of perfectoid $\Z_{p}^{\cyc}$-algebras which satisfy conditions $(a)$ and $(b)$ is closed under the formation of products.
\end{remark}

\begin{remark}
For any $p$-complete commutative ring $R$, Theorem \ref{theorem:syntomic-chern-class-isomorphism}
and Remark \ref{remark:Tate-module-appearance} supply an isomorphism
$$ T_p(R^{\times} ) \xrightarrow{ \widehat{c}_{1}^{\Syn} } \mathrm{H}^{0}_{\Syn}( \Spf(R), \Z_p(1) ).$$
If $R$ satisfies the hypotheses of Proposition \ref{proposition:compute-syntomic}, then $\RGamma_{\Syn}( \Spf(R), \Z_p(1) )$ is concentrated
in cohomological degree zero, and can therefore be identified with the Tate module $T_p(R^{\times})$.
\end{remark}

\begin{proposition}\label{proposition:comparison-degree-zero}
Let $R$ be $p$-complete animated commutative ring. Then the comparison map
$$ \gamma_{\Syn}^{\mathet}\{0\}: \RGamma_{\Syn}( \Spf(R), \Z_p) \rightarrow \RGamma_{\mathet}( \Spec(R[1/p]), \Z_p)$$
of Theorem \ref{theorem:etale-comparison} is obtained by composing the
isomorphism $$\alpha_{R}: \RGamma_{\Syn}( \Spf(R), \Z_p) \simeq \RGamma_{\mathet}( \Spec(R), \Z_p)$$ of
Theorem \ref{theorem:syntomic-complex-degree-zero} with the tautological map
$\RGamma_{\mathet}( \Spec(R), \Z_p) \rightarrow \RGamma_{\mathet}( \Spec(R[1/p]), \Z_p)$ (up to a homotopy
depending functorially on $R$).
\end{proposition}

\begin{proof}
Let $\calC$ be the category of perfectoid commutative rings $R$ for which there exists a ring homomorphism
$\Z_p^{\cyc} \rightarrow R$ satisfying conditions $(a)$ and $(b)$ of Proposition \ref{proposition:compute-syntomic}.
By virtue of Example \ref{example:spherical-complete} and Remark \ref{remark:product-closure},
the category $\calC$ contains every commutative ring of the form $\prod_{i \in I} \calO_{ C_i }$, where
$\{ C_i \}_{i \in I}$ is a collection of algebraically closed fields of characteristic zero which are
complete and spherically complete with respect to non-archimedean absolute values of residue characteristic $p$.
In particular, $\calC$ forms a basis for the $\arc_{p}$ topology on the category of commutative rings
(see Definition~6.14 of \cite{arcs}). Using Corollary~6.17 of \cite{arcs}, we see that the functor
$\RGamma_{\mathet}( \Spec(R[1/p]), \Z_p)$ is a right Kan extension of its restriction to $\calC$.
It will therefore suffice to prove Proposition \ref{proposition:comparison-degree-zero} in the special case
where $R$ belongs to $\calC$. In this case, the complex $\RGamma_{\Syn}( \Spf(R), \Z_p )$
is concentrated in cohomological degree zero (Corollary \ref{corollary:0th-syntomic}); it will therefore suffice to prove the commutativity of the diagram of rings
\begin{equation}
\begin{gathered}\label{equation:silly-diagram}
\xymatrix@R=50pt@C=50pt{ \mathrm{H}^{0}_{\Syn}( \Spf(R), \Z_p) \ar[rr]^{\alpha_R} \ar[dr]_{ \gamma_{\Syn}^{\mathet}\{0\}} & &
\mathrm{H}^{0}_{\mathet}( \Spec(R), \Z_p) \ar[dl] \\
& \mathrm{H}^{0}_{\mathet}( \Spec(R[1/p]), \Z_p). & }
\end{gathered}
\end{equation}
Here $\mathrm{H}^{0}_{\mathet}( \Spec(R[1/p]), \Z_p)$ can be identified with the ring of continuous
$\Z_p$-valued functions on the topological space $\Spec( R[1/p] )$. To show that two such functions
coincide, it suffices to check this at each point of the topological space $\Spec(R[1/p])$. By functoriality,
we are reduced to checking the commutativity of the diagram (\ref{equation:silly-diagram}) in the case
where $R = \calO_{C}$ is the valuation ring of an algebraically closed field $C$ of characteristic zero which is complete and spherically complete having residue characteristic $p$. In this case, $\mathrm{H}^{0}_{\Syn}( \Spf(R), \Z_p) \simeq \Z_p$
is the initial object of the category of $p$-complete commutative rings, so the commutativity of the diagram is automatic.
\end{proof}

\begin{proof}[Proof of Theorem \ref{theorem:etale-comparison}]
Let $\calC$ be as in the proof of Proposition \ref{proposition:comparison-degree-zero}. It follows from
Corollary~6.17 of \cite{arcs} (and Remark \ref{remark:ignore-higher-homotopy-twisted}) that, for every integer $n$, the functor
$$R \mapsto \RGamma_{\mathet}( \Spec(R[1/p]), \Z_p(n) )$$ on the $\infty$-category of $p$-complete animated commutative rings is a right Kan extension of its restriction to $\calC$.
Consequently, to construct the comparison map
$$ (\bigoplus_{n \in \Z} \gamma_{\Syn}^{\mathet}\{n\}): \bigoplus_{n \in \Z} \RGamma_{\Syn}( \Spf(R), \Z_p(n)) \rightarrow
\bigoplus_{n \in \Z} \RGamma_{\mathet}( \Spec(R[1/p]), \Z_p(n) )$$
and the homotopy $\gamma_{\Syn}^{\mathet}\{1\} \circ c_{1}^{\Prism} \simeq c_{1}^{\mathet}$, we may assume without loss of generality that
$R$ belongs to $\calC$. In this case, Proposition \ref{proposition:compute-syntomic} guarantees that 
$\bigoplus_{n \in \Z} \RGamma_{\Syn}( \Spf(R), \Z_p(n))$ is isomorphic to a polynomial ring in one variable over
the ring $\mathrm{H}^{0}_{\mathet}( \Spec(R); \Z_p )$ of continuous $\Z_p$-valued functions on $\Spec(R)$;
in particular, it is concentrated in cohomological degree zero. Since the cohomology groups of
the complex $\bigoplus_{n \in \Z} \RGamma_{\mathet}( \Spec(R[1/p]), \Z_p(n) )$ are concentrated in degrees $\geq 0$,
it can be replaced by its cohomological truncation
$$ \tau^{\leq 0} (\bigoplus_{n \in \Z} \RGamma_{\mathet}( \Spec(R[1/p]), \Z_p(n) )) \simeq
\bigoplus_{n \in \Z} \mathrm{H}^{0}_{\mathet}( \Spec(R[1/p]), \Z_p(n) ).$$
Moreover, the proof of Proposition \ref{proposition:comparison-degree-zero} shows that,
for the morphism $\gamma_{\Syn}^{\mathet}\{0\}$ to depend functorially on $R$, it must coincide
with the restriction map $\mathrm{H}^{0}_{\mathet}( \Spec(R), \Z_p) \rightarrow \mathrm{H}^{0}_{\mathet}( \Spec(R[1/p]), \Z_p)$.
Let $\rho: T_p( R[1/p]^{\times} ) \rightarrow \mathrm{H}^0_{\mathet}( \Spec(R[1/p]), \Z_p(1) )$ be the isomorphism
of Example \ref{example:Zp1-etale}. Theorem \ref{theorem:etale-comparison} then reduces to the following more concrete assertion:

\begin{itemize}
\item[$(\ast)$] For $R \in \calC$, the restriction map $\mathrm{H}^{0}_{\mathet}( \Spec(R), \Z_p) \rightarrow \mathrm{H}^{0}_{\mathet}( \Spec(R[1/p]), \Z_p)$
admits a unique extension to a homomorphism of graded rings
$$ ( \bigoplus_{n \in \Z} \gamma_{\Syn}^{\mathet}\{n\}): \bigoplus_{n \in \Z} \mathrm{H}^{0}_{\Syn}( \Spf(R), \Z_p(n)) \rightarrow \bigoplus_{n \in \Z} \mathrm{H}^{0}_{\mathet}( \Spec(R[1/p]), \Z_p(n) )$$
for which the diagram of abelian groups
$$ \xymatrix@R=50pt@C=50pt{ T_p( R^{\times} ) \ar[d]^{ \log_{\Prism} }_{\sim} \ar[r] & T_p( R[1/p]^{\times} ) \ar[d]^{\rho}_{\sim} \\
\mathrm{H}^{0}_{\Syn}( \Spf(R), \Z_p(1) ) \ar[r]^-{ \gamma_{\Syn}^{\mathet}\{1\} } & \mathrm{H}^{0}_{\mathet}( \Spec(R[1/p]), \Z_p(1)) }$$
commutes. Moreover, each of the maps $\gamma_{\Syn}^{\mathet}\{n\}$ is natural in $R$.
\end{itemize}

To prove $(\ast)$, choose a homomorphism $\Z_p^{\cyc} \rightarrow R$ satisfying the hypotheses of Proposition \ref{proposition:compute-syntomic}, which we identify with an element $\epsilon \in T_p( R^{\times} )$. Proposition \ref{proposition:compute-syntomic} then guarantees that $\bigoplus_{n \geq 0} \RGamma_{\Syn}( \Spf(R), \Z_p(n))$ is isomorphic to a polynomial algebra over the ring $\mathrm{H}^{0}_{\Syn}( \Spf(R), \Z_p) \simeq \mathrm{H}^{0}_{\mathet}( \Spec(R), \Z_p)$ on the generator $\log_{\Prism}(\epsilon)$. 
It follows that there is a unique homomorphism of graded rings
$$ \beta: \bigoplus_{n \in \Z} \RGamma_{\Syn}( \Spf(R), \Z_p(n)) \rightarrow \bigoplus_{n \in \Z} \mathrm{H}^{0}_{\mathet}( \Spec(R[1/p]), \Z_p(n) )$$
which agrees with the restriction map in degree zero and satisfies $\beta( \log_{\Prism}( \epsilon) ) = \rho(\epsilon)$
We will complete the proof by showing that $\beta( \log_{\Prism}(u) ) = \rho(u)$ for every element $u \in T_p( R^{\times} )$
(in particular, this shows that the homomorphism $\beta$ does not depend on the choice of homomorphism $\Z_p^{\cyc} \rightarrow R$, and therefore depends functorially on $R$). To prove this, we note that $\mathrm{H}^{0}( \Spec(R[1/p]), \Z_p(1) )$ is a 
free $\mathrm{H}^{0}_{\mathet}( \Spec(R[1/p]), \Z_p)$-module of rank $1$, generated by the element $\rho(\epsilon)$ (Example \ref{example:trivialize-roots}). We therefore have
$$ \beta( \log_{\Prism}(u) ) = a \rho(\epsilon) \quad \quad \rho(u) = b \rho(\epsilon)$$
for unique elements $a,b \in \mathrm{H}^{0}_{\mathet}( \Spec(R[1/p]), \Z_p )$, which we can identify with continuous $\Z_p$-valued functions on the topological space $\Spec( R[1/p] )$. To show complete the proof, it will suffice to show that the functions $a$ and $b$ coincide at each point of $\Spec(R[1/p])$. As in the proof of 
Proposition \ref{proposition:comparison-degree-zero}, we can reduce to the case where $R = \calO_{C}$ is the
valuation ring of an algebraically closed field $C$ of characteristic zero which is complete and spherically complete
with respect to a non-archimedean absolute value (of residue characteristic $p$). In this case, $T_p(R^{\times} )$ is the
free $\Z_p$-module generated by $\epsilon$. We may therefore replace $u$ by $\epsilon$, so that the desired result follows by construction.
\end{proof}

\subsection{Syntomic Cohomology of Schemes}\label{subsection:syntomic-complexes-general}

We now use the \'{e}tale comparison morphism of Theorem \ref{theorem:etale-comparison}
to produce a refinement of Construction \ref{construction:syntomic-complex}.

\begin{construction}\label{construction:syntomic-complexes-general}
Let $R$ be an animated commutative ring having $p$-completion $\widehat{R}$. For every integer $n$, we let $\RGamma_{\Syn}(\Spec(R), \Z_p(n) )$
denote the pullback of the diagram
$$ \RGamma_{\Syn}( \Spf(\widehat{R}), \Z_p(n) )
\xrightarrow{ \gamma_{\Syn}^{\mathet}\{n\} } \RGamma_{\mathet}( \Spec( \widehat{R}[1/p]), \Z_p(n) )
\leftarrow \RGamma_{\mathet}( \Spec(R[1/p]), \Z_p(n)),$$
formed in the $\infty$-category $\widehat{\calD}(\Z_p)$. We will refer to the complexes $\{ \RGamma_{\Syn}( \Spec(R), \Z_p(n)) \}_{n \in \Z}$ as the {\it syntomic complexes of $R$}.
\end{construction}

\begin{remark}\label{remark:arithmetic-square-for-syntomic}
Let $R$ be an animated commutative ring and let $n$ be an integer. By construction, there is a pullback diagram
\begin{equation}
\begin{gathered}\label{equation:arithmetic-square-syntomic}
\xymatrix@R=50pt@C=50pt{ \RGamma_{\Syn}( \Spec(R), \Z_p(n)) \ar[r] \ar[d] &  \RGamma_{\mathet}( \Spec( R[1/p] ), \Z_p(n) ) \ar[d] \\
\RGamma_{\Syn}( \Spf( \widehat{R}), \Z_p(n) ) \ar[r]^-{ \gamma_{\Syn}^{\mathet}\{n\} } \ar[r] &  \RGamma_{\mathet}( \Spec( \widehat{R}[1/p] ), \Z_p(n) ), }
\end{gathered}
\end{equation}
where the upper left corner is the syntomic complex of Construction \ref{construction:syntomic-complexes-general}
and the lower left corner is the syntomic complex of Construction \ref{construction:syntomic-complex}. By
a slight abuse of terminology, we will denote the upper horizontal map by
$$\gamma_{\Syn}^{\mathet}\{n\}: \RGamma_{\Syn}( \Spec(R), \Z_p(n)) \rightarrow   \RGamma_{\mathet}( \Spec( R[1/p] ), \Z_p(n) )$$
and refer to it as the {\it \'{e}tale comparison morphism}. Note that: 

\begin{itemize}
\item If the animated commutative ring $R$ is $p$-complete, then the right vertical map in the diagram (\ref{equation:arithmetic-square-syntomic}) is invertible. It follows that the left vertical map is also invertible:
in other words, the syntomic complex $\RGamma_{\Syn}( \Spec(R), \Z_p(n))$ of Construction \ref{construction:syntomic-complexes-general} agrees with
the syntomic complex $\RGamma_{\Syn}( \Spf(R), \Z_p(n))$ of Construction \ref{construction:syntomic-complex}.

\item If $p$ is invertible in $R$, then the lower horizontal map in the diagram (\ref{equation:arithmetic-square-syntomic}) is invertible (since the source and target are both zero). It follows that the upper horizontal map is also invertible: in other words,
the \'{e}tale comparison morphism $$\gamma_{\Syn}^{\mathet}\{n\}: \RGamma_{\Syn}( \Spec(R), \Z_p(n)) \rightarrow \RGamma_{\mathet}( \Spec(R), \Z_p(n) )$$ is an isomorphism.
\end{itemize}
Consequently, we can identify (\ref{equation:arithmetic-square-syntomic}) with the diagram
$$ \xymatrix@R=50pt@C=50pt{ \RGamma_{\Syn}( \Spec(R), \Z_p(n)) \ar[r] \ar[d] & \RGamma_{\Syn}( \Spec(R[1/p]), \Z_p(n)) \ar[d] \\
\RGamma_{\Syn}( \Spec( \widehat{R}), \Z_p(n) ) \ar[r] & \RGamma_{\Syn}( \Spec( \widehat{R}[1/p]), \Z_p(n) ) . }$$
\end{remark}

\begin{remark}\label{remark:selmer-K}
For every animated commutative ring $R$, we write $K^{\mathrm{Sel}}(R, \Z_p)$ for the $p$-completed Selmer $K$-theory spectrum of $R$ (or equivalently of the
$\infty$-category $\Perf_{R}$ of perfect $R$-modules) introduced in \cite{clausen2017ktheoretic}. The spectrum $K^{\mathrm{Sel}}(R, \Z_p)$ is equipped with a decreasing
filtration $\Fil^{\bullet}_{\mot} K^{\mathrm{Sel}}(R, \Z_p)$ (which is convergent under some mild assumptions), and for every integer $n$ there is a canonical isomorphism $\gr^{n} K^{\mathrm{Sel}}(R, \Z_p) \simeq \RGamma_{\Syn}( \Spec(R), \Z_p(n))[2n]$ (where we abuse notation by identifying the complex $\RGamma_{\Syn}( \Spec(R), \Z_p(n)) \in \calD(\Z)$ with the associated generalized Eilenberg-MacLane spectrum).
\end{remark}

\begin{remark}\label{remark:extension-by-zero-comparison}
Let $R$ be an animated commutative ring and let $j: \Spec( R[1/p] ) \hookrightarrow \Spec(R)$ be the associated open immersion.
For every integer $n$ and every $k \geq 0$, we let $j_{!} (\Z / p^k \Z)(n)$ denote the extension by zero of the locally constant
sheaf $(\Z / p^{k} \Z)(n)$ on the \'{e}tale site of $\Spec( R[1/p] )$, and we write
$\RGamma_{\mathet}( \Spec(R), j_{!} \Z_p(n) )$ for the inverse limit $\varprojlim_{k} \RGamma_{\mathet}( \Spec(R), (j_{!} \Z / p^{k} \Z)(n) )$.
It follows from Theorem~6.4 of \cite{arcs} that the diagram
$$ \xymatrix@R=50pt@C=50pt{ \RGamma_{\mathet}( \Spec(R), j_{!} \Z_p(n) ) \ar[r] \ar[d] & \RGamma_{\mathet}( \Spec(R[1/p]), \Z_p(n) ) \ar[d] \\
\RGamma_{\mathet}( \Spec( \widehat{R} ), j_{!} \Z_p(n) ) \ar[r] & \RGamma_{\mathet}( \Spec(\widehat{R}[1/p]), \Z_p(n) ) }$$
is a pullback square, and from Gabber's affine proper base change theorem (\cite{gabberaffine}, \cite[Corollary 1.18]{arcs}) that the lower left corner vanishes.
Consequently, $\RGamma_{\mathet}( \Spec(R), j_{!} \Z_p(n) )$ can be identified with the fiber of the restriction map
$\RGamma_{\mathet}( \Spec(R[1/p]), \Z_p(n) )  \rightarrow \RGamma_{\mathet}( \Spec(\widehat{R}[1/p]), \Z_p(n) )$.
Applying this observation to the pullback square (\ref{equation:arithmetic-square-syntomic}), we obtain a fiber sequence
$$ \RGamma_{\mathet}( \Spec(R), j_{!} \Z_p(n) ) \rightarrow \RGamma_{\Syn}( \Spec(R), \Z_p(n))
\rightarrow \RGamma_{\Syn}( \Spf(R), \Z_p(n) )$$
in the $\infty$-category $\widehat{\calD}(\Z_p)$.
\end{remark}

\begin{example}
Let $R$ be an animated commutative ring and let $n$ be a negative integer. Then the syntomic complex 
$\RGamma_{\Syn}( \Spf(R), \Z_p(n) )$ vanishes (Corollary \ref{corollary:vanishing-negative}). Applying Remark \ref{remark:extension-by-zero-comparison}, we obtain an isomorphism $\RGamma_{\mathet}( \Spec(R), j_{!} \Z_p(n) ) \xrightarrow{\sim} \RGamma_{\Syn}( \Spec(R), \Z_p(n))$, where $j: \Spec(R[1/p]) \hookrightarrow \Spec(R)$ denotes the inclusion map.
\end{example}

\begin{proposition}\label{rhoax}
For every integer $n$, the construction $R \mapsto \RGamma_{\Syn}( \Spec(R), \Z_p(n))$ satisfies descent for the fpqc
topology.
\end{proposition}

\begin{proof}
Remark \ref{remark:extension-by-zero-comparison} supplies a fiber sequence
$$ \RGamma_{\mathet}( \Spec(R), j_{!} \Z_p(n) ) \rightarrow \RGamma_{\Syn}( \Spec(R), \Z_p(n))
\rightarrow \RGamma_{\Syn}( \Spf(R), \Z_p(n) ).$$
Here the first term satisfies arc descent (Theorem~5.4 of \cite{arcs}), and the third term
satisfies $p$-complete fpqc descent (Proposition \ref{proposition:syntomic-sheaf-descent}).
In particular, both the first and third term satisfy fpqc descent, so the middle term also satisfies fpqc descent.
\end{proof}

\begin{variant}[Globalization]\label{variant:syntomic-sheaf-integrally-globalized}
Let $X$ be a scheme, formal scheme, or algebraic stack. For every integer $n$, we let
$\RGamma_{\Syn}( X, \Z_p(n) )$ denote the limit
$$\varprojlim_{ \eta: \Spec(R) \rightarrow X } \RGamma_{\Syn}( \Spec(R), \Z_p(n)),$$
where the limit is indexed by category of pairs $(R, \eta)$ where $R$ is a commutative ring and
$\eta: \Spec(R) \rightarrow X$ is a morphism. It follows from Remark \ref{remark:arithmetic-square-for-syntomic}
that this definition is consistent with the construction of Variant \ref{variant:syntomic-sheaf-globalized}
in the case where $X$ is a bounded $p$-adic formal scheme. We denote the cohomology groups of the complex $\RGamma_{\Syn}(X, \Z_p(n) )$ by
$\mathrm{H}^{\ast}_{\Syn}(X, \Z_p(n) )$, which we refer to as {\it syntomic cohomology groups} of $X$.

When restricted to the category of schemes, the construction $X \mapsto \RGamma_{\Syn}(X, \Z_p(n) )$ is characterized the following properties:
\begin{itemize}
\item The functor $X \mapsto \RGamma_{\Syn}(X, \Z_p(n) )$ satisfies descent for the fpqc topology
(see Proposition \ref{rhoax}); in particular, it satisfies Zariski descent and is therefore determined by its values on affine schemes.
\item When $X = \Spec(R)$ is an affine scheme, $\RGamma_{\Syn}(X, \Z_p(n) )$ agrees with the syntomic complex
$\RGamma_{\Syn}( \Spec(R), \Z_p(n) )$ introduced in Construction \ref{construction:syntomic-complexes-general}.
\end{itemize}
\end{variant}

\begin{remark}\label{remark:globalized-fiber-sequence}
Let $X$ be a scheme and let $j: U \hookrightarrow X$ be the inclusion of the open subscheme $U = \Spec( \Z[1/p]) \times X$.
For every integer $n$, we have a canonical fiber sequence
$$ \RGamma_{\mathet}(X, j_{!} \Z_p(n) ) \rightarrow
\RGamma_{\Syn}( X, \Z_p(n) ) \rightarrow \fib( \varphi\{n\} - \id: \Fil^{n}_{\Nyg} \RGamma_{\Prism}(X)\{n\} \rightarrow \RGamma_{\Prism}(X)\{n\} ),$$
which specializes to the fiber sequence of Remark \ref{remark:extension-by-zero-comparison} in the case where $X = \Spec(R)$ is affine.
If the structure sheaf $\calO_{X}$ has bounded $p$-power torsion, then we can identify the third term in this fiber sequence
with the syntomic complex $\RGamma_{\Syn}( \mathfrak{X}, \Z_p(n) )$ of Variant \ref{variant:syntomic-sheaf-globalized},
where $\mathfrak{X} = \Spf(\Z_p) \times X$ is the formal completion of $X$ along the vanishing locus of $p$.
\end{remark}

\begin{example}\label{example:proper-gabber-application}
Let $R$ be a commutative ring and let $X$ be a proper $R$-scheme. Assume that the structure sheaf $\calO_{X}$ has bounded $p$-power torsion,
and let $\mathfrak{X} = \Spf(\Z_p) \times X$ be the formal completion of $X$ along the vanishing locus of $p$. If $R$ is Henselian along the ideal $(p) \subseteq R$, 
the proper base change theorem for \'{e}tale cohomology (and Gabber's affine analogue thereof) guarantee that the complex
$\RGamma_{\mathet}(X, j_{!} \Z_p(n) )$ vanishes for every integer $n$ (where $j: U \hookrightarrow X$ is the open immersion of Remark \ref{remark:globalized-fiber-sequence}).
It follows that the inclusion $\mathfrak{X} \hookrightarrow X$ induces an isomorphism $\RGamma_{\Syn}( X, \Z_p(n) ) \rightarrow \RGamma_{\Syn}( \mathfrak{X}, \Z_p(n) )$ 
in the $\infty$-category $\widehat{\calD}(\Z_p)$. In particular, we have a canonical fiber sequence
$$ \RGamma_{\Syn}(X, \Z_p(n) ) \rightarrow \Fil^{n}_{\Nyg} \RGamma_{\Prism}(X)\{n\} \xrightarrow{ \varphi\{n\} - 1} 
\RGamma_{\Prism}(X)\{n\}.$$
\end{example}

\begin{proposition}\label{proposition:syntomic-LKE}
For every integer $n$, the functor 
$$ \CAlg^{\anim} \rightarrow \widehat{\calD}(\Z_p) \quad \quad R \mapsto \RGamma_{\Syn}( \Spec(R), \Z_p(n))$$
is a left Kan extension of its restriction to the category $\CAlg_{\Z}^{\mathrm{sm}}$ of smooth $\Z$-algebras.
\end{proposition}

\begin{proof}
For every animated commutative ring $R$, let $j: \Spec(R[1/p]) \hookrightarrow \Spec(R)$ be the inclusion map
and let $F(R)$ denote the complex $\RGamma_{\mathet}( \Spec(R), j_{!} \Z_p(n) )$. Remark \ref{remark:extension-by-zero-comparison} supplies
a fiber sequence $F(R) \rightarrow \RGamma_{\Syn}( \Spec(R), \Z_p(n)) \rightarrow \RGamma_{\Syn}( \Spf(R), \Z_p(n))$ depending functorially on $R$.
By virtue of Proposition \ref{proposition:sifted-colimit-syntomic}, the functor $R \mapsto \RGamma_{\Syn}( \Spf(R), \Z_p(n))$ commutes
with sifted colimits and is therefore a left Kan extension of its restriction to the category $\Poly_{\Z}$ of finitely generated polynomial algebras over $\Z$;
in particular, it is a left Kan extension of its restriction to $\CAlg_{\Z}^{\mathrm{sm}}$. We will complete the proof by showing that the functor
$F$ is also a left Kan extension of its restriction to $\CAlg_{\Z}^{\mathrm{sm}}$.

Let $\calC$ denote the category of commutative rings which can be obtained as a filtered colimit of smooth $\Z$-algebras.
Since the functor $F$ commutes with filtered colimits, the restriction $F|_{\calC}$ is a left Kan extension of its
restriction to $\CAlg_{\Z}^{\mathrm{sm}}$. It will therefore suffice to show that the functor $F$ is a left Kan extension of $F|_{\calC}$.
Fix an animated commutative ring $R$, and let $\calC_{/R}$ denote the full subcategory of $\CAlg^{\anim}_{/R}$ spanned by those maps
$A \rightarrow R$, where $A$ is an ordinary commutative ring which can be written as a filtered colimit of smooth $\Z$-algebras. We wish to show that
the tautological map $\varinjlim_{A \in \calC_{/R} } F(A) \rightarrow F(R)$ is an isomorphism in the $\infty$-category $\widehat{\calD}(\Z_p)$.

Let $\calE \subseteq \calC_{/R}$ denote the full subcategory spanned by those maps $A \rightarrow R$ which are surjective on connected components.
Fixing one such map $A \rightarrow R$, the construction $B \mapsto A \otimes_{\Z} B$ determines a functor
$G: \calC_{/R} \rightarrow \calE$ which is equipped with a natural transformation $\id_{\calC_{/R}} \rightarrow G$. It follows
that the inclusion $\calE \hookrightarrow \calC_{/R}$ is cofinal. We are therefore reduced to showing that the natural map
$\varinjlim_{A \in \calE} F(A) \rightarrow F(R)$ is an isomorphism in $\widehat{\calD}(\Z_p)$.

Let $\calE' \subseteq \calE$ denote the full subcategory spanned by those commutative rings $A \in \calC$ equipped
with a morphism $f: A \rightarrow R$ with the property that $A$ is Henselian along the kernel ideal $I = \ker( A \rightarrow \pi_0(R) )$.
Note that the inclusion $\calE' \hookrightarrow \calE$ admits a left adjoint (given by the formation of Henselization along $I$),
and is therefore cofinal. It will therefore suffice to show that the natural map $\varinjlim_{A \in \calE} F(A) \rightarrow F(R)$ is an isomorphism in $\widehat{\calD}(\Z_p)$.
For each object $A \in \calE'$, Gabber's affine proper base change theorem (\cite{gabberaffine}) guarantees that the map $F(A) \rightarrow F(R)$ is an isomorphism.
We conclude by observing that the $\infty$-category $\calE'$ is weakly contractible (since it is nonempty and admits pairwise coproducts). 
\end{proof}

\begin{corollary}\label{corollary:filtered-colimit-syn}
For every integer $n$, the functor 
$$ \CAlg^{\anim} \rightarrow \widehat{\calD}(\Z_p) \quad \quad R \mapsto \RGamma_{\Syn}( \Spec(R), \Z_p(n))$$
commutes with filtered colimits.
\end{corollary}

\begin{warning}
Let $n$ be an integer. The functor 
$$ \CAlg^{\anim} \rightarrow \widehat{\calD}(\Z_p) \quad \quad R \mapsto \RGamma_{\Syn}( \Spec(R), \Z_p(n))$$
does not commute with sifted colimits (and is therefore not a left Kan extension of its restriction to finitely generated polynomial algebras over $\Z$). For example, suppose that $k$ is an algebraically closed field of characteristic zero and let $R$ be a $k$-algebra. Then $R$ can be realized as a sifted colimit of polynomial algebras of the form $P = k[x_1, \cdots, x_d]$. However,
the syntomic complex $\RGamma_{\Syn}( \Spec(R), \Z_p(n))$ generally cannot be realized as a sifted colimit of complexes
of the form $\RGamma_{\Syn}( \Spec(P), \Z_p(n))$, since the \'{e}tale cohomology groups $\mathrm{H}^{m}_{\mathet}( \Spec(P), \Z_p(n) )$ vanish for $m > 0$
while the cohomology groups $\mathrm{H}^{m}_{\mathet}( \Spec(R), \Z_p(n) )$ are generally nonzero.
\end{warning}

\begin{proposition}\label{proposition:keron}
For every animated commutative ring $R$, there is a canonical isomorphism
$$ \RGamma_{\mathet}( \Spec(R), \Z_p) \simeq \RGamma_{\Syn}( \Spec(R), \Z_p),$$
of commutative algebra objects of the $\infty$-category $\widehat{\calD}(\Z_p)$, which is determined (up to homotopy) by the requirement that it depends functorially on $R$.
\end{proposition}

\begin{proof}
Since the functor $R \mapsto \RGamma_{\Syn}( \Spec(R), \Z_p)$ satisfies \'{e}tale descent (Proposition \ref{rhoax}), there is an essentially unique comparison map
$\beta_{R}: \RGamma_{\mathet}( \Spec(R), \Z_p) \rightarrow \RGamma_{\Syn}( \Spec(R), \Z_p)$ which depends functorially on $R$. We wish to show that $\beta_{R}$ is
an isomorphism in the $\infty$-category $\widehat{\calD}(\Z_p)$. By virtue of Remark \ref{remark:arithmetic-square-for-syntomic}, it will suffice to prove the following:
\begin{itemize}
\item[$(a)$] The diagram 
$$ \xymatrix@R=50pt@C=50pt{ \RGamma_{\mathet}( \Spec(R), \Z_p) \ar[r] \ar[d] & \RGamma_{\mathet}( \Spec(R[1/p]), \Z_p) \ar[d] \\
\RGamma_{\mathet}( \Spec(\widehat{R}), \Z_p) \ar[r] & \RGamma_{\mathet}( \Spec(\widehat{R}[1/p]), \Z_p) }$$
is a pullback square in $\widehat{\calD}(\Z_p)$.
\item[$(b)$] The map $\beta_{R}$ is an isomorphism when $p$ is invertible in $R$.
\item[$(c)$] The map $\beta_{R}$ is an isomorphism when $R$ is $p$-complete.
\end{itemize}
Assertion $(a)$ follows from Theorem~5.13 of \cite{arcs}, assertion $(b)$ is immediate from the definitions, and assertion $(c)$ follows from Theorem \ref{theorem:syntomic-complex-degree-zero}.
\end{proof}

Let $R$ be an animated commutative ring. It follows from the defining property of the \'{e}tale comparison map 
$$\gamma_{\Syn}^{\mathet}\{1\}: \RGamma_{\Syn}( \Spf( \widehat{R}), \Z_p(1)) \to \RGamma_{\mathet}( \Spec( \widehat{R}[1/p] ), \Z_p(1) )$$ that the diagram
$$ \xymatrix@R=50pt@C=50pt{ \RGamma_{\mathet}( \Spec(R), \mathbf{G}_m)[-1] \ar[d]^{c_{1}^{\Syn}} \ar[r]^-{c_{1}^{\mathet}} & \RGamma_{\mathet}( \Spec(R[1/p]), \Z_p(1) ) \ar[d] \\
\RGamma_{\Syn}( \Spf( \widehat{R}), \Z_p(1)) \ar[r]^-{ \gamma_{\Syn}^{\mathet}\{1\} } & 
 \RGamma_{\mathet}( \Spec(\widehat{R}[1/p]), \Z_p(1) ) }$$
commutes (up to a homotopy depending functorially on $R$). We therefore obtain a map $c_{1}^{\Syn}: \RGamma_{\mathet}( \Spec(R), \mathbf{G}_m)[-1] \rightarrow \RGamma_{\Syn}( \Spec(R), \Z_p(1))$
depending functorially on $R$, which recovers the syntomic first Chern class of Proposition \ref{proposition:prismatic-chern-class-construction} in the case where $R$ is $p$-complete, and recovers the \'{e}tale first Chern class of Construction \ref{construction:etale-chern} in the case where $p$ is invertible in $R$. We will abuse terminology by referring to
$c_{1}^{\Syn}$ as the {\it syntomic first Chern class}.

\begin{proposition}\label{vrox}
Let $R$ be an animated commutative ring. Then the syntomic first Chern class $$c_{1}^{\Syn}: \RGamma_{\mathet}( \Spec(R), \mathbf{G}_m)[-1] \rightarrow \RGamma_{\Syn}( \Spec(R), \Z_p(1) )$$
exhibits $\RGamma_{\Syn}( \Spec(R), \Z_p(1) )$ as the $p$-completion of the complex
$\RGamma_{\mathet}( \Spec(R), \mathbf{G}_m)[-1]$.
\end{proposition}

\begin{proof}
Since $\RGamma_{\Syn}( \Spec(R), \Z_p(1) )$ is $p$-complete, the morphism $c_{1}^{\Syn}$
admits an essentially unique factorization
$$ \RGamma_{\mathet}( \Spec(R), \mathbf{G}_m)[-1] 
\rightarrow\RGamma_{\mathet}( \Spec(R), \mathbf{G}_m)^{\wedge}[-1] 
\xrightarrow{ \widehat{c}_{1}^{\Syn} } \RGamma_{\Syn}( \Spec(R), \Z_p(1) ).$$
By virtue of Theorem \ref{theorem:syntomic-chern-class-isomorphism}, the morphism
$\widehat{c}_{1}^{\Syn}$ is an isomorphism when $R$ is $p$-complete.
It is also an isomorphism when $p$ is invertible in $R$ (Example \ref{example:Zp1-etale}). To complete the proof, it will suffice to show that the diagram of complexes
$$ \xymatrix@R=50pt@C=50pt{ \RGamma_{\mathet}( \Spec(R), \mathbf{G}_m)^{\wedge}[-1] \ar[r] \ar[d] & \RGamma_{\mathet}( \Spec(R[1/p]), \mathbf{G}_m)^{\wedge}[-1] \ar[d] \\
\RGamma_{\mathet}( \Spec(\widehat{R}), \mathbf{G}_m)^{\wedge}[-1] \ar[r] & \RGamma_{\mathet}( \Spec(\widehat{R}[1/p]), \mathbf{G}_m)^{\wedge}[-1] }$$
is a pullback square. This can be checked after reduction modulo $p$, in which case it is a special case of Lemma~5.4.2 of \cite{cesnavicius2019purity}; alternately, one can use Proposition~\ref{proposition:p-adic-continuity} as well as the Fujiwara-Gabber theorem \cite[Corollary 6.10]{arcs} to replace $\widehat{R}$ in the bottom row of the square with the $p$-henselization $R^h$ of $R$; after this replacement, the desired property is a direct consequence of pro-Nisnevich descent of \'etale cohomology.
\end{proof}

\begin{variant}[Globalization to Schemes]\label{variant:globalization-of-Chern-class-scheme}
Let $X$ be a scheme, formal scheme, or algebraic stack. Then Proposition \ref{vrox} supplies a tautological map
$$ c_{1}^{\Syn}: \RGamma_{\mathet}( X, \mathbf{G}_m )[-1] \rightarrow \RGamma_{\Syn}( X, \Z_p(1) ),$$
which exhibits $\RGamma_{\Syn}(X, \Z_p(1) )$ as the $p$-completion of $\RGamma_{\mathet}(X, \mathbf{G}_m)[-1]$ in the derived $\infty$-category $\calD(\Z)$. Passing to cohomology in degree $2$, we obtain a homomorphism of abelian groups
$\Pic( X) \rightarrow \mathrm{H}^{2}_{\Syn}( X, \Z_p(1) )$ which we will also denote by $c_{1}^{\Syn}$ and refer to as the {\it syntomic first Chern class}.
\end{variant}

\begin{remark}[Comparison with Derived Hodge Cohomology]
Let $R$ be an animated commutative ring. We let $\gamma_{\Syn}^{\Hodge}$ denote the composite map
$$ \RGamma_{\Syn}(\Spec(R), \Z_p(n)) \rightarrow \RGamma_{\Syn}( \Spf(R), \Z_p(n) )
\rightarrow \Fil^{n}_{\Nyg} \Prism_{R}\{n\} \xrightarrow{ \gamma_{\Prism}^{\Hodge} } \widehat{\Omega}^{n}_{R}\{n\},$$
where $\gamma_{\Prism}^{\Hodge}$ is the prismatic-to-Hodge comparison morphism of Notation \ref{notation:prismatic-to-Hodge}. 
This construction depends functorially on $R$, and therefore globalizes to give a comparison map
$$ \gamma_{\Syn}^{\Hodge}: \RGamma_{\Syn}(X, \Z_p(n) ) \rightarrow \RGamma(X, L \widehat{\Omega}^{n}_{X})[-n]$$
whenever $X$ is a scheme, formal, scheme, or algebraic stack. It follows from Theorem \ref{theorem:dR-comparison-c1}
that the diagram
$$ \xymatrix@R=50pt@C=50pt{ \RGamma_{\mathet}( X, \mathbf{G}_m)[-1] \ar[r]^-{ c_{1}^{\Syn} } \ar[d]^{c_{1}^{\Hodge}} & \RGamma_{\Syn}( X, \Z_p(1) ) \ar[d]^{ \gamma_{\Syn}^{\Hodge} }\\
\RGamma(X, L \Omega^{1}_{X} )[-1] \ar[r] & \RGamma(X, L \widehat{\Omega}^{1}_{X} ) }$$
commutes (up to a canonical homotopy depending functorially on $X$), where $c_{1}^{\Hodge}$ denotes the first Chern class in Hodge cohomology
(Notation \ref{notation:hodge-c1}).
\end{remark}

\begin{variant}[The Relative First Chern Class]\label{variant:relative-first-Chern-class}
Let $f: U \rightarrow X$ be a morphism between schemes, formal schemes, or algebraic stacks.
For each integer $n$, we let $\mathrm{H}^{\ast}_{\Syn}( X,U, \Z_p(n) )$ denote the cohomology
of the fiber of the restriction map $f^{\ast}: \RGamma_{\Syn}(X, \Z_p(n) ) \rightarrow \RGamma_{\Syn}(U, \Z_p(n) )$;
we will refer to $\mathrm{H}^{\ast}_{\Syn}( X,U, \Z_p(n) )$ as the {\em syntomic cohomology of $X$ relative to $U$}.
Let $\mathscr{L}$ be a line bundle on $X$, and let $\alpha$ be a trivialization of the restriction $\mathscr{L}|_{U}$.
Then the pair $(\mathscr{L}, \alpha)$ determines a $1$-cocycle in the fiber of the restriction map
$\RGamma_{\mathet}( X, \mathbf{G}_m) \rightarrow \RGamma_{\mathet}(U, \mathbf{G}_m )$.
Applying the map $c_1^{\Syn}$ of Proposition \ref{vrox}, we obtain a relative cohomology class
$$ c_{1}^{\Syn}( \calL, \alpha) \in \mathrm{H}^{2}_{\Syn}(X, U, \Z_p(1) )$$
which we will refer to as the {\em relative syntomic first Chern class} of the pair $(\calL, \alpha )$.
\end{variant}

\subsection{Localization to the Generic Fiber}\label{subsection:localization-generic-fiber}

Throughout this section, we let $\Z[ \zeta_{p^{\infty}} ]$ denote the quotient ring $\Z[ q^{1/p^{\infty}} ] / ( [p]_q )$, and we let
$\epsilon \in T_p( \Z[ \zeta_{p^{\infty}} ]^{\times} )$ denote the compatible system of $p$th power roots of unity given by
$(1, \zeta_p, \zeta_p^2, \cdots ) = ( q^{p}, q, q^{1/p}, \cdots )$. By virtue of Proposition \ref{vrox} (and Remark \ref{remark:Tate-module-appearance}),
we can identify $\epsilon$ with an element of the abelian group $\mathrm{H}^{0}_{\Syn}( \Spec( \Z_p[ \zeta_{p^{\infty}}), \Z_p(1) )$.
For every animated commutative $\Z[ \zeta_{p^{\infty}} ]$-algebra $R$, we will abuse notation by identifying $\epsilon$
with its image in $\mathrm{H}^{0}_{\Syn}( \Spec(R), \Z_p(1) )$. Note that, in this case, the image of $\epsilon$ in
$\mathrm{H}^{0}_{\Syn}( \Spec(R[1/p]), \Z_p(1) )$ is an invertible element of the graded ring
$$\bigoplus_{n \in \Z} \mathrm{H}^{0}_{\Syn}( \Spec(R[1/p]), \Z_p(n)) \simeq \bigoplus_{n \in \Z} \mathrm{H}^{0}_{\mathet}( \Spec( R[1/p]), \Z_p(n) ).$$
We therefore obtain a comparison map
$$ (\bigoplus_{ n \in \Z} \RGamma_{\Syn}( \Spec(R), \Z_p(n)))[ 1 / \epsilon ] \rightarrow \bigoplus_{n \in \Z}
\RGamma( \Spec(R[1/p]), \Z_p(n) ).$$
Our goal is to prove the following:

\begin{theorem}\label{theorem:localization-to-generic-fiber}
Let $R$ be an animated $\Z[ \zeta_{p^{\infty}}]$-algebra. Then comparison map 
$$ (\bigoplus_{ n \in \Z} \RGamma_{\Syn}( \Spec(R), \Z_p(n)))[ 1 / \epsilon ] \rightarrow \bigoplus_{n \in \Z}
\RGamma( \Spec(R[1/p]), \Z_p(n) ).$$
becomes an isomorphism after $p$-completion (in the $\infty$-category $\widehat{\calD}(\Z_p)$).
\end{theorem}

\begin{remark}
Let $R$ be an animated $\Z[ \zeta_{p^{\infty}} ]$-algebra. The choice of compatible $p$th power roots of unity in $R$
determines a ``Bott element'' $\beta \in \pi_{2}( \Fil^{1} K^{\mathrm{Sel}}(R, \Z_p) )$, whose image in
$\gr^{1} K^{\mathrm{Sel}}(R, \Z_p) ) \simeq \SynSheaf{1}{R}[2]$ can be identified with $\epsilon$
(see Remark \ref{remark:selmer-K}). It follows from the main result of \cite{MR4110725}
that the $p$-completed localization $K^{\mathrm{Sel}}(R, \Z_p)[ \beta^{-1}]^{\wedge}$ can be identified
with the $p$-completed Selmer $K$-theory spectrum $K^{\mathrm{Sel}}( R[1/p], \Z_p)$ (see Corollary~2.22 of
\cite{MR4110725}). Theorem \ref{theorem:localization-to-generic-fiber} asserts the existence of a similar isomorphism
at the associated graded level for the motivic filtration of Remark \ref{remark:selmer-K}.
\end{remark}

As we will see below, Theorem \ref{theorem:localization-to-generic-fiber} is essentially a reformulation of the \'{e}tale comparison theorem for relative prismatic cohomology (Theorem~9.1 of \cite{prisms}), specialized to the perfectoid base ring $\Z_p^{\cyc}$. To carry out the reduction, we will use the following observation:

\begin{proposition}\label{proposition:rustle}
Let $(A,I)$ be the perfected $q$-de Rham prism (so that $A$ is the $(p,q-1)$-completion of $\Z[q^{1/p^{\infty}}]$), and let
$R$ be an animated commutative algebra over the quotient ring $\Z_p^{\cyc} = A/I$. Then the tautological map
$$ \theta: \varinjlim_{n} \Fil^{n}_{\Nyg} \Prism_R\{n\} \rightarrow \varinjlim_{n} \Prism_R\{n\}$$
is an isomorphism in the $p$-complete derived $\infty$-category $\widehat{\calD}(\Z_p)$; here the transition maps in both direct systems are given
by multiplication by the element $\log_{\Prism}( q^{p} ) \in \Fil^{1}_{\Nyg} A\{1\}$.
\end{proposition}

\begin{proof}
Note that $\theta$ can be realized as a filtered colimit of maps
$$ \theta_{m}:  \varinjlim_{n} \Fil^{n}_{\Nyg} \Prism_R\{n\} \rightarrow \varinjlim_{n} \Fil^{n-m}_{\Nyg} \Prism_R\{n\}.$$
(where we have $\Fil^{k}_{\Nyg} \Prism_{R}\{n\} = \Prism_{R}\{n\}$ for $k \leq 0$). It will therefore
suffice to show that each $\theta_m$ is an isomorphism in the $p$-complete derived $\infty$-category $\widehat{\calD}(\Z_p)$.
Proceeding by induction on $m$, we are reduced to showing that the $p$-completed colimit
$\varinjlim_{n} \gr^{n-m}_{\Nyg} \Prism_{R} \{n\}$ vanishes for every integer $m$ (where the transition maps
are given by multiplication by $\log_{\Prism}(q^p)$). Note that we can regard 
$\{ \gr^{n-m}_{\Nyg} \Prism_{R}\{n\} \}$ as a direct system in the derived $\infty$-category of the commutative ring $\gr^{0}_{\Nyg}(A)
\simeq A / ( [p]_{q^{1/p} } )$. This is a valuation ring whose $p$-adic topology coincides with the $(q^{1/p}-1)$-adic topology.
It will therefore suffice to show that each of the transition maps in the direct system
$\{ \gr^{n-m}_{\Nyg} \Prism_{R}\{n\} \}$ is divisible by the element $(q^{1/p} - 1) \in A$. This follows from the observation that
the ratio $\frac{ \log_{ \Prism }( q^{p} ) }{ q^{1/p} - 1}$ belongs to $\Fil^{1}_{\Nyg} A\{1\}$.
\end{proof}

\begin{remark}\label{remark:syntomic-colimit}
Let $R$ be a $p$-complete animated $\Z_p^{\cyc}$-algebra. Then the Frobenius maps
$\varphi\{n\}: \Fil^{n}_{\Nyg} \Prism_{R}\{n\} \rightarrow \Prism_{R}\{n\}$ of 
Notation \ref{notation:Frobenius-on-absolute}
$$ \varphi_{\infty}: \varinjlim_{n} \Fil^{n}_{\Nyg} \Prism_{R}\{n\} \rightarrow \varinjlim_{n} \Prism_{R}\{n\},$$
where the transition maps on both sides are given by multiplication by the element $\log_{\Prism}(\epsilon) \in \Fil^{1}_{\Nyg} \Prism_{\Z_p^{\cyc}}\{1\}$
(and colimits computed in the $p$-complete derived $\infty$-category $\widehat{\calD}(\Z_p)$). The map $\varphi_{\infty}$ admits an essentially unique factorization as a composition
$$ \varinjlim_{n} \Fil^{n}_{\Nyg} \Prism_{R}\{n\} \xrightarrow{\theta} \varinjlim_{n} \Prism_{R}\{n\} \xrightarrow{ \widetilde{\varphi}_{\infty}} \Prism_R\{n\},$$
where $\theta$ is the isomorphism of Proposition \ref{proposition:rustle}. We therefore obtain a fiber sequence
$$ \varinjlim_{n} \RGamma_{\Syn}( \Spf(R), \Z_p(n) ) \rightarrow \varinjlim_{n} \Prism_{R}\{n\} \xrightarrow{ \widetilde{\varphi}_{\infty} - 1} \varinjlim_{n} \Prism_{R}\{n\}$$
in the $p$-complete derived $\infty$-category $\widehat{\calD}(\Z_p)$.
\end{remark}

\begin{proof}[Proof of Theorem \ref{theorem:localization-to-generic-fiber}]
Without loss of generality, we may assume that $R$ is $p$-complete, and therefore has the structure of an animated $\Z_p^{\cyc}$-algebra. Write $\Z_p^{\cyc} = A/I$, where $(A,I)$ is the perfected $q$-de Rham prism. Theorem \ref{theorem:etale-comparison} supplies a comparison map
$$ (\bigoplus \gamma_{\Syn}^{\mathet}\{n\}): \bigoplus_{n} \RGamma_{\Syn}( \Spec(R), \Z_p(n)) \rightarrow \bigoplus_{n} \RGamma_{\mathet}( \Spec(R[1/p]), \Z_p(n) ).$$
of graded algebra objects of $\calD(\Z)$. Inverting $\epsilon$ and passing to degree zero, we obtain a comparison map
$$ \rho: \varinjlim_{n} \RGamma_{\Syn}( \Spf(R), \Z_p(n)) \rightarrow \RGamma_{\mathet}( \Spec(R[1/p]), \Z_p),$$
where the transitition maps on the left hand side are given by multiplication by $\log_{\Prism}(q^p)$ and
the colimit is formed in the $p$-complete $\infty$-category $\widehat{\calD}(\Z_p)$. We wish to show that $\rho$ is an isomorphism. 
Using Remark \ref{remark:syntomic-colimit} and Proposition \ref{proposition:twist-in-q-de-Rham-case}, we can identify $\varinjlim_{n} \RGamma_{\Syn}( \Spf(R), \Z_p(n))$
with the fiber of the map
$$ (\varphi - 1): A[ \frac{1}{q-1} ] \widehat{\otimes}_A^{L} \Prism_{R} \rightarrow A[ \frac{1}{q-1} ] \widehat{\otimes}_A^{L} \Prism_{R},$$
where the derived tensor products on each side are completed with respect to $p$ (but not with respect to $[p]_q$). Since the
ideals $(q-1)$ and $( [p]_q )$ generate the same topology on the quotient ring $A/pA$, we can rewrite this map as
$$ (\varphi -1 ): \varinjlim_{m} I^{-m} \otimes_{A}^{L} \Prism_{R} \rightarrow \varinjlim_{m} I^{-m} \otimes_{A}^{L} \Prism_{R}.$$
The proof of Theorem~9.1 of \cite{prisms} now shows that there is a unique isomorphism
$$\gamma: \RGamma_{\mathet}( \Spec(R[1/p]), \Z_p) \xrightarrow{\sim} \varinjlim_{n} \RGamma_{\Syn}( \Spf(R), \Z_p(n))$$
of commutative algebra objects of $\widehat{\calD}(\Z_p)$ which depends functorially on $R$. Note that
the composition $\rho \circ \gamma$ must be the identity on $\RGamma_{\mathet}( \Spec(R[1/p]), \Z_p)$ (since
this can be checked locally with respect to the $\arc_{p}$-topology of \cite{arcs}), so that $\rho$ is also an isomorphism.
\end{proof}

\begin{variant}\label{variant:p-quasi-stack}
Let $X$ be an algebraic stack over $\Z[ \zeta_{p^{\infty}}]$ which is quasi-compact, quasi-separated, and $p$-quasisyntomic.
For each integer $n$, let $\RGamma_{\Syn}( \mathfrak{X}, \Z_p(n) )$ denote the fiber of the canonical map
$$ (\varphi\{n\} - \id): \Fil^{n}_{\Nyg} \RGamma_{\Prism}(X)\{n\} \rightarrow \RGamma_{\Prism}(X)\{n\}.$$
Then there is a canonical pullback square
$$ \xymatrix@R=50pt@C=50pt{ \RGamma_{\Syn}(X, \Z_p(n) ) \ar[r] \ar[d] & \RGamma_{\Syn}( \mathfrak{X}, \Z_p(n) ) \ar[d] \\
\RGamma_{\mathet}( \Spec(\Z[1/p]) \times X, \Z_p(n) ) \ar[r] & \varinjlim_{m} \RGamma_{\Syn}( \mathfrak{X}, \Z_p(n+m) ); }$$
here the colimit in the lower right hand square is $p$-completed and the transition maps are given by multiplication by 
$\epsilon \in \mathrm{H}^{0}_{\Syn}( X, \Z_p(1) )$.
\end{variant}

\begin{remark}
In the statement of Variant \ref{variant:p-quasi-stack}, one can think of $\RGamma_{\Syn}( \mathfrak{X}, \Z_p(n) )$
as the syntomic cohomology of the formal completion $\mathfrak{X} = \Spf(\Z_p) \times X$, and
the colimit $\varinjlim_{m} \RGamma_{\Syn}( \mathfrak{X}, \Z_p(n+m) )$ as the \'{e}tale cohomology of its
generic fiber.
\end{remark}

\begin{proof}[Proof of Variant \ref{variant:p-quasi-stack}]
Since $X$ is quasi-compact and quasi-separated, it admits a smooth hypercovering $U_{\bullet} = \Spec( A^{\bullet} )$,
where $A^{\bullet}$ is a cosimplicial $\Z[ \zeta_{p^{\infty}} ]$-algebra. For each integer $k \geq 0$,
Theorem \ref{theorem:localization-to-generic-fiber} supplies an identification of
$\RGamma_{\mathet}( \Spec( \widehat{A}^{k}[1/p] ), \Z_p(n) )$ with the $p$-completed direct limit
$\varinjlim_{m} \RGamma_{\Syn}( \Spf( \widehat{A}^{k} ), \Z_p(n+m) )$, so that we have a pullback diagram $\sigma^{k}$:
$$ \xymatrix@R=50pt@C=50pt{
\RGamma_{\Syn}( \Spec(A^k), \Z_p(n) ) \ar[r] \ar[d] & \RGamma_{\Syn}( \Spf( \widehat{A}^{k} ), \Z_p(n) ) \ar[d] \\
\RGamma_{\mathet}( \Spec(A^k[1/p]), \Z_p(n) ) \ar[r] & \varinjlim_{m} \RGamma_{\Syn}( \Spf( \widehat{A}^{k} ), \Z_p(n+m) ), }$$
depending functorially on $k$. Our assumption on $X$ guarantees that each of the commutative rings
$A^{k}$ is $p$-quasisyntomic, so that each of the complexes appearing in the diagram $\sigma^{k}$ is coconnective.
The desired result now follows by passing to the totalization of the cosimplicial diagram $\sigma^{\bullet}$.
\end{proof}

\begin{corollary}[Purity]
\label{cor:purity-general}
Let $X$ be qcqs $p$-$p$-quasisyntomic scheme. Let $Z \subset X$ be a constructible closed subset contained in $X_{p=0} = X \times_{\Spec(\Z)} \Spec(\Z/p)$. Assume there exists an integer $d \geq 0$ such that for every affine open $\Spec(R) \subset X$, we have $\RGamma_Z(R) \in D^{\geq d}$. Then  $\RGamma_{Z}(X, \mu_{p^n}) \in D^{\geq d}$ for all $n \geq 1$.
\end{corollary}

As Akhil Mathew observed, the corollary and its proof apply to all the syntomic complexes $\RGamma_{Z,\mathrm{syn}}(X, \mathbf{Z}_p(m))/p^n$, with the version above corresponding $m=1$. 

\begin{proof}
It is enough to prove the statement for $n=1$. The hypotheses and the conclusion are stable under replacing the pair $(X,Z)$ with its base change to $\Z[ \zeta_{p^{\infty}}]$, so we may assume $X$ is equipped with the structure of a  $\Z[ \zeta_{p^{\infty}}]$-scheme. In particular, the $p$-completion $\mathfrak{X}$ is a $p$-adic formal scheme over the perfectoid ring $\Z_p^{\cyc}$, so we can identify the absolute prismatic cohomology of $\mathfrak{X}$ with the relative prismatic cohomology of $\mathfrak{X}$ over the perfect prism $(A,I)$ attached to $\Z_p^{\cyc}$. We now use the pullback square of Variant~\ref{variant:p-quasi-stack} to study $\RGamma_Z(X,\mu_{p}) \simeq \RGamma_{Z,\Syn}(X, \mathbf{Z}_p(1))/p$. As $Z \subset X_{p=0}$, we have $\RGamma_{Z,\mathet}( \Spec(\Z[1/p]) \times X, -) = 0$, so the pullback square collapses to a fibre sequence
\[ \RGamma_Z(X,\mu_{p}) \simeq \RGamma_{Z,\Syn}(X, \mathbf{Z}_p(1))/p \to \RGamma_{Z,\Syn}(\mathfrak{X}, \mathbf{Z}_p(1))/p \to  \varinjlim_{m} \RGamma_{Z,\Syn}( \mathfrak{X}, \Z_p(1+m))/p. \]
Unwrapping the definition of the second map appearing above, we learn that $\RGamma_Z(X,\mu_{p})$ is computed by applying $\RGamma_Z(\mathfrak{X},-)$ to the total fibre of the diagram
\[ \xymatrix@R=50pt@C=50pt{ \Fil^1_{\Nyg} \Prism_{\mathfrak{X}}\{1\}/p \ar[r]^-{\varphi-1} \ar[d] & \Prism_{\mathfrak{X}}\{1\}/p \ar[d] \\
\varinjlim_m \Prism_{\mathfrak{X}}\{m\}/p \ar[r]^-{\widetilde{\varphi}_\infty-1} & \varinjlim_m \Prism_{\mathfrak{X}}\{m\} }\]
in $D(\mathfrak{X},A/p)$; here the vertical maps come from multiplication by $\log_\Prism(q^p)$ as in Proposition~\ref{proposition:rustle}, and the bottom horizontal map is the one from Remark~\ref{remark:syntomic-colimit}. Trivializing the Breuil-Kisin twist via the generator $\frac{\log_\Prism(q^p)}{q-1} \in A\{1\}$ (Proposition~\ref{proposition:twist-in-q-de-Rham-case}) and tracing through the identifications, the previous diagram can be identified with the diagram
\[ \xymatrix@R=50pt@C=50pt{ \Fil^1_{\Nyg} \Prism_{\mathfrak{X}}/p \ar[r]^-{\frac{\varphi}{[p]_q}-1} \ar[d]^{\mathrm{can}} & \Prism_{\mathfrak{X}}\{1\}/p \ar[d]^{\mathrm{can}} \\
\left(\Fil^1_{\Nyg} \Prism_{\mathfrak{X}}/p\right)[\frac{1}{q-1}] \ar[r]^-{\frac{\varphi}{[p]_q}-1}  & \Prism_{\mathfrak{X}}\{1\}/p[\frac{1}{q-1}]. }\]
Thus, we learn that
\[ \RGamma_Z(X,\mu_p)\simeq \mathrm{fib}\left(\RGamma_Z(\mathfrak{X}, \RGamma_{q-1}(\Fil^1_{\Nyg} \Prism_{\mathfrak{X}}/p)) \xrightarrow{\frac{\varphi}{[p]_q}-1} \RGamma_Z(\mathfrak{X}, \RGamma_{q-1}(\Prism_{\mathfrak{X}}/p))\right).\]
Using the identification $\gr^0_{\Nyg} \Prism_{\mathfrak{X}} \simeq \mathcal{O}_{\mathfrak{X}}$, it suffices to check the following:
\[  \RGamma_{Z}(\mathfrak{X}, \RGamma_{q-1}(\Prism_{\mathfrak{X}}/p)) \in D^{\geq d} \quad \text{and} \quad \RGamma_{Z}(\mathfrak{X}, \mathcal{O}_{\mathfrak{X}}[-1])/p \in D^{\geq d}.\]
These assertions are local on $X$, so we can assume $X=\Spec(R)$ is affine. The hypothesis on $X$ implies that $\RGamma_Z(M) \simeq M \otimes_R^L \RGamma_Z(R)$ lies in $D^{\geq d}$ for any $p$-completely flat $R$-module $M$, which immediately gives the second vanishing above. For the first, observe that $[p]_q = (q-1)^{p-1} \in A/p$. As $[p]_q$ annihilates the Hodge-Tate complex $\overline{\Prism}_{\mathfrak{X}}/p$, we obtain a fibre sequence
\[ \overline{\Prism}_{\mathfrak{X}}/p[-1] \to \RGamma_{q-1}(\Prism_{\mathfrak{X}}/p) \xrightarrow{[p]_q = (q-1)^{p-1}} \RGamma_{q-1}(\Prism_{\mathfrak{X}}/p)\]
in $D(\mathfrak{X}, A/p)$. Applying the reasoning used above for the containment $\RGamma_{Z}(\mathfrak{X}, \mathcal{O}_{\mathfrak{X}}[-1])/p \in D^{\geq d}$ to the conjugate filtration then shows that $\RGamma_Z(\mathfrak{X},  \overline{\Prism}_{\mathfrak{X}}/p[-1]) \in D^{\geq d}$. By the previous fibre sequence, this implies that $(q-1)^{p-1}$ acts injectively on $\mathrm{H}^{<d}(\RGamma_{Z}(\mathfrak{X}, \RGamma_{q-1}(\Prism_{\mathfrak{X}}/p)))$. But these $A/p$-modules are also $(q-1)$-nilpotent as $\RGamma_{q-1}(\Prism_{\mathfrak{X}}/p))$ is itself $(q-1)$-nilpotent. It follows that $\mathrm{H}^{<d}(\RGamma_{Z}(\mathfrak{X}, \RGamma_{q-1}(\Prism_{\mathfrak{X}}/p))) = 0$, which gives the desired vanishing.
\end{proof}

Corollary~\ref{cor:purity-general} specializes to yield the following result: 

\begin{theorem}[Gabber's conjectures on Picard and Brauer groups, \cite{cesnavicius2019purity}]
\label{thm:cspurity}
Let $(R,\mathfrak{m})$ be a noetherian local ring of residue characteristic $p > 0$ that  is lci. In particular, if $U_R = \mathrm{Spec}(R)-\{\mathfrak{m}\}$, denotes the punctured spectrum, we have:
\begin{enumerate}
\item If $\dim(R) \geq 3$, then $\mathrm{Pic}(U_R)$ is $p$-torsionfree.
\item If $\dim(R) \geq 4$, then $\mathrm{Br}(R) \simeq \mathrm{Br}(U_R)$.
\end{enumerate}
\end{theorem}
\begin{proof}
(1): Via the Kummer sequence, it suffices to show that the map $\mathrm{H}^1(U_R. \mu_p) \to \mathrm{Pic}(U_R)$ coming from the inclusion $\mu_p \subset \mathbf{G}_m$ is the $0$ map. Consider the commutative square
\[ \xymatrix@R=50pt@C=50pt{ \mathrm{H}^1(\mathrm{Spec}(R), \mu_p) \ar[r] \ar[d] & \mathrm{Pic}(\mathrm{Spec}(R)) \ar[d] \\
\mathrm{H}^1(U_R,\mu_p) \ar[r] & \mathrm{Pic}(U_R) }\]
Now $\mathrm{Pic}(\mathrm{Spec}(R)) = 0$ as $R$ is local, so it suffices to show that the left vertical map is bijective. Let $Z=\{\mathfrak{m}\} \subset \mathrm{Spec}(R)$ be the closed point, so we have a triangle 
\[ R\Gamma_Z(\mathrm{Spec}(R),\mu_p) \to R\Gamma(\mathrm{Spec}(R),\mu_p) \to R\Gamma(U_R,\mu_p).\]
By the long exact sequence, it is enough to show that $\mathrm{H}^{<3}_Z(\mathrm{Spec}(R),\mu_p) = 0$. But this follows from Corollary~\ref{cor:purity-general} as $R$ is lci (and hence Cohen-Macaulay) of dimension $\geq 3$.

(2): The claim is classical for prime-to-$p$ torsion classes, so we focus on the $p$-primary torsion. Fix an integer $n \geq 1$. Consider the commutative diagram
\[ \xymatrix@R=50pt@C=50pt{ \mathrm{H}^2(\mathrm{Spec}(R),\mu_{p^n}) \ar[r] \ar[d] & \mathrm{H}^2(U_R,\mu_{p^n}) \ar[d] \\
     \mathrm{H}^2(\mathrm{Spec}(R),\mathbf{G}_m)[p^n] = \mathrm{Br}(R)[p^n] \ar[r] & \mathrm{H}^2(U_R,\mathbf{G}_m)[p^n] = \mathrm{Br}(U_R)[p^n] }\]
     where the horizontal maps come from functoriality while the vertical maps come from the inclusion $\mu_{p^n} \subset \mathbf{G}_m$. By the Kummer sequence, both vertical maps are sujrective. Moreover, as $\mathrm{H}^1(\mathrm{Spec}(R),\mathbf{G}_m) = 0$, the left vertical map is bijective; similarly, it is classically known that $\mathrm{H}^1(U_R,\mathbf{G}_m) = 0$ under the hypotheses in (2) (see \cite[Theorem XI.3.1.13 (ii)]{SGA2}, so the right vertical map is also bijective. Finally, applying the same analysis used in part (1) one degree higher (possible as $\dim(R) \geq 4$) shows that the top horizontal map is bijective. Consequently, the bottom horizontal map is also bijective. As this is true for all $n$, we conclude that $\mathrm{Br}(R) \to \mathrm{Br}(U_R)$ is an isomorphism on $p$-primary components.
 \end{proof}

Theorem~\ref{thm:cspurity} is one of the main results of \cite{cesnavicius2019purity}, powering the proof of Gabber's conjectures on Picard and Brauer groups. The proof given above, at its core, is not significantly different from the one in \cite{cesnavicius2019purity}: both rely ultimately on the ability to access $\mu_p$-cohomology via coherent cohomology through perfectoid/prismatic techniques. The main difference is that our work on syntomic cohomology of schemes allows us to streamline the presentation by avoiding the need to translate Theorem~\ref{thm:cspurity} to a statement about perfectoid rings.

\newpage \section{Calculations in Syntomic Cohomology}\label{section:calculate-syntomic}

Let $X$ be a scheme. To every line bundle $\mathscr{L}$ on $X$, Variant \ref{variant:globalization-of-Chern-class-scheme}
associates a syntomic cohomology class $c_{1}^{\Syn}( \mathscr{L} ) \in \mathrm{H}^{2}_{\Syn}( X, \Z_p(1) )$.
Our goal in this section is to extend the theory of syntomic Chern classes to vector bundles of higher rank.
Our main result can be summarized as follows:

\begin{theorem}[Syntomic Chern Classes]\label{theorem:syntomic-chern-existence}
To every scheme $X$ and every vector bundle $\mathscr{E}$ on $X$, one can assign a system of syntomic cohomology classes
\[ c_{n}^{\Syn}(\mathscr{E} ) \in \mathrm{H}^{2n}_{\Syn}( X, \Z_p(n) ).\]
This assignment is uniquely determined by the following requirements:
\begin{enumerate}
\item {\bf Functoriality:} Let $f: Y \rightarrow X$ be a morphism of schemes and let $\mathscr{E}$ be a vector bundle on $X$.
Then, for every integer $n \geq 0$, we have an equality
\[ c_{n}^{\Syn}( f^{\ast} \mathscr{E} ) = f^{\ast}( c_{n}^{\Syn}( \mathscr{E} ) ) \in \mathrm{H}^{2n}_{\Syn}( Y, \Z_p(n) ).\]

\item {\bf Normalization:} Let $\mathscr{E}$ be a vector bundle on a scheme $X$. Then $c_{0}^{\Syn}( \mathscr{E} ) = 1$,
and $c_{n}^{\Syn}( \mathscr{E} )$ vanishes for $n > \mathrm{rank}(\mathscr{E} )$. Moreover, if $\mathscr{E}$
is a line bundle, then $c_{1}^{\Syn}(\mathscr{E} ) \in \mathrm{H}^{2}_{\Syn}(X, \Z_p(1) )$ is the syntomic cohomology class defined 
Variant \ref{variant:globalization-of-Chern-class-scheme}.

\item {\bf Additivity:} Let $0 \rightarrow \mathscr{E}' \rightarrow \mathscr{E} \rightarrow \mathscr{E}'' \rightarrow 0$ be a short
exact sequence of vector bundles on $X$. Then, for every integer $n \geq 0$, we have
\[ c_{n}^{\Syn}( \mathscr{E} ) = \sum_{n = n' + n''} c_{n'}^{\Syn}( \mathscr{E}' ) \cdot c_{n''}^{\Syn}( \mathscr{E}'' ).\]
\end{enumerate}
\end{theorem}

The proof of Theorem \ref{theorem:syntomic-chern-existence} follows a standard pattern. We begin in
\S\ref{subsection:projective-bundle-scheme} by showing that syntomic cohomology satisfies a projective bundle formula.
More precisely, we show that if $\mathscr{E}$ is a vector bundle of rank $r$ over a scheme $X$, then the syntomic cohomology
of the projective bundle $\mathbf{P}( \mathscr{E} )$ is described by the formula
$$ \mathrm{H}^{\ast}_{\Syn}( \mathbf{P}(\mathscr{E}), \Z_p(n) )
\simeq \bigoplus_{0 \leq i < r}  \mathrm{H}^{\ast-2i}_{\Syn}(X, \Z_p(n-i) ) c_{1}^{\Syn}( \calO(1) )^{i}$$
(see Theorem \ref{theorem:projective-bundle-formula}).

From the projective bundle formula, we can immediately deduce the uniqueness assertion of Theorem \ref{theorem:syntomic-chern-existence}.
Let $\mathscr{E}$ be a vector bundle of rank $r$ on a scheme $X$, and let $\pi: \mathrm{Flag}(\mathscr{E} ) \rightarrow X$ be the scheme parametrizing
complete flags in $\mathscr{E}$, so that the pullback $\pi^{\ast}( \mathscr{E} )$ admits a filtration
$$ \mathscr{F}_{1} \hookrightarrow \mathscr{F}_2 \hookrightarrow \mathscr{F}_3 \hookrightarrow
\cdots \hookrightarrow \mathscr{F}_r = \pi^{\ast}(\mathscr{E})$$
where each quotient $\mathscr{L}_{i} = \mathscr{F}_{i} / \mathscr{F}_{i-1}$ is a line bundle on
$\mathrm{Flag}( \mathscr{E} )$. Using the projective bundle formula repeatedly, we deduce that the pullback map
$$ \pi^{\ast}: \mathrm{H}^{\ast}_{\Syn}( X, \Z_p(n) ) \rightarrow \mathrm{H}^{\ast}_{\Syn}( \mathrm{Flag}(\mathscr{E}), \Z_p(n) )$$
is a monomorphism. Consequently, the axioms of Theorem \ref{theorem:syntomic-chern-existence} allow us to characterize
the syntomic Chern class $c_{n}^{\Syn}( \mathscr{E} )$ as the unique element whose image in $\mathrm{H}^{2n}_{\Syn}( \mathrm{Flag}(\mathscr{E}), \Z_p(n) )$ is equal to the $n$th elementary
symmetric function of the Chern classes $\{ c_1^{\Syn}( \mathscr{L}_{i} ) \}_{0 \leq i \leq r}$.

To prove the existence assertion of Theorem \ref{theorem:syntomic-chern-existence}, we need to work a bit harder.
Let $\mathscr{E}$ be a vector bundle of rank $r$ on a scheme $X$, and let $t = c_{1}^{\Syn}( \calO(1) ) \in \mathrm{H}^{2}_{\Syn}( \mathbf{P}( \mathscr{E} ), \Z_p(1) )$.
It follows from the projective bundle formula that $t$ satisfies a unique monic polynomial of degree $r$
$$ t^{r} + c_1 \cdot t^{r-1} + c_2 \cdot t^{r-2} + \cdots + c_{r-1} \cdot t + c_{r} = 0,$$
where each coefficient $c_{n}$ belongs to $\mathrm{H}^{2n}_{\Syn}( X, \Z_p(n) )$. We will
define $c_{n}^{\Syn}( \mathscr{E} )$ to be the coefficient $c_{n}$ which appears in this identity (Construction \ref{construction:higher-chern-classes}).
It follows immediately that this definition satisfies the first two axioms of Theorem \ref{theorem:syntomic-chern-existence} follow immediately. However, the additivity formula for syntomic Chern classes requires some additional arguments, which we carry out in \S \ref{subsection:higher-chern} (Theorem \ref{theorem:additivity-chern}).

The syntomic Chern classes $\{ c_{n}^{\Syn}( \mathscr{E} ) \}_{n \geq 0}$ are essentially
the {\em only} invariants in syntomic cohomology which can be associated to an algebraic vector bundle $\mathscr{E}$. To make this precise,
it is convenient to consider the universal case where $X = \BGL_{r}$ is the classifying stack of the general linear group $\GL_{r}$,
and $\mathscr{E} = \mathscr{E}_{\mathrm{univ} }$ is the tautological rank $r$ vector bundle on $X$. In \S\ref{subsection:cohomology-classifying-stack},
we show that the the syntomic cohomology classes $\{ c_{i}^{\Syn}( \mathscr{E} ) \}_{1 \leq i \leq r}$ form polynomial generators
for bigraded cohomology ring
$$ \bigoplus_{n \in \Z} \mathrm{H}^{\ast}_{\Syn}( \BGL_r, \Z_p(n) )$$
as an algebra over $\bigoplus_{n \in \Z} \mathrm{H}^{\ast}_{\Syn}( \Spec(\Z), \Z_p(n) )$ (for a more general statement, see Theorem \ref{theorem:syntomic-cohomology-of-BGL}).

We close this section by discussing a companion to the projective bundle formula of \S\ref{subsection:projective-bundle-scheme}. Let
$X$ be a scheme, let $\widetilde{X}$ be the blowup of $X$ along a closed subscheme $Y \subseteq X$, and let $D \subseteq \widetilde{X}$
be the exceptional divisor. Assume that $Y$ can be described locally as the vanishing locus of a regular sequence.
In \S\ref{subsection:blowup-formula}, we show that the commutative diagram of schemes
$$ \xymatrix@R=50pt@C=50pt{ D \ar[r] \ar[d] & \widetilde{X} \ar[d] \\
Y \ar[r] & X }$$
determines a pullback diagram after passing to syntomic complexes (Theorem \ref{theorem:syntomic-blow-up}). Under some mild additional assumptions,
we can combine this with the projective bundle formula (applied to the map $D \rightarrow Y$) to obtain a simple description of the
syntomic cohomology groups of $\widetilde{X}$ (in terms of those of $X$ and $Y$); see Corollary \ref{corollary:blowup-application}.

\subsection{The Projective Bundle Formula}\label{subsection:projective-bundle-scheme}

Let $X$ be a scheme, formal scheme, or algebraic stack, and let $\mathscr{E}$ be a vector bundle on $X$.
We write $\mathbf{P}( \mathscr{E} )$ for the projectivization of $\mathscr{E}$, so that $\mathbf{P}( \mathscr{E} )$
is equipped with a proper smooth morphism $\pi: \mathbf{P}( \mathscr{E} ) \rightarrow X$. We adopt the convention that $\mathbf{P}(\calE)$ parametrizes subbundles of $\mathscr{E}$ of rank $1$, so that the pullback $\pi^{\ast}(\calE)$ is equipped with a tautological line subbundle $\calO(-1) \subseteq \pi^{\ast}(\mathscr{E})$, whose inverse we denote by $\calO(1)$. This line bundle 
$\calO(1)$ has a syntomic first Chern class (see Variant \ref{variant:globalization-of-Chern-class-scheme}), which we denote by $c_{1}^{\Syn}( \calO(1) ) \in \mathrm{H}^{2}_{\Syn}( \mathbf{P}( \mathscr{E} ), \Z_p(1) )$. The goal of this section is to prove the following:

\begin{theorem}[Projective Bundle Formula]\label{theorem:projective-bundle-formula}
Let $X$ be a scheme, formal scheme, or algebraic stack. Let $\mathscr{E}$ be a vector bundle of rank $r$ on $X$.
Then, for every integer $n$, the syntomic cohomology classes $\{ c_{1}^{\Syn}( \calO(1) )^{i} \}_{0 \leq i < r}$ induce an isomorphism 
$$ \bigoplus_{0 \leq i < r} \RGamma_{\Syn}(X, \Z_p(n-i) )[ -2i] \rightarrow \RGamma_{\Syn}( \mathbf{P}(\mathscr{E}), \Z_p(n))$$
in the derived $\infty$-category $\widehat{\calD}(\Z_p)$.
\end{theorem}

We will deduce Theorem \ref{theorem:projective-bundle-formula} from the following classical calculation:

\begin{lemma}\label{lemma:classical-hodge-calculation}
Let $R$ be a commutative ring and let $\mathbf{P}^{k}$ denote projective space of dimension $k$ over $\Spec(R)$.
For $0 \leq d \leq k$, the cohomology group
$\mathrm{H}^{d}( \mathbf{P}^{k}_{R}, \Omega^{d}_{ \mathbf{P}^{k} / R} )$ is a free $R$-module of rank $1$, generated by the image
of $c_{1}^{\Hodge}( \calO(1) )^{d}$. Moreover, the cohomology groups $\mathrm{H}^{\ast}( \mathbf{P}^{k}_{R}, \Omega^{d}_{ \mathbf{P}^{k} / R} )$ vanish for $\ast \neq d$.
\end{lemma}

\begin{proof}
See \cite[\href{https://stacks.math.columbia.edu/tag/0FMI}{Tag 0FMI}]{stacks-project}
\end{proof}

We now prove a variant of Lemma \ref{lemma:classical-hodge-calculation} for derived Hodge cohomology.

\begin{lemma}\label{lemma:straight}
Let $f: X \rightarrow Y$ be a morphism of schemes, let $\mathbf{P}^{r-1}_{X}$ be an $(r-1)$-dimensional projective space over $X$, and let
$t$ denote the first Chern class of the line bundle $\calO(1)$ in the derived Hodge cohomology group $\mathrm{H}^{1}( \mathbf{P}^{r-1}_{X}, L \Omega^{m}_{ \mathbf{P}^{r-1}_{X} / Y} )$. Then, for every integer $d$, the map
$$ \bigoplus_{0 \leq i < r} \RGamma(X, L \Omega^{d-i}_{X/Y} )[ -i ] \xrightarrow{ (1, t, \cdots, t^{r-1})} \RGamma( \mathbf{P}^{k}_{X}, L \Omega^{d}_{ \mathbf{P}^{k}_{X}/Y} )$$
is an isomorphism in $\calD(\Z)$.
\end{lemma}

\begin{proof}
Without loss of generality, we may assume that $Y = \Spec(R)$ is affine and that $X$ is quasi-compact and quasi-separated. In this case,
we can use the K\"{u}nneth formula of Remark \ref{remark:Kunneth-non-affine} to reduce to the case $X = Y$, in which case the desired
result reduces to Lemma \ref{lemma:classical-hodge-calculation}.
\end{proof}

\begin{lemma}\label{lemma:Hodge-basis-diffracted-scheme}
Let $X$ be a bounded $p$-adic formal scheme, let $\mathbf{P}^{r-1}_{X}$ be projective space of dimension $r-1$ over $X$,
and let $t$ denote the syntomic first Chern class $c_{1}^{\Syn}( \calO(1) )$. Then:
\begin{itemize}
\item[$(1)$] For every integer $d$, the elements $\{ t^i \}_{0 \leq i < r}$ induce an isomorphism
$$ \bigoplus_{0 \leq i < r} \Fil_{d-i}^{\conj} \RGamma(X, \widehat{\Omega}^{\DHod}_{X})[-2i]
\xrightarrow{} \Fil_{d}^{\conj} \RGamma( \mathbf{P}^{k}_{X},  \widehat{\Omega}^{\DHod}_{\mathbf{P}^{r-1}_{X}}).$$

\item[$(2)$] The elements $\{ t^i \}_{0 \leq i < r}$ induce an isomorphism
$$  \bigoplus_{0 \leq i < r} \RGamma(X, \widehat{\Omega}^{\DHod}_{X})[-2i] \rightarrow
\RGamma( \mathbf{P}^{k}_{X},  \widehat{\Omega}^{\DHod}_{\mathbf{P}^{r-1}_{X}}).$$

\item[$(3)$] For every integer $n$, the elements $\{ t^i \}_{0 \leq i < r}$ induce an isomorphism
$$ \bigoplus_{0 \leq i < r} \RGamma_{ \overline{\Prism}}(X)\{n-i\}[-2i] \rightarrow  \RGamma_{\overline{\Prism}}( \mathbf{P}^{r-1}_{X})\{n\}.$$

\item[$(4)$] For every integer $n$, the elements $\{ t^i \}_{0 \leq i < r}$ induce an isomorphism
$$ \bigoplus_{0 \leq i < r} \RGamma_{ \Prism}(X)\{n-i\}[-2i] \rightarrow \RGamma_{\Prism}( \mathbf{P}^{r-1}_{X})\{n\}.$$

\item[$(5)$] For every pair of integers $m$ and $n$, the elements $\{ t^i \}_{0 \leq i < r}$ induce an isomorphism
$$ \bigoplus_{0 \leq i < r} \Fil^{m-i}_{\Nyg} \RGamma_{ \Prism}(X)\{n-i\}[-2i] \rightarrow  \Fil^{m}_{\Nyg} \RGamma_{\Prism}( \mathbf{P}^{r-1}_{X})\{n\}.$$

\item[$(6)$] For every integer $n$, the elements $\{ t^i \}_{0 \leq i < r}$ induce an isomorphism
$$ \bigoplus_{0 \leq i < r} \RGamma_{\Syn}( X, \Z_p(n-i) )[-2i] \rightarrow \RGamma_{\Syn}(\mathbf{P}^{r-1}_{X}, \Z_p(n) ).$$
\end{itemize}
\end{lemma}

\begin{proof}
Without loss of generality we may assume that $X$ is a scheme (with $p$ nilpotent in the structure sheaf $\calO_X$).
To prove $(1)$, we proceed by induction on $d$. The case $d < 0$ is vacuous. It will therefore suffice to show that
the elements $\{ t^i \}_{0 \leq i < r}$ induce an isomorphism
$$ \bigoplus_{0 \leq i < r} \gr_{d-i}^{\conj} \RGamma(X, \widehat{\Omega}^{\DHod}_{X})[-2i]
\rightarrow \gr_{d}^{\conj} \RGamma( \mathbf{P}^{r-1}_{X},  \widehat{\Omega}^{\DHod}_{\mathbf{P}^{r-1}_{X}})$$
which is a special case of Lemma \ref{lemma:straight} (applied to the ground scheme $Y = \Spec(\Z)$).

Assertion $(2)$ follows from $(1)$ by passing to the colimit over $d$, and assertion $(3)$ follows by combining $(2)$ with Remark \ref{remark:diffracted-vs-prismatic2}. We now prove $(4)$. For each integer $k \geq 0$, let $\RGamma_{\Prism}^{[k]}(X)\{n\}$ denote the limit
$$ \varprojlim_{\Spec(R) \rightarrow X} \Prism_{R}^{[k]}\{n\},$$ and define $\RGamma_{\Prism}^{[k]}( \mathbf{P}^{r-1}_{X} )\{n\}$ similarly. It will then suffice to show that, for each $k \geq 0$, the elements $\{ t^i \}_{0 \leq i < r}$ induce an isomorphism
$$ \bigoplus_{0 \leq i < r} (\RGamma_{\Prism}( X)\{n-i\} / \RGamma_{\Prism}^{[k]}( X)\{n-i\})[-2i] 
\rightarrow \RGamma_{\Prism}( \mathbf{P}^{r-1}_{X})\{n\} / \RGamma_{\Prism}^{[k]}( \mathbf{P}^{r-1}_{X})\{n\}.$$
This follows by induction on $k$, using assertion $(3)$.

To prove $(5)$, we proceed by induction on $m$. For $m \leq 0$, the desired result follows from $(4)$. To handle the inductive step, it will suffice to show that the elements $\{ t^i \}_{0 \leq i < r}$ induce an isomorphism
$$ \bigoplus_{0 \leq i < r} \gr^{m}_{\Nyg} \RGamma_{\Prism}( X)\{n-i\}[-2i] \rightarrow \gr^{m}_{\Nyg} \RGamma_{\Prism}( \mathbf{P}^{r-1}_{X})\{n\}.$$
This follows by combining $(1)$ with the fiber sequence of Remark \ref{remark:Nygaard-associated-graded}.

Assertion $(6)$ is an immediate consequence of $(5)$.
\end{proof}

We will also need the analogue of Lemma \ref{lemma:Hodge-basis-diffracted-scheme} for \'{e}tale cohomology:

\begin{lemma}\label{lemma:etale-calculation}
Let $C$ be a separably closed field of characteristic $\neq p$ and let $\mathbf{P}^{r-1}_{C}$ denote projective space of dimension $(r-1)$ over $C$.
Choose a primitive $p$th root of unity in $C$, so that the \'{e}tale Chern class $c_{1}^{\mathet}( \calO(1 ) )$ determines an element
$t \in \mathrm{H}^{2}_{\mathet}( \mathbf{P}^{r-1}_{C}, \F_p)$. Then the elements $\{ t^{i} \}_{0 \leq i < r}$ form a basis for the \'{e}tale
cohomology ring $\mathrm{H}^{\ast}_{\mathet}( \mathbf{P}^{r-1}_{C}, \F_p)$ as a vector space over $\F_p$.
\end{lemma}

\begin{proof}
This is a standard calculation in \'{e}tale cohomology. However, it is amusing to note that it can also be deduced from Lemma \ref{lemma:Hodge-basis-diffracted-scheme}. Since the projection map $\pi: \mathbf{P}^{r-1}_{ \Z} \rightarrow \Spec(\Z)$ is smooth and proper, the (derived) direct image
$\pi_{\ast} \underline{\F_p}$ is locally constant when restricted to $\Spec( \Z[1/p] )$ (see \cite{SGAIV}, Expos\'{e} XIV, Corollaire~2.2).
Consequently, the conclusion of Lemma \ref{lemma:etale-calculation} does not depend on the choice of $C$. We may therefore assume without loss of generality
that $C$ is an algebraically closed field of characteristic zero which is complete with respect to nonarchimedean absolute value having residue characteristic $p$.
Let $\calO_{C}$ be valuation ring of $C$, set $X = \Spf(\calO_{C} )$, and let $\epsilon$ be a generator of $\mathrm{H}^{0}_{\Syn}( X, \Z_p(1) ) \simeq T_p( \calO_C^{\times} )$ (corresponding to a compatible system of $p^n$th roots of unity in $C$). Applying Lemma \ref{lemma:Hodge-basis-diffracted-scheme},
we deduce that the powers of $c_{1}^{\Syn}( \calO(1) )$ induce an isomorphism
$$ \bigoplus_{n \in \Z} \bigoplus_{0 \leq i < r} \RGamma_{\Syn}( X, \Z_p(n-i) )[-2i] \rightarrow \bigoplus_{n \in \Z} \RGamma_{\Syn}( \mathbf{P}^{r-1}_{X}, \Z_p(n) ).$$
The conclusion of Lemma \ref{lemma:etale-calculation} follows by reducing modulo $p$, inverting the element $\epsilon$, and applying Theorem \ref{theorem:localization-to-generic-fiber}.
\end{proof}

\begin{proof}[Proof of Theorem \ref{theorem:projective-bundle-formula}]
Without loss of generality, we may assume that $X = \Spec(R)$ is an affine scheme and that $\mathscr{E} = \calO_{X}^{r}$ is a trivial vector bundle of rank
on $X$, so that $\mathbf{P}(\mathscr{E})$ can be identified with the projective space $\mathbf{P}^{r-1}_{X}$. Set $$t = c_{1}^{\Syn}( \calO(1) ) \in \mathrm{H}^{2}_{\Syn}( \mathbf{P}^{r-1}_{X}, \Z_p(1) ).$$
We wish to show that, for every integer $n$, the syntomic cohomology classes $\{ t^i \}_{0 \leq i < r}$ induce an isomorphism
$$ \bigoplus_{0 \leq i < r} \RGamma_{\Syn}(X, \Z_p(n-i) )[-2i] \rightarrow \RGamma_{\Syn}( \mathbf{P}^{r-1}_{X}, \Z_p(n) ).$$
Writing $R$ as a filtered colimit of finitely generated subrings (and using Corollary \ref{corollary:filtered-colimit-syn}), we may further assume that
the ring $R$ is finitely generated and therefore has bounded $p$-power torsion. Let $\mathfrak{X} = \Spf( \widehat{R} )$ be the formal completion of $X$ along the vanishing locus of $p$, and let $U  = \Spec( R[1/p] )$ be the open subscheme of $X$ where $p$ is invertible, so that we have a pullback diagram
of schemes
$$ \xymatrix@R=50pt@C=50pt{ \mathbf{P}^{r-1}_{U} \ar[r]^-{ \overline{j} } \ar[d] & \mathbf{P}^{r-1}_{X} \ar[d]^{\pi} \\
U \ar[r]^-{j} & X }$$
where the horizontal maps are open immersions. 

Using Remark \ref{remark:globalized-fiber-sequence}, we obtain a commutative diagram of fiber sequences
$$ \xymatrix@R=50pt@C=50pt{ \bigoplus_{0 \leq i < r} \RGamma_{\mathet}(X, j_{!} \Z_p(n-i) )[-2i] 
\ar[r]^-{\theta'} \ar[d] & \RGamma( \mathbf{P}^{r-1}_{X}, \overline{j}_{!} \Z_p(n) ) \ar[d] \\
 \bigoplus_{0 \leq i < r} \RGamma_{\Syn}(X, \Z_p(n-i) )[-2i] \ar[r]^-{\theta} \ar[d] & 
\RGamma_{\Syn}( \mathbf{P}^{r-1}_{X}, \Z_p(n) ) \ar[d] \\
 \bigoplus_{0 \leq i < r} \RGamma_{\Syn}(\mathfrak{X}, \Z_p(n-i) )[-2i]\ar[r]^-{\theta''} & 
\RGamma_{\Syn}( \mathbf{P}^{r-1}_{\mathfrak{X}}, \Z_p(n) ), }$$
where $\theta''$ is an isomorphism by virtue of Lemma \ref{lemma:Hodge-basis-diffracted-scheme}. It will
therefore suffice to show that $\theta'$ is an isomorphism. Unwinding the definitions, we see that the reduction of
$\theta'$ modulo $p$ can be obtained by applying the functor $\mathscr{F} \mapsto \RGamma_{\mathet}(X, \mathscr{F} )$
to a morphism
$$ \rho: j_{!}( \bigoplus_{0 \leq i < r} \F_p(n-i)[-2i] ) \rightarrow R\pi_{\ast} \overline{j}_{!}( \underline{\F_p}(n) )$$
of complexes of \'{e}tale sheaves on $X$. To complete the proof, it will suffice to show that $\rho$ induces
an isomorphism after taking stalks at each geometric point $\Spec(k) \rightarrow X$. Using the proper base change theorem for \'{e}tale cohomology, we can replace $X$ by $\Spec(k)$ and thereby reduce to the case where $X$ is the spectrum of a separably closed field $k$. We may further assume that $k$ has characteristic $\neq p$ (otherwise, the open set $U$ is empty and there is nothing to prove). In this case, the desired result follows from Lemma \ref{lemma:etale-calculation}.
\end{proof}

Theorem \ref{theorem:projective-bundle-formula} has analogues for other cohomology theories studied in this paper, which can be proved using the same arguments. For example:

\begin{variant}\label{variant:projectivelineBKtwist}
Let $(A,I)$ be a bounded prism, let $\mathfrak{X}$ be a bounded $p$-adic formal scheme over the quotient ring $A/I$,
and let $\mathscr{E}$ be a vector bundle of rank $r$ on $\mathfrak{X}$. Then the prismatic cohomology classes
$\{ c_{1}^{\Prism}( \calO(1) )^{i} \}_{0 \leq i < r }$ induce an isomorphism
$$ \bigoplus_{0 \leq i < r } \mathrm{H}^{\ast-2i}_{\Prism}( \mathfrak{X}/A)\{-i\} \rightarrow
\mathrm{H}^{\ast}_{\Prism}( \mathbf{P}( \mathscr{E} ) / A).$$
In particular, the invertible $A$-module $A\{-1\}$ can be identified with the relative prismatic cohomology group
$\mathrm{H}^{2}_{\Prism}( \mathbf{P}^{1}_{\overline{A}} / A )$.
\end{variant}

\subsection{Higher Chern Classes}\label{subsection:higher-chern}

We now use projective bundle formula of Theorem \ref{theorem:projective-bundle-formula}
to introduce higher Chern classes in the setting of syntomic cohomology, using a formal procedure
outlined by Grothendieck in \cite{grothendieckchern}.

\begin{construction}[Chern Classes]\label{construction:higher-chern-classes}
Let $X$ be a scheme, formal scheme, or algebraic stack. Let $\mathscr{E}$ be a vector bundle
of rank $r$ on $X$, and let $\pi: \mathbf{P}( \mathscr{E} ) \rightarrow X$ denote the associated projective bundle.
By virtue of Theorem \ref{theorem:projective-bundle-formula}, the syntomic cohomology classes
$$\{ c_{1}^{\Syn}( \calO(1) )^{i} \in \mathrm{H}^{2i}_{\Syn}( \mathbf{P}(\mathscr{E}), \Z_p(i) ) \}_{0 \leq i < r}$$ determine an isomorphism of abelian groups
$$ \bigoplus_{0 \leq i < r } \mathrm{H}^{2r-2i}_{\Syn}(X, \Z_p(r-i) ) \rightarrow
\mathrm{H}^{2r}_{\Syn}( \mathbf{P}(\mathscr{E}), \Z_p(r) ).$$
It follows that there exist unique elements $\{ c_j^{\Syn}( \mathscr{E} ) \in \mathrm{H}^{2j}( X, \Z_p(j) ) \}_{0 \leq j \leq r}$
satisfying the identities
$$ c_0^{\Syn}( \mathscr{E} ) = 1 \quad \quad \sum_{0 \leq j \leq r} \pi^{\ast}( c^{\Syn}_{j}(\mathscr{E})) c^{\Syn}_1(\calO(1) )^{r-j} = 0.$$
We will refer to $c_{j}^{\Syn}( \mathscr{E} )$ as the {\it syntomic $j$th Chern class of $\mathscr{E}$}.
By convention, we set $c_{j}^{\Syn}( \mathscr{E} ) = 0$ if $j$ is an integer which does not belong to the set $\{ 0, 1, \cdots, r \}$.
\end{construction}

\begin{example}\label{example:first-chern-class}
Let $X$ be a scheme, formal scheme, or algebraic stack. Let $\mathscr{L}$ be a line bundle on $X$.
Then the projection map $\pi: \mathbf{P}( \mathscr{L} ) \rightarrow X$ is an isomorphism, and
the line bundle $\calO(1)$ is (by our convention) the pullback $\pi^{\ast}( \mathscr{L} )^{-1}$.
It follows that the syntomic first Chern class $c_{1}^{\Syn}( \mathscr{L} )$ of Construction \ref{construction:higher-chern-classes}
agrees with the syntomic first Chern class $c_{1}^{\Syn}( \mathscr{L} )$ introduced in Variant \ref{variant:globalization-of-Chern-class-scheme}.
\end{example}

\begin{variant}
In the situation of Construction \ref{construction:higher-chern-classes}, suppose that the vector bundle $\mathscr{E}$
is not assumed to have fixed rank. In that case, we can write $X$ as a disjoint union of closed and open subsets $\{ X(r) \}_{r \geq 0}$, where each restriction $\mathscr{E}|_{ X(r) }$ is a vector bundle of rank $r$. For every integer $i$, we let
$$ c_{i}^{\Syn}( \mathscr{E} ) \in \mathrm{H}^{2i}_{\Syn}(X, \Z_p(i) ) \simeq \prod_{r \geq 0} \mathrm{H}^{2i}_{\Syn}(X(r), \Z_p(i) )$$
denote the unique element whose restriction to each subset $X(r)$ is the syntomic Chern class
$c_i^{\Syn}( \mathscr{E}|_{ X(r) } )$ of Construction \ref{construction:higher-chern-classes}.
\end{variant}

\begin{remark}[Functoriality]\label{remark:functoriality-of-chern-class}
Let $f: X \rightarrow Y$ be a morphism of schemes, formal schemes, or algebraic stacks, and let
$\mathscr{E}$ be a vector bundle on $Y$. For every integer $i$, we have an equality
$$ c_{i}^{\Syn}( f^{\ast}(\mathscr{E} ) ) = f^{\ast} c_i^{\Syn}( \mathscr{E} )$$
in the abelian group $\mathrm{H}^{2i}_{\Syn}(X, \Z_p(i) )$.
\end{remark}

\begin{remark}[Comparison with \'{E}tale Cohomology]\label{remark:etale-comparison}
Let $X$ be a scheme, let $\mathscr{E}$ be a vector bundle on $X$, and let
let $U = \Spec( \Z[1/p] ) \times X$ be the open subscheme where $p$ is invertible. For every integer $i$,
the \'{e}tale comparison morphism
$$ \gamma_{\Syn}^{\mathet}\{i\}: \mathrm{H}^{2i}_{\Syn}(X, \Z_p(i) )
\rightarrow \mathrm{H}^{2i}_{\Syn}(U, \Z_p(i) ) \simeq \mathrm{H}^{2i}_{\mathet}(U, \Z_p(i) )$$
carries the syntomic Chern class $c_{i}^{\Syn}(\mathscr{E} )$ to the usual \'{e}tale Chern class
$c_{i}^{\mathet}( \mathscr{E}|_{U} )$. This follows from the compatibility of the \'{e}tale comparison morphism
with first Chern classes, which is immediate from the construction.
\end{remark}

\begin{remark}[Comparison with Crystalline Cohomology]
Let $k$ be a perfect field of characteristic $p$, let $X$ be a smooth $k$-scheme, and let $\mathscr{E}$ be a vector bundle on $X$. For every integer $i$, the crystalline comparison morphism
$$ \gamma_{\Syn}^{\crys}: \mathrm{H}^{2i}_{\Syn}(X, \Z_p(i))
\rightarrow \mathrm{H}^{2i}_{\Prism}(X)\{i\} \rightarrow \mathrm{H}^{2i}_{\crys}(X / \Z_p)$$
carries the syntomic first Chern class $c_i^{\Syn}( \mathscr{E} )$ to the usual $i$th Chern class
$c_i^{\crys}( \mathscr{E} )$ in crystalline cohomology. As in Remark \ref{remark:etale-comparison},
this follows immediately from the analogous statement for line bundles, which is a special case of Proposition \ref{proposition:crystalline-agreement}.
\end{remark}

\begin{theorem}[Additivity Formula]\label{theorem:additivity-chern}
Let $X$ be a scheme, formal scheme, or algebraic stack. Suppose we are given a short exact sequence
$$ 0 \rightarrow \mathscr{E}' \rightarrow \mathscr{E} \rightarrow \mathscr{E}'' \rightarrow 0$$
of vector bundles on $X$. Then, for every integer $n$, we have an equality
$$ c_{n}^{\Syn}( \mathscr{E} ) = \sum_{i+j = n} c_{i}^{\Syn}( \mathscr{E}') \cdot c_j^{\Syn}( \mathscr{E}'' )$$
in the abelian group $\mathrm{H}^{2n}_{\Syn}( X, \Z_p(n) )$.
\end{theorem}

We begin by proving Theorem \ref{theorem:additivity-chern} in the special case of a split extension.

\begin{lemma}\label{lemma:rore}
Let $X$ be a scheme, formal scheme, or algebraic stack, and let $\mathscr{E}'$ and $\mathscr{E}''$ be
vector bundles on $X$. Then, for every integer $n$, we have an equality 
$$ c_{n}^{\Syn}( \mathscr{E}' \oplus \mathscr{E}'' ) = \sum_{i+j = n} c_{i}^{\Syn}( \mathscr{E}') \cdot c_j^{\Syn}( \mathscr{E}'' )$$
in the abelian group $\mathrm{H}^{2n}_{\Syn}( X, \Z_p(n) )$.
\end{lemma}

\begin{proof}
Without loss of generality, we may assume that the vector bundles $\mathscr{E}'$ and $\mathscr{E}''$ have
fixed ranks $r'$ and $r''$. Set $\mathscr{E} = \mathscr{E}' \oplus \mathscr{E}''$, let
$\pi: \mathbf{P}( \mathscr{E} ) \rightarrow X$ denote the projection map, and let $\calO(-1)$
be the tautological subbundle of $\pi^{\ast} \mathscr{E}$. Let $U \subseteq \mathbf{P}( \mathscr{E} )$
be the open subset for which composite map
$$ \calO(-1) \hookrightarrow \calE' \oplus \calE'' \twoheadrightarrow \calE'$$
is the inclusion of a subbundle. It follows that sum
$$ \eta' = \sum_{0 \leq i \leq r'} \pi^{\ast}( c_{i}^{\Syn}( \mathscr{E}') ) c_1^{\Syn}( \calO(1) )^{i}$$
vanishes when restricted to $U$, and can therefore be promoted to a class in the relative
syntomic cohomology group $\mathrm{H}^{2r'}_{\Syn}( \mathbf{P}(\mathscr{E}), U, \Z_p(r') )$.
Let $V \subseteq \mathbf{P}( \mathscr{E} )$ be the open subset on which the composite map
$$ \calO(-1) \hookrightarrow \calE' \oplus \calE'' \twoheadrightarrow \calE''$$
is the inclusion of a subbundle; the same argument shows that the sum
$$ \eta'' = \sum_{0 \leq j \leq r''} \pi^{\ast}( c_{j}^{\Syn}( \mathscr{E}'') ) c_1^{\Syn}( \calO(1) )^{j}$$
can be promoted to an element of $\mathrm{H}^{2r''}_{\Syn}( \mathbf{P}(\mathscr{E}), V, \Z_p(r') )$.
It follows that the product class $\eta' \cdot \eta''$ can be promoted to an element of the relative syntomic cohomology group
$$ \mathrm{H}^{2r'+2r''}_{\Syn}( \mathbf{P}(\mathscr{E}), U \cup V, \Z_p(r'+r'') ),$$
which is trivial (since $U$ and $V$ are an open covering of $\mathbf{P}( \mathscr{E} )$). We therefore
obtain the identity
\begin{eqnarray*}
0 & = & \eta' \cdot \eta'' \\
& = & \sum_{n} \sum_{i+j = n} \pi^{\ast}( c_i^{\Syn}(\mathscr{E}') \cdot c_j^{\Syn}(\mathscr{E}'') ) c_{1}^{\Syn}( \calO(1) )^{n},
\end{eqnarray*}
in the cohomology group $\mathrm{H}^{2r'+2r''}_{\Syn}( \mathbf{P}(\mathscr{E}), \Z_p(r'+r'') )$, which is a restatement
of Lemma \ref{lemma:rore}.
\end{proof}

To reduce Theorem \ref{theorem:additivity-chern} to Lemma \ref{lemma:rore}, we will use the following general principle:

\begin{proposition}\label{proposition:iso-on-syntomic}
Let $f: X \rightarrow Y$ be a morphism of algebraic stacks satisfying the following assumptions:
\begin{itemize}
\item[$(1)$] The algebraic stacks $X$ and $Y$ are quasi-compact, quasi-separated, and $p$-quasisyntomic.
\item[$(2)$] For every morphism $\Spec(R) \rightarrow Y$, where $R$ is a strictly Henselian ring of
residue chararacteristic $\neq p$, the unit map
$$ \F_p \rightarrow \RGamma_{\mathet}( \Spec(R) \times_{Y} X, \F_p)$$
is an isomorphism in $\widehat{\calD}(\Z_p)$.
\item[$(3)$] For every integer $d$, the restriction map
$$ f^{\ast}: \RGamma(Y, L \widehat{\Omega}^{d}_{Y} ) \rightarrow \RGamma(X, L \widehat{\Omega}^{d}_{X} )$$
is an isomorphism in $\widehat{\calD}(\Z_p)$.
\end{itemize}
Then, for every integer $n$, the restriction map $f^{\ast}:  \RGamma_{\Syn}( Y, \Z_p(n) ) \rightarrow \RGamma_{\Syn}(X,\Z_p(n) )$
is an isomorphism in $\widehat{\calD}(\Z_p)$.
\end{proposition}

\begin{proof}
We will prove more generally that for every $p$-quasisyntomic commutative ring $A$ and every integer $n$, the restriction map
$$ \RGamma_{\Syn}( \Spec(A) \times Y, \Z_p(n) ) \rightarrow \RGamma_{\Syn}(\Spec(A) \times X,\Z_p(n) )$$
is an isomorphism in $\widehat{\calD}(\Z_p)$. By virtue of Proposition \ref{rhoax}, we may assume without loss of
generality that $A$ is an algebra over the cyclotomic number ring $\Z[ \zeta_{p^{\infty}}]$ of \S\ref{subsection:localization-generic-fiber}. Replacing $f$ by the map $(\id \times f): \Spec(A) \times X \rightarrow \Spec(A) \times Y$, we are reduced to verifying
Proposition \ref{proposition:iso-on-syntomic} in the special case where $Y$ is defined over the ring $\Z[ \zeta_{p^{\infty}} ]$.
Let $\epsilon$ denote the element $(1, \zeta_p, \zeta_{p^2}, \cdots ) \in T_p( \Z[\zeta_{p^{\infty}}]^{\times} )$, 
which we identify with its image in $\mathrm{H}^{0}_{\Syn}( Y, \Z_p(1) )$.

For every integer $n$, let $\RGamma_{\Syn}( \mathfrak{Y}, \Z_p(n) )$
denote the fiber of the map $$(\varphi\{n\} - 1): \Fil^{n}_{\Nyg} \RGamma_{\Prism}(Y)\{n\} \rightarrow \RGamma_{\Prism}(Y)\{n\},$$
and define $\RGamma_{\Syn}( \mathfrak{X}, \Z_p(n) )$ similarly. Using assumption $(1)$ and Variant \ref{variant:p-quasi-stack}, we see that
the restriction map
$\RGamma_{\Syn}( Y, \Z_p(n) ) \rightarrow \RGamma_{\Syn}(X,\Z_p(n) )$ can be obtained from a commutative diagram
$$ \xymatrix@R=50pt@C=30pt{ \RGamma_{\Syn}(\mathfrak{Y}, \Z_p(n) ) \ar[r] \ar[d] &
\varinjlim_{d} \RGamma_{\Syn}(\mathfrak{Y}, \Z_p(n+d)) \ar[d] & \RGamma_{\mathet}( \Spec(\Z[1/p]) \times Y, \Z_p(n)) \ar[l] \ar[d] \\
\RGamma_{\Syn}(\mathfrak{X}, \Z_p(n) ) \ar[r] & \varinjlim_{d} \RGamma_{\Syn}(\mathfrak{X}, \Z_p(n+d))
& \RGamma_{\mathet}( \Spec(\Z[1/p]) \times X, \Z_p(n)) \ar[l] }$$
by passing to inverse limits along the rows (where the colimits are $p$-completed). It
will therefore suffice to prove the following pair of assertions:
\begin{itemize}
\item[$(2')$] For every integer $n$, the restriction map 
$$ \RGamma_{\mathet}(\Spec(\Z[1/p]) \times Y, \Z_p(n) ) \rightarrow \RGamma_{\mathet}( \Spec(\Z[1/p]) \times X, \Z_p(n))$$
is an isomorphism in $\widehat{\calD}(\Z_p)$.

\item[$(3')$] For every integer $n$, the restriction map $\RGamma_{\Syn}( \mathfrak{Y}, \Z_p(n) )
\rightarrow \RGamma_{\Syn}( \mathfrak{X}, \Z_p(n) )$ is an isomorphism in $\widehat{\calD}(\Z_p)$.
\end{itemize}

We first prove $(2')$. In fact, we prove a stronger assertion: for every morphism of algebraic stacks $Y' \rightarrow Y$
where $Y'$ is defined over $\Z[1/p]$, the induced map $\RGamma_{\mathet}(Y', \Z_p(n) ) \rightarrow \RGamma_{\mathet}( Y' \times_{Y} X, \Z_p(n))$ is an isomorphism. To prove this, we may assume without loss of generality that $Y' = \Spec(R)$
is affine. Set $X' = Y' \times_{Y} X$ and let $f': X' \rightarrow Y'$ be the projection map. It will then
suffice to show that the unit map $\underline{\F_p}_{Y'} \rightarrow R f'_{\ast} \underline{\F_p}_{X'}$ is an isomorphism
in the derived category of \'{e}tale sheaves on $\Spec(Y')$. This can be checked stalkwise, and assumption
$(1)$ guarantees that the formation of the sheaf $R f'_{\ast} \underline{\F_p}_{X'}$ is compatible with filtered colimits
in $R$. We can therefore assume without loss of generality that $R$ is strictly Henselian, in which case the
desired result reduces to assumption $(2)$. 

To deduce $(3')$, we proceed in several steps:
\begin{itemize}
\item[$(a)$] For every integer $d$, the restriction map
$$ \RGamma(Y, \Fil_{d}^{\conj} \widehat{\Omega}^{\DHod}_{Y} ) \rightarrow \RGamma(Y, \Fil_{d}^{\conj} \widehat{\Omega}^{\DHod}_{X} )$$
is an isomorphism. This follows from assumption $(3)$ using induction on $d$.

\item[$(b)$] For every pair of integers $m$ and $n$, the restriction map
$$ \gr^{m}_{\Nyg} \RGamma_{\Prism}(Y)\{n\} \rightarrow \gr^{m}_{\Nyg} \RGamma_{\Prism}(X)\{n\}$$
is an isomorphism. This follows by combining $(a)$ with the fiber sequence of Remark \ref{remark:Nygaard-associated-graded}.

\item[$(c)$] For every triple of integers $m$, $m'$, and $n$ with $m \leq m'$, the restriction map
$$ \Fil^{m}_{\Nyg} \RGamma_{\Prism}(Y)\{n\} / \Fil^{m'}_{\Nyg} \RGamma_{\Prism}(Y)\{n\}
\rightarrow \Fil^{m}_{\Nyg} \RGamma_{\Prism}(X)\{n\} / \Fil^{m'}_{\Nyg} \RGamma_{\Prism}(X)\{n\}$$
is an isomorphism. This follows from $(b)$, using induction on the difference $m' - m$.

\item[$(d)$] For every pair of integers $m$ and $n$, the restriction map
$$\Fil^{m}_{\Nyg} \RGamma_{\widehat{\Prism} }(Y)\{n\} \rightarrow \Fil^{m}_{\Nyg} \RGamma_{\widehat{\Prism} }(X)\{n\}$$
is an isomorphism. This follows from $(c)$ by passing to the limit over $m'$.
\end{itemize}
To prove $(3')$, we observe that Proposition \ref{proposition:syntomic-Nygaard-complete} supplies
a commutative diagram of fiber sequences
$$ \xymatrix@R=50pt@C=50pt{ \RGamma_{\Syn}( \mathfrak{Y}, \Z_p(n) ) \ar[r] \ar[d] & \Fil^{n}_{\Nyg} \RGamma_{\widehat{\Prism} }(Y)\{n\} \ar[r]^-{\varphi\{n\}-\id} \ar[d] & \RGamma_{\widehat{\Prism}}(Y)\{n\} \ar[d] \\
\RGamma_{\Syn}( \mathfrak{X}, \Z_p(n) ) \ar[r] & \Fil^{n}_{\Nyg} \RGamma_{\widehat{\Prism} }(X)\{n\} \ar[r]^-{\varphi\{n\}-\id} & \RGamma_{\widehat{\Prism}}(X)\{n\}. }$$
Assertion $(d)$ guarantees that the middle and right vertical maps are isomorphisms, so that the left
vertical map is also an isomorphism.
\end{proof}

\begin{corollary}\label{corollary:restriction-to-Levi}
Let $R$ be a $p$-quasisyntomic commutative ring, let $G$ be a reductive algebraic group over $R$,
let $P \subseteq G$ be a parabolic subgroup and let $L \subseteq P$ be a Levi splitting. Then
the induced map of classifying stacks $\BL \rightarrow \BP$ induces an isomorphism
$$ \theta: \RGamma_{\Syn}( \BP, \Z_p(n) ) \rightarrow \RGamma_{\Syn}( \BL, \Z_p(n) )$$
for every integer $n$.
\end{corollary}

\begin{proof}
The assertion is local with respect to the \'{e}tale topology on $\Spec(R)$. We may therefore assume without loss
of generality that there is a short exact sequence 
$$ 0 \rightarrow \mathbf{G}_m \xrightarrow{\lambda} L \twoheadrightarrow L' \rightarrow 0,$$
where the cocharacter $\lambda$ is positive on the roots of $P$ which do not belong to $L$.

For each integer $k \geq 0$, let $Z_{k}$ denote the $(k+1)$-fold fiber power of $\BL$ over $\BP$.
Then $Z_{\bullet}$ forms a simplicial algebraic stack, and we can identify $\theta$ with the restriction map
$$ \Tot( \RGamma_{\Syn}( Z_{\bullet}, \Z_p(n) ) ) \rightarrow \RGamma_{\Syn}( Z_0, \Z_p(n) ) = \RGamma_{\Syn}( \BL, \Z_p(n) ).$$
To show that this map is an isomorphism, it will suffice to show that each of the projection maps $\pi: Z_{k} \rightarrow \BL$
induces an isomorphism on syntomic cohomology.

Let us henceforth regard $k$ as fixed. We will prove more generally that if $Y$ is a $p$-quasisyntomic algebraic stack
equipped with a map $Y \rightarrow \BL'$, then the projection map 
$$Y \times_{\BL'} Z_k \rightarrow Y \times_{ \BL' } \BL$$ induces an isomorphism on syntomic cohomology. The assertion is local on $Y$. We may assume without loss of generality that $Y = \Spec(R')$ and that the map $Y \rightarrow \BL'$ classifies the trivial $L'$-torsor on $Y$. Replacing $R$ by $R'$, we are reduced to proving that the projection map
$\pi_{\lambda}: \BGm \times_{\BL} Z_k \rightarrow \BGm$ induces an isomorphism on syntomic cohomology.

Let $U \subseteq P$ be the unipotent radical of $P$ and let $\calO_{U}$ denote its coordinate ring, which
we regard as a smooth $R$-algebra. The group $\mathbf{G}_m$ acts on $U$ via the cocharacter $\lambda$, which determines a grading $\calO_{U} = \bigoplus_{d \in \Z} \calO_{U,d}$. Since $\lambda$ is positive on the roots of $P$ which do not belong to $L$, the $R$-module $\calO_{U,d}$ vanishes for $d< 0$, and the unit map $R \rightarrow \calO_{U,0}$ is an isomorphism.
Unwinding the definitions, we see that $\BGm \times_{\BL} Z_k$ can be identified with the stack-theoretic quotient
$[\Spec(\calO^{\otimes k}_{U}) / \mathbf{G}_m]$. Applying Proposition \ref{proposition:graded-Hodge-cohomology},
we see that the projection map $\pi_{\lambda}$ induces an isomorphism on Hodge cohomology.
Moreover, if we are given any map $\Spec(A) \rightarrow \BL$, then
the projection map $\Spec(A) \times_{ \BL } Z_k \rightarrow \Spec(A)$ is an affine space bundle over $\Spec(A)$,
and therefore induces an isomorphism $\RGamma_{\mathet}( \Spec(A), \F_p )
\rightarrow \RGamma_{\mathet}( \Spec(A) \times_{\BL} Z_k, \F_p)$ if $p$ is invertible in $A$.
Applying the criterion of Proposition \ref{proposition:iso-on-syntomic}, we conclude that
$\pi_{\lambda}$ is an isomorphism on syntomic cohomology, as desired.
\end{proof}

\begin{proof}[Proof of Theorem \ref{theorem:additivity-chern}]
Suppose we are given a short exact sequence
\begin{equation}\label{equation:ses-for-additivity}
0 \rightarrow \mathscr{E}' \rightarrow \mathscr{E} \rightarrow \mathscr{E}'' \rightarrow 0
\end{equation}
of vector bundles on $X$. We wish to prove that, for every integer $n$, we have an equality
$$ c_{n}^{\Syn}( \mathscr{E} ) = \sum_{i+j = n} c_{i}^{\Syn}( \mathscr{E}') \cdot c_j^{\Syn}( \mathscr{E}'' )$$
in the syntomic cohomology group $\mathrm{H}^{2n}_{\Syn}( X, \Z_p(n) )$. Without loss of generality, we may assume that $\mathscr{E}'$ and $\mathscr{E}''$ are vector bundles of fixed ranks $r'$ and $r''$, respectively. 
Set $r = r' + r''$, and let $P \subseteq \GL_{r}$ be the standard parabolic subgroup with Levi factor
$\GL_{r'} \times \GL_{r''} \subseteq P$, which we regard as an affine group scheme over $\Z$.
The short exact sequence (\ref{equation:ses-for-additivity}) is classified by a map $u: X \rightarrow \BP$.
It will therefore suffice to treat the universal case where $X = \BP$ and $u$ is the identity map.
Corollary \ref{corollary:restriction-to-Levi} guarantees that the restriction map
$$ \mathrm{H}^{\ast}_{\Syn}( \BP, \Z_p(n) ) \rightarrow \mathrm{H}^{\ast}_{\Syn}( \BGL_{r'} \times \BGL_{r''}, \Z_p(n))$$
is an isomorphism, it suffices to prove the relevant identity in the syntomic cohomology group
$\mathrm{H}^{2n}_{\Syn}( \BGL_{r'} \times \BGL_{r''}, \Z_p(n))$. We are therefore reduced
to proving Theorem \ref{theorem:additivity-chern} in the special case where the sequence
(\ref{equation:ses-for-additivity}) splits, which follows from Lemma \ref{lemma:rore}.
\end{proof}

\subsection{The Syntomic Cohomology of \texorpdfstring{$\BGL_m$}{BGLm}}\label{subsection:cohomology-classifying-stack}

Fix an integer $m \geq 0$, and let $\BGL_m$ denote the classifying stack of the group scheme $\GL_m$;
we regard $\BGL_m$ as an algebraic stack which is smooth over $\Spec(\Z)$.
Let $\mathscr{E}_{\mathrm{univ}}$ denote the vector bundle on $\BGL_m$ corresponding to the tautological representation of $\GL_m$. Applying  Construction \ref{construction:higher-chern-classes}, we obtain syntomic Chern classes
$$ \{ c_{i}^{\Syn}( \mathscr{E}_{\mathrm{univ} } ) \in \mathrm{H}^{2i}_{\Syn}( \BGL_m, \Z_p(i) ) \}_{1 \leq i \leq m}.$$
Our goal in this section is to prove the following:

\begin{theorem}\label{theorem:syntomic-cohomology-of-BGL}
Let $X$ be a $p$-quasisyntomic algebraic stack which is quasi-compact and quasi-separated. Then the syntomic cohomology ring
$$ \bigoplus_{n \in \Z} \mathrm{H}^{\ast}_{\Syn}( \BGL_m \times X, \Z_p(n) )$$
is a polynomial algebra on generators $\{ c_{i}^{\Syn}( \mathscr{E}_{\mathrm{univ} } ) \}_{1 \leq i \leq m}$
over the ring $\bigoplus_{n \in \Z} \mathrm{H}^{\ast}_{\Syn}( X, \Z_p(n) )$.
\end{theorem}

We begin by treating the special case $m=1$.

\begin{lemma}\label{lemma:isovar}
Let $X$ be an algebraic stack which is $p$-quasisyntomic, and let $\mathscr{L}$ denote
the tautological line bundle on the product $\BGm \times X$. Then the syntomic cohomology ring
$$ \bigoplus_{n \in \Z} \mathrm{H}^{\ast}_{\Syn}( \BGm \times X, \Z_p(n) )$$
is a polynomial algebra over $\bigoplus_{n \in \Z} \mathrm{H}^{\ast}_{\Syn}(X, \Z_p(n) )$ on the generator
$c_{1}^{\Syn}( \mathscr{L} )$.
\end{lemma}

\begin{proof}
Let $\mathscr{E}$ denote the direct sum $\mathscr{L} \oplus \calO$,
which we regard as a vector bundle of rank $2$ over the product $\BG_m \times X$. 
Let $\pi: \mathbf{P}(\mathscr{E} ) \rightarrow \BGm \times X$ denote the associated projective bundle
and let $\calO(-1) \subseteq \pi^{\ast} \mathscr{E}$ be the tautological subbundle. Let
$U \subseteq \mathbf{P}( \mathscr{E} )$ be the open subset over which the composite map
$$ \calO(-1) \hookrightarrow \pi^{\ast} \mathscr{E} \twoheadrightarrow \pi^{\ast} \mathscr{L}$$
is an isomorphism, and let $V \subseteq \mathbf{P}( \mathscr{E} )$ be the open subset over which
the composite map 
$$ \calO(-1) \hookrightarrow \pi^{\ast} \mathscr{E} \twoheadrightarrow \calO_{ \mathbf{P}(\mathscr{E} )}$$
is an isomorphism. More concretely, $\mathscr{E}$ can be described as the product of $X$ with the
stack-theoretic quotient $[\mathbf{P}^{1} / \mathbf{G}_m]$, and $U$ and $V$ with the products of $X$
with the stack-theoretic quotients $[\mathbf{A}^{1}_{+} / \mathbf{G}_m]$ and $[\mathbf{A}^{1}_{-} / \mathbf{G}_m]$,
where $\mathbf{A}^{1}_{+}$ and $\mathbf{A}^{1}_{-}$ comprise the usual open covering of the projective line $\mathbf{P}^{1}$
by affine spaces. In particular, the composite map
$$ U \cap V \hookrightarrow \mathbf{P}(\mathscr{E} ) \xrightarrow{\pi} \BGm \times X \rightarrow X$$
is an isomorphism.

By construction, the line bundle $\calO(-1)$ is equipped with isomorphisms
$$ \alpha: \calO(-1)|_{U} \simeq (\pi^{\ast} \mathscr{L})|_{U} \quad \quad
\beta: \calO(-1)|_{V} \simeq \calO_{V}.$$
Composing $\alpha$ with the inverse of $\beta$, we obtain a trivialization $\gamma$ of the line bundle
$(\pi^{\ast} \mathscr{L})|_{U \cap V}$. The pair $( \mathscr{L}, \gamma)$ determines a relative Chern class
$$c_{1}^{\Syn}( \mathscr{L}, \gamma) \in \mathrm{H}^{2}_{\Syn}( \BGm \times X, U \cap V, \Z_p(1))$$
(see Variant \ref{variant:relative-first-Chern-class}). We will prove the following:
\begin{itemize}
\item[$(\ast)$] For every integer $n$, multiplication by the element $c_{1}^{\Syn}( \mathscr{L}, \gamma)$ induces an isomorphism
$$ \RGamma_{\Syn}( \BGm \times X, \Z_p(n-1) )[-2] \rightarrow \RGamma_{\Syn}( \BGm \times X, U \cap V, \Z_p(n) ).$$
\end{itemize}
Assume $(\ast)$ for the moment. Using the identification $U \cap V \simeq X$, we can reformulate $(\ast)$ as the assertion that for every pair of integers $m$ and $n$, the tautological map
$$ \mathrm{H}^{m}_{\Syn}( X, \Z_p(n) ) \oplus \mathrm{H}^{m-2}_{\Syn}( \BG_m \times X, \Z_p(n-1) )
\rightarrow \mathrm{H}^{m}_{\Syn}(\BGm \times X, \Z_p(n) )$$
is an isomorphism (where the map on the second factor is given by multiplication with $c_{1}^{\Syn}( \mathscr{L} )$.
Lemma \ref{lemma:isovar} now follows by induction on $m$, using the fact that the groups
$\mathrm{H}^{m}_{\Syn}(X, \Z_p(n) )$ and $\mathrm{H}^{m}_{\Syn}( \BGm \times X, \Z_p(n) )$ vanish for $m < 0$
(by virtue of our assumption that $X$ is $p$-quasisyntomic).

It remains to prove $(\ast)$. The assertion is local on $X$, so we may assume without loss of generality that
$X = \Spec(R)$ for some $p$-quasisyntomic ring $R$. In this case, we can identify $U$ with the
stack-theoretic quotient $[\Spec(R[t]) / \mathbf{G}_m]$ (where $R[t]$ is equipped with the usual grading).
Applying Propositions \ref{proposition:graded-Hodge-cohomology} and \ref{proposition:iso-on-syntomic},
we see that the composite map $U \hookrightarrow \mathbf{P}(\mathscr{E} ) \xrightarrow{\pi} \BGm \times X$
induces an isomorphism on syntomic cohomology. We can therefore reformulate $(\ast)$ as follows:
\begin{itemize}
\item[$(\ast')$] For every integer $n$, multiplication by the relative Chern class $c_{1}( \calO(-1), \beta )$ induces an isomorphism
$$ \RGamma_{\Syn}( \BGm \times X, \Z_p(n-1))[-2] \rightarrow \RGamma_{\Syn}( U, U \cap V, \Z_p(n) ).$$
\end{itemize}
Since $U$ and $V$ are an open covering of the projective bundle $\mathbf{P}( \mathscr{E} )$, this is equivalent
to the following:
\begin{itemize}
\item[$(\ast'')$] For every integer $n$, multiplication by the relative Chern class $c_{1}( \calO(-1), \beta )$ induces an isomorphism
$$ \RGamma_{\Syn}( \BGm \times X, \Z_p(n-1))[-2] \rightarrow \RGamma_{\Syn}( \mathbf{P}^{1}(\mathscr{E}), V, \Z_p(n) ).$$
\end{itemize}
Arguing as above, we see that the composite map $V \hookrightarrow \mathbf{P}^{1}(\mathscr{E} ) \xrightarrow{\pi} \BGm \times X$
is an isomorphism on syntomic cohomology. It follows that the natural map
$$ \RGamma_{\Syn}( \mathbf{P}^{1}(\mathscr{E}), V, \Z_p(n) ) \oplus
\RGamma_{\Syn}( \BGm \times X, \Z_p(n) ) \rightarrow \RGamma_{\Syn}( \mathbf{P}^{1}(\mathscr{E}), \Z_p(n) )$$
is an isomorphism, so that $(\ast'')$ can be reformulated as follows:
\begin{itemize}
\item[$(\ast''')$] For every integer $n$, the map
$$ \RGamma_{\Syn}( \BGm \times X, \Z_p(n) ) \oplus \RGamma_{\Syn}( \BGm \times X, \Z_p(n-1))[-2]
\xrightarrow{(1, c_{1}^{\Syn}(\calO(-)) )} \RGamma_{\Syn}( \mathbf{P}(\mathscr{E} ), \Z_p(n-1))[-2]$$
is an isomorphism.
\end{itemize}
This is a special case of the projective bundle formula (Theorem \ref{theorem:projective-bundle-formula}).
\end{proof}

\begin{proof}[Proof of Theorem \ref{theorem:syntomic-cohomology-of-BGL}]
Let $X$ be an algebraic stack which is quasi-compact, quasi-separated, and $p$-quasisyntomic, let
$A$ denote the syntomic cohomology ring $\bigoplus_{n \in \Z} \mathrm{H}^{\ast}_{\Syn}( X, \Z_p(n) )$,
and let $B$ denote the syntomic cohomology ring $\bigoplus_{n \in \Z} \mathrm{H}^{\ast}_{\Syn}( \BGL_m \times X, \Z_p(n) )$.
We have a canonical map of bigraded rings $\rho: A[ c_1, \cdots, c_m ] \rightarrow B$, carrying each of the formal
variable $c_{i}$ to the syntomic Chern class  $c_{i}^{\Syn}( \mathscr{E}_{\mathrm{univ}} ) \in \mathrm{H}^{2i}_{\Syn}(\BGL_m \times X, \Z_p(i) )$. We wish to show that $\rho$ is an isomorphism.

The proof proceeds by induction on $m$, the case $m=0$ being trivial. Assume that $m > 0$, and let
$\mathbf{P}( \mathscr{E}_{\mathrm{univ}} )$ denote the projectivization of the vector bundle $\mathscr{E}_{\mathrm{univ}}$,
so that we have a projection map $\pi: \mathbf{P}( \mathscr{E}_{\mathrm{univ}} ) \rightarrow \BGL_m \times X$
and a short exact sequence of vector bundles
$$ 0 \rightarrow \calO(-1) \rightarrow \pi^{\ast} \mathscr{E}_{\mathrm{univ} } \rightarrow \mathscr{E}' \rightarrow 0.$$
Let $C$ denote the syntomic cohomology ring $\bigoplus_{n \in \Z} \mathrm{H}^{\ast}_{\Syn}( \mathbf{P}( \mathscr{E}_{\mathrm{univ}}), \Z_p(n) )$. It follows from the definition of the syntomic Chern classes $c_{i}^{\Syn}( \mathscr{E}_{\mathrm{univ}} )$
that $\rho$ extends to a ring homomorphism
$$ \rho^{+}: A[ c_1, \cdots, c_m, t ] / ( t^{m} + c_1 t^{m-1} + \cdots +c_m ) \rightarrow C,$$
carrying the formal variable $t$ to the syntomic Chern class $c_{1}^{\Syn}( \calO(1) ) \in \mathrm{H}^{2}_{\Syn}( \mathbf{P}( \mathscr{E}_{\mathrm{univ} }), \Z_p(1) )$. Note that the domain of $\rho^{+}$ can be identified with the polynomial ring
$A[ c'_1, c'_2, \cdots, c'_{m-1}, t]$, where $c'_{i} = c_{i} + t c_{i-1} + \cdots + t^{i}$. Moreover,
it follows from Theorem \ref{theorem:additivity-chern} that the map $\rho^{+}$ carries each $c'_{i}$ to the
syntomic Chern class $c_{i}^{\Syn}( \mathscr{E}' ) \in \mathrm{H}^{2i}_{\Syn}( \mathbf{P}( \mathscr{E}_{\mathrm{univ}}, \Z_p(i) )$.

Note that the projective bundle $\mathbf{P}( \mathscr{E}_{\mathrm{univ}} )$ can be identified with the product
$\BP \times X$, where $P \subseteq \GL_m$ is the standard parabolic subgroup with Levi factor
$\GL_{m-1} \times \GL_{1} \subseteq \GL_{m}$. Combining our inductive hypothesis with Lemma \ref{lemma:isovar},
we deduce that the composite map
$$ A[ c'_1, \cdots, c'_{m-1}, t] \xrightarrow{\rho^{+}} C \rightarrow
\bigoplus_{n} \mathrm{H}^{\ast}_{\Syn}( \BGL_1 \times \BGL_{m-1} \times X, \Z_p(n) )$$
is an isomorphism of bigraded rings. Together with Corollary \ref{corollary:restriction-to-Levi}, this
guarantees that $\rho^{+}$ is an isomorphism of commutative rings.

Note that the quotient ring $A[ c_1, \cdots, c_m, t ] / ( t^{m} + c_1 t^{m-1} + \cdots +c_m )$
is a free module over $A[c_1,\cdots, c_m]$, with basis given by $\{ t^{i} \}_{0 \leq i < m}$.
Consequently, to show that $\rho$ is an isomorphism, it will suffice to show that the images of the elements
$\{ t^{i} \}_{0 \leq i < m}$ form a basis for $C$ as a $B$-module. This follows from the projective bundle
formula (Theorem \ref{theorem:projective-bundle-formula}).
\end{proof}

\subsection{The Blowup Formula for Syntomic Cohomology}\label{subsection:blowup-formula}

Throughout this section, we fix a scheme $X$ together with a closed subscheme $Y \subseteq X$ 
for which the inclusion $Y \hookrightarrow X$ is a regular immersion: that is, locally on $X$, the closed subscheme $Y$ can be realized as the
vanishing locus of a regular sequence of sections of the structure sheaf $\calO_{X}$. We write $\widetilde{X}$ for the blowup of $X$ along $Y$
and $D \subseteq \widetilde{X}$ for the exceptional divisor, so that we have a commutative diagram of schemes
\begin{equation}
\begin{gathered}\label{equation:blowup-diagram}
\xymatrix@R=50pt@C=50pt{ D \ar[r] \ar[d] & \widetilde{X} \ar[d] \\
Y \ar[r] & X }
\end{gathered}
\end{equation}

Our goal is to prove the following:
 
\begin{theorem}[Blowup Formula]\label{theorem:syntomic-blow-up}
For every integer $n$, the diagram (\ref{equation:blowup-diagram}) determines a pullback square of syntomic complexes
\begin{equation}
\begin{gathered}\label{equation:syntomic-complexes-blowup-formula}
\xymatrix@R=50pt@C=50pt{ \RGamma_{\Syn}(D, \Z_p(n) ) & \RGamma_{\Syn}( \widetilde{X}, \Z_p(n) ) \ar[l] \\
\RGamma_{\Syn}(Y, \Z_p(n) ) \ar[u] & \RGamma_{\Syn}(X, \Z_p(n) ) \ar[l] \ar[u] }
\end{gathered}
\end{equation}
in the $\infty$-category $\widehat{\calD}(\Z_p)$.
\end{theorem}

\begin{corollary}\label{corollary:blowup-application}
Let $f: X \rightarrow Y$ be a smooth morphism of schemes of relative dimension $r$, let $s: Y \rightarrow X$
be a section of $f$, and let $\widetilde{X}$ be the blowup of $X$ along the image of $s$, and let $D \subseteq \widetilde{X}$ be the exceptional divisor. Then, for every integer $n$, the syntomic cohomology classes $\{ c_{1}^{\Syn}( \calO(D) )^{i} \}_{0 < i < r}$ induce an isomorphism
$$ \RGamma_{\Syn}( X, \Z_p(n) )
\oplus \bigoplus_{0 < i < r} \RGamma_{\Syn}(Y, \Z_p(n-i) )[-2i] 
\rightarrow \RGamma_{\Syn}( \widetilde{X}, \Z_p(n) )$$
in the derived $\infty$-category $\widehat{\calD}(\Z_p)$.
\end{corollary}

\begin{proof}
Combine Theorem \ref{theorem:syntomic-blow-up} with Theorem \ref{theorem:projective-bundle-formula}.
\end{proof}

The proof of Theorem \ref{theorem:syntomic-blow-up} will proceed in several steps. We begin by proving the analogue of
Theorem \ref{theorem:syntomic-blow-up} for (derived) Hodge cohomology (see \cite{rao2019hodge} for a closely related statement):

\begin{lemma}\label{lemma:hodge-blowup}
The diagram (\ref{equation:blowup-diagram}) determines a pullback square of derived Hodge complexes
\begin{equation}
\begin{gathered}\label{equation:blowup-diagram-Hodge1}
\xymatrix@R=50pt@C=50pt{ \bigoplus_{d \in \Z} \RGamma(D, L \Omega^{d}_{D} )[-d] & \bigoplus_{d \in \Z} \RGamma( \widetilde{X}, L \Omega^{d}_{\widetilde{X}} )[-d] \ar[l] \\
\bigoplus_{d \in \Z} \RGamma(Y, L \Omega^{d}_{Y}  )[-d] \ar[u] & \bigoplus_{d \in \Z} \RGamma(X, L \Omega^{d}_{X} )[-d]. \ar[l] \ar[u] }
\end{gathered}
\end{equation}
\end{lemma}

\begin{proof}
The assertion is local on $X$. We may therefore assume without loss of generality that $X = \Spec(R)$ is affine
and that $Y$ is the vanishing locus of a regular sequence $f_1, f_2, \cdots, f_n \in R$, which determines a morphism $u$ from
$X$ to the affine space $\mathbf{A}^{n} = \Spec( \Z[x_1, \cdots, x_n] )$. Let $\widetilde{ \mathbf{A} }^{n}$ denote the blowup of
$\mathbf{A}^{n}$ at the origin, and let $D' \subseteq \widetilde{ \mathbf{A} }^{n}$ be the exceptional divisor. Note that pullback
along $u$ induces a map of derived Hodge complexes
$$ u^{\ast}: \bigoplus_{d \in \Z} \Omega^{d}_{ \Z[x_1, \cdots, x_n] }[-d] \rightarrow \bigoplus_{d \in \Z} L \Omega^{d}_{R}[-d].$$
Applying Remark \ref{remark:Kunneth-non-affine}, we see that the diagram (\ref{equation:blowup-diagram-Hodge1}) can be obtained from the simpler diagram
\begin{equation}
\begin{gathered}\label{equation:blowup-diagram-Hodge2}
\xymatrix@R=50pt@C=50pt{  \RGamma(D', \Omega^{\ast}_{D'} ) & \RGamma( \widetilde{ \mathbf{A} }^{n}, \Omega^{\ast}_{ \widetilde{ \mathbf{A} }^{n} } ) \ar[l] \\
\RGamma( \Spec(\Z), \Omega^{\ast}(\Spec(\Z)) ) \ar[u] & \RGamma( \mathbf{A}^{n}, \Omega^{\ast}_{ \mathbf{A}^{n} }). \ar[l] \ar[u] }
\end{gathered}
\end{equation}
via derived extension of scalars along $u$. It will therefore suffice to show that (\ref{equation:blowup-diagram-Hodge2}) is a pullback square.

Let us regard $\mathbf{A}^{n}$ as an open subscheme of the projective space $\mathbf{P}^{n}$ (over $\Z$), and let
$\widetilde{\mathbf{P}}^{n}$ denote the blowup of $\mathbf{P}^{n}$ at the origin. We can then expand (\ref{equation:blowup-diagram-Hodge2}) to a commutative
diagram
\begin{equation}
\begin{gathered}\label{equation-fat-diagram}
\xymatrix@R=50pt@C=50pt{ \RGamma(D', \Omega^{\ast}_{D'} ) & \RGamma( \widetilde{ \mathbf{A}}^{n}, \Omega^{\ast}_{ \widetilde{\mathbf{A}}^{n}} ) \ar[l] & \RGamma( \widetilde{\mathbf{P}}^{n}, \Omega^{d}_{ \widetilde{\mathbf{P}}^{n} }) \ar[l] \\
\RGamma( \Spec(\Z), \Omega^{\ast}_{\Spec(\Z)} ) \ar[u] & \RGamma(\mathbf{A}^{n}, \Omega^{\ast}_{ \mathbf{A}^{n} } ) \ar[l] \ar[u] & \RGamma( \mathbf{P}^{n}, \Omega^{\ast}_{\mathbf{P}^{n} }), \ar[u] \ar[l] }
\end{gathered}
\end{equation}
where the square on the right is a pullback (by Zariski descent). It will therefore suffice to show that the outer rectangle is also a pullback diagram.

We will assume that $n \geq 1$ (otherwise, there is nothing to prove). Let $\pi: \widetilde{\mathbf{P}}^{n} \rightarrow \mathbf{P}^{n}$ denote the projection map
and set $t = c_1^{\Hodge}( \pi^{\ast} \calO(1)) \in \mathrm{H}^{1}( \widetilde{\mathbf{P}}^{n}, \Omega^{1}_{\widetilde{\mathbf{P}}^{n}} )$.
Note that there is a canonical map $\pi': \widetilde{\mathbf{P}}^{n} \rightarrow \mathbf{P}^{n-1}$, which exhibits $\widetilde{\mathbf{P}}^{n}$
as the projectivization of the rank $2$ vector bundle $\calO(-1) \oplus \calO$. Set 
$t' = c_1^{\Hodge}( \pi'^{\ast} \calO(1) ) \in \mathrm{H}^{1}( \widetilde{\mathbf{P}}^{n}, \Omega^{1}_{\widetilde{\mathbf{P}}^{n}} )$.
Note that we have a subbundle inclusion $\pi^{\ast}( \calO(-1) ) \hookrightarrow \pi'^{\ast}( \calO(-1) \oplus \calO )$.
In particular, $\widetilde{\mathbf{P}}^{n}$ can be covered by open subschemes $U = \widetilde{\mathbf{A}}^{n}$ and $V$ for which the line bundle
$\pi^{\ast} \calO(-1)$ is trivial when restricted to $U$ and isomorphic to $\pi'^{\ast}( \calO(-1) )$ when restricted to $V$.
It follows that $t|_{U} = 0$ and $t|_{V} = t'|_{V}$, so that the product $t \cdot (t - t')$ vanishes: that is, we have
an equality $t^2 = tt'$ in the cohomology group $\mathrm{H}^{2}( \widetilde{\mathbf{P}}^{n}, \Omega^{2}_{ \widetilde{\mathbf{P}}^{n} } )$. 

Applying Lemma \ref{lemma:straight}, we deduce that the map
$$ \RGamma( \mathbf{P}^{n-1}, \Omega^{\ast}_{ \mathbf{P}^{n-1} }) \oplus
\RGamma( \mathbf{P}^{n-1}, \Omega^{\ast-1}_{ \mathbf{P}^{n-1} })[-1] \xrightarrow{(1,t)} \RGamma( \widetilde{\mathbf{P}}^{n}, \Omega^{\ast}_{ \widetilde{\mathbf{P}}^{n} })$$
is an isomorphism. Since $\pi'|_{D'}$ is an isomorphism, it follows that the restriction map 
$$ u: \mathrm{H}^{\ast}( \widetilde{\mathbf{P}}^{n}, \Omega^{\ast}_{ \widetilde{\mathbf{P}}^{n} }) \rightarrow
\mathrm{H}^{\ast}( D', \Omega^{\ast}_{ D'})$$
is a surjection, whose kernel can be identified (via multiplication by $t$) with the bigraded abelian group
$\mathrm{H}^{\ast-1}( \mathbf{P}^{n-1}, \Omega^{\ast-1}_{ \mathbf{P}^{n-1} })$. Applying Lemma \ref{lemma:classical-hodge-calculation} (to the projective space $\mathbf{P}^{n-1}$), we conclude that the kernel $\ker(u)$ is a free abelian group
generated by the cohomology classes $t \cdot t'^{d-1} = t^{d}$ for $0 < d \leq r$.
Applying Lemma \ref{lemma:classical-hodge-calculation} to the projective space $\mathbf{P}^{n}$, we conclude that the restriction map 
$$ v: \mathrm{H}^{\ast}( \mathbf{P}^{n}, \Omega^{\ast}_{\mathbf{P}^{n} }) \rightarrow
\mathrm{H}^{\ast}( \Spec(\Z), \Omega^{\ast}_{\Spec(\Z)})$$
is also a surjection, whose kernel is the free abelian group generated by $c_{1}( \calO(1) )^{d}$ if $0 < d \leq r$.
It follows that pullback along $\pi$ induces an isomorphism of abelian groups $\ker(v) \rightarrow \ker(u)$, so that the outer rectangle
(\ref{equation-fat-diagram}) is a pullback diagram, as desired.
\end{proof}

\begin{lemma}\label{lemma:diffracted-blowup}
For every integer $d$, the diagram (\ref{equation:blowup-diagram}) determines a pullback square
$$ \xymatrix@R=50pt@C=50pt{
 \RGamma(D, \Fil_{d}^{\conj} \Omega^{\DHod}_{D} ) & \RGamma( \widetilde{X}, \Fil_{d}^{\conj} \Omega^{\DHod}_{\widetilde{X}} ) \ar[l] \\
\RGamma(Y, \Fil_{d}^{\conj} \Omega^{\DHod}_{Y}  ) \ar[u] &  \RGamma(X, \Fil_{d}^{\conj} \Omega^{\DHod}_{X} ). \ar[l] \ar[u] }
$$
\end{lemma}

\begin{proof}
For $d < 0$, the assertion is vacuous. The general case follows by induction on $d$, using Lemma \ref{lemma:hodge-blowup}.
\end{proof}

\begin{lemma}\label{lemma:HT-blowup}
For every integer $n$, the diagram (\ref{equation:blowup-diagram}) determines a pullback square of absolute Hodge-Tate complexes
$$ \xymatrix@R=50pt@C=50pt{
 \RGamma_{\overline{\Prism}}(D)\{n\} & \RGamma_{\overline{\Prism}}( \widetilde{X})\{n\} \ar[l] \\
\RGamma_{\overline{\Prism}}(Y)\{n\} \ar[u] &  \RGamma_{\overline{\Prism}}(X)\{n\}. \ar[l] \ar[u] }
$$
\end{lemma}

\begin{proof}
The assertion is local on $X$, so we may assume without loss of generality that $X$ is quasi-compact and quasi-separated.
Applying Lemma \ref{lemma:diffracted-blowup} and passing to the colimit over $d$, we deduce that the diagram of
diffracted Hodge complexes
$$ \xymatrix@R=50pt@C=50pt{
 \RGamma(D, \Omega^{\DHod}_{D} ) & \RGamma( \widetilde{X}, \Omega^{\DHod}_{\widetilde{X}} ) \ar[l] \\
\RGamma(Y,\Omega^{\DHod}_{Y}  ) \ar[u] &  \RGamma(X, \Omega^{\DHod}_{X} ). \ar[l] \ar[u] }
$$
is a pullback square. Lemma \ref{lemma:HT-blowup} now follows by $p$-completing and 
passing to eigenspaces of the Sen operator (see Remark \ref{remark:diffracted-vs-prismatic2}).
\end{proof}

\begin{lemma}\label{lemma:prismatic-blowup}
For every pair of integers $m$ and $n$, the diagram (\ref{equation:blowup-diagram}) determines a pullback square
$$ \xymatrix@R=50pt@C=50pt{
\Fil^{m}_{\Nyg} \RGamma_{\Prism}(D)\{n\} & \Fil^{m}_{\Nyg} \RGamma_{\Prism}( \widetilde{X})\{n\} \ar[l] \\
\Fil^{m}_{\Nyg} \RGamma_{\overline{\Prism}}(Y)\{n\} \ar[u] & \Fil^{m}_{\Nyg} \RGamma_{\overline{\Prism}}(X)\{n\}. \ar[l] \ar[u] }
$$
\end{lemma}

\begin{proof}
For $m \leq 0$, this follows formally from Lemma \ref{lemma:HT-blowup} (see Remark \ref{remark:diffracted-vs-prismatic1}).
Proceeding by induction on $m$, we are reduced to showing that the diagram
$$ \xymatrix@R=50pt@C=50pt{
\gr^{m}_{\Nyg} \RGamma_{\Prism}(D)\{n\} & \gr^{m}_{\Nyg} \RGamma_{\Prism}( \widetilde{X})\{n\} \ar[l] \\
\gr^{m}_{\Nyg} \RGamma_{\overline{\Prism}}(Y)\{n\} \ar[u] & \gr^{m}_{\Nyg} \RGamma_{\overline{\Prism}}(X)\{n\}. \ar[l] \ar[u] }
$$
is a pullback square. This follows by combining Lemma \ref{lemma:diffracted-blowup} with the fiber sequence of Remark \ref{remark:Nygaard-associated-graded}.
\end{proof}

\begin{proof}[Proof of Theorem \ref{theorem:syntomic-blow-up}]
Fix an integer $n$. For every scheme $Z$, let $j_{Z}: \Spec(\Z[1/p]) \times Z \hookrightarrow Z$ be the inclusion map, so that
Remark \ref{remark:globalized-fiber-sequence} supplies a fiber sequence
$$ \RGamma_{\mathet}( Z, j_{Z!} \Z_p(n) ) \rightarrow \RGamma_{\Syn}( Z, \Z_p(n) )
\rightarrow \fib( \varphi\{n\} - \id: \Fil^{n}_{\Nyg} \RGamma_{\Prism}( Z)\{n\} \rightarrow \RGamma_{\Prism}(Z)\{n\}).$$
By virtue of Lemma \ref{lemma:prismatic-blowup}, to show that the diagram (\ref{equation:syntomic-complexes-blowup-formula}) is a pullback square,
it will suffice to show that the analogous diagram
\begin{equation}
\begin{gathered}\label{equation:etale-blowup-square}
\xymatrix@R=50pt@C=50pt{ \RGamma_{\mathet}(D, j_{D!} \F_p(n) ) & \RGamma_{\mathet}( \widetilde{X}, j_{\widetilde{X}!} \Z_p(n) ) \ar[l] \\
\RGamma_{\mathet}(Y, j_{Y!} \Z_p(n) ) \ar[u] & \RGamma_{\mathet}(X, j_{X!} \Z_p(n) ) \ar[l] \ar[u] }
\end{gathered}
\end{equation}
Let $\pi: \widetilde{X} \rightarrow X$ denote the projection map, and let $i: Y \hookrightarrow X$ denote the inclusion.
Unwinding the definitions, we see that after reducing modulo $p$, the diagram (\ref{equation:etale-blowup-square})
is obtained by applying the functor $\mathscr{F} \mapsto \RGamma_{\mathet}( X, \mathscr{F} )$ to a diagram
\begin{equation}
\begin{gathered}\label{equation:etale-blowup-square-local}
 \xymatrix@R=50pt@C=50pt{ (\pi|_{D})_{\ast} j_{D!} \underline{\F_p}(n) & \pi_{\ast} j_{\widetilde{X}!} \underline{\F_p}(n) \ar[l] \\
i_{\ast} j_{Y!} \underline{\F_p}(n) \ar[u] & j_{X!} \underline{\F_p}(n) \ar[u] \ar[l] }
\end{gathered}
\end{equation}
of complexes of \'{e}tale sheaves on $X$. It will therefore suffice to show that the diagram (\ref{equation:etale-blowup-square-local})
is a pullback square. This is clear: the horizontal maps in the diagram (\ref{equation:etale-blowup-square-local}) induce isomorphisms on
stalks at each geometric point of $Y$ (by definition for the bottom row, and by proper base change for the top row), and the vertical maps induce isomorphisms of stalks at each geometric point of the complement $X \setminus Y$ (as $\pi$ is an isomorphism over $U$).
\end{proof}

\renewcommand{\thesubsection}{\Alph{subsection}}

\appendix

\newpage \section*{Appendix}

\subsection{Animated Commutative Rings}
\label{ss:CAlgAnim}

In this appendix, we briefly review the formalism of animated commutative rings, which is used freely throughout this paper.

\begin{definition}\label{definition:animated-commutative-ring}
Let $k$ be a commutative ring. We let $\Poly_{k}$ denote the category whose objects are the polynomial rings $k[x_1, \cdots, x_n]$ and whose morphisms are 
$k$-algebra homomorphisms. Let $\SSet$ denote the $\infty$-category of spaces. An {\it animated commutative $k$-algebra} is a functor $$\Poly_{k}^{\op} \rightarrow \SSet$$ which preserves finite products (that is, it carries tensor products
of polynomial algebras to products of spaces). We let $\CAlg^{\anim}_{k}$ denote the full subcategory of $\Fun( \Poly_{k}^{\op}, \SSet)$ spanned by the animated commutative $k$-algebras. We will be primarily interested in the special case where $k = \Z$ is the ring of integers. In this case, we denote the $\infty$-category
$\CAlg^{\anim}_{k}$ by $\CAlg^{\anim}$, and refer to its objects as {\it animated commutative rings}. 
\end{definition}

\begin{example}[$k$-Algebras as Animated $k$-Algebras]
Let $R$ be a commutative $k$-algebra. Then the functor
$$(P \in \Poly_{k}^{\op} ) \mapsto \{ \textnormal{$k$-algebra homomorphisms $P \rightarrow R$} \}$$
is an animated commutative $k$-algebra. This construction determines a fully faithful embedding
from the ordinary category of commutative $k$-algebras to the $\infty$-category $\CAlg^{\anim}_{k}$ of animated commutative $k$-algebras.
The essential image of this embedding consists of those animated commutative $k$-algebras $X: \Poly_{k}^{\op} \rightarrow \SSet$ which are {\em discrete},
in the sense that the higher homotopy groups of the space $X(P)$ vanish for every polynomial algebra $P \in \Poly_{k}$.
We will systematically abuse notation by identifying a commutative $k$-algebra $R$ with its image in the $\infty$-category $\CAlg^{\anim}_{k}$.
\end{example}

\begin{notation}[The Underlying $k$-Algebra]\label{notation:pi-0-animated}
Let $R$ be an animated commutative $k$-algebra. Then the functor
$$ \Poly_{k}^{\op} \xrightarrow{ R } \SSet \xrightarrow{ \pi_0 } \Set$$
is also an animated commutative $k$-algebra, which we will denote by $\pi_0(R)$. By construction, $\pi_0(R)$ is discrete, and can therefore be viewed
as a commutative $k$-algebra in the usual sense. 
\end{notation}

\begin{remark}
\label{rmk:CAlgModel}
Let $\mathbf{A}$ denote the ordinary category of simplicial commutative $k$-algebras. Then $\mathbf{A}$ has the structure of a simplicial model category,
whose underlying $\infty$-category is canonically equivalent to the $\infty$-category $\CAlg^{\anim}_{k}$ of animated commutative $k$-algebras.
In other words, the $\infty$-category $\CAlg^{\anim}_{k}$ can be obtained from the ordinary category $\mathbf{A}$ by formally adjoining inverses to weak homotopy equivalences, or more explicitly as the homotopy coherent nerve of the category of fibrant-cofibrant objects of $\mathbf{A}$. See Corollary~5.5.9.3 of \cite{HTT}.
\end{remark}

The $\infty$-category $\CAlg^{\anim}_{k}$ can be characterized by a universal mapping property: it can be obtained from the ordinary category
of finitely generated polynomial algebras over $k$ by freely adjoining sifted colimits. More precisely, we have the following:

\begin{proposition}\label{proposition:universal-of-animated}
Let $k$ be a commutative ring and let $\calC$ be an $\infty$-category which admits small sifted colimits.
Then every functor $F: \Poly_{k} \rightarrow \calC$ admits a canonical extension $LF: \CAlg^{\anim}_{k} \rightarrow \calC$,
which is uniquely determined (up to isomorphism) by the requirement that it commutes with sifted colimits.
\end{proposition}

\begin{proof}
See Proposition~5.5.8.15 of \cite{HTT}.
\end{proof}

\begin{remark}[Derived Functors as Kan Extensions]
\label{rmk:LKECAlgAnim}
In the situation of Proposition \ref{proposition:universal-of-animated}, we will refer to $LF$ as the {\it nonabelian left
derived functor} of $F$. It is characterized (up to isomorphism) by the property that it is a left Kan extension of $F$.
In particular, if $\calD \subseteq \CAlg^{\anim}_{k}$ is a full subcategory containing $\Poly_{k}$ and $G: \calD \rightarrow \calC$
is any functor, then every natural transformation $\alpha_0: F \rightarrow G|_{ \Poly_{k} }$ admits an essentially unique extension to a natural transformation $\alpha: (LF)|_{\calD} \rightarrow G$.
\end{remark}

\begin{remark}
\label{AnimationColimits}
In the situation of Proposition \ref{proposition:universal-of-animated}, suppose that the $\infty$-category $\calC$ admits small colimits.
Then the functor $LF: \CAlg^{\anim}_{k} \rightarrow \calC$ preserves small colimits if and only if the original functor
$F: \Poly_{k} \rightarrow \calC$ preserves finite coproducts: that is, it carries tensor products of polynomial algebras
to coproducts in the $\infty$-category $\calC$. See Corollary~5.5.8.17 of \cite{HTT}.
\end{remark}

\subsection{Derived Hodge Cohomology}

We begin by reviewing a standard application of Proposition \ref{proposition:universal-of-animated}.

\begin{construction}[Exterior Powers of the Cotangent Complex]\label{construction:exterior-of-cotangent}
Let $A$ be a commutative ring. For every commutative $A$-algebra $B$, we let $\Omega^{1}_{B/A}$ denote the module of K\"{a}hler differentials of $B$
over $A$, and we let $\Omega^{d}_{B/A}$ denote its $d$th exterior power in the abelian category of $B$-modules; by convention, we have
$\Omega^{d}_{B/A} \simeq 0$ for $d < 0$. Applying Proposition \ref{proposition:universal-of-animated}, we see that there
is an essentially unique functor 
$$ \CAlg^{\anim}_{A} \rightarrow \calD(A) \quad \quad B \mapsto L\Omega^{d}_{B/A}$$
with the following properties:
\begin{itemize}
\item The functor $B \mapsto L \Omega^{d}_{B/A}$ commutes with sifted colimits.
\item When $B$ is a finitely generated polynomial algebra over $A$, there is canonical isomorphism $L \Omega^{n}_{B/A} \simeq \Omega^{n}_{B/A}$.
\end{itemize}
We will refer to $L \Omega^{1}_{B/A}$ as the {\it cotangent complex of $B$ over $A$} and to $L \Omega^{d}_{B/A}$ as the (derived) $d$th exterior power
of $L \Omega^{1}_{B/A}$. We let $L \widehat{\Omega}^{d}_{B/A}$ denote the $p$-completion of $L \Omega^{d}_{B/A}$.
\end{construction}

\begin{remark}
Let $A$ be a commutative ring. For every commutative $A$-algebra $B$, we can regard the direct sum
$$ \bigoplus_{d \in \Z} \Omega^{d}_{B/A}[-d]$$
as a commutative differential graded algebra over $A$, where the differential is identically zero. Applying Proposition \ref{proposition:universal-of-animated}
to the functor $B \mapsto  \bigoplus_{d \in \Z} \Omega^{d}_{B/A}[-d]$, we deduce that for every animated commutative $A$-algebra $B$,
the direct sum
$$ \bigoplus_{d \in \Z} L \Omega^{d}_{B/A}[-d]$$
can be regarded as a graded commutative algebra object of the $\infty$-category $\calD(A)$. In particular, each $L \Omega^{d}_{B/A}$ has
the structure of a module over $B \simeq L\Omega^{0}_{B/A}$, and can therefore be regarded as an object of the $\infty$-category $\calD(B)$.
\end{remark}

\begin{remark}\label{remark:kunneth-for-derived-Hodge}
Let $A$ be a commutative ring. Then the functor
$$ \CAlg^{\anim}_{A} \rightarrow \CAlg( \calD(A) ) \quad \quad B \mapsto \bigoplus_{d \in \Z} L\Omega^{d}_{B/A}[-d]$$
commutes with {\em all} small colimits (not just sifted colimits). This follows formally from the observation that the functor
$$ \Poly_{A} \rightarrow \CAlg( \calD(A) ) \quad \quad P \mapsto \bigoplus_{d \in \Z} \Omega^{d}_{P/A}[-d]$$
commutes with finite coproducts. 
\end{remark}

\begin{remark}
Let $f: A \rightarrow B$ be a morphism of commutative rings. For every integer $d$, there is a canonical map
$L \Omega^{d}_{B/A} \rightarrow \Omega^{d}_{B/A}$, which identifies $\Omega^{d}_{B/A}$ with the $0$th cohomology of
the complex $L \Omega^{d}_{B/A}$. Moreover, the cohomology groups of the complex $L \Omega^{d}_{B/A}$ vanish in degrees $> 0$.
\end{remark}

\begin{remark}\label{remark:LOmega-is-Omega}
Let $f: A \rightarrow B$ be a morphism of commutative rings. Suppose that the cotangent complex $L\Omega^{1}_{B/A}$ is a flat
$B$-module concentrated in cohomological degree zero (this condition is satisfied, for example, if $B$ is smooth over $A$). Then, for every integer $d$, the canonical map $L \Omega^{d}_{B/A} \rightarrow \Omega^{d}_{B/A}$ is an isomorphism (in the derived $\infty$-category $\calD(B)$).
\end{remark}

\begin{remark}\label{remark:descent-for-LOmega}
Let $A$ be a commutative ring. For every nonnegative integer $n$, the functor
$$ \CAlg_{A}^{\anim} \rightarrow \calD(A) \quad \quad B \mapsto L \Omega^{n}_{B/A}$$
satisfies descent for the flat topology. See \cite[Remark 2.8]{bhatt-completions} or \cite[Theorem 3.1]{BMS2}.
\end{remark}

\begin{construction}[The Derived Hodge Complex]
Let $A$ be a commutative ring and let $X$ be a scheme, formal scheme, or algebraic stack defined over $\Spec(A)$.
For every integer $d$, we let $\RGamma(X, L \Omega^{d}_{X/A} )$ denote the limit
$$ \varprojlim_{ \Spec(B) \rightarrow X} L \Omega^{d}_{B/A},$$
formed in the $\infty$-category $\calD(A)$, and we denote its $p$-completion by $\RGamma(X, L \widehat{\Omega}^{d}_{X/A} )$.
We will refer to the direct sum
$$ \bigoplus_{d \in \Z} \RGamma(X, L \Omega^{d}_{X/A} )[-d]$$
as the {\it derived Hodge complex of $X$ relative to $A$}, and to its cohomology as the {\it derived Hodge complex of $X$ relative to $A$}.

We will be primarily interested in the case where $A = \Z$; in this case, we denote $\RGamma(X, L \Omega^{d}_{X/A})$ by $\RGamma(X, L \Omega^{d}_{X})$
and its $p$-completion by $\RGamma(X, \widehat{\Omega}^{d}_{X} )$.
\end{construction}

\begin{remark}[K\"{u}nneth Formula]\label{remark:Kunneth-non-affine}
Let $f: A \rightarrow B$ be a morphism of commutative rings, let $X$ and $Y$ be $B$-schemes which are quasi-compact and quasi-separated,
and let $Z = X \times_{\Spec(B)} Y$ denote their product. Suppose that $X$ and $Y$ are $\Tor$-independent over $B$ (this assumption holds,
for example, if either $X$ or $Y$ is flat over $B$). Then the diagram of derived Hodge complexes
$$ \xymatrix@R=50pt@C=50pt{ \bigoplus_{d \in \Z} L \Omega^{d}_{B/A}[-d] \ar[r] \ar[d] & \bigoplus_{d \in \Z} \RGamma(X, L \Omega^{d}_{X/A} )[-d] \ar[d] \\
\bigoplus_{d \in \Z} \RGamma( Y, L \Omega^{d}_{Y/A})[-d] \ar[r] & \bigoplus_{d \in \Z} \RGamma( Z, L \Omega^{d}_{Z/A})[-d] }$$
is a pushout square in $\CAlg( \calD(A) )$. To prove this, we can work locally to reduce to the case where $X$ and $Y$
are affine, in which case the desired result follows from Remark \ref{remark:kunneth-for-derived-Hodge}.
\end{remark}

We will need the following result, whose proof is left to the reader:

\begin{proposition}\label{proposition:graded-Hodge-cohomology}
Let $A = \bigoplus_{n \geq 0} A_{n}$ be nonnegatively graded ring, so that the inclusion $A_0 \hookrightarrow A$ induces a map of quotient stacks
$$ f: [ \Spec(A) / \mathbf{G}_m ] \rightarrow [ \Spec(A_0) / \mathbf{G}_m ] \simeq \Spec(A_0) \times \BGm.$$
Then, for every integer $d$, the induced map
$$ \RGamma( [ \Spec(A_0) / \mathbf{G}_m], L \Omega^{d}_{  [ \Spec(A_0) / \mathbf{G}_m ] } )
\rightarrow \RGamma( [ \Spec(A) / \mathbf{G}_m], L \Omega^{d}_{  [ \Spec(A) / \mathbf{G}_m ] } )$$
is an isomorphism in the derived $\infty$-category $\calD(\Z)$.
\end{proposition}

\subsection{Quasisyntomic Morphisms}

\begin{definition}\label{definition:complete-flatness}
Let $R$ be a commutative ring and let $I \subseteq R$ be a finitely generated ideal. We say that an object $M \in \calD(R)$
is {\it $I$-completely flat} if the derived tensor product $(R/I) \otimes^{L}_{R} M$ is a flat $(R/I)$-module (regarded as a chain complex
concentrated in cohomological degree zero). We say that $M$ is {\it $p$-completely flat} if it is $I$-completely flat for the principal ideal $I = (p)$.
\end{definition}

\begin{warning}
In the situation of Definition \ref{definition:complete-flatness}, we do not require that $M$ is $I$-complete. However,
if $\widehat{M}$ is the $I$-completion of $M$, then $M$ is $I$-completely flat if and only if $\widehat{M}$ is $I$-completely flat. 
\end{warning}

\begin{example}
Let $R$ be a commutative ring and let $M$ be a flat $R$-module. Then $M$ is $I$-completely flat for every finitely generated ideal $I \subseteq R$.
\end{example}

\begin{remark}
Let $R$ be a Noetherian ring, let $I \subseteq R$ be an ideal, and let $M$ be an object of $\calD(R)$. Then $M$ is $I$-completely flat if and only
if the $I$-completion $\widehat{M}$ is a flat $R$-module, concentrated in cohomological degree zero.
\end{remark}

\begin{definition}
Let $f: R \rightarrow S$ be a homomorphism of commutative rings and let $I \subseteq R$ be a finitely generated ideal. We will say that
$f$ is {\it $I$-completely flat} if the commutative ring $S$ is $I$-completely flat when viewed as an $R$-module (via restriction of scalars along $f$).
Note that, in this case, the induced map of quotient rings $R/I \rightarrow S/IS$ is flat. We say that $f$ is {\it $I$-completely faithfully flat}
if it is $I$-completely flat and the map $R/I \rightarrow S/IS$ is faithfully flat.
\end{definition}

\begin{definition}\label{definition:prequasisyntomic}
Let $R$ be a commutative ring. We will say that $R$ is {\it $p$-quasisyntomic} if it satisfies the following conditions:
\begin{itemize}
\item[$(1)$] The commutative ring $R$ has bounded $p$-power torsion (in particular, the $p$-completion
$\widehat{R}$ coincides with the separated $p$-completion $\varprojlim_{m} R/p^{m}R$; see Remark \ref{remark:bounded-p-torsion}).

\item[$(2)$] The derived tensor product $(R/pR) \otimes^{L}_{R} L\Omega^{1}_{R} \in \calD(R/pR)$ has $\Tor$-amplitude contained in $[-1,0]$.
\end{itemize}
We let $\CAlg^{\QSyn}$ denote the category whose objects are $p$-quasisyntomic commutative rings (and whose morphisms are ring homomorphisms).
\end{definition}

\begin{warning}
Definition \ref{definition:prequasisyntomic} is a slight variant of Definition~1.7 of \cite{BMS2}. A commutative ring $R$ is
{\it quasisyntomic} in the sense of \cite{BMS2} if and only if if it is both $p$-quasisyntomic (in the sense of Definition \ref{definition:prequasisyntomic})
and $p$-complete. Conversely, a commutative ring $R$ is $p$-quasisyntomic if and only if its $p$-completion $\widehat{R}$ is a quasisyntomic
commutative ring (concentrated in cohomological degree zero).
\end{warning}

\begin{example}
Every regular Noetherian ring is $p$-quasisyntomic. In particular, every polynomial ring $\Z[x_1, x_2, \cdots, x_n]$ is $p$-quasisyntomic.
\end{example}

We will need a relative version of Definition \ref{definition:prequasisyntomic}.

\begin{definition}[The $p$-Quasisyntomic Topology]\label{definition:qsyn-topology}
Let $R$ be a commutative ring. We will say that a ring homomorphism $f: R \rightarrow S$ is {\it $p$-quasisyntomic} if it is $p$-completely flat and the complex
$$(S/pS) \otimes_{S}^{L} L\Omega^{1}_{S/R} \in \calD(S)$$
has $\Tor$-amplitude contained in $[0,-1]$. In this case, we say that $S$ is a {\it $p$-quasisyntomic $R$-algebra}. We let
$\CAlg^{\QSyn}_{R}$ denote the category of $p$-quasisyntomic $R$-algebras.

We say that a ring homomorphism $f: R \rightarrow S$ is a {\it $p$-quasisyntomic cover} if $f$ is $p$-quasisyntomic and $p$-completely faithfully flat.
The $p$-quasisyntomic coverings determine a Grothendieck topology on (the opposite of) the category of
$p$-quasisyntomic commutative rings, where a collection of morphisms $\{ R \rightarrow S_i \}_{i \in I}$ generate a covering sieve if and only if
there exists a finite subset $I_0 \subseteq I$ and a ring homomorphism $\prod_{i \in I_0} S_i \rightarrow S$ for which
the composite map $$ R \rightarrow \prod_{i \in I_0} S_i \rightarrow S$$
is a $p$-quasisyntomic covering (see Lemma~4.17 of \cite{BMS2}). We will refer to this Grothendieck topology as the {\it $p$-quasisyntomic topology} on
$( \CAlg^{\QSyn})^{\op}$.
\end{definition}

\begin{warning}
Let $f: R \rightarrow S$ be a ring homomorphism. If $R$ is $p$-quasisyntomic and $f$ is $p$-quasisyntomic, then $S$ is also $p$-quasisyntomic.
Conversely, if $S$ is $p$-quasisyntomic and $f$ is a $p$-quasisyntomic covering, then $R$ is $p$-quasisyntomic (see Lemma~4.15 of \cite{BMS2}).

In particular, if $S$ is a commutative ring for which the unit map $u: \Z \rightarrow S$ is $p$-quasisyntomic (in the sense of Definition \ref{definition:qsyn-topology}),
then $S$ is $p$-quasisyntomic (in the sense of Definition \ref{definition:prequasisyntomic}). Beware that the converse is false: the finite field $\F_p$
is a $p$-quasisyntomic commutative ring, but the quotient map $\Z \twoheadrightarrow \F_p$ is not a $p$-quasisyntomic ring homomorphism (since it is not $p$-completely flat).
It follows that the inclusion $\CAlg_{\Z}^{\QSyn} \subset \CAlg^{\QSyn}$ is strict.
\end{warning}

\begin{example}
Let $R$ be a commutative $\F_p$-algebra. The following conditions are equivalent:
\begin{itemize}
\item The commutative ring $R$ is $p$-quasisyntomic, in the sense of Definition \ref{definition:prequasisyntomic}.
\item The unit map $u: \F_p \rightarrow R$ is $p$-quasisyntomic, in the sense of Definition \ref{definition:qsyn-topology}.
\item The relative cotangent complex $L \Omega^{1}_{R/\F_p}$ has $\Tor$-amplitude contained in $[-1,0]$. 
\end{itemize}
If these conditions are satisfied, we will say that $R$ is a {\it quasisyntomic $\F_p$-algebra}. We say that
an $\F_p$-scheme $X$ is {\it quasisyntomic} if, for every affine open subset $U \subseteq X$, the coordinate
ring $\RGamma(U, \calO_U )$ is a quasisyntomic $\F_p$-algebra.
\end{example}

\begin{example}
Let $R$ be a commutative ring having bounded $p$-power torsion, and let $\widehat{R}$ denote its $p$-completion. Then the canonical map
$R \rightarrow \widehat{R}$ is a $p$-quasisyntomic covering.
\end{example}

\begin{example}
Every smooth ring homomorphism is $p$-quasisyntomic. In particular, every \'{e}tale ring homomorphism is quasisyntomic.
\end{example}

\begin{remark}
Let $\calC$ be an $\infty$-category which admits small limits, and let $\mathscr{F}: \CAlg^{\QSyn} \rightarrow \calC$ be a functor. 
Then $\mathscr{F}$ satisfies descent with respect to the $p$-quasisyntomic topology of Definition \ref{definition:qsyn-topology} if and only if
it satisfies the following conditions:
\begin{itemize}
\item[$(1)$] For every $p$-quasisyntomic commutative ring $R$, the canonical map $\mathscr{F}(R) \rightarrow \mathscr{F}( \widehat{R} )$ is
an isomorphism.

\item[$(2)$] The functor $\mathscr{F}$ commutes with finite products.

\item[$(3)$] Let $f: R \rightarrow R^{0}$ be a $p$-quasisyntomic covering between commutative rings which are $p$-quasisyntomic and $p$-complete,
and let $R^{\bullet}$ be the cosimplicial $R$-algebra obtained by $p$-completing the tensor powers of $R^{0}$ over $R$. Then the tautological map
$$ \mathscr{F}(R) \rightarrow \Tot( \mathscr{F}(R^{\bullet} ) )$$
is an isomorphism in $\calC$.
\end{itemize}
\end{remark}

\begin{definition}[\cite{BMS2}, Definition~4.20]\label{definition:qrsp}
Let $S$ be a commutative ring. We say that $S$ is {\it quasiregular semiperfectoid} if it satisfies the following conditions:
\begin{itemize}
\item The commutative ring $S$ is $p$-quasisyntomic and $p$-complete.
\item There exists a ring homomorphism $R \rightarrow S$, where $R$ is perfectoid.
\item The quotient ring $S/pS$ is semiperfect. That is, every element $x \in S$ can be written as $y^p + pz$ for some elements $y,z \in S$.
\end{itemize}
We let $\CAlg^{\qrsp}$ denote the category whose objects are quasiregular semiperfectoid rings and whose morphisms are ring homomorphisms.
\end{definition}

\begin{remark}
Let $R$ be a commutative ring. Then $R$ is $p$-quasisyntomic if and only if there exists a $p$-quasisyntomic covering $f: R \rightarrow S$, where
$S$ is quasiregular semiperfectoid (see Lemma~4.28 of \cite{BMS2}). In particular, the full subcategory $\CAlg^{\qrsp} \subseteq \CAlg^{\QSyn}$ forms a basis for the quasisyntomic topology of Definition~\ref{definition:qsyn-topology}.
\end{remark}

\begin{example}
Let $R$ be an $\F_p$-algebra. The following conditions are equivalent:
\begin{itemize}
\item The commutative ring $R$ is quasiregular semiperfectoid (Definition \ref{definition:qrsp}).
\item The $\F_p$-algebra $R$ is quasiregular semiperfect: that is, the Frobenius map $$ \varphi_{R}: R \rightarrow R \quad \quad x \mapsto x^{p}$$ 
is surjective and the shifted cotangent complex $L \Omega^{1}_{R/\F_p}[-1]$ is
a flat $R$-module (concentrated in cohomological degree zero).
\end{itemize}
We let $\CAlg^{\qrsp}_{\F_p}$ denote the category of $\F_p$-algebras which satisfy these conditions.
\end{example}

\subsection{Filtered Complexes}

In this appendix, we review the homological algebra of filtered complexes, expressed in a language which is convenient for our applications in this paper.

\begin{definition}\label{definition:filtered-derived}
Let $R$ be a commutative ring. We let $\DFilt(R)$ denote the $\infty$-category $\Fun( (\Z, \geq)^{\op}, \calD(R) )$ of functors; here $\calD(R)$ denotes the derived
$\infty$-category of $R$-modules, and $(\Z, \geq )$ denotes the linearly ordered set of integers. We refer to $\DFilt(R)$ as the {\it filtered derived $\infty$-category} of $R$. 
We denote objects of $\DFilt(R)$ by $\Fil^{\bullet}(M)$, which we view as diagrams 
$$ \cdots \rightarrow \Fil^{2}(M) \rightarrow \Fil^{1}(M) \rightarrow \Fil^{0}(M) \rightarrow \Fil^{-1}(M) \rightarrow \Fil^{-2}(M) \rightarrow \cdots$$
in the $\infty$-category $\calD(R)$. For each integer $n$, we let $\gr^{n}(M)$ denote the cofiber of the transition map $\Fil^{n+1}(M) \rightarrow \Fil^{n}(M)$,
formed in the $\infty$-category $\calD(R)$.
\end{definition}

\begin{warning}
Definition \ref{definition:filtered-derived} is not standard. There is a fully faithful embedding $\iota: \calD(R) \hookrightarrow \DFilt(R)$, which carries each complex $M \in \calD(R)$
to the constant diagram
$$ \cdots \rightarrow M \xrightarrow{\id} M \xrightarrow{\id} M \xrightarrow{\id} M \xrightarrow{\id} M \rightarrow \cdots;$$
the essential image of this embedding consists of those filtered complexes $\Fil^{\bullet}(M)$ for which the associated graded complex $\gr^{\bullet}(N)$ vanishes.
Most authors define the filtered derived $\infty$-category of $R$ to be (the homotopy category of) the Verdier quotient $\DFilt(R) / \calD(R)$, which can be identified
with the full subcategory $\DFiltComp(R) \subseteq \DFilt(R)$ introduced in Definition \ref{definition:filtration-complete} below.
\end{warning}

\begin{definition}[The Beilinson t-Structure]\label{definition:Beilinson-t}
Let $R$ be a commutative ring. We define a pair of full subcategories
$\DFilt(R)^{\geq 0}, \DFilt(R)^{\leq 0} \subseteq \DFilt(R)$ as follows:
\begin{itemize}
\item An object $\Fil^{\bullet}(M) \in \DFilt(R)$ belongs to $\DFilt(R)^{\geq 0}$ if and only if, for every integer $n$, the cohomology groups of the complex
$\Fil^{n}(M)$ are concentrated in degrees $\geq n$.

\item An object $\Fil^{\bullet}(M) \in \DFilt(R)$ belongs to $\DFilt(R)^{\leq 0}$ if and only if, for every integer $n$, the cohomology groups of the complex
$\gr^{n}(M) = \Fil^{n}(M) / \Fil^{n+1}(M)$ are concentrated in degrees $\leq n$.
\end{itemize}
\end{definition}

\begin{proposition}\label{proposition:Beilinson-exists}
Let $R$ be a commutative ring. Then the subcategories $( \DFilt(R)^{\geq 0}, \DFilt(R)^{\leq 0})$ determine a t-structure on the stable $\infty$-category $\DFilt(R)$.
\end{proposition}

\begin{proof}
See Theorem~5.4 of \cite{BMS2} (or the appendix to \cite{MR923133} in the case of bounded filtrations).
\end{proof}

\begin{definition}\label{definition:filtration-complete}
Let $R$ be a commutative ring and let $\Fil^{\bullet}(M)$ be an object of $\DFilt(R)$. We say that $\Fil^{\bullet}(M)$ is {\it filtration-complete}
if the limit $\varprojlim_{n} \Fil^{n}(M)$ vanishes in the $\infty$-category $\calD(R)$. We let $\DFiltComp(R)$ denote the full subcategory of $\DFilt(R)$ spanned by the filtration-complete objects. 
\end{definition}

We refer to the t-structure of Proposition \ref{proposition:Beilinson-exists} as the {\it Beilinson t-structure} on $\DFilt(R)$.
If $\Fil^{\bullet}(M)$ is an object of $\DFilt(R)$, we denote its truncations for the Beilinson t-structure by
$\tau^{\leq n}_{\Beil}(\Fil^{\bullet}(M))$ and $\tau^{\geq n}_{\Beil}( \Fil^{\bullet}(M) )$. 

\begin{remark}[The Decalage]\label{remark:decalage}
Let $\Fil^{\bullet}(M)$ be an object of $\DFilt(R)$. For every integer $n$, let us write $\tau^{\leq n}_{\Beil}(M)$ for the underlying complex of
filtered complex $\tau^{\leq n}_{\Beil}( \Fil^{\bullet}(M) )$. The diagram
$$ \cdots \rightarrow \tau^{\leq -1}_{\Beil}(M) \rightarrow \tau^{\leq 0}_{\Beil}(M) \rightarrow \tau^{\leq 1}_{\Beil}(M) \rightarrow \cdots$$
determines a new object of $\DFilt(R)$. At the level of filtered cochain complexes, the passage from $\Fil^{\bullet}(M)$ to
$\tau^{\leq -\bullet}_{\Beil}(M)$ is implemented by the {\it decalage construction} introduced by Deligne in \cite{Hodge2}.
\end{remark}

\begin{remark}\label{remark:bullet}
Let $R$ be a commutative ring and let $\Fil^{\bullet}(M)$ be an object of $\DFilt(R)$. Then $\Fil^{\bullet}(M)$ belongs to $\DFilt(R)^{\geq 0}$ if and only if it satisfies both of the following conditions:
\begin{itemize}
\item For every integer $n$, the cohomology groups of the complex $\gr^{n}(M)$ are concentrated in degrees $\geq n$.
\item The filtration $\Fil^{\bullet}(M)$ is complete: that is, the limit $\varprojlim_{n} \Fil^{n}(M)$ vanishes in the $\infty$-category $\calD(R)$.
\end{itemize}
\end{remark}

\begin{example}
\label{example:Beilinson-t-heart}
Let $(M^{\ast}, \partial)$ be a cochain complex of $R$-modules. For each integer $n$, let $(M^{\geq n}, \partial)$ denote the subcomplex of $(M^{\ast}, \partial)$
which coincides with $M^{\ast}$ in cohomological degrees $\geq n$ and vanishes in degrees $< n$. The construction
$n \mapsto (M^{\geq n}, \partial)$ determines an object of $\DFilt(R)$ which lies in the heart $\DFilt(R)^{\heartsuit} = \DFilt(R)^{\geq 0}(R) \cap \DFilt(R)^{\leq 0}(R)$
of the Beilinson t-structure. Moreover, the construction $M^{\ast} \mapsto M^{\geq \bullet}$ induces an equivalence of abelian categories
$$ \{ \text{Chain complexes of $R$-modules} \} \rightarrow \DFilt(R)^{\heartsuit}.$$
\end{example}

\begin{example}\label{example:I-adic-filtration-complete}
Let $R$ be a commutative ring, let $I \subseteq R$ be an invertible ideal, and let $M$ be an object of $\calD(R)$. Let $I^{\bullet} M$
denote the filtered complex corresponding to the diagram
For each element $f \in R$, let $I^{\bullet} M$ denote the filtered complex corresponding to the diagram
$$ \cdots \rightarrow I^{2} \otimes_{R}^{L} M \rightarrow I^{1} \otimes_{R}^{L} M 
\rightarrow I^{0} \otimes_{R}^{L} M \rightarrow I^{-1} \otimes_{R}^{L} M \rightarrow
I^{-2} \otimes_{R}^{L} M \rightarrow \cdots$$
Then $M$ is $I$-complete (in the sense of Definition \ref{definition:I-complete}) if and only if $I^{\bullet} M$ is filtration-complete.
\end{example}

\begin{remark}
\label{remark:Beilinson-connective-characterize}
Let $R$ be a commutative ring and let $\Fil^{\bullet}(M)$ be an object of $\DFilt(R)$. For every integer $n$, we have a fiber sequence
$$ \tau^{\leq n}_{\Beil}( \Fil^{\bullet}(M) ) \rightarrow \Fil^{\bullet}(M) \rightarrow \tau^{\geq n+1}_{\Beil}( \Fil^{\bullet}(M) ),$$
where the third term is automatically filtration-complete (Remark \ref{remark:bullet}). 
Consequently, $\Fil^{\bullet}(M)$ is filtration-complete if and only if $\tau^{\leq n}_{\Beil}( \Fil^{\bullet}(M) )$ is filtration-complete.
In particular, the subcategories
$$ \DFiltComp(R)^{\leq 0} = \DFilt(R)^{\leq 0} \cap \DFiltComp(R) \quad \quad \DFiltComp(R)^{\geq 0} = \DFilt(R)^{\geq 0} \cap \DFiltComp(R)$$
determine a t-structure on the complete filtered derived $\infty$-category $\DFiltComp(R)$, which we will also refer to as the {\it Beilinson t-structure}.
\end{remark}

\begin{remark}\label{remark:characterize-Beilinson-connective}
Let $R$ be a commutative ring and let $u: \Fil^{\bullet}(M') \rightarrow \Fil^{\bullet}(M)$ be a morphism in $\DFilt(R)$. 
Assume that $\Fil^{\bullet}(M)$ is filtration-complete. The following conditions are equivalent:
\begin{itemize}
\item The morphism $u$ exhibits $\Fil^{\bullet}(M')$ as a connective cover of $\Fil^{\bullet}(M)$ with respect to the Beilinson t-structure:
that is, it induces an isomorphism $\Fil^{\bullet}(M') \simeq \tau^{\leq 0}_{\Beil}( \Fil^{\bullet}(M) )$.

\item The filtered complex $\Fil^{\bullet}(M')$ is filtration-complete and the morphism $u$ induces isomorphisms of abelian groups
$$ \mathrm{H}^{m}( \gr^{n}(M') ) \xrightarrow{\sim} \begin{cases} \mathrm{H}^{m}( \gr^{n}(M) ) & \text{ if } m \leq n \\
0 & \text{ otherwise. } \end{cases}$$
\end{itemize}
\end{remark}

We will often need to contemplate a different completeness condition on filtered complexes.

\begin{notation}\label{notation:filtered-derived-I-complete}
Let $R$ be a commutative ring, let $I \subseteq R$ be a finitely generated ideal, and let
$\Fil^{\bullet}(M)$ be an object of the filtered derived $\infty$-category $\DFilt(R)$. We will say that
$\Fil^{\bullet}(M)$ is {\it $I$-complete} if each of the complexes $\Fil^{n}(M)$ is $I$-complete. We let $\DFiltI(R)$ denote the full subcategory of
$\DFilt(R)$ spanned by the $I$-complete filtered complexes: that is, the $\infty$-category of
functors from the linearly ordered set $(\Z, \geq)^{\op}$ to the complete derived $\infty$-category $\widehat{\calD}(R)$
of Notation \ref{notation:complete-derived}. We let $\DFiltIComp(R) = \DFiltI(R) \cap \DFiltComp(R)$ denote the full subcategory of $\DFilt(R)$ spanned by those filtered complexes which are both $I$-complete and filtration-complete.
\end{notation}

\begin{warning}
Notation \ref{notation:filtered-derived-I-complete} is potentially ambiguous: the $\infty$-category $\DFiltI(R)$ depends not only on $R$,
but also on the finitely generated ideal $I \subseteq R$. However, we will use this notation only in three situations:
\begin{itemize}
\item If $(A,I)$ is a prism, we write $\DFiltI(A)$ for the filtered derived $\infty$-category of $(p,I)$-complete filtered complexes over $A$.
\item If $(A,I)$ is a prism, we write $\DFiltI( \overline{A} )$ for the filtered derived $\infty$-category of $p$-complete filtered complexes
over the quotient ring $\overline{A} = A/I$.
\item We write $\DFiltI(\Z_p)$ for the filtered derived $\infty$-category of $p$-complete complexes over the ring $\Z_p$ of $p$-adic integers
(or equivalently over $\Z$).
\end{itemize}
\end{warning}

\begin{proposition}\label{proposition:check-truncations}
Let $R$ be a commutative ring, let $I \subseteq R$ be a finitely generated ideal, and let
$\Fil^{\bullet}(M) \in \DFilt(R)$ be an $I$-complete filtered complex. Then, for every integer $n$,
the Beilinson truncations $\tau^{\leq n}_{\mathrm{Dec}} \Fil^{\bullet}(M)$ and $\tau^{\geq n}_{\mathrm{Dec}} \Fil^{\bullet}(M)$
are also $I$-complete.
\end{proposition}

\begin{proof}
We will prove that the truncations $\tau^{\geq n}_{\mathrm{Dec}} \Fil^{\bullet}(M)$ are $I$-complete; the analogous
assertions for $\tau^{\leq n}_{\mathrm{Dec}} \Fil^{\bullet}(M)$ follows from the fiber sequence
$$ \tau^{\leq n}_{\mathrm{Dec}} \Fil^{\bullet}(M) \rightarrow \Fil^{\bullet}(M) \rightarrow
\tau^{\geq n+1}_{\mathrm{Dec}} \Fil^{\bullet}(M).$$
Without loss of generality, we may assume that $n=0$. Let us denote the truncation $\tau^{\geq 0}_{\mathrm{Dec}} \Fil^{\bullet}(M)$ by $\Fil^{\bullet}(N)$. We wish to show that, for every integer $m$, the complex $\Fil^{m}(N) \in \calD(R)$ is $I$-complete.

For each integer $m$, let $\Fil^{m}(N)^{\wedge}_{I}$ denote the $I$-completion of $\Fil^{m}(N)$, and let 
$K(m)$ denote the fiber of the natural map $\Fil^{m}( N) \rightarrow \Fil^{m}( N )^{\wedge}_{I}$.
We wish to show that each of the complexes $K(m)$ vanishes (as an object of $\calD(R)$). 
It follows from the definition of the Beilinson t-structure that each $\Fil^{m}(N)$ is concentrated in cohomological degrees $\geq m$.
Choose an integer $d \geq 0$ such that the ideal $I$ is generated by $d$ elements,
so that each $\Fil^{m}( N )^{\wedge}_{I}$ has cohomology concentrated in degrees $\geq m-d$ (see \cite[Tag 091V]{stacks-project}),
and therefore $K(m)$ has cohomology concentrated in degrees $\geq m-d$. We will complete the proof by showing that each
of the transition maps $K(m+1) \rightarrow K(m)$ is an isomorphism in $\calD(R)$. Note that we have a canonical isomorphism
$$ \cofib( K(m+1) \rightarrow K(m) ) \simeq \fib( \gr^{m}(N) \rightarrow \gr^{m}(N)^{\wedge}_{I} ).$$
We are therefore reduced to proving that the complex $\gr^{m}(N)$ is $I$-complete. This is clear,
since since $\gr^{m}(N)$ can be obtained by applying the cohomological truncation functor $\tau^{\geq n}$ to
the $I$-complete complex $\gr^{m}(M) = \Fil^{m}(M) / \Fil^{m+1}(M)$.
\end{proof}

\begin{corollary}\label{corollary:Beilinson-descends}
Let $R$ be a commutative ring and let $I \subseteq R$ be a finitely generated ideal. Then
the subcategories
$$ \DFiltI(R)^{\leq 0} = \DFilt(R)^{\leq 0} \cap \DFiltI(R) \quad \quad \DFiltI(R)^{\geq 0} = \DFilt(R)^{\geq 0} \cap \DFiltI(R)$$
determine a t-structure on $\DFiltI(R)$, and the subcategories
$$ \DFiltIComp(R)^{\leq 0} = \DFilt(R)^{\leq 0} \cap \DFiltIComp(R) \quad \quad \DFiltIComp(R)^{\geq 0} = \DFilt(R)^{\geq 0} \cap \DFiltIComp(R)$$
determine a t-structure on $\DFiltIComp(R)$.
\end{corollary}

We will abuse terminology by referring to the t-structures of Corollary \ref{corollary:Beilinson-descends} as the
{\it Beilinson t-structures} on $\DFiltI(R)$ and $\DFiltIComp(R)$.

\begin{remark}
Let $R$ be a commutative ring, let $I \subseteq R$ be a finitely generated ideal, and let
$\Fil^{\bullet}(M)$ be an object of $\DFilt(R)$. Then the filtered complex $\Fil^{\bullet}(M)$ is $I$-complete
if and only if each of truncations $\Fil^{\bullet}(M^{\leq n} ) = \tau^{\leq n}_{\mathrm{Dec}} \Fil^{\bullet}(M)$ is $I$-complete.
The ``only if'' direction follows from Proposition \ref{proposition:check-truncations}. To prove the converse, it
suffices to observe that for every triple of integers $k$, $m$, and $n$, the map of filtered complexes
$\Fil^{\bullet}( M^{\leq n} ) \rightarrow \Fil^{\bullet}(M)$ induces $R$-module homomorphisms $\mathrm{H}^{k}( \Fil^{m}( M^{\leq n} ) \rightarrow \mathrm{H}^{k}( \Fil^{m}( M) )$ which are isomorphisms for $k \leq m+n$. 
\end{remark}

\begin{warning}
Let $R$ be a commutative ring, let $I \subseteq R$ be a finitely generated ideal, and let
$\Fil^{\bullet}(M)$ be an object of $\DFilt(R)$. If $\Fil^{\bullet}(M)$ is filtration-complete and
each truncation $\tau^{\geq n}_{\mathrm{Dec}} \Fil^{\bullet}(M)$ is $I$-complete, then
$\Fil^{\bullet}(M)$ is also $I$-complete. Beware that the assumption that $\Fil^{\bullet}(M)$ is filtration-complete
cannot be omitted.
\end{warning}

\subsection{De Rham Complexes}

Throughout this section, we fix a commutative ring $A$.

\begin{notation}[The Classical de Rham Complex]\label{notation:classical-de-Rham}
Let $R$ be an $A$-algebra. For every integer $n \geq 0$, we let $\Omega^{n}_{R/A}$ denote the $n$th exterior power of the
module of K\"{a}hler differentials $\Omega^{1}_{R/A}$. Allowing $n$ to vary, we obtain a chain complex
$$ R = \Omega^{0}_{R/A} \xrightarrow{d} \Omega^{1}_{ R/ A} \xrightarrow{d} \Omega^{2}_{R/A} \rightarrow \cdots$$
We will denote this $A$-linear complex by $(\Omega^{\ast}_{R/A}, d)$ and refer to it as the {\it algebraic de Rham complex of $R$ relative to $A$}.

For every integer $n \geq 0$, we let $\widehat{\Omega}^{n}_{R / A}$ denote the $p$-completion of
$\Omega^{n}_{R/ A}$ in the abelian category of $R$-modules. Beware that $\widehat{\Omega}^{n}_{R / A}$ is not
necessarily $p$-adically separated (see Warning \ref{warning:compare-completions}). Allowing $n$ to vary, we obtain a chain complex
$$ \widehat{\Omega}^{0}_{R/A} \xrightarrow{d} \widehat{\Omega}^{1}_{ R/ A} \xrightarrow{d} \widehat{\Omega}^{2}_{R/A} \rightarrow \cdots$$
We will denote this $A$-linear complex by $(\widehat{\Omega}^{\ast}_{R/A}, d)$ and refer to it as the {\it $p$-complete de Rham complex of $R$ relative to $A$}.

More generally, if $R$ is an animated commutative $A$-algebra, we define $(\Omega^{\ast}_{R/A}, d)$ and
$(\widehat{\Omega}^{\ast}_{R/A}, d)$ to be the de Rham complexes 
$(\Omega^{\ast}_{\pi_0(R)/A}, d)$ and $(\widehat{\Omega}^{\ast}_{\pi_0(R)/A}, d)$, respectively, where $\pi_0(R)$ denotes the underlying commutative $A$-algebra of $R$.
\end{notation}
We now review the definition of derived de Rham cohomology, following \cite{bhatt2012padic}.

\begin{construction}[The derived de Rham complex and its Hodge filtration]
\label{construction:derived-de-Rham}
The construction carrying a polynomial algebra $R$ over $A$ to the filtered complex
\[\Fil^\bullet_{\Hodge} \Omega_{R/A}  := \left( \cdots \to (\Omega^{\geq 3}_{R/A},d) \to (\Omega^{\geq 2}_{R/A}, d) \to (\Omega^{\geq 1}_{R/A},d) \to (\Omega^\ast_{R/A},d)\right)\]
determines a functor $\Fil^\bullet_{\Hodge} \Omega_{-/A}$ from the category $\Poly_{A}$ to the filtered derived $\infty$-category $\DFilt(A)$;
here we set $\Fil^i_{\Hodge} \Omega_{-/A} = \Omega_{-/A}$ for $i \leq 0$. Applying Proposition \ref{proposition:universal-of-animated}, we see that this construction admits an essentially unique extension to a
functor of $\infty$-categories $\CAlg_{A}^{\anim} \rightarrow \DFilt(A)$ which commutes with sifted colimits. We denote this extension by $R \mapsto \Fil^{\bullet}_{\Hodge} {\dR}_{R/A}$. We also denote $\Fil^{0}_{\Hodge} {\dR}_{R/A}$ by ${\dR}_{R/A}$ and refer to it as the {\it derived de Rham complex of $R$ relative to $A$}. The filtered complex $\Fil^\bullet_{\Hodge} \dR_{R/A}$ can be viewed as a decreasing $\mathbf{N}$-indexed filtration on $\dR_{R/A} \in \calD(A)$; we refer to this filtration as the  {\it Hodge filtration} on $\dR_{R/A}$. By construction, we have a natural isomorphism
\[ \gr^n_{\Hodge} \dR_{R/A} \simeq L\Omega^n_{R/A}[-n]\]
for each $n \geq 0$. 
\end{construction}

Beware that derived de Rham complex $\dR_{R/A}$ of Construction \ref{construction:derived-de-Rham} has somewhat pathological behavior rational behavior:
for example, if $A$ is a $\Q$-algebra, then the unit map $A \rightarrow \dR_{R/A}$ is an isomorphism for every (animated) commutative $A$-algebra $R$
(see \cite[Corollary 2.5]{bhatt2012padic}). 

\begin{construction}[The $p$-completed derived de Rham complex]
Let $R$ be an animated commutative $A$-algebra. For every integer $n$, we let $\Fil^n_{\Hodge} \dR_{R/A}$
denote the $p$-completion of the complex $\Fil^{n}_{\Hodge} \dR_{R/A}$. We will refer to $\widehat{\dR}_{R/A} := \Fil^0_{\Hodge} \widehat{\dR}_{R/A}$ as {\em the $p$-completed derived de Rham complex}, and to $\Fil^\bullet_{\Hodge} \widehat{\dR}_{R/A}$ as the {\em Hodge filtration} on $\widehat{\dR}_{R/A}$. By construction, we have a natural isomorphism \[ \gr^n_{\Hodge} \widehat{\dR}_{R/A} \simeq L\Omega^n_{R/A}[-n]^{\wedge}_p\]
for each $n \geq 0$.
\end{construction}

\begin{remark}
\label{remark:p-comp-dR}
Let $R$ be an animated $A$-algebra with $p$-completion $\widehat{R}$. The natural map 
$$\Fil^\bullet_{\Hodge} \widehat{\dR}_{R/A} \to \Fil^\bullet_{\Hodge} \widehat{\dR}_{\widehat{R}/A}$$ is an isomorphism. In other words, the $p$-completed derived de Rham complex (and its Hodge filtration) are insensitive to replacing the animated $A$-algebra with its $p$-completion.
\end{remark}

\begin{remark}\label{remark:gr-derived-deRham}
LLet $R$ be an animated commutative $A$-algebra, and let $\widehat{\dR}_{R/A}$ denote the $p$-adic derived de Rham complex
of Construction \ref{construction:derived-de-Rham}. Then, for every integer $d$, we have a canonical isomorphism
$$\gr^{d}_{\Hodge} \widehat{\dR}_{B/A} \simeq L \widehat{\Omega}^{d}_{B/A}[-d]$$
in the $p$-complete derived $\infty$-category $\widehat{\calD}(B)$.
\end{remark}

\begin{notation}[The de Rham Augmentation]\label{notation:dra}
Let $R$ be an animated commutative $A$-algebra. We let
$\epsilon_{\dR}: \widehat{\dR}_{R/A} \rightarrow \widehat{R}$ denote the composition
$$ \widehat{\dR}_{R/A} = \Fil^{0}_{\Hodge} \widehat{\dR}_{R/A}  \rightarrow \gr^{0}_{\Hodge} \widehat{\dR}_{R/A} \simeq 
L \widehat{\Omega}^{0}_{R/A} \simeq \widehat{R}.$$
We will refer to $\epsilon_{\dR}$ as the {\it de Rham augmentation}.
\end{notation}

\begin{construction}[The Conjugate Filtration]\label{construction:confilt}
Let $A$ be an $\F_p$-algebra. For every integer $n \geq 0$, Proposition \ref{proposition:universal-of-animated} implies that there is an essentially unique functor
$$ \Fil_{n}^{\conj} \dR_{-/A}: \CAlg_{A}^{\anim} \rightarrow \calD(A)$$
which commutes with sifted colimits and satisfies $\Fil_{n}^{\conj} \dR_{R/A} = \tau^{\leq n} \dR_{R/A}$ for $R \in \Poly_{A}$.
If $R$ is an animated commutative $A$-algebra, then the construction $n \mapsto \Fil_{n}^{\conj} \dR_{R/A}$ determines a diagram
$$ \Fil_{0}^{\conj} \dR_{R/A} \rightarrow \Fil_{1}^{\conj} \dR_{R/A} \rightarrow \Fil_{2}^{\conj} \dR_{R/A} \rightarrow \cdots$$
whose colimit is the derived de Rham complex $\dR_{R/A}$. We will refer to this diagram as the {\it conjugate filtration} on $\dR_{R/A}$.
\end{construction}

\begin{remark}[The Cartier Isomorphism]\label{cons:conjugate-filt}
Let $A$ be an $\F_p$-algebra, let $R$ be a flat $A$-algebra, and let $R^{(1)}$ denote its Frobenius twist $R^{(1)} := A \otimes_{\varphi,A} R$ denote the 
animated commutative $A$-algebra obtained from $R$ by extending scalars along the Frobenius map. Note that the differential on the de Rham
complex $(\Omega^{\bullet}_{R/A},d)$ is $R^{(1)}$-linear, where $R^{(1)}$ acts via the relative Frobenius map $\varphi_{R/A}: R^{(1)} \rightarrow R$.
Moreover, there is a natural map
\[ \mathrm{Cart}^{-1}: \Omega^*_{R^{(1)}/A} \to \mathrm{H}^*(\Omega_{R/A})\]
of graded $R^{(1)}$-algebras, which is uniquely determined by the formula
$\mathrm{Cart}^{-1}(d(1 \otimes y)) = y^{p-1} dy \in \mathrm{H}^1(\Omega_{R/A})$ for $y \in R$. If $R$ is a smooth $A$-algebra, then
$\mathrm{Cart}^{-1}$ is an isomorphism (see \cite[Theorem 1.2]{MR894379}); we will refer to $\Cart^{-1}$ as the {\it Cartier isomorphism}.
In particular, $\Cart^{-1}$ is an isomorphism if $R$ is a polynomial algebra over $A$. 

Now suppose that $R$ is an animated commutative $A$-algebra, and let $R^{(1)}$ denote the derived tensor product  $A \otimes_{\varphi,A}^{L} R$.
For each $n \geq 0$, Remark \ref{rmk:LKECAlgAnim} supplies an isomorphism
$$ \Cart^{-1}: L \Omega^{n}_{ R^{(1)} / A }[-n] \rightarrow \gr_{n}^{\conj} \dR_{R/A},$$
which is uniquely determined by the requirement that it depends functorially on $R$ and agrees with the Cartier isomorphism when $R$ is a polynomial
algebra over $A$. By slight abuse of terminology, we will refer to this extension also as the {\it Cartier isomorphism}. Note that
if $R$ is a smooth $A$-algebra, then the existence of this isomorphism guarantees that $\gr_{n}^{\conj} \dR_{R/A}$ is concentrated in cohomological
degree $n$. It follows that the conjugate filtration on $\dR_{R/A}$ coincides with the Postnikov filtration: that is, we have
$\Fil_{d}^{\conj} \dR_{R/A} \simeq \tau^{\leq d} \dR_{R/A}$ for every integer $d$.
\end{remark}

\begin{remark}\label{remark:derived-deRham-product}
For every commutative ring $A$, the functor
$$ \CAlg_{A}^{\anim} \rightarrow \widehat{\calD}(A) \quad \quad R \mapsto \dR_{R/A}$$
commutes with finite products. To prove this, we can assume without loss of generality that $A$ is an $\F_p$-algebra.
It will then suffice to show that the functor $R \mapsto \gr_{n}^{\conj} \dR_{R/A}$ commutes with finite products
for each integer $n \geq 0$. This follows by induction on $n$, using the isomorphisms
$\gr_{n}^{\conj} \dR_{R/A} \simeq L\Omega^{n}_{ R^{(1)}/A}[-n]$ of Remark \ref{cons:conjugate-filt}.
\end{remark}

For every animated commutative ${A}$-algebra $R$, there is a filtered comparison map $\Fil^\bullet(\rho_{R}): \Fil^\bullet_{\Hodge} \widehat{\dR}_{R/A} \rightarrow ( \widehat{\Omega}^{\geq \bullet}_{R/A}, d)$ in the $\infty$-category $\widehat{\DFilt}(A)$ which is characterized (up to essentially unique homotopy) by the requirement that it depends functorially on $R$ and restricts to the identity when $R$ is a finitely generated polynomial algebra over $A$. 
The map $\Fil^\bullet(\rho_{R})$ is also an isomorphism when the commutative ring $R$ is $p$-completely smooth over $A$ (Corollary~3.10 of \cite{bhatt2012padic}). 
We now formulate a generalization of this result.

\begin{definition}\label{definition:Cartier-smooth}
Let $f: A \rightarrow R$ be a homomorphism of $\F_p$-algebras. We say that $f$ is {\it Cartier smooth} if it satisfies the following conditions:
\begin{itemize}
\item[$(1)$] The morphism $f$ is flat.
\item[$(2)$] The relative cotangent complex $L\Omega^{1}_{R/A}$ is a flat $R$-module.
\item[$(3)$] The Cartier operator $\Cart^{-1}: \Omega^*_{R^{(1)}/A} \to \mathrm{H}^*(\Omega_{R/A},d)$ is an isomorphism in every degree.
\end{itemize}
\end{definition}

\begin{remark}
Every smooth (or ind-smooth) morphism of $\F_p$-algebras is Cartier smooth. In particular, every regular noetherian $\F_p$-algebra is Cartier smooth over $\F_p$ (this follows from Popescu's theorem; see \cite[\S 9.5]{DRW} for a more direct proof). Gabber has also shown (unpublished) that any valuation ring of characteristic $p$ is Cartier smooth over $\F_p$. We refer to \cite[\S 2]{KellyMorrowValuation} for further discussion.
\end{remark}

\begin{proposition}\label{proposition:derived-to-classical-de-Rham}
Let $A$ be a commutative ring with bounded $p$-power torsion, and let $R$ be a $p$-completely flat $A$-algebra. If
the underlying morphism of $\F_p$-algebras $A/pA \rightarrow R/pR $ is Cartier smooth, then the map
$$ \Fil^{\bullet}( \rho_{R} ):  \Fil^\bullet_{\Hodge} \widehat{\dR}_{R/A} \rightarrow ( \widehat{\Omega}^{\geq \bullet}_{R/A}, d)$$
is an isomorphism (in the filtered derived $\infty$-category $\widehat{\DFilt}(A)$).
\end{proposition}

\begin{proof}
As $R$ is $p$-completely flat over $A$, we have $R \otimes_A^L A/pA \simeq R/pR$. The Cartier smoothness hypothesis and base change for the cotangent complex then show that $L\Omega^1_{R/A}$ is $p$-completely flat over $R$. It follows that the natural map $L\Omega^1_{R/A} \to \Omega^1_{R/A}$ becomes an isomorphism after $p$-completion. Since $A$ has bounded $p$-torsion and $R$ is $p$-completely flat over $A$, it follows that $\widehat{\Omega}^1_{R/A}$ is a $p$-completely flat $R$-module and that $(L\Omega^n_{R/A})^{\wedge}_p \to \widehat{\Omega}^n_{R/A}$ is an isomorphism of $p$-completely flat $R$-modules for all $n \geq 0$. In other words, for each $n$, the map 
\[ \gr^n(\rho_R):  L\Omega^n_{R/A}[-n]^{\wedge}_p \to \widehat{\Omega}^n_{R/A}[-n]\]
is an isomorphism. It is therefore enough to check that the map $\rho_R:\widehat{\dR}_{R/A} \to (\Omega^{\ast}_{R/A},d)^{\wedge}_p$ is an isomorphism in $\widehat{\calD}(A)$.
Since the source and target of $\rho_{R}$ are $p$-complete, it will suffice to check this after extending scalars along the map $A \twoheadrightarrow A/pA$. Since
Our assumptions on $A$ and $R$ (and the description of the target in terms of the cotangent complex) guarantee that the formation of $\rho_{R}$
is compatible with base change in $A$. We may therefore assume without loss of generality that $A$ is an $\mathbf{F}_p$-algebra. 

Let $\calC \subseteq \CAlg_{A}^{\anim}$ denote the category of Cartier smooth $A$-algebras. To complete the proof, it will suffice to show that the functor
$$ \calC \rightarrow \calD(A) \quad \quad R \mapsto ( \Omega^{\ast}_{R/A}, d)$$
is a left Kan extension of its restriction to the full subcategory $\Poly_{A} \subseteq \calC$. Writing the de Rham complex $( \Omega^{\ast}_{R/A}, d)$ as a colimit
of its truncations, we are reduced to showing that for every integer $n$, the functor
$$ \calC \rightarrow \calD(A) \quad \quad R \mapsto \mathrm{H}^{n}( \Omega^{\ast}_{R/A}, d)$$
is a left Kan extension of its restriction to $\Poly_{A}$. Invoking the definition of Cartier smoothness, we can rewrite this functor
as $R \mapsto \Omega^{n}_{ R^{(1)}/ A}$, in which case the desired result follows from Remark \ref{remark:LOmega-is-Omega}.
\end{proof}

\begin{warning}
In the formulation of Proposition \ref{proposition:derived-to-classical-de-Rham}, condition $(3)$ of Definition \ref{definition:Cartier-smooth}
cannot be omitted. Gabber has given examples of non-reduced $\mathbf{F}_p$-algebras $R$ with $L_{R/\F_p} = 0$ (see \cite{GabberImperfect}). For such algebras, the map $\rho_R:\dR_{R/\F_p} \to \Omega_{R/\F_p}$ identifies with the Frobenius endomorphism of $R$, which is not an isomorphism.
\end{warning}

Beware that the derived de Rham complex $\widehat{\dR}_{R/A}$ is not necessarily complete with respect to its Hodge filtration: that is, the limit $\varprojlim_{n} \Fil^{n}_{\Hodge} \widehat{\dR}_{R/A}$ is generally nonzero. For instance, this happens if $A = \mathbf{F}_p$ and $R = A[x]/(x^2)$. For this reason, it is sometimes useful
to contemplate a different completion of derived de Rham cohomology.

\begin{construction}[The Hodge-completed derived de Rham complex]\label{construction:Hodge-complete-dR}
Let $R$ be an animated commutative $A$-algebra. We let $\dR^{\hc}_{R/A}$ denote the completion of the derived de Rham complex $\dR_{R/A}$
with respect to its Hodge filtration: that is, the homotopy limit of the tower
$$ \cdots \rightarrow \dR_{R/A} / \Fil^{2}_{\Hodge} \dR_{R/A} \rightarrow \dR_{R/A} / \Fil^{1}_{\Hodge} \dR_{R/A} \rightarrow 
\dR_{R/A} / \Fil^{0}_{\Hodge} \dR_{R/A} = 0.$$
We will refer to ${\dR}^{\hc}_{R/A}$ as the {\em Hodge-completed derived de Rham complex} of $R$. It is equipped with
a Hodge filtration $\Fil^{\bullet}_{\Hodge} \dR^{\hc}_{R/A}$ (see Construction~\ref{construction:derived-de-Rham}), with
associated graded given by
$$ \gr^{n}_{ \Hodge} \dR^{\hc}_{R/A} \xleftarrow{\sim} \gr^{n}_{\Hodge} \dR_{R/A} \simeq L\Omega^{n}_{R/A}[-n].$$
\end{construction}

\begin{remark}
If $A=\mathbf{C}$ is the field of complex numbers and $R/A$ is a finitely generated algebra corresponding to the affine variety $X=\Spec(R)$, 
then \cite{bhatt-completions} constructs a canonical isomorphism $\dR^{\hc}_{R/A} \simeq \RGamma(X(\mathbf{C}), \mathbf{C})$: that is,
the Hodge-completed derived de Rham complex $\dR^{\hc}{R/A}$ computes the Betti cohomology of the underlying complex analytic space $X(\mathbf{C})$.
Beware that this isomorphism does not carry the Hodge filtration on $\dR^{\hc}_{R/A}$ to the Hodge filtration on $\RGamma(X(\mathbf{C}),\mathbf{C})$ 
introduced by Deligne in \cite{Hodge2} and \cite{Hodge3}.
\end{remark}

\begin{proposition}[Descent for derived de Rham complexes]
\label{proposition:qs-de-Rham}
Let $A$ be a commutative ring. Then the $\DFiltComp(A)$-valued functor $R \mapsto \dR^{\hc}_{R/A}$ satisfies descent for the flat topology on $\CAlg_A^{\anim}$.
\end{proposition}

\begin{proof}
Since isomorphisms in $\DFiltComp(A)$ can be detected at the associated graded level, it suffices to observe that for every integer $n$, the functor
$R \mapsto \gr^{n}_{\Hodge} \dR^{\hc}_{R/A} \simeq L\Omega^{n}_{R/A}[-n]$ satisfies descent for the flat topology 
(see Remark \ref{remark:descent-for-LOmega}).
\end{proof}

\begin{variant}\label{variant:qs-de-Rham}
The functor
$$ \CAlg^{\QSyn} \rightarrow \calD(\Z) \quad \quad R \mapsto \widehat{\dR}_{ R / \Z}$$
satisfies descent for the quasisyntomic topology.
\end{variant}

\begin{proof}
By virtue of Remark \ref{remark:derived-deRham-product} and Remark~\ref{remark:p-comp-dR}, it will suffice to show that
for every $p$-quasisyntomic cover $R \rightarrow R^0$ between $p$-complete objects of $\CAlg^{\QSyn}$ having
$p$-completed \v{C}ech nerve $R^{\bullet}$, the tautological map $\widehat{\dR}_{R/\Z} \rightarrow \Tot( \widehat{\dR}_{ R^{\bullet} / \Z} )$
is an isomorphism in $\widehat{\calD}(\Z_p)$. Note that this assertion can be checked after extending scalars to $\F_p$. Setting
$\overline{R} = \F_p \otimes^{L} R$ and $\overline{R}^{\bullet} = \F_p \otimes^{L} R^{\bullet}$, we are reduced to showing that the comparison map
$\theta: \dR_{ \overline{R} / \F_p} \rightarrow \Tot( \dR_{ \overline{R}^{\bullet} / \F_p } )$ is an isomorphism. Note that we can regard
$\theta$ as a morphism of filtered complexes (using the conjugate filtration of Construction \ref{construction:confilt}) and that
it is an isomorphism at the associated graded level (by combining the Cartier isomorphism of Remark \ref{cons:conjugate-filt} with
Remark \ref{remark:descent-for-LOmega}). To complete the proof, it will suffice to show that the conjugate filtration on
$\Tot( \dR_{ \overline{R}^{\bullet} / \F_p } )$ is exhaustive: that is, that the tautological map
\[ \varinjlim_{n} \Fil_n^{\conj} \dR_{ \overline{R}^{\bullet} /\F_p} \simeq \Tot \varinjlim_n \Fil_n^{\conj} \dR_{\overline{R}^{\bullet}/\F_p}\]
is an isomorphism. This follows from the observation the formation of totalizations in $\calD^{\geq -1}(\F_p)$ commutes with filtered colimits,
since each of the complexes $\Fil_n^{\conj} \dR_{\overline{R}^{\bullet}/\F_p}$ is concentrated in cohomological degrees $\geq -1$
(by virtue of our assumptions on $R$ and $R^{0}$).
\end{proof}

\begin{proposition}\label{proposition:de-Rham-qrsp}
Let $R$ be a quasiregular semiperfectoid $\Z_p$-algebra which is $p$-torsion-free. Then $\widehat{\dR}_{R/\Z_p}$ is concentrated in degree $0$ and $p$-torsion-free. 
\end{proposition}


\begin{proof}
It will suffice to show that the derived tensor product $\widehat{\dR}_{R/\Z_p} \otimes_{\Z_p}^L \F_p$ is concentrated in cohomological degree $0$. Setting $\overline{R} = R/pR$
(and invoking the compatibility of derived de Rham cohomology with base change), we are reduced to checking that $\dR_{ \overline{R} / \F_p}$ is concentrated in degree zero.
Our assumption that $R$ is quasiregular semiperfectoid (and $p$-torsion-free) guarantees that $\overline{R}$ is quasiregular semiperfect. In particular,
the cotangent complex $L\Omega^{1}_{ \overline{R} / \F_p}[-1]$ is a flat $R$-module. It follows that each of the complexes
$L \Omega^{n}_{ \overline{R} / \F_p}[-n] \simeq \gr_n^{\conj} \dR_{ \overline{R} / \F_p }$ is concentrated in cohomological degree zero (Remark \ref{cons:conjugate-filt}).
Proceeding by induction on $n$, we see that each $\Fil_{n}^{\conj} \dR_{\overline{R}/\F_p}$ is also concentrated in cohomological degree zero, so that
$\dR_{ \overline{R} / \F_p} = \varinjlim_{n} \Fil_{n}^{\conj} \dR_{ \overline{R}/ \F_p}$ has the same property.
\end{proof}

\subsection{Crystalline Cohomology}

In this section, we briefly review the construction of crystalline cohomology in a form that is convenient for our applications in this paper.

\begin{notation}\label{notation:big-crystalline-site}
Let $X$ be an $\F_p$-scheme. We define a category $\Crys(X / \Z_p)$ as follows:
\begin{itemize}
\item The objects of $\Crys(X / \Z_p)$ are triples $(A,I, v)$, where $(A,I)$ is a divided power algebra over $(\Z_p, (p) )$,
the commutative ring $A$ is $p$-adically separated and complete, and $v: \Spec(A/I) \rightarrow X$ is a morphism of $\F_p$-schemes.

\item A morphism from $(A,I,v)$ to $(A', I', v')$ in $\Crys(X / \Z_p)$ is a homomorphism of divided power rings $(A,I) \rightarrow (A', I')$
for which the diagram
$$ \xymatrix@R=50pt@C=50pt{ \Spec(A'/I') \ar[rr] \ar[dr]_{v'} & & \Spec(A/I) \ar[dl]^{v} \\
& X & }$$
commutes.
\end{itemize}
We say that a morphism $(A,I,v) \rightarrow (A', I', v')$ is a {\it flat cover} if the underlying ring homomorphism
$A \rightarrow A'$ is $p$-completely faithfully flat and $I' = IA'$ (in which case the system of divided powers on $I'$
is uniquely determined by the system of divided powers on $I$). The collection of flat coverings determines an extensive
Grothendieck topology on the category $\Crys(X/\Z_p)^{\op}$, which we will refer to as the {\it flat topology}. We refer to the category
$\Crys(X/\Z_p)$ (with the flat topology) as the {\it big crystalline site of $X$}. We will be primarily interested in the case where 
$X = \Spec(R)$ is affine; in this case, we abuse notation by writing $\Crys(R/\Z_p)$ for the category $\Crys(X/\Z_p)$.
\end{notation}

\begin{construction}[Crystalline Cohomology]\label{construction:crystalline-cochain-complex}
Let $X$ be an $\F_p$-scheme. The assignment $(A,I,v) \mapsto A$ determines a presheaf of commutative rings on
the category $\Crys(X/\Z_p)^{\op}$, which is a sheaf with respect to the flat topology of Notation \ref{notation:big-crystalline-site}.
We denote the derived global sections of this sheaf by $\RGamma_{\crys}(X/\Z_p)$, which we regard as a commutative
algebra object of the $\infty$-category $\widehat{\calD}(\Z_p)$. Phrased differently, we let $\RGamma_{\crys}(X/\Z_p)$
denote the homotopy inverse limit $\varprojlim_{ (A,I,v) \in \Crys(X/\Z_p) } A$, formed in the $\infty$-category $\widehat{\calD}(\Z_p)$.
We will refer to $\RGamma_{\crys}(X/\Z_p)$ as the {\it crystalline cochain complex} of $X$. We denote the cohomology groups
of $\RGamma_{\crys}(X/\Z_p)$ by $\mathrm{H}^{\ast}_{\crys}(X/\Z_p)$, and refer to them as the {\it crystalline cohomology groups of $X$}.
In the special case where $X = \Spec(R)$ is an affine $\F_p$-scheme, we denote the crystalline cochain complex $\RGamma_{\crys}(X/\Z_p)$ by
$\RGamma_{\crys}(R/\Z_p)$, and its cohomology groups by $\mathrm{H}^{\ast}_{\crys}(R/\Z_p)$.
\end{construction}

\begin{notation}\label{notation:crystalline-augmentation}
Let $X$ be an $\F_p$-scheme. For every open immersion $j: \Spec(R) \hookrightarrow X$, the triple $(R, (0), j)$ can be regarded
as an object of the big crystalline site $\Crys(X/\Z_p)$, so that we have an evaluation map
$$ \RGamma_{\crys}(X/\Z_p) = \varprojlim_{ (A,I,v) \in \Crys(X/\Z_p) } A \rightarrow R.$$
Passing the the inverse limit over affine open subsets of $X$, we obtain a map $\epsilon_{\crys}: \RGamma_{\crys}(X/\Z_p) \rightarrow
\RGamma(X, \calO_X)$, which we will refer to as the {\it crystalline augmentation morphism}. In the special case where $X = \Spec(R)$ is
an affine $\F_p$-scheme, we will denote $\epsilon_{\crys}$ by $\epsilon_{\crys}: \RGamma_{\crys}(R/\Z_p) \rightarrow R$.
\end{notation}

The crystalline cochain functor $\RGamma_{\crys}(\bullet/\Z_p)$ enjoys the following descent property:

\begin{proposition}\label{proposition:descend-to-semiperfect}
Let $P = \F_p[ \{ x_s \} ]$ be a polynomial algebra over $\F_p$ (possibly on infinitely many generators), let $P^{0} = \F_p[ \{ x_s^{1/p^{\infty}} \} ]$ be the perfection of
$P$, and let $P^{\bullet}$ be the cosimplicial $P$-algebra given by the tensor powers of $P^{0}$ relative to $P$. For every $P$-algebra $R$, the tautological map
$$ \RGamma_{\crys}(R / \Z_p) \rightarrow \Tot( \RGamma_{\crys}( R \otimes_{P} P^{\bullet} / \Z_p ) )$$
is an isomorphism in the $\infty$-category $\widehat{\calD}(\Z_p)$.
\end{proposition}

\begin{proof}
Let $\Crys(R/\Z_p)$ be the big crystalline site of $R$. Unwinding the definitions, we see that the totalization
$\Tot( \RGamma_{\crys}( R \otimes_{P} P^{\bullet} / \Z_p ) )$ can be identified with the homotopy inverse limit
$\varprojlim_{ (A,I,v) \in \calC } A$ (formed in the $\infty$-category $\widehat{\calD}(\Z_p)$), where $\calC \subseteq \Crys(R/\Z_p)$
is the full subcategory spanned by those objects $(A,I,v)$ for which the composite map $P \rightarrow R \xrightarrow{v} A/I$
factors through $P^{0}$. Since the functor $(A,I,v) \mapsto A$ is a $\widehat{\calD}(\Z_p)$-valued sheaf with respect to the
flat topology on $\Crys(R/\Z_p)^{\op}$, it will suffice to show that the subcategory $\calC^{\op} \subseteq \Crys(R/\Z_p)^{\op}$ is a
covering sieve with respect to the flat topology. Fix an object $(A,I,v) \in \Crys(R/\Z_p)$, and let $\overline{f}$ denote the composite map
$\F_p[ \{ x_s \} ] =  P \rightarrow R \xrightarrow{v} A/I$.
Then $\overline{f}$ can be lifted to a ring homomorphism $f: \Z_p[ \{ x_s \} ] \rightarrow A$. Let $A'$ denote the $p$-completed tensor product
$\Z_p[ \{ x_s^{1/p^{\infty}} \} ] \widehat{\otimes}_{ \Z_p[ \{ x_s \} ] } A$.
Then $A'$ is $p$-completely faithfully flat over $A$, so there is a unique system of divided powers on the ideal $I' = IA'$ which are compatible with those on $I$.
Let $v'$ denote the composite map $R \xrightarrow{v} A/I \rightarrow A'/I'$. Then the tautological map $(A,I,v) \rightarrow (A', I', v')$ is a flat covering,
and $(A', I', v')$ belongs to the subcategory $\Crys(R/\Z_p)$.
\end{proof}

Note that, for any commutative $\F_p$-algebra $R$, we can choose a polynomial ring $P = \F_p[ \{ x_s \} ]$ and a surjection $P \twoheadrightarrow R$.
In this case, each term of the cosimplicial commutative ring $R \otimes_{P} P^{\bullet}$ appearing in Proposition \ref{proposition:descend-to-semiperfect}
is a semiperfect $\F_p$-algebra. In the semiperfect case, the crystalline cochain complex admits a very concrete description.

\begin{construction}\label{construction:Acrys}
Let $R$ be an $\F_p$-algebra which is {\it semiperfect}: that is, every element of $R$ is a $p$th power. Let $R^{\flat}$ denote the inverse limit perfection
$$ \varprojlim( \cdots \rightarrow R \xrightarrow{ x \mapsto x^{p} } R \xrightarrow{ x \mapsto x^{p} } R \xrightarrow{ x \mapsto x^{p} } R ),$$
and let $I^{+}$ denote the kernel of the map $W( R^{\flat} ) \twoheadrightarrow R^{\flat} \twoheadrightarrow R$. We let
$( A_{\crys}(R), I_{\crys}(R) )$ denote the $p$-adically separated and complete divided power envelope of $( W(R^{\flat}), I^{+} )$ relative to $(\Z_p, (p) )$.
\end{construction}

\begin{remark}\label{remark:senty}
Let $R$ be a semiperfect $\F_p$-algebra and let $I^{+} \subseteq W(R^{\flat} )$ be as in Construction \ref{construction:Acrys}.
If $(B,J)$ is any $p$-adically separated and complete divided power algebra over $(\Z_p, (p))$, then we have canonical bijections
$$ \xymatrix@R=50pt@C=30pt{ \{ \textnormal{Divided power morphisms $(A_{\crys}(R), I_{\crys}(R)) \rightarrow (B,J)$} \} \ar[d]^{\sim} \\
\{ \textnormal{Morphisms of pairs $(W(R^{\flat}), I^{+}) \rightarrow (B,J)$} \} \ar[d]^{\sim} \\
\{ \textnormal{Ring homomorphisms $R \rightarrow B/J$} \}. }$$
Consequently, if we let $v$ denote the composite isomorphism $$R \simeq W(R^{\flat}) / I^{+} \simeq A_{\crys}(R) / I_{\crys}(R),$$ then
the triple $( A_{\crys}(R), I_{\crys}(R), v)$ is an initial object of the big crystalline site $\Crys(R/\Z_p)$. It follows
that the crystalline cochain complex $\RGamma_{\crys}( R/ \Z_p)$ can be identified with $A_{\crys}(R)$, regarded as an abelian group
concentrated in cohomological degree zero. Under this identification, the crystalline augmentation $\epsilon_{\crys}$ of Notation \ref{notation:crystalline-augmentation}
corresponds to the composition $A_{\crys}(R) \twoheadrightarrow A_{\crys}(R) / I_{\crys}(R) \xrightarrow{ v^{-1} } R$.
\end{remark}

\begin{remark}\label{remark:xi-0}
Let $R$ be a semiperfect $\F_p$-algebra and let $I$ denote the kernel of the quotient map $R^{\flat} \twoheadrightarrow R$. 
Then the quotient ring $A_{\crys}(R) / p A_{\crys}(R)$ can be identified with the divided power envelope of $(R^{\flat}, I)$
relative to $\F_p$. In particular, for each element $x \in I$, the image of $x^{p}$ vanishes in $A_{\crys}(R) / p A_{\crys}(R)$.
It follows that the composite map
$$ R^{\flat} \xrightarrow{ \varphi_{ R^{\flat} } } R^{\flat} \rightarrow A_{\crys}(R) / p A_{\crys}(R)$$
annihilates the ideal $I$, and therefore factors through a homomorphism
$\xi_0: R \rightarrow A_{\crys}(R) / p A_{\crys}(R)$. By construction, the composite map
$$ R \xrightarrow{\xi_0} A_{\crys}(R) / p A_{\crys}(R) \twoheadrightarrow R$$
is equal to the Frobenius map $\varphi_{R}: R \twoheadrightarrow R$.
\end{remark}

\begin{construction}\label{construction:conjugate-filtration}
Let $R$ be a semiperfect $\F_p$-algebra, let $I$ denote the kernel of the surjection $R^{\flat} \twoheadrightarrow R$,
and regard the quotient ring $A_{\crys}(R) / p A_{\crys}(R)$ as an $R$-algebra via the homomorphism $\xi_0$ of Remark \ref{remark:xi-0}.

For every integer $m > 0$, we let $\gamma_m: I \rightarrow A_{\crys}(R) / p A_{\crys}(R)$ denote the map determined by
the system of divided powers on the ideal $I_{\crys} / pA_{\crys}(R) \subseteq A_{\crys}(R) / p A_{\crys}(R)$.
For each $n \geq 0$, we let 
$$\Fil_{n}^{\conj}( A_{\crys}(R) / p A_{\crys}(R) ) \subseteq A_{\crys}(R) / p A_{\crys}(R)$$
denote the $R$-submodule of $A_{\crys}(R) / p A_{\crys}(R)$ generated by elements of the form
$$\gamma_{pm_1}(x_1) \cdots \gamma_{pm_k}(x_k),$$ 
where $x_1, x_2, \cdots, x_k$ are elements of $I$ and $m_1, m_2, \cdots, m_k$ are nonnegative integers satisfying
$m_1 + \cdots + m_k \leq n$. This construction determines an increasing filtration
$$ \Fil^{\conj}_0( A_{\crys}(R) / p A_{\crys}(R) ) \subseteq \Fil^{\conj}_1( A_{\crys}(R) / p A_{\crys}(R) ) \subseteq \Fil^{\conj}_2( A_{\crys}(R) / p A_{\crys}(R) ) \subseteq \cdots$$
which we will refer to as the {\it conjugate filtration}. 

Note that for every pair of elements 
$$y \in \Fil^{\conj}_m( A_{\crys}(R) / p A_{\crys}(R) ) \quad \quad z \in \Fil^{\conj}_n( A_{\crys}(R) / p A_{\crys}(R) ),$$ the product $y \cdot z$
belongs to $\Fil^{\conj}_{m+n}(  A_{\crys}(R) / p A_{\crys}(R) )$. It follows that the associated graded module
$$ \gr^{\conj}_{\ast}( A_{\crys}(R) / p A_{\crys}(R) ) =  \Fil^{\conj}_{\ast}( A_{\crys}(R) / p A_{\crys}(R) ) /  \Fil^{\conj}_{\ast-1}( A_{\crys}(R) / p A_{\crys}(R) )$$
inherits the structure of a graded $R$-algebra.
\end{construction}

\begin{proposition}\label{proposition:describe-gr-conjugate}
Let $R$ be a semiperfect $\F_p$-algebra, let $I$ denote the kernel of the map $R^{\flat} \twoheadrightarrow R$, and let
$\Gamma_{R}^{\ast}( I / I^2 )$ denote the divided power algebra on the $R$-module $I/I^2$. Then there is a unique homomorphism of graded rings
$$ \xi: \Gamma_{R}^{\ast}( I / I^2 ) \rightarrow \gr^{\conj}_{\ast}( A_{\crys}(R) / p A_{\crys}(R) )$$
with the following properties:
\begin{itemize}
\item[$(1)$] In degree zero, $\xi$ coincides with the ring homomorphism 
$$\xi_0: R \twoheadrightarrow \Fil_{0}^{\conj}( A_{\crys}(R) / p A_{\crys}(R) ) \subseteq A_{\crys}(R) / p A_{\crys}(R)$$
of Remark \ref{remark:xi-0}.

\item[$(2)$] For each element $x \in I$ having image $\overline{x} \in I/I^2$, the image
$$\xi( \gamma_n( \overline{x}) ) \in \gr_{n}^{\conj}( A_{\crys}(R) / p A_{\crys}(R) )$$ is represented
by the element $\gamma_{pn}(x) \in \Fil_{n}^{\conj}( A_{\crys}(R) / p A_{\crys}(R) )$.
\end{itemize}
Moreover, the homomorphism $\xi$ is surjective.
\end{proposition}

\begin{proof}
We will show that there exists a unique homomorphism of graded rings $\widetilde{\xi}: \Gamma^{\ast}_{ R^{\flat} }(I) \rightarrow \gr^{\conj}_{\ast}( A_{\crys}(R) / p A_{\crys}(R) )$
satisfying the following analogues of $(1)$ and $(2)$:
\begin{itemize}
\item[$(1')$] In degree zero, $\widetilde{\xi}$ agrees with the map
$$ R^{\flat} \xrightarrow{ \varphi_{ R^{\flat} }} R^{\flat} \rightarrow A_{\crys}(R) / p A_{\crys}(R)$$
(whose image is the $0$th step of the conjugate filtration).

\item[$(2')$] For each element $x \in I$, the image
$\widetilde{\xi}( \gamma_n(x) ) \in \gr_{n}^{\conj}( A_{\crys}(R) / p A_{\crys}(R) )$ is represented
by the element $\gamma_{pn}(x) \in \Fil_{n}^{\conj}( A_{\crys}(R) / p A_{\crys}(R) )$.
\end{itemize}
The homomorphism $\xi$ of Proposition \ref{proposition:describe-gr-conjugate} can then be obtained from $\widetilde{\xi}$ by extending
scalars along the quotient map $R^{\flat} \twoheadrightarrow R$ (using Remark \ref{remark:xi-0}). It follows immediately from the definitions that
if the homomorphism $\widetilde{\xi}$ exists, then it is uniquely determined by requirements $(1')$ and $(2')$ and is surjective
(from which it follows that $\xi$ is also surjective). To prove existence, we recall that the divided power algebra
$\Gamma^{\ast}_{ R^{\flat} }(I)$ can be described as the $R^{\flat}$-algebra generated by symbols $\{ \gamma_{n}(x) \}_{x \in I, n \geq 0}$, which are subject to the relations
$$ \gamma_{n}(x+y) = \sum_{i=0}^{n} \gamma_i(x) \gamma_{n-i}(y) \quad \quad \gamma_n( \lambda x) = \lambda^{n} \gamma_n(x)$$
$$ \gamma_{0}(x) = 1 \quad \quad \gamma_{m}(x) \gamma_{n}(x) = \frac{ (m+n)!}{m! \cdot n!} \gamma_{m+n}(x).$$
We are therefore reduced to proving the following:
\begin{itemize}
\item For every pair of elements $x,y \in I$ and every integer $n \geq 0$, we have
$$ \gamma_{pn}(x+y) \equiv \sum_{i=0}^{n} \gamma_{pi}(x) \gamma_{pn-pi}(y) \pmod{ \Fil^{\conj}_{n-1}( A_{\crys} / p A_{\crys} ) }.$$
Writing $\gamma_{pn}(x+y)$ as a sum $\sum_{j + k = pn} \gamma_j(x) \gamma_{k}(y)$, we are reduced to showing that
if $j$ and $k$ are nonnegative integers which are {\em not} divisible by $p$ satisfying $j+k = pn$, then 
$\gamma_{j}(x) \gamma_{k}(y)$ belongs to $\Fil^{\conj}_{n-1}( A_{\crys} / p A_{\crys} )$. Writing
$j = p j_0 + r$ and $k = pk_0 + (p-r)$ with $0< r < p$, we observe that 
$\gamma_{j}(x) \gamma_{k}(y)$ is a unit multiple of $x^{r} y^{p-r} \gamma_{pj_0}(x) \gamma_{pk_0}(y)$,
and therefore belongs to $\Fil^{\conj}_{n-1}( A_{\crys} / p A_{\crys} )$ by virtue of the identity $j_0 + k_0 = n-1$. 

\item For every element $x \in I$ and every element $\lambda \in R^{\flat}$, we have
$$\gamma_{pn}( \lambda x) = \lambda^{pn} \gamma_{pn}(x) = \varphi_{R^{\flat}}(\lambda)^{n} \gamma_{pn}(x)$$
in $A_{\crys}(R) / p A_{\crys}(R)$

\item For every element $x \in I$, we have $\gamma_{0}(x) = 1$ in $A_{\crys}(R) / p A_{\crys}(R)$.

\item For every element $x \in I$ and every pair of nonnegative integers $m$ and $n$, we have
$$ \gamma_{pm}(x) \gamma_{pn}(x) = \frac{ (pm+pn)!}{(pm)! \cdot (pn)!} \gamma_{p(m+n)}(x)
= \frac{ (m+n)!}{m! \cdot n!} \gamma_{p(m+n)}(x)$$
in $A_{\crys}(R) / p A_{\crys}(R)$ by virtue of the congruence
$\frac{ (pm+pn)!}{(pm)! \cdot (pn)!} \equiv  \frac{ (m+n)!}{m! \cdot n!} \pmod{ p}$.
\end{itemize}
\end{proof}

\newpage

\bibliography{bibliography}{}
\bibliographystyle{plain}

\end{document}